\documentclass[11pt]{amsbook}
\usepackage{amsfonts, amsmath, amssymb, amsthm, stmaryrd, color, enumerate, relsize}
\usepackage[pdftex]{graphicx}
\usepackage{mathrsfs,array}
\usepackage{xy}
\usepackage{hyperref}
\usepackage{tikz-cd}
\input xy
\xyoption{all}
\usepackage{scalerel, stackengine}
\stackMath
\newcommand\widecheck[1]{%
\savestack{\tmpbox}{\stretchto{%
  \scaleto{%
    \scalerel*[\widthof{\ensuremath{#1}}]{\kern-.6pt\bigwedge\kern-.6pt}%
    {\rule[-\textheight/2]{1ex}{\textheight}}%WIDTH-LIMITED BIG WEDGE
  }{\textheight}% 
}{0.5ex}}%
\stackon[1pt]{#1}{\scalebox{-1}{\tmpbox}}%
}
\DeclareMathOperator{\Spa}{Spa}
\DeclareMathOperator{\Set}{Set}

\DeclareMathOperator{\Spf}{Spf}

\DeclareMathOperator{\Spd}{Spd}

\DeclareMathOperator{\Proj}{Proj}

\DeclareMathOperator{\GL}{GL}
\DeclareMathOperator{\PGL}{PGL}
\DeclareMathOperator{\Fun}{Fun}
\DeclareMathOperator{\SL}{SL}
\DeclareMathOperator{\Gal}{Gal}
\DeclareMathOperator{\Hom}{Hom}
\DeclareMathOperator{\End}{End}
\DeclareMathOperator{\Ext}{Ext}
\DeclareMathOperator{\sHom}{\!{\mathscr{H}\! om}}

\DeclareMathOperator{\Aut}{Aut}
\DeclareMathOperator{\Map}{Map}
\DeclareMathOperator{\Mod}{Mod}
\DeclareMathOperator{\Rep}{Rep}

\DeclareMathOperator{\Frob}{Frob}
\DeclareMathOperator{\Spec}{Spec}
\DeclareMathOperator{\Bun}{Bun}
\DeclareMathOperator{\FilBun}{FilBun}
\DeclareMathOperator{\GrBun}{GrBun}
\DeclareMathOperator{\Sht}{Sht}
\DeclareMathOperator{\Perf}{Perf}

\DeclareMathOperator{\op}{op}
\newcommand{\cont}{\mathrm{cont}}
\DeclareMathOperator{\Lie}{Lie}

\DeclareMathOperator{\Coh}{Coh}
\DeclareMathOperator{\Gr}{Gr}
\DeclareMathOperator{\BC}{\mathcal{B}\mathcal{C}}
\DeclareMathOperator{\coker}{coker}
\DeclareMathOperator{\Fil}{Fil}

\DeclareMathOperator{\Isoc}{Isoc}
\DeclareMathOperator{\Pic}{Pic}

\DeclareMathOperator{\dimtrg}{dim.trg}

\DeclareMathOperator{\Sat}{Sat}
\DeclareMathOperator{\Perv}{Perv}
\newcommand{\dotimes}{\otimes^{\mathbb L}}
\newcommand{\soliddotimes}{\stackrel{\mathsmaller{\blacksquare}}{\otimes}\!\!{}^\mathbb L}
\newcommand{\soliddotimesLambda}{\stackrel{\mathsmaller{\blacksquare}}{\otimes}\!\!{}^\mathbb L_\Lambda}
\def\OO{\mathcal{O}}
\newcommand{\et}{{\mathrm{\acute{e}t}}}

\newcommand{\proet}{{\mathrm{pro\acute{e}t}}}
\newcommand{\qproet}{{\mathrm{qpro\acute{e}t}}}

\newcommand{\Pro}{\mathrm{Pro}}
\renewcommand{\op}{\mathrm{op}}
\renewcommand{\int}{\mathrm{int}}

\newcommand{\lc}{\mathrm{lc}}

\newcommand{\perf}{\mathrm{perf}}

\newcommand{\sm}{\mathrm{sm}}
\newcommand{\id}{\mathrm{id}}
\newcommand{\tw}{\mathrm{tw}}

\newcommand{\ULA}{\mathrm{ULA}}
\DeclareMathOperator{\Ind}{Ind}
\DeclareMathOperator{\colim}{colim}
\DeclareMathOperator{\LocSys}{LocSys}
\DeclareMathOperator{\Hglob}{Hck}
\DeclareMathOperator{\Hloc}{{\mathcal Hck}}
\DeclareMathOperator{\bigast}{{\mathlarger{\mathlarger{\ast}}}}
\DeclareMathOperator{\bigboxtimes}{{\mathlarger{\mathlarger{\boxtimes}}}}

\def\Fl{{\mathcal{F}}\!\ell}

\newcommand{\dR}{\mathrm{dR}}

\DeclareMathOperator{\LT}{LT}
\renewcommand*{\diamond}{\diamondsuit}
\renewcommand*{\hat}{\widehat}
\renewcommand*{\tilde}{\widetilde}
\newcommand{\lis}{\mathrm{lis}}
\newcommand{\solid}{{\mathsmaller{\blacksquare}}}

\def\Z{\mathbb{Z}}
\def\Q{\mathbb{Q}}

\def\Fq{\mathbb{F}_q}
\def\O{\mathcal{O}}
\def\E{\mathcal{E}}
\def\iso{\xrightarrow{\sim}}
\def\Div{\mathrm{Div}}
\def\lto{\longrightarrow}
\renewcommand*{\phi}{\varphi}

\newcommand{\powerseries}[1]{[\![ #1 ]\!]}
\newcommand{\laurentseries}[1]{(\!( #1 )\!)}
\def\Qlb{\overline{\mathbb{Q}}_\ell}

\numberwithin{equation}{section}
\newtheorem{theorem}{Theorem}
\numberwithin{theorem}{section}
\newtheorem{lemma}[theorem]{Lemma}
\newtheorem{corollary}[theorem]{Corollary}
\newtheorem{proposition}[theorem]{Proposition}
\newtheorem{conj}[theorem]{Conjecture}
\newtheorem{problem}[theorem]{Problem}

\newtheorem{definition}[theorem]{Definition}
\newtheorem{defprop}[theorem]{Definition/Proposition}
\newtheorem{question}[theorem]{Question}
\newtheorem{hypothesis}[theorem]{Hypothesis}

\theoremstyle{definition}
\newtheorem{remark}[theorem]{Remark}
\newtheorem{example}[theorem]{Example}

\newtheorem{convention}[theorem]{Convention}

\newtheorem{warning}[theorem]{Warning}

\newenvironment{altenumerate}
   {\begin{list}
      {\textup{(\theenumi)} }
      {\usecounter{enumi}
       \setlength{\labelwidth}{0pt}
       \setlength{\labelsep}{2pt}
       \setlength{\leftmargin}{0pt}
       \setlength{\itemsep}{\the\smallskipamount}
       \renewcommand{\theenumi}{\roman{enumi}}
      }}
   {\end{list}}
\newenvironment{altitemize}
   {\begin{list}
      {$\bullet$ }
      {\setlength{\labelwidth}{0pt}
       \setlength{\labelsep}{2pt}
       \setlength{\leftmargin}{0pt}
       \setlength{\itemsep}{\the\smallskipamount}
      }}
   {\end{list}}

\renewcommand{\thechapter}{\Roman{chapter}}

\date{\today}

\title{Geometrization of the local Langlands correspondence}
\author{Laurent Fargues and Peter Scholze}
\begin{document}

\begin{abstract} Following the idea of \cite{FarguesGeom}, we develop the foundations of the geometric Langlands program on the Fargues--Fontaine curve. In particular, we define a category of $\ell$-adic sheaves on the stack $\Bun_G$ of $G$-bundles on the Fargues--Fontaine curve, prove a geometric Satake equivalence over the Fargues--Fontaine curve, and study the stack of $L$-parameters. As applications, we prove finiteness results for the cohomology of local Shimura varieties and general moduli spaces of local shtukas, and define $L$-parameters associated with irreducible smooth representations of $G(E)$, a map from the spectral Bernstein center to the Bernstein center, and the spectral action of the category of perfect complexes on the stack of $L$-parameters on the category of $\ell$-adic sheaves on $\Bun_G$.
\end{abstract}

\maketitle

\tableofcontents

\renewcommand{\thesection}{\thechapter.\arabic{section}}

\chapter{Introduction}

\section{The local Langlands correspondence}

The local Langlands correspondence aims at a description of the irreducible smooth representations $\pi$ of $G(E)$, for a reductive group $G$ over a local field $E$. Until further notice, we will simplify our life by assuming that $G$ is split; the main text treats the case of general reductive $G$, requiring only minor changes.

The case where $E$ is archimedean, i.e.~$E=\mathbb R$ or $E=\mathbb C$, is the subject of Langlands' classical work \cite{LanglandsArchimedean}. Based on the work of Harish-Chandra, cf.~e.g.~\cite{HarishChandra}, Langlands associates to each $\pi$ an $L$-parameter, that is a continuous homomorphism
\[
\varphi_\pi: W_E\to \hat{G}(\mathbb C)
\]
where $W_E$ is the Weil group of $E=\mathbb R,\mathbb C$ (given by $W_{\mathbb C}=\mathbb C^\times$ resp.~a nonsplit extension $1\to W_{\mathbb C}\to W_{\mathbb R}\to \mathrm{Gal}(\mathbb C/\mathbb R)\to 1$), and $\hat{G}$ is the Langlands dual group. This is the split reductive group over $\mathbb Z$ whose root datum is dual to the root datum of $G$. The map $\pi\mapsto \varphi_\pi$ has finite fibres, and a lot of work has been done on making the fibres, the so-called $L$-packets, explicit. If $G=\GL_n$, the map $\pi\mapsto \varphi_\pi$ is essentially a bijection.

Throughout this paper, we assume that $E$ is nonarchimedean, of residue characteristic $p>0$, with residue field $\mathbb F_q$. Langlands has conjectured that one can still naturally associate an $L$-parameter
\[
\varphi_\pi: W_E\to \hat{G}(\mathbb C)
\]
to any irreducible smooth representation $\pi$ of $G(E)$. In the nonarchimedean case, $W_E$ is the dense subgroup of the absolute Galois group $\mathrm{Gal}(\overline{E}|E)$, given by the preimage of $\mathbb Z\subset \mathrm{Gal}(\overline{\mathbb F}_q|\Fq)$ generated by the Frobenius $x\mapsto x^q$. This raises the question where such a parameter should come from. In particular,
\begin{altenumerate}
\item[{\rm (1)}] How does the Weil group $W_E$ relate to the representation theory of $G(E)$?
\item[{\rm (2)}] How does the Langlands dual group $\hat{G}$ arise?
\end{altenumerate}

The goal of this paper is to give a natural construction of a parameter $\varphi_\pi$ (only depending on a choice of isomorphism $\mathbb C\cong \overline{\mathbb Q}_\ell$), and in the process answer questions (1) and (2).

\section{The big picture}

In algebraic geometry, to any ring $A$ corresponds a space $\Spec A$. The starting point of our investigations is a careful reflection on the space $\Spec E$ associated with $E$.\footnote{Needless to say, the following presentation bears no relation to the historical developments of the ideas, which as usual followed a far more circuitous route. We will discuss some of our original motivation in Section~\ref{sec:introorigin} below.} Note that the group $G(E)$ is the automorphism group of the trivial $G$-torsor over $\Spec E$, while the Weil group of $E$ is essentially the absolute Galois group of $E$, that is the (\'etale) fundamental group of $\Spec E$. Thus, $G(E)$ relates to coherent information (especially $G$-torsors) on $\Spec E$, while $W_E$ relates to \'etale information on $\Spec E$. Moreover, the perspective of $G$-torsors is a good one: Namely, for general groups $G$ there can be nontrivial $G$-torsors $\mathcal E$ on $\Spec E$, whose automorphism groups are then the so-called pure inner forms of Vogan \cite{VoganPureInner}. Vogan realized that from the perspective of the local Langlands correspondence, and in particular the parametrization of the fibres of $\pi\mapsto \varphi_\pi$, it is profitable to consider all pure inner forms together; in particular, he was able to formulate a precise form of the local Langlands conjecture (taking into account the fibres of $\pi\mapsto \varphi_\pi$) for pure inner forms of (quasi)split groups. All pure inner forms together arise by looking at the groupoid of all $G$-bundles on $\Spec E$: This is given by
\[
[\ast/G](\Spec E) = \bigsqcup_{[\alpha]\in H^1_\et(\Spec E,G)} [\ast/G_\alpha(E)],
\]
where $H^1_\et(\Spec E,G)$ is the set of $G$-torsors on $\Spec E$ up to isomorphism, and $G_\alpha$ the corresponding pure inner form of $G$. Also, we already note that representations of $G(E)$ are equivalent to sheaves on $[\ast/G(E)]$ (this is a tautology if $G(E)$ were a discrete group; in the present context of smooth representations, it is also true for the correct notion of ``sheaf''), and hence sheaves on
\[
[\ast/G](\Spec E) = \bigsqcup_{[\alpha]\in H^1_\et(\Spec E,G)} [\ast/G_\alpha(E)],
\]
are equivalent to tuples $(\pi_\alpha)_{[\alpha]}$ of representations of $G_\alpha(E)$.\footnote{The point of replacing $[\ast/G(E)]$ by $[\ast/G](\Spec E)$ was also stressed by Bernstein.}

Looking at the \'etale side of the correspondence, we observe that the local Langlands correspondence makes the Weil group $W_E$ of $E$ appear, not its absolute Galois group $\mathrm{Gal}(\overline{E}|E)$. Recall that $W_E\subset \mathrm{Gal}(\overline{E}|E)$ is the dense subgroup given as the preimage of the inclusion $\mathbb Z\subset \mathrm{Gal}(\overline{\mathbb F}_q|\Fq)\cong \hat{\mathbb Z}$, where $\mathrm{Gal}(\overline{\mathbb F}_q|\Fq)$ is generated by its Frobenius morphism $x\mapsto x^q$. On the level of geometry, this change corresponds to replacing a scheme $X$ over $\mathbb F_q$ with the (formal) quotient $X_{\overline{\mathbb F}_q}/\mathrm{Frob}$.

In the function field case $E=\Fq\laurentseries{t}$, we are thus led to replace $\Spec E$ by $\Spec \breve{E}/\phi^\Z$ where $\breve{E} = \overline{\mathbb F}_q\laurentseries{t}$. We can actually proceed similarly in general, taking $\breve{E}$ to be the completion of the maximal unramified extension of $E$. For a natural definition of $\pi_1$, one then has $\pi_1(\Spec(\breve{E})/\phi^\Z)=W_E$ --- for example, $\Spec \overline{\breve{E}}\to \Spec(\breve{E})/\phi^\Z$ is a $W_E$-torsor, where $\overline{\breve{E}}$ is a separable closure.

Let us analyze what this replacement entails on the other side of the correspondence: Looking at the coherent theory of $\Spec \breve{E}/\phi^\Z$, one is led to study $\breve{E}$-vector spaces $V$ equipped with $\phi$-linear automorphisms $\sigma$. This is known as the category of isocrystals $\mathrm{Isoc}_E$. The category of isocrystals is much richer than the category of $E$-vector spaces, which it contains fully faithfully. Namely, by the Dieudonn\'e--Manin classification, the category $\mathrm{Isoc}_E$ is semisimple, with one simple object $V_\lambda$ for each rational number $\lambda\in \mathbb Q$. The endomorphism algebra of $V_\lambda$ is given by the central simple $E$-algebra $D_\lambda$ of Brauer invariant $\overline{\lambda}\in \mathbb Q/\mathbb Z$. Thus, there is an equivalence of categories
\[
\mathrm{Isoc}_E\cong \bigoplus_{\lambda\in \mathbb Q} \mathrm{Vect}_{D_\lambda}\otimes V_\lambda.
\]
Here, if one writes $\lambda=\frac sr$ with coprime integers $s,r$, $r>0$, then $V_\lambda\cong \breve{E}^r$ is of rank $r$ with $\sigma$ given by the matrix
\[
\left(\begin{array}{ccccc} 0 & 1 & 0 & \ldots & \\
& 0 & 1 & 0 & \\
& & \ddots & \ddots & \ddots \\
0 & & & 0 & 1\\
\pi^s & 0 & & & 0
\end{array}\right)
\]
where $\pi\in E$ is a uniformizer.

Above, we were considering $G$-torsors on $\Spec E$, thus we should now look at $G$-torsors in $\Isoc_E$. These are known as $G$-isocrystals and have been extensively studied by Kottwitz \cite{KottwitzIsocrystals}, \cite{KottwitzIsocrystals2}. Their study has originally been motivated by the relation of isocrystals to $p$-divisible groups and accordingly a relation of $G$-isocrystals to the special fibre of Shimura varieties (parametrizing abelian varieties with $G$-structure, and thus $p$-divisible groups with $G$-structure). Traditionally, the set of $G$-isocrystals is denoted $B(E,G)$, and for $b\in B(E,G)$ we write $\mathcal E_b$ for the corresponding $G$-isocrystal. In particular, Kottwitz has isolated the class of basic $G$-isocrystals; for $G=\GL_n$, a $G$-isocrystal is just a rank $n$ isocrystal, and it is basic precisely when it has only one slope $\lambda$. There is an injection $H^1_\et(\Spec E,G)\hookrightarrow B(E,G)$ as any $G$-torsor on $\Spec E$ ``pulls back'' to a $G$-torsor in $\Isoc_E$; the image lands in $B(E,G)_{\mathrm{basic}}$. For any $b\in B(E,G)_{\mathrm{basic}}$, the automorphism group of $\mathcal E_b$ is an inner form $G_b$ of $G$; the set of such inner forms of $G$ is known as the extended pure inner forms of $G$. Note that for $G=\GL_n$, there are no nontrivial pure inner forms of $G$, but all inner forms of $G$ are extended pure inner forms, precisely by the occurence of all central simple $E$-algebras as $H_\lambda$ for some slope $\lambda$. More generally, if the center of $G$ is connected, then all inner forms of $G$ can be realized as extended pure inner forms. Kaletha, \cite{KalethaLocalLanglands}, has extended Vogan's results on pure inner forms to extended pure inner forms, giving a precise form of the local Langlands correspondence (describing the fibres of $\pi\mapsto \varphi_\pi$) for all extended pure inner forms and thereby showing that $G$-isocrystals are profitable from a purely representation-theoretic point of view. We will actually argue below that it is best to include $G_b$ for all $b\in B(E,G)$, not only the basic $b$; the resulting automorphism groups $G_b$ are then inner forms of Levi subgroups of $G$. Thus, we are led to consider the groupoid of $G$-torsors in $\mathrm{Isoc}_E$,
\[
G\text-\mathrm{Isoc}\cong \bigsqcup_{[b]\in B(E,G)} [\ast/G_b(E)].
\]
Sheaves on this are then tuples of representations $(\pi_b)_{[b]\in B(E,G)}$ of $G_b(E)$. The local Langlands conjecture, including its expected functorial behaviour with respect to passage to inner forms and Levi subgroups, then still predicts that for any irreducible sheaf $\mathcal F$ --- necessarily given by an irreducible representation $\pi_b$ of $G_b(E)$ for some $b\in B(E,G)$ --- one can associate an $L$-parameter $\varphi_{\mathcal F}: W_E\to \hat{G}(\mathbb C)$.

To go further, we need to bring geometry into the picture: Indeed, it will be via geometry that (sheaves on the groupoid of) $G$-torsors on $\Spec \breve{E}/\phi^\Z$ will be related to the fundamental group $W_E$ of $\Spec \breve{E}/\phi^\Z$. The key idea is to study a moduli stack of $G$-torsors on $\Spec \breve{E}/\phi^\Z$.

There are several ways to try to define such a moduli stack. The most naive may be the following. The category $\mathrm{Isoc}_E$ is an $E$-linear category. We may thus, for any $E$-algebra $A$, consider $G$-torsors in $\mathrm{Isoc}_E\otimes_E A$. The resulting moduli stack will then actually be
\[
\bigsqcup_{b\in B(E,G)} [\ast/G_b],
\]
an Artin stack over $E$, given by a disjoint union of classifying stacks for the algebraic groups $G_b$. This perspective is actually instrumental in defining the $G_b$ as algebraic groups. However, it is not helpful for the goal of further geometrizing the situation. Namely, sheaves on $[\ast/G_b]$ are representations of the algebraic group $G_b$, while we are interested in representations of the locally profinite group $G_b(E)$.

A better perspective is to treat the choice of $\overline{\mathbb F}_q$ as auxiliary, and replace it by a general $\overline{\mathbb F}_q$-algebra $R$. In the equal characteristic case, we can then replace $\breve{E} = \overline{\mathbb F}_q\laurentseries{t}$ with $R\laurentseries{t}$. This carries a Frobenius $\phi=\phi_R$ acting on $R$. To pass to the quotient $\Spec R\laurentseries{t}/\phi^\Z$, we need to assume that the Frobenius of $R$ is an automorphism, i.e.~that $R$ is perfect. (The restriction to perfect $R$ will become even more critical in the mixed characteristic case. For the purpose of considering $\ell$-adic sheaves, the passage to perfect rings is inconsequential, as \'etale sheaves on a scheme $X$ and on its perfection are naturally equivalent.) We are thus led to the moduli stack on perfect $\overline{\mathbb F}_q$-algebras
\[
G\text-\mathcal{I}\mathrm{soc}: \{\mathrm{perfect}\ \overline{\mathbb F}_q\text-\mathrm{algebras}\}\to \{\mathrm{groupoids}\}: R\mapsto \{G\text-\mathrm{torsors\ on\ } \Spec R\laurentseries{t}/\phi^\Z\}.
\]
These are also known as families of $G$-isocrystals over the perfect scheme $\Spec R$. (Note the curly $\mathcal I$ in $G\text-\mathcal{I}\mathrm{soc}$, to distinguish it from the groupoid $G\text-\mathrm{Isoc}$.)

This definition can be extended to the case of mixed characteristic. Indeed, if $R$ is a perfect $\Fq$-algebra, the analogue of $R\powerseries{t}$ is the unique $\pi$-adically complete flat $\mathcal O_E$-algebra $\tilde{R}$ with $\tilde{R}/\pi=R$; explicitly, $\tilde{R}=W_{\mathcal O_E}(R)=W(R)\otimes_{W(\Fq)} \mathcal O_E$ in terms of the $p$-typical Witt vectors $W(R)$ or the ramified Witt vectors $W_{\mathcal O_E}(R)$. Thus, if $E$ is of mixed characteristic, we define
\[
G\text-\mathcal{I}\mathrm{soc}: \{\mathrm{perfect}\ \overline{\mathbb F}_q\text-\mathrm{algebras}\}\to \{\mathrm{groupoids}\}: R\mapsto \{G\text-\mathrm{torsors\ on\ } \Spec (W_{\mathcal O_E}(R)[\tfrac 1\pi])/\phi^\Z\}.
\]

We will not use the stack $G\text-\mathcal{I}\mathrm{soc}$ in this paper. However, it has been highlighted recently among others implicitly by Genestier--V.~Lafforgue, \cite{GenestierLafforgue}, and explicitly by Gaitsgory, \cite[Section 4.2]{GaitsgoryLGsigma}, and Zhu, \cite{ZhuConjectures}, and one can hope that the results of this paper have a parallel expression in terms of $G\text-\mathcal{I}\mathrm{soc}$, so let us analyze it further in this introduction. It is often defined in the following slightly different form. Namely, v-locally on $R$, any $G$-torsor over $R\laurentseries{t}$ resp.~$W_{\mathcal O_E}(R)[\tfrac 1\pi]$ is trivial by a recent result of Ansch\"utz \cite{AnschuetzParahoric}. Choosing such a trivialization, a family of $G$-isocrystals is given by some element of $LG(R)$, where we define the loop group
\[
LG(R) = G(R\laurentseries{t})\ (\mathrm{resp.\ } LG(R) = G(W_{\mathcal O_E}(R)[\tfrac 1\pi])).
\]
Changing the trivialization of the $G$-torsor amounts to $\sigma$-conjugation on $LG$, so as v-stacks
\[
G\text-\mathcal{I}\mathrm{soc} = LG/_{\mathrm{Ad,\sigma}} LG
\]
is the quotient of $LG$ under $\sigma$-conjugation by $LG$.

The stack $G\text-\mathcal{I}\mathrm{soc}$ can be analyzed. More precisely, we have the following result.\footnote{This result seems to be well-known to experts, but we are not aware of a full reference. For the v-descent (even arc-descent), see \cite[Lemma 5.9]{IvanovArcDescent}. The stratification is essentially constructed in \cite{RapoportRichartz}; the local constancy of the Kottwitz map is proved in general in Corollary~\ref{cor:RapoportRichartz}. The identification of the strata in some cases is in \cite[Proposition 4.3.13]{CaraianiScholze}, and in general in \cite[Theorem 1.4]{HamacherKim}.}

\begin{theorem}\label{thm:GIsoc} The prestack $G\text-\mathcal{I}\mathrm{soc}$ is a stack for the v-topology on perfect $\overline{\mathbb F}_q$-algebras. It admits a stratification into locally closed substacks
\[
G\text-\mathcal{I}\mathrm{soc}^b\subset G\text-\mathcal{I}\mathrm{soc}
\]
for $b\in B(E,G)$, consisting of the locus where at each geometric point, the $G$-isocrystal is isomorphic to $\mathcal E_b$. Moreover, each stratum
\[
G\text-\mathcal{I}\mathrm{soc}^b\cong [\ast/G_b(E)]
\]
is a classifying stack for the locally profinite group $G_b(E)$.
\end{theorem}

The loop group $LG$ is an ind-(infinite dimensional perfect scheme), so the presentation
\[
G\text-\mathcal{I}\mathrm{soc}=LG/_{\mathrm{Ad},\sigma} LG
\]
is of extremely infinite nature. We expect that this is not an issue with the presentation, but that the stack itself has no good finiteness properties; in particular note that all strata appear to be of the same dimension $0$, while admitting nontrivial specialization relations. Xiao--Zhu (see \cite{XiaoZhu}, \cite{ZhuConjectures}) have nonetheless been able to define a category $D(G\text-\mathcal{I}\mathrm{soc},\overline{\mathbb Q}_\ell)$ of $\ell$-adic sheaves on $G\text-\mathcal{I}\mathrm{soc}$, admitting a semi-orthogonal decomposition into the various $D(G\text-\mathcal{I}\mathrm{soc}^b,\overline{\mathbb Q}_\ell)$. Each $D(G\text-\mathcal{I}\mathrm{soc}^b,\overline{\mathbb Q}_\ell)\cong D([\ast/G_b(E)],\overline{\mathbb Q}_\ell)$ is equivalent to the derived category of the category of smooth representations of $G_b(E)$ (on $\overline{\mathbb Q}_\ell$-vector spaces). Here, as usual, we have to fix an auxiliary prime $\ell\neq p$ and an isomorphism $\mathbb C\cong \overline{\mathbb Q}_\ell$.

At this point we have defined a stack $G\text-\mathcal{I}\mathrm{soc}$, with a closed immersion
\[
i: [\ast/G(E)]\cong G\text-\mathcal{I}\mathrm{soc}^1\subset G\text-\mathcal{I}\mathrm{soc},
\]
thus realizing a fully faithful embedding
\[
i_\ast: D(G(E),\overline{\mathbb Q}_\ell)\hookrightarrow D(G\text-\mathcal{I}\mathrm{soc},\overline{\mathbb Q}_\ell)
\]
of the derived category of smooth representations of $G(E)$ into the derived category of $\overline{\mathbb Q}_\ell$-sheaves on $G\text-\mathcal{I}\mathrm{soc}$. It is in this way that we ``geometrize the representation theory of $G(E)$''.

The key additional structure that we need are the Hecke operators: These will simultaneously make the Weil group $W_E$ (i.e.~$\pi_1(\Spec \breve{E}/\phi^\Z)$) and, by a careful study, also the Langlands dual group $\hat{G}$ appear. Recall that Hecke operators are related to modifications of $G$-torsors, and are parametrized by a point $x$ of the curve where the modification happens, and the type of the modification at $x$ (which can be combinatorially encoded in terms of a cocharacter of $G$ --- this eventually leads to the appearance of $\hat{G}$). Often, the effect of Hecke operators is locally constant for varying $x$. In that case, letting $x$ vary amounts to an action of $\pi_1(X)$, where $X$ is the relevant curve; thus, the curve should now be $\Spec \breve{E}/\phi^\Z$.

Thus, if we want to define Hecke operators, we need to be able to consider modifications of $G$-isocrystals. These modifications ought to happen at a section of $\Spec R\laurentseries{t}\to \Spec R$ (resp.~a non-existent map $\Spec(W(R)\otimes_{W(\Fq)} E)\to \Spec R$). Unfortunately, the map $R\to R\laurentseries{t}$ does not admit any sections. In fact, we would certainly want to consider continuous sections; such continuous sections would then be in bijection with maps $\overline{\mathbb F}_q\laurentseries{t}=\breve{E}\to R$. In other words, in agreement with the motivation from the previous paragraph, the relevant curve should be $\Spec \breve{E}$, or really $\Spec \breve{E}$ modulo Frobenius --- so we can naturally hope to get actions of $\pi_1(\Spec \breve{E}/\phi^\Z)$ by the above recipe.

However, in order for this picture to be realized we need to be in a situation where we have continuous maps $\Fq\laurentseries{t}\to R$. In other words, we can only hope for sections if we put ourselves into a setting where $R$ is itself some kind of Banach ring.

This finally brings us to the setting considered in this paper. Namely, we replace the category of perfect $\overline{\mathbb F}_q$-schemes with the category of perfectoid spaces $\Perf = \Perf_{\overline{\mathbb F}_q}$ over $\overline{\mathbb F}_q$. Locally any $S\in \Perf$ is of the form $S=\Spa(R,R^+)$ where $R$ is a perfectoid Tate $\overline{\mathbb F}_q$-algebra: This means that $R$ is a perfect topological algebra that admits a topologically nilpotent unit $\varpi\in R$ (called a pseudouniformizer) making it a Banach algebra over $\overline{\mathbb F}_q\laurentseries{\varpi}$. Moreover, $R^+\subset R$ is an open and integrally closed subring of powerbounded elements. Often $R^+=R^\circ$ is the subring of powerbounded elements, and we consequently use the abbreviation $\Spa R = \Spa (R,R^\circ)$. The geometric (rank $1$) points of $S$ are given by $\Spa C$ for complete algebraically closed nonarchimedean fields $C$, and as usual understanding geometric points is a key first step. We refer to \cite{Berkeley} for an introduction to adic and perfectoid spaces. 

For any $S=\Spa(R,R^+)$, we need to define the analogue of $\Spec R\laurentseries{t}/\phi^\Z$, taking the topology of $R$ into account. Note that for discrete $R'$, we have
\[
\Spa R'\laurentseries{t} = \Spa R'\times_{\Spa \Fq} \Spa \Fq\laurentseries{t},
\]
and we are always free to replace $\Spec R'\laurentseries{t}/\phi^\Z$ by $\Spa R'\laurentseries{t}/\phi^\Z$ as they have the same category of vector bundles. This suggests that the analogue of $\Spec R'\laurentseries{t}$ is
\[
\Spa(R,R^+)\times_{\Spa \Fq} \Spa \Fq\laurentseries{t} = \mathbb D^\ast_{\Spa(R,R^+)},
\]
a punctured open unit disc over $\Spa(R,R^+)$, with coordinate $t$. Note that
\[
\Spa(R,R^+)\times_{\Spa \Fq} \Spa \Fq\laurentseries{t}\subset \Spa R^+\times_{\Spa \Fq} \Spa \Fq\powerseries{t} = \Spa R^+\powerseries{t}
\]
is the locus where $t$ and $\varpi\in R^+$ are invertible, where $\varpi$ is a topologically nilpotent unit of $R$. The latter definition can be extended to mixed characteristic: We let
\[
\Spa(R,R^+)\dot\times_{\Spa \Fq} \Spa E\subset \Spa R^+\dot\times_{\Spa\Fq} \Spa \mathcal O_E := \Spa W_{\mathcal O_E}(R^+)
\]
be the open subset where $\pi$ and $[\varpi]\in W_{\mathcal O_E}(R^+)$ are invertible. This space is independent of the choice of $\varpi$ as for any other such $\varpi'$, one has $\varpi|\varpi'^n$ and $\varpi'|\varpi^n$ for some $n\geq 1$, and then the same happens for their Teichm\"uller representatives. We note that the symbol $\dot\times$ is purely symbolic: There is of course no map of adic spaces $\Spa E\to \Spa \Fq$ along which a fibre product could be taken.

\begin{definition} The Fargues--Fontaine curve (for the local field $E$, over $S=\Spa(R,R^+)\in \Perf$) is the adic space over $E$ defined by
\[
X_S = X_{S,E} = \left(\Spa(R,R^+)\times_{\Spa \Fq} \Spa \Fq\laurentseries{t}\right)/\phi^\Z,
\]
respectively
\[
X_S = X_{S,E} = \left(\Spa(R,R^+)\dot\times_{\Spa \Fq} \Spa E\right)/\phi^\Z,
\]
where the Frobenius $\phi$ acts on $(R,R^+)$.
\end{definition}

A novel feature, compared to the discussion of $G$-isocrystals, is that the action of $\phi$ is free and totally discontinuous, so the quotient by $\phi$ is well-defined in the category of adic spaces. In fact, on $Y_S=\Spa(R,R^+)\dot\times_{\Spa \Fq} \Spa E\subset \Spa W_{\mathcal O_E}(R^+)$ one can compare the absolute values of $\pi$ and $[\varpi]$. As both are topologically nilpotent units, the ratio
\[
\mathrm{rad}=\log(|[\varpi]|)/\log(|\pi|): |Y_S|\to (0,\infty)
\]
gives a well-defined continuous map. The Frobenius on $|Y_S|$ multiplies $\mathrm{rad}$ by $q$, proving that the action is free and totally discontinuous.

We note that in the function field case $E=\Fq\laurentseries{t}$, the space
\[
Y_S= S\times_{\Spa \Fq} \Spa \Fq\laurentseries{t} = \mathbb D^\ast_S
\]
is precisely a punctured open unit disc over $S$. In this picture, the radius function measures the distance to the origin: Close to the origin, the radius map is close to $0$, while close to the boundary of the open unit disc it is close to $\infty$. The quotient by $\phi$ is however not an adic space over $S$ anymore, as $\phi$ acts on $S$. Thus,
\[
X_S=Y_S/\phi^\Z = \mathbb D^\ast_S/\phi^\Z
\]
is locally an adic space of finite type over $S$, but not globally so. This space, for $S=\Spa C$ a geometric point, has been first studied by Hartl--Pink \cite{HartlPink}.

If $S=\Spa C$ is a geometric point but $E$ is general, this curve (or rather a closely related schematic version) has been extensively studied by Fargues--Fontaine \cite{FarguesFontaine}, where it was shown that it plays a central role in $p$-adic Hodge theory. From the perspective of adic spaces, it has been studied by Kedlaya--Liu \cite{KedlayaLiu1}. In particular, in this case where $S$ is a point, $X_S$ is indeed a curve: It is a strongly noetherian adic space whose connected affinoid subsets are spectra of principal ideal domains. In particular, in this situation there is a well-behaved notion of ``classical points'', referring to those points that locally correspond to maximal ideals. These can be classified. In the equal characteristic case, the description of
\[
Y_S = \mathbb D^\ast_S = S\times_{\Spa \Fq} \Spa \Fq\laurentseries{t}
\]
shows that the closed points are in bijection with maps $S\to \Spa \Fq\laurentseries{t}$ up to Frobenius; where now one has to take the quotient under $t\mapsto t^q$. In mixed characteristic, the situation is more subtle, and brings us to the tilting construction for perfectoid spaces.

\begin{proposition} If $E$ is of mixed characteristic and $S=\Spa C$ is a geometric point, the classical points of $X_C$ are in bijection with untilts $C^\sharp|E$ of $C$, up to the action of Frobenius.
\end{proposition}

Here, we recall that for any complete algebraically closed field $C'|E$, or more generally any perfectoid Tate ring $R$, one can form the tilt
\[
R^\flat = \varprojlim_{x\mapsto x^p} R,
\]
where the addition is defined on the ring of integral elements in terms of the bijection $R^{\flat +}=\varprojlim_{x\mapsto x^p} R^+\cong \varprojlim_{x\mapsto x^p} R^+/\pi$, where now $x\mapsto x^p$ is compatible with addition on $R^+/\pi$. Then $R^\flat$ is a perfectoid Tate algebra of characteristic $p$. Geometrically, sending $\Spa(R,R^+)$ to $\Spa(R^\flat,R^{\flat +})$ defines a tilting functor on perfectoid spaces $T\mapsto T^\flat$, preserving the underlying topological space and the \'etale site, cf.~\cite{Berkeley}.

One sees that the classical points of $X_S$, for $S=\Spa C$ a geometric point, are in bijection with untilts $S^\sharp$ of $S$ together with a map $S\to \Spa E$, modulo the action of Frobenius. Recall from \cite{Berkeley} that for any adic space $Z$ over $W(\overline{\mathbb F}_q)$, one defines a functor
\[
Z^\diamond: \Perf\to \mathrm{Sets}: S\mapsto \{S^\sharp, f: S^\sharp\to Z\}
\]
sending a perfectoid space $S$ over $\overline{\mathbb F}_q$ to pairs $S^\sharp$ of an untilt of $S$, and a map $S^\sharp\to Z$. If $Z$ is an analytic adic space, then $Z^\diamond$ is a diamond, that is a quotient of a perfectoid space by a pro-\'etale equivalence relation. Then the classical points of $X_S$ are in bijection with the $S$-valued points of the diamond
\[
(\Spa \breve{E})^\diamond/\phi^\Z.
\]
More generally, for any $S\in \Perf$, maps $S\to (\Spa \breve{E})^\diamond/\phi^\Z$ are in bijection with degree $1$ Cartier divisors $D_S\subset X_S$, so we define
\[
\mathrm{Div}^1 = (\Spa \breve{E})^\diamond/\phi^\Z.
\]
We warn the reader the action of the Frobenius here is a geometric Frobenius. In particular, it only exists on $(\Spa \breve{E})^\diamond$, not on $\Spa \breve{E}$, in case $E$ is of mixed characteristic. However, one still has $\pi_1(\mathrm{Div}^1) = W_E$.

This ends our long stream of thoughts on the geometry of $\Spec E$: We have arrived at the Fargues--Fontaine curve, in its various incarnations. To orient the reader, we recall them here:

\begin{altenumerate}
\item For any complete algebraically closed nonarchimedean field $C|\overline{\mathbb F}_q$, the curve $X_C=X_{C,E}$, a strongly noetherian adic space over $E$, locally the adic spectrum of a principal ideal domain. One can also construct a schematic version $X_C^{\mathrm{alg}}$, with the same classical points and the same category of vector bundles. The classical points are in bijection with untilts $C^\sharp|E$ of $C$, up to Frobenius.
\item More generally, for any perfectoid space $S\in \Perf$, the ``family of curves'' $X_S$, again an adic space over $E$, but no longer strongly noetherian. If $S$ is affinoid, there is a schematic version $X_S^{\mathrm{alg}}$, with the same category of vector bundles. 
\item The ``mirror curve'' $\mathrm{Div}^1 = (\Spa \breve{E})^\diamond/\phi^\Z$, which is only a diamond. For any $S\in \Perf$, this parametrizes ``degree $1$ Cartier divisors on $X_S$''.
\end{altenumerate}

A peculiar phenomenon here is that there is no ``absolute curve'' of which all the others are the base change. Another peculiar feature is that the space of degree $1$ Cartier divisors is not the curve itself.

Again, it is time to study $G$-torsors. This leads to the following definition.

\begin{definition} Let
\[
\Bun_G: \Perf\to \{\mathrm{groupoids}\}: S\mapsto \{G\text-\mathrm{torsors}\ \mathrm{on}\ X_S\}
\]
be the moduli stack of $G$-torsors on the Fargues--Fontaine curve.
\end{definition}

\begin{remark} Let us stress here that while ``the Fargues--Fontaine curve'' is not really a well-defined notion, ``the moduli stack of $G$-torsors on the Fargues--Fontaine curve'' is.
\end{remark}

As $X_S$ maps towards $\Spa \breve{E}/\phi^\Z$, there is a natural pullback functor $G\text-\mathrm{Isoc}\to \Bun_G(S)$. The following result is in most cases due to Fargues \cite{FarguesGBun}, completed by Ansch\"utz, \cite{AnschuetzGBun}.

\begin{theorem}\label{thm:introFarguesBG} If $S=\Spa C$ is a geometric point, the map
\[
B(G)\to \Bun_G(S)/\cong
\]
is a bijection. In particular, any vector bundle on $X_S$ is a direct sum of vector bundles $\mathcal O_{X_S}(\lambda)$ associated to $D_{-\lambda}$, $\lambda\in \mathbb Q$.

Under this bijection, $b\in B(G)$ is basic if and only if the corresponding $G$-torsor $\mathcal E_b$ on $X_S$ is semistable in the sense of Atiyah--Bott \cite{AtiyahBott}.
\end{theorem}

However, it is no longer true that the automorphism groups are the same. On the level of the stack, we have the following result.

\begin{theorem} The prestack $\Bun_G$ is a v-stack. It admits a stratification into locally closed substacks
\[
i^b: \Bun_G^b\subset \Bun_G
\]
for $b\in B(G)$ consisting of the locus where at each geometric point, the $G$-torsor is isomorphic to $\mathcal E_b$. Moreover, each stratum
\[
\Bun_G^b\cong [\ast/\widetilde{G}_b]
\]
is a classifying space for a group $\widetilde{G}_b$ that is an extension of the locally profinite group $G_b(E)$ by a ``unipotent group diamond''.

The semistable locus $\Bun_G^{\mathrm{ss}}\subset \Bun_G$ is an open substack, and
\[
\Bun_G^{\mathrm{ss}} = \bigsqcup_{b\in B(G)_{\mathrm{basic}}} [\ast/G_b(E)].
\]
\end{theorem}

\begin{remark} The theorem looks formally extremely similar to Theorem~\ref{thm:GIsoc}. However, there is a critical difference, namely the closure relations are reversed: For $\Bun_G$, the inclusion of $\Bun_G^b$ for $b\in B(G)$ basic is an open immersion while it was a closed immersion in Theorem~\ref{thm:GIsoc}. Note that basic $b\in B(G)$ correspond to semistable $G$-bundles, and one would indeed expect the semistable locus to be an open substack. Generally, $\Bun_G$ behaves much like the stack of $G$-bundles on the projective line.
\end{remark}

\begin{remark} We define a notion of Artin stacks in this perfectoid setting, and indeed $\Bun_G$ is an Artin stack; we refer to Section~\ref{sec:introBunG0} for a more detailed description of our geometric results on $\Bun_G$. This shows that $\Bun_G$ has much better finiteness properties than $G\text-\mathcal{I}\mathrm{soc}$, even if it is defined on more exotic test objects.
\end{remark}

We can define a derived category of $\ell$-adic sheaves
\[
D(\Bun_G,\overline{\mathbb Q}_\ell)
\]
on $\Bun_G$. This admits a semi-orthogonal decomposition into all $D(\Bun_G^b,\overline{\mathbb Q}_\ell)$, and
\[
D(\Bun_G^b,\overline{\mathbb Q}_\ell)\cong D([\ast/G_b(E)],\overline{\mathbb Q}_\ell)\cong D(G_b(E),\overline{\mathbb Q}_\ell)
\]
is equivalent to the derived category of smooth $G_b(E)$-representations.

\begin{remark} It is reasonable to expect that this category is equivalent to the category $D(G\text-\mathcal{I}\mathrm{soc},\overline{\mathbb Q}_\ell)$ defined by Xiao--Zhu. However, we do not pursue this comparison here.
\end{remark}

Finally, we can define the Hecke stack that will bring all key players together. Consider the global Hecke stack $\Hglob_G$ parametrizing pairs $(\mathcal E,\mathcal E')$ of $G$-bundles on $X_S$, together with a map $S\to \mathrm{Div}^1$ giving rise to a degree $1$ Cartier divisor $D_S\subset X_S$, and an isomorphism
\[
f:\mathcal E|_{X_S\setminus D_S}\cong \mathcal E'|_{X_S\setminus D_S}
\]
that is meromorphic along $D_S$. This gives a correspondence
\[
\Bun_G\xleftarrow{h_1}\Hglob_G\xrightarrow{h_2} \Bun_G\times \Div^1.
\]
To define the Hecke operators, we need to bound the modification, i.e.~bound the poles of $f$ along $D_S$. This is described by the local Hecke stack $\Hloc_G$, parametrizing pairs of $G$-torsors on the completion of $X_S$ along $D_S$, together with an isomorphism away from $D_S$; thus, there is a natural map $\Hglob_G\to \Hloc_G$ from the global to the local Hecke stack. Geometrically, $\Hloc_G$ admits a Schubert stratification in terms of the conjugacy classes of cocharacters of $G$; in particular, there are closed Schubert cells $\Hloc_{G,\leq \mu}$ for each conjugacy class $\mu: \mathbb G_m\to G$. By pullback, this defines a correspondence
\[
\Bun_G\xleftarrow{h_{1,\leq\mu}}\Hglob_{G,\leq\mu}\xrightarrow{h_{2,\leq\mu}} \Bun_G\times \Div^1
\]
where now $h_{1,\leq\mu}$ and $h_{2,\leq\mu}$ are proper. One can then consider Hecke operators
\[
Rh_{2,\leq\mu,\ast} h_{1,\leq\mu}^\ast: D(\Bun_G,\Lambda)\to D(\Bun_G\times \Div^1,\Lambda).
\]

The following theorem ensures that Hecke operators are necessarily locally constant as one varies the point of $\Div^1$, and hence give rise to representations of $\pi_1(\Div^1) = W_E$. In the following, we are somewhat cavalier about the precise definition of $D(-,\overline{\mathbb Q}_\ell)$ employed, and the notion of $W_E$-equivariant objects: The fine print is addressed in the main text.

\begin{theorem} Pullback along the map $\Div^1\to [\ast/W_E]$ induces an equivalence
\[
D(\Bun_G\times \Div^1,\overline{\mathbb Q}_\ell)\cong D(\Bun_G\times [\ast/W_E],\overline{\mathbb Q}_\ell)\cong D(\Bun_G,\overline{\mathbb Q}_\ell)^{BW_E}.
\]
\end{theorem}

Thus, Hecke operators produce $W_E$-equivariant objects in $D(\Bun_G,\overline{\mathbb Q}_\ell)$, making the Weil group appear naturally.

One also wants to understand how Hecke operators compose. This naturally leads to the study of $D(\Hloc_G,\overline{\mathbb Q}_\ell)$ as a monoidal category, under convolution. Here, we have the geometric Satake equivalence. In the setting of usual smooth projective curves (over $\mathbb C$), this was established in the papers of Lusztig \cite{LusztigKatoLusztig}, Ginzburg \cite{GinzburgSatake} and Mirkovi\'c--Vilonen \cite{MirkovicVilonen}. The theorem below is a first approximation; we will actually prove a more precise version with $\mathbb Z_\ell$-coefficients, describing all perverse sheaves on $\Hloc_G$, and applying to the Beilinson--Drinfeld Grassmannians in the spirit of Gaitsgory's paper \cite{GaitsgoryBeilinsonDrinfeldSatake}.

\begin{theorem}\label{thm:geometricsatakeintro} There is a natural monoidal functor from $\Rep \hat{G}$ to $D(\Hloc_G,\overline{\mathbb Q}_\ell)$.
\end{theorem}

\begin{remark} Our proof of Theorem~\ref{thm:geometricsatakeintro} follows the strategy of Mirkovi\'c--Vilonen's proof, and in particular defines a natural symmetric monoidal structure on the category of perverse sheaves by using the fusion product. This requires one to work over several copies of the base curve, and let the points collide. It is a priori very surprising that this can be done in mixed characteristic, as it requires a space like $\Spa \mathbb Q_p\dot\times_{\Spa \mathbb F_p} \Spa \mathbb Q_p$. Spaces of this type do however exist as diamonds, and this was one of the main innovations of \cite{Berkeley}.
\end{remark}

\begin{remark} Using a degeneration of the local Hecke stack, which is essentially the $B_{\mathrm{dR}}^+$-affine Grassmannian of \cite{Berkeley}, to the Witt vector affine Grassmannian, Theorem~\ref{thm:geometricsatakeintro} gives a new proof of Zhu's geometric Satake equivalence for the Witt vector affine Grassmannian \cite{ZhuWitt}. In fact, we even prove a version with $\mathbb Z_\ell$-coefficients, thus also recovering the result of Yu \cite{YuWittSatake}.
\end{remark}

\begin{remark} Regarding the formalism of $\ell$-adic sheaves, we warn the reader that we are cheating slightly in the formulation of Theorem~\ref{thm:geometricsatakeintro}; the definition of $D(\Bun_G,\overline{\mathbb Q}_\ell)$ implicit above is not the same as the one implicit in Theorem~\ref{thm:geometricsatakeintro}. With torsion coefficients, the problem would disappear, and in any case the problems are essentially of technical nature.
\end{remark}

Thus, this also makes the Langlands dual group $\hat{G}$ appear naturally. For any representation $V$ of $\hat{G}$, we get a Hecke operator
\[
T_V: D(\Bun_G,\overline{\mathbb Q}_\ell)\to D(\Bun_G,\overline{\mathbb Q}_\ell)^{BW_E}.
\]
Moreover, the Hecke operators commute and
\[
T_{V\otimes W}\cong T_V\circ T_W|_{\Delta(W_E)}
\]
where we note that $T_V\circ T_W$ naturally takes values in $W_E\times W_E$-equivariant objects; the restriction on the right means the restriction to the action of the diagonal copy $\Delta(W_E)\subset W_E\times W_E$.

At this point, the representation theory of $G(E)$ (which sits fully faithfully in $D(\Bun_G,\overline{\mathbb Q}_\ell)$), the Weil group $W_E$, and the dual group $\hat{G}$, all interact with each other naturally. It turns out that this categorical structure is precisely what is needed to construct $L$-parameters for (Schur-)irreducible objects $A\in D(\Bun_G,\overline{\mathbb Q}_\ell)$, and in particular for irreducible smooth representations of $G(E)$. We will discuss the construction of $L$-parameters below in Section~\ref{sec:introLparameters}.

We note that the whole situation is exactly parallel to the Betti geometric Langlands situation considered by Nadler--Yun \cite{NadlerYun}, and indeed the whole strategy can be described as ``the geometric Langlands program on the Fargues--Fontaine curve''. It is curious that our quest was to understand the local Langlands correspondence in an arithmetic setting, for potentially very ramified representations, and eventually we solved it by relating it to the global Langlands correspondence in a geometric setting, in the everywhere unramified setting.

In the rest of this introduction, we give a more detailed overview of various aspects of this picture:

\begin{altenumerate}
\item The Fargues--Fontaine curve (Section~\ref{sec:introFF});
\item The geometry of the stack $\Bun_G$ (Section~\ref{sec:introBunG0});
\item The derived category of $\ell$-adic sheaves on $\Bun_G$ (Section~\ref{sec:introBunGladic});
\item The geometric Satake equivalence (Section~\ref{sec:introGeomSatake});
\item Finiteness of the cohomology of Rapoport--Zink spaces, local Shimura varieties, and more general moduli spaces of shtukas (Section~\ref{sec:introRZspaces});
\item The stack of $L$-parameters (Section~\ref{sec:introstackLparameters});
\item The construction of $L$-parameters (Section~\ref{sec:introLparameters});
\item The spectral action (Section~\ref{sec:introspectral});
\item The origin of the ideas fleshed out in this paper (Section~\ref{sec:introorigin}).
\end{altenumerate}

These items largely mirror the chapters of this paper, and each chapter begins with a reprise of these introductions.

\section{The Fargues--Fontaine curve}\label{sec:introFF}

The Fargues--Fontaine curve has been studied extensively in the book of Fargues--Fontaine \cite{FarguesFontaine} and further results, especially in the relative situation, have been obtained by Kedlaya--Liu \cite{KedlayaLiu1}. In the first chapter, we reprove these foundational results, thereby also collecting and unifying certain results (proved often only for $E=\mathbb Q_p$).

The first results concern the Fargues--Fontaine curve $X_C=X_S$ when $S=\Spa C$ for some complete algebraically closed nonarchimedean field $C|\Fq$. We define a notion of classical points of $X_C$ in that case; they form a subset of $|X_C|$. The basic finiteness properties of $X_C$ are summarized in the following result.

\begin{theorem}\label{thm:introFFcurve0} The adic space $X_C$ is locally the adic spectrum $\Spa(B,B^+)$ where $B$ is a principal ideal domain; the classical points of $\Spa(B,B^+)\subset X_C$ are in bijection with the maximal ideals of $B$. For each classical point $x\in X_C$, the residue field of $x$ is an untilt $C^\sharp$ of $C$ over $E$, and this induces a bijection of the classical points of $X_C$ with untilts $C^\sharp$ of $C$ over $E$, taken up to the action of Frobenius.
\end{theorem}

In the equal characteristic case, Theorem~\ref{thm:introFFcurve0} is an immediate consequence of the presentation $X_C = \mathbb D^\ast_C/\phi^\Z$ and classical results in rigid-analytic geometry. In the $p$-adic case, we use tilting to reduce to the equal characteristic case. At one key turn, in order to understand Zariski closed subsets of $X_C$, we use the result that Zariski closed implies strongly Zariski closed \cite{BhattScholzePrism}. Using these ideas, we are able to give an essentially computation-free proof.

A key result is the classification of vector bundles.

\begin{theorem}\label{thm:introFFclassification} The functor from $\Isoc_E$ to vector bundles on $X_C$ induces a bijection on isomorphism classes. In particular, there is a unique stable vector bundle $\mathcal O_{X_C}(\lambda)$ of any slope $\lambda\in \mathbb Q$, and any vector bundle $\mathcal E$ can be written as a direct sum of stable bundles.
\end{theorem}

We give a new self-contained proof of Theorem~\ref{thm:introFFclassification}, making critical use of the v-descent results for vector bundles obtained in \cite{ECoD} and \cite{Berkeley}, and basic results on the geometry of Banach--Colmez spaces established here. The proof in the equal characteristic case by Hartl--Pink \cite{HartlPink} and the proof of Kedlaya in the $p$-adic case \cite{KedlayaLocalMonodromy} relied on heavy computations, while the proof of Fargues--Fontaine \cite{FarguesFontaine} relied on the description of the Lubin--Tate and Drinfeld moduli spaces of $\pi$-divisible $\mathcal O$-modules. Our proof is related to the arguments of Colmez in \cite{ColmezEspaces}.

Allowing general $S\in \Perf_{\overline{\mathbb F}_q}$, we define the moduli space of degree $1$ Cartier divisors as $\Div^1 = \Spd \breve{E}/\phi^\Z$. Given a map $S\to \Div^1$, one can define an associated closed Cartier divisor $D_S\subset X_S$; locally, this is given by an untilt $D_S=S^\sharp\subset X_S$ of $S$ over $E$, and this embeds $\Div^1$ into the space of closed Cartier divisors on $X_S$ (justifying the name). Another important result is the following ampleness result, cf.~\cite[Proposition 6.2.4]{KedlayaLiu1}, which implies that one can define an algebraic version of the curve, admitting the same theory of vector bundles.

\begin{theorem}\label{thm:introFFGAGA} Assume that $S\in \Perf$ is affinoid. For any vector bundle $\mathcal E$ on $X_S$, the twist $\mathcal E(n)$ is globally generated and has no higher cohomology for all $n\gg 0$. Defining the graded ring
\[
P=\bigoplus_{n\geq 0} H^0(X_S,\mathcal O_{X_S}(n))
\]
and the scheme $X_S^{\mathrm{alg}} = \Proj P$, there is a natural map of locally ringed spaces $X_S\to X_S^{\mathrm{alg}}$, pullback along which defines an equivalence of categories of vector bundles, preserving cohomology.

If $S=\Spa C$ for some complete algebraically closed nonarchimedean field $C$, then $X_C^{\mathrm{alg}}$ is a regular noetherian scheme of Krull dimension $1$, locally the spectrum of a principal ideal domain, and its closed points are in bijection with the classical points of $X_C$.
\end{theorem}

We also need to understand families of vector bundles, i.e.~vector bundles $\mathcal E$ on $X_S$ for general $S$. Here, the main result is the following.

\begin{theorem}\label{thm:introFFKL} Let $S\in \Perf$ and let $\mathcal E$ be a vector bundle on $X_S$. Then the function taking a point $s\in S$ to the Harder--Narasimhan polygon of $\mathcal E|_{X_s}$ defines a semicontinuous function on $S$. If it is constant, then $\mathcal E$ admits a global Harder--Narasimhan stratification, and pro-\'etale locally on $S$ one can find an isomorphism with a direct sum of $\mathcal O_{X_S}(\lambda)$'s.

In particular, if $\mathcal E$ is everywhere semistable of slope $0$, then $\mathcal E$ is pro-\'etale locally trivial, and the category of such $\mathcal E$ is equivalent to the category of pro-\'etale $\underline{E}$-local systems on $S$.
\end{theorem}

The key to proving Theorem~\ref{thm:introFFKL} is the construction of certain global sections of $\mathcal E$. To achieve this, we use v-descent techniques, and an analysis of the spaces of global sections of $\mathcal E$; these are known as Banach--Colmez spaces, and were first introduced (in slightly different terms) in \cite{ColmezEspaces}.

\begin{definition}\label{def:introBanachColmez} Let $\mathcal E$ be a vector bundle on $X_S$. The Banach--Colmez space $\BC(\mathcal E)$ associated with $\mathcal E$ is the locally spatial diamond over $S$ whose $T$-valued points, for $T\in \Perf_S$, are given by
\[
\BC(\mathcal E)(T) = H^0(X_T,\mathcal E|_{X_T}).
\]
Similarly, if $\mathcal E$ is everywhere of only negative Harder--Narasimhan slopes, the negative Banach--Colmez space $\BC(\mathcal E[1])$ is the locally spatial diamond over $S$ whose $T$-valued points are
\[
\BC(\mathcal E[1])(T) = H^1(X_T,\mathcal E|_{X_T}).
\]
\end{definition}

Implicit here is that this functor actually defines a locally spatial diamond. For this, we calculate some key examples of Banach--Colmez spaces. For example, if $\mathcal E=\mathcal O_{X_S}(\lambda)$ with $0<\lambda\leq [E:\mathbb Q_p]$ (resp.~all positive $\lambda$ if $E$ is of equal characteristic), then $\BC(\mathcal E)$ is representable by a perfectoid open unit disc (of dimension given by the numerator of $\lambda$). A special case of this is the identification of $\BC(\mathcal O_{X_S}(1))$ with the universal cover of a Lubin--Tate formal group law, yielding a very close relation between Lubin--Tate theory, and thus local class field theory, and the Fargues--Fontaine curve; see also \cite{FarguesClassFieldTheory}. On the other hand, for larger $\lambda$, or negative $\lambda$, Banach--Colmez spaces are more exotic objects; for example, the negative Banach--Colmez space
\[
\BC(\mathcal O_{X_C}(-1)[1])\cong (\mathbb A^1_{C^\sharp})^\diamond/\underline{E}
\]
is the quotient of the affine line by the translation action of $\underline{E}\subset \mathbb A^1_{C^\sharp}$. We remark that our proof of the classification theorem, Theorem~\ref{thm:introFFclassification}, ultimately relies on the negative result that $\BC(\mathcal O_{X_C}(-1)[1])$ is not representable by a perfectoid space!\footnote{Actually, we only know this for sure if $E$ is $p$-adic; in the function field case, we supply a small extra argument circumventing the issue.}

For the proof of Theorem~\ref{thm:introFFKL}, a key result is that projectivized Banach--Colmez spaces
\[
(\BC(\mathcal E)\setminus \{0\})/\underline{E^\times}
\]
are proper --- they are the relevant analogues of ``families of projective spaces over $S$''. In particular, their image in $S$ is a closed subset, and if the image is all of $S$, then we can find a nowhere vanishing section of $\mathcal E$ after a v-cover, as then the projectivized Banach--Colmez space is a v-cover of $S$. From here, Theorem~\ref{thm:introFFKL} follows easily.

\section{The geometry of $\mathrm{Bun}_G$}\label{sec:introBunG0}

Let us discuss the geometry of $\Bun_G$. Here, $G$ can be any reductive group over a nonarchimedean local field $E$, with residue field $\Fq$ of characteristic $p$. Recall that Kottwitz' set $B(G)=B(E,G)$ of $G$-isocrystals can be described combinatorially, by two discrete invariants. The first is the Newton point
\[
\nu: B(G)\to (X_\ast(T)_{\mathbb Q}^+)^\Gamma,
\]
where $T$ is the universal Cartan of $G$ and $\Gamma=\mathrm{Gal}(\overline{E}|E)$. More precisely, any $G$-isocrystal $\mathcal E$ defines a slope morphism $\mathbb D\to G_{\breve{E}}$ where $\mathbb D$ is the diagonalizable group with cocharacter group $\mathbb Q$; its definition reduces to the case of $\GL_n$, where it amounts to the slope decomposition of isocrystals. Isomorphisms of $G$-isocrystals lead to conjugate slope morphisms, and this defines the map $\nu$.

The other map is the Kottwitz invariant
\[
\kappa: B(G)\to \pi_1(G_{\overline{E}})_\Gamma.
\]
Its definition is indirect, starting from tori, passing to the case of $G$ with simply connected derived group, and finally to the general case by z-extensions. Then Kottwitz shows that
\[
(\nu,\kappa): B(G)\to (X_\ast(T)_{\mathbb Q}^+)^\Gamma\times \pi_1(G_{\overline{E}})_\Gamma
\]
is injective. Moreover, $\kappa$ induces a bijection between $B(G)_{\mathrm{basic}}$ and $\pi_1(G_{\overline{E}})_\Gamma$. The non-basic elements can be described in terms of Levi subgroups.

Using $\nu$ and $\kappa$, one can define a partial order on $B(G)$ by declaring $b\leq b'$ if $\kappa(b)=\kappa(b')$ and $\nu_b\leq \nu_{b'}$ with respect to the dominance order.

Up to sign, one can think of $\nu$, resp.~$\kappa$, as the Harder--Narasimhan polygon, resp.~first Chern class, of a $G$-bundle.

\begin{theorem}\label{thm:introBunG} The prestack $\Bun_G$ satisfies the following properties.
\begin{altenumerate}
\item[{\rm (i)}] The prestack $\Bun_G$ is a stack for the v-topology.
\item[{\rm (ii)}] The points $|\Bun_G|$ are naturally in bijection with Kottwitz' set $B(G)$ of $G$-isocrystals.
\item[{\rm (iii)}] The map
\[
\nu: |\Bun_G|\to B(G)\to (X_\ast(T)_{\mathbb Q}^+)^\Gamma
\]
is semicontinuous, and
\[
\kappa: |\Bun_G|\to B(G)\to \pi_1(G_{\overline{E}})_\Gamma
\]
is locally constant. Equivalently, the map $|\Bun_G|\to B(G)$ is continuous when $B(G)$ is equipped with the order topology.
\item[{\rm (iv)}] For any $b\in B(G)$, the corresponding subfunctor
\[
i^b: \Bun_G^b=\Bun_G\times_{|\Bun_G|} \{b\}\subset \Bun_G
\]
is locally closed, and isomorphic to $[\ast/\tilde{G}_b]$, where $\tilde{G}_b$ is a v-sheaf of groups such that $\tilde{G}_b\to \ast$ is representable in locally spatial diamonds with $\pi_0 \tilde{G}_b=G_b(E)$. The connected component $\tilde{G}_b^\circ\subset \tilde{G}_b$ of the identity is cohomologically smooth of dimension $\langle 2\rho,\nu_b\rangle$.
\item[{\rm (v)}] In particular, the semistable locus $\Bun_G^{\mathrm{ss}}\subset \Bun_G$ is open, and given by
\[
\Bun_G^{\mathrm{ss}}\cong \bigsqcup_{b\in B(G)_{\mathrm{basic}}} [\ast/\underline{G_b(E)}].
\]
\item[{\rm (vi)}] For any $b\in B(G)$, there is a map
\[
\pi_b: \mathcal M_b\to \Bun_G
\]
that is representable in locally spatial diamonds, partially proper and cohomologically smooth, where $\mathcal M_b$ parametrizes $G$-bundles $\mathcal E$ together with an increasing $\mathbb Q$-filtration whose associated graded is, at all geometric points, isomorphic to $\mathcal E_b$ with its slope grading. The v-stack $\mathcal M_b$ is representable in locally spatial diamonds, partially proper and cohomologically smooth over $[\ast/\underline{G_b(E)}]$.
\item[{\rm (vii)}] The v-stack $\Bun_G$ is a cohomologically smooth Artin stack of dimension $0$.
\end{altenumerate}
\end{theorem}

As examples, let us analyze the case of $\GL_1$ and $\GL_2$. For $\GL_1$, and general tori, everything is semistable, so
\[
\mathrm{Pic}:=\Bun_{\GL_1}\cong \bigsqcup_{\mathbb Z} [\ast/\underline{E^\times}].
\]
For $\GL_2$, the Kottwitz invariant gives a decomposition
\[
\Bun_{\GL_2} = \bigsqcup_{\alpha\in \tfrac 12 \mathbb Z} \Bun_{\GL_2}^\alpha.
\]
Each connected component has a unique semistable point, given by the basic element $b\in B(\GL_2)_{\mathrm{basic}}$ with $\kappa(b)=\alpha$. For $b\in B(\GL_2)_{\mathrm{basic}}\cong \tfrac 12\mathbb Z$, the corresponding group $G_b(E)$ is given by $\GL_2(E)$ when $b\in \mathbb Z$, and by $D^\times$ when $b\in \tfrac 12\mathbb Z\setminus \mathbb Z$, where $D|E$ is the quaternion algebra. 

The non-semistable points of $\Bun_{\GL_2}$ are given by extensions of line bundles, which are of the form $\mathcal O(i)\oplus \mathcal O(j)$ for some $i,j\in \mathbb Z$, with $2\alpha=i+j$. Let us understand the simplest degeneration inside $\Bun_{\GL_2}$, which is from $\mathcal O(\tfrac 12)$ to $\mathcal O\oplus \mathcal O(1)$. The individual strata here are
\[
[\ast/\underline{D^\times}],\ [\ast/\underline{\mathrm{Aut}}(\mathcal O\oplus \mathcal O(1))].
\]
Here
\[
\underline{\mathrm{Aut}}(\mathcal O\oplus \mathcal O(1)) = \left(\begin{array}{cc} \underline{E^\times} & \BC(\mathcal O(1))\\ 0 & \underline{E^\times}\end{array}\right).
\]
Here $\BC(\mathcal O(1))$ is representable by a perfectoid open unit disc $\Spd \overline{\mathbb F}_q\powerseries{t^{1/p^\infty}}$.

In this case, the local chart $\mathcal M_b$ for $\Bun_{\GL_2}$ parametrizes rank $2$ bundles $\mathcal E$ written as an extension
\[
0\to \mathcal L\to \mathcal E\to \mathcal L'\to 0
\]
such that at all geometric points, $\mathcal L\cong \mathcal O$ and $\mathcal L'\cong \mathcal O(1)$. Fixing such isomorphisms defines a $\underline{E^\times}\times \underline{E^\times}$-torsor
\[
\tilde{\mathcal M}_b\to \mathcal M_b
\]
with $\tilde{\mathcal M}_b = \BC(\mathcal O(-1)[1])$ a ``negative Banach--Colmez space''. This local chart shows that the local structure of $\Bun_G$ is closely related to the structure of negative Banach--Colmez spaces. It also shows that while the geometry of $\Bun_G$ is quite nonstandard, it is still fundamentally a finite-dimensional and ``smooth'' situation.

For general $G$, we still get a decomposition into connected components
\[
\Bun_G = \bigsqcup_{\alpha\in \pi_1(G)_\Gamma} \Bun_G^\alpha
\]
and each connected component $\Bun_G^\alpha$ admits a unique semistable point.

By a recent result of Viehmann \cite{ViehmannOberwolfach}, the map $|\Bun_G|\to B(G)$ is a homeomorphism. This had previously been proved for $G=\GL_n$ by Hansen \cite{HansenGLn} based on \cite{Hansenetal}; that argument was extended to some classical groups in unpublished work of Hamann.

Let us say some words about the proof of Theorem~\ref{thm:introBunG}. Part (i) has essentially been proved in \cite{Berkeley}, and part (ii) follows from the result of Fargues and Ansch\"utz, Theorem~\ref{thm:introFarguesBG}. In part (iii), the statement about $\nu$ reduces to $\GL_n$ by an argument of Rapoport--Richartz \cite{RapoportRichartz}, where it is Theorem~\ref{thm:introFFKL}. The statement about $\kappa$ requires more work, at least in the general case: If the derived group of $G$ is simply connected, one can reduce to tori, which are not hard to handle. In general, one approach is to use z-extensions $\tilde{G}\to G$ to reduce to the case of simply connected derived group. For this, one needs that $\Bun_{\tilde{G}}\to \Bun_G$ is a surjective map of v-stacks; we prove this using Beauville--Laszlo uniformization. Alternatively, one can use the abelianized Kottwitz set of Borovoi \cite{BorovoiAbelian}, which we prove to behave well relatively over a perfectoid space $S$. Part (iv) is a also consequence of Theorem~\ref{thm:introFFKL}. Part (v) is a consequence of parts (iii) and (iv). The key point is then part (vi), which will imply (vii) formally. (One can, and we do, also prove part (vii) directly. Part (vi) is however critical for the other results we prove below.) The properties of $\mathcal M_b$ itself are easy to establish --- the analysis for $\GL_2$ above generalizes easily to show that $\tilde{\mathcal M}_b$ is a successive extension of negative Banach--Colmez spaces. The key difficulty is to prove that
\[
\pi_b: \mathcal M_b\to \Bun_G
\]
is cohomologically smooth. Note that as we are working with perfectoid spaces, there are no tangent spaces, and we cannot hope to prove smoothness via deformation theory. To attack this problem, we nonetheless prove a general ``Jacobian criterion of cohomological smoothness''. The setup here is the following.

Let $S$ be a perfectoid space, and let $Z\to X_S$ be a smooth map of (sousperfectoid) adic spaces; this means that $Z$ is an adic space that is locally \'etale over a finite-dimensional ball over $X_S$. In this situation, we can define a v-sheaf $\mathcal M_Z\to S$ parametrizing sections of $Z\to X_S$, i.e.~the $S'$-valued points, for $S'/S$ a perfectoid space, are given by the maps $s: X_{S'}\to Z$ lifting $X_{S'}\to X_S$. For each such section, we get the vector bundle $s^\ast T_{Z/X_S}$ on $S'$, where $T_{Z/X_S}$ is the tangent bundle. Naively, deformations of $S'\to \mathcal M_Z$, i.e.~of $X_{S'}\to Z$ over $X_{S'}\to X_S$, should correspond to global sections $H^0(X_{S'},s^\ast T_{Z/X_S})$, and obstructions to $H^1(X_{S'},s^\ast T_{Z/X_S})$. If $s^\ast T_{Z/X_S}$ has everywhere only positive Harder--Narasimhan slopes, then this vanishes locally on $S'$. By analogy with the classical situation, we would thus expect the open subspace
\[
\mathcal M_Z^{\mathrm{sm}}\subset \mathcal M_Z,
\]
where $s^\ast T_{Z/X_S}$ has positive Harder--Narasimhan slopes, to be (cohomologically) smooth over $S$. Our key geometric result confirms this, at least if $Z\to X_S$ is quasiprojective.

\begin{theorem}\label{thm:introJacobianCriterion} Assume that $Z\to X_S$ can, locally on $S$, be embedded as a Zariski closed subset of an open subset of (the adic space) $\mathbb P^n_{X_S}$. Then $\mathcal M_Z\to S$ is representable in locally spatial diamonds, compactifiable, and of locally finite $\dimtrg$. Moreover, the open subset $\mathcal M_Z^{\mathrm{sm}}\subset \mathcal M_Z$ is cohomologically smooth over $S$.
\end{theorem}

In the application, the space $Z\to X_S$ will be the flag variety parametrizing $\mathbb Q$-filtrations on a given $G$-torsor $\mathcal E$ on $X_S$. Then $\mathcal M_b$ will be an open subset of $\mathcal M_Z^{\mathrm{sm}}$.

The proof of Theorem~\ref{thm:introJacobianCriterion} requires several innovations. The first is a notion of formal smoothness, in which infinitesimal thickenings (that are not available in this perfectoid setting) are replaced by small \'etale neighborhoods. This leads to a notion with a close relation to the notion of absolute neighborhood retracts \cite{BorsukRetractLivre} in classical topology. We prove that virtually all examples of cohomologically smooth maps are also formally smooth, including Banach--Colmez spaces and $\Bun_G$. We also prove that $\mathcal M_Z^{\mathrm{sm}}\to S$ is formally smooth, which amounts to some delicate estimates, spreading sections $X_{T_0}\to Z$ into small neighborhoods of $T_0\subset T$, for any Zariski closed immersion $T_0\subset T$ of affinoid perfectoid spaces --- here we crucially use the assumption that all Harder--Narasimhan slopes are positive. Coupled with the theorem that Zariski closed implies strongly Zariski closed \cite{BhattScholzePrism} this makes it possible to write $\mathcal M_Z^{\mathrm{sm}}$, up to (cohomologically and formally) smooth maps, as a retract of a space that is \'etale over a ball over $S$. Certainly in classical topology, this is not enough to ensure cohomological smoothness --- a coordinate cross is a retract of $\mathbb R^2$ --- but it does imply that the constant sheaf $\mathbb F_\ell$ is universally locally acyclic over $S$. For this reason, and other applications to sheaves on $\Bun_G$ as well as geometric Satake, we thus also develop a general theory of universally locally acyclic sheaves in our setting. To finish the proof, we use a deformation to the normal cone argument to show that the dualizing complex is ``the same'' as the one for the Banach--Colmez space $\BC(s^\ast T_{Z/X_S})$.

\section{$\ell$-adic sheaves on $\mathrm{Bun}_G$}\label{sec:introBunGladic}

For our results, we need to define the category of $\ell$-adic sheaves on $\Bun_G$. More precisely, we will define for each $\mathbb Z_\ell$-algebra $\Lambda$ a category
\[
D(\Bun_G,\Lambda)
\]
of sheaves of $\Lambda$-modules on $\Bun_G$. If $\Lambda$ is killed by some power of $\ell$, such a definition is the main achievement of \cite{ECoD}. Our main interest is however the case $\Lambda=\overline{\mathbb Q}_\ell$. In the case of schemes (of finite type over an algebraically closed field), the passage from torsion coefficients to $\overline{\mathbb Q}_\ell$-coefficients is largely formal: Roughly,
\[
D(X,\overline{\mathbb Q}_\ell) = \Ind(\varprojlim_n D^b_c(X,\mathbb Z/\ell^n\mathbb Z)\otimes_{\mathbb Z_\ell} \overline{\mathbb Q}_\ell).
\]
Behind this definition are however strong finiteness results for constructible sheaves; in particular, the morphism spaces between constructible sheaves are finite. For $\Bun_G$, or for the category of smooth representations, there are still compact objects (given by compactly induced representations in the case of smooth representations), but their endomorphism algebras are Hecke algebras, which are infinite-dimensional. A definition along the same lines would then replace all Hecke algebras by their $\ell$-adic completions, which would drastically change the category of representations.

Our definition of $D(\Bun_G,\Lambda)$ in general involves some new ideas, employing the idea of solid modules developed by Clausen--Scholze \cite{Condensed} in the context of the pro-\'etale (or v-)site; in the end, $D(\Bun_G,\Lambda)$ is defined as a certain full subcategory
\[
D_\lis(\Bun_G,\Lambda)\subset D_\solid(\Bun_G,\Lambda)
\]
of the category $D_\solid(\Bun_G,\Lambda)$ of solid complexes of $\Lambda$-modules on the v-site of $\Bun_G$. The formalism of solid sheaves, whose idea is due to Clausen and the second author, is developed in Chapter~\ref{ch:solid}. It presents some interesting surprises; in particular, there is always a left adjoint $f_\natural$ to pullback $f^\ast$, satisfying base change and a projection formula. (In return, $Rf_!$ fails to exist in general.)

\begin{theorem} Let $\Lambda$ be any $\mathbb Z_\ell$-algebra.
\begin{altenumerate}
\item[{\rm (i)}] Via excision triangles, there is an infinite semiorthogonal decomposition of $D(\Bun_G,\Lambda)$ into the various $D(\Bun_G^b,\Lambda)$ for $b\in B(G)$.
\item[{\rm (ii)}] For each $b\in B(G)$, pullback along
\[
\Bun_G^b\cong [\ast/\tilde{G}_b]\to [\ast/\underline{G_b(E)}]
\]
gives an equivalence
\[
D([\ast/\underline{G_b(E)}],\Lambda)\cong D(\Bun_G^b,\Lambda),
\]
and $D([\ast/\underline{G_b(E)}],\Lambda)\cong D(G_b(E),\Lambda)$ is equivalent to the derived category of the category of smooth representations of $G_b(E)$ on $\Lambda$-modules.
\item[{\rm (iii)}] The category $D(\Bun_G,\Lambda)$ is compactly generated, and a complex $A\in D(\Bun_G,\Lambda)$ is compact if and only if for all $b\in B(G)$, the restriction
\[
i^{b\ast} A\in D(\Bun_G^b,\Lambda)\cong D(G_b(E),\Lambda)
\]
is compact, and zero for almost all $b$. Here, compactness in $D(G_b(E),\Lambda)$ is equivalent to lying in the thick triangulated subcategory generated by $c\text-\mathrm{Ind}_K^{G_b(E)} \Lambda$ as $K$ runs over open pro-$p$-subgroups of $G_b(E)$.
\item[{\rm (iv)}] On the subcategory $D(\Bun_G,\Lambda)^\omega\subset D(\Bun_G,\Lambda)$ of compact objects, there is a Bernstein--Zelevinsky duality functor
\[
\mathbb D_{BZ}: (D(\Bun_G,\Lambda)^\omega)^{\mathrm{op}}\to D(\Bun_G,\Lambda)^\omega
\]
with a functorial identification
\[
R\Hom(A,B)\cong \pi_\natural(\mathbb D_{BZ}(A)\dotimes_\Lambda B)
\]
for $B\in D(\Bun_G,\Lambda)$, where $\pi: \Bun_G\to \ast$ is the projection. The functor $\mathbb D_{BZ}$ is an equivalence, and $\mathbb D_{BZ}^2$ is naturally equivalent to the identity. It is compatible with usual Bernstein--Zelevinsky duality on $D(G_b(E),\Lambda)$ for basic $b\in B(G)$.
\item[{\rm (v)}] An object $A\in D(\Bun_G,\Lambda)$ is universally locally acyclic (with respect to $\Bun_G\to \ast$) if and only if for all $b\in B(G)$, the restriction
\[
i^{b\ast} A\in D(\Bun_G^b,\Lambda)\cong D(G_b(E),\Lambda)
\]
is admissible, i.e.~for all pro-$p$ open subgroups $K\subset G_b(E)$, the complex $(i^{b\ast} A)^K$ is perfect. Universally locally acyclic complexes are preserved by Verdier duality, and satisfy Verdier biduality.
\end{altenumerate}
\end{theorem}

This theorem extends many basic notions from representation theory --- finitely presented objects, admissible representations, Bernstein--Zelevinsky duality, smooth duality --- to the setting of $D(\Bun_G,\Lambda)$.

Parts (i) and (ii) are easy when $\Lambda$ is $\ell$-power torsion. In general, their proofs invoke the precise definition of $D(\Bun_G,\Lambda)=D_\lis(\Bun_G,\Lambda)$ and are somewhat subtle. Part (iii) uses that $i^{b\ast}$ admits a left adjoint, which will then automatically preserve compact objects (inducing compact generators). Using the diagram
\[
[\ast/\underline{G_b(E)}]\xleftarrow{q_b} \mathcal M_b\xrightarrow{\pi_b}\Bun_G,
\]
this left adjoint is defined as $\pi_{b\natural} q_b^\ast$. The verification that this is indeed a left adjoint amounts in some sense to the assertion that $\mathcal M_b$ is ``strictly local'' along the closed subspace $[\ast/\underline{G_b(E)}]\subset \mathcal M_b$ in the sense that for all $A\in D(\mathcal M_b,\Lambda)$, the restriction
\[
R\Gamma(\mathcal M_b,A)\to R\Gamma([\ast/\underline{G_b(E)}],A)
\]
is an isomorphism. This builds on a detailed analysis of the topological nature of $\mathcal M_b$, in particular that $\tilde{\mathcal M}_b\setminus \ast$ is a spatial diamond, and Theorem~\ref{thm:intropartialcompactsupport} below. For part (iv), the constructions in (iii) imply the existence of $\mathbb D_{BZ}(A)$ on a class of generators, thus in general, and similar arguments to the ones in (iii) prove the biduality. Finally, part (v) is essentially a formal consequence.

The key cohomological result for the proof is the following result, applied to $\tilde{\mathcal M}_b\setminus \ast$ (or quotients of it). It plays on the subtle point that the point $\ast$ is not quasiseparated.

\begin{theorem}\label{thm:intropartialcompactsupport} Let $X$ be a spatial diamond such that $f: X\to \ast$ is partially proper, and of finite $\dimtrg$. Then for any affinoid perfectoid space $S$, the base change $X_S=X\times S$ naturally admits two ends. Taking compactly supported cohomology with respect to one end (but no support condition at the other end), one has
\[
R\Gamma_{\partial\text-c}(X_S,A)=0
\]
for all $A\in D_\solid^+(X,\mathbb Z_\ell)$ (resp.~all $A\in D_\solid(X,\mathbb Z_\ell)$ if $f$ is $\ell$-cohomologically smooth).
\end{theorem}

As an example, if $X=\Spa \overline{\mathbb F}_q\laurentseries{t}$, then $X_S=\mathbb D^\ast_S$ is an open unit disc over $S$, whose two ends are the origin and the boundary, and one has
\[
R\Gamma_{\partial\text-c}(\mathbb D^\ast_S,\mathbb Z_\ell)=0.
\]
In particular, the cohomology of $\Spa \overline{\mathbb F}_q\powerseries{t}$ agrees with sections on the closed point, showing that $\Spa \overline{\mathbb F}_q\powerseries{t}$ is ``strictly local''. The same phenomenon is at work for $\mathcal M_b$.

\section{The geometric Satake equivalence}\label{sec:introGeomSatake}

In order to define the Hecke operators, we need to prove the geometric Satake equivalence, taking representations of the dual group $\hat{G}$ to sheaves on the local Hecke stack. In order to analyze compositions of Hecke operators, it will in fact be necessary to analyze modifications at several points.

Thus, for any finite set $I$, we consider the moduli space $(\Div^1)^I$ parametrizing degree $1$ Cartier divisors $D_i\subset X_S$, $i\in I$. Locally on $S$, each $D_i$ defines an untilt $S_i^\sharp$ of $S$ over $E$, and one can form the completion $B^+$ of $\mathcal O_{X_S}$ along the union of the $D_i$. Inverting the $D_i$ defines a localization $B$ of $B^+$. One can then define a positive loop group $L^+_{(\Div^1)^I} G$ and loop group $L_{(\Div^1)^I} G$, with values given by $G(B^+)$ resp.~$G(B)$; for brevity, we will simply write $L^+G$ and $LG$ here. One can then define the local Hecke stack
\[
\Hloc_G^I = [L^+ G\backslash LG / L^+ G]\to (\Div^1)^I.
\]
For $d=|I|$, this is in fact already defined over the moduli space $\Div^d = (\Div^1)^d/\Sigma_d$ of degree $d$ Cartier divisors. We will often break symmetry, and first take the quotient on the right to define the Beilinson--Drinfeld Grassmannian
\[
\Gr_G^I = LG / L^+ G\to (\Div^1)^I
\]
so that
\[
\Hloc_G^I = L^+ G\backslash \Gr_G^I.
\]

The Beilinson--Drinfeld Grassmannian $\Gr_G^I\to (\Div^1)^I$ is a small v-sheaf that can be written as an increasing union of closed subsheaves that are proper and representable in spatial diamonds, by bounding the relative position; this is one main result of \cite{Berkeley}. On the other hand, $L^+ G$ can be written as an inverse limit of truncated positive loop groups, which are representable in locally spatial diamonds and cohomologically smooth; moreover, on each bounded subset, it acts through such a finite-dimensional quotient. This essentially reduces the study of all bounded subsets of $\Hloc_G^I$ to Artin stacks.

In particular, one can write the local Hecke stack as an increasing union of closed substacks that are quasicompact over $(\Div^1)^I$, by bounding the relative position. In the following, we assume that the coefficients $\Lambda$ are killed by some power of $\ell$, so that we can use the theory from \cite{ECoD}. Let
\[
D_\et(\Hloc_G^I,\Lambda)^{\mathrm{bd}}\subset D_\et(\Hloc_G^I,\Lambda)
\]
be the full subcategory of all objects with quasicompact support over $(\Div^1)^I$. This is a monoidal category under convolution $\star$. Here, we use the convolution diagram
\[
\Hloc_G^I\times_{(\Div^1)^I} \Hloc_G^I\xleftarrow{(p_1,p_2)} L^+ G\backslash L G\times^{L^+ G} L G / L^+ G\xrightarrow{m} \Hloc_G^I
\]
and define
\[
A\star B = Rm_\ast(p_1^\ast A\dotimes_\Lambda p_2^\ast B).
\]
The map $m$ is ind-proper (its fibres are $\Gr_G^I$), and in particular proper on any bounded subset; thus, proper base change ensures that this defines an associative monoidal structure.

On $D_\et(\Hloc_G^I,\Lambda)^{\mathrm{bd}}$, one can define a relative perverse $t$-structure (where an object is perverse if and only if it is perverse over any geometric fibre of $(\Div^1)^I$). For this $t$-structure, the convolution $\star$ is left $t$-exact (and $t$-exactness only fails for issues related to non-flatness over $\Lambda$). To prove that there is a well-defined $t$-structure, and the preservation of perversity under convolution, we adapt Braden's theorems \cite{Braden} on hyperbolic localization, and a degeneration to the Witt vector affine Grassmannian \cite{ZhuWitt}, \cite{BhattScholzeWitt}. We will discuss hyperbolic localization further below.

We remark that there is no general theory of perverse sheaves in $p$-adic geometry, the issue being that it is difficult to unambiguously assign a dimension to a point of an adic space (cf.~\cite{TemkinTranscendence} for what is known about topological transcendence degrees of points, and the subtleties especially in characteristic $p$). In particular, we would not know how to define a notion of perverse sheaf on $(\Div^1)^I$ in general, which is the reason we revert to asking perversity only in the fibres. Here, we use that all geometric fibres of the stack $\Hloc_G^I\to (\Div^1)^I$ have only countably many points enumerated explicitly in terms of dominant cocharacters $\mu_i$, and one can assign by hand the dimension $\sum_i \langle 2\rho,\mu_i\rangle$ of the corresponding open Schubert cells.

\begin{remark} Inspired by this, we realized that for any map $f: X\to S$ locally of finite type between schemes, one can define a relative perverse $t$-structure, with relative perversity equivalent to perversity on all geometric fibres, cf.~\cite{HansenScholze}.
\end{remark}

Moreover, one can restrict to the complexes $A\in D_\et(\Hloc_G^I,\Lambda)^{\mathrm{bd}}$ that are universally locally acyclic over $(\Div^1)^I$. This condition is also preserved under convolution.

\begin{definition} The Satake category
\[
\Sat_G^I(\Lambda)\subset D_\et(\Hloc_G^I,\Lambda)^{\mathrm{bd}}
\]
is the category of all $A\in D_\et(\Hloc_G^I,\Lambda)^{\mathrm{bd}}$ that are perverse, flat over $\Lambda$ (i.e., for all $\Lambda$-modules $M$, also $A\dotimes_\Lambda M$ is perverse), and universally locally acyclic over $(\Div^1)^I$.
\end{definition}

Intuitively, $\Sat_G^I(\Lambda)$ are the ``flat families of perverse sheaves on $\Hloc_G^I\to (\Div^1)^I$'', where flatness refers both to the geometric aspect of flatness over $(\Div^1)^I$ (encoded in universal local acyclicity) and the algebraic aspect of flatness in the coefficients $\Lambda$. The Satake category $\Sat_G^I(\Lambda)$ is a monoidal category under convolution. Moreover, it is covariantly functorial in $I$.

In fact, the monoidal structure naturally upgrades to a symmetric monoidal structure. This relies on the fusion product, for which it is critical to allow general finite sets $I$. Namely, given finite sets $I_1,\ldots,I_k$, letting $I=I_1\sqcup\ldots\sqcup I_k$, one has an isomorphism
\[
\Hloc_G^I\times_{(\Div^1)^I} (\Div^1)^{I;I_1,\ldots,I_k}\cong \prod_{j=1}^k \Hloc_G^{I_j}\times_{(\Div^1)^I} (\Div^1)^{I;I_1,\ldots,I_k}
\]
where $(\Div^1)^{I;I_1,\ldots,I_k}\subset (\Div^1)^I$ is the open subset where $x_i\neq x_{i'}$ whenever $i,i'\in I$ lie in different $I_j$'s. The exterior tensor product then defines a functor
\[
\boxtimes_{j=1}^k: \prod_{j=1}^k \Sat_G^{I_j}(\Lambda)\to \Sat_G^{I;I_1,\ldots,I_k}(\Lambda)
\]
where $\Sat_G^{I;I_1,\ldots,I_k}(\Lambda)$ is the variant of $\Sat_G^I(\Lambda)$ for $\Hloc_G^I\times_{(\Div^1)^I} (\Div^1)^{I;I_1,\ldots,I_k}$. However, the restriction functor
\[
\Sat_G^I(\Lambda)\to \Sat_G^{I;I_1,\ldots,I_k}(\Lambda)
\]
is fully faithful, and the essential image of the exterior product lands in its essential image. Thus, we get a natural functor
\[
\bigast_{j=1}^k: \prod_{j=1}^k \Sat_G^{I_j}(\Lambda)\to \Sat_G^I(\Lambda),
\]
independent of the ordering of the $I_j$. In particular, for any $I$, we get a functor
\[
\Sat_G^I(\Lambda)\times \Sat_G^I(\Lambda)\to \Sat_G^{I\sqcup I}(\Lambda)\to \Sat_G^I(\Lambda),
\]
using functoriality of $\Sat_G^J(\Lambda)$ in $J$, which defines a symmetric monoidal structure $\ast$ on $\Sat_G^I(\Lambda)$, commuting with $\star$. This is called the fusion product. In general, for any symmetric monoidal category $(\mathcal C,\ast)$ with a commuting monoidal structure $\star$, the monoidal structure $\star$ necessarily agrees with $\ast$; thus, the fusion product refines the convolution product. (As usual in geometric Satake, we actually need to change $\ast$ slightly by introducing certain signs into the commutativity constraint, depending on the parity of the support of the perverse sheaves.)

Moreover, restricting $A\in \Sat_G^I(\Lambda)$ to $\Gr_G^I$ and taking the pushforward to $(\Div^1)^I$, all cohomology sheaves are local systems of $\Lambda$-modules on $(\Div^1)^I$. By a version of Drinfeld's lemma, these are equivalent to representations of $W_E^I$ on $\Lambda$-modules. This defines a symmetric monoidal fibre functor
\[
F^I: \Sat_G^I(\Lambda)\to \Rep_{W_E^I}(\Lambda),
\]
where $\Rep_{W_E^I}(\Lambda)$ is the category of continuous representations of $W_E^I$ on finite projective $\Lambda$-modules. Using a version of Tannaka duality, one can then build a Hopf algebra in the $\Ind$-category of $\Rep_{W_E^I}(\Lambda)$ so that $\Sat_G^I(\Lambda)$ is given by its category of representations (internal in $\Rep_{W_E^I}(\Lambda)$). For any finite set $I$, this is given by the tensor product of $I$ copies of the corresponding Hopf algebra for $I=\{\ast\}$, which in turn is given by some affine group scheme $\widecheck{G}$ over $\Lambda$ with $W_E$-action.

\begin{theorem} There is a canonical isomorphism $\widecheck{G}\cong \hat{G}_\Lambda$ with the Langlands dual group, under which the action of $W_E$ on $\widecheck{G}$ agrees with the usual action of $W_E$ on $\hat{G}$ up to an explicit cyclotomic twist. If $\sqrt{q}\in \Lambda$, the cyclotomic twist can be trivialized, and $\Sat_G^I(\Lambda)$ is naturally equivalent to the category of $(\hat{G}\rtimes W_E)^I$-representations on finite projective $\Lambda$-modules.
\end{theorem}

This theorem is thus a version of the theorem of Mirkovi\'c--Vilonen \cite{MirkovicVilonen}, coupled with the refinements of Gaitsgory \cite{GaitsgoryBeilinsonDrinfeldSatake} for general $I$. (We remark that we formulate a theorem valid for any $\Lambda$, not necessarily regular; such a formulation does not seem to be in the literature. Also, we give a purely local proof: Most proofs require a globalization on a (usual) curve.) Contrary to Mirkovi\'c--Vilonen, we actually construct an explicit pinning of $\widecheck{G}$. For the proof, one can restrict to $\Lambda=\mathbb Z/\ell^n\mathbb Z$; passing to a limit over $n$, one can actually build a group scheme over $\mathbb Z_\ell$. Its generic fibre is reductive, as the Satake category with $\mathbb Q_\ell$-coefficients is (geometrically) semisimple: For this, we again use the degeneration to the Witt vector affine Grassmannian and the decomposition theorem for schemes. To identify the reductive group, we argue first for tori, and then for rank $1$ groups, where everything reduces to $G=\mathrm{PGL}_2$ which is easy to analyze by using the minuscule Schubert cell. Here, the pinning includes a cyclotomic twist as of course the cohomology of the minuscule Schubert variety $\mathbb P^1$ of $\Gr_{\mathrm{\PGL}_2}$ contains a cyclotomic twist. Afterwards, we apply hyperbolic localization in order to construct symmetric monoidal functors $\Sat_G\to \Sat_M$ for any Levi $M$ of $G$, inducing dually maps $\widecheck{M}\to \widecheck{G}$. This produces many Levi subgroups of $\widecheck{G}_{\mathbb Q_\ell}$ from which it is easy to get the isomorphism with $\hat{G}_{\mathbb Q_\ell}$, including a pinning. As these maps $\widecheck{M}\to \widecheck{G}$ are even defined integrally, and $\hat{G}(\mathbb Z_\ell)\subset \hat{G}(\mathbb Q_\ell)$ is a maximal compact open subgroup by Bruhat--Tits theory, generated by the rank $1$ Levi subgroups, one can then deduce that $\widecheck{G}\cong \hat{G}$ integrally, again with an explicit (cyclotomic) pinning.

We will also need the following addendum regarding a natural involution. Namely, the local Hecke stack $\Hloc_G$ has a natural involution $\mathrm{sw}$ given by reversing the roles of the two $G$-torsors; in the presentation in terms of $LG$, this is induced by the inversion on $LG$. Then $\mathrm{sw}^\ast$ induces naturally an involution of $\Sat_G(\Lambda)$, and this involution can be upgraded to a symmetric monoidal functor commuting with the fibre functor, thus realizing a $W_E$-equivariant automorphism of $\widecheck{G}\cong \hat{G}$.

\begin{proposition}\label{prop:chevalleyinvolutionintro} The action of $\mathrm{sw}^\ast$ on $\mathrm{Sat}_G$ induces the automorphism of $\hat{G}$ that is the Chevalley involution of the split group $\hat{G}$, conjugated by $\hat{\rho}(-1)$.
\end{proposition}

Critical to all of our arguments is the hyperbolic localization functor. In the setting of the Beilinson--Drinfeld Grassmannian, assume that $P^+,P^-\subset G$ are two opposite parabolics, with common Levi $M$. We get a diagram
\[\xymatrix{
& \Gr_{P^+}^I\ar[dl]_{q^+}\ar[dr]^{p^+} & \\
\Gr_G^I && \Gr_M^I.\\
& \Gr_{P^-}^I\ar[ul]_{q^-}\ar[ur]^{p^-} &
}\]
We get two ``constant term'' functors
\[
\mathrm{CT}^+ = R(p^+)_! (q^+)^\ast, \mathrm{CT}^- = R(p^-)_\ast R(q^-)^!: D_\et(\Gr_G^I,\Lambda)^{\mathrm{bd}}\to D_\et(\Gr_M^I,\Lambda)^{\mathrm{bd}},
\]
and one can construct a natural transformation $\mathrm{CT}^-\to \mathrm{CT}^+$. The functor $\mathrm{CT}^+$ corresponds classically to the Satake transform, of integrating along orbits under the unipotent radical of $U^+$. Hyperbolic localization claims that the transformation $\mathrm{CT}^-\to \mathrm{CT}^+$ is an equivalence when restricted to $L^+ G$-equivariant objects. This has many consequences; note that $\mathrm{CT}^+$ is built from left adjoint functors while $\mathrm{CT}^-$ is built from right adjoint functors, so if they are isomorphic, hyperbolic localization has the best of both worlds. In particular, hyperbolic localization commutes with all colimits and all limits, preserves (relative) perversity, universal local acyclicity, commutes with any base change, etc.~.

This is in fact a special case of the following more general assertion. Let $S$ be any small v-stack, and $f: X\to S$ be a proper map that is representable in spatial diamonds with $\dimtrg f<\infty$. Assume that there is an action of $\mathbb G_m$ on $X/S$, where $\mathbb G_m(R,R^+)= R^\times$. The fixed points $X^0\subset X$ of the $\mathbb G_m$-action form a closed substack. We assume that one can define an attractor locus $X^+\subset X$ and a repeller locus $X^-\subset X$, given by disjoint unions of locally closed subspaces, on which the $t\in \mathbb G_m$-action admits a limit as $t\to 0$ (resp.~$t\to \infty$). We get a diagram
\[\xymatrix{
& X^+\ar[dl]_{q^+}\ar[dr]^{p^+} & \\
X && X^0,\\
& X^-\ar[ul]_{q^-}\ar[ur]^{p^-} &
}\]
generalizing (bounded parts of) the above diagram if one chooses a cocharacter $\mu: \mathbb G_m\to G$ whose dynamic parabolics are $P^+,P^-$. One can define
\[
L^+ = R(p^+)_! (q^+)^\ast, L^- = R(p^-)_\ast R(q^-)^!: D_\et(X,\Lambda)\to D_\et(X^0,\Lambda)
\]
and a natural transformation $L^-\to L^+$. The following is our version of Braden's theorem \cite{Braden}, cf.~also \cite{RicharzHL}.

\begin{theorem} The transformation $L^-\to L^+$ is an equivalence when restricted to the essential image of $D_\et(X/\mathbb G_m,\Lambda)\to D_\et(X,\Lambda)$.
\end{theorem}

The proof makes use of the following principle: If $Y\to S$ is partially proper with a $\mathbb G_m$-action such that the quotient stack $Y/\mathbb G_m$ is qcqs over $S$, then again $Y$ admits two ends, and the partially compactly supported cohomology of $Y$ with coefficients in any $A\in D_\et(Y/\mathbb G_m,\Lambda)$ vanishes identically.

\section{Cohomology of moduli spaces of shtuka}\label{sec:introRZspaces}

At this point, we have defined
\[
D(\Bun_G,\Lambda),
\]
and using the geometric Satake equivalence and the diagram
\[\xymatrix{
& \Hglob_G^I\ar[r]^q\ar[dl]^{h_1}\ar[dr]^{h_2} & \Hloc_G^I\\
\Bun_G&& \Bun_G\times (\Div^1)^I
}\]
one can define the Hecke operator
\[
T_V = Rh_{2\ast}(h_1^\ast \dotimes_\Lambda q^\ast \mathcal S_V): D(\Bun_G,\Lambda)\to D(\Bun_G\times (\Div^1)^I,\Lambda)
\]
for any $V\in \Sat_G^I(\Lambda)$, where $\mathcal S_V$ is the corresponding sheaf on $\Hloc_G^I$. This works at least if $\Lambda$ is killed by some power of $\ell$. We can in fact extend this functor to all $\mathbb Z_\ell$-algebras $\Lambda$. Moreover, its image lies in the full subcategory of those objects that are locally constant in the direction of $(\Div^1)^I$, thereby giving a functor
\[
T_V: D(\Bun_G,\Lambda)\to D(\Bun_G,\Lambda)^{BW_E^I}
\]
to the category of $W_E^I$-equivariant objects in $D(\Bun_G,\Lambda)$. The proof is surprisingly formal: One reduces to $I=\{\ast\}$ by an inductive argument, and then uses that $\Div^1 = \Spd \breve{E}/\phi^\Z$ is still just a point. More precisely, one uses that
\[
D(\Bun_G,\Lambda)\to D(\Bun_G\times \Spd \hat{\overline{E}},\Lambda)
\]
is an equivalence.

\begin{remark} To define $D(\Bun_G,\Lambda)^{BW_E^I}$, we need to upgrade $D(\Bun_G,\Lambda)$ to a condensed $\infty$-category; then it is the notion of $W_E^I$-equivariant objects for the condensed group $W_E^I$.
\end{remark}

A first consequence of our results is that $T_V$, forgetting the $W_E^I$-equivariance, preserves finiteness properties. Note that $T_V\circ T_W\cong T_{V\otimes W}$ as the geometric Satake equivalence is monoidal. This formally implies that $T_V$ is left and right adjoint to $T_{V^\ast}$. From here, it is not hard to prove the following result.

\begin{theorem} The functor $T_V: D(\Bun_G,\Lambda)\to D(\Bun_G,\Lambda)$ preserves compact objects and universally locally acyclic objects. Moreover, it commutes with Bernstein--Zelevinsky and Verdier duality in the sense that there are natural isomorphisms $\mathbb D_{\mathrm{BZ}}(T_V(A))\cong T_{\mathrm{sw}^\ast V^\vee}(\mathbb D_{\mathrm{BZ}}(A))$ and $R\sHom(T_V(A),\Lambda)\cong T_{(\mathrm{sw}^\ast V)^\vee} R\sHom(A,\Lambda)$.
\end{theorem}

Here $\mathrm{sw}^\ast$ is the involution of $\Sat_G^I$ which by Proposition~\ref{prop:chevalleyinvolutionintro} is induced by the Chevalley involution of $\hat{G}$, conjugated by $\hat{\rho}(-1)$.

This theorem has concrete consequences for the cohomology of moduli spaces of shtukas. For simplicity, we formulate it here with coefficients in a $\Lambda$-algebra that is killed by $\ell^n$ for some $n$; for the general formulation, we would need to discuss more precisely the foundational issues surrounding the derived categories. In \cite[Lecture XXIII]{Berkeley}, for any collection $\{\mu_i\}_i$ of conjugacy classes of cocharacters with fields of definition $E_i/E$ and $b\in B(G)$, there is defined a tower of moduli spaces of local shtukas
\[
f_K: (\Sht_{(G,b,\mu_\bullet),K})_{K\subset G(E)}\to \prod_{i\in I} \Spd \breve{E}_i
\]
as $K$ ranges over compact open subgroups of $G(E)$, equipped with compatible \'etale period maps
\[
\pi_K: \Sht_{(G,b,\mu_\bullet),K}\to \Gr^\tw_{G,\prod_{i\in I} \Spd \breve{E}_i,\leq \mu_\bullet}.
\]
Here, $\Gr^\tw_{G,\prod_{i\in I} \Spd \breve{E}_i}\to \prod_{i\in I} \Spd \breve{E}$ is a certain twisted form of the convolution affine Grassmannian, cf.~\cite[Section 23.5]{Berkeley}. Let $W$ be the exterior tensor product $\boxtimes_{i\in I} V_{\mu_i}$ of highest weight representations, and $\mathcal S_W$ the corresponding sheaf on $\Gr^\tw_{G,\prod_{i\in I} \Spd \breve{E}_i}$. We continue to write $\mathcal S_W$ for its pullback to $\Sht_{(G,b,\mu_\bullet),K}$.

\begin{corollary}\label{cor:finitenessmain} The sheaf
\[
Rf_{K!} \mathcal S_W\in D([\ast/\underline{G_b(E)}]\times \prod_{i\in I} \Spd \breve{E}_i,\Lambda)
\]
is equipped with partial Frobenii, thus descends to an object of
\[
D([\ast/\underline{G_b(E)}]\times \prod_{i\in I} \Spd \breve{E}_i/\varphi_i^{\mathbb Z},\Lambda).
\]
This object lives in the full subcategory
\[
D(G_b(E),\Lambda)^{B\prod_{i\in I} W_{E_i}}\subset D([\ast/\underline{G_b(E))}]\times \prod_{i\in I} \Spd \breve{E}_i/\varphi_i^{\mathbb Z},\Lambda),
\]
and its restriction to $D(G_b(E),\Lambda)$ is compact. In particular, for any admissible representation $\rho$ of $G_b(E)$, the object
\[
R\Hom_{G_b(E)}(Rf_{K!} \mathcal S_W,\rho)\in D(\Lambda)^{B\prod_{i\in I} W_{E_i}}
\]
is a representation of $\prod_{i\in I} W_{E_i}$ on a perfect complex of $\Lambda$-modules. Taking the colimit over $K$, this gives rise to a complex of admissible $G(E)$-representations
\[
\varinjlim_K R\Hom_{G_b(E)}(Rf_{K!} \mathcal S_W,\rho)
\]
equipped with a $\prod_{i\in I} W_{E_i}$-action.

If $\rho$ is compact, then so is
\[
\varinjlim_K R\Hom_{G_b(E)}(Rf_{K!} \mathcal S_W,\rho)
\]
as a complex of $G(E)$-representations.
\end{corollary}

Specializing to $I=\{\ast\}$ and $\mu$ minuscule, we get local Shimura varieties, and this proves the finiteness properties of \cite[Proposition 6.1]{RapoportViehmann} unconditionally, as well as \cite[Remark 6.2 (iii)]{RapoportViehmann}. We note that those properties seem inaccessible using only the definition of the moduli spaces of shtukas, i.e.~without the use of $\Bun_G$.

\section{The stack of $L$-parameters}\label{sec:introstackLparameters}

Let us discuss the other side of the Langlands correspondence, namely (the stack of) $L$-parameters. This has been previously done by Dat--Helm--Kurinczuk--Moss \cite{DatHelmKurinczukMoss} and Zhu \cite{ZhuConjectures}. One wants to define a scheme whose $\Lambda$-valued points, for a $\mathbb Z_\ell$-algebra $\Lambda$, are the continuous $1$-cocycles
\[
\varphi: W_E\to \hat{G}(\Lambda).
\]
(Here, we endow $\hat{G}$ with its usual $W_E$-action, that factors over a finite quotient $Q$ of $W_E$. As discussed above, the difference between the two actions disappears over $\mathbb Z_\ell[\sqrt{q}]$, and we find it much more convenient to use the standard normalization here, so that we can sometimes make use of the algebraic group $\hat{G}\rtimes Q$.)

There seems to be a mismatch here, in asking for an algebraic stack, but continuous cocycles. Interestingly, there is a way to phrase the continuity condition that produces a scheme. Namely, we consider $\Lambda$ as a condensed $\mathbb Z_\ell$-algebra that is ``relatively discrete over $\mathbb Z_\ell$''. Abstract $\mathbb Z_\ell$-modules $M$ embed fully faithfully into condensed $\mathbb Z_\ell$-modules, via sending $M$ to $M_{\mathrm{disc}}\otimes_{\mathbb Z_{\ell,\mathrm{disc}}} \mathbb Z_\ell$.

\begin{theorem} There is a scheme $Z^1(W_E,\hat{G})$ over $\mathbb Z_\ell$ whose $\Lambda$-valued points, for a $\mathbb Z_\ell$-algebra $\Lambda$, are the condensed $1$-cocycles
\[
\varphi: W_E\to \hat{G}(\Lambda),
\]
where we regard $\Lambda$ as a relatively discrete condensed $\mathbb Z_\ell$-algebra. The scheme $Z^1(W_E,\hat{G})$ is a union of open and closed affine subschemes $Z^1(W_E/P,\hat{G})$ as $P$ runs through open subgroups of the wild inertia subgroup of $W_E$, and each $Z^1(W_E/P,\hat{G})$ is a flat local complete intersection over $\mathbb Z_\ell$ of dimension $\dim G$.
\end{theorem}

The point here is that the inertia subgroup of $W_E$ has a $\mathbb Z_\ell$-factor, and this can map in interesting ways to $\Lambda$ when making this definition. To prove the theorem, following \cite{DatHelmKurinczukMoss} and \cite{ZhuConjectures} we define discrete dense subgroups $W\subset W_E/P$ by discretizing the tame inertia, and the restriction $Z^1(W_E/P,\hat{G})\to Z^1(W,\hat{G})$ is an isomorphism, where the latter is clearly an affine scheme.

We can also prove further results about the $\hat{G}$-action on $Z^1(W_E,\hat{G})$, or more precisely each $Z^1(W_E/P,\hat{G})$. For this result, we need to make a very minor assumption on $\ell$.

\begin{theorem}\label{thm:verygoodprimeintro} Assume that $\ell$ does not divide the order of $\pi_0 Z(G)$ (equivalently, $\ell$ does not divide the order of $\pi_1(\widehat{G})_{\mathrm{tor}}$). Then $H^i(\hat{G},\mathcal O(Z^1(W_E/P,\hat{G})))=0$ for $i>0$ and the formation of the invariants $\mathcal O(Z^1(W_E/P,\hat{G}))^{\hat{G}}$ commutes with any base change. The algebra $\mathcal O(Z^1(W_E/P,\hat{G}))^{\hat{G}}$ admits an explicit presentation in terms of excursion operators,
\[
\mathcal O(Z^1(W_E/P,\hat{G}))^{\hat{G}} = \colim_{(n,F_n\to W)} \mathcal O(Z^1(F_n,\hat{G}))^{\hat{G}}
\]
where the colimit runs over all maps from a free group $F_n$ to $W\subset W_E/P$, and $Z^1(F_n,\hat{G})\cong \hat{G}^n$ with the simultaneous twisted $\hat{G}$-conjugation.

Moreover, the $\infty$-category $\Perf(Z^1(W_E/P,\hat{G})/\hat{G})$ is generated under cones and retracts by the image of $\Rep(\hat{G})\to \Perf(Z^1(W_E/P,\hat{G})/\hat{G})$, and $\Ind\Perf(Z^1(W_E/P,\hat{G}))$ is equivalent to the $\infty$-category of modules over $\mathcal O(Z^1(W_E/P,\hat{G}))$ in $\Ind\Perf(\ast/\hat{G})$.

All of these results also hold with $\mathbb Q_\ell$-coefficients, without the assumption on $\ell$.
\end{theorem}

With $\mathbb Q_\ell$-coefficients, these results are simple, as the representation theory of $\hat{G}$ is semisimple. However, with $\mathbb Z_\ell$-coefficients, these results are quite subtle, and we need to dive into modular representation theory of reductive groups. In fact, we give a new perspective on (and provide some new examples of) the phenomenon that restriction along an embedding $H\subset G$ of reductive groups preserves representations admitting a good filtration.

\begin{theorem} Let $G$ be a reductive group over an algebraically closed field $L$ of characteristic $\ell$. Let $P$ be a finite solvable group of order prime to $\ell$ acting on $G$. The fixed point group $H=G^P$ is a smooth linear algebraic group with $H^\circ$ reductive, and with $\pi_0 H$ of order prime to $\ell$.

In this situation, for any representation $V$ of $G$ admitting a good $G$-filtration, also $V|_{H^\circ}$ admits a good $H^\circ$-filtration.
\end{theorem}

The case of Levi subgroups is a classical theorem, while the case $P=\mathbb Z/2\mathbb Z$ was known as Brundan's conjecture and proved by exhaustive case-by-case analysis in \cite{Brundan}, \cite{vanderKallenBrundan}. We give a new proof that works uniformly in all cases.

We also prove that in the situation of the theorem, the image of $\Perf(\ast/G)\to \Perf(\ast/H)$ generates the whole category under cones and retracts. In the first version, we proved this by a very explicit (and exhausting) analysis of all possible cases, but there is now a uniform proof.

\section{Construction of $L$-parameters}\label{sec:introLparameters}

Finally, we can discuss the construction of $L$-parameters. Assume first for simplicity that $\Lambda=\overline{\mathbb Q}_\ell$ with fixed $\sqrt{q}\in \overline{\mathbb Q}_\ell$, and let $A\in D(\Bun_G,\overline{\mathbb Q}_\ell)$ be any Schur-irreducible object, i.e.~$\End(A) = \overline{\mathbb Q}_\ell$. For example, $A$ could correspond to an irreducible smooth representation of $G(E)$, taking the extension by zero along $[\ast/G(E)]\hookrightarrow \Bun_G$. Then, following V.~Lafforgue \cite{VLafforgue1}, we can define excursion operators as follows. For any representation $V$ of $(\hat{G}\rtimes W_E)^I$ over $\overline{\mathbb Q}_\ell$, together with maps
\[
\alpha: \overline{\mathbb Q}_\ell\to V|_{\hat{G}}, \beta: V|_{\hat{G}}\to \overline{\mathbb Q}_\ell
\]
when restricted to the action of the diagonal copy $\hat{G}\subset (\hat{G}\rtimes W_E)^I$, and elements $\gamma_i\in W_E$ for $i\in I$, we can define the endomorphism
\[
A\xrightarrow{T_\alpha} T_V(A)\xrightarrow{(\gamma_i)_{i\in I}} T_V(A)\xrightarrow{T_\beta} A
\]
of $A$, defining an element of $\overline{\mathbb Q}_\ell$. With all the formalism in place, the following result is essentially due to V.~Lafforgue \cite[Proposition 11.7]{VLafforgue1}.

\begin{proposition} There is a unique continuous semisimple $L$-parameter
\[
\varphi_A: W_E\to \hat{G}(\overline{\mathbb Q}_\ell)
\]
such that for all $(I,V,\alpha,\beta,(\gamma_i)_{i\in I})$ as above, the excursion operator
\[
A\xrightarrow{T_\alpha} T_V(A)\xrightarrow{(\gamma_i)_{i\in I}} T_V(A)\xrightarrow{T_\beta} A
\]
is given by multiplication with the scalar
\[
\overline{\mathbb Q}_\ell\xrightarrow{\alpha} V\xrightarrow{(\varphi_A(\gamma_i)_{i\in I})} V\xrightarrow{\beta} \overline{\mathbb Q}_\ell.
\]
\end{proposition}

Note that in fact, the excursion operators define elements in the Bernstein center of $G(E)$, as they define endomorphisms of the identity functor. From this perspective, let us make the following definition.

\begin{definition}\label{def:bernsteincenterintro}\leavevmode
\begin{altenumerate}
\item[{\rm (i)}] The Bernstein center of $G(E)$ is
\[
\mathcal Z(G(E),\Lambda) = \pi_0\mathrm{End}(\mathrm{id}_{\mathcal D(G(E),\Lambda)}) = \varprojlim_{K\subset G(E)} \mathcal Z(\Lambda[K\backslash G(E)/K])
\]
where $K$ runs over open pro-$p$ subgroups of $G(E)$, and $\Lambda[K\backslash G(E)/K]=\mathrm{End}_{G(E)}(c\text-\mathrm{Ind}_K^{G(E)} \Lambda)$ is the Hecke algebra of level $K$.
\item[{\rm (ii)}] The geometric Bernstein center of $G$ is
\[
\mathcal Z^{\mathrm{geom}}(G,\Lambda) = \pi_0\mathrm{End}(\mathrm{id}_{\mathcal D_\lis(\Bun_G,\Lambda)}).
\]
Inside $\mathcal Z^{\mathrm{geom}}(G,\Lambda)$, we let $\mathcal Z^{\mathrm{geom}}_{\mathrm{Hecke}}(G,\Lambda)$ be the subring of all endomorphisms $f: \mathrm{id}\to \mathrm{id}$ commuting with Hecke operators, in the sense that for all $V\in \mathrm{Rep}(\hat{G}^I)$ and $A\in \mathcal D_\lis(\Bun_G,\Lambda)$, one has $T_V(f(A))=f(T_V(A))\in \mathrm{End}(T_V(A))$.
\item[{\rm (iii)}] The spectral Bernstein center of $G$ is
\[
\mathcal Z^{\mathrm{spec}}(G,\Lambda) = \mathcal O(Z^1(W_E,\hat{G})_\Lambda)^{\hat{G}},
\]
the ring of global functions on the quotient stack $Z^1(W_E,\hat{G})_\Lambda/\hat{G}$.
\end{altenumerate}
\end{definition}

The inclusion $\mathcal D(G(E),\Lambda)\hookrightarrow \mathcal D_\lis(\Bun_G,\Lambda)$ induces a map of algebra $\mathcal Z^{\mathrm{geom}}(G,\Lambda)\to \mathcal Z(G(E),\Lambda)$.

Now the construction of excursion operators, together with Theorem~\ref{thm:verygoodprimeintro} imply the following. Here $\Lambda$ is a $\mathbb Z_\ell[\sqrt{q}]$-algebra such that the order of $\pi_0 Z(G)$ is invertible in $\Lambda$.

\begin{proposition} There is a canonical map
\[
\mathcal Z^{\mathrm{spec}}(G,\Lambda)\to \mathcal Z^{\mathrm{geom}}_{\mathrm{Hecke}}(G,\Lambda)\subset \mathcal Z^{\mathrm{geom}}(G,\Lambda),
\]
and in particular a map
\[
\Psi_G: \mathcal Z^{\mathrm{spec}}(G,\Lambda)\to\mathcal Z(G(E),\Lambda).
\]
\end{proposition}

The construction of $L$-parameters above is then a consequence of this map on Bernstein centers. The existence of such an integral map is due to Helm--Moss \cite{HelmMoss} in the case $G=\GL_n$.

\begin{remark} In the function-field case, a similar construction has been given by Genestier--Lafforgue \cite{GenestierLafforgue}. Li-Huerta \cite{LiHuertaFSGL} has proved that these constructions agree.
\end{remark}

We make the following conjecture regarding independence of $\ell$. For its formulation, we note that there is a natural $\mathbb Q$-algebra $\mathcal Z^{\mathrm{spec}}(G,\mathbb Q)$ whose base change to $\mathbb Q_\ell$ is $\mathcal Z^{\mathrm{spec}}(G,\mathbb Q_\ell)$ for any $\ell\neq p$; in fact, one can take the global functions on the stack of $L$-parameters that are continuous for the discrete topology (i.e.~trivial on an open subgroup of $W_E$); see also \cite{DatHelmKurinczukMoss}.

\begin{conj}\label{conj:independenceofell} There is a (necessarily unique) map $\mathcal Z^{\mathrm{spec}}(G,\mathbb Q(\sqrt{q}))\to \mathcal Z(G(E),\mathbb Q(\sqrt{q}))$ that after base extension to any $\mathbb Q_\ell$ for $\ell\neq p$ recovers the composite
\[
\mathcal Z^{\mathrm{spec}}(G,\mathbb Q_\ell(\sqrt{q}))\to \mathcal Z^{\mathrm{geom}}(G,\mathbb Q_\ell(\sqrt{q}))\to \mathcal Z(G(E),\mathbb Q_\ell(\sqrt{q})).
\]
\end{conj}

This would ensure that the $L$-parameters we construct are independent of $\ell$ in the relevant sense. Further conjectures about this map and its relation to the stable Bernstein center have been formulated by Haines \cite{HainesStable} (see also \cite{BezrukavnikovKazhdanVarshavsky}, \cite[Section 6]{ScholzeShin}). In particular, it is conjectured that for $G$ quasisplit, the map $\Psi_G$ is injective, and its image can be characterized as those elements of the Bernstein center of $G(E)$ whose corresponding distribution is invariant under stable conjugation.

One can also construct the map to the Bernstein center in terms of moduli spaces of local shtukas, as follows. For simplicity, we discuss this again only if $\Lambda$ is a $\mathbb Z/\ell^n$-algebra for some $n$. Given $I$ and $V$ as above, we can consider a variant $\Sht_{(G,1,V),K}$ of the spaces $\Sht_{(G,b,\leq \mu_\bullet),K}$ considered above, where the bound is given by the support of $V$ and we fix the element $b=1$. They come with an \'etale period map
\[
\pi_K: \Sht_{(G,1,V),K}\to \Gr^\tw_{G,\prod_{i=1}^n \Spd \breve{E}}
\]
and a perverse sheaf $\mathcal S_V$. When restricted to the geometric diagonal
\[
\overline{x}: \Spd \hat{\overline{E}}\to \prod_{i=1}^n \Spd \breve{E},
\]
they become a corresponding moduli space of shtukas with one leg
\[
f_K^\Delta: \Sht_{(G,1,V|_{\hat{G}})}\to \Spd \hat{\overline{E}}
\]
with the sheaf $\mathcal S_{V|_{\hat{G}}}$. The sheaf $\mathcal S_{V|_{\hat{G}}}$ admits maps $\alpha$ (resp.~$\beta$) from (resp.~to) the sheaf $i_\ast \Lambda$, where $i: G(E)/K=\Sht_{(G,1,\overline{\mathbb Q}_\ell),K}\hookrightarrow \Sht_{(G,1,V|_{\hat{G}}),K}$ is the subspace of shtukas with no legs. This produces an endomorphism
\[
c\text-\mathrm{Ind}_K^{G(E)} \Lambda\xrightarrow{\alpha} Rf^\Delta_{K!} \mathcal S_{V|_{\hat{G}}} = (Rf_{K!} \mathcal S_V)_{\overline{x}}\xrightarrow{(\gamma_i)_{i\in I}} (Rf_{K!} \mathcal S_V)_{\overline{x}}= Rf^\Delta_{K!} \mathcal S_{V|_{\hat{G}}}\xrightarrow{\beta} c\text-\mathrm{Ind}_K^{G(E)} \Lambda.
\]
Here, the action of $(\gamma_i)_{i\in I}$ is defined by Corollary~\ref{cor:finitenessmain}. It follows from the definitions that this is precisely the previous construction applied to the representation $c\text-\mathrm{Ind}_K^{G(E)} \Lambda$. Note that these endomorphisms are $G(E)$-equivariant, so define elements in the Hecke algebra
\[
\Lambda[K\backslash G(E)/K] = \mathrm{End}_{G(E)}(c\text-\mathrm{Ind}_K^{G(E)}\Lambda);
\]
in fact, these elements are central (as follows by comparison to the previous construction). Taking the inverse limit over $K$, one gets the elements in the Bernstein center of $G(E)$.\footnote{When the second author gave his Berkeley lectures \cite{Berkeley}, this was the construction of excursion operators that we envisaged. Note that a key step here is that the cohomology of moduli spaces of local shtukas defines a local system on $(\Div^1)^I$. It is however not clear how to prove this purely in terms of moduli spaces of shtukas. In the global function field case, this result has been obtained by Xue \cite{XueSmoothness}.}

Concerning the $L$-parameters we construct, we can prove the following basic results.

\begin{theorem}\leavevmode
\begin{altenumerate}
\item[{\rm (i)}] If $G=T$ is a torus, then $\pi\mapsto \varphi_\pi$ is the usual Langlands correspondence.
\item[{\rm (ii)}] The correspondence $\pi\mapsto \varphi_\pi$ is compatible with twisting.
\item[{\rm (iii)}] The correspondence $\pi\mapsto \varphi_\pi$ is compatible with central characters (cf.~\cite[10.1]{Borel79}).
\item[{\rm (iv)}] The correspondence $\pi\mapsto \varphi_\pi$ is compatible with passage to congradients (cf.~\cite{AdamsVoganContragredient}).
\item[{\rm (v)}] If $G'\to G$ is a map of reductive groups inducing an isomorphism of adjoint groups, $\pi$ is an irreducible smooth representation of $G(E)$ and $\pi'$ is an irreducible constitutent of $\pi|_{G'(E)}$, then $\varphi_{\pi'}$ is the image of $\varphi_\pi$ under the induced map $\hat{G}\to \hat{G'}$.
\item[{\rm (vi)}] If $G=G_1\times G_2$ is a product of two groups and $\pi$ is an irreducible smooth representation of $G(E)$, then $\pi=\pi_1\boxtimes \pi_2$ for irreducible smooth representations $\pi_i$ of $G_i(E)$, and $\varphi_\pi=\varphi_{\pi_1}\times \varphi_{\pi_2}$ under $\hat{G}=\hat{G}_1\times \hat{G}_2$.
\item[{\rm (vii)}] If $G=\mathrm{Res}_{E'|E} G'$ is the Weil restriction of scalars of a reductive group $G'$ over some finite separable extension $E'|E$, so that $G(E)=G'(E')$, then $L$-parameters for $G|E$ agree with $L$-parameters for $G'|E'$.
\item[{\rm (viii)}] The correspondence $\pi\mapsto \varphi_\pi$ is compatible with parabolic induction.
\item[{\rm (ix)}] For $G=\GL_n$ and supercuspidal $\pi$, the correspondence $\pi\mapsto \varphi_\pi$ agrees with the usual local Langlands correspondence \cite{LaumonRapoportStuhler}, \cite{HarrisTaylor}, \cite{HenniartGLn}.
\end{altenumerate}
\end{theorem}

Note that parts (viii) and (ix) together say that for $\GL_n$ and general $\pi$, the $L$-parameter $\varphi_\pi$ is what is usually called the semisimple $L$-parameter.

\section{The spectral action}\label{sec:introspectral}

The categorical structure we have constructed actually produces something better. Let $\Lambda$ be the ring of integers in a finite extension of $\mathbb Q_\ell(\sqrt{q})$. We have the stable $\infty$-category $\mathcal C=\mathcal D_\lis(\Bun_G,\Lambda)^\omega$ of compact objects, which is linear over $\Lambda$, and functorially in the finite set $I$ an exact monoidal functor $\Rep_\Lambda(\hat{G}\rtimes Q)^I\to \End_\Lambda(\mathcal C)^{BW_E^I}$ that is linear over $\Rep_\Lambda(Q^I)$; here, $\End_\Lambda(\mathcal C)$ denotes the stable $\infty$-category of $\Lambda$-linear endofunctors of $\mathcal C$, and we regard it as being enriched in condensed $\Lambda$-modules via regarding $\mathcal C$ as enriched in relatively discrete condensed $\Lambda$-modules. A first version of the following theorem is due to Nadler--Yun \cite{NadlerYun} in the context of Betti geometric Langlands, and a more general version appeared in the work of Gaitsgory--Kazhdan--Rozenblyum--Varshavsky \cite{GaitsgoryKazhdanRozenblyumVarshavsky}. Both references, however, effectively assume that $G$ is split, work only with characteristic $0$ coefficients, and work with a discrete group in place of $W_E$. At least the extension to $\mathbb Z_\ell$-coefficients is a nontrivial matter.

Note that $Z^1(W_E,\hat{G})$ is not quasicompact, as it has infinitely many connected components; it can be written as the increasing union of open and closed quasicompact subschemes $Z^1(W_E/P,\hat{G})$. We say that an action of $\Perf(Z^1(W_E,\hat{G})/\hat{G})$ on a stable $\infty$-category $\mathcal C$ is compactly supported if for all $X\in \mathcal C$ the functor $\Perf(Z^1(W_E,\hat{G})/\hat{G})\to \mathcal C$ (induced by acting on $X$) factors over some $\Perf(Z^1(W_E/P,\hat{G})/\hat{G})$.

\begin{theorem} Assume that $\ell$ does not divide the order of $\pi_0 Z(G)$. Let $\mathcal C$ be a small $\Lambda$-linear stable $\infty$-category. Then giving, functorially in the finite set $I$, an exact $\Rep_\Lambda(Q^I)$-linear monoidal functor
\[
\Rep(\hat{G}\rtimes Q)^I\to \End_\Lambda(\mathcal C)^{BW_E^I}
\]
is equivalent to giving a compactly supported $\Lambda$-linear action of
\[
\Perf(Z^1(W_E,\hat{G})_\Lambda/\hat{G}).
\]
Here, given a compactly supported $\Lambda$-linear action of $\Perf(Z^1(W_E,\hat{G})_\Lambda/\hat{G})$, one can produce such an exact $\Rep_\Lambda(Q^I)$-linear monoidal functor
\[
\Rep_\Lambda(\hat{G}\rtimes Q)^I\to \End_\Lambda(\mathcal C)^{BW_E^I}
\]
functorially in $I$ by composing the exact $\Rep_\Lambda(Q^I)$-linear symmetric monoidal functor
\[
\Rep(\hat{G}\rtimes Q)^I\to \Perf(Z^1(W_E,\hat{G})_\Lambda/\hat{G})^{BW_E^I}
\]
with the action of $\Perf(Z^1(W_E,\hat{G})_\Lambda/\hat{G})$.

The same result holds true with $\Lambda$ a field over $\mathbb Q_\ell$, for any prime $\ell$.
\end{theorem}

Here, the exact $\Rep_\Lambda(Q^I)$-linear symmetric monoidal functor
\[
\Rep_\Lambda(\hat{G}\rtimes Q)^I\to \Perf(Z^1(W_E,\hat{G})_\Lambda/\hat{G})^{BW_E^I}
\]
is induced by tensor products and the exact $\Rep_\Lambda(Q)$-linear symmetric monoidal functor
\[
\Rep_\Lambda(\hat{G}\rtimes Q)\to \Perf(Z^1(W_E,\hat{G})_\Lambda/\hat{G})^{BW_E}
\]
corresponding to the universal $\hat{G}\rtimes Q$-torsor, with the universal $W_E$-equivariance as parametrized by $Z^1(W_E,\hat{G})/\hat{G}$.

The key part of the proof is actually the final part of Theorem~\ref{thm:verygoodprimeintro} above, which effectively describes $\Perf(Z^1(W_E/P,\hat{G})/\hat{G})$ in terms of generators and relations, as does the present theorem.

In particular, we get an action of $\Perf(Z^1(W_E,\hat{G})_\Lambda/\hat{G})$ on $\mathcal D_\lis(\Bun_G,\Lambda)$, suitably compatible with the Hecke action.

With everything in place, it is now obvious that the main conjecture is the following, cf.~\cite{ArinkinGaitsgory}, \cite{BenZviChenHelmNadler}, \cite{ZhuConjectures}, \cite{HellmannLLC}:\footnote{The previous version of this manuscript made a seemingly less precise conjecture by asking for the existence of a functor instead of noting that it must necessarily be realized as a right adjoint and hence is unique if it exists. This uniqueness was pointed out to us in particular by Hansen \cite{HansenBeijing}.}

\begin{conj} Assume that $G$ is quasisplit and choose Whittaker data consisting of a Borel $B\subset G$ and generic character $\psi: U(E)\to \mathcal O_L^\times$ of the unipotent radical $U\subset B$, where $L/\mathbb Q_\ell$ is some algebraic extension; also fix $\sqrt{q}\in \mathcal O_L$. Let $n$ be the order of $\pi_0 Z(G)$ and let $\Lambda=\mathcal O_L[\tfrac 1n]$. Let
\[
\mathcal W_\psi\in \mathcal D_\lis(\Bun_G,\Lambda)
\]
be the Whittaker sheaf, which is the sheaf concentrated on $\Bun_G^1$ corresponding to the Whittaker representation $c\text-\mathrm{Ind}_{U(E)}^{G(E)} \psi$, and let
\[
\Ind\Perf^{\mathrm{qc}}(Z^1(W_E,\hat{G})_\Lambda/\hat{G})\to \mathcal D_\lis(\Bun_G,\Lambda): M\mapsto \mathrm{Act}_M(\mathcal W_\psi)
\]
be defined as the colimit-preserving extension of the spectral action $\mathrm{Act}$ on $\mathcal W_\psi$. Then the corresponding right adjoint functor is fully faithful when restricted to the compact objects, and induces an equivalence of ($\Perf(Z^1(W_E,\hat{G})_\Lambda/\hat{G})$-linear small stable) $\infty$-categories
\[
\mathcal D(\Bun_G,\Lambda)^\omega\cong \mathcal D^{b,\mathrm{qc}}_{\mathrm{coh},\mathrm{Nilp}}(Z^1(W_E,\hat{G})_\Lambda/\hat{G}).
\]
\end{conj}

We inverted the order $n$ of $\pi_0 Z(G)$ here, because only then the spectral action has been constructed. We are not sure what to expect without inverting $n$. (In fact, we would not be surprised if the notion of ``nilpotent singular support'' that we use has to be modified at bad primes.)

Here, we use the notion of complexes of coherent sheaves with nilpotent singular support, see \cite{ArinkinGaitsgory}. More precisely, $\mathcal D^{b,\mathrm{qc}}_{\mathrm{coh},\mathrm{Nilp}}$ is the $\infty$-category of bounded complexes with quasicompact support, coherent cohomology, and nilpotent singular support. With characteristic $0$ coefficients, or at banal primes $\ell$, the condition of nilpotent singular support is actually automatic.

If $\mathcal{W}_{\psi}$ is the Whittaker sheaf and we note $\ast$  the spectral action, the conjecture thus says that 
\begin{align*}
\Perf^{\mathrm{qc}}(Z^1(W_E,\hat{G})_\Lambda/\hat{G}) & \lto  \mathcal D(\Bun_G,\Lambda)\\ 
M & \longmapsto M\ast \mathcal{W}_{\psi} 
\end{align*}
is fully faithful and extends to an equivalence of stable $\infty$-categories
$$
\mathcal D^{b,\mathrm{qc}}_{\mathrm{coh},\mathrm{Nilp}}(Z^1(W_E,\hat{G})_\Lambda/ \hat{G})\cong \mathcal D(\Bun_G,\Lambda)^\omega.
$$

Recall that the right-hand side contains $\mathcal D(G(E),\Lambda)^\omega$ fully faithfully, so in particular this $\infty$-category should embed fully faithfully into the left-hand side. This has been conjectured by Hellmann in \cite{HellmannLLC} and Ben-Zvi--Chen--Helm--Nadler \cite{BenZviChenHelmNadler} have proved parts of this (they use $\overline{\mathbb Q}_\ell$-coefficients, and work with split groups and the Bernstein component corresponding to representations with Iwahori fixed vector).

\begin{remark} Consider the conjecture with coefficients in $\overline{\mathbb Q}_\ell$. Ideally, the conjecture should also include a comparison of $t$-structures. Unfortunately, we did not immediately see a good candidate for matching $t$-structures. Ideally, this would compare the perverse $t$-structure on the left (which is well-defined, for abstract reasons, and appears at least implicitly in \cite{CaraianiScholze}, \cite{CaraianiScholze2}; it seems to be the ``correct'' $t$-structure for questions of local-global compatibility) with some ``perverse-coherent'' $t$-structure on the right. If so, the equivalence would also yield a bijection between irreducible objects in the abelian hearts. On the left-hand side, these irreducible objects would then be enumerated by pairs $(b,\pi_b)$ of an element $b\in B(G)$ and an irreducible smooth representation $\pi_b$ of $G_b(E)$, by using intermediate extensions. On the right-hand side, they would likely correspond to a Frobenius-semisimple $L$-parameter $\phi: W_E\to \hat{G}(\overline{\mathbb Q}_\ell)$ together with an irreducible representation of the centralizer $S_\phi$ of $\phi$. Independently of the categorical conjecture, one can wonder whether these two sets are in fact canonically in bijection.\footnote{This question has been answered affirmatively by Bertoloni Meli--Oi \cite{BertoloniMeliOi}. Hansen \cite{HansenBeijing} has moreover made progress in understanding possibly matching $t$-structures, by introducing the hadal $t$-structure on $\Bun_G$.}
\end{remark}

\section{The origin of the ideas}\label{sec:introorigin}

Finally, let us give some account of the historical developments of these ideas, from our own biased perspective. Let us first recall some of our early work in the direction of local Langlands correspondences. Fargues \cite{FarguesMantovan} has proved that in the cohomology of basic Rapoport--Zink spaces for $\GL_n$ (and $U(3)$) and general minuscule cocharacters, an appropriate version of the local Langlands correspondence is realized. Moreover, Fargues \cite{FarguesIsomorphism} has proved the duality isomorphism between the Lubin--Tate and Drinfeld tower. Already at this point Fargues thought of this as an attempt to geometrize the Jacquet-Langlands correspondence, see \cite[Theorem 2 of the Pr\'eambule]{FarguesIsomorphism}. On the other hand, Scholze \cite{ScholzeLLC} has given a new proof of the local Langlands correspondence for $\GL_n$. His results pointed to the idea that there ought to exist certain sheaves on the moduli stack of $p$-divisible groups (which, when restricted to perfect schemes, can be regarded as a ``part'' of the stack $\GL_n\text-\mathcal{I}\mathrm{soc}$ considered above), giving a certain geometrization of the local Langlands correspondence, then formulated as a certain character sheaf property (inspired by the character formulas in \cite{ScholzeLLC}). Related observations were also made by Boyer (cf.~e.g.~\cite{Boyer}) and in unpublished work of Dat. However, Scholze was always uneasy with the very bad geometric properties of the stack of $p$-divisible groups.

At this point, both of us had essentially left behind local Langlands to study other questions. Fargues found the fundamental curve of $p$-adic Hodge theory in his work with Fontaine \cite{FarguesFontaine}; an initial critical motivation for Fargues was a development of ``$p$-adic Hodge theory without Galois actions'', i.e.~for fields like $\mathbb C_p$. Indeed, this was required in some of his work on Rapoport--Zink spaces. On the other hand, Scholze developed perfectoid spaces \cite{ScholzePerfectoidSpaces}, motivated by the weight-monodromy conjecture. After his talk at a conference in Princeton in March 2011, Weinstein gave a talk about his results on the Lubin--Tate tower at infinite level, which made it clear that it is in fact a perfectoid space. Scholze at the time was already eager to understand the isomorphism between Lubin--Tate and Drinfeld tower, and it now became clear that it should really be an isomorphism of perfectoid spaces. This was worked out in \cite{ScholzeWeinstein}. At the time of writing of \cite{ScholzeWeinstein}, the perspective of the Fargues--Fontaine curve had already become central, and we realized that the isomorphism of the towers simply amounts to two dual descriptions of the space of minuscule modifications $\mathcal O_X^n\to \mathcal O_X(\tfrac 1n)$ on the Fargues--Fontaine curve, depending on which bundle is fixed and which one is the modification. This was the first clear connection between local Langlands (as encoded in the cohomology of Lubin--Tate and Drinfeld space) and the theory of vector bundles on the Fargues--Fontaine curve, which Scholze had however not taken seriously enough. Moreover, Fargues had noted in \cite{FarguesFontaine}, in the proof of ``weakly admissible implies admissible", that modifications of vector bundles were playing an important role: the Hodge filtration of a filtered $\phi$-module allows one to define a new vector bundle by modifying the vector bundle associated to an isocrystal i.e.~by ``applying a Hecke correspondence'' as he said in the talk \cite{FarguesExposeFontaine} at the conference in honor of Jean-Marc Fontaine.

This duality perspective also put the two dual period morphisms into the center of attention: The Hodge--de Rham period mapping, and the Hodge--Tate period mapping (which are swapped under the duality isomorphism). Thinking about the Lubin--Tate tower as part of the moduli space of elliptic curves, Scholze then realized that the Hodge--Tate period map even exists globally on the moduli space of elliptic curves with infinite level (on the level of Berkovich topological spaces, this had also been observed by Fargues before). Moreover, Scholze realized that the Hodge--Tate period map gives a substitute for the map from the moduli space of elliptic curves to the moduli space of $p$-divisible groups, and that the sheaves he sought for a geometric interpretation of \cite{ScholzeLLC} have a better chance of existing on the target of the Hodge-Tate period map, which is simply a projective space over $\mathbb C_p$; he sketched these ideas in an MSRI talk \cite{ScholzeMSRISheaves}. (Again, Dat has had similar ideas.) Eventually, this perspective was used in his work with Caraiani \cite{CaraianiScholze}, \cite{CaraianiScholze2} to study torsion in the cohomology of Shimura varieties. The work with Caraiani required the classification of $G$-torsors on the Fargues--Fontaine curve, which was proved by Fargues \cite{FarguesGBun}.

Increasingly taking the perspective of studying all geometric objects by mapping only perfectoid spaces in, the idea of diamonds emerged quickly, including the possibility of getting several copies of $\Spec \mathbb Q_p$ (the earliest published incarnation of this idea is \cite{WeinsteinGalois}), and of defining general moduli spaces of $p$-adic shtukas. These ideas were laid out in Scholze's Berkeley course \cite{Berkeley} during the MSRI trimester in Fall 2014. The eventual goal was always to adapt V.~Lafforgue's work \cite{VLafforgue1} to the case of $p$-adic fields; the original strategy was to define the desired excursion operators via the cohomology of moduli spaces of local shtukas. At the beginning of the trimester, Scholze was still very wary about the geometric Langlands program, as it did not seem to be able to incorporate the subtle arithmetic properties of supercuspidal representations of $p$-adic groups. It was thus a completely unexpected conceptual leap that in fact the best perspective for the whole subject is to view the local Langlands correspondence as a geometric Langlands correspondence on the Fargues--Fontaine curve, which Fargues suggested over a coffee break at MSRI (partly inspired by having thought intensely about the space of $G$-bundles on the curve in relation to \cite{FarguesGBun}). Fargues was taking the perspective of Hecke eigensheaves then, seeking to construct for any (discrete) $L$-parameter $\phi$ an associated Hecke eigensheaf $A_\phi$ on $\Bun_G$ with eigenvalue $\phi$. This should define a functor $\phi\mapsto A_\phi$, and thus carry an action of the centralizer group $S_\phi\subset \hat{G}$ of $\phi$, and the corresponding $S_\phi$-isotypic decomposition of $A_\phi$ should realize the internal structures of the $L$-packets. Moreover, the Hecke eigensheaf property should imply the Kottwitz conjecture \cite[Conjecture 7.3]{RapoportViehmann} on the cohomology of local Shimura varieties. This made everything come together. In particular, it gave a compelling geometric origin for the internal structure of $L$-packets, and also matched the recent work of Kaletha \cite{KalethaLocalLanglands} who used basic $G$-isocrystals for the fine study of $L$-packets.

Unfortunately, the conjecture was formulated on extremely shaky grounds: It presumed that one could work with the moduli stack $\Bun_G$ as if it were an object of usual algebraic geometry. Of course, it also presumed that there is a version of geometric Satake, etc.pp. On the other hand, we realized that once we could merely \emph{formulate} Fargues' conjecture, enough machinery is available to apply Lafforgue's ideas \cite{VLafforgue1} to get the ``automorphic-to-Galois'' direction and define (semisimple) $L$-parameters (as Genestier--Lafforgue \cite{GenestierLafforgue} did in equal characteristic).

Since then, it has been a long and very painful process. The first step was to give a good definition of the category of geometric objects relevant to this picture, i.e.~diamonds. In particular, one had to prove that the relevant affine Grassmannians have this property. This was the main result of the Berkeley course \cite{Berkeley}. For the proof, the concept of v-sheaves was introduced, which has since taken on a life of its own also in algebraic geometry (cf.~\cite{BhattMathew}). (Generally, v-descent turned out to be an extremely powerful proof technique. We use it here to reprove the basic theorems about the Fargues--Fontaine curve, recovering the main theorems of \cite{FarguesFontaine} and \cite{KedlayaLiu1} with little effort.) Next, one had to develop a $6$-functor formalism for the \'etale cohomology of diamonds, which was achieved in \cite{ECoD}, at least with torsion coefficients. The passage to $\mathbb Q_\ell$-coefficients requires more effort than for schemes, and we will comment on it below. A central technique of \cite{ECoD} is pro-\'etale descent, and more generally v-descent. In fact, virtually all theorems of \cite{ECoD} are proved using such descent techniques, essentially reducing them to profinite collections of geometric points. It came as a surprise to Scholze that this process of disassembling smooth spaces into profinite sets has any power in proving geometric results, and this realization gave a big impetus to the development of condensed mathematics (which in turn fueled back into the present project).

At this point, it became possible to contemplate Fargues' conjecture. In this respect, the first result that had to be established is that $D_\et(\Bun_G,\mathbb Z/\ell^n\mathbb Z)$ is well-behaved, for example satisfies Verdier biduality for ``admissible'' sheaves. We found a proof, contingent on the cohomological smoothness of a certain ``chart'' $\pi_b: \mathcal M_b\to \Bun_G$ for $\Bun_G$ near any $b\in B(G)$; this was explained in Scholze's IH\'ES course \cite{ScholzeIHES2}. While for $G=\GL_n$, the cohomological smoothness of $\pi_b$ could be proved by a direct attack, in general we could only formulate it as a special case of a general ``Jacobian criterion of smoothness'' for spaces parametrizing sections of $Z\to X_S$ for some smooth adic space $Z$ over the Fargues--Fontaine curve. Proving this Jacobian criterion required three further key ideas. The first is the notion of ``formal smoothness'', where liftings to infinitesimal thickenings (that do not exist in perfectoid geometry) are replaced by liftings to actual small open (or \'etale) neighborhoods. The resulting notion is closely related to the notion of absolute neighborhood retracts in classical topology \cite{BorsukRetractLivre}. Through some actual ``analysis'', it is not hard to prove that the space of sections is formally smooth. Unfortunately, this does not seem to be enough to guarantee cohomological smoothness. The first issue is that formal smoothness does not imply any finite-dimensionality. Here, the second key idea comes in, which is Bhatt's realization \cite{BhattScholzePrism} that Zariski closed immersions are strongly Zariski closed in the sense of \cite[Section II.2]{ScholzeTorsion} (contrary to a claim made by Scholze there). At this point, it would be enough to show that spaces that are formally smooth and Zariski closed in a finite-dimensional perfectoid ball are cohomologically smooth. Unfortunately, despite many tries, we are still unable to prove that even the different notions of dimension of \cite{ECoD} (Krull dimension, $\dimtrg$, cohomological dimension) agree for such spaces. This may well be the most important foundational open problem in the theory:

\begin{problem} Let $X\subset \tilde{\mathbb B}^n_C$ be Zariski closed, where $\tilde{\mathbb B}^n$ is a perfectoid ball. Show that $X$ has a well-behaved dimension.
\end{problem}

In fact, we find it crazy that we are able to prove all sorts of nontrivial geometric results without ever being able to unambiguously talk about dimensions!

Our attacks on this failing, a third key idea comes in: Namely, the notion of universally locally acyclic sheaves, that we also developed independently in order to prove geometric Satake. It is easy to see that formal smoothness plus finite-dimensionality implies that the constant sheaf is universally locally acyclic; it remains to see that the dualizing sheaf is invertible. This can be proved by a deformation to the normal cone (using universal local acyclicity to spread the result on the normal cone to a neighborhood). We found this argument at a conference in Luminy in July 2018; an inspiration to use a deformation to the normal cone may have been Clausen's use in the proof of the ``linearization hypothesis''.

These results are enough to show that $D_\et(\Bun_G,\mathbb Z/\ell^n\mathbb Z)$ is well-behaved, and are already enough to prove new finiteness results on the cohomology of Rapoport--Zink spaces (with torsion coefficients). Our next emphasis was on geometric Satake. This essentially required the theory of universally locally acyclic sheaves, and a version of Braden's hyperbolic localization theorem \cite{Braden}. We were able to find substitutes for both. Regarding universally locally acyclic sheaves, we were able to prove analogues of most basic theorems, however we failed to prove that in general they are preserved under relative Verdier duality (even while we could check it by hand in all relevant cases). Lu--Zheng \cite{LuZhengDuality} then found a new characterization of universally locally acyclic sheaves, making stability under relative Verdier duality immediate. Their arguments immediately transport to our setting. Eventually we used a slightly different characterization, but in spirit the argument is still the same as theirs. Regarding hyperbolic localization, we could not follow Braden's arguments that rely on nice coordinate choices. Instead, we reduce all arguments to the following (simple to prove) principle: If $X$ is a (partially proper) space with a $\mathbb G_m$-action such that $[X/\mathbb G_m]$ is qcqs, and $A\in D_\et([X/\mathbb G_m],\Lambda)$, then the partially compactly supported cohomology of $X$ with coefficients in $A$ vanishes. The idea here is that the $\mathbb G_m$-action contracts $X$ towards one of the ends. Afterwards, the proof of geometric Satake largely follows the lines of \cite{MirkovicVilonen}, although there are certain improvements in the argument; in particular, we give a simple reduction to groups of rank $1$, and pin the isomorphism with the dual group.

Using these results, one has all ingredients in place, but only working with torsion coefficients. One can formally pass to $\ell$-adically complete sheaves, but this leads to studying representations on Banach $\mathbb Q_\ell$-vector spaces, which is very unnatural. During this time, Clausen came to Bonn, and Clausen and Scholze started to develop condensed mathematics, and the theory of solid modules \cite{Condensed}. They realized that one could also define solid $\mathbb Z_\ell$-sheaves on schemes or diamonds, and that this makes it possible to study representations on discrete $\mathbb Q_\ell$- or $\overline{\mathbb Q}_\ell$-vector spaces, as desired. We take this up here, and first define solid $\mathbb Z_\ell$-sheaves on any small v-stack, together with some $5$-functor formalism (involving relative homology in place of compactly supported cohomology; its right adjoint is then pullback, so there are only $5$ functors), and afterwards pass to a certain subcategory of ``lisse-\'etale'' sheaves to define the desired category $D_\lis(\Bun_G,\overline{\mathbb Q}_\ell)$, with exactly the desired properties.

In the meantime, there was related work in the geometric Langlands program by Nadler--Yun \cite{NadlerYun} and Gaitsgory--Kazhdan--Rozenblyum--Varshavsky \cite{GaitsgoryKazhdanRozenblyumVarshavsky} that implied that the categorical structures we have now constructed --- $D(\Bun_G,\mathbb Q_\ell)$ together with the action of Hecke operators --- formally induce an action of the category of perfect complexes on the stack of $L$-parameters on $D(\Bun_G,\mathbb Q_\ell)$, giving a categorical upgrade to the construction of $L$-parameters based on excursion operators. (We were aware of some weak form of this, when restricted to elliptic parameters; this was discussed in the last lecture of \cite{ScholzeIHES2}, based on some unpublished results of Ansch\"utz.) Here, we make the effort of proving a result with $\mathbb Z_\ell$-coefficients.

\section{Acknowledgments}

We apologize for the long delay in the preparation of the manuscript. Already in April 2016, an Oberwolfach Arbeitsgemeinschaft gave an introduction to these ideas. The required foundations on \'etale cohomology of diamonds (with torsion coefficients) were essentially worked out in the winter term 2016/17, where also the ARGOS seminar in Bonn studied the evolving manuscript \cite{ECoD}, which was an enormous help. Next, we gave the Hadamard course at the IH\'ES in Spring 2017, \cite{ScholzeIHES2}. At the time we had essentially obtained enough results to achieve the construction of $L$-parameters; the missing ingredient was the Jacobian criterion. In particular, we understood that one has to use the local charts $\mathcal M_b$ to study $\Bun_G$. In more detail, these results were discussed in the ARGOS seminar in the summer term of 2017, and the proof of geometric Satake in the ARGOS seminar in the summer term of 2018. We heartily thank all the participants of these seminars for working through these manuscripts, and their valuable feedback. In the meantime, in the spring of 2018, we revised the Berkeley notes \cite{Berkeley} and used \cite{ECoD} to give a full construction of moduli spaces of $p$-adic shtukas. Then we eventually found a full proof of the Jacobian criterion in July 2018 at a conference in Luminy. Our early attempts at writing up the results ran into the issue that we wanted to first prove basic sanity results on dimensions. Later, we got sidetracked by the development of condensed mathematics with Dustin Clausen, which eventually made it possible to use $\overline{\mathbb Q}_\ell$-coefficients in this paper, finally giving the most natural results. In fact, we always wanted to develop the foundations so that papers such as the paper of Hansen--Kaletha--Weinstein \cite{HansenKalethaWeinstein} on the Kottwitz conjecture can build on a proper foundation, and prove results about general smooth $\overline{\mathbb Q}_\ell$-representations; we believe we have finally achieved this. In the winter term 2020/21, Scholze gave an (online) course in Bonn about these results, and wants to thank all the participants for their valuable feedback. We thank the many mathematicians that we discussed with about these results, including Johannes Ansch\"utz, David Ben-Zvi, Joseph Bernstein, Bhargav Bhatt, Arthur-C\'esar le Bras, Ana Caraiani, Dustin Clausen, Pierre Colmez, Jean-Fran\c{c}ois Dat, Vladimir Drinfeld, Jessica Fintzen, Dennis Gaitsgory, Toby Gee, David Hansen, Eugen Hellmann, Alex Ivanov, Tasho Kaletha, Kiran Kedlaya, Dmitry Kubrak, Vincent Lafforgue, Judith Ludwig, Andreas Mihatsch, Sophie Morel, Jan Nekov\'a\v{r}, Wies\l awa Nizio\l, Vincent Pilloni, Michael Rapoport, Timo Richarz, Will Sawin, Tony Scholl, Antoine Touz\'e, Roman Travkin, Wilberd van der Kallen, Yakov Varshavsky, Jared Weinstein, Zhiyou Wu, Xinwen Zhu, Konrad Zou, and undoubtedly many others. Special thanks go to David Hansen for many comments on a preliminary draft. He also observed that the categorical equivalence is unique if it exists, an observation that we had not made in the original version of this manuscript. Konrad Zou observed that at bad primes, the nilpotent cone used in the formulation of the nilpotent singular support condition may not behave correctly. Moreover, we thank Toby Gee, David Helm, Sug Woo Shin and especially Zhiyu Zhang for comments and suggestions on a first version. Finally, we heartily thank the referees for a very careful reading, and many helpful suggestions.

During the preparation of this manuscript, Fargues was supported by the ANR grant ANR-14-CE25 PerCoLaTor and the ERC Advanced Grant 742608 GeoLocLang, and Scholze was supported by a DFG Leibniz Prize, and by the DFG under the Excellence Strategy – EXC-2047/1 – 390685813.

\section{Notation}

Throughout most of this paper, $E$ denotes a nonarchimedean local field with residue field $\mathbb F_q$ of characteristic $p>0$, and we fix an algebraic closure $k=\overline{\mathbb F}_q$ of $\mathbb F_q$. Then $\breve{E}$ is the completed unramified extension of $E$ with residue field $k$. We also fix a separable closure $\overline{E}$ of $E$, with absolute Galois group $\Gamma=\mathrm{Gal}(\overline{E}|E)$, containing the Weil group $W_E$, inertia subgroup $I_E$, and wild inertia $P_E$. The letter $P$ usually denotes open subgroups of $P_E$, but is occasionally also used to denote a finite $p$-group (or more generally finite solvable group of order prime to $\ell$).

The group $G$ is usually a reductive group over $E$; reductive groups are always assumed to be connected.

For any topological space $X$, we denote by $\underline{X}$ the sheaf taking any $S$ (in the relevant test category, usually a perfectoid space) to the continuous maps from $|S|$ to $X$. This is in the spirit of the passage from topological spaces to condensed sets, see \cite{Condensed}. We make occasional use of the condensed language, but do not make use of any nontrivial results from \cite{Condensed}. In particular, our discussion of solid $\ell$-adic sheaves is self-contained.

We will occasionally use the ``animated'' terminology, see \cite{Condensed}, \cite{CesnaviciusScholze}. In particular, we use the term anima for what is variously called spaces in \cite{LurieHTT}, $\infty$-groupoids, or homotopy types, and for any ring $A$, the $\infty$-category of animated $A$-algebras is the $\infty$-category obtained from simplicial $A$-algebras by inverting weak equivalences. Thus, animated $A$-algebras are freely generated under sifted colimits by polynomial algebras $A[X_1,\ldots,X_n]$.

If $\mathcal C$ is an ($\infty$-)category equipped with an action of a group $G$, we write $\mathcal C^{BG}$ for the ($\infty$-)category of $G$-equivariant objects in $\mathcal C$. Note that the data here is really a functor $BG\to \mathrm{Cat}_\infty$, and $\mathcal C^{BG}$ is by definition the limit of this diagram. (It would be more customary to write $\mathcal C^G$, but this leads to inconsistent notation.) Also, we often write classifying stacks as $\ast/G$ instead of $BG$ as the letter $B$ also denotes Borel subgroups (and we strongly prefer $\ast/B$ to $BB$), and appears in Kottwitz' set $B(G)$.

\chapter{The Fargues--Fontaine curve and vector bundles}

The goal of this chapter is to define the Fargues--Fontaine curve, in its various incarnations, and the category of vector bundles on the Fargues--Fontaine curve. Throughout this chapter, we fix a nonarchimedean local field $E$ with residue field $\Fq$ of characteristic $p$. We let $\mathcal O_E\subset E$ be the ring of integers, and $\pi$ a uniformizing element in $E$.

For any perfectoid space $S$ over $\Fq$, we introduce a curve $\mathcal Y_S$, to be thought of as the hypothetical product $S\times_{\Spa \Fq} \Spa \mathcal O_E$, together with an open subset $Y_S\subset \mathcal Y_S$ given by the locus where $\pi\neq 0$. This carries a Frobenius $\phi$ induced from the Frobenius on $S$, and $X_S$ is the quotient $Y_S/\phi^\Z$.

The first results concern the Fargues--Fontaine curve $X_C=X_S$ when $S=\Spa C$ for some complete algebraically closed nonarchimedean field $C|\Fq$. We define a notion of classical points of $X_C$ in that case; they form a subset of $|X_C|$. The basic finiteness properties of $X_C$ are summarized in the following result.

\begin{theorem}[Proposition~\ref{prop:classicalpointsmaximalspectrum}, Corollary~\ref{cor:principalidealdomains}, Definition/Proposition~\ref{def:classicalpointsFF}] The adic space $\mathcal Y_C$ is locally the adic spectrum $\Spa(B,B^+)$ where $B$ is a principal ideal domain; the classical points of $\Spa(B,B^+)\subset \mathcal Y_C$ are in bijection with the maximal ideals of $B$. For each classical point $x\in \mathcal Y_C$, the residue field of $x$ is an untilt $C^\sharp$ of $C$ over $\mathcal O_E$, and this induces a bijection of the classical points of $\mathcal Y_C$ with untilts $C^\sharp$ of $C$ over $\mathcal O_E$. A similar result holds true for $Y_C\subset \mathcal Y_C$, and the quotient $X_C=Y_C/\phi^\Z$.
\end{theorem}

In the equal characteristic case, this is an immediate consequence of $\mathcal Y_C = \mathbb D_C$ and classical results in rigid-analytic geometry. In the $p$-adic case, we use tilting to reduce to the equal characteristic case. More precisely, if $E$ is $p$-adic and $E_\infty$ is the completion of $E(\pi^{1/p^\infty})$, then $\mathcal Y_C\times_{\Spa \mathcal O_E} \Spa \mathcal O_{E_\infty}$ is perfectoid, with tilt given by a perfectoid open unit disc $\tilde{\mathbb D}_C$. The corresponding map $|\tilde{\mathbb D}_C|\to |\mathcal Y_C|$ induces a surjective map on classical points, see Proposition~\ref{prop:classicalpointstilting}. At one key turn, in order to understand Zariski closed subsets of $\mathcal Y_C$, we use the result that Zariski closed subspaces are invariant under tilting, to reduce to $\tilde{\mathbb D}_C$. More precisely, we recall the following result.

\begin{proposition}[{\cite[Section II.2]{ScholzeTorsion}, \cite[Remark 7.5]{BhattScholzePrism}, \cite[Definition 5.7, Theorem 5.8]{ECoD}}]\label{prop:zariskiclosedstronglyzariskiclosed} Let $S=\Spa(R,R^+)$ be an affinoid perfectoid space with tilt $S^\flat=\Spa(R^\flat,R^{\flat +})$. Then a closed subspace $Z\subset |S|$ is the vanishing locus of an ideal $I\subset R$ if and only if $Z\subset |S|\cong |S^\flat|$ is the vanishing locus of an ideal $J\subset R^\flat$. In that case, there is a universal perfectoid space $S_Z\to S$ such that $|S_Z|\to |S|$ factors over $Z$, and $S_Z=\Spa(T,T^+)$ is affinoid perfectoid with $|S_Z|\to Z$ a homeomorphism, $R\to T$ surjective, $R^+\to T^+$ almost surjective, and $T^+$ is the integral closure of $R^+$ in $T$.
\end{proposition}

A key result is the classification of vector bundles.

\begin{theorem}[Theorem~\ref{thm:classificationvectorbundles}] The functor from $\Isoc_E$ to vector bundles on $X_C$ induces a bijection on isomorphism classes. In particular, there is a unique stable vector bundle $\mathcal O_{X_C}(\lambda)$ of any slope $\lambda\in \mathbb Q$, and any vector bundle $\mathcal E$ can be written as a direct sum of stable bundles.
\end{theorem}

We give a new self-contained proof of the classification theorem, making critical use of the v-descent results for vector bundles obtained in \cite{ECoD} and \cite{Berkeley}, and basic results on the geometry of Banach--Colmez spaces established here.

Allowing general $S\in \Perf_{\overline{\mathbb F}_q}$, we define the moduli space of degree $1$ Cartier divisors as $\Div^1 = \Spd \breve{E}/\phi^\Z$. Given a map $S\to \Div^1$, one can define an associated closed Cartier divisor $D_S\subset X_S$; locally, this is given by an untilt $D_S=S^\sharp\subset X_S$ of $S$ over $E$, and this embeds $\Div^1$ into the space of closed Cartier divisors on $X_S$. Another important result is the following ampleness result, cf.~\cite[Proposition 6.2.4]{KedlayaLiu1}, which implies that one can define an algebraic version of the curve, admitting the same theory of vector bundles.

\begin{theorem}[Theorem~\ref{thm:O1ample}, Proposition~\ref{prop:ampleGAGA}, Proposition~\ref{prop:classicalpointsschematiccurve}] Assume that $S\in \Perf$ is affinoid. For any vector bundle $\mathcal E$ on $X_S$, the twist $\mathcal E(n)$ is globally generated and has no higher cohomology for all $n\gg 0$. Defining the graded ring
\[
P=\bigoplus_{n\geq 0} H^0(X_S,\mathcal O_{X_S}(n))
\]
and the scheme $X_S^{\mathrm{alg}} = \Proj P$, there is a natural map of locally ringed spaces $X_S\to X_S^{\mathrm{alg}}$, pullback along which defines an equivalence of categories of vector bundles, preserving cohomology.

If $S=\Spa C$ for some complete algebraically closed nonarchimedean field $C$, then $X_C^{\mathrm{alg}}$ is a regular noetherian scheme of Krull dimension $1$, locally the spectrum of a principal ideal domain, and its closed points are in bijection with the classical points of $X_C$.
\end{theorem}

We also need to understand families of vector bundles, i.e.~vector bundles $\mathcal E$ on $X_S$ for general $S$. Here, the main result is the following, which is originally due to Kedlaya--Liu \cite{KedlayaLiu1}.

\begin{theorem}[Theorem~\ref{thm:kedlayaliu}, Corollary~\ref{cor:Elocalsystems}] Let $S\in \Perf$ and let $\mathcal E$ be a vector bundle on $X_S$. Then the function taking a point $s\in S$ to the Harder--Narasimhan polygon of $\mathcal E|_{X_s}$ defines a semicontinuous function on $S$. If it is constant, then $\mathcal E$ admits a global Harder--Narasimhan stratification, and pro-\'etale locally on $S$ one can find an isomorphism with a direct sum of $\mathcal O_{X_S}(\lambda)$'s.

In particular, if $\mathcal E$ is everywhere semistable of slope $0$, then $\mathcal E$ is pro-\'etale locally trivial, and the category of such $\mathcal E$ is equivalent to the category of pro-\'etale $\underline{E}$-local systems on $S$.
\end{theorem}

The key to proving this theorem is the construction of certain global sections of $\mathcal E$. To achieve this, we use v-descent techniques, and an analysis of the spaces of global sections of $\mathcal E$; these are known as Banach--Colmez spaces, and were first introduced (in slightly different terms) by Colmez \cite{ColmezEspaces}; see also le Bras' thesis \cite{leBras}.

\begin{definition} Let $\mathcal E$ be a vector bundle on $X_S$. The Banach--Colmez space $\BC(\mathcal E)$ associated with $\mathcal E$ is the locally spatial diamond over $S$ whose $T$-valued points, for $T\in \Perf_S$, are given by
\[
\BC(\mathcal E)(T) = H^0(X_T,\mathcal E|_{X_T}).
\]
Similarly, if $\mathcal E$ is everywhere of only negative Harder--Narasimhan slopes, the negative Banach--Colmez space $\BC(\mathcal E[1])$ is the locally spatial diamond over $S$ whose $T$-valued points are
\[
\BC(\mathcal E[1])(T) = H^1(X_T,\mathcal E|_{X_T}).
\]
\end{definition}

Implicit here is that this functor actually defines a locally spatial diamond. For this, we calculate some key examples of Banach--Colmez spaces. For example, if $\mathcal E=\mathcal O_{X_S}(\lambda)$ with $0<\lambda\leq [E:\mathbb Q_p]$ (resp.~all positive $\lambda$ if $E$ is of equal characteristic), then $\BC(\mathcal E)$ is representable by a perfectoid open unit disc (of dimension given by the numerator of $\lambda$). A special case of this is the identification of $\BC(\mathcal O_{X_S}(1))$ with the universal cover of a Lubin--Tate formal group law, yielding a very close relation between Lubin--Tate theory, and thus local class field theory, and the Fargues--Fontaine curve. This case actually plays a special role in getting some of the theory started, and we recall it explicitly in Section~\ref{sec:LubinTate}. On the other hand, for larger $\lambda$, or negative $\lambda$, Banach--Colmez spaces are more exotic objects; for example, the negative Banach--Colmez space
\[
\BC(\mathcal O_{X_C}(-1)[1])\cong (\mathbb A^1_{C^\sharp})^\diamond/\underline{E}
\]
is the quotient of the affine line by the translation action of $\underline{E}\subset \mathbb A^1_{C^\sharp}$.

A key result is Proposition~\ref{prop:relativebanachcolmez}, stating in particular that projectivized Banach--Colmez spaces
\[
(\BC(\mathcal E)\setminus \{0\})/\underline{E^\times}
\]
are proper --- they are the relevant analogues of ``families of projective spaces over $S$''. In particular, their image in $S$ is a closed subset, and if the image is all of $S$, then we can find a nowhere vanishing section of $\mathcal E$ after a v-cover, as then the projectivized Banach--Colmez space is a v-cover of $S$.

\section{The Fargues--Fontaine curve}
\subsection{The curve $\mathcal Y_C$}
\label{sec:preFarguesFontaine}

Recall that for any perfect $\Fq$-algebra $R$, there is a unique $\pi$-adically complete flat $\mathcal O_E$-algebra $\tilde{R}$ such that $\tilde{R}=R/\pi$. There is a unique multiplicative lift $[\cdot]: R\to \tilde{R}$ of the identity $R\to R$, called the Teichm\"uller lift. Explicitly, one can take
\[
\tilde{R}=W_{\mathcal O_E}(R)=W(R)\hat{\otimes}_{W(\mathbb F_q)} \mathcal O_E
\]
in terms of the ramified Witt vectors; here the completion is the $\pi$-adic completion. (In the case $E=\mathbb F_q((\pi))$ is of equal characteristic, this becomes simply $R[[\pi]]$.)

The construction of the Fargues--Fontaine curve is based on this construction on the level of perfectoid spaces $S$ over $\Fq$. Its construction is done in three steps. First, one constructs a curve $\mathcal Y_S$, an adic space over $\mathcal O_E$, which carries a Frobenius action $\phi$. Passing to the locus $Y_S=\mathcal Y_S\setminus \{\pi=0\}$, i.e.~the base change to $E$, the action of $\phi$ is free and totally discontinuous, so that one can pass to the quotient $X_S=Y_S/\phi^{\mathbb Z}$, which will be the Fargues--Fontaine curve.

We start by constructing $\mathcal Y_S$ in the affinoid case. More precisely, if $S=\Spa(R,R^+)$ is an affinoid perfectoid space over $\Fq$, and $\varpi\in R^+$ is a pseudouniformizer (i.e.~a topologically nilpotent unit of $R$), we let
\[
\mathcal Y_S = \Spa W_{\mathcal O_E}(R^+)\setminus V([\varpi]).
\]
Here $W_{\mathcal O_E}(R^+)$ has the $(\pi,[\varpi])$-adic topology. These objects do not depend on the choice of $\varpi$, as for any choice of $\varpi,\varpi'\in R$, one has $\varpi|\varpi^{\prime n}$, $\varpi'|\varpi^n$ for some $n>0$. The $q$-th power Frobenius of $R^+$ induces an automorphism $\phi$ of $\mathcal Y_S$. To construct the Fargues--Fontaine curve, we will eventually remove $V(\pi)$ from $\mathcal Y_S$ and quotient by $\phi$, but for now we recall some properties of $\mathcal Y_S$.

\begin{proposition}\label{prop:YSadicspace} The above defines an analytic adic space $\mathcal Y_S$ over $\mathcal O_E$. Letting $E_\infty$ be the completion of $E(\pi^{1/p^\infty})$, the base change
\[
\mathcal Y_S\times_{\Spa \mathcal O_E} \Spa \mathcal O_{E_\infty}
\]
is a perfectoid space, with tilt given by
\[
S\times_{\Fq} \Spa \Fq\powerseries{t^{1/p^\infty}} = \mathbb D_{S,\perf},
\]
a perfectoid open unit disc over $S$.
\end{proposition}

\begin{proof} One can cover $\mathcal Y_S$ by the subsets $\mathcal Y_{S,[0,n]}:=\{|\pi|^n\leq |[\varpi]|\neq 0\}\subset \mathcal Y_S$, which are rational subsets of $\Spa W_{\mathcal O_E}(R^+)$, where $n>0$ is some integer that we assume to be a power of $p$ for simplicity. Then
\[
\mathcal Y_{S,[0,n]} = \Spa(B_{S,[0,n]},B_{S,[0,n]}^+)
\]
where
\[
B_{S,[0,n]} = W_{\mathcal O_E}(R^+)\langle \frac{\pi^n}{[\varpi]}\rangle [\tfrac 1{[\varpi]}]
\]
and $B_{S,[0,n]}^+\subset B_{S,[0,n]}$ is the integral closure of $W_{\mathcal O_E}(R^+)\langle \frac{\pi^n}{[\varpi]}\rangle$. To see that $\mathcal Y_S$ is an adic space (i.e.~the structure presheaf is a sheaf) and $\mathcal Y_S\times_{\Spa \mathcal O_E} \Spa \mathcal O_{E_\infty}$ is perfectoid, it is enough to prove that $B_{S,[0,n]}\hat{\otimes}_{\mathcal O_E} \mathcal O_{E_\infty}$ is a perfectoid Tate algebra. Indeed, the algebra $B_{S,[0,n]}$ splits off $B_{S,[0,n]}\hat{\otimes}_{\mathcal O_E} \mathcal O_{E_\infty}$ as a direct factor as topological $B_{S,[0,n]}$-module, and hence the sheaf property for perfectoid spaces gives the result for $\mathcal Y_{S,[0,n]}$ and thus all of $\mathcal Y_S$ (cf.~the sousperfectoid property of \cite{HansenKedlaya}, \cite[Section 6.3]{Berkeley}). Using the Frobenius automorphism of $(R,R^+)$, one can in fact assume that $n=1$.

Let us abbreviate
\[
A=B_{S,[0,1]}\hat{\otimes}_{\mathcal O_E} \mathcal O_{E_\infty}
\]
and $A^+\subset A$ the integral closure of $B_{S,[0,1]}^+\hat{\otimes}_{\mathcal O_E} \mathcal O_{E_\infty}$. In particular
\[
A_0^+=(W_{\mathcal O_E}(R^+)\hat{\otimes}_{\mathcal O_E} \mathcal O_{E_\infty})[(\tfrac{\pi}{[\varpi]})^{1/p^\infty}]^\wedge_{[\varpi]}\subset A^+,
\]
and $A=A_0^+[\tfrac 1{[\varpi]}]$. Note that
\[
A_0^+/[\varpi]=(R^+/\varpi\otimes_{\Fq} \mathcal O_{E_\infty}/\pi)[t_1^{1/p^\infty}]/(\pi^{1/p^m}-[\varpi]^{1/p^m}t_1^{1/p^m})\cong R^+/\varpi[t_1^{1/p^\infty}].
\]
This implies already that $A_0^+$ is integral perfectoid by \cite[Lemma 3.10 (ii)]{BMS1}, and thus necessarily (cf.~\cite[Lemma 3.21]{BMS1}) $A_0^+\to A^+$ is an almost isomorphism and $A_0^+[\tfrac 1{[\varpi]}]\cong A$ is perfectoid. Moreover, one can see that the tilt of $A$ is given by $R\langle t_1^{1/p^\infty}\rangle$, where $t_1^\sharp = \tfrac \pi{[\varpi]}$, which corresponds to the subset
\[
\{|t|\leq |\varpi|\neq 0\}\subset S\times_{\Fq} \Spa \Fq\powerseries{t^{1/p^\infty}} = \mathbb D_{S,\perf}.\qedhere
\]
\end{proof}

\begin{proposition}\label{prop:diamondequation} For any perfectoid space $T$ over $\Fq$, giving an untilt $T^\sharp$ of $T$ together with a map $T^\sharp\to \mathcal Y_S$ of analytic adic spaces is equivalent to giving an untilt $T^\sharp$ together with a map $T^\sharp\to \Spa \mathcal O_E$, and a map $T\to S$. In other words, there is a natural isomorphism
\[
\mathcal Y_S^\diamond\cong \Spd \mathcal O_E\times S.
\]
\end{proposition}

\begin{proof} Changing notation, we need to see that for any perfectoid space $T$ over $\mathcal O_E$, giving a map $T\to \mathcal Y_S$ is equivalent to giving a map $T^\flat\to S$. Without loss of generality, assume that $T=\Spa(A,A^+)$ is affinoid. Giving a map $T\to \mathcal Y_S$ is equivalent to giving a map $W_{\mathcal O_E}(R^+)\to A^+$ such that the image of $[\varpi]$ in $A$ is invertible. By the universal property of $W_{\mathcal O_E}(R^+)$ in case $R^+$ is perfect, this is equivalent to giving a map $R^+\to (A^+)^\flat$ such that the image of $\varpi$ in $A^\flat$ is invertible. But this is precisely a map $T^\flat=\Spa(A^\flat,A^{\flat +})\to S=\Spa(R,R^+)$.
\end{proof}

In particular, there is a natural map
\[
|\mathcal Y_S|\cong |\mathcal Y_S^\diamond|\cong |\Spd \mathcal O_E\times S|\to |S|.
\]

The following proposition ensures that we may glue $\mathcal Y_S$ for general $S$, i.e.~for any perfectoid space $S$ there is an analytic adic space $\mathcal Y_S$ equipped with an isomorphism
\[
\mathcal Y_S^\diamond\cong \Spd \mathcal O_E\times S
\]
(and in particular a map $|\mathcal Y_S|\to |S|$) such that for $U=\Spa(R,R^+)\subset S$ an affinoid subset, the corresponding pullback of $\mathcal Y_S$ is given by $\mathcal Y_U$.

\begin{proposition}\label{prop:basechangecurve} If $S'\subset S$ is an affinoid subset, then $\mathcal Y_{S'}\to \mathcal Y_S$ is an open immersion, with
\[\xymatrix{
|\mathcal Y_{S'}|\ar[r]\ar[d] & |\mathcal Y_S|\ar[d] \\
|S'|\ar[r] & |S|
}\]
cartesian.
\end{proposition}

\begin{proof} Let $Z\subset \mathcal Y_S$ be the open subset corresponding to $|\mathcal Y_S|\times_{|S|} |S'|\subset |\mathcal Y_S|$. Then by functoriality of the constructions, we get a natural map of adic spaces $\mathcal Y_{S'}\to Z$. To see that it is an isomorphism, we can check after base change to $\mathcal O_{E_\infty}$ (as the maps on structure sheaves are naturally split injective). The base change of $\mathcal Y_{S'}$ and $Z$ become perfectoid, and hence it suffices to see that one gets an isomorphism after passing to diamonds, where it follows from Proposition~\ref{prop:diamondequation}.
\end{proof}

Next, we recall the ``sections of $\mathcal Y_S\to S$''.

\begin{proposition}[{\cite[Proposition 11.3.1]{Berkeley}}]\label{prop:degreeonecartier} Let $S$ be a perfectoid space over $\Fq$. The following objects are in natural bijection.
\begin{altenumerate}
\item[{\rm (i)}] Sections of $\mathcal Y_S^\diamond\to S$;
\item[{\rm (ii)}] Morphisms $S\to \Spd \mathcal O_E$;
\item[{\rm (iii)}] Untilts $S^\sharp$ over $\mathcal O_E$ of $S$.
\end{altenumerate}

Moreover, given an untilt $S^\sharp$ over $\mathcal O_E$ of $S$, there is a natural closed immersion of adic spaces
\[
S^\sharp\hookrightarrow \mathcal Y_S
\]
that presents $S^\sharp$ as a closed Cartier divisor in $\mathcal Y_S$.
\end{proposition}

\begin{proof} The equivalence of (i), (ii) and (iii) is a direct consequence of Proposition~\ref{prop:diamondequation}. Thus, let $S^\sharp$ be an untilt of $S$ over $\mathcal O_E$. We may work locally, so assume $S=\Spa(R,R^+)$ is affinoid. Then $S^\sharp=\Spa(R^\sharp,R^{\sharp +})$ is affinoid perfectoid as well, and
\[
R^{\sharp +}=W_{\mathcal O_E}(R^+)/\xi
\]
for some nonzerodivisor $\xi\in W_{\mathcal O_E}(R^+)$ that can be chosen to be of the form $\pi-a[\varpi]$ for some $a\in W_{\mathcal O_E}(R^+)$ and suitable topologically nilpotent $\varpi\in R$ (choose $\varpi \in R^+$ a pseudouniformizing element such that $\varpi^\sharp |\pi$, and write $\pi = \varpi^\sharp \theta(a)$ for some $a$). To see that $S^\sharp$ defines a closed Cartier divisor in $\mathcal Y_S$, that is to say the sequence 
\[
0 \to \mathcal O_{\mathcal Y_S}\xrightarrow{\xi} \mathcal O_{\mathcal Y_S} \to i_\ast \mathcal O_{S^\sharp}\to 0
\]
is exact with $i: S^\sharp \hookrightarrow \mathcal Y_S$, we need to see that for any open affinoid $U=\Spa(A,A^+)\subset \mathcal Y_S$ with affinoid perfectoid pullback $V=\Spa(B,B^+)\subset S^\sharp$, the sequence
\[
0\to A\xrightarrow{\xi} A\to B\to 0
\]
is exact. To see this, we are free to localize near $S^\sharp=V(\xi)\subset \mathcal Y_S$. In particular, replacing $S$ by $V^\flat$, we can assume that $V=S^\sharp$. In that case, any neighborhood of $S^\sharp=V(\xi)$ in $\mathcal Y_S$ contains $\{|\xi|\leq |[\varpi]|^n\}$ for some $n>0$, so we can assume that $U$ is of this form.

Endow $A$ with the spectral norm, where we normalize the norm on each completed residue field of $\mathcal Y_S$ by $|[\varpi]|=\tfrac 1q$. We claim that with this choice of norm, one has
\[
|\xi a|\geq q^{-n} |a|
\]
for all $a\in A$. In particular, this implies that $\xi: A\to A$ is injective, and has closed image (as the preimage of any Cauchy sequence in the image is a Cauchy image). On the other hand, $R^\sharp$ is the separated completion of $A/\xi$, so $B=A/\xi$.

To verify the claimed inequality, it is enough to see that the norm of $|a|$ is equal to the supremum norm over $\{|\xi|=|[\varpi]|^n\}$. In fact, it is enough to consider the points in the Shilov boundary, i.e.~those points $\Spa(C,\mathcal O_C)\to U$ that admit a specialization $\Spa(C,C^+)\to \mathcal Y_S$ whose image is not contained in $U$; any such is necessarily contained in $\{|\xi|=|[\varpi]|^n\}$. This will in fact hold for all functions on $U\times_{\Spa \mathcal O_E} \Spa \mathcal O_{E_\infty}$, for which the claim reduces to the tilt, which is an affinoid subset of $\mathbb D_{S,\perf}$.
 By approximation, it then reduces to the case of affinoid subsets of $\mathbb D_S$, where it is well-known that the maximum is taken on the Shilov boundary. (Note that this question immediately reduces to the case that $S$ is a geometric point.)
\end{proof}

\begin{remark}
The preceding Cartier divisor satisfies the stronger property of being a ``relative Cartier divisor'' in the sense that for all $s\in S$ its pullback to $\mathcal Y_{\Spa(K(s),K(s)^+)}$ is a Cartier divisor.
\end{remark}

Now let us analyze the case $S=\Spa C$ for some complete algebraically closed nonarchimedean field over $\Fq$.

\begin{example} Assume that $E=\mathbb F_q\laurentseries{t}$ is of equal characteristic. Then $\mathcal Y_C = \mathbb D_C$ is an open unit disc over $C$, with coordinate $t$. In particular, inside $|\mathcal Y_C|$, we have the subset of classical points $|\mathcal Y_C|^{\mathrm{cl}}\subset |\mathcal Y_C|$, which can be identified as
\[
|\mathcal Y_C|^{\mathrm{cl}}=\{x\in C\mid |x|<1\}.
\]
Note that these classical points are in bijection with maps $\mathcal O_E\to C$ (over $\mathbb F_q$), i.e.~with ``untilts of $C$ over $\mathcal O_E$''.
\end{example}

With suitable modifications, the same picture exists also when $E$ is of mixed characteristic.

\begin{defprop}\label{prop:classicalpointsY} Any untilt $C^\sharp$ of $C$ over $\mathcal O_E$ defines a closed Cartier divisor $\Spa C^\sharp\hookrightarrow \mathcal Y_S$, and in particular a closed point of $|\mathcal Y_C|$. This induces an injection from the set of such untilts to $|\mathcal Y_C|$.

The set of classical points $|\mathcal Y_C|^{\mathrm{cl}}\subset |\mathcal Y_C|$ is defined to be the set of such points.
\end{defprop}

\begin{proof} We have seen that any untilt $C^\sharp$ defines such a map $\Spa C^\sharp\hookrightarrow \mathcal Y_S$. As it is a closed Cartier divisor, the corresponding point is closed in $|\mathcal Y_C|$. One can recover $C^\sharp$ as the completed residue field at the point, together with the map $W_{\mathcal O_E}(\mathcal O_C)\to \mathcal O_{C^\sharp}$, which induces the isomorphism $\mathcal O_C\cong \mathcal O_{C^\sharp}^\flat$ and thus $C\cong (C^\sharp)^\flat$, giving the untilt structure on $C^\sharp$; this shows that the map is injective.
\end{proof}

Recall that $\mathcal Y_C$ is preperfectoid. In fact, if one picks a uniformizer $\pi\in E$ and lets $E_\infty$ be the completion of $E(\pi^{1/p^\infty})$, then $\mathcal Y_C\times_{\mathcal O_E} \mathcal O_{E_\infty}$ is perfectoid, and its tilt is given by
\[
\Spa C\times \Spa \mathcal O_{E_\infty}^\flat\cong \Spa C\times \Spa \mathbb F_q\powerseries{t^{1/p^\infty}}.
\]
Thus, we get a map
\[
|\mathbb D_C|=|\Spa C\times \Spa \mathbb F_q\powerseries{t}|\cong |\Spa C\times \Spa \mathbb F_q\powerseries{t^{1/p^\infty}}|\cong |\mathcal Y_C\times_{\mathcal O_E}\mathcal O_{E_\infty}|\to |\mathcal Y_C|.
\]

\begin{proposition}\label{prop:classicalpointstilting} Under this map, the classical points $|\mathbb D_C|^{\mathrm{cl}}=\{x\in C\mid |x|<1\}\subset |\mathbb D_C|$ are exactly the preimage of the classical points $|\mathcal Y_C|^{\mathrm{cl}}\subset |\mathcal Y_C|$.
\end{proposition}

Unraveling the definitions, one sees that the map
\[
\{x\in C\mid |x|<1\}=|\mathbb D_C|^{\mathrm{cl}}\to |\mathcal Y_C|^{\mathrm{cl}}
\]
sends any $x\in C$ with $|x|<1$ to the closed point defined by the ideal $(\pi-[x])$. In particular, the proposition shows that any classical point of $\mathcal Y_C$ can be written in this form.

\begin{proof} This is clear as classical points are defined in terms of maps of diamonds, which are compatible with this tilting construction on topological spaces.
\end{proof}

The formation of classical points is also compatible with changing $C$ in the following sense.

\begin{proposition}\label{prop:classicalpointschangeoffield} Let $C'|C$ be an extension of complete algebraically closed nonarchimedean fields over $\mathbb F_q$, inducing the map $\mathcal Y_{C'}\to \mathcal Y_C$. A point $x\in |\mathcal Y_C|$ is classical if and only if its preimage in $|\mathcal Y_{C'}|$ is a classical point. Moreover, if $x\in |\mathcal Y_C|$ is a rank-$1$-point that is not classical, then there is some $C'|C$ such that the preimage of $x$ contains a nonempty open subset of $|\mathcal Y_{C'}|$.
\end{proposition}

In other words, one can recognize classical points as those points that actually stay points after any base change; all other rank $1$ points actually contain whole open subsets after some base change.

\begin{proof} It is clear that if $x$ is classical, then its preimage is a classical point. Conversely, if $x\in |\mathcal{Y}_C|$ is a rank $1$ point, and $S=\Spd K(x)$, the point $x$ is given by a morphism $S\to \Spa C\times \Spd \O_E$. If the preimage of $x$ is a classical point,  the induced morphism $S\to \Spa C$ becomes an isomorphism after pullback via $\Spa C'\to \Spa C$. Since $S$ is a v-sheaf (\cite[Proposition 11.9]{ECoD}) and $\Spa C'\to \Spa C$ a v-cover, the morphism $S\to \Spa C$ is an isomorphism, and thus $x$ is a classical point.

Now assume that $x$ is nonclassical rank-$1$-point; we want to find $C'|C$ such that the preimage of $x$ contains an open subset of $|\mathcal Y_{C'}|$. By Proposition~\ref{prop:classicalpointstilting}, it is enough to prove the similar result for $\mathbb D_C$, using that $|\mathbb D_{C'}|\to |\mathcal Y_{C'}|$ is open.\footnote{Any quasicompact open of $|\mathbb D_{C'}|$ is the base change of a quasicompact open of $|\mathcal Y_{C'}\times_{\Spa \mathcal O_E} \Spa \mathcal O_{E'}|$ for a finite extension $E'|E$. Passing to the Galois hull of $E'$ and taking the orbit of the open subset under the Galois group, the openness of the image follows from the map being a quotient map, as is any surjective quasicompact map of analytic adic spaces.} Thus, assume $x\in |\mathbb D_C|$ is a non-classical point. Let $C'$ be a completed algebraic closure of the corresponding residue field. Then the preimage of $x$ in $|\mathbb D_{C'}|$ has a tautological section $\tilde{x}\in \mathbb D_{C'}(C')$ which is a classical point, and the preimage of $x$ contains a small disc $\mathbb B(\tilde{x},r)\subset \mathbb D_{C'}$ for some $r>0$. Indeed, this follows from the description of the rank $1$ points of $\mathbb D_C$ as being either the Gauss norm for some disc $\mathbb B(x,r_0)\subset \mathbb D_C$ of radius $r_0>0$, or the infimum of such over a decreasing sequence of balls (but with radii not converging to zero). See Lemma \ref{lem:image reciproque point non classique ouvert}.
\end{proof}

\begin{lemma}\label{lem:image reciproque point non classique ouvert}
Let $x\in \mathbb{D}_C(C)$, $\rho \in (0,1]$, and $x_\rho\in |\mathbb{D}_C|$ be the Gauss norm with radius $\rho$ centered at $x$. The preimage of $x_\rho$ in $|\mathbb{D}_{C(x_\rho)}|$ contains the open disk with radius $\rho$ centered at $x_\rho\in \mathbb D_{C(x_\rho)}(C(x_\rho))$.
\end{lemma}

\begin{proof}
We can suppose $x=0$. The point $x_\rho$ is given by the morphism $C\langle T\rangle \to C(x_\rho)$ that sends $T$ to $t$. Let $y\in |\mathbb{D}_{C(x_\rho)}|$. This corresponds to a morphism 
$C(x_\rho)\langle T\rangle \to C(x_\rho)(y)$.
 Let us note $u\in C(x_\rho)(y)$ the image of $T$ via the preceding map. Suppose $y$ lies in the open disk with radius $\rho$ centered at $x_\rho$. This means $|u-t|<\rho = |t|$. 
Let us remark that this implies that for any $n\geq 1$,
\[
|u^n-t^n|=|u-t||u^{n-1}+\ldots+t^{n-1}|\leq |u-t|\rho^{n-1}<\rho^n.
\]
For $f=\sum_{n\geq 1} a_n T^n\in  C\langle T\rangle$, one then has
\[
|\sum_{n\geq 1} a_n (u^n-t^n)|< \sup_{n\geq 1} |a_n| \rho^n = |\sum_{n\geq 1} a_n t^n|.
\]
We deduce that
\[
|\sum_{n\geq 1} a_n u^n|= |\sum_{n\geq 1} a_n t^n|=|f(x_\rho)|.\qedhere
\]
\end{proof}

There is in fact another characterization of the classical points in terms of maximal ideals.

\begin{proposition}\label{prop:classicalpointsmaximalspectrum} Let $U=\Spa(B,B^+)\subset \mathcal Y_C$ be an affinoid subset. Then for any maximal ideal $\mathfrak m\subset B$, the quotient $B/\mathfrak m$ is a nonarchimedean field, inducing an injection $\mathrm{Spm}(B)\hookrightarrow |U|$. This gives a bijection between $\mathrm{Spm}(B)$ and $|U|^{\mathrm{cl}}:=|U|\cap |\mathcal Y_C|^{\mathrm{cl}}\subset |\mathcal Y_C|$.
\end{proposition}

\begin{proof} First, if $x\in |U|^{\mathrm{cl}}$, then it corresponds to a closed Cartier divisor $\Spa C^\sharp\hookrightarrow U\subset \mathcal Y_C$, and thus defines a maximal ideal of $B$, yielding an injection $|U|^{\mathrm{cl}}\hookrightarrow \mathrm{Spm}(B)$. We need to see that this is a bijection.

Note that using the tilting map $|\mathbb D_C|\to |\mathcal Y_C|$, one sees that the preimage of $U$ in $|\mathbb D_C|$ has only finitely many connected components  (any quasicompact open subset of $|\mathbb D_C|$ has finitely many connected components); we can thus assume that $U$ is connected. In that case, we claim that any nonzero element $f\in B$ vanishes only at classical points of $|U|$. By Proposition~\ref{prop:classicalpointschangeoffield}, it suffices to see that for any nonempty open subset $U'\subset U$, the map $\mathcal O(U)\to \mathcal O(U')$ is injective. In fact, if $V(f)$ contains a nonclassical point, it also contains a nonclassical rank $1$ point as $V(f)$ is generalizing, then after base changing to some $C'|C$, $V(f)$ contains an open subset $U'$, and this is impossible if $\O (U)\hookrightarrow \O(U')$. For this it suffices to prove that  $\O (V) \hookrightarrow \O (V')$ where $V$ is a connected component of $U\hat{\otimes}_{\O_E} \O_{E_\infty}$, and $V'$ the intersection of $U'\hat{\otimes}_{\O_E}\O_{E_\infty}$ with this connected component. Now for any $g\in \O (V) \setminus \{0\}$, $V(g)\neq V$, as perfectoid spaces are uniform (and hence vanishing at all points implies vanishing).  We thus have to prove that for any Zariski closed subset $Z\subsetneq V$, $V'\not\subset Z$.

By Proposition~\ref{prop:zariskiclosedstronglyzariskiclosed}, it suffices to prove the similar property for open subsets $V'\subset V\subset \mathbb D_{C,\perf}$, with $V$ connected. But then $V'=W'_\perf$ and $V=W_\perf$ for $W'\subset W\subset \mathbb D_C$, and $\mathcal O(V)\to \mathcal O(V')$ is topologically free (with basis $t^i$, $i\in [0,1)\cap \mathbb Z[p^{-1}]$) over the corresponding map $\mathcal O(W)\to \mathcal O(W')$ of classical Tate algebras over $C$, for which injectivity is classical.
\end{proof}

The previous proposition implies that, once $U$ is connected, the rings $B$ are principal ideal domains (cf.~\cite{KedlayaNoetherianProperties}).

\begin{corollary}\label{cor:principalidealdomains} Let $U=\Spa(B,B^+)\subset \mathcal Y_C$ be an affinoid subset. Then $U$ has finitely many connected components. Assuming that $U$ is connected, the ring $B$ is a principal ideal domain.
\end{corollary}

\begin{proof} We have already seen in the preceding proof that $U$ has finitely many connected components. Passing to one component, we can assume that $U$ is connected. Each maximal ideal of $B$ is principal, as it comes from a closed Cartier divisor on $U$. Now take any nonzero $f\in B$. We have seen (in the preceding proof) that the vanishing locus of $f$ is contained in $|U|^{\mathrm{cl}}$, and it is also closed in $|U|$. It is thus a spectral space with no nontrivial specializations, and therefore a profinite set. We claim that it is in fact discrete. For this, let $x\in V(f)$ be any point. We get a generator $\xi_x\in B$ for the corresponding maximal ideal. We claim that there is some $n\geq 1$ such that $f=\xi_x^n g$ where $g$ does not vanish at $x$. Assume otherwise. Note that the spectral norm on $U$ is given by the supremum over finitely many points, the Shilov boundary of $U$ (cf.~proof of Proposition~\ref{prop:degreeonecartier}). We may normalize $\xi_x$ so that its norm at all of these finitely many points is $\geq 1$. Then for any $n$, if $f=\xi_x^n g_n$, one has $||g_n||\leq ||f||$. But inside the open neighborhood $U_x=\{|\xi_x|\leq |[\varpi]|\}$ of $x$, this implies that $||f||_{U_x}\leq |[\varpi]|^n ||f||$ for all $n$, and thus $||f||_{U_x}=0$ as $n\to \infty$. Thus, $f$ vanishes on all of $U_x$, which is a contradiction.

By the above, we can write $f=\xi_x^n g$ where $g$ does not vanish at $x$. But then $g$ does not vanish in a neighborhood of $x$, and therefore $x\in V(f)$ is an isolated point, and hence $V(f)$ is profinite and discrete, and thus finite. Enumerating these points $x_1,\ldots,x_m$, we can thus write $f=\xi_{x_1}^{n_1}\cdots \xi_{x_m}^{n_m} g$ where $g$ does not vanish at $x_1,\ldots,x_m$, and thus vanishes nowhere, and hence is a unit. This finishes the proof.
\end{proof}

\begin{remark}
The main new ingredient compared to \cite{KedlayaNoetherianProperties} or \cite[Theorem 2.5.1]{FarguesFontaine} that allows us to shorten the proof is Proposition~\ref{prop:zariskiclosedstronglyzariskiclosed}, i.e.~the use of the fact (proved in \cite{BhattScholzePrism}) that ``Zariski closed implies strongly Zariski closed'' in the terminology of \cite[Section II.2]{ScholzeTorsion}. 
\end{remark}

Later (cf.~Proposition~\ref{prop:drinfeldlemmalocallyconstant}), we will also need the following lemma about non-classical points of $Y_C=\mathcal Y_C\times_{\Spa \mathcal O_E} \Spa E$.

\begin{lemma}\label{lem:genericpointcurve} There is a point $x\in |Y_C|$, with completed residue field $K(x)$, such that the induced map $\Gal(\overline{K(x)}|K(x))\to I_E$ is surjective, where $I_E$ is the inertia subgroup of the absolute Galois group of $E$. 
\end{lemma}

Note that a priori we have a map $\Gal(\overline{K(x)}|K(x))\to \mathrm{Gal}(\overline{E}|E)$, but it is clear that its image is contained in $I_E$, as $K(x)$ contains $\breve E$.

\begin{proof} In fact, we can be explicit: Looking at the surjection
\[
|\mathbb D^\ast_C|\to |Y_C|
\]
from the tilting construction, the image of any Gau\ss point (corresponding to a disc of radius $r$, $0<r<1$, around the origin) will have the desired property. This follows from the observation that this locus of Gau\ss\ points lifts uniquely to $|Y_C\times_{\Spa \breve E}\Spa E'|$ for any finite extension $E'|\breve E$. In fact, this cover admits a similar surjection from a punctured open unit disc over $C$, and there is again one Gau\ss\ point for each radius (i.e.~the set of Gau\ss\ points maps isomorphically to $(0,\infty)$ via $\mathrm{rad}: |Y_C|\to (0,\infty)$).
\end{proof}

\subsection{The Fargues--Fontaine curve}
\label{subsec:FarguesFontaine}

Now we can define the Fargues--Fontaine curve.

\begin{definition}\label{def:farguesfontaine} For any perfectoid space $S$ over $\Fq$, the relative Fargues--Fontaine curve is
\[
X_S = Y_S/\phi^{\mathbb Z}
\]
where
\[
Y_S=\mathcal Y_S\times_{\Spa \mathcal O_E}\Spa E = \mathcal Y_S\setminus V(\pi),
\]
which for affinoid $S=\Spa(R,R^+)$ with pseudouniformizer $\varpi$ is given by
\[
Y_S = \Spa W_{\mathcal O_E}(R^+)\setminus V(\pi[\varpi]).
\]
\end{definition}

To see that this is well-formed, we note the following proposition, cf.~\cite[Lecture 12]{Berkeley}.

\begin{proposition}\label{prop:actiontotallydiscontinuous} The action of $\phi$ on $Y_S$ is free and totally discontinuous. In fact, if $S=\Spa(R,R^+)$ is affinoid and $\varpi\in R$ is a pseudouniformizer, one can define a map
\[
\mathrm{rad}: |Y_S|\longrightarrow (0,\infty)
\]
taking any point $x\in Y_S$ with rank-$1$-generalization $\tilde{x}$ to $\log |[\varpi](\tilde{x})|/\log |\pi(\tilde{x})|$. This factorizes through the Berkovich space quotient of $|Y_S|$ and satisfies $\mathrm{rad}\circ \phi = q\cdot \mathrm{rad}$.

For any interval $I=[a,b]\subset (0,\infty)$ with rational ends (possibly with $a=b$), there is the open subset
\[
Y_{S,I}=\{|\pi|^b\leq |[\varpi]|\leq |\pi|^a\}\subset \mathrm{rad}^{-1}(I)\subset Y_S
\]
which is in fact a rational open subset of $\Spa W_{\mathcal O_E}(R^+)$ and thus affinoid,
\[
Y_{S,I} = \Spa(B_{S,I},B_{S,I}^+),
\]
and one can form $X_S$ as the quotient of $Y_{S,[1,q]}$ via the identification $\phi: Y_{S,[1,1]}\cong Y_{S,[q,q]}$. In particular, $X_S$ is qcqs in case $S$ is affinoid.
\end{proposition}

\begin{proof} This follows directly from the definitions.
\end{proof}

In terms of the preceding radius function, the end $0$ corresponds to the boundary divisor $(\pi)$, and $\infty$ to the boundary divisor $([\varpi])$.

For each $s\in S$ corresponding to a map $\Spa(K(s),K(s)^+)\to S$, functoriality defines a morphism $X_{K(s),K(s)^+}\rightarrow X_S$. We way think of $X_S$ as the collection of curves $(X_{K(s),K(s)^+})_{s\in S}$, the one defined and studied in \cite{FarguesFontaine}, merged in a ``family of curves''. Although $X_S$ does not sit over $S$, 
the absolute Frobenius $\phi\times \phi$ of $S\times \Spd (E)$ acts trivially on the topological space and one has 
$$
|X_S|\cong |X_S^\diamond|\cong |S\times \Spd(E) / \phi^\Z\times \mathrm{id}|\cong | S\times \Spd (E) / \mathrm{id} \times \phi^\Z |\longrightarrow |S|.
$$
Thus the topological space $|X_S|$ sits over $|S|$, and for all $S$ the map $|X_S|\to |S|$ is qcqs. Here, we used the following identification of the diamond.

\begin{proposition}\label{prop:diamondequation2} There is a natural isomorphism
\[
Y_S^\diamond\cong S\times \Spd(E),
\]
descending to an isomorphism
\[
X_S^\diamond\cong (S\times \Spd(E)) / \phi^\Z\times \mathrm{id}.
\]
\end{proposition}

\begin{proof} This is immediate from Proposition~\ref{prop:diamondequation}.
\end{proof}

Moreover, we have the following version of Proposition~\ref{prop:degreeonecartier}.

\begin{proposition}\label{prop:degreeonecartier1} The following objects are naturally in bijection.
\begin{altenumerate}
\item[{\rm (i)}] Sections of $Y_S^\diamond\to S$;
\item[{\rm (ii)}] Maps $S\to \Spd(E)$;
\item[{\rm (iii)}] Untilts $S^\sharp$ over $E$ of $S$.
\end{altenumerate}

Given such a datum, in particular an untilt $S^\sharp$ over $E$ of $S$, there is a natural closed immersion $S^\sharp\hookrightarrow Y_S$ presenting $S^\sharp$ as a closed Cartier divisor in $Y_S$. The composite map $S^\sharp\to Y_S\to X_S$ is still a closed Cartier divisor, and depends only on the composite $S\to \Spd(E)\to \Spd(E)/\phi^\Z$. In this way, any map $S\to \Spd(E)/\phi^\Z$ defines a closed Cartier divisor $D\subset X_S$; this gives an injection of $\Spd(E)/\phi^\Z$ into the space of closed Cartier divisors on $X_S$.
\end{proposition}

\begin{proof} This is immediate from Proposition~\ref{prop:degreeonecartier}.
\end{proof}

\begin{definition}\label{def:degree1cartier} A closed Cartier divisor of degree $1$ on $X_S$ is a closed Cartier divisor $D\subset X_S$ that arises from a map $S\to \Spd(E)/\phi^\Z$. Equivalently, it arises locally on $S$ from an untilt $S^\sharp$ over $E$ of $S$.
\end{definition}

The quotient $\Spd(E)/\phi^\Z$ that occurs here is the quotient in the category of v-sheaves; but we note that it agrees with the quotient computed in the category of sheaves on $\mathrm{Perf}_{\mathbb F_q}$ for the topology of open covers. In particular, ``locally on $S$'' in the preceding definition can be taken to mean v-locally, or on open subsets of $|S|$.

In particular, we see that the moduli space $\mathrm{Div}^1$ of degree $1$ closed Cartier divisors is given by
\[
\mathrm{Div}^1 = \Spd (E) / \phi^\Z.
\]
Note that something strange is happening in the formalism here: Usually the curve itself would represent the moduli space of degree $1$ Cartier divisors!

\begin{remark}
In \cite[D\'efinition 2.6]{FarguesClassFieldTheory} Fargues gives a definition of a Cartier divisor of degree $1$ on $X_S$ equivalent to the preceding one, similar to the definition of a relative Cartier divisor in classical algebraic geometry.
\end{remark}

In the next proposition and elsewhere, we write $\ast$ for the v-sheaf taking any $S\in \mathrm{Perf}_{\mathbb F_q}$ to a point $\ast$; one could also write $\ast=\mathrm{Spd}(\mathbb F_q)$.

\begin{proposition}\label{prop:div1propersmooth} The map $\mathrm{Div}^1\to \ast$ is proper, representable in spatial diamonds, and cohomologically smooth.
\end{proposition}

\begin{proof} First, $\Spd(E)\to \ast$ is representable in locally spatial diamonds and cohomologically smooth by \cite[Proposition 24.5]{ECoD} (for $E=\mathbb Q_p$, which formally implies the case of $E$ finite over $\mathbb Q_p$, and the equal characteristic case is handled in the proof). As $|\Spd (E)\times S|\cong |Y_S|\to |S|$, we see that $\phi^\Z$ acts totally discontinuously with quotient $|\Spd(E)/\phi^\Z\times S|\cong |X_S|\to |S|$ being qcqs in case $|S|$ is qcqs; thus, $\Spd(E)/\phi^\Z\to \ast$ is representable in spatial diamonds, in particular qcqs. Then being proper follows from the valuative criterion \cite[Proposition 18.3]{ECoD}.
\end{proof}

In particular, the map
$$
|X_S|=|\mathrm{Div}^1\times S|\longrightarrow |S|
$$
is open and closed. {\it We can thus picture $X_S$ as being ``a proper and smooth family over $S$''.}

Further motivation for Definition~\ref{def:degree1cartier} is given by the following.

\begin{defprop}\label{def:classicalpointsFF} The classical points of $X_C$ are $|X_C|^{\mathrm{cl}} := |Y_C|^{\mathrm{cl}}/\phi^\Z\subset |X_C|=|Y_C|/\phi^\Z$. They are in bijection with $(\Spd(E)/\phi^\Z)(C)=\mathrm{Div}^1(C)$, i.e.~are given untilts of $C$ over $E$ up to Frobenius, or by degree $1$ closed Cartier divisors on $X_C$. For any affinoid open subset $U=\Spa(B,B^+)\subset X_C$, the maximal ideals of $B$ are in bijection with $|U|^{\mathrm{cl}}=|U|\cap |X_C|^{\mathrm{cl}}$. Any such $U$ has only finitely many connected components, and if $U$ is connected, then $B$ is a Dedekind domain.\footnote{The results on the Picard group of $X_C$ proved below actually imply that $B$ is a principal ideal domain -- the map $\mathrm{Pic}(X_C)\to \mathrm{Pic}(B)$ is surjective, the source is $\mathbb Z$ and generated by any classical point outside of $U$, so the map is zero and hence $\mathrm{Pic}(B)=0$.}
\end{defprop}

\begin{proof} This follows immediately from Proposition~\ref{prop:classicalpointsmaximalspectrum} and Corollary~\ref{cor:principalidealdomains} if $U$ lifts to $Y_C$. In general, $Y_C\to X_C$ is locally split, so the result is true locally on $U$; and then it easily follows by gluing in general.
\end{proof}

\section{Vector bundles on the Fargues--Fontaine curve}
\label{sec:vectorbundlesFarguesFontaine}

Let us recall a few basic facts about the cohomology of vector bundles. Suppose $S=\Spa(R,R^+)$ is affinoid perfectoid. Then $Y_S$ is "Stein", one has $Y_S=\bigcup_{I\subset (0,\infty)} Y_{(R,R^+),I}$ where
\begin{altenumerate}
\item as before $I$ is a compact interval with rational ends
\item $Y_{(R,R^+),I}$ is affinoid sous-perfectoid
\item for $I_1\subset I_2$, the restriction morphism $\O(Y_{(R,R^+),I_2})\to \O(Y_{(R,R^+),I_1})$ has dense image.
\end{altenumerate}
Let $\mathcal{F}$ be a vector bundle on $Y_S$. Point (2) implies that $H^i (Y_{S,I},\mathcal{F}_{|Y_{(R,R^+),I}})=0$ when $i>0$. Point (3) implies that $R^1\varprojlim_I \Gamma (Y_{(R,R^+),I},\mathcal{F})=0$ (\cite[0.13.2.4]{EGAIII1}). We thus have $H^i (Y_S,\mathcal{F})=0$ when $i>0$. 

Thus, if $\E$ is a vector bundle on $X_S$, one has 
$$
R\Gamma (X_S,\E)= \big [ H^0(Y_S,\E_{|Y_S})\xrightarrow{\ \phi-1 \ } H^0(Y_S,\E_{|Y_S}) ].
$$ 
In particular, this vanishes in degree $>1$.

Moreover, one has the following important (cohomological) descent result.

\begin{proposition}\label{prop:relativebanachcolmezvsheaf} Let $S$ be a perfectoid space over $\mathbb F_q$ and $\mathcal E$ a vector bundle on $X_S$. The functor taking any $T\in \Perf_S$ to
\[
R\Gamma(X_T,\mathcal E|_{X_T})
\]
is a v-sheaf of complexes. In fact, the functor taking any $T\in \Perf_{S}$ to $H^0(Y_T,\mathcal E|_{Y_T})$ is a v-sheaf, whose cohomology vanishes in case $T$ is affinoid.

Moreover, sending $S$ to the groupoid of vector bundles on $X_S$ defines a v-stack.
\end{proposition}

\begin{proof} By the displayed formula for $R\Gamma(X_S,\mathcal E)$ as Frobenius fixed points, it suffices to prove the result about $Y_T$. We can assume that $S=\Spa(R,R^+)$ is affinoid, pick a pseudouniformizer $\varpi\in R$, and one can further reduce to the similar claim for $Y_{T,I}$ for any compact interval $I$ with rational ends. Then $\mathcal E|_{Y_{T,I}}$ is a retract of $\mathcal O^n_{Y_{T,I}}$, so we can reduce to the structure sheaf. We need to see that for any v-cover $T=\Spa(R',R^{\prime +})$, the corresponding \v{C}ech complex
\[
0\to \mathcal O(Y_{S,I})\to \mathcal O(Y_{T,I})\to \mathcal O(Y_{T\times_S T,I})\to \ldots
\]
of $E$-Banach spaces is exact. This can be checked after taking a completed tensor product with $E_\infty = E(\pi^{1/p^\infty})^\wedge$. In that case, all algebras become perfectoid, and $Y_{T,I}\times_E E_\infty\to Y_{S,I}\times_E E_\infty$ is a v-cover of affinoid perfectoid spaces, so the result follows from \cite[Theorem 8.7, Proposition 8.8]{ECoD}.

Similarly, one proves v-descent for the groupoid of vector bundles, cf.~\cite[Lemma 17.1.8, Proposition 19.5.3]{Berkeley}.
\end{proof}

If $[\mathcal E_1\to \mathcal E_0]$ is a complex of vector bundles on $X_S$ sitting in homological degrees $[0,1]$, such that $H^0(X_T,\mathcal E_1|_{X_T})=0$ for all $T\in \Perf_{S}$, we let
\[
\BC([\mathcal E_1\to \mathcal E_0]): T\mapsto \mathbb H^0(X_T,[\mathcal E_1\to \mathcal E_0]|_{X_T})
\]
be the corresponding v-sheaf on $\Perf_{S}$. We refer to this as {\it the Banach--Colmez space associated with $[\mathcal E_1\to \mathcal E_0]$}. We will usually apply this only when either of $\mathcal E_1$ and $\mathcal E_0$ is zero.

Let us also recall the basic examples of vector bundles. Already here it is useful to fix an algebraically closed field $k|\Fq$, e.g.~$k=\overline{\mathbb F}_q$. Let $\breve E = W_{\mathcal O_E}(k)[\tfrac 1 \pi]$, the complete unramified extension of $E$ with residue field $k$, equipped with its Frobenius automorphism $\sigma$. Recall that, functorially in $S\in \Perf_k$, there is a natural exact $\otimes$-functor
\begin{align*}
\mathrm{Isoc}_k &\longrightarrow \mathrm{Bun}(X_S) \\
 (D,\phi) & \longmapsto \mathcal E(D,\phi)
\end{align*}
from the category of isocrystals (of a finite-dimensional $\breve E$-vector space $D$ equipped with a $\sigma$-linear automorphism $\phi: D\iso D$) to the category of vector bundles on $X_S$, defined via descending $D\otimes_{\breve E} \mathcal O_{Y_S}$ to $X_S$ via $\phi\otimes \phi$. We denote by $\mathcal O_{X_S}(n)$ the image of $(\breve E,\pi^{-n}\sigma)$ (note the change of sign --- the functor $\mathcal E$ reverses slopes); more generally, if $(D_\lambda,\phi_\lambda)$ is the simple isocrystal of slope $\lambda\in \mathbb Q$ in the Dieudonn\'e--Manin classification, we let $\mathcal O_{X_S}(-\lambda)=\mathcal E(D_\lambda,\phi_\lambda)$.\\

\subsection{Lubin--Tate formal groups}
\label{sec:LubinTate}

The claim of this paper is that the Fargues--Fontaine curve enables a geometrization of the local Langlands correspondence. As a warm-up, let us recall the relation between $\mathcal O_{X_S}(1)$ and local class field theory in the form of Lubin--Tate theory.

Up to isomorphism, there is a unique $1$-dimensional formal group $G$ over $\mathcal O_{\breve E}$ with action by $\mathcal O_E$, such that the two induced actions on $\mathrm{Lie}\ G$ coincide; this is ``the'' Lubin--Tate formal group $G=G_{\mathrm{LT}}$ of $E$. Fixing a uniformizer $\pi\in E$, we normalize this as follows. First, any Lubin--Tate formal group law over $\mathcal O_E$ is the unique (up to unique isomorphism) lift of a $1$-dimensional formal group over $k$ whose Lie algebra has the correct $\mathcal O_E$-action. Now, if $E$ is $p$-adic then $G_k$ is classified by Dieudonn\'e theory by a finite projective $W_{\mathcal O_E}(k)$-module $M$ equipped with a $\sigma$-linear isomorphism $F: M[\tfrac 1\pi]\cong M[\tfrac 1\pi]$ such that $M\subset F(M)\subset \tfrac 1\pi M$.\footnote{As in \cite[p. 99]{Berkeley}, we renormalize usual covariant Dieudonn\'e theory for $p$-divisible groups by dividing $F$ by $p$; and then in the case of $\pi$-divisible $\mathcal O_E$-modules as here, we base change along $W(k)\otimes_{\mathbb Z_p}\mathcal O_E\to W_{\mathcal O_E}(k)$.} Here, we take $M=W_{\mathcal O_E}(k)$ with $F=\tfrac 1\pi \sigma$. One can similarly define $G$ in equal characteristic, but actually we will explain a different way to pin down the choice just below; under our normalization, $G$ is already defined over $\mathcal O_E$.

After passing to the generic fibre, $G_E$ is isomorphic to the additive group $\mathbb G_a$, compatibly with the $\mathcal O_E$-action, and one can choose a coordinate on $G\cong \Spf \mathcal O_E\powerseries{X}$ so that explicitly, the logarithm map is given by
\[
\mathrm{log}_G: G_E\to \mathbb G_{a,E} : X\mapsto X+\tfrac 1{\pi} X^q+\tfrac 1{\pi^2} X^{q^2} + \ldots + \tfrac 1{\pi^n} X^{q^n} + \ldots .
\]
Regarding the convergence of $\mathrm{log}_G$, we note that in fact it defines a map of rigid-analytic varieties (i.e.~adic spaces locally of finite type over $E$)
\[
\mathrm{log}_G: G_E^{\mathrm{ad}}\cong \mathbb D_E\to \mathbb G_{a,E}^{\mathrm{ad}}
\]
from the open unit disc
\[
G^{\mathrm{ad}}_E\cong \Spa \mathcal O_E\powerseries{X}\times_{\Spa \mathcal O_E} \Spa E
\]
to the adic space corresponding to $\mathbb G_a$. From the formula, one sees that in small enough discs it defines an isomorphism, and via rescaling by powers of $\pi$ (which on the level of $G^{\mathrm{ad}}_E$ defines finite \'etale covers of degree $q$, while it is an isomorphism on $\mathbb G_{a,E}^{\mathrm{ad}}$), one sees that one has an exact sequence
\[
0\to G^{\mathrm{ad}}_E[\pi^\infty]\to G^{\mathrm{ad}}_E\to \mathbb G_{a,E}^{\mathrm{ad}}\to 0
\]
on the big \'etale site of adic spaces over $\Spa E$, where $G^{\mathrm{ad}}_E[\pi^\infty]\subset G^{\mathrm{ad}}_E$ is the torsion subgroup. This is, in fact, the generic fibre of $G[\pi^\infty]=\bigcup_n G[\pi^n]$ over $\Spa \mathcal O_E$, and each $G[\pi^n] = \Spa A_n$ is represented by some finite $\mathcal O_E$-algebra $A_n$ of degree $q^n$. Inductively, $G[\pi^{n-1}]\subset G[\pi^n]$ giving a map $A_n\to A_{n-1}$; after inverting $\pi$, this is split, and the other factor is a totally ramified extension $E_n|E$. Then
\[
G^{\mathrm{ad}}_E[\pi^\infty] = \bigcup_n \Spa A_n[\tfrac 1\pi] = \bigsqcup_n \Spa E_n.
\]

We also need the ``universal cover'' of $G$, defined as
\[
\tilde{G} = \varprojlim_{\times \pi} G\cong \mathrm{Spf} \mathcal O_E\powerseries{\tilde{X}^{1/p^\infty}},
\]
where the inverse limit is over the multiplication by $\pi$ maps. The isomorphism with $\mathrm{Spf} \mathcal O_E\powerseries{\tilde{X}^{1/p^\infty}}$ is evident modulo $\pi$, but as this gives a perfect algebra, we see that in fact the isomorphism lifts uniquely to $\mathcal O_E$. Explicitly, the coordinate $\tilde{X}$ is given by
\[
\tilde{X} = \lim_{n\to \infty} X_n^{q^n}
\]
where $X_n$ is the coordinate on the $n$-th copy of $G$ in the formula $\tilde{G}=\varprojlim_{\times \pi} G$; in fact, $\tilde{X}\equiv X_n^{q^n}$ modulo $\pi^n$. In particular, the logarithm map
\[
\log_G: \tilde{G}_E\to G_E\to \mathbb G_{a,E}
\]
is given by the series
\[
\sum_{i\in \mathbb Z} \pi^i \tilde{X}^{q^{-i}} = \mathrm{lim}_{n\to \infty} \pi^n \mathrm{log}_G(X_n) = \mathrm{lim}_{n\to \infty} \mathrm{log}_G([\pi^n]_G(X_n)).
\]

Note that for any $\pi$-adically complete $\mathcal O_E$-algebra $A$, one has
\[
\tilde{G}(A)\cong \tilde{G}(A/\pi) = \mathrm{Hom}_{\mathcal O_E}(E/\mathcal O_E,G(A/\pi))[\tfrac 1\pi].
\]
Indeed, the first equality follows from $\mathcal O_E\powerseries{\tilde{X}^{1/p^\infty}}$ being relatively perfect over $\mathcal O_E$, and the second equality by noting that any element of $G(A/\pi)$ is $\pi^n$-torsion for some $n$. A different description based on $\tilde{G} = \Spf \mathcal O_E\powerseries{\tilde{X}^{1/p^\infty}}$ is
\[
\tilde{G}(A) = \varprojlim_{x\mapsto x^p} A^{\circ\circ} = A^{\flat,\circ\circ}\subset A^\flat,
\]
the subset of topologically nilpotent elements of the tilt.

This is related to the line bundle $\mathcal O_{X_S}(1)$ as follows.

\begin{proposition}\label{prop:lubintateO1} Let $S=\Spa(R,R^+)$ be an affinoid perfectoid space over $\Fq$ and let $S^\sharp=\Spa(R^\sharp,R^{\sharp +})$ be an untilt of $S$ over $E$, giving rise to the closed immersion $S^\sharp\hookrightarrow X_S$. Let $\mathcal O_{X_S}(1)$ be the line bundle on $X_S$ corresponding to the isocrystal $(E,\pi^{-1})$. Then the map
\[
\tilde{G}(R^{\sharp +})\cong R^{\circ\circ}\to H^0(Y_S,\mathcal O_{Y_S}): X\mapsto \sum_{i\in \mathbb Z} \pi^i [X^{q^{-i}}]
\]
defines a natural isomorphism
\[
\tilde{G}(R^{\sharp +})\cong H^0(X_S,\mathcal O_{X_S}(1)) = H^0(Y_S,\mathcal O_{Y_S})^{\phi=\pi}.
\]
Under this isomorphism, the map
\[
H^0(X_S,\mathcal O_{X_S}(1))\to H^0(S^\sharp,\mathcal O_{S^\sharp})=R^\sharp
\]
of evaluation at $S^\sharp$ is given by the logarithm map
\[
\mathrm{log}_G: \tilde{G}(R^{\sharp +})\to G(R^{\sharp +})\to R^\sharp.
\]
\end{proposition}

\begin{proof} The compatibility with the logarithm map is clear from the explicit formulas. Assume first that $E$ is of characteristic $p$. Then $H^0(Y_S,\mathcal O_{Y_S})$, where $Y=\mathbb D^\ast_S$ is a punctured open unit disc over $S$, can be explicitly understood as certain power series $\sum_{i\in \mathbb Z} r_i \pi^i$ with coefficients $r_i\in R$ (subject to convergence conditions as $i\to \pm \infty$). Then
\[
H^0(X_S,\mathcal O_{X_S}(1)) = H^0(Y_S,\mathcal O_{Y_S})^{\phi=\pi}
\]
amounts to those series such that $r_i=r_{i+1}^q$ for all $i\in \mathbb Z$. Thus, all $r_i$ are determined by $r_0$, which in turn can be any topologically nilpotent element of $R$. This gives the desired isomorphism
\[
H^0(X_S,\mathcal O_{X_S}(1))\cong R^{\circ\circ} = \tilde{G}(R^+) = \tilde{G}(R^{\sharp +})
\]
(as $\tilde{G}=\mathrm{Spa} \mathcal O_E\powerseries{\tilde{X}^{1/p^\infty}}$ and $R^\sharp = R$).

If $E$ is $p$-adic, then we argue as follows. First, as in the proof of Proposition~\ref{prop:standardbanachcolmez} below, one can rewrite $H^0(X_S,\mathcal O_{X_S}(1))$ as $B_{R,[1,\infty]}^{\phi=\pi}$ where
\[
B_{R,[1,\infty]} = \mathcal O(Y_{[1,\infty]}),\ \mathrm{for}\ Y_{[1,\infty]} = \{|[\varpi]|\leq |\pi|\neq 0\}\subset \Spa W_{\mathcal O_E}(R^+).
\]
By the contracting property of Frobenius, one can also replace $B_{R,[1,\infty]}$ with the crystalline period ring $B_{\mathrm{crys}}^+$ of $R^{\sharp +}/\pi$ here, and then \cite[Theorem A]{ScholzeWeinstein} gives the desired
\[
B_{R,[1,\infty]}^{\phi=\pi} = \mathrm{Hom}_{\mathcal O_E}(E/\mathcal O_E,G(R^{\sharp +}/\pi))[\tfrac 1\pi] = \tilde{G}(R^{\sharp +}/\pi) = \tilde{G}(R^{\sharp +}).
\]
That this agrees with the explicit formula follows from \cite[Lemma 3.5.1]{ScholzeWeinstein}.
\end{proof}

Recall also that the field $E_\infty$ obtained as the completion of the union of all $E_n$ is perfectoid --- in fact, one has a closed immersion $\Spf \mathcal O_{E_\infty}\hookrightarrow \tilde{G} = \Spf \mathcal O_E\powerseries{\tilde{X}^{1/p^\infty}}$, which induces an isomorphism $\Spf \mathcal O_{E^\infty}^\flat\cong \Spf \Fq\powerseries{X^{1/p^\infty}}$. Over $E_\infty$, we have an isomorphism $\mathcal O_E\cong (T_\pi G)(\mathcal O_{E_\infty})\subset \tilde{G}(\mathcal O_{E_\infty})$. By the last proposition, if $S^\sharp$ lives over $E_\infty$, we get a nonzero section of $\mathcal O_{X_S}(1)$, vanishing at $S^\sharp\subset X_S$.

\begin{proposition}\label{prop:lubintateuntilt} For any perfectoid space $S$ with untilt $S^\sharp$ over $E_\infty$, the above construction defines an exact sequence
\[
0\to \mathcal O_{X_S}\to \mathcal O_{X_S}(1)\to \mathcal O_{S^\sharp}\to 0
\]
of $\mathcal O_{X_S}$-modules.
\end{proposition}

\begin{proof} The above constructions show that one has a map $\mathcal O_{X_S}\to \mathcal I(1)$ where $\mathcal I\subset \mathcal O_{X_S}$ is the ideal sheaf of $S^\sharp$, which by Proposition~\ref{prop:degreeonecartier1} is a line bundle. To see that this map is an isomorphism, it suffices to check on geometric points, so we can assume that $S=\Spa C$ for some complete algebraically closed extension $C$ of $\Fq$. We have now fixed some nonzero global section of $\mathcal O_{X_S}(1)$, which by Proposition~\ref{prop:lubintateO1} corresponds to some nonzero topologically nilpotent $\tilde{X}\in C$; explicitly this section is given by
\[
f=\sum_{i\in \Z} \pi^i [\tilde{X}^{q^{-i}}]\in H^0(Y_C,\mathcal O_{Y_C})^{\phi=\pi}.
\]
This is the base change of the function
\[
\sum_{i\in \Z} \pi^i \tilde{X}^{q^{-i}}\in \mathcal O((\Spa \mathcal O_E\powerseries{\tilde{X}^{1/p^\infty}})_E\setminus V(\tilde{X}))
\]
under the induced map
\[
Y_C\to \Spa \mathcal O_E\powerseries{\tilde{X}^{1/p^\infty}})_E\setminus V(\tilde{X}),
\]
so it is enough to determine the vanishing locus of this function. But note that under the identification $\tilde{G} = \Spf \mathcal O_E\powerseries{\tilde{X}^{1/p^\infty}}$, this is precisely the logarithm function
\[
\mathrm{log}_G: \tilde{G}^{\mathrm{ad}}_E\setminus \{0\}\to \mathbb G_{a,E}^{\mathrm{ad}};
\]
thus, it is enough to determine the vanishing locus of the logarithm function. But this is precisely
\[
\bigsqcup_n \Spa E_\infty\subset \tilde{G}^{\mathrm{ad}}_E\setminus \{0\},
\]
with a simple zero at each of these points. This gives exactly the claimed statement.
\end{proof}

\begin{corollary}[{\cite[Proposition 2.12]{FarguesClassFieldTheory}}]\label{cor:descriptiondiv1} There is a well-defined map $\BC(\mathcal O(1))\setminus \{0\}\to \mathrm{Div}^1$ sending a nonzero section $f\in H^0(X_S,\mathcal O_{X_S}(1))$ to the closed Cartier divisor given by $V(f)$. This descends to an isomorphism
\[
(\BC(\mathcal O(1))\setminus \{0\})/\underline{E^\times}\cong \mathrm{Div}^1.
\]
\end{corollary}

\begin{proof} Note that $\BC(\mathcal O(1))\cong \Spd \Fq\powerseries{X^{1/p^\infty}}$ by Proposition~\ref{prop:lubintateO1}, and hence $\BC(\mathcal O(1))\setminus \{0\}\cong \Spa \Fq\laurentseries{X^{1/p^\infty}}$ is representable by a perfectoid space. In fact, it is naturally isomorphic to $\Spd E_\infty=\Spa E_\infty^\flat$, and the previous proposition ensures that the map to $\mathrm{Div}^1$ is well-defined and corresponds to the projection $\Spd E_\infty\to \Spd E\to \Spd E/\phi^\Z = \mathrm{Div}^1$. Here, the first map $\Spd E_\infty\to \Spd E$ is a quotient under $\underline{\mathcal O_E^\times}$, and the second map $\Spd E\to \Spd E/\phi^\Z$ then corresponds to the quotient by $\pi^\Z$, as $\phi=\pi$ on $\BC(\mathcal O(1))$.
\end{proof}

In particular, if one works on $\Perf_k$, then $\mathrm{Div}^1 = \Spd \breve E/\phi^\Z$, whose $\pi_1^{\et}$ is given by the absolute Galois group of $E$. On the other hand, the preceding gives a canonical $\underline{E^\times}$-torsor, giving a natural map from the absolute Galois group of $E$ to the profinite completion of $E^\times$. By comparison with Lubin--Tate theory, this is the usual {\it Artin reciprocity map}, see \cite[Section 2.3]{FarguesClassFieldTheory} for more details.

\subsection{Absolute Banach--Colmez spaces}
\label{sec:absoluteBanachColmez}

In this section, we analyze the Banach--Colmez spaces in the case $\mathcal E=\mathcal E(D)$ for some isocrystal $D=(D,\phi)$. We then sometimes write $\BC(D)$ and $\BC(D[1])$ for the corresponding functors on $\Perf_k$; or also $\BC(\mathcal O(\lambda))$, $\BC(\mathcal O(\lambda)[1])$ for $\lambda\in\mathbb Q$ when $D=D_{-\lambda}$. These are in fact already defined for all $S\in \Perf_{\Fq}$.

\begin{proposition}\label{prop:standardbanachcolmez} Let $\lambda\in \mathbb Q$.
\begin{altenumerate}
\item[{\rm (i)}] If $\lambda<0$, then $H^0(X_S,\mathcal O_{X_S}(\lambda))=0$ for all $S\in \Perf_{\Fq}$. Moreover, the projection from
\[
\BC(\mathcal O(\lambda)[1]): S\mapsto H^1(X_S,\mathcal O_{X_S}(\lambda))
\]
to the point $\ast$ is relatively representable in locally spatial diamonds, partially proper, and cohomologically smooth.
\item[{\rm (ii)}] For $\lambda=0$, the map
\[
\underline{E}\to \BC(\mathcal O)
\]
is an isomorphism of pro-\'etale sheaves, and the pro-\'etale sheafification of $S\mapsto H^1(X_S,\mathcal O_{X_S})$ vanishes. In particular, for all $S$ one gets an isomorphism
\[
R\Gamma_\proet(S,\underline{E})\to R\Gamma(X_S,\mathcal O_{X_S}).
\]
\item[{\rm (iii)}] For $\lambda>0$, one has $H^1(X_S,\mathcal O_{X_S}(\lambda))=0$ for all affinoid $S\in \Perf_{\Fq}$, and the projection from
\[
\BC(\mathcal O(\lambda)): S\mapsto H^0(X_S,\mathcal O_{X_S}(\lambda))
\]
to the point $\ast$ is relatively representable in locally spatial diamonds, partially proper, and cohomologically smooth.
\item[{\rm (iv)}] If $0<\lambda\leq [E:\mathbb Q_p]$ (resp.~for all positive $\lambda$ if $E$ is of equal characteristic), there is an isomorphism
\[
\BC(\mathcal O(\lambda))\cong \Spd k\powerseries{x_1^{1/p^\infty},\ldots,x_r^{1/p^\infty}}
\]
where $\lambda=r/s$ with coprime integers $r,s>0$.
\end{altenumerate}
\end{proposition}

\begin{proof} For all statements, we can reduce to the case $\lambda=n\in \mathbb Z$ by replacing $E$ by its unramified extension of degree $s$. Regarding the vanishing of $H^1(X_S,\mathcal O_{X_S}(n))$ for $n>0$ and $S=\Spa(R,R^+)$ affinoid, pick a pseudouniformizer $\varpi\in R$. In terms of the presentation of $X_S$ as gluing $Y_{S,[1,q]}$ along $\phi: Y_{S,[1,1]}\cong Y_{S,[q,q]}$, it suffices to see that
\[
\phi-\pi^n: B_{R,[1,q]}\to B_{R,[1,1]}
\]
is surjective. Any element of $B_{R,[1,1]}$ can be written as the sum of an element of $B_{R,[0,1]}[\tfrac 1{\pi}]$ and an element of $[\varpi] B_{R,[1,\infty]}$. Here the rings $B_{R,[0,1]} = \mathcal O(Y_{S,[0,1]})$ and $B_{R,[1,\infty]} = \mathcal O(Y_{S,[1,\infty]})$ correspond to the affinoid subsets
\[
Y_{S,[0,1]} = \{|\pi|\leq |[\varpi]|\neq 0\}\subset \Spa W_{\mathcal O_E}(R^+)
\]
resp.~
\[
Y_{S,[1,\infty]} = \{|[\varpi]|\leq |\pi|\neq 0\}\subset \Spa W_{\mathcal O_E}(R^+).
\]
(We warn the reader that $Y_{S,[0,1]}$ and $Y_{S,[1,\infty]}$ are not contained in $Y_S=Y_{S,(0,\infty)}$; we hope this clash of notation will not cause confusion.) If $f\in B_{R,[0,1]}$, then the series
\[
g=\phi^{-1}(f) + \pi^n \phi^{-2}(f) + \pi^{2n} \phi^{-3}(f) + \ldots
\]
converges in $B_{R,[0,q]}$ (with its evident definition) and thus in $B_{R,[1,q]}$, and $f=\phi(g) - \pi^n g$. The same then applies to elements of $B_{R,[0,1]}[\tfrac 1{\pi}]$. On the other hand, if $f\in [\varpi] B_{R,[1,\infty]}$, then the series
\[
g=-\pi^{-n} f - \pi^{-2n} \phi(f) - \pi^{-3n} \phi^2(f) - \ldots
\]
converges in $B_{R,[1,q]}$, and $f=\phi(g) - \pi^n g$.

In fact, the same arguments prove that the map
\[
[B_{R,[1,\infty]}\xrightarrow{\phi-\pi^n} B_{R,[1,\infty]}]\to [B_{R,[1,q]}\xrightarrow{\phi-\pi^n} B_{R,[1,1]}]
\]
is a quasi-isomorphism. Indeed, we have a short exact sequence
\[
0\to W_{\mathcal O_E}(R^+)[\tfrac 1{\pi}]\to B_{R,[1,\infty]}\oplus B_{R,[0,q]}[\tfrac 1{\pi}]\to B_{R,[1,q]}\to 0
\]
(obtained from sheafyness of $W_{\mathcal O_E}(R^+)[\tfrac 1{\pi}]$ when endowed with the $\pi$-adic topology on $W_{\mathcal O_E}(R^+)$), and similarly
\[
0\to W_{\mathcal O_E}(R^+)[\tfrac 1{\pi}]\to B_{R,[1,\infty]}\oplus B_{R,[0,1]}[\tfrac 1{\pi}]\to B_{R,[1,1]}\to 0.
\]
Therefore, it suffices to see that the maps
\[
B_{R,[0,q]}[\tfrac 1{\pi}]\xrightarrow {\phi-\pi^n} B_{R,[0,1]}[\tfrac 1{\pi}]
\]
and
\[
W_{\mathcal O_E}(R^+)[\tfrac 1{\pi}]\xrightarrow {\phi-\pi^n} W_{\mathcal O_E}(R^+)[\tfrac 1{\pi}]
\]
are isomorphisms. In both cases, this follows from convergence of $\phi^{-1}+\pi^n \phi^{-2} + \pi^{2n} \phi^{-3} + \ldots$ on these algebras, giving an explicit inverse.

For part (iv), note that in equal characteristic one can describe $\mathcal O(Y_{S,I}) = B_{R,I}$, for $S=\Spa(R,R^+)$ affinoid, explicitly as power series $\sum_{i\in \mathbb Z} r_i \pi^i$ with $r_i\in R$, satisfying some convergence conditions as $i\to \pm \infty$. Taking the part where $\phi=\pi^n$, we require $\phi(r_i) = r_{i+n}$, and we see that we can freely choose $r_1,\ldots,r_n$. The required convergence holds precisely when all $r_i$ are topologically nilpotent, giving the isomorphism in that case. If $E$ is $p$-adic, we can reduce to $E=\mathbb Q_p$ (but now $\lambda$ rational, $0<\lambda\leq 1$), taking a pushforward of the sheaf along $X_{S,E} = X_{S,\mathbb Q_p}\times_{\mathbb Q_p} E\to X_{S,\mathbb Q_p}$. In that case, the result follows from the equality $H^0(X_S,\mathcal O_{X_S}(\lambda)) = B_{R,[1,\infty]}^{\phi^r=p^s}$ proved above, and \cite[Theorem A, Proposition 3.1.3 (iii)]{ScholzeWeinstein}.

In particular, for affinoid $S$ we can choose a fibrewise nonzero map $\mathcal O_{X_S}\to \mathcal O_{X_S}(1)$, by taking a map $S\to \BC(\mathcal O(1))\cong \Spd \Fq\powerseries{x^{1/p^\infty}}$ sending $x$ to a pseudouniformizer. By Proposition~\ref{prop:lubintateuntilt}, for any $n\in \mathbb Z$, we get an exact sequence
\[
0\to \mathcal O_{X_S}(n)\to \mathcal O_{X_S}(n+1)\to \mathcal O_{S^\sharp}\to 0.
\]
Applying this for $n>0$, we get inductively an exact sequence
\[
0\to \BC(\mathcal O(n))|_S\to \BC(\mathcal O(n+1))|_S\to (\mathbb A^1_{S^\sharp})^\diamond\to 0.
\]
Starting with the base case $n=1$ already handled, this allows one to prove part (iii) by induction, using \cite[Proposition 23.13]{ECoD}.

Now for part (ii), we use the sequence for $n=0$. In that case, for $S=\Spa(R,R^+)$, we get an exact sequence
\[
0\to H^0(X_S,\mathcal O_{X_S})\to H^0(X_S,\mathcal O_{X_S}(1))\to R^\sharp\to H^1(X_S,\mathcal O_{X_S})\to 0
\]
where the map in the middle can be identified with the logarithm map of the universal cover of the Lubin--Tate formal group. This is pro-\'etale locally surjective, with kernel given by $\underline{E}$, proving (ii).

Finally, for part (i), we first treat the case $n=-1$, where we get an exact sequence
\[
0\to \underline{E}\to (\mathbb A^1_{S^\sharp})^\diamond\to \BC(\mathcal O(-1)[1])|_S\to 0
\]
showing in particular the vanishing of $H^0(X_S,\mathcal O_{X_S}(-1))=0$. As $\underline{E}\to (\mathbb A^1_{S^\sharp})^\diamond$ is a closed immersion, the result follows from \cite[Proposition 24.2]{ECoD}. Now for $n<-1$, the result follows by induction from the sequence
\[
0\to (\mathbb A^1_{S^\sharp})^\diamond\to \BC(\mathcal O(-n)[1])|_S\to \BC(\mathcal O(-n+1)[1])|_S\to 0
\]
and \cite[Proposition 23.13]{ECoD}.
\end{proof}

\subsection{The algebraic curve}

We recall the following important ampleness result.

\begin{theorem}[{\cite[Proposition 6.2.4]{KedlayaLiu1}}]\label{thm:O1ample} Let $S=\Spa(R,R^+)$ be an affinoid perfectoid space over $\mathbb F_q$ and let $\mathcal E$ be any vector bundle on $X_S$. Then there is an integer $n_0$ such that for all $n\geq n_0$, the vector bundle $\mathcal E(n)$ is globally generated, i.e.~there is a surjective map
\[
\mathcal O_{X_S}^m\to \mathcal E(n)
\]
for some $m\geq 0$, and moreover $H^1(X_S,\mathcal E(n))=0$.
\end{theorem}

\begin{proof} Pick a pseudouniformizer $\varpi\in R$, thus defining a radius function on $Y_S$. Write $X_S$ as the quotient of $Y_{S,[1,q]}$ along the isomorphism $\phi: Y_{S,[1,1]}\cong Y_{S,[q,q]}$. Correspondingly, $\mathcal E$ is given by some finite projective $B_{R,[1,q]}$-module $M_{[1,q]}$, with base changes $M_{[1,1]}$ and $M_{[q,q]}$ to $B_{R,[1,1]}$ and $B_{R,[q,q]}$, and an isomorphism $\phi_M: M_{[q,q]}\cong M_{[1,1]}$, linear over $\phi: B_{R,[q,q]}\cong B_{R,[1,1]}$.

For convenience, we first reduce to the case that $M_{[1,q]}$ is free (cf.~\cite[Corollary 1.5.3]{KedlayaLiu1}). Indeed, pick a surjection $\psi: F_{[1,q]}:=B_{R,[1,q]}^m\to M_{[1,q]}$. We want to endow the source with a similar $\phi$-module structure $\phi_F: F_{[q,q]}\cong F_{[1,1]}$ (with obvious notation), making $\psi$ equivariant. For this, we would like to find a lift
\[\xymatrix{
F_{[q,q]}\ar@{-->}[r]^{\phi_F}\ar[d]^\psi & F_{[1,1]}\ar[d]^\psi\\
M_{[q,q]}\ar[r]^{\phi_M} & M_{[1,1]}
}\]
such that $\phi_F$ is an isomorphism. Let $N_{[1,q]}=\ker(\psi)$, with base change $N_{[q,q]}$, $N_{[1,1]}$. Choosing a splitting $F_{[1,q]}\cong M_{[1,q]}\oplus N_{[1,q]}$, we see that we could find $\phi_F$ if there is an isomorphism $\phi^\ast N_{[q,q]}\cong N_{[1,1]}$ of $B_{R,[1,1]}$-modules. But in the Grothendieck group of finite projective $B_{R,[1,1]}$-modules, both are given by $[B_{R,[1,1]}^m] - [M_{[1,1]}]$. But equality in the Grothendieck group is the same thing as stable isomorphism; thus, after possibly adding a free module (i.e.~increasing $m$), they are isomorphic, giving the claim.

Thus, we can assume that $M_{[1,q]}\cong B_{R,[1,q]}^m$ is a free $B_{R,[1,q]}$-module, and then
\[
\phi_M=A^{-1}\phi
\]
for some matrix $A\in \GL_m(B_{R,[1,1]})$. Actually, repeating the above argument starting with the presentation of $X_S$ as the quotient of $Y_{S,[1,q^2]}$ via identifying $Y_{S,[q,q^2]}$ with $Y_{S,[1,q]}$, one can ensure that
\[
A\in \GL_m(B_{R,[1,q]}).
\]
Twisting by $\mathcal O_{X_S}(n)$ amounts to replacing $A$ by $A\pi^n$. Let us choose integers $N$ and $N'$ such that
\begin{itemize}
\item the matrix $A$ has entries in $\pi^N W_{\mathcal O_E}(R^+)\langle (\tfrac{[\varpi]}{\pi})^{\pm 1}\rangle$ 
\item the matrix $A^{-1}$ has entries in $\pi^{-N'} W_{\mathcal O_E}(R^+)\langle \tfrac{\pi}{[\varpi]^{1/q}},\tfrac{[\varpi]}{\pi}\rangle$.
\end{itemize}
By twisting, we can replace $N$ and $N'$ by $N+n$ and $N'+n$; we can thus arrange that $qN>N'$, $N>0$.

Fix some rational $r$ such that $1<r\leq q$. We will now show that there are $m$ elements
\[
v_1,\ldots,v_m\in (B_{R,[1,q]}^m)^{\phi=A}=H^0(X_S,\mathcal E)
\]
that form a basis of $B_{R,[r,q]}^m$. Repeating the above analysis for different strips (and different choices of pseudouniformizers $\varpi\in R$ to get overlapping strips), we can then get global generation of $\mathcal E$.

In fact, we will choose $v_i$ to be of the form $[\varpi]^M e_i - v_i'$, for some positive integer $M$ chosen later, where $e_i\in B_{R,[1,q]}^m$ is the $i$-th basis vector and $v_i'$ is such that
\[
||v_i'||_{B_{R,[r,q]}}\leq ||[\varpi]^{M+1}||_{B_{R,[r,q]}} = q^{-M-1}.
\]
Here, we endow all $B_{R,I}$ with the spectral norm, normalizing the norms on all completed residue fields via $||[\varpi]||=\tfrac 1q$. These $v_1,\ldots,v_m$ restrict to a basis of $B_{R,[r,q]}^m$ since the base change matrix from the canonical basis is given by an element of
\[
[\varpi^M] (\mathrm{Id} + [\varpi] M_m (B_{R,[r,q]}^\circ))\subset \GL_m (B_{R,[r,q]}).
\]

In order to find the $v_i'$, it suffices to prove that the map
\[
\phi - A: B_{R,[1,q]}^m\to B_{R,[1,1]}^m
\]
is surjective (yielding $H^1(X_S,\mathcal E)=0$), in the following quantitative way: If, for some positive integer $M$ chosen later,
\[
w\in \pi^M W_{\mathcal O_E}(R^+)\langle (\tfrac{[\varpi]}{\pi})^{\pm 1}\rangle^m\subset B_{R,[1,1]}^m,
\]
then there is some
\[
v\in B_{R,[1,q]}^m
\]
such that
\begin{equation}\label{eq:estimee1}
(\phi-A)v=w\ \mathrm{and}\ ||v||_{B_{R,[r,q]}}\leq q^{-M-1}.
\end{equation}
Indeed, we can then apply this to $w_i=(\phi-A)([\varpi]^M e_i)$ (since $N>0$ and thus $A$ has entries in $W_{\O_E}(R^+)\langle (\tfrac{[\varpi]}{\pi})^{\pm 1}\rangle$), getting some $v_i'$ with $w_i=(\phi-A)(v_i')$ and
\[
||v_i'||_{B_{R,[r,q]}}\leq q^{-M-1},
\]
as desired.

Thus, take any
\[
w\in \pi^M W_{\mathcal O_E}(R^+)\langle (\tfrac{[\varpi]}{\pi})^{\pm 1}\rangle^m.
\]
We can write
\[
w=w_1+w_2\ \mathrm{where}\ w_1\in [\varpi]^{N-1} \pi^{M-N+1} W_{\mathcal O_E}(R^+)\langle \tfrac{\pi}{[\varpi]}\rangle^m,\ w_2\in [\varpi]^N \pi^{M-N} W_{\mathcal O_E}(R^+)\langle \tfrac{[\varpi]}{\pi}\rangle^m .
\]
Let
\[
v=\phi^{-1}(w_1)-A^{-1}w_2\in B_{R,[1,q]}^m\ \mathrm{so}\ \mathrm{that}\ w':=w-\phi(v)+Av=\phi(A^{-1}w_2)+A\phi^{-1}(w_1).
\]
Note that (as $N>0$)
\[
A\phi^{-1}(w_1)\in \pi^N W_{\mathcal O_E}(R^+)\langle (\tfrac{[\varpi]}{\pi})^{\pm 1}\rangle\cdot [\varpi]^{(N-1)/q} \pi^{M-N+1} W_{\mathcal O_E}(R^+)\langle \tfrac{\pi}{[\varpi]^{1/q}}\rangle^m\subset \pi^{M+1} W_{\mathcal O_E}(R^+)\langle (\tfrac{[\varpi]}{\pi})^{\pm 1}\rangle^m
\]
and also (as $qN>N'$)
\[
\phi(A^{-1}w_2)\in \pi^{-N'} [\varpi]^{Nq} \pi^M W_{\mathcal O_E}(R^+)\langle \tfrac{[\varpi]^q}{\pi}\rangle^m\subset \pi^{M+1} W_{\mathcal O_E}(R^+)\langle (\tfrac{[\varpi]}{\pi})^{\pm 1}\rangle^m
\]
so that
\[
w'\in \pi^{M+1} W_{\mathcal O_E}(R^+)\langle (\frac{[\varpi]}{\pi})^{\pm 1}\rangle^m.
\]
If one can thus prove the required bounds on $v$, this process will converge and prove the desired statement. It remains to estimate $v$. On the one hand, its norm is clearly bounded in terms of the norm of $w$ (as both $w_1$ and $w_2$ are, and $\phi^{-1}$ and $A^{-1}$ are bounded operators), and thus, since when one iterates $w$ goes to zero, $v$ goes to zero, and the process converges by summing to obtain some $v$ such that $(\varphi-A)v=w$. But we need an improved estimate over $B_{R,[r,q]}$ to obtain \eqref{eq:estimee1}. Note that the norm of $\phi^{-1}(w_1)$ is bounded above by the norm of $[\varpi]^{(N-1)/q} \pi^{M-N+1}$, which in $B_{R,[r,q]}$ is given by $q^{-(N-1)/q-rM+rN-r}$. This is at most $q^{-M-1}$ once $M$ is large enough. On the other hand, $w_2\in \pi^M W_{\mathcal O_E}(R^+)\langle \frac{[\varpi]}{\pi}\rangle^m$ and so the norm of $A^{-1} w_2$ is bounded by the norm of $\pi^{-N'}\pi^M$, which in $B_{R,[r,q]}$ is given by $q^{rN'-rM}$. Again, this is at most $q^{-M-1}$ once $M$ is large enough. Thus, taking $M$ large enough (depending only on $N$, $N'$ and $r>1$), the process above converges, giving the desired result.
\end{proof}

We have the following general {\it GAGA theorem}. Its proof is an axiomatization of \cite[Theorem 6.3.9]{KedlayaLiu1}.

\begin{proposition}[GAGA]\label{prop:ampleGAGA} Let $(X,\mathcal O_X)$ be a locally ringed spectral space equipped with a line bundle $\mathcal O_X(1)$ such that for any vector bundle $\mathcal E$ on $X$, there is some $n_0$ such that for all $n\geq n_0$, the bundle $\mathcal E(n)$ is globally generated. Moreover, assume that for $i>0$, the cohomology group $H^i(X,\mathcal E(n))=0$ vanishes for all sufficiently large $n$.

Let $P=\bigoplus_{n\geq 0} H^0(X,\mathcal O_X(n))$ be the graded ring and $X^{\mathrm{alg}} = \Proj(P)$. There is a natural map $(X,\mathcal O_X)\to X^{\mathrm{alg}}$ of locally ringed topological spaces, and pullback along this map induces an equivalence of categories between vector bundles on $X^{\mathrm{alg}}$ and vector bundles on $(X,\mathcal O_X)$. Moreover, for any vector bundle $\mathcal E^{\mathrm{alg}}$ on $X^{\mathrm{alg}}$ with pullback $\mathcal E$ to $X$, the map
\[
H^i(X^\mathrm{alg},\mathcal E^{\mathrm{alg}})\to H^i(X,\mathcal E)
\]
is an isomorphism for all $i\geq 0$.
\end{proposition}

Recall that for any graded ring $P=\bigoplus_{n\geq 0} P_n$, one can define a separated scheme $\Proj(P)$ by gluing $\Spec P[f^{-1}]_0$ for all $f\in P_n$, $n>0$, where $P[f^{-1}]_0 = \varinjlim_i f^{-i} P_{in}$ is the degree $0$ part of $P[f^{-1}]$. In our situation, if $n$ is large enough so that $\mathcal O_X(n)$ is globally generated, then it is enough to consider only $f\in P_n$ for this given $n$, and in fact only a finite set of them (as $X$ is quasicompact); in particular, $\Proj(P)$ is quasicompact. Moreover, one sees that there is a tautological line bundle $\mathcal O_{\Proj(P)}(n)$ for all sufficiently large $n$, compatible with tensor products; thus, there is also a tautological line bundle $\mathcal O_{\Proj(P)}(1)$, which is an ample line bundle on $\Proj(P)$. The pullback of $\mathcal O_{\Proj(P)}(1)$ is then given by $\mathcal O_X(1)$.

\begin{proof} The construction of the map $f: (X,\mathcal O_X)\to X^{\mathrm{alg}}$ is formal (and does not rely on any assumptions): if $g\in P_n$, then on the non-vanishing locus $U=D(g)\subset X$, there is an isomorphism $g_{|U}:\O_U\iso \O_U(n)$. Now, for $x=\frac{a}{g^k} \in P[g^{-1}]_0$, $g_{|U}^{-k}\circ a\in \O(U)$, and this defines a morphism of rings $P[g^{-1}]_0\to \O(U)$. One deduces a morphism of locally ringed spaces $U\to D^+(g)$, and those glue when $g$ varies to a morphism of locally ringed spaces $(X,\O_X)\to X^{\mathrm{alg}}$.

We consider the functor taking any vector bundle $\mathcal E$ on $X$ to the quasicoherent $\mathcal O_{X^{\mathrm{alg}}}$-module $\overline{\mathcal E}$ associated to the graded $P$-module $\bigoplus_{n\geq 0} H^0(X,\mathcal E(n))$. This functor is exact as $H^1(X,\mathcal E(n))=0$ for all sufficiently large $n$, and it commutes with twisting by $\mathcal O(1)$. We claim that it takes values in vector bundles on $X^{\mathrm{alg}}$. To see this, take a surjection $\mathcal O_X^m\to \mathcal E(n)$ with kernel $\mathcal F$, again a vector bundle. The map $\mathcal O_X^m\to \mathcal E(n)$ splits after twisting, i.e.~for any $f\in P_{n'}$ with $n'$ large enough, there is a map $\mathcal E(n-n')\to \mathcal O_X^m$ such that $\mathcal E(n-n')\to \mathcal O_X^m\to \mathcal E(n)$ is multiplication by $f$. Indeed, the obstruction to such a splitting is a class in $H^1(X,\sHom(\mathcal E,\mathcal F)(n'))$ which vanishes for $n'$ large enough. This implies that $\overline{\mathcal E}$ is a vector bundle on $\Spec P[f^{-1}]_0$ for any such $f$, and these cover $X^{\mathrm{alg}}$.

There is a natural map $f^\ast \overline{\mathcal E}\to \mathcal E$, and the preceding arguments show that this is an isomorphism (on the preimage of any $\Spec P[f^{-1}]_0$, and thus globally). It now remains to show that if $\mathcal E^{\mathrm{alg}}$ is any vector bundle on $X^{\mathrm{alg}}$, the map
\[
H^i(X^\mathrm{alg},\mathcal E^{\mathrm{alg}})\to H^i(X,\mathcal E)
\]
is an isomorphism for all $i\geq 0$. (Indeed, for $i=0$ this implies, by passing to internal Hom's, that $\mathcal E^{\mathrm{alg}}\mapsto \mathcal E$ is fully faithful, and we have just seen that this functor is essentially surjective.) By ampleness of $\mathcal O_{X^{\mathrm{alg}}}(1)$, there is some surjection $\mathcal O_{X^{\mathrm{alg}}}(-n)^m\to \mathcal (\mathcal E^{\mathrm{alg}})^\vee$, with kernel a vector bundle $\mathcal F$. Dualizing, we get an injection $\mathcal E^{\mathrm{alg}}\to \mathcal O_{X^{\mathrm{alg}}}(n)^m$ with cokernel a vector bundle. This already gives injectivity on $H^0$ by reduction to $\mathcal O_{X^{\mathrm{alg}}}(n)$ where it is clear. Applying this injectivity also for $\mathcal F$, we then get bijectivity on $H^0$. This already implies that we get an equivalence of categories (exact in both directions). Finally, picking $f_1,\ldots,f_m\in P_n$ so that the $\Spec P[f_i^{-1}]_0$ cover $X^{\mathrm{alg}}$, we can look at the corresponding \v{C}ech complex. Each term is a filtered colimit of global sections of vector bundles $\mathcal E(n)$ along multiplication by products of powers of $f_i$'s. This reduces the assertion to the case of $H^0$ and the vanishing of $H^i(X,\mathcal E(n))$ for $n$ large enough.
\end{proof}

\begin{remark}
One can check that $X^{\mathrm{alg}}$ is up to canonical isomorphism independent of the choice of a line bundle $\O_X(1)$ satisfying the preceding properties. 
\end{remark}

In particular, for any affinoid perfectoid space $S$ over $\mathbb F_q$, we can define the algebraic curve
\[
X_S^{\mathrm{alg}} = \mathrm{Proj} \bigoplus_{n\geq 0} H^0(X_S,\mathcal O_{X_S}(n)).
\]
There is a well-defined map $X_S\to X_S^{\mathrm{alg}}$ of locally ringed spectral spaces, pullback along which defines an equivalence of categories of vector bundles, and is compatible with cohomology.

Notably, this connects the present discussion to the original definition of the Fargues--Fontaine curve as given in \cite{FarguesFontaine}, where the case $S=\Spa(F,\mathcal O_F)$ is considered, for a perfectoid field $F$ of characteristic $p$. We will restrict ourselves, as above, to the case that $F=C$ is algebraically closed.

\begin{proposition}\label{prop:classicalpointsschematiccurve} Let $C$ be a complete algebraically closed nonarchimedean field over $\mathbb F_q$. Then $X_C^{\mathrm{alg}}$ is a connected regular noetherian scheme of Krull dimension $1$, and the map $|X_C|\to |X_C^{\mathrm{alg}}|$ induces a bijection between $|X_C|^{\mathrm{cl}}$ and the closed points of $|X_C^{\mathrm{alg}}|$. Moreover, for any classical point $x\in |X_C|$, the complement $X_C^{\mathrm{alg}}\setminus \{x\}$ is the spectrum of a principal ideal domain.
\end{proposition}

\begin{proof} Let $x\in |X_C|^{\mathrm{cl}}$ be any classical point, corresponding to some untilt $C^\sharp$ over $E$ of $C$. Using Lubin--Tate formal groups, we see that there is an exact sequence
\[
0\to \mathcal O_{X_C}\to \mathcal O_{X_C}(1)\to \mathcal O_{C^\sharp}\to 0
\]
on $X_C$. The corresponding section $f\in H^0(X_C,\mathcal O_{X_C}(1))$ defines its vanishing locus in $X_C^{\mathrm{alg}}$, which is then also given by $\Spec C^\sharp$. Indeed, this vanishing locus is affine as it is Zariski closed in the affine scheme $D^+(g)$ for any $g\in H^0(X_C,\mathcal O_{X_C}(1))$ that does not vanish at $x$; and one can compute the global sections via Proposition~\ref{prop:lubintateO1}. In particular, $x$ defines a closed point of $|X_C^{\mathrm{alg}}|$. Now we want to show that $P[f^{-1}]_0$ is a principal ideal domain. Thus, take any nonzero $g\in H^0(X_C,\mathcal O_{X_C}(n))$. This has finitely many zeroes on $X_C$, all at classical points $x_1,\ldots,x_m$. For each $x_i$, we have a section $f_{x_i}\in H^0(X_C,\mathcal O_{X_C}(1))$ as before, and then $g=f_1^{n_1}\cdots f_m^{n_m} h$ for some $n_i\geq 1$, and some $h\in H^0(X_C,\mathcal O_{X_C}(n'))$ that is everywhere nonzero. In particular, $h$ defines an isomorphism $\mathcal O_{X_C}\to \mathcal O_{X_C}(n')$, whence $n'=0$, and $h\in E^\times$. This decomposition implies easily that $P[f^{-1}]_0$ is indeed a principal ideal domain, and it shows that all maximal ideals arise from classical points of $|X_C|$, finishing the proof.
\end{proof}

\subsection{Classification of vector bundles}

At this point, we can recall the classification of vector bundles over $X_C$; so here we take $S=\Spa C$ for a complete algebraically closed nonarchimedean field $C$ over $\Fq$. First, one classifies line bundles.

\begin{proposition}\label{prop:picardFF} The map $\mathbb Z\to \mathrm{Pic}(X_C)$, $n\mapsto \mathcal O_{X_C}(n)$, is an isomorphism.
\end{proposition}

\begin{proof} By Proposition~\ref{prop:classicalpointsschematiccurve}, any line bundle becomes trivial after removing one closed point $x\in X_C^{\mathrm{alg}}$. As the local rings of $X_C^{\mathrm{alg}}$ are discrete valuation rings, this implies that any line bundle is of the form $\mathcal O_{X_C}(n[x])$ for some $n\in \mathbb Z$. But $\mathcal O_{X_C}([x])\cong \mathcal O_{X_C}(1)$ by Proposition~\ref{prop:lubintateO1}, so the result follows.
\end{proof}

In particular, one can define the degree of any vector bundle $\mathcal E$ on $X_C$ via
\[
\mathrm{deg}(\mathcal E) = \mathrm{deg}(\mathrm{det}(\mathcal E))\in \mathbb Z
\]
where $\mathrm{det}(\mathcal E)$ is the determinant line bundle, and $\mathrm{deg}: \mathrm{Pic}(X_C)\cong \mathbb Z$ is the isomorphism from the proposition. Of course, one can also define the rank $\mathrm{rk}(\mathcal E)$ of any vector bundle, and thus for any nonzero vector bundle its slope
\[
\mu(\mathcal E)=\frac{\mathrm{deg}(\mathcal E)}{\mathrm{rk}(\mathcal E)}\in \mathbb Q.
\]

It is easy to see that this satisfies the Harder--Narasimhan axiomatics \cite[5.5.1]{FarguesFontaine} (for example, rank and degree are additive in short exact sequences). In particular, one can define semistable vector bundles as those vector bundles $\mathcal E$ such that for all proper nonzero $\mathcal F\subset \mathcal E$, one has $\mu(\mathcal F)\leq \mu(\mathcal E)$. One says that $\mathcal E$ is stable if in fact $\mu(\mathcal F)<\mu(\mathcal E)$ for all such $\mathcal F$.

\begin{example} For any $\lambda\in \mathbb Q$, the bundle $\mathcal O_{X_C}(\lambda)$ is stable of slope $\lambda$. Indeed, assume that $0\neq \mathcal F\subsetneq \mathcal O_{X_C}(\lambda)$ is a proper nonzero subbundle, and let $r=\mathrm{rk}(\mathcal F)$, $s=\mathrm{deg}(\mathcal F)$. Passing to $r$-th wedge powers, we get an injection
\[
\det(\mathcal F)\cong \mathcal O_{X_C}(s)\hookrightarrow \mathcal O_{X_C}(r\lambda)^m,
\]
using that $\bigwedge^r \mathcal O_{X_C}(\lambda)$ is a direct sum of copies of $\mathcal O_{X_C}(r\lambda)$. This implies that $s\leq r\lambda$. Moreover, if we have equality, then $r$ is at least the denominator of $\lambda$, which is the rank of $\mathcal O_{X_C}(\lambda)$, i.e.~$\mathcal F$ has the same rank as $\mathcal O_{X_C}(\lambda)$. Thus, $\mathcal O_{X_C}(\lambda)$ is stable.
\end{example}

\begin{proposition}\label{prop:hardernarasimhan} Any vector bundle $\mathcal E$ on $X_C$ admits a unique exhaustive separating $\mathbb Q$-indexed filtration by saturated subbundles $\mathcal E^{\geq \lambda}\subset \mathcal E$, called the Harder--Narasimhan filtration, such that
\[
\mathcal E^\lambda := \mathcal E^{\geq \lambda}/\mathcal E^{>\lambda},\ \mathrm{where}\ \mathcal E^{>\lambda}=\bigcup_{\lambda'>\lambda} \mathcal E^{\geq \lambda'},
\]
is semistable of slope $\lambda$. The formation of the Harder--Narasimhan filtration is functorial in $\mathcal E$.$\hfill \Box$
\end{proposition}

As a preparation for the next theorem, we note that the Harder--Narasimhan filtration is also compatible with change of $C$.

\begin{proposition}\label{prop:hardernarasimhanbasechange} Let $\mathcal E$ be a vector bundle on $X_C$, and let $C'|C$ be an extension of complete algebraically closed nonarchimedean fields, with pullback $\mathcal E'$ of $\mathcal E$ to $X_{C'}$. Then $(\mathcal E')^{\geq \lambda}$ is the pullback of $\mathcal E^{\geq \lambda}$.

Similarly, if $E'|E$ is a finite separable extension of degree $r$, and $\mathcal E'$ is the pullback of $\mathcal E$ along $X_{C,E'} = X_{C,E}\otimes_E E'\to X_{C,E}=X_C$, then $(\mathcal E')^{\geq \lambda}$ is the pullback of $\mathcal E^{\geq \lambda/r}$.
\end{proposition}

\begin{proof} Consider the case of $C'|C$. By uniqueness of the Harder--Narasimhan filtration, it suffices to see that pullbacks of semistable vector bundles remain semistable. Thus, assume that $\mathcal E$ is semistable, and assume by way of contradiction that $\mathcal E'$ is not semistable. By induction on the rank, we can assume that the formation of the Harder--Narasimhan filtration of $\mathcal E'$ is compatible with any base change. Consider the first nontrivial piece of the Harder--Narasimhan filtration $0\neq \mathcal F\subsetneq \mathcal E'$. This is a vector bundle on $X_{C'}$ with $\mu(\mathcal F)>\mu(\mathcal E')$. We claim that $\mathcal F$ descends to $X_C$. By Proposition~\ref{prop:relativebanachcolmezvsheaf}, it suffices to see that the two pullbacks of $\mathcal F$ to $X_{C'\hat{\otimes}_C C'}$ agree. This is true as there are no nonzero maps from $\mathcal F$ to $\mathcal E'/\mathcal F$ after base change to $X_R$ for any perfectoid $C'$-algebra $R$: If there were such a nonzero map, there would also be a nonzero map for some choice of $R=C''$ a complete algebraically closed nonarchimedean field. But then $\mathcal F$ is still semistable and all pieces of the Harder--Narasimhan filtration of $\mathcal E'/\mathcal F$ are of smaller slope, so such maps do not exist.

For an extension $E'|E$, the similar arguments work, using Galois descent instead (noting that one may assume that $E'|E$ is Galois by passing to Galois hulls). Note that the pullback of $\mathcal O_{X_{C,E}}(1)$ is $\mathcal O_{X_{C,E'}}(r)$, causing the mismatch in slopes.
\end{proof}

The main theorem on the classification of vector bundles is the following. Our proof follows the arguments of Hartl-Pink, \cite{HartlPink}, to reduce to Lemma~\ref{lem:keylemma} below. However, we give a new and direct proof of this key lemma, which avoids any hard computations by using the geometry of diamonds and v-descent. We thus get a new proof of the classification theorem.\footnote{First, it has been proven for $E$ of equal characteristic in \cite{HartlPink} and for $p$-adic $E$ by Kedlaya in \cite{KedlayaLocalMonodromy}; both of these proofs used heavy computations to prove Lemma~\ref{lem:keylemma}. A more elegant proof was given by Fargues--Fontaine \cite{FarguesFontaine} (for all $E$) by reducing to the description of the Lubin--Tate and Drinfeld moduli spaces of $\pi$-divisible $\mathcal O_E$-modules, and their Grothendieck--Messing period morphisms (which arguably also involve some nontrivial computations). Finally, for $p$-adic $E$ a proof is implicit in Colmez' work \cite{ColmezEspaces} on Banach--Colmez spaces.}

\begin{theorem}\label{thm:classificationvectorbundles} Any vector bundle $\mathcal E$ on $X_C$ is isomorphic to a direct sum of vector bundles of the form $\mathcal O_{X_C}(\lambda)$ with $\lambda\in \mathbb Q$. If $\mathcal E$ is semistable of slope $\lambda$, then $\mathcal E\cong \mathcal O_{X_C}(\lambda)^m$ for some $m\geq 0$.
\end{theorem}

\begin{proof} We argue by induction on the rank $n$ of $\mathcal E$, so assume the theorem in rank $\leq n-1$ (and for all choices of $E$); the case $n=1$ has been handled already. By the vanishing of $H^1(X_C,\mathcal O_{X_C}(\lambda))=0$ for $\lambda>0$, the theorem follows for $\mathcal E$ if $\mathcal E$ is not semistable. Thus assume $\mathcal E$ is semistable of slope $\lambda=\frac sr$ with $s\in \mathbb Z$ and $r>0$ coprime. It suffices to find a nonzero map $\mathcal O_{X_C}(\lambda)\to \mathcal E$: Indeed, by stability of $\mathcal O_{X_C}(\lambda)$, the map is necessarily injective (the category of semi-stable vector bundles of slope $\lambda$ is abelian with simple objects the stable vector bundles of slope $\lambda$), and the quotient will then again be semistable of slope $\lambda$, and thus by induction isomorphic to $\mathcal O_{X_C}(\lambda)^{m-1}$. One finishes by observing that $\mathrm{Ext}^1_{X_C}(\mathcal O_{X_C}(\lambda),\mathcal O_{X_C}(\lambda))=0$ by Proposition~\ref{prop:standardbanachcolmez}~(ii).

Thus, it suffices to find a nonzero map $\mathcal O_{X_C}(\lambda)\to \mathcal E$. Let $E'|E$ be the unramified extension of degree $r$, and consider the covering $f: X_{C,E'} = X_{C,E}\otimes_E E'\to X_{C,E}=X_C$. Then $\mathcal O_{X_C}(\lambda) = f_\ast \mathcal O_{X_{C,E'}}(s)$, and so it suffices to find a nonzero map $\mathcal O_{X_{C,E'}}(s)\to f^\ast \mathcal E$. In other words, up to changing $E$, we can assume that $\lambda\in \mathbb Z$. Then by twisting, we can assume $\lambda=0$.

Next, we observe that we are free to replace $C$ by an extension. Indeed, consider the v-sheaf sending $S\in \Perf_{C}$ to the isomorphisms $\mathcal E\cong \mathcal O_{X_C}^n$. This is a v-quasitorsor under $\underline{\GL_n(E)}$ (using Proposition~\ref{prop:standardbanachcolmez}~(ii)). If there is some extension of $C$ where we can find a nonzero section of $\mathcal E$ (and thus also trivialize $\mathcal E$), then it is a v-torsor under $\underline{\GL_n(E)}$. By v-descent of $\underline{\GL_n(E)}$-torsors, cf.~\cite[Lemma 10.13]{ECoD}, it is then representable by a space pro-\'etale over $\Spa C$, and thus admits a section.

Let $d\geq 0$ be minimal such that there is an injection $\mathcal O_{X_C}(-d)\hookrightarrow \mathcal E$, possibly after base enlarging $C$; by Theorem~\ref{thm:O1ample} some such $d$ exists. We want to see that $d=0$, so assume $d>0$ by way of contradiction. By minimality of $d$, the quotient $\mathcal F=\mathcal E/\mathcal O_{X_C}(-d)$ is a vector bundle, and by induction the classification theorem holds true for $\mathcal F$.

If $d\geq 2$, then we can by induction find an injection $\mathcal O_{X_C}(-d+2)\hookrightarrow \mathcal F$; taking the pullback defines an extension
\[
\mathcal O_{X_C}(-d)\to \mathcal G\to \mathcal O_{X_C}(-d+2).
\]
so by twisting
\[
\mathcal O_{X_C}(-1)\to \mathcal G(d-1)\to \mathcal O_{X_C}(1).
\]
By the key lemma, Lemma~\ref{lem:keylemma} below, we would, possibly after enlarging $C$, get an injection $\mathcal O_{X_C}\hookrightarrow \mathcal G(d-1)$, and hence an injection $\mathcal O_{X_C}(-d+1)\hookrightarrow \mathcal G\hookrightarrow \mathcal E$, contradicting our choice of $d$.

Thus, we may assume that $d=1$. If $\mathcal F$ is not semistable, then it admits a subbundle $\mathcal F'\subset \mathcal F$ of degree $\geq 1$ and rank $\leq n-2$. Applying the classification theorem to the pullback
\[
0\to \mathcal O_{X_C}(-1)\to \mathcal E'\to \mathcal F'\to 0
\]
of $\mathcal F'$, which is of slope $\geq 0$, we then get that $\mathcal E'\subset \mathcal E$ has a global section.

It remains the case that $d=1$ and that $\mathcal F$ is semistable, thus necessarily isomorphic to $\mathcal O_{X_C}(\tfrac 1{n-1})$. This is the content of the next lemma.
\end{proof}

\begin{lemma}\label{lem:keylemma} Let
\[
0\to \mathcal O_{X_C}(-1)\to \mathcal E\to \mathcal O_{X_C}(\tfrac 1n)\to 0
\]
be an extension of vector bundles on $X_C$, for some $n\geq 1$. Then there is some extension $C'|C$ of complete algebraically closed nonarchimedean fields such that $H^0(X_{C'},\mathcal E|_{X_{C'}})\neq 0$.
\end{lemma}

\begin{proof} Assume the contrary. Passing to Banach--Colmez spaces, we find an injection of v-sheaves
\[
f: \BC(\mathcal O_{X_C}(\tfrac 1n))\hookrightarrow \BC(\mathcal O_{X_C}(-1)[1]).
\]
The image cannot be contained in the classical points, i.e. the $C$-points (as these form a totally disconnected subset while the source is connected and not reduced to a point), so the image contains some non-classical point. After base change to some $C'|C$, we thus find that the image contains some nonempty open subset of $\BC(\mathcal O_{X_C}(-1)[1])$, as follows from the presentation
\[
\BC(\mathcal O_{X_C}(-1)[1]) = (\mathbb A^1_{C^\sharp})^\diamond/\underline{E}
\]
and the similar behaviour of non-classical points of $\mathbb A^1_{C^\sharp}$, cf.~proof of Proposition~\ref{prop:classicalpointsmaximalspectrum}. Translating this nonempty open subset to the origin, we find that the image of $f$ contains an open neighborhood of $0$, and then by rescaling by the contracting action of $E^\times$, we find that the map $f$ must be surjective, and thus an isomorphism.

In particular, this would mean that $\BC(\mathcal O_{X_C}(-1)[1])$ is a perfectoid space. This is patently absurd if $E$ is $p$-adic, as then the given presentation shows that $(\mathbb A^1_{C^\sharp})^\diamond$ is pro-\'etale over a perfectoid space and thus itself a perfectoid space, but $\mathbb A^1_{C^\sharp}$ is clearly not a perfectoid space.\footnote{We believe that also when $E$ is of equal characteristic, $\BC(\mathcal O_{X_C}(-1)[1])$ is not a perfectoid space, but we were not able to settle this easily.}

In general, we can argue as follows. There is a nonzero map
\[
\BC(\mathcal O_{X_C}(\tfrac 1n))\to (\mathbb A^1_{C^\sharp})^\diamond
\]
as $H^0(\mathcal O_{X_C}(\tfrac 1n))$ maps nontrivially to its fibre at the chosen untilt $\Spa C^\sharp\hookrightarrow X_C$. If $f$ is an isomorphism, we would then get a nonzero map
\[
(\mathbb A^1_{C^\sharp})^\diamond/\underline{E}\cong \BC(\mathcal O_{X_C}(-1)[1])\xrightarrow{f^{-1}} \BC(\mathcal O_{X_C}(\tfrac 1n))\to (\mathbb A^1_{C^\sharp})^\diamond.
\]
On the other hand, one can classify all $E$-linear maps $(\mathbb A^1_{C^\sharp})^\diamond\to (\mathbb A^1_{C^\sharp})^\diamond$. The latter are the same as maps $\mathbb A^1_{C^\sharp}\to \mathbb A^1_{C^\sharp}$ if $E$ is $p$-adic (by \cite[Proposition 10.2.3]{Berkeley}), respectively maps $\mathbb A^1_{C,\perf}\to \mathbb A^1_{C,\perf}$ if $E$ is of equal characteristic. Thus, they are given by some convergent power series $g(X)$ that is additive, i.e.~$g(X+Y)=g(X)+g(Y)$, and satisfies $g(aX)=ag(X)$ for all $a\in E$. (If $E$ is of characteristic $p$, then $g$ may a priori involve fractional powers $X^{1/p^i}$.) The equation $g(\pi X)=\pi g(X)$ alone in fact shows that only the linear coefficient of $g$ may be nonzero, so $g(X)=cX$ for some $c\in C^\sharp$, and thus $g$ is either an isomorphism or zero. But our given map is nonzero with nontrivial kernel, giving a contradiction.
\end{proof}

\subsection{Families of vector bundles}

Using the ampleness of $\mathcal O(1)$, we can now prove the following result on relative Banach--Colmez spaces.

\begin{proposition}\label{prop:relativebanachcolmez} Let $S$ be a perfectoid space over $\mathbb F_q$. Let $\mathcal E$ be a vector bundle on $X_S$. Then the Banach--Colmez space
\[
\BC(\mathcal E): T\mapsto H^0(X_T,\mathcal E|_{X_T})
\]
is a locally spatial diamond, partially proper over $S$. Moreover, the projectivized Banach--Colmez space
\[
(\BC(\mathcal E)\setminus \{0\})/\underline{E^\times}
\]
is a locally spatial diamond, proper over $S$.
\end{proposition}

\begin{proof} Using Theorem~\ref{thm:O1ample}, choose a presentation $\mathcal O_{X_S}(-n')^{m'}\to \mathcal O_{X_S}(-n)^m\to \mathcal E^\vee$ with $n,n'>0$. Dualizing, we get an exact sequence
\[
0\to \mathcal E\to \mathcal O_{X_S}(n)^m\to \mathcal O_{X_S}(n')^{m'}.
\]
This implies that $\BC(\mathcal E)\subset \BC(\mathcal O_{X/S}(n))^m$ is a closed subspace, so the first part follows from Proposition~\ref{prop:standardbanachcolmez}~(iii). For the second part, we may assume that $S$ is qcqs. It is also enough to prove the similar result for $(\BC(\mathcal E)\setminus \{0\})/\pi^\Z$ as the $\underline{\mathcal O_E^\times}$-action is free (so one can apply the last part of \cite[Proposition 11.24]{ECoD}). This follows from the following general lemma about contracting group actions on locally spectral spaces, noting that checking the conditions formally reduces to the case of $\BC(\mathcal O_{X_S}(n)^m)$ and from there to $\mathbb A^1_{S^\sharp}$ by evaluating sections at some collection of untilts.
\end{proof}

\begin{lemma}\label{lem:contractinglocspec} Let $X$ be a taut locally spectral space such that for any $x\in X$, the set $X_x\subset X$ of generalizations of $x$ is a totally ordered chain under specialization. Let $\gamma: X\iso X$ be an automorphism of $X$ such that the subset $X_0\subset X$ of fixed points is a spectral space. Moreover, assume that
\begin{altenumerate}
\item[{\rm (i)}] for all $x\in X$, the sequence $\gamma^n(x)$ for $n\to \infty$ converges towards $X_0$, i.e.~for all open neighborhoods $U$ of $X_0$, one has $\gamma^n(x)\in U$ for all sufficiently positive $n$;
\item[{\rm (ii)}] for all $x\in X\setminus X_0$, the sequence $\gamma^n(x)$ for $n\to -\infty$ diverges, i.e.~for all quasicompact open subspaces $U\subset X$, one has $\gamma^n(x)\not\in U$ for all sufficiently negative $n$.
\end{altenumerate}

Then $X_0\subset X$ is a closed subspace, the action of $\gamma$ on $X\setminus X_0$ is free and totally discontinuous (i.e.~the action map $(X\setminus X_0)\times \mathbb Z\to (X\setminus X_0) \times (X\setminus X_0)$ is a closed immersion), and the quotient $(X\setminus X_0)/\gamma^{\mathbb Z}$ is a spectral space.
\end{lemma}

\begin{remark} For applications of this lemma, we recall the following facts:
\begin{altenumerate}
\item[{\rm (i)}] If $X$ is any locally spatial diamond, then $|X|$ is a locally spectral space such that all for all $x\in |X|$, the set of generalizations of $x$ in $|X|$ is a totally ordered chain under specialization. Indeed, this follows from \cite[Proposition 11.19]{ECoD} and the similar property for analytic adic spaces.
\item[{\rm (ii)}] If $X$ is in addition partially proper over a spatial diamond, then $|X|$ is taut by \cite[Proposition 18.10]{ECoD}.
\end{altenumerate}
This means that the first sentence of the lemma is practically always satisfied.
\end{remark}

\begin{proof} Let $U\subset X$ be some quasicompact open neighborhood of $X_0$. First, we claim that one can arrange that $\gamma(U)\subset U$. Indeed, one has
\[
U\subset \gamma^{-1}(U)\cup \gamma^{-2}(U)\cup \ldots \cup \gamma^{-n}(U)\cup \ldots,
\]
as for any $x\in U\subset X$, also $\gamma^n(x)\in U$ for all sufficiently large $n$ by assumption, and so $x\in \gamma^{-n}(U)$ for some $n>0$. By quasicompacity of $U$, this implies that $U\subset \gamma^{-1}(U)\cup \ldots \cup \gamma^{-n}(U)$ for some $n$, and then $U^\prime = U\cup \gamma^{-1}(U)\cup \ldots\cup \gamma^{-n+1}(U)$ is a quasicompact open neighborhood of $X_0$ with $\gamma(U^\prime)\subset U^\prime$.

Now fix a quasicompact open neighborhood $U$ of $X_0$ with $\gamma(U)\subset U$. We claim that
\[
X_0=\bigcap_{n\geq 0} \gamma^n(U).
\]
Indeed, if $x\in X\setminus X_0$, then by assumption there is some positive $n$ such that $\gamma^{-n}(x)\not\in U$, giving the result.

In particular, for any other quasicompact open neighborhood $V$ of $X_0$, there is some $n$ such that $\gamma^n(U)\subset V$. Indeed, the sequence of spaces $\gamma^n(U)\setminus V$ is a decreasing sequence of spectral spaces with empty inverse limit, and so one of the terms is empty.

Consider the closure $\overline{U}\subset X$ of $U$ in $X$. As $X$ is taut, this is still quasicompact. Repeating the above argument, we see that for some $n>0$, one has
\[
\overline{U}\subset \gamma^{-1}(U)\cup \ldots\cup \gamma^{-n}(U) = \gamma^{-n}(U).
\]
This implies that the sequences $\{\gamma^n(U)\}_{n\geq 0}$ and $\{\gamma^n(\overline{U})\}_{n\geq 0}$ are cofinal. In particular,
\[
X_0 = \bigcap_{n\geq 0} \gamma^n(U) = \bigcap_{n\geq 0} \gamma^n(\overline{U})
\]
is a closed subset of $X$.

Next, we check that any point $x\in X\setminus X_0$ has an open neighborhood $V$ such that $\{\gamma^n(V)\}_{n\in \mathbb Z}$ are pairwise disjoint; for this it suffices to arrange that $V\cap \gamma^i(V)=\emptyset$ for all $i>0$. For this, note that if $n$ is chosen such that $\gamma^n(\overline{U})\subset U$, then up to rescaling by a power of $\gamma$, we can assume that $x\in U\setminus \gamma^{n+1}(\overline{U})$. Let $V\subset U\setminus \gamma^{n+1}(\overline{U})$ be a quasicompact open neighborhood of $x$. Then $\gamma^i(V)\cap V=\emptyset$ as soon as $i\geq n+1$. For the finitely many $i=1,\ldots,n$, we can use a quasicompacity argument, and reduce to proving that if $X_x$ is the localization of $X$ at $x$ (i.e., the set of all generalizations of $x$), then $X_x\cap \gamma^i(X_x)=\emptyset$ for $i=1,\ldots,n$. By our assumption on $X$, the space $X_x$ has a unique generic point $\eta\in X_x$ ($X_x$ is pro-constructible in a spectral space thus spectral and by our hypothesis $X_x$ is irreducible), which must then also be the unique generic point of $\gamma^i(X_x)$ if $X_x\cap\gamma^i(X_x)\neq \emptyset$. Thus, if $X_x\cap\gamma^i(X_x)\neq \emptyset$, then $\gamma^i(\eta)=\eta$, so $\eta\in X_0$. But $X_0$ is closed, so that $x\in X_0$, which is a contradiction.

In particular, the action of $\gamma$ on $X\setminus X_0$ is free and totally discontinuous, and the quotient $\overline{X}=(X\setminus X_0)/\gamma^{\mathbb Z}$ is a locally spectral space which is locally isomorphic to $X\setminus X_0$. A basis of open neighborhoods of $\overline{X}$ is given by the image of quasicompact open subsets $V\subset X\setminus X_0$ for which $\{\gamma^n(V)\}_{n\in \mathbb Z}$ are pairwise disjoint; it follows that these are quasicompact open subsets of $\overline{X}$. Also, the intersection of two such subsets is of the same form, so the quotient $\overline{X}$ is quasiseparated. Finally, note that $U\setminus \gamma(U)\to \overline{X}$ is a bijective continuous map, and the source is a spectral space (as $\gamma(U)\subset U$ is a quasicompact open subspace of the spectral space $U$), and in particular quasicompact, and so $\overline{X}$ is quasicompact.
\end{proof}

The result on properness of the projectivized Banach--Colmez space enables us to give quick proofs of the main results of \cite{KedlayaLiu1} (including an extension to the case of general $E$, in particular of equal characteristic).

\begin{theorem}[{\cite[Theorem 7.4.5, Theorem 7.4.9, Theorem 7.3.7, Proposition 7.3.6]{KedlayaLiu1}}]\label{thm:kedlayaliu} Let $S$ be a perfectoid space over $\Fq$ and let $\mathcal E$ be a vector bundle over $X_S$ of constant rank $n$.
\begin{altenumerate}
\item[{\rm (i)}] The function taking a geometric point $\Spa C\to S$ of $S$ to the Harder--Narasimhan polygon of $\mathcal E|_{X_C}$ is upper semicontinuous.
\item[{\rm (ii)}] Assume that the Harder--Narasimhan polygon of $\mathcal E$ is constant. Then there exists a global (separated exhaustive decreasing) Harder--Narasimhan filtration
\[
\mathcal E^{\geq \lambda}\subset \mathcal E
\]
specializing to the Harder--Narasimhan filtration at each point. Moreover, after replacing $S$ by a pro-\'etale cover, the Harder--Narasimhan filtration can be split, and there are isomorphisms
\[
\mathcal E^\lambda\cong \mathcal O_{X_S}(\lambda)^{n_\lambda}
\]
for some integers $n_\lambda\geq 0$.
\end{altenumerate}
\end{theorem}

\begin{proof} Note that the Harder--Narasimhan polygon can be described as the convex hull of the points $(i,d_i)$ for $i=0,\ldots,n$, where $d_i$ is the maximal integer such that $H^0(X_C,(\wedge^i \mathcal E)(-d_i)|_{X_C})\neq 0$. To prove part (i), it therefore suffices to show that for any vector bundle $\mathcal F$ on $X_S$, the locus of all geometric points $\Spa C\to S$ for which $H^0(X_C,\mathcal F|_{X_C})\neq 0$ is closed in $S$. But note that this is precisely the image of
\[
(\BC(\mathcal F)\setminus \{0\})/\underline{E^\times}\to S.
\]
As this map is proper by Proposition~\ref{prop:relativebanachcolmez}, its image is closed. To see that the endpoint of the Harder--Narasimhan polygon is locally constant, apply the preceding also to the dual of the determinant of $\mathcal E$.

For part (ii), it is enough to prove that v-locally on $S$, there exists an isomorphism $\mathcal E\cong \bigoplus_\lambda \mathcal O_{X_S}(\lambda)^{n_\lambda}$. Indeed, the desired global Harder--Narasimhan filtration will then exist v-locally, and it necessarily descends. The trivialization of each $\mathcal E^\lambda$ amounts to a torsor under some locally profinite group, and can thus be done after a pro-\'etale cover by \cite[Lemma 10.13]{ECoD}. Then the ability to split the filtration follows from Proposition~\ref{prop:standardbanachcolmez}~(iii).

We argue by induction on the rank of $\mathcal E$. Let $\lambda$ be the maximal slope of $\mathcal E$. We claim that v-locally on $S$, there is a map $\mathcal O_{X_S}(\lambda)\to \mathcal E$ that is nonzero in each fibre. Indeed, finding such a map is equivalent to finding a fibrewise nonzero map $\mathcal O_{X_S}\to \mathcal F=\sHom(\mathcal O_{X_S}(\lambda),\mathcal E)$. But then
\[
\BC(\mathcal F)\setminus \{0\}\to (\BC(\mathcal F)\setminus \{0\})/\underline{E^\times}\to S
\]
is a v-cover over which such a map exists: The first map is an $\underline{E^\times}$-torsor and thus a v-cover, while the second map is proper and surjective on geometric points, thus surjective by \cite[Lemma 12.11]{ECoD}. The dual map $\mathcal E^\vee\to \mathcal O_{X_S}(-\lambda)$ is surjective as can be checked over geometric points (using that $\mathcal O_{X_C}(-\lambda)$ is stable), thus the cokernel of $\mathcal O_{X_S}(\lambda)\to \mathcal E$ is a vector bundle $\mathcal E'$, that again has constant Harder--Narasimhan polygon. By induction, one can find an isomorphism
\[
\mathcal E'\cong \bigoplus_{\lambda'\leq \lambda} \mathcal O_{X_S}(\lambda')^{n_{\lambda'}'}.
\]
By Proposition~\ref{prop:standardbanachcolmez}~(i)--(ii), the extension
\[
0\to \mathcal O_{X_C}(\lambda)\to \mathcal E\to \bigoplus_{\lambda'\leq \lambda} \mathcal O_{X_S}(\lambda')^{n_{\lambda'}'}\to 0
\]
can be split after a further pro-\'etale cover, finishing the proof.
\end{proof}

Let us explicitly note the following corollary.

\begin{corollary}[{\cite[Theorem 8.5.12]{KedlayaLiu1}}]\label{cor:Elocalsystems} Let $S$ be a perfectoid space. The category of pro-\'etale $\underline{E}$-local systems $\mathbb L$ is equivalent to the category of vector bundles on $X_S$ whose Harder--Narasimhan polygon is constant $0$, via $\mathbb L\mapsto \mathbb L\otimes_{\underline{E}} \mathcal O_{X_S}$.
\end{corollary}

\begin{proof} First, the functor is fully faithful, as we can see by pro-\'etale descent (to assume $\mathbb L$ is trivial) and Proposition~\ref{prop:standardbanachcolmez}. Now essential surjectivity follows from Theorem~\ref{thm:kedlayaliu}.
\end{proof}

\section{Further results on Banach--Colmez spaces}

We include some further results on Banach--Colmez spaces.

\subsection{Cohomology of families of vector bundles}

First, we generalize the vanishing results of Proposition~\ref{prop:standardbanachcolmez} to families of vector bundles. A key tool is given by the following result, which is a small strengthening of \cite[Lemma 8.8.13]{KedlayaLiu1}.

\begin{proposition}\label{prop:localresolutionnonnegative} Let $S$ be a perfectoid space over $\Fq$, and let $\mathcal E$ be a vector bundle on $X_S$ such that all Harder--Narasimhan slopes of $\mathcal E$ at all geometric points are $\geq \tfrac 1r$, for some $r\geq 1$. Then locally (in the analytic topology) on $S$, there is an exact sequence
\[
0\to \mathcal O_{X_S}^m\to \mathcal F\to \mathcal E\to 0
\]
where $\mathcal F$ is semistable of degree $\tfrac 1r$ at all geometric points.
\end{proposition}

\begin{proof} We may assume that the rank and degree of $\mathcal E$ is constant, given by some $n\geq 0$ and $d\geq \tfrac nr$. We set $m=dr-n$. We can assume $S=\Spa(R,R^+)$ is affinoid perfectoid and pick $m$ untilts $S_i^\sharp=\Spa(R_i^\sharp,R_i^{\sharp +})$ over $E$, $i=1,\ldots,m$, such that $S_1^\sharp,\ldots,S_m^\sharp\subset X_S$ are pairwise disjoint; more precisely, choose $m$ maps $S\to \BC(\mathcal O(1))\setminus \{0\}$. (The disjointness can be ensured by defining these maps through suitable fractional powers of a pseudouniformizer so that each $S_i^\sharp$ has a fixed image under the radius map $Y_S\to (0,\infty)$.) Let $W_i$ be the fibre of $\mathcal E$ over $R_i^\sharp$, which is a finite projective $R^\sharp_i$-module. For any rank $1$ quotients $W_i\to R_i^\sharp$, we can pull back the sequence
\[
0\to \mathcal O_{X_S}^m\to \mathcal O_{X_S}(1)^m\to \bigoplus_{i=1}^m \mathcal O_{S_i^\sharp}\to 0,
\]
obtained from Proposition~\ref{prop:lubintateuntilt}, along
\[
\mathcal E\to \bigoplus_{i=1}^m \mathcal E\otimes_{\mathcal O_{X_S}} \mathcal O_{S_i^\sharp} = \bigoplus_{i=1}^m W_i\otimes_{R_i^\sharp} \mathcal O_{S_i^\sharp}\to \bigoplus_{i=1}^m \mathcal O_{S_i^\sharp}
\]
to get an extension
\[
0\to \mathcal O_{X_S}^m\to \mathcal E'\to \mathcal E\to 0.
\]

We claim that one can choose, locally on $S$, the rank $1$ quotients so that $\mathcal E'$ is semistable (necessarily of degree $\tfrac 1r$). For this, we argue by induction on $i=1,\ldots,m$ that one can choose (locally on $S$) the quotients of $W_1,\ldots,W_i$ so that the modification
\[
0\to \mathcal O_{X_S}^i\to \mathcal E_i\to \mathcal E\to 0,
\]
invoking only the quotients of $W_1,\ldots,W_i$, has the property that the Harder--Narasimhan slopes of $\mathcal E_i$ at all geometric points are $\geq \tfrac 1r$. We treat the case $i=1$; the general case uses the same argument, by looking at the exact sequence
\[
0\to \mathcal O_{X_S}\to \mathcal E_{i+1}\to \mathcal E_i\to 0
\]
presenting $\mathcal E_{i+1}$ as a modification of $\mathcal E_i$ (which, by induction, still has Harder--Narasimhan slopes $\geq \tfrac 1r$ everywhere) at $S_{i+1}^\sharp$, using the quotient $W_{i+1}\to R_{i+1}^\sharp$. (Note that by disjointness of the untilts, the fibres of $\mathcal E$ and $\mathcal E_i$ agree at $S_{i+1}^\sharp$.)

First, we handle the case that $S=\Spa(K,K^+)$ for a perfectoid field $K$ (not necessarily algebraically closed). Then $\mathcal E$ has a Harder--Narasimhan filtration, and look at the subbundle $\mathcal E^{\geq \lambda}\subset \mathcal E$ of maximal slope, where necessarily $\lambda> \tfrac 1r$. Also recall that any nonsplit extension
\[
0\to \mathcal O_{X_S}\to \mathcal G\to \mathcal E^{\geq \lambda}\to 0
\]
necessarily has Harder--Narasimhan slopes $\geq \tfrac 1r$. Indeed, if $\mathcal G\to \overline{\mathcal G}$ is a quotient of slope $\tfrac{d'}{r'}<\tfrac 1r$, then the quotient of $\overline{\mathcal G}$ by $\mathcal O_{X_S}$ has degree $\leq \tfrac{d'}{r'-1}$, but we still have $\tfrac{d'}{r'-1}\leq \tfrac 1r$ and $\mathcal E^{\geq \lambda}$ has no quotients of slope $\leq \frac 1r$. Thus, it is enough to ensure that the pullback of the extension to $\mathcal E^{\geq \lambda}$ is nonsplit. But if it splits, then the given map $\mathcal E^{\geq \lambda}\to \mathcal O_{S_1^\sharp}$ lifts to a map $\mathcal E^{\geq \lambda}\to \mathcal O_{X_S}$; by consideration of slopes, this map is necessarily trivial. Thus, if we let $W'_1\subset W_1$ be the fibre of $\E^{\geq \lambda}\subset \E$ at $S_1^\sharp$, it suffices to pick a quotient $W_1\to K_1^\sharp$ whose restriction to $W'_1$ is nonzero.

Going back to general affinoid $S$, pick any point $s\in S$. By the preceding argument, we can locally on $S$ find a quotient $W_1\to R_1^\sharp$ such that the corresponding extension $\mathcal E_1$ has the property that the Harder--Narasimhan slopes at $s$ are still $\geq \tfrac 1r$. By Theorem~\ref{thm:kedlayaliu}, the same is true in an open neighborhood, finishing the proof.
\end{proof}

In applications, it is often more useful to have the following variant, switching which of the two bundles is trivialized, at the expense of assuming slopes strictly larger than $\tfrac 1r$ (and allowing \'etale localizations in place of analytic localizations --- this is probably unnecessary).

\begin{proposition}\label{prop:localresolutionpositive} Let $S$ be a perfectoid space over $\Fq$ and let $\mathcal E$ be a vector bundle on $X_S$ such that at all geometric points of $S$, all Harder--Narasimhan slopes of $\mathcal E$ are $>\tfrac 1r$, for some $r\geq 1$. Then \'etale locally on $S$, there is a short exact sequence
\[
0\to \mathcal G\to \mathcal O_{X_S}(\tfrac 1r)^m\to \mathcal E\to 0
\]
where $\mathcal G$ is semistable of slope $0$ at all geometric points.
\end{proposition}

\begin{proof} We can assume that $\mathcal E$ has constant rank $n$ and degree $d$; we set $m=dr$. Let $\mathcal E'$ be the bundle of maps $\mathcal O_{X_S}(\tfrac 1r)^m\to \mathcal E$. Inside $\BC(\mathcal E')$, we can look at the locus $U\subset \BC(\mathcal E')$ of those maps $\mathcal O_{X_S}(\tfrac 1r)^m\to \mathcal E$ that are surjective and whose kernel is semistable of slope $0$. This is an open subdiamond of the Banach--Colmez space $\BC(\mathcal E')$: This is clear for the condition of surjectivity (say, as the cokernel of the universal map $\mathcal O_{X_T}(\tfrac 1r)^m\to \mathcal E|_{X_T}$ is supported on a closed subset of $X_T$, whose image is then closed in $T$), and then the locus where the kernel is semistable of slope $0$ is open by Theorem~\ref{thm:kedlayaliu}. By Proposition~\ref{prop:localresolutionnonnegative}, we see moreover that all geometric fibres of $U\to S$ are nonempty. It thus suffices to prove that for any geometric point $T=\Spa(C,C^+)\to S$, given as a cofiltered inverse limit of \'etale maps $S_i=\Spa(R_i,R_i^+)\to S$, and any section $s\in \BC(\mathcal E')(T)$, one can find a sequence of $i$'s and sections $s_i\in \BC(\mathcal E')(S_i)$ such that $s_i|_T\to s$ as $i\to \infty$. Indeed, one of the $s_i$ will then lie in $U(S_i)$, giving the desired short exact sequence.

To prove that one can approximate $s$, we argue in a way similar to the proof of Theorem~\ref{thm:O1ample}. To facilitate the estimates, it is useful to assume that all Harder--Narasimhan slopes of $\mathcal E'$ at $T$ are integral; this can always be achieved through pulling back $\mathcal E'$ to a cover $f: X_{S,E'}=X_S\otimes_E E'\to X_S$ for some unramified extension $E'|E$ (as this pullback multiplies slopes by $[E':E]$), noting that $\mathcal E'$ is a direct factor of $f_\ast f^\ast \mathcal E'=\mathcal E'\otimes_E E'$. We analyze $\mathcal E'$ in terms of its pullback to $Y_{S,[1,q]}$ and the isomorphism over $Y_{S,[q,q]}\cong Y_{S,[1,1]}$. Note that as $B_{C,[1,q]}$ is a principal ideal domain, the pullback of $\mathcal E'$ to $Y_{C,[1,q]}$ is necessarily free, and by approximation we can already find a basis over some $Y_{S_i,[1,q]}$; replacing $S$ by $S_i$ we can then assume that the pullback of $\mathcal E'$ to $Y_{S,[1,q]}$ is free. The descent datum is then given by $A^{-1}\phi$ for some matrix $A\in \GL_n(B_{R,[1,1]})$. After pullback to $T$, by Theorem~\ref{thm:classificationvectorbundles} and the assumption of integral slopes, one can in fact choose a basis so that $A$ is a diagonal matrix $D$ with positive powers of $\pi$ along the diagonal. Approximating this basis, we can assume that $A-D\in \pi^N B_{(R,R^+),[1,1]}^+$ for any chosen $N>0$.

Now the map
\[
\phi-D:B_{R,[1,q]}^r\to B_{R,[1,1]}^r
\]
is surjective by Proposition~\ref{prop:standardbanachcolmez}, and in fact there is some $M$ (depending only on $D$, not on $R$) such that for any $x\in (B_{(R,R^+),[1,1]}^+)^r$, there is some $y\in \pi^{-M} (B_{(R,R^+),[1,q]}^+)^r$ with $x=(\phi-D)(y)$. There are two ways to see the existence of $M$: Either by an explicit reading of the proof of Proposition~\ref{prop:standardbanachcolmez}, or as follows. Assume no such $M$ exists; then we can find perfectoid algebras $R_0,R_1,\ldots$ with integral elements $R_i^+\subset R_i$ and pseudouniformizers $\varpi_i\in R_i$, and sections $x_i\in (B_{(R_i,R_i^+),[1,1]}^+)^r$ such that there is no $y_i\in \pi^{-2i} (B_{(R_i,R_i^+),[1,q]}^+)^r$ with $x_i=(\phi-D)(y_i)$. Let $R^+$ be the product of all $R_i^+$, and $R=R^+[\tfrac 1{\varpi}]$ where $\varpi=(\varpi_i)_i\in R^+=\prod_i R_i^+$. Then all $x_i$ define elements of $(B_{(R,R^+),[1,1]}^+)^r$, and $x=x_0+\pi x_1+\pi^2 x_2+\ldots$ another element. By surjectivity of $\phi-D$, there is some $y\in (B_{(R,R^+),[1,q]})^r$ with $x=(\phi-D)(y)$. But then $y\in \pi^{-i} (B_{(R,R^+),[1,q]}^+)^r$ for some $i$, and then projecting along $(R,R^+)\to (R_i,R_i^+)$ contradicts the choice of $x_i$.

Taking $N>M$ above, one sees that also
\[
\phi-A=\phi-D+(D-A): B_{R,[1,q]}^r\to B_{R,[1,1]}^r
\]
is surjective, with the same bound (in particular independent of $R$). But now the section $s$ of $\mathcal E'$ over $X_T$ can be approximated by a section $s_i'$ of $B_{R_i,[1,q]}$, so that its image under $\phi-A$ will be small. By the preceding surjectivity, we can then replace $s_i'$ by $s_i=s_i'+\epsilon_i$ for some still small $\epsilon_i$ such that
\[
s_i\in B_{R_i,[1,q]}^{\phi=A} = H^0(X_{S_i},\mathcal E'|_{X_{S_i}}).
\]
This gives the desired conclusion.
\end{proof}

One can prove the following variants.

\begin{corollary}\label{cor:niceshortexseq} Let $S$ be a perfectoid space over $\Fq$ and let $\mathcal E$ be a vector bundle on $X_S$.
\begin{altenumerate}
\item[{\rm (i)}] Assume that all Harder--Narasimhan slopes of $\mathcal E$ are $\geq \tfrac 1r$. Then locally on $S$, for some $m\geq 0$ there is a short exact sequence
\[
0\lto \mathcal O_{X_S}^m\lto \mathcal F\lto \mathcal E\lto 0,
\]
where $\mathcal F$ is fibrewise semistable of slope $\tfrac 1r$.
\item[{\rm (ii)}] Assume that all Harder--Narasimhan slopes of $\mathcal E$ are $\geq \tfrac 1r$. Then locally on $S$, for some $m\geq 0$ there is a short exact sequence
\[
0\lto \mathcal O_{X_S}(\tfrac 1{2r})^m\lto \mathcal F\lto \mathcal E'\lto 0,
\]
where $\mathcal F$ is fibrewise semistable of slope $\tfrac 1r$, and $\mathcal E$ is a direct summand of $\mathcal E'$.
\item[{\rm (iii)}] Assume that all Harder--Narasimhan slopes of $\mathcal E$ are $>\tfrac 1r$. Then \'etale locally on $S$, for some $m\geq 0$ there is a short exact sequence
\[
0\lto \mathcal G\lto \mathcal O_{X_S}(\tfrac 1r)^m\lto \mathcal E\lto 0,
\]
where $\mathcal G$ is fibrewise semistable of slope $0$.
\item[{\rm (iv)}] Assume that all Harder--Narasimhan slopes of $\mathcal E$ are $>\tfrac 1r$. Then \'etale locally on $S$, for some $m\geq 0$ there is a short exact sequence
\[
0\lto \mathcal G\lto \mathcal O_{X_S}(\tfrac 1r)^m\lto \mathcal E'\lto 0,
\]
where $\mathcal G$ is fibrewise semistable of slope $\tfrac 1{2r}$, and $\mathcal E$ is a direct summand of $\mathcal E'$.
\end{altenumerate}
\end{corollary}

In parts (ii) and (iv), we believe that it is not really necessary to pass from $\mathcal E$ to $\mathcal E'$, but the argument below does not immediately prove this strengthening. For our applications, the passage from $\mathcal E$ to $\mathcal E'$ will not matter, so we content ourselves with the stated version.

\begin{proof} We can suppose $S$ is affinoid. We have already seen parts (i) and (iii). For part (ii), let $\pi_{2r}: X_{S,2r} = Y_S/\phi^{2r\mathbb Z}\to X_S$ be the finite \'etale cover $X_{S,E_{2r}}=X_{S,E}\otimes_E E_{2r} \to X_{S,E}=X_S$, where $E_{2r}$ is the unramified extension of degree $2r$ of $E$. We apply Proposition~\ref{prop:localresolutionnonnegative} to $(\pi_{2r}^\ast \mathcal E)(-1)$. We get, locally on $S$, a short exact sequence
\[
0\to \mathcal O_{X_{S,E_{2r}}}(1)^{m'}\to \mathcal F'\to \pi_{2r}^\ast \mathcal E\to 0
\]
where $\mathcal F'$ is fiberwise semistable of slope $2$. Thus, applying $\pi_{2r\ast}$, we get a short exact sequence
\[
0\to \mathcal O_{X_S}(\tfrac 1{2r})^{m'}\to \pi_{2r\ast} \mathcal F'\to \pi_{2r\ast} \pi_{2r}^\ast \mathcal E\to 0.
\]
Here $\pi_{2r\ast} \mathcal F'$ is fiberwise semistable of slope $\tfrac 1r$, and $\mathcal E$ is a direct summand of $\pi_{2r\ast} \pi_{2r}^\ast \mathcal E=\mathcal E\otimes_E E_{2r}$. One argues similarly for part (iv).
\end{proof}

\begin{proposition}\label{prop:h1loczero} Let $S\in \Perf_{\mathbb{F}_q}$, and let $\mathcal E$ be a vector bundle on $X_S$.
\begin{altenumerate}
\item If at all geometric points of $S$, all slopes of $\mathcal E$ are negative, then $H^0(X_S,\mathcal E)=0$.
\item If at all geometric points of $S$, all slopes of $\mathcal E$ are nonnegative, then there is a pro-\'etale cover $\tilde{S}\to S$ such that
\[
H^1(X_{\tilde{S}},\mathcal E|_{X_{\tilde{S}}})=0.
\]
\item If at all geometric points of $S$, all slopes of $\mathcal E$ are positive, then there is an \'etale cover $S'\to S$ such that for any affinoid perfectoid $T$ over $S'$, one has $H^1(X_T,\mathcal E|_{X_T})=0$.
\end{altenumerate}
\end{proposition}

\begin{proof} Part (i) can be checked on geometric points, where it follows from Theorem~\ref{thm:classificationvectorbundles} and Proposition~\ref{prop:standardbanachcolmez}~(i). For part (ii), we use Proposition~\ref{prop:localresolutionnonnegative} (applied to $\mathcal E(1)$) to produce locally on $S$ an exact sequence
\[
0\to \mathcal O_{X_S}(-1)^d\to \mathcal E'\to \mathcal E\to 0
\]
where $\mathcal E'$ is everywhere semistable of degree $0$. By Theorem~\ref{thm:kedlayaliu} we can find a pro-\'etale cover of $S$ over which $\mathcal E'\cong \mathcal O_{X_S}^r$. By the vanishing of $H^2$, this induces a surjection from $H^1(X_S,\OO_{X_S})^r$ onto $H^1(X,\mathcal E)$. Since $H^1(X_S,\O_{X_S})=H^1_{\proet}(S,\underline{E})$ by Proposition~\ref{prop:standardbanachcolmez}~(ii), this vanishes pro-\'etale locally on $S$ (e.g., when $S$ is strictly totally disconnected, by Lemma~\ref{lemma:trivial loc sys} below).

For part (iii), we use Corollary~\ref{cor:niceshortexseq}~(iii) to produce an \'etale cover of $S$ over which there is an exact sequence
\[
0\to \mathcal G\to \mathcal O_{X_S}(\tfrac 1r)^m\to \mathcal E\to 0.
\]
For any affinoid $T|S$, this induces a surjection from $H^1(X_T,\mathcal O_{X_T}(\tfrac 1r)^m)$ onto $H^1(X_T,\mathcal E|_{X_T})$, so we conclude by Proposition~\ref{prop:standardbanachcolmez}~(iii).
\end{proof}

\subsection{Families of Banach--Colmez spaces}

We can now prove the following strengthening of Proposition~\ref{prop:relativebanachcolmez}.

\begin{proposition}\label{prop:relativebanachcolmez1} Let $S$ be a perfectoid space over $\mathbb F_q$. Let $[\mathcal E_1\to \mathcal E_0]$ be a map of vector bundles on $X_S$ such that at all geometric points of $S$, the bundle $\mathcal E_1$ has only negative Harder--Narasimhan slopes.
\begin{altenumerate}
\item[{\rm (i)}] The Banach--Colmez space
\[
\BC([\mathcal E_1\to \mathcal E_0]): T\mapsto \mathbb H^0(X_T,[\mathcal E_1\to \mathcal E_0]|_{X_T})
\]
is a locally spatial diamond, partially proper over $S$.
\item[{\rm (ii)}] The projectivized Banach--Colmez space
\[
(\BC([\mathcal E_1\to \mathcal E_0])\setminus \{0\})/\underline{E^\times}
\]
is a locally spatial diamond, proper over $S$.
\item[{\rm (iii)}] Assume that all Harder--Narasimhan slopes of $\mathcal E_0$ at all geometric points are positive. Then
\[
\BC([\mathcal E_1\to \mathcal E_0])\to S
\]
is cohomologically smooth.
\end{altenumerate}
\end{proposition}

\begin{proof} All assertions are \'etale local (in fact v-local) on $S$. For parts (i) and (ii), let us first simplify the form of the complex $[\mathcal E_1\to \mathcal E_0]$. By Theorem~\ref{thm:O1ample}, we can find (for $S$ affinoid) some $d>0$ and a surjection
\[
\mathcal O_{X_S}(-d)^m\to \mathcal E_0.
\]
Let $\mathcal E_1'$ be the kernel of $\mathcal E_1\oplus \mathcal O_{X_S}(-d)^m\to \mathcal E_0$. Then we find a quasi-isomorphism
\[
[\mathcal E_1'\to \mathcal O_{X_S}(-d)^m]\to [\mathcal E_1\to \mathcal E_0].
\]
Note also that $\mathcal E_1'$ still satisfies the assumption on negative slopes. We get an exact sequence
\[
0\to \BC([\mathcal E_1'\to \mathcal O_{X_S}(-d)^m])\to \BC(\mathcal E_1'[1])\to \BC(\mathcal O_{X_S}(-d)^m[1]).
\]
As $\BC(\mathcal O_{X_S}(-d)^m[1])$ is separated by Proposition~\ref{prop:standardbanachcolmez}~(i), we see that parts (i) and (ii) reduce to the case of $\BC(\mathcal E_1'[1])$. Now applying Corollary~\ref{cor:niceshortexseq}~(iv) to the dual of $\mathcal E_1'$, we get (\'etale locally on $S$, and after replacing $\mathcal E_1'$ by the direct sum with another bundle) an exact sequence
\[
0\to \BC(\mathcal E_1'[1])\to \BC(\mathcal O_{X_S}(-\tfrac 1r)^m[1])\to \BC(\mathcal G[1])
\]
where $\mathcal G$ is semistable of slope $\tfrac 1{2r}$ everywhere. In particular, $\BC(\mathcal G[1])$ is separated over $S$ by pro-\'etale descent and Proposition~\ref{prop:standardbanachcolmez}~(i). Thus, $\BC(\mathcal E_1'[1])\subset \BC(\mathcal O_{X_S}(-\tfrac 1r)^m[1])$ is a closed subfunctor, finishing the proof of part (i) by applying Proposition~\ref{prop:standardbanachcolmez}~(i) again. Part (ii) is then reduced to the similar assertion for $\BC(\mathcal O_{X_S}(-\tfrac 1r)^m[1])$. Replacing $E$ by its unramified extension of degree $r$, this reduces to $\BC(\mathcal O_{X_S}(-1)^m[1])$. Now, as in the proof of Proposition~\ref{prop:relativebanachcolmez}, this follows from Lemma~\ref{lem:contractinglocspec}, where one checks the required contracting property of multiplication by $\pi$ by using the presentation
\[
\BC(\mathcal O_{X_S}(-1)[1]) = (\mathbb A^1_{S^\sharp})^\diamond/\underline{E}
\]
for an untilt $S^\sharp$ of $S$ over $E$.

It remains to prove part (iii). Note that one has a short exact sequence
\[
0\to \BC(\mathcal E_0)\to \BC([\mathcal E_1\to \mathcal E_0])\to \BC(\mathcal E_1[1])\to 0;
\]
by \cite[Proposition 23.13]{ECoD}, we can thus handle $\BC(\mathcal E_1[1])$ and $\BC(\mathcal E_0)$ individually. For the case of $\BC(\mathcal E_0)$, we use Corollary~\ref{cor:niceshortexseq}~(iv) to get, pro-\'etale locally on $S$ (and after replacing $\mathcal E_0$ by the direct sum with another bundle), an exact sequence
\[
0\to \mathcal O_{X_S}(\tfrac 1{2r})^{m'}\to \mathcal O_{X_S}(\tfrac 1r)^m\to \mathcal E_0\to 0,
\]
inducing a similar sequence of Banach--Colmez spaces. Then the result follows from \cite[Proposition 23.13]{ECoD}. For the case of $\BC(\mathcal E_1[1])$, choose a surjection $\mathcal O_{X_S}(-d)^m\to \mathcal E_1^\vee$ for some $d>0$; we get an exact sequence
\[
0\to \mathcal E_1\to \mathcal O_{X_S}(d)^m\to \mathcal F\to 0
\]
where necessarily all Harder--Narasimhan slopes of $\mathcal F$ are positive everywhere. This gives an exact sequence
\[
0\to \BC(\mathcal O_{X_S}(d)^m)\to \BC(\mathcal F)\to \BC(\mathcal E_1[1])\to 0,
\]
so the result follows from \cite[Proposition 23.13]{ECoD} and the case of positive slopes already established.
\end{proof}

\subsection{Punctured absolute Banach--Colmez spaces}

Finally, we analyze punctured absolute Banach--Colmez spaces. Recall that, in the situation of Proposition~\ref{prop:standardbanachcolmez}~(iv), one has
\[
\BC (\O(d))\cong \Spd (k\powerseries{x_1^{1/p^\infty},\dots,x_d^{1/p^\infty}}),
\]
so the v-sheaf $\BC(\O(d))$ fails to be a perfectoid space, or even a diamond, as it contains the non-analytic point $\Spd k$. However, passing to the punctured Banach--Colmez space
\[
\BC (\O(d))\setminus \{0\}\cong \Spa(k\powerseries{x_1^{1/p^\infty},\dots,x_d^{1/p^\infty}})_{\mathrm{an}}
\]
identifies with the analytic points, which form a perfectoid space; in fact, a qcqs perfectoid space. These objects first showed up in \cite{FarguesClassFieldTheory} in the case of positive slopes. It was remarked in \cite{FarguesClassFieldTheory} that the punctured version $\BC (\O(d))\setminus \{0\}$ is a diamond for all $d\geq 1$, that is moreover simply connected when $d>2$. This plays a key role in \cite{FarguesClassFieldTheory} since after base changing from $\Spd k$ to $\Spa C$ this is not simply connected anymore. In the above example,
$$
\BC ( \O(d))\setminus \{0\}= \Spa (k\powerseries{x_1^{1/p^\infty},\dots,x_d^{1/p^\infty}}) \setminus V(x_1,\cdots,x_d)
$$
is a qcqs perfectoid space that is simply connected when $d>1$. After base changing to $\Spa (C)$ this is a punctured $n$-dimensional open ball over $\Spa(C)$ that is not quasicompact anymore, and not simply connected. Thus, some new interesting phenomena appear when we consider absolute Banach--Colmez spaces.

Let us first continue the discussion with the case of $\O(d)$ for $d\geq 1$. In that case, there is a relation to Cartier divisors. Recall that any closed Cartier divisor $D\subset X_S$ is given by a line bundle $\mathcal I$ on $X_S$ together with an injection $\mathcal I\hookrightarrow \mathcal O_{X_S}$ with closed image. We will only consider the case of relative Cartier divisors, so that this map stays injective after base change to any geometric point. Now Theorem~\ref{thm:kedlayaliu} implies that after replacing $S$ by an open and closed cover, $\mathcal I$ is of degree $-d$ for some integer $d\geq 0$, and that there is an $\underline{E^\times}$-torsor of isomorphisms $\mathcal I\cong \mathcal O_{X_S}(-d)$. This shows that the v-sheaf $\mathrm{Div}$ sending any $S$ to the closed relative Cartier divisors is given by
\[
\mathrm{Div} = \bigsqcup_{d\geq 0} \mathrm{Div}^d,\ \mathrm{Div}^d\cong (\BC(\mathcal O(d))\setminus \{0\})/\underline{E^\times}.
\]
Note that we are implicitly using a different definition of $\mathrm{Div}^1$ here, but Corollary~\ref{cor:descriptiondiv1} shows that they agree.

In particular, the moduli space $\mathrm{Div}^d$ of degree $d$ Cartier divisors is given by the projectivized Banach--Colmez space for $\mathcal O(d)$. On the other hand, in terms of divisors we can see the following proposition. Recall that one can take sums of Cartier divisors (by tensoring their ideal sheaves).

\begin{proposition}\label{prop:divdsymmetricpower} For any $d\geq 1$, the sum map
\[
(\mathrm{Div}^1)^d\to \mathrm{Div}^d: (D_1,\ldots,D_d)\mapsto D_1+\ldots+D_d
\]
is a quasi-pro-\'etale cover, identifying
\[
\mathrm{Div}^d = (\mathrm{Div}^1)^d/\Sigma_d,
\]
where $\Sigma_d$ is the symmetric group. In particular, $\mathrm{Div}^d$ is a diamond.
\end{proposition}

\begin{proof} By Proposition~\ref{prop:relativebanachcolmez}~(ii), all occuring spaces are proper over $\ast$. In particular, the sum map is proper. To check surjectivity as v-sheaves, we can then check on geometric points, where it follows from Proposition~\ref{prop:classicalpointsschematiccurve} (in whose proof we checked that any element of $P_d$ is a product of elements of $P_1$). In fact, we even get bijectivity up to the $\Sigma_d$-action, and thus the isomorphism
\[
\mathrm{Div}^d = (\mathrm{Div}^1)^d/\Sigma_d
\]
as v-sheaves. But the projection $(\mathrm{Div}^1)^d\to (\mathrm{Div}^1)^d/\Sigma_d$ is quasi-pro-\'etale by \cite[Lemma 7.19, Definition 10.1 (i)]{ECoD}. As $\mathrm{Div}^1=\Spd E/\phi^\Z$ is a diamond, it follows that $\mathrm{Div}^d$ is a diamond by \cite[Proposition 11.4, Proposition 11.6]{ECoD}.
\end{proof}

Now we can analyze the case of general absolute Banach--Colmez spaces. Here, we abbreviate $\BC(D) = \BC(\mathcal E(D))$.

\begin{proposition}\label{prop:absolutebanachcolmez} Let $D$ be an isocrystal with only negative slopes (resp.~with only positive slopes), and work on $\Perf_k$.
\begin{altenumerate}
\item[{\rm (i)}] The punctured Banach--Colmez space $\BC(D)\setminus \{0\}$ (resp.~$\BC(D[1])\setminus \{0\}$) is a spatial diamond.
\item[{\rm (ii)}] The quotient
\[
\big (\BC(D)\setminus \{0\}\big )/\underline{E^\times}\lto \ast\ \big (\mathrm{resp.~} \big (\BC(D[1])\setminus \{0\}\big )/\underline{E^\times}\lto \ast\big )
\]
is proper, representable in spatial diamonds, and cohomologically smooth.
\end{altenumerate}
\end{proposition}

\begin{proof} Part (ii) follows from Proposition~\ref{prop:relativebanachcolmez1} and (for the cohomological smoothness after taking the quotient by $\underline{E^\times}$) \cite[Proposition 24.2]{ECoD}.

For part (i), we are going to apply Lemma \ref{lemma:critere spatial v faisceau}, so we first want to see that $\BC(D)\setminus \{0\}$ is a spatial v-sheaf. By the Dieudonn\'e--Manin classification, we can find a basis for $D$ so that $\phi$ is $E$-rational and $U:=\phi^N$ is a diagonal matrix with entries powers of $\pi$ for some $N>0$; this essentially means that $V$ is decent in the sense of \cite[Definition 1.8]{RapoportZink}. Then $\BC(D)$ (resp.~$\BC(D[1])$) is already defined on $\Perf_{\mathbb F_q}$. In particular, the absolute $q$-power Frobenius $\mathrm{Frob}$ acts on them (as it acts on $\Perf_{\mathbb F_q}$), but also $U$ is an endomorphism of $D$ commuting with $\phi$ and hence acts. Moreover, as $U=\phi^N$, the action of $U$ agrees with the action of $\Frob^N$. Moreover, the action of $U^{-1}$ (resp.~$U$) on $|\BC(D)\times_{\Fq} \Spa \Fq\laurentseries{t^{1/p^\infty}}|$ (resp.~$|\BC(D[1])\times_{\Fq} \Spa \Fq\laurentseries{t^{1/p^\infty}}|$) still satisfies the hypotheses of Lemma~\ref{lem:contractinglocspec}. This implies that
\[
(\BC(D)\setminus \{0\})/\phi^N\times_{\Fq} \Spa \Fq\laurentseries{t^{1/p^\infty}}
\]
is a spatial diamond, which can be translated into
\[
(\BC(D)\setminus \{0\})\times_{\Fq} \Spa \Fq\laurentseries{t^{1/p^\infty}}/\phi^N
\]
being a spatial diamond, as the absolute Frobenius acts trivially on the topological space. But $\Spa \Fq\laurentseries{t^{1/p^\infty}}/\phi^N\to \ast$ is qcqs, even proper, and cohomologically smooth. We can thus apply point (i) of Lemma \ref{lemma:critere spatial v faisceau} to conclude that $\BC(D)\setminus \{0\}$ (resp.~$\BC(D[1])\setminus \{0\}$, for which the same argument applies) is spatial.

It remains to see that it is a diamond. One easily reduces to the case that $D$ is simple, and allowing ourselves to replace $E$ by a finite unramified extension, to $D$ of rank $1$. The case of positive Banach--Colmez spaces now follows from Proposition~\ref{prop:divdsymmetricpower}, as it is an $\underline{E^\times}$-torsor over a diamond (so \cite[Proposition 11.7]{ECoD} applies). It remains to prove that this is a diamond in the case of a negative absolute Banach--Colmez spaces, i.e.~for $D=(E,\pi^n\phi)$ with $n>0$. Then $\mathcal{D}:=\BC(D[1])\setminus \{0\}$ classifies extensions
\[
0\lto \mathcal O_{X_S}(-n)\lto \mathcal E\lto \mathcal O_{X_S}\lto 0
\]
that are geometrically fiberwise non split on $S$ (remark that those extensions are rigid). We now apply point (ii) of Lemma \ref{lemma:critere spatial v faisceau} using the Harder--Narasimhan stratification of $\mathcal{D}$ defined by $\E$. We can pass to the subsheaf of $\mathcal{D}$ where $\mathcal E$ is, at each geometric point, isomorphic to a given rank $2$ bundle, necessarily of the form $\mathcal O_{X_S}(-n+i)\oplus \mathcal O_{X_S}(-i)$ for some $0<i\leq \frac n2$ or to $\mathcal O_{X_S}(-\frac n2)$. 

 On such a stratum $\mathcal{D}_\alpha\subset \mathcal{D}$ there is a global Harder--Narasimhan filtration by Theorem~\ref{thm:kedlayaliu}, and trivializing the graded piece of lowest slope defines a pro-\'etale morphism $\tilde{\mathcal{D}}_\alpha \to \mathcal{D}_\alpha$. For $S\to \tilde{\mathcal{D}}_\alpha$ there is a morphism from $\O_{X_S}(-n)$ to this quotient of $\E$ by composing with the inclusion $\O_{X_S}(-n)\hookrightarrow \E$. Since the extension is non-split geometrically fiberwise on $S$, this morphism is non-zero geometrically fiberwise. 
 This defines a morphism $\tilde{\mathcal{D}}_\alpha\to X$ from $\tilde{\mathcal{D}}_\alpha$ to a punctured positive absolute Banach--Colmez space $X$, which is a diamond by what we have already proved. Thus
\[
\tilde{\mathcal{D}}_\alpha\subset \tilde{\mathcal{D}}_\alpha\times X
\]
where the latter is a diamond as $\tilde{\mathcal D}_\alpha\to \ast$ is representable in diamonds, so \cite[Proposition 11.10]{ECoD} shows that $\tilde{\mathcal D}_\alpha$ is a diamond.
\end{proof}

\begin{lemma}\label{lemma:critere spatial v faisceau}
Let $\mathcal{F}$ be a small v-sheaf.
\begin{altenumerate}
\item[{\rm (i)}] Suppose there exists a surjective qcqs cohomologically smooth morphism $D\to \mathcal{F}$ where $D$ is a spatial diamond. Then $\mathcal{F}$ is a spatial v-sheaf. 
\item[{\rm (ii)}] Suppose moreover there is a family of locally closed generalizing subsets $(X_\alpha)_\alpha$, $X_\alpha\subset |\mathcal{F}|$, such that for each $\alpha$ the associated subsheaf of $\mathcal{F}$ is a diamond.  Then $\mathcal{F}$ is a spatial diamond.
\end{altenumerate}
\end{lemma}
\begin{proof}
For point (i), since $D$ is qcqs and $D\to \mathcal{F}$ qcqs surjective, $\mathcal{F}$ is qcqs. Since cohomologically smooth implies universally open we can apply \cite[Lemma 2.10]{ECoD} to conclude it is spatial. For point (ii) we apply \cite[Theorem 12.18]{ECoD}. Let $\mathcal{G}_\alpha\subset \mathcal{F}$ be associated to $X_\alpha$. From \cite[Lemma 7.6]{ECoD}  we deduce that $\mathcal{G}_\alpha\hookrightarrow \mathcal{F}$ is quasi-pro-\'etale. This implies the result.
\end{proof}

\begin{remark}
The proof of Proposition \ref{prop:absolutebanachcolmez} for negative absolute Banach--Colmez space goes the same way as the proof of the fact that $\mathrm{Gr}_{\leq \mu}$ is a spatial diamond, \cite[Theorem 19.2.4]{Berkeley}. One first proves this is a spatial v-sheaf and then one stratifies it by locally closed generalizing subsets that are diamonds.
\end{remark}

\begin{remark}
One has to be careful that although the absolute $\BC (\O(d))\setminus \{0\}$ is a spatial diamond, $(\BC (\O(d))\setminus \{0\})/\pi^\Z$ is not spatial anymore since not quasiseparated, \cite[Remarque 2.15]{FarguesClassFieldTheory}. In this context the good object is the morphism $(\BC (\O(d))\setminus \{0\})/\pi^\Z\to \ast$ that is representable in spatial diamonds.
\end{remark}

\begin{remark}
In the equal characteristic case, $E=\mathbb{F}_q\laurentseries{\pi}$, the structure of punctured positive absolute Banach--Colmez spaces is much simpler since they are perfectoid spaces. Nevertheless the structure of the punctured negative one is not, they are only spatial diamonds. We will see below that $\BC(\mathcal O(-1)[1])$ is stratified into the open part, which can be written as the quotient of $\Spd k\laurentseries{t}$ by the action of a profinite group, and a point $\Spd k$. However, the degeneration to the point happens at the boundary of the open unit disc as $|t|\to 1$, not as $|t|\to 0$ as in $\Spa k\powerseries{t}$. Thus $\BC(\mathcal O(-1)[1])$ is a rather strange geometric object.\footnote{This example was critical in convincing us to not try to develop a version of the theory of diamonds that would allow non-analytic test objects like $\Spa k\powerseries{t^{1/p^\infty}}$ and would thus make $\BC(\mathcal O(1))$ itself representable: After all, in the context of absolute Banach--Colmez spaces, $\BC(\mathcal O(1))=\Spd k\powerseries{t^{1/p^\infty}}$ and $\BC(\mathcal O(-1)[1])$ play very similar roles, so the formalism should also treat them similarly.}
\end{remark}

\begin{example}\label{ex:BCO-1} The absolute $\BC(\O(-1)[1])\setminus \{0\}$ classifies extensions
\[
0\lto \mathcal O_{X_S}(-1)\lto \mathcal E\lto \mathcal O_{X_S}\to 0
\]
that are non-split fiberwise on $S$. Any such extension is, at each geometric point, isomorphic to $\mathcal O_{X_S}(-\tfrac 12)$. Parametrizing isomorphisms $\mathcal E\cong \mathcal O_{X_S}(-\tfrac 12)$ defines a $\underline{D^\times}$-torsor, where $D$ is the quaternion algebra over $E$; here we use Theorem~\ref{thm:kedlayaliu}. Remark that if $0\to \O_{X_S}(-1)\to \O_{X_S}(-\tfrac{1}{2})\to \mathcal{L}\to 0$ is an extension, then taking the determinant automatically fixes an isomorphism $\mathcal{L}\cong \O_{X_S}$, and thus the $\underline{E^\times}$-torsor of isomorphisms between $\O_{X_S}$ and $\mathcal{L}$ is trivial. From this we deduce that 
$$
\BC (\O(-1)[1])\setminus \{0\} \simeq (\BC (\O(\tfrac{1}{2}))\setminus \{0\})/ \underline{\mathrm{SL}_1(D)}
$$
with 
$$
\BC (\O(\tfrac{1}{2})) \setminus \{0\} \simeq \Spa (k\laurentseries{x^{1/p^\infty}})
$$
the punctured universal cover of a $1$-dimensional formal $\pi$-divisible $\O_E$-module of height $2$.

Let us compare this with our previous description of $\BC (\O(-1)[1])$ after pullback to $\Spa (C)$, fixing an untilt $C^\sharp$ over $E$ and $t\in H^0(X_C,\O_{X_C}(1))\setminus \{0\}$: the exact sequence 
$$
0\to \O_{X_C}(-1)\to \O_{X_C}\to \mathcal O_{C^\sharp}\to 0
$$
induces an isomorphism 
$$
\BC (\O(-1)[1])\times_k \Spa (C) \cong (\mathbb A^1_{C^\sharp})^\diamond /\underline{E}.
$$
We thus have 
$$
\BC (\O(-1)[1])\setminus \{0\}\times_k \Spa (C) \cong (\Omega_{C^\sharp})^\diamond / \underline{E}
$$
where $\Omega = \mathbb A^1_E\setminus \underline{E}$ is Drinfeld's upper half plane over $E$.
\\

We deduce an isomorphism 
$$
\big ( ( \BC(\O(\tfrac{1}{2}))\setminus \{0\} ) \times_k \Spa (C) \big )/ \underline{\mathrm{SL}_1 (D)}\cong (\Omega_{C^\sharp})^\diamond /\underline{E}.
$$
This isomorphism is in fact deduced from the isomorphism between Lubin--Tate and Drinfeld towers. In fact, \cite{ScholzeWeinstein}, the Lubin--Tate tower in infinite level, $\LT_\infty$ over $\Spd(C^\sharp)$, is the moduli of modifications $\O_{X_S}^2\hookrightarrow \O_{X_S} (\tfrac{1}{2})$ at the point of the curve defined by the untilt $C^\sharp$. From this one deduces a $\underline{D^\times}$-equivariant isomorphism 
$$
\LT_\infty^\flat \Big / \left ( \begin{matrix}
1 & \underline{E} \\
0 & \underline{E^\times}
\end{matrix} \right ) \cong  (\BC ( \O (\tfrac{1}{2} ))\setminus \{0\} )\times_k \Spa (C).
$$
Dividing this isomorphism by $\underline{\mathrm{SL}_1(D)}$ one obtains the preceding isomorphism.
\end{example}

\begin{example} The absolute $\BC (\O(-2)[1])\setminus \{0\}$ classifies extensions $$0\lto \O_{X_S}(-1)\lto \E\lto \O_{X_S}(1)\lto 0$$ that are non-split fiberwise on $S$. There is only one Harder--Narasimhan stratum and geometrically fiberwise on $S$, $\E$ is a trivial vector bundle. The moduli of surjections $\O_{X_S}^2\twoheadrightarrow \O_{X_S}(1)$ is the open subset 
$$
U\subset (\BC(\O(1))\setminus \{0\} )\times_k (\BC(\O(1))\setminus \{0\} )
$$
equal to 
$$
U=  (\BC(\O(1))\setminus \{0\} )^2 \setminus (\underline{E}^\times\times 1).\Delta
$$
where $\Delta$ is the diagonal of $(\BC(\O(1))\setminus \{0\} )^2$, that is to say couples $(x,y)$ of sections of $H^0(X_S,\O_{X_S}(1))$ that are fiberwise$/S$ non-zero and linearly independent over $E$. Here $(\underline{E}^\times\times 1).\Delta$ is a locally profinite union of copies of $\Delta$.

Again, by consideration of determinants, $\ker ( \O_{X_S}^2\twoheadrightarrow \O_{X_S}(1))$ is canonically identified with $\O_{X_S}(-1)$. This implies that 
$$
\BC (\O(-2)[1])\setminus \{0\} = U/\underline{\SL_2(E)}.
$$
\end{example}

\chapter{$\mathrm{Bun}_G$}

Throughout this chapter, we fix a reductive group $G$ over the nonarchimedean local field $E$. As it will be important to study $\Bun_G$ over a geometric base point, we fix from now on an algebraical closure $k$ of $\Fq$ and work with perfectoid spaces $S$ over $\Spd k$; write $\Perf_k$ for the category.

\begin{definition} Let $\mathrm{Bun}_G$ be the prestack taking a perfectoid space $S\in \Perf_k$ to the groupoid of $G$-bundles on $X_S$.
\end{definition}

The main results of this chapter are summarized in the following theorem.

\begin{theorem}[Proposition~\ref{prop:BunGsmall}; Theorem~\ref{thm:classifyGbun}; Theorem~\ref{thm:newtonsemicont} and Theorem~\ref{thm:kappaconstant}; Theorem~\ref{cor:descrsemistable}; Proposition~\ref{prop:BunGb}] The prestack $\Bun_G$ satisfies the following properties.
\begin{altenumerate}
\item[{\rm (i)}] The prestack $\Bun_G$ is a small v-stack.
\item[{\rm (ii)}] The points $|\Bun_G|$ are naturally in bijection with Kottwitz' set $B(G)$ of isomorphism classes of $G$-isocrystals.
\item[{\rm (iii)}] The map
\[
\nu: |\Bun_G|\to B(G)\to (X_\ast(T)_{\mathbb Q}^+)^\Gamma
\]
is semicontinuous, and
\[
\kappa: |\Bun_G|\to B(G)\to \pi_1(G_{\overline{E}})_\Gamma
\]
is locally constant. Equivalently, the map $|\Bun_G|\to B(G)$ is continuous when $B(G)$ is equipped with the order topology.
\item[{\rm (iv)}] The semistable locus $\Bun_G^{\mathrm{ss}}\subset \Bun_G$ is open, and given by
\[
\Bun_G^{\mathrm{ss}}\cong \bigsqcup_{b\in B(G)_{\mathrm{basic}}} [\ast/\underline{G_b(E)}].
\]
\item[{\rm (v)}] For any $b\in B(G)$, the corresponding subfunctor
\[
i^b: \Bun_G^b=\Bun_G\times_{|\Bun_G|} \{b\}\subset \Bun_G
\]
is locally closed, and isomorphic to $[\ast/\tilde{G}_b]$, where $\tilde{G}_b$ is a v-sheaf of groups such that $\tilde{G}_b\to \ast$ is representable in locally spatial diamonds with $\pi_0 \tilde{G}_b=G_b(E)$. The connected component $\tilde{G}_b^\circ\subset \tilde{G}_b$ of the identity is cohomologically smooth of dimension $\langle 2\rho,\nu_b\rangle$.
\end{altenumerate}
\end{theorem}

The hardest part of this theorem is that $\kappa$ is locally constant. We give two proofs of this fact. If the derived group of $G$ is simply connected, one can reduce to tori, which are not hard to handle. In general, one approach is to use z-extensions $\tilde{G}\to G$ to reduce to the case of simply connected derived group. For this, one needs that $\Bun_{\tilde{G}}\to \Bun_G$ is a surjective map of v-stacks; we prove this using Beauville--Laszlo uniformization. Alternatively, at least for $p$-adic $E$, one can use the abelianized Kottwitz set of Borovoi \cite{BorovoiAbelian}, which we prove to behave well relatively over a perfectoid space $S$.

\section{Generalities}
\label{sec:introBunG}

There is a good notion of $G$-torsors in $p$-adic geometry:

\begin{defprop}[{\cite[Proposition 19.5.1]{Berkeley}\footnote{The reference applies in the case of $\mathbb Z_p$, but it extends verbatim to $\mathcal O_E$.}}] Let $X$ be a sousperfectoid space over $E$. The following categories are naturally equivalent.
\begin{altenumerate}
\item[{\rm (i)}] The category of adic spaces $T\to X$ with a $G$-action such that \'etale locally on $X$, there is a $G$-equivariant isomorphism $T\cong G\times X$.
\item[{\rm (ii)}] The category of \'etale sheaves $\mathcal Q$ on $X$ equipped with an action of $G$ such that \'etale locally, $\mathcal Q\cong G$.
\item[{\rm (iii)}] The category of exact $\otimes$-functors
\[
\mathrm{Rep}_E G\to \Bun(X)
\]
to the category of vector bundles on $X$.
\end{altenumerate}

A $G$-bundle on $X$ is an exact $\otimes$-functor
\[
\mathrm{Rep}_E G\to \Bun(X);
\]
by the preceding, it can equivalently be considered in a geometric or cohomological manner.
\end{defprop}

In particular, $G$-torsors up to isomorphism are classified by $H^1_\et(X,G)$. By Proposition~\ref{prop:relativebanachcolmezvsheaf}, the following defines a v-stack.

\begin{definition} Let $\mathrm{Bun}_G$ be the v-stack taking a perfectoid space $S\in \Perf_k$ to the groupoid of $G$-bundles on $X_S$.
\end{definition}

By the GAGA results from the previous chapter, we are free to replace $X_S$ by $X_S^{\mathrm{alg}}$ here, when $S$ is affinoid.

Our goal in this chapter is to analyze this v-stack. Before going on, let us quickly observe that it is small, i.e~there are perfectoid spaces $S$, $R$ with a v-surjection $S\to \Bun_G$ and a v-surjection $R\to S\times_{\Bun_G} S$.

\begin{proposition}\label{prop:BunGsmall} The v-stack $\Bun_G$ is small.
\end{proposition}

\begin{proof} It is enough to prove that if $S_i=\Spa(R_i,R_i^+)$, $i\in I$, is an $\omega_1$-cofiltered inverse system of affinoid perfectoid spaces with inverse limit $S=\Spa(R,R^+)$, then
\[
\Bun_G(S) = \varinjlim \Bun_G(S_i).
\]
Indeed, then any section of $\Bun_G$ over an affinoid perfectoid space $S=\Spa(R,R^+)$ factors over $S'=\Spa(R',R^{\prime +})$ for some topologically countably generated perfectoid algebra $R'$. But there is only a set worth of such $R'$ up to isomorphism, and then taking the disjoint union $T=\bigsqcup_{S',\alpha\in \Bun_G(S')} S'$ gives a perfectoid space that surjects onto $\Bun_G$. Moreover, the equivalence relation $T\times_{\Bun_G} T$ satisfies the same limit property, and hence also admits a similar surjection.

To see the claim, note first that $R=\varinjlim R_i$ as any Cauchy sequence already lies in some $R_i$. The same applies to $B_{R,I}$ for any interval $I$, and hence one sees that
\[
\Bun(X_S) = \varinjlim \Bun(X_{S_i}).
\]
Now the definition of $G$-torsors gives the claim.
\end{proof}

This proof uses virtually no knowledge about $\Bun_G$ and shows that any reasonable v-stack is small.

\section{The topological space $|\mathrm{Bun}_G|$}
\subsection{Points}

As a first step, we recall the classification of $G$-bundles on the Fargues--Fontaine curve over geometric points. This is based on the following definition of Kottwitz, \cite{KottwitzIsocrystals}.

\begin{definition} A $G$-isocrystal is an exact $\otimes$-functor
\[
\Rep_E G\to \Isoc_E.
\]
The set of isomorphism classes of $G$-isocrystals is denoted $B(G)$.
\end{definition}

By Steinberg's theorem, the underlying fibre functor to $\breve E$-vector spaces is isomorphic to the standard fibre functor; this shows one can identify $B(G)$ with the quotient of $G(\breve E)$ under $\sigma$-conjugation.

Composing with the exact $\otimes$-functor
\[
\Isoc_E\to \Bun(X_S): D\mapsto \mathcal E(D)
\]
any $G$-isocrystal defines a $G$-bundle on $X_S$, for any $S\in \Perf_k$. In particular, for any $b\in B(G)$, we denote by $\mathcal E_b$ the corresponding $G$-bundle on $X_S$.

\begin{theorem}[\cite{FarguesGBun}, \cite{AnschuetzGBun}]\label{thm:classifyGbun} For any complete algebraically closed nonarchimedean field $C$ over $k$, the construction above defines a bijection
\[
B(G)\to \mathrm{Bun}_G(C)/\! \cong\ : b\mapsto \mathcal E_b.
\]
\end{theorem}

\begin{proof} For the convenience of the reader, and as some of the constructions will resurface later, we give a sketch of the proof in \cite{AnschuetzGBun}. Any $G$-bundle on $X_C$ has its Harder--Narasimhan filtration, and the formation of the Harder--Narasimhan filtration is compatible with tensor products. This implies that any exact $\otimes$-functor $\Rep_E G\to \Bun(X_C)$ lifts canonically to an exact $\otimes$-functor $\Rep_E G\to \FilBun(X_C)$ to $\mathbb Q$-filtered vector bundles. To check exactness, note that if $E$ is $p$-adic, the category $\Rep_E G$ is semisimple and thus exactness reduces to additivity, which is clear. If $E$ is of equal characteristic, one needs to argue more carefully, and we refer to the proof of \cite[Theorem 3.11]{AnschuetzGBun}.

We can now project $\Rep_E G\to \FilBun(X_C)$ to the category $\GrBun(X_C)$ of $\mathbb Q$-graded vector bundles, and note that the essential image of this functor is landing in the category of bundles $\bigoplus_\lambda \mathcal E^\lambda$ such that each $\mathcal E^\lambda$ is semistable of slope $\lambda$. This category is in fact equivalent to $\Isoc_E$ by Theorem~\ref{thm:classificationvectorbundles} and Proposition~\ref{prop:standardbanachcolmez}~(ii). Thus, it suffices to see that the filtration on the exact $\otimes$-functor $\Rep_E G\to \FilBun(X_C)$ can be split. By GAGA, we are free to work with $X_C^{\mathrm{alg}}$ in place of $X_C$.

Looking at splittings (fpqc, but also Zariski) locally on $X_C^{\mathrm{alg}}$, they exist, and form a torsor under a unipotent group scheme $U$ over $X_C^{\mathrm{alg}}$, where $U$ is parametrizing automorphisms of the filtered fibre functor that are trivial on the graded pieces. One can then filter $U$ by vector bundles of positive slopes, and using the vanishing of their $H^1$, we get the desired splitting.
\end{proof}

In particular, using \cite[Proposition 12.7]{ECoD} it follows that the map $$B(G)\longrightarrow |\mathrm{Bun}_G|$$ is a bijection.

\subsection{Harder--Narasimhan stratification}

Now we need to recall Kottwitz's description of $B(G)$, following \cite{KottwitzIsocrystals2}. This relies on two invariants there, {\it the Newton point and the Kottwitz point}. Let $\overline{E}$ be a separable closure of $E$ and fix a maximal torus inside a Borel subgroup $T\subset B\subset G_{\overline{E}}$; the set of dominant cocharacters $X_\ast(T)^+$ is naturally independent of the choice of $T$ and $B$, and acquires an action of $\Gamma=\mathrm{Gal}(\overline{E}|E)$ via its identification with
$$
\Hom ( \mathbb{G}_{m \overline{E}},G_{\overline{E}})/G(\overline{E})\rm{-conjugacy}.
$$
 The Newton point is a map
\begin{align*}
\nu: B(G) &\longrightarrow (X_\ast(T)^+_{\mathbb Q})^\Gamma \\
 b &\longmapsto \nu_b.
\end{align*}
When $G=\mathrm{GL}_n$, then $X_\ast(T)\cong \mathbb Z^n$ and the target is the set of nonincreasing sequences of rational numbers, which are the slopes of the Newton polygon of the corresponding isocrystal. The Kottwitz point is a map
\begin{align*}
\kappa: B(G) &\longrightarrow \pi_1(G)_\Gamma \\
 b & \longmapsto \kappa(b),
\end{align*}
where $\pi_1(G):=\pi_1(G_{\overline{E}}) = X_\ast(T)/(\mathrm{coroot\ lattice})$ is the Borovoi fundamental group. For $G=\mathrm{GL}_n$, this is naturally isomorphic to $\mathbb Z$, and in this case $\kappa(b)$ is the endpoint of the Newton polygon. In general, this compatibility is expressed by saying that the images of $\kappa(b)$ and $\nu_b$ in
\[
\pi_1(G)_{\mathbb Q}^\Gamma
\]
agree (using an averaging operation for $\kappa(b)$). However, this means that in general $\kappa(b)$ is not determined by $\nu_b$, as $\pi_1(G)_\Gamma$ may contain torsion.

The definition of $\kappa$ is done in steps. First, one defines it for tori, where it is actually a bijection. Then one defines it for $G$ whose derived group is simply connected; in that case, it is simply done via passage to the torus $G/G_{\mathrm{der}}$ which does not change $\pi_1$. In general, one uses a $z$-extension $\tilde{G}\to G$ such that $\tilde{G}$ has simply connected derived group, observing that $B(\tilde{G})\to B(G)$ is surjective.

Borovoi, \cite{BorovoiAbelian}, gave a more canonical construction of $\kappa$ as an abelianization map that does not use the choice of a $z$-extension, at least in the case of $p$-adic $E$. We will recall the construction in Section~\ref{sec:Second proof local constancy}.

Finally, recall that $(\nu,\kappa): B(G)\to (X^*(T)_\mathbb{Q}^+)^\Gamma\times \pi_1 (G)_\Gamma$ is injective.

It is possible to define the Harder--Narasimhan polygon and first Chern class of a $G$-bundle on $X_C$ (at least for $G=\GL_n$). This matches the invariants $(\nu,\kappa)$ up to sign (which results from $D\mapsto \mathcal E(D)$ reversing slopes): Let $v\mapsto v^*= w_0(-v)$ be the involution of the positive Weyl chamber $X_*(T)_\Q^+$ where $w_0$ is the longest element of the Weyl group. Then the Harder--Narasimhan polygon of $\E_b$ is $\nu_b^*$, and its first Chern class is $-\kappa(b)$. For general $G$, we may take this as the definition of the Harder--Narasimhan polygon of $\mathcal E_b$, and its first Chern class.

We need to understand how $\nu$ and $\kappa$ vary on $B(G)$. The following result follows from Theorem~\ref{thm:kedlayaliu} and \cite[Lemma 2.2]{RapoportRichartz}.

\begin{theorem}[{\cite[Corollary 22.5.1]{Berkeley}}]\label{thm:newtonsemicont} The map
\[
\nu^\ast: |\mathrm{Bun}_G|\cong B(G)\to (X_\ast(T)^+_{\mathbb Q})^\Gamma
\]
is upper semicontinuous.
\end{theorem}

We will later prove in Theorem~\ref{thm:kappaconstant} that $\kappa$ is locally constant on $\Bun_G$.

\subsection{Geometrically fiberwise trivial $G$-bundles}

Let 
$$
[\ast/\underline{G(E)}]
$$
be the classifying stack of pro-\'etale $\underline{G(E)}$-torsors, and 
$$
\Bun_G^1\subset \Bun_G
$$
be the substack of geometrically fiberwise trivial $G$-bundles. One has $H^0(X_S,\O_{X_S})=\underline{E}(S)$ and thus $\underline{G(E)}$ acts on the trivial $G$-bundle. From this we deduce a morphism
$$
[\ast/\underline{G(E)}] \longrightarrow \Bun_G^1.
$$
We are going to prove that this is an isomorphism. Let us note that, although this is an isomorphism at the level of geometric points, we can not apply \cite[Lemma 12.5]{Berkeley} since it is not clear that it is qcqs.
\\

Theorems~\ref{thm:newtonsemicont} and~\ref{thm:kappaconstant} (to follow) taken together imply that the locus $$\mathrm{Bun}_G^1\subset \mathrm{Bun}_G$$ is an open substack. One of our proofs of Theorem~\ref{thm:kappaconstant} will however require this statement as an input. Of course, when  $\pi_1(G)_\Gamma$ is torsion free, that is to say $H^1 (E,G)=\{1\}$, Theorem~\ref{thm:newtonsemicont} is enough to obtain the openness.

\begin{theorem}\label{thm:trivialGtorsor} The substack $\mathrm{Bun}_G^1\subset \mathrm{Bun}_G$ is open, and the map
\[
[\ast/\underline{G(E)}]\iso \mathrm{Bun}_G^1
\]
defined above is an isomorphism.
\end{theorem}

\begin{proof} Let $S\in \Perf_k$ be qcqs with a map to $\mathrm{Bun}_G$. We need to see that the subset of $|S|$ over which this map is trivial at any geometric point is open; and that if this is all of $S$, then the data is equivalent to a pro-\'etale $\underline{G(E)}$-torsor.

Let us check the openness assertion. If $T\to S$ is surjective with $T$ qcqs then $|T|\to |S|$ is a quotient map. 
We can thus assume that $S$ is strictly totally disconnected. The locus where the Newton point is identically zero is an open subset of $S$ by Theorem~\ref{thm:newtonsemicont}, so passing to this open subset, we can assume that the Newton point is zero. In that case, for any algebraic representation $\rho: G\to \mathrm{GL}_n$, the corresponding rank $n$ vector bundle on $X_S$ is trivial (by functoriality of the Newton map). Now, geometrically fiberwise on $S$ trivial vector bundles on $X_S$ are equivalent to $\underline{E}$-local systems on $S$ by Corollary~\ref{cor:Elocalsystems}. On the other hand, as $S$ is strictly totally disconnected, all $\underline{E}$-local systems on $S$ are trivial (Lemma~\ref{lemma:trivial loc sys} for $H=\GL_n(E)$), and their category is equivalent to the category of finite free modules over $\mathcal{C}^0(|S|,E)$. Thus, the preceding discussion defines a fibre functor on $\Rep_{E} (G)$ with values in $\mathcal{C}^0(|S|,E)=\mathcal{C}^0(\pi_0 (S),E)$, i.e.~a $G$-torsor over $\Spec(\mathcal{C}^0(\pi_0 (S),E))$. Note that for all $s\in \pi_0 (S)$, the local ring $\varinjlim_{U\ni s} \mathcal{C}^0(|U|,E)$ is henselian as the local ring of the analytic adic space $\underline{\pi_0(S)}\times \Spa (E)$ at $s$. This implies that if this $G$-torsor is trivial at some point of $S$, then it is trivial in a neighborhood. This concludes the openness assertion.

Moreover, the preceding argument shows that the map $\ast\to \mathrm{Bun}_G^1$ is a pro-\'etale cover. As $\ast\times_{\mathrm{Bun}_G^1} \ast=\underline{G(E)}$, as automorphisms of the trivial $G$-torsor are given by $\underline{G(E)}$, we thus get the desired isomorphism.
\end{proof}

\begin{remark}[{\cite[Lemma 10.13]{ECoD}}]\label{remark:torseur representable I}
If $S$ is a perfectoid space, and $T\rightarrow S$ is a pro-\'etale $\underline{G(E)}$-torsor then $T$ is representable by a perfectoid space. In fact, $T=\varprojlim_K \underline{K} \backslash T$ where $K$ goes through the set of compact open subgroups of $G(E)$. By descent of \'etale separated morphisms (\cite[Proposition 9.7]{ECoD}), for each such $K$, $\underline{K}\backslash T$ is represented by a separated \'etale perfectoid space over $S$. The transition morphisms in the preceding limit are finite \'etale.
\end{remark}

\begin{lemma}\label{lemma:trivial loc sys}
Let $S$ be a strictly totally disconnected perfectoid space, and let $H$ be a first-countable locally profinite group. Then any pro-\'etale $\underline{H}$-torsor on $S$ is trivial.
\end{lemma}

\begin{proof}
Let $T\to S$ be such a torsor. Fix a compact open subgroup $K\subset H$; by first-countability, this is a countable limit of finite groups. Since $\underline{K}\backslash T\to S$ is an \'etale cover of perfectoid spaces it has a section and we can assume $T\to S$ is in fact a $\underline{K}$-torsor. Now, $T=\varprojlim_n  \underline{K_n}\backslash T$ where $K_n\subset K$ is a cofinal system of open subgroups (with $K_0=K$). Each map $\underline{K_{n+1}}\backslash T\to \underline{K_n}\backslash T$ is a finite \'etale cover and hence split. Inductively choosing splittings, we get the result.
\end{proof}

\subsection{Local constancy of the Kottwitz invariant}

A central result is the following.

\begin{theorem}\label{thm:kappaconstant} The map
\[
\kappa: |\mathrm{Bun}_G|\cong B(G)\to \pi_1(G)_\Gamma
\]
is locally constant.
\end{theorem}

Let us note the following corollary. We give a new proof (and slight strengthening) of a result of Rapoport--Richartz (when $p \mid |\pi_1(G)|$ the original proof used $p$-adic nearby cycles and relied on a finite type hypothesis).

\begin{corollary}[{\cite[Corollary 3.11]{RapoportRichartz}}]\label{cor:RapoportRichartz}
Let $S$ be an $\mathbb{F}_q$-scheme and $\E$ an $G$-isocrystal on $S$. The map $|S|\to \pi_1(G)_\Gamma$ that sends a geometric point $\bar{s}\to S$ to $\kappa (\E_{\bar{s}})$ is locally constant.
\end{corollary}

\begin{proof}
We can suppose $S=\Spec(R)$ is affine and defined over $k$. We get a small v-sheaf $\Spd(R,R)$, and $\E$ defines a morphism $\Spd(R,R)\to \Bun_G$. The induced map $\kappa: |\Spd(R,R)|\to |\Bun_G|\to \pi_1(G)_\Gamma$ is locally constant by Theorem~\ref{thm:kappaconstant}. As open and closed subsets of $\Spd(R,R)$ are in bijection with open and closed subschemes of $\Spec(R)$ (by \cite[Proposition 18.3.1]{Berkeley} applied to morphisms to $\ast\sqcup \ast$), we can thus assume that $\kappa: |\Spd(R,R)|\to |\Bun_G|\to \pi_1(G)_\Gamma$ is constant. But now for any geometric point $\bar{s}\to S$, the element $\kappa(\E_{\bar{s}})\in \pi_1(G)_\Gamma$ agrees with the image of $\ast=|\Spd(\bar{s},\bar{s})|\to |\Spd(R,R)|\to |\Bun_G|\to \pi_1(G)_\Gamma$, giving the desired result.
\end{proof}

\begin{remark} There is a natural map $|\Spd(R,R)|\to |\Spa(R,R)|$, the latter of which admits two natural maps to $|\Spec(R)|$ (given by the support of the valuation, or the prime ideal of all elements of norm $<1$). However, with either choice of map $|\Spa(R,R)|\to |\Spec(R)|$, the resulting map $|\Spd(R,R)|\to |\Spec(R)|$ does not commute with the $\kappa$ maps. Still, there is also the map $|\Spec(R)|\to |\Spd(R,R)|$, used in the proof, and this commutes with the $\kappa$ maps.
\end{remark}

We give two different proofs of Theorem~\ref{thm:kappaconstant}.
\\

\subsubsection{First proof}

For the first proof of Theorem~\ref{thm:kappaconstant}, we also need the following lemma that we will prove in the next section.

\begin{lemma}\label{lem:zextension} Let $\tilde{G}\to G$ be a central extension with kernel a torus. Then
\[
\mathrm{Bun}_{\tilde{G}}\to \mathrm{Bun}_G
\]
is a surjective map of v-stacks.
\end{lemma}

In fact, up to correctly interpreting all the relevant structure, if $Z\subset \tilde{G}$ is the kernel, then $\mathrm{Bun}_Z$ is a Picard stack (as for commutative $Z$ one can tensor $Z$-bundles) which acts on $\mathrm{Bun}_{\tilde{G}}$, and $\mathrm{Bun}_G$ is the quotient stack. It is in fact clear that it is a quasitorsor, and the lemma ensures surjectivity.

\begin{proof}[First Proof of Theorem~\ref{thm:kappaconstant}] Picking a $z$-extension, we can by Lemma~\ref{lem:zextension} reduce to the case that $G$ has simply connected derived group. Then we may replace $G$ by $G/G_{\mathrm{der}}$, and so reduce to the case that $G$ is a torus. By a further application of Lemma~\ref{lem:zextension}, we can reduce to the case that $G$ is an induced torus. In that case $\pi_1(G)_\Gamma$ is torsion-free, and so the Kottwitz map is determined by the Newton map, so the result follows from Theorem~\ref{thm:newtonsemicont}, noting that in the case of tori there are no nontrivial order relations so semicontinuity means local constancy.
\end{proof}

It remains to prove Lemma~\ref{lem:zextension}. This will be done in the next section, using Beauville--Laszlo uniformization.

\subsubsection{Second proof}
\label{sec:Second proof local constancy}

For this proof, we assume that $E$ is $p$-adic (otherwise certain non-\'etale finite flat group schemes may appear). We define
$$
B_{ab} (G) =  H^1 (W_E, [G_{sc} (\overline{\breve{E}})\rightarrow G(\overline{\breve{E}})]),
$$
the abelianized Kottwitz set (cohomology with coefficient in a crossed module, see \cite{BorovoiAbelian} and \cite[Appendix B]{LabesseBaseChange}). There is an abelianization map 
$$
B(G) \longrightarrow B_{ab} (G)
$$
deduced from the morphism $[1\rightarrow G]\to [G_{sc}\rightarrow G]$. If $T$ is a maximal torus in $G$ with reciprocal image $T_{sc}$ in $G_{sc}$ then
$$
[T_{sc}\to T]\longrightarrow [G_{sc}\to G]
$$
is a homotopy equivalence. If $Z$, resp.~$Z_{sc}$, is the center of $G$, resp.~$G_{sc}$, there is a homotopy equivalence
$$
[Z_{sc}\to Z]\longrightarrow [G_{sc}\to G].
$$

\begin{lemma}
There is an identification
$$
B_{ab} (G) = \pi_1(G)_{\Gamma}
$$
through which Kottwitz map $\kappa$ is identified with the abelianization map $B(G)\rightarrow B_{ab}(G)$.
\end{lemma}
\begin{proof}
Choose a maximal torus $T$ in $G$. One has 
\begin{align*}
B_{ab} (G) &= H^1 (W_E, [T_{sc} (\overline{\breve{E}})\to T(\overline{\breve{E}})]) \\
&= \rm{coker}\big ( \it {B(T_{sc}) \to B(T) }\big )
\end{align*}
since $H^2 ( W_{E}, T_{sc} ( \overline{\breve{E}}))=0$ (use \cite[Chapter II.3.3 example (c)]{SerreCohoGal} and \cite[Chapter III.2.3 Theorem 1']{SerreCohoGal}). The result is deduced using Kottwitz description of $B(T_{sc})$ and $B(T)=X_\ast(T)_\Gamma$.
\end{proof}

For $S\in \Perf_k$ there is a morphism of sites 
$$
\tau:(X_S)_{\et}\longrightarrow S_{\et}
$$
deduced from the identifications
$$
(X_S)_{\et} = (X_S^\diamond)_{\et} = (\Div^1_S)_{\et}
$$
and the projection $\Div^1_S=\Div^1\times S\to S$. Equivalently, $\tau^\ast$ takes any \'etale $T\to S$ to $X_T\to X_S$, which is again \'etale.
\\

We now interpret some \'etale cohomology groups of the curve as Galois cohomology groups, as in \cite{FarguesGBun} where this type of computation was done for the schematical curve attached to an algebraically closed perfectoid field. Below, we abbreviate $X_{C,C^+} := X_{\mathrm{Spa}(C,C^+)}$.

\begin{proposition}\label{prop:relative etale cohomology of curve}
Let $S\in \Perf_k$. 
\begin{altenumerate}
\item 
Let $\mathcal{F}$ be a locally constant sheaf of finite abelian groups on $\Spa (E)_{\et}$. One has
$$
R\tau_* \mathcal{F}_{|X_S} = \underline{R\Gamma_\et (\Spa (E),\mathcal{F})}
$$
as a constant complex on $S_{\et}$.
\item If $D$ is a diagonalizable algebraic group over $E$, the pro-\'etale sheaf associated to 
$$
T/S\longmapsto H^1_\et (X_T,D)
$$
is the constant sheaf with value $H^1(W_E,D(\overline{\breve{E}}))$.
\end{altenumerate}
\end{proposition}
\begin{proof}
Let us note $\mathcal{G}=\mathcal{F}|_{X_S}$.
There is a natural morphism $\underline{R\Gamma_\et (\Spa (E),\mathcal{F})}\to R\tau_* \mathcal{G}$. The morphism $\Div^1_S\to S$ is proper and applying \cite[Corollary 16.10 (ii)]{ECoD}, we are reduced to prove that 
$$
H^\bullet_{\et} (\Spa (E),\mathcal{F}) \iso H^\bullet_{\et} (X_{C,C^+}, \mathcal{G})
$$
when $C$ is an algebraically closed field. Since $X_{C,C^+}$ is quasicompact quasiseparated 
$$
H^\bullet_{\et} (X_{C,C^+}\hat\otimes_E \hat{\overline{E}}, \mathcal{G}) = \varinjlim_{E'|E \text{ finite}} 
H^\bullet (X_{C,C^+}\otimes_E E',\mathcal{G}),
$$
and, using Galois descent, it thus suffices to prove that the left member vanishes in positive degrees, and equals $\mathcal F|_{\overline{E}}$ in degree $0$.

Let $K=\Fq\laurentseries{T}$ and $X_{(C,C^+),K}$ the equal characteristic Fargues--Fontaine curve over $\Spa (K)$. Identifying $\hat{\overline{E}}^\flat$ with $\hat{K^{\mathrm{sep}}}$, one has $$(X_{(C,C^+),E}\hat\otimes_E \hat{\overline{E}})^\flat = X_{(C,C^+),K}\hat\otimes \hat{K^{\mathrm{sep}}}.$$ Using this we are reduced to prove that for any prime number $n$, for $i>0$ $$H^i_{\et} ( X_{(C,C^+),K}\hat\otimes_K \hat{K^{\mathrm{sep}}},\mathbb{Z}/n\mathbb{Z})=0.$$
This is reduced, as above, to prove that any class in $H^i_{\et} (X_{(C,C^+),K},\mathbb{Z}/n\mathbb{Z})$, $i>0$, is killed by pullback to a finite separable extension of $K$.

 When $n\neq p$ one has 
$$
R\Gamma_{\et} ( X_{(C,C^+),K},\mathbb{Z}/n\mathbb{Z}) = R\Gamma (\phi^\mathbb{Z}, R\Gamma_{\et} ( \mathbb{D}^*_{C,C^+}, \mathbb{Z}/n\mathbb{Z}))
$$
where $T$ is the coordinate on the open punctured disk $$  \mathbb{D}^*_{C,C^+} = \Spa (C,C^+)\times_{\Spa(\mathbb F_q)} \Spa (K).$$
One has $H^k (\mathbb{D}^*_{C,C^+}, \mathbb{Z}/n\mathbb{Z})=0$ for $k>1$, and this is equal, via Kummer theory, to $\mathbb{Z}/n\mathbb{Z} (-1)$ for $k=1$. The Kummer covering of $\mathbb{D}^*_{C,C^+}$ induced by $T\mapsto T^n$  kills any class in $H^1 (\mathbb{D}^*_{C,C^+}, \mathbb{Z}/n\mathbb{Z})$. Also $H^0(\mathbb D^*_{C,C^+},\mathbb Z/n\mathbb Z)= \mathbb Z/n\mathbb Z$ and the class in $H^1(\phi^{\mathbb Z},\mathbb Z/n\mathbb Z)=\mathbb Z/n\mathbb Z$ is killed by passing up along an unramified extension of $K$ of degree $n$.

When $n=p$ we use Artin-Schreier theory. Since $C$ is an algebraically closed field we have $H^i (X_{(C,C^+),K},\O)=0$ when $i>0$. Since the adic space $X_{(C,C^+),K}$ is noetherian we deduce that $H^i_\et (X_{(C,C^+),K},\O)=0$ for $i>0$. Thus, $H^i_{\et} (X_{(C,C^+),K},\mathbb{Z}/n\mathbb{Z})$ is $0$ for $i>1$ and $\coker (K\xrightarrow{F-Id} K)$ when $i=1$, which is killed by pullback to an Artin-Schreier extension of $K$. This finishes the proof of point (1).

For point (2). 
There is a natural morphism $$H^1 (W_E, D(\overline{\breve{E}}))\to H^1(X_S,D)$$ (see just after the proof of this proposition). 
Suppose first that $D$ is a torus. Then point (2) is the computation of the coarse moduli space of $\Bun_D$ as a pro-\'etale stack. This itself is a consequence of Theorem~\ref{thm:trivialGtorsor} using a translation argument from $1$ to any $[b]\in B(D)$ (use the Picard stack structure on $\Bun_D$).

For any $D$ we use the exact sequence 
$$
1\longrightarrow D^0\longrightarrow D\longrightarrow \underline{\pi}_0 (D)\longrightarrow 1,
$$
where $D^0$ is a torus as $E$ is of characteristic $0$. For $T/S$ there is a diagram 
$$
\begin{tikzcd}
H^0 (E,\underline{\pi}_0 (D)) \ar[r] \ar[d] & 
B(D^0) \ar[r]\ar[d] & H^1 (W_E,D (\overline{\breve{E}}))\ar[r]\ar[d] & H^1 (E,\underline{\pi}_0 (D)) \ar[d] \ar[r] & 0 \\
H^0(X_T,\underline{\pi}_0(D)) \ar[r] & 
H^1_{\et} ( X_T, D^0) \ar[r] &  H^1_{\et} (X_T,D) \ar[r] & H^1_{\et} (X_T,\underline{\pi}_0(D)) 
\end{tikzcd}
$$
since $H^2 (W_E, D^0 ( \overline{\breve{E}}))=0$ and $H^\bullet (W_E, \underline{\pi}_0(D) (\overline{\breve{E}}))=H^\bullet (W_E, \underline{\pi}_0(D)(\overline{E}))$. The result is then deduced from part (1) and the torus case.
\end{proof}

For $S\in \Perf_k$ there is a natural morphism of groups
$$
B_{ab} (G)\longrightarrow H^1_{\et} (X_S, [G_{sc}\to G]).
$$
This is deduced from the natural continuous morphism of sites 
$$
(X_S)_{\et} \longrightarrow \{ \text{discrete }W_E\text{-sets} \}.
$$

\begin{proposition}
\label{prop:relative abelianized coho}
For $S\in \Perf_k$, the pro-\'etale sheaf on $S$ associated with 
$$
T/S\longmapsto H^1_{\et} ( X_T, [G_{sc}\to G]) 
$$
is the constant sheaf with value $B_{ab} (G)$.
\end{proposition}
\begin{proof}
We use the homotopy equivalence $[Z_{sc}\to Z]\to [G_{sc}\to G]$. There is a diagram
\[\xymatrix{
H^1 (W_E, Z_{sc} (\overline{\breve{E}})) \ar[r]\ar[d]^\simeq &
H^1 (W_E, Z (\overline{\breve{E}})) \ar[r]\ar[d]^\simeq & H^1 (W_E,[G_{sc} ( \overline{\breve{E}})\to G(\overline{\breve{E}})]) \ar[r]\ar[d] & H^2 (W_E, Z_{sc}(\overline{\breve{E}})) \ar[d]^\simeq \\
H^1_{\et}(X_T, Z_{sc})\ar[r] &
 H^1_{\et}(X_T, Z)\ar[r] & H^1_{\et} ( X_T, [Z_{sc}\to Z]) \ar[r] & H^2_{\et} (X_T, Z_{sc}) \\
\ar[r] & H^2 (W_E,Z(\overline{\breve{E}}))\Big(\ar[r]^\simeq\ar@{^(->}[d] & H^2 (E,\underline{\pi}_0(Z))\Big) \ar[d]^\simeq \\
\ar[r] & H^2_\et ( X_T,Z)\Big(\ar[r] & H^2_\et (X_T,\underline{\pi}_0 (Z))\Big).
}\]
Using Proposition~\ref{prop:relative etale cohomology of curve} and some diagram chasing we conclude.
\end{proof}

\begin{proof}[Second Proof of Theorem~\ref{thm:kappaconstant}]
The theorem is now deduced from the preceding Proposition~\ref{prop:relative abelianized coho} and the abelianization map $H^1_{\et} (X_S,G)\to H^1_{\et} ( X_S, [G_{sc}\to G])$.
\end{proof}

\begin{remark}
Let $X_C^{\mathrm{alg}}$ be the schematical curve associated to $C|\mathbb{F}_q$ algebraically closed.
The results of \cite{FarguesGBun} for the \'etale cohomology of torsion local systems, \cite[Theorem 3.7]{FarguesGBun} and the vanishing of the $H^2(X_C,T)$ for a torus $T$, \cite[Theorem 2.7]{FarguesGBun}, can be stated in a more uniform way; if $D$ is a diagonalizable group over $E$ then $H^i(W_E, D(\overline{\breve{E}}))\iso H^i_{\et} (X_C,D)$ for $0\leq i\leq 2$. {\it Weil cohomology of 
$E$ is the natural cohomology theory that corresponds to \'etale cohomology of the curve.} For example Theorem~\ref{thm:classifyGbun} can be restated as $H^1(W_E, G(\overline{\breve{E}}))\iso H^1_{\et} (X_C,G)$ for a reductive group $G$.
\end{remark}

\subsection{The explicit description of $|\Bun_G|$}

Theorem~\ref{thm:newtonsemicont} and Theorem~\ref{thm:kappaconstant} imply that the map
\[
|\Bun_G|\to B(G)
\]
is continuous when the target is endowed with the topology induced by the order on $(X_\ast(T)^+_{\mathbb Q})^\Gamma$ and the discrete topology on $\pi_1(G)_\Gamma$.

\begin{conj} The map $|\Bun_G|\to B(G)$ is a homeomorphism.
\end{conj}

In other words, whenever $b,b'\in B(G)$ such that $b>b'$, there should be a specialization from $b$ to $b'$ in $|\Bun_G|$.

The conjecture is known for $G=\GL_n$ by work of Hansen, \cite{HansenGLn}, based on \cite{Hansenetal}. The argument has been extended to some other classical groups in unpublished work of Hamann. While finishing our manuscript, a proof for general $G$ has been given by Viehmann \cite{ViehmannOberwolfach}.

We will later prove some weak form of the conjecture in Corollary~\ref{cor:pi0bunG}, determining the connected components of $\Bun_G$ using simple geometric considerations.

\section{Beauville--Laszlo uniformization}

Recall from \cite[Lecture XIX]{Berkeley} the $B_{\mathrm{dR}}^+$-affine Grassmannian $$\mathrm{Gr}_G$$ of $G$ over $\Spd E$, sending an affinoid perfectoid $S=\Spa(R,R^+)$ over $\Spd E$ to the $G$-torsors over $\Spec (B_{\mathrm{dR}}^+(R^\sharp))$ with a trivialization over $B_{\mathrm{dR}}(R^\sharp)$; here $R^\sharp/E$ is the untilt of $R$ given by $S\to \Spd (E)$. In \cite{Berkeley} this was considered over $\Spa (C^\flat)$ for some $C|E$ algebraically closed but we want now to consider it in a more ``absolute'' way over $\Spd (E)$.

Since any $G$-torsor over $\Spec (B^+_{\mathrm{dR}} (R^\sharp))$ is trivial locally on $\Spa (R,R^+)_{\et}$, this coincides with the \'etale sheaf associated to the presheaf $(R,R^+)\mapsto G(B_{\mathrm{dR}} (R^\sharp))/G(B^+_{\mathrm{dR}} (R^\sharp))$. 

{\it This has an interpretation as a Beilinson--Drinfeld type affine Grassmannian.} If $\E, \E'\in \Bun_G (S)$ and $D\in \Div^1 (S)$,  a modification between $\E$ and $\E'$ at $D$ is an isomorphism
$$
\E|_{X_S\setminus D} \xrightarrow{\ \sim\ } \E'|_{X_S\setminus D}
$$
that is meromorphic along $D$. The latter means that for any representation in $\Rep_E (G)$, the associated isomorphism between vector bundles, $\mathcal{F}|_{X_S\setminus D} \xrightarrow{\ \sim\ } \mathcal{F}'|_{X_S\setminus D}$ extends to a morphism $\mathcal{F}\to \mathcal{F}'(kD)$ for $k\gg 0$ via $\mathcal{F}'\hookrightarrow \mathcal{F}' (kD)$. Beauville--Laszlo gluing then identifies $$\Gr_G /\phi^\Z\longrightarrow \Div^1$$ with the moduli of $D\in \Div^1 (S)$, $\E\in \Bun_G (X_S)$, and a modification between the trivial $G$-bundle and $\E$ at $D$,
cf.~\cite[Proposition 19.1.2]{Berkeley}. This defines a morphism of v-stacks
$$
\Gr_G\longrightarrow \Bun_G.
$$

\begin{proposition}\label{prop:beauvillelaszlosurj} The Beauville--Laszlo morphism
\[
\mathrm{Gr}_G\longrightarrow \mathrm{Bun}_G
\]
is a surjective map of v-stacks; in fact, of pro-\'etale stacks.
\end{proposition}

\begin{proof} Pick any $S=\Spa(R,R^+)\in \Perf_k$ affinoid perfectoid with a map to $\mathrm{Bun}_G$, given by some $G$-bundle $\mathcal E$ on $X_S$. Fix an untilt $S^\sharp$ of $S$ over $\Spa(E)$. 
To prove surjectivity as pro-\'etale stacks, we can assume that $S$ is strictly totally disconnected. By \cite[Th\'eor\`eme 7.1]{FarguesGBun} (in case $G$ quasisplit) and \cite[Theorem 6.5]{AnschuetzGBun} (for general $G$), for any connected component $\Spa(C,C^+)\subset S$ of $S$, the map $\mathrm{Gr}_G(C)\to \mathrm{Bun}_G(C)$ is surjective, so in particular, we can pick a modification $\mathcal E'_C$ of $\mathcal E|_{X_{C,C^+}}$ at $\Spa (C^\sharp,C^{\sharp,+})\hookrightarrow X_{C,C^+}$ such that $\mathcal E'_C$ is trivial.

Now, since $S$ is strictly totally disconnected, we can trivialize $\mathcal E$ at the completion at $S^\sharp$; as $\mathrm{Gr}_G(S)\to \mathrm{Gr}_G(C)$ is surjective (Lemma~\ref{lemma:surjectivite BdR}), we can lift $\mathcal E'_C$ to a modification $\mathcal E'$ of $\mathcal E$. Now Theorem~\ref{thm:trivialGtorsor} implies that $\mathcal E'$ is trivial in a neighborhood of the given point, and as $S$ is strictly totally disconnected, the corresponding $\underline{G(E)}$-torsor is trivial (Lemma~\ref{lemma:trivial loc sys}), so we can trivialize $\mathcal E'$ in a neighborhood of the given point. This shows that locally on $S$, the bundle $\mathcal E$ is in the image of $\mathrm{Gr}_G\to \mathrm{Bun}_G$, as desired.
\end{proof}

\begin{lemma}\label{lemma:surjectivite BdR}
For $S=\Spa (R,R^+)$ a strictly totally disconnected perfectoid space over $\Spa (E)$, and $s\in S$, the map $\mathrm{Gr}_G (R)\to \mathrm{Gr}_G (K(s))$ is  surjective.
\end{lemma}
\begin{proof}
Set $C=K(s)$. First note that $R\to C$ is surjective, as any connected component of $S$ is an intersection of open and closed subsets, and thus Zariski closed.

Since $C$ is algebraically closed, $\mathrm{Gr}_G (C)= G(B_{\mathrm{dR}} (C))/G(B^+_{\mathrm{dR}}(C))$. For any finite degree extension $E'|E$ in $C$, since $S\otimes_E E'\to S$ is \'etale and $S$ strictly totally disconnected, $E'$ can be embedded in $R$. Fix such an $E'\subset R$ that splits $G$ and a pair $T\subset B$ inside $G_{E'}$.  Let $\xi\in W_{\O_E}(R^{\flat,+})$ be a generator of the kernel of $\theta$. The Cartan decomposition $$G(B_{\mathrm{dR}}(C))=\coprod_{\mu\in X^*(T)^+} G(B^+_{\mathrm{dR}} (C)) \mu (\xi) G(B_{\mathrm{dR}}^+(C))$$ shows then that we only need to prove that $G(B^+_{\mathrm{dR}}(R))\to G(B^+_{\mathrm{dR}}(C))$ is surjective. 

Let us first remark that $G(R)\to G(C)$ is surjective. In fact, since $S$ is totally disconnected it suffices to check that $G(\O_{S,s})\to G(C)$ is surjective. But this is a consequence of the smoothness of $G$ and the fact that $\O_{S,s}$ is Henselian (with residue field $C$).

Moreover $G(B^+_{\mathrm{dR}}(R))=\varprojlim_{n\geq 1} G(B^+_{\mathrm{dR}}(R)/\Fil^n)$. Using the surjectivity of $\mathrm{Lie} (G)\otimes R\to \mathrm{Lie} (G)\otimes C$ the result is then deduced by an approximation argument.
\end{proof}

\begin{remark}
In the ``classical case'' of the moduli of $G$-bundles over a proper smooth algebraic curve over a field $k$,  $G/k$, Proposition~\ref{prop:beauvillelaszlosurj} is true only when $G$ is semi-simple (\cite{DrinfeldSimpson}). Typically this is false for $\GL_n$ in general. The main reason why it is true in our situation is that $\mathrm{Pic}^0 (X_{C,C^+})$ is trivial, equivalently that $X_C^{\mathrm{alg}}\setminus \{x\}$ is the spectrum of a principal ideal domain in Proposition~\ref{prop:classicalpointsschematiccurve}.
\end{remark}

\begin{lemma}\label{lemma:Affine Grassmanian Inductive limit}
If $G$ is split then 
$$
\mathrm{Gr}_G =\varinjlim_{\mu\in X_*(T)^+} \mathrm{Gr}_{G,\leq \mu}
$$
as a v-sheaf, where the index set is a partially ordered set according to the dominance order ($\mu\leq \mu'$ if $\mu'-\mu$ is a nonnegative integral sum of positive coroots).
\end{lemma}

\begin{proof}
Consider a morphism $S\to \mathrm{Gr}_G$ with $S$ quasicompact quasiseparated. Fix an embedding $G\hookrightarrow \GL_n$ such that the image of $T$ lies in the standard maximal torus of $\GL_n$. This induces an embedding $X_*(T)^+\hookrightarrow \mathbb{Z}^n$. One checks easily that the image of $|S|\to |\mathrm{Gr}_{\GL_n}|$ lies in a finite union of affine Schubert cells. Since the fibers of $X_*(T)^+ \to \mathbb{Z}^n/\mathfrak{S}_n$ are finite we deduce that there is a finite collection $(\mu_i)_{i}$, $\mu_i\in X_* (T)^+$, such that the image of $|S|\to |\mathrm{Gr}_G|$ lies in $\cup_i |\mathrm{Gr}_{G,\leq \mu_i}|$. By \cite[Proposition 19.2.3]{Berkeley}, $S_i:=S\times_{\mathrm{Gr}_G} \mathrm{Gr}_{G,\leq \mu_i}$ is closed in $S$ and thus quasicompact. Since the morphism $\coprod_i S_i\to S$ is surjective at the level of points with quasicompact source it is quasicompact and thus a v-cover. This allows us to conclude.
\end{proof}

\begin{lemma}\label{lemma:Surjectivity Gr extension}
Suppose $\tilde{G}\to G$ is a central extension with kernel a torus. Then
$$
\mathrm{Gr}_{\tilde{G}}\to \mathrm{Gr}_G
$$
is a surjective map of v-sheaves.
\end{lemma}

\begin{proof}
Up to replacing $E$ by a finite degree extension we can suppose $G$ and $\tilde{G}$ are split. Fix $\tilde{T}\to T$ inside $\tilde{G}\to G$. According to Lemma~\ref{lemma:Affine Grassmanian Inductive limit} it is enough to prove that for any $\tilde{\mu}\in X_*(\tilde{T})^+$, if $\mu\in X_*(T)^+$ is its image in $G$, then 
$\mathrm{Gr}_{\tilde{G},\leq \tilde{\mu}}\to \mathrm{Gr}_{G,\leq \mu}$ is surjective. This is clearly surjective at the level of points since, if $D$ is the kernel of $\tilde{G}\to G$, then $H^1_{\et} ( \Spec (B_{\mathrm{dR}}(C)), D)=0$, and thus $\tilde{G} ( B_{\mathrm{dR}}(C))\to G(B_{\mathrm{dR}}(C))$ is surjective. Since $\mathrm{Gr}_{\tilde{G},\leq \tilde{\mu}}$ is quasicompact over $\Spd (E)$ and $\mathrm{Gr}_{G,\leq \mu}$ quasiseparated over $\Spd (E)$ (both are proper according to \cite[Proposition 19.2.3]{Berkeley}), $\mathrm{Gr}_{\tilde{G},\leq \tilde{\mu}}\to \mathrm{Gr}_{G,\leq \mu}$ is quasicompact and thus a v-cover by \cite[Lemma 12.11]{ECoD}.
\end{proof}

{\it Using Proposition~\ref{prop:beauvillelaszlosurj}, we thus have now a proof of Lemma~\ref{lem:zextension}.}
\\

Let us record a few facts we can deduce from the preceding results. Here, $\mu^\sharp\in \pi_1(G)$ denotes the image of $\mu\in X_\ast(T)^+$ under the quotient map $X_\ast(T)\to \pi_1(G)$.

\begin{proposition}
Suppose $G$ is split. 
\begin{altenumerate}
\item There is a locally constant map $|\mathrm{Gr}_G|\to \pi_1(G)$ inducing a decomposition in open/closed subsheaves $$\mathrm{Gr}_G=\coprod_{\alpha\in \pi_1(G)} \mathrm{Gr}_G^\alpha$$ 
characterized by $\mathrm{Gr}_{G,\mu}\subset \mathrm{Gr}_G^{\mu^\sharp}$.
\item The composite 
$$
|\mathrm{Gr}_G|\xrightarrow{\text{Beauville--Laszlo}} |\Bun_G | \xrightarrow{\ \kappa\ } \pi_1(G)
$$
is the opposite of the preceding map.
\item For each $\alpha\in \pi_1(G)$,
\[
\mathrm{Gr}_G^\alpha=\varinjlim_{\mu \in X_*(T)^+, \mu^\sharp= \alpha} \mathrm{Gr}_{G,\leq \mu}^\alpha
\]
as a filtered colimit of v-sheaves.
\end{altenumerate}
\end{proposition}
\begin{proof}
Point (1) is reduced to the case when $G_{\mathrm{der}}$ is simply connected using Lemma~\ref{lemma:Surjectivity Gr extension} and a $z$-extension. Now, if $G_{\mathrm{der}}$ is simply connected, for $\mu_1,\mu_2\in X_*(T)^+$, $\mu_1\leq \mu_2$ implies $\mu_1^\sharp=\mu_2^\sharp$. The result is then deduced from the fact that $\mathrm{Gr}_{G,\leq \mu}\subset \mathrm{Gr}_G$ is closed for any $\mu$.

Point (2) is can be similarly reduced first to the case that $G_{\mathrm{der}}$ is simply connected by passage to a z-extension; then to the case of a torus by taking the quotient by $G_{\mathrm{der}}$; then to the case of an induced torus by another z-extension; and then to $\mathbb G_m$ by changing $E$. In that case, it follows from Proposition~\ref{prop:lubintateuntilt}.

Point (3) is deduced from Lemma~\ref{lemma:Affine Grassmanian Inductive limit} and the fact that for $\alpha\in \pi_1(G)$, $\{\mu\in X_*(T)^+\ | \ \mu^\sharp =\alpha\}$ is a filtered ordered set. (Indeed, to see that this is filtered, note that if $\mu_1$ and $\mu_2$ satisfy $\mu_1^\sharp=\alpha=\mu_2^\sharp$, then $\mu_1-\mu_2$ lies in the coroot lattice, so is a sum of coroots with integer coefficients. Rearranging terms, it follows that one can add sums of positive coroots to $\mu_1$ and to $\mu_2$ so that they become equal, giving a common majorization.)
\end{proof}

When $G$ is not split, choosing $E'|E$ Galois of finite degree splitting $G$, using the formula $\mathrm{Gr}_G\times_{\Spd (E)} \Spd (E')= \mathrm{Gr}_{G_{E'}}$, one deduces:
\begin{altenumerate}
\item There is a decomposition $\mathrm{Gr}_G = \coprod_{\bar{\alpha}\in \Gamma\backslash \pi_1(G)} \mathrm{Gr}_G^{\bar{\alpha}}$ and a formula $\mathrm{Gr}_G=\varinjlim_{\bar{\mu}\in \Gamma\backslash X^* (T)^+} \mathrm{Gr}_{G,\leq \bar{\mu}}$ such that $\mathrm{Gr}_{G,\leq \bar{\mu}}\subset \mathrm{Gr}_G^{\bar{\mu^\sharp}}$.
\item The composite $|\mathrm{Gr}_G|\xrightarrow{\rm{BL}} |\Bun_G|\xrightarrow{\kappa} \pi_1 (G)_\Gamma$ is induced by the opposite of $\Gamma\backslash \pi_1(G) \to \pi_1 (G)_\Gamma$ and the preceding decomposition.
\end{altenumerate}

This description of the Kottwitz map, together with Proposition~\ref{prop:beauvillelaszlosurj}, in fact gives another proof of Theorem~\ref{thm:kappaconstant}.

\section{The semistable locus}

\subsection{Pure inner twisting}

Recall the following. In the particular case of non abelian group cohomology this is called ``torsion au moyen d'un cocycle'' in \cite[I.5.3]{SerreCohoGal}.

\begin{proposition}\label{prop:def inner twisting}
Let $X$ be a topos, $H$ a group in $X$ and $T$ an $H$-torsor. Let $H^T=\underline{\Aut} (T)$ as a group in $X$. Then:
\begin{altenumerate}
\item $H^T$ is the ``pure inner twisting'' of $H$ by $T$, $H^T=H\overset{H}{\wedge} T$ where $H$ acts by conjugation on $H$. In particular $[H^T]\in H^1(X,H_{\mathrm{ad}})$ is the image of $[T]$ via $H^1(X,H)\to H^1(X,H_{\mathrm{ad}})$.
\item The morphism of stacks on $X$, $[\ast /H]\to [\ast /{H^T}]$, that sends an $H$-torsor $S$ to $\underline{\mathrm{Isom}} (S,T)$, is an equivalence.
\end{altenumerate}
\end{proposition}

In the following we use the cohomological description of $G$-bundles on the curve as $G$-torsors on the \'etale site of the sous-perfectoid space $X_S$ (here $G$ is seen as an $E$-adic group, for $(R,R^+)$ a sous-perfectoid $E$-algebra its $\Spa (R,R^+)$-points being $G(R)$).

For the next proposition, recall (cf.~\cite[3.3]{KottwitzIsocrystals2}, where $G_b$ is denoted by $J$) that for any $b\in B(G)$, the automorphism group of the corresponding $G$-isocrystal defines a reductive group $G_b$ over $E$, via
\[
G_b(R)=\{g\in G(R\otimes_E \breve{E})\mid gb = b\sigma(g)\}.
\]
If $G$ is quasisplit, then $G_b$ is an inner form of a Levi subgroup of $G$. More generally $G_b$ is an inner form of a Levi subgroup of the quasisplit inner form of $G$. It is an inner form of $G$ precisely when $b$ is basic, i.e.~the Newton point is central. Recall from \cite{FarguesGBun} that a $G$-bundle over $X_{C,C^+}$ is semistable if and only if it corresponds to some basic element of $B(G)$; the reader may also take this as the definition of semistability for $G$-bundles.

\begin{proposition}\label{prop:pure inner twisting}
Let $S\in \Perf_k$, $b\in B(G)$ basic and $\E_b\to X_S$ the associated \'etale $G$-torsor. Then the \'etale sheaf of groups $G_b\times_{\Spa(E)} X_S$ over $X_S$ is the pure inner twisting of $G\times_{\Spa(E)} X_S$ by $\E_b$. 
\end{proposition}

\begin{proof}
One has 
$$
\E_b=  ( G_{\breve{E}}\times_{\Spa (\breve{E})} Y_S ) / ( (b\sigma)\times \phi )^\mathbb{Z} \longrightarrow X_S
$$
where $b\sigma$ acts on $G_{\breve{E}}$ by translation on the right. The $G\times_{\Spa(E)} X_S=(G_{\breve{E}}\times_{\Spa (\breve{E})} Y_S)/ (Id\times \phi)^\Z$-torsor structure is given by multiplication on the left on $G_{\breve{E}}$. The group $G_b\times_{\Spa(E)} X_S = (G_b\times Y_S)/(Id\times \phi)^\Z$ acts on this torsor on the right via the morphism $G_b\to G_{\breve{E}}$, which gives a morphism $$G_b\times_{\Spa(E)} X_S \to \underline{\Aut} (\E_b).$$ After pullback via the \'etale cover $Y_S\to X_S$ and evaluation on $T=\Spa (R,R^+)\to Y_S$ affinoid sous-perfectoid, this is identified with the map
$$
G_b(R)=\{g\in G(\breve{E}\otimes_E R)\ |\ g \cdot(b\sigma\otimes 1 )= (b\sigma \otimes 1)\cdot g\} \lto G(R)
$$ 
deduced from the $\breve{E}$-algebra structure of $R$. But as $b$ is basic, the natural map
\[
G_b\times_E \breve{E}\to G\times_E \breve{E}
\]
is an isomorphism, cf.~\cite[Corollary 1.14]{RapoportZink}.
\end{proof}

Thus, {\it extended pure inner forms, as defined by Kottwitz, become pure inner forms, as defined by Vogan, when pulled back to the curve.}

\begin{corollary}\label{cor:basictwistBunG}
For $b$ basic there is an isomorphism of v-stacks
$$
\Bun_G\simeq \Bun_{G_b}
$$
that induces an isomorphism $\Bun_G^b\simeq \Bun_{G_b}^1$.
\end{corollary}

\begin{example}
Take $G=\GL_n$ and $(D,\phi)$ an isoclinic isocrystal of height $n$. Let $B=\mathrm{End} (D,\phi)$ be the associated simple algebra over $E$. Since $(D,\phi)$ is isoclinic the action of $B$ on $\E(D,\phi)$ induces an isomorphism $B\otimes_E \O_{X_S} \iso \underline{\mathrm{End}} (\E(D,\phi))$ for any $S\in \Perf_{k}$. The stack $\Bun_{B^\times}$ is identified 
with the stack $S\mapsto \{\text{rank }1 \text{ locally free } B\otimes_E \O_{X_S}\text{-modules}\}$.
There is then an isomorphism (Morita equivalence)
\begin{align*}
\Bun_{\GL_n} & \xrightarrow{\ \sim\ } \Bun_{B^\times} \\
\E & \longmapsto \underline{\Hom}_{\O_{X_S}} (\E, \E(D,\phi)).
\end{align*}
\end{example}

\subsection{Description of the semi-stable locus}

As recalled in the previous section, a $G$-bundle on $X_{C,C^+}$ is semistable if the corresponding Newton point is central. A family of $G$-bundles is semistable if all of its geometric fibres are.

\begin{theorem}\label{cor:descrsemistable} The semistable locus
\[
\mathrm{Bun}_G^{\mathrm{ss}}\subset \mathrm{Bun}_G
\]
is open, and there is a canonical decomposition as open/closed substacks
\[
\mathrm{Bun}_G^{\mathrm{ss}}= \coprod_{b\in B(G)_{\mathrm{basic}}}  \mathrm{Bun}_G^b.
\]
For $b$ basic there is an isomorphism
\[
 [\ast/\underline{G_b(E)}]\iso \mathrm{Bun}_G^b.
\]
\end{theorem}

\begin{proof} Theorem~\ref{thm:newtonsemicont} implies that $\mathrm{Bun}_G^{\mathrm{ss}}$ is open, using that the condition that $\mu$ is central is a minimality condition in the dominance order. Recall that the basic elements of $B(G)$ map isomorphically to $\pi_1(G)_\Gamma$ via the Kottwitz map \cite[4.9, (4.4.1)]{KottwitzIsocrystals2}. Thus Theorem~\ref{thm:kappaconstant} gives a disjoint decomposition
\[
\mathrm{Bun}_G^{\mathrm{ss}} = \coprod_{b\in B(G)_{\mathrm{basic}}} \mathrm{Bun}_G^b .
\]
The result is then a consequence of Proposition~\ref{prop:pure inner twisting} and Theorem~\ref{thm:trivialGtorsor}.
\end{proof}

\begin{example}
For a torus $T$, $\Bun_T=\Bun_T^{ss}$ and there is an exact sequence of Picard stacks
$$
0\longrightarrow [\ast /\underline{T(E)}]\longrightarrow \Bun_T\longrightarrow \underline{X_*(T)_\Gamma}\longrightarrow 0,
$$
where we recall that $X_*(T)_\Gamma=B(T)$. The fiber of $\Bun_T\to \underline{B(T)}$ over $\beta$ is a gerbe banded by $\underline{T(E)}$ over $\ast$ via the action $\Bun_T^1$. This gerbe is neutralized after choosing some $b$ such that $[b]=\beta$. In case there is a section of $T(\breve{E})\twoheadrightarrow B(T)$, for example if $B(T)$ is torsion free, then $\Bun_T \simeq [\ast / \underline{T(E)}]\times \underline{X_*(T)_\Gamma}$ as a Picard stack.
\end{example}

\subsection{Splittings of the Harder--Narasimhan filtration}

We can also consider the following moduli problem, parametrizing $G$-bundles with a splitting of their Harder--Narasimhan filtration.

\begin{proposition}\label{prop:GBunHNsplit} Consider the functor $\Bun_G^{\mathrm{HN}\text-\mathrm{split}}$ taking each $S\in \Perf_k$ to the groupoid of exact $\otimes$-functors from $\Rep_E G$ to the category of $\mathbb Q$-graded vector bundles $\mathcal E=\bigoplus_\lambda \mathcal E^\lambda$ on $X_S$ such that $\mathcal E^\lambda$ is everywhere semistable of slope $\lambda$ for all $\lambda\in \mathbb Q$. For any $b\in B(G)$, the bundle $\E_b$ naturally refines to a $\mathbb Q$-graded bundle $\mathcal E_b^{\mathrm{gr}}$, using the $\mathbb Q$-grading on isocrystals, and for $S$ affinoid the natural map
\[
G_b\times_E X_S^{\mathrm{alg}}\to \underline{\Aut}(\mathcal E_b^{\mathrm{gr}})
\]
of group schemes over $X_S^{\mathrm{alg}}$ is an isomorphism. In particular, we get a natural map
\[
\bigsqcup_{b\in B(G)} [\ast/\underline{G_b(E)}]\to \Bun_G^{\mathrm{HN}\text-\mathrm{split}},
\]
and this is an isomorphism.
\end{proposition}

\begin{proof} Recall that the natural map $G_b\times_E \breve{E}\to G\times_E \breve{E}$, recording the map of underlying $\breve E$-vector spaces, is a closed immersion identifying $G_b\times_E \breve E$ with the centralizer of the slope homomorphism $\nu_b: \mathbb D\to G\times_E \breve{E}$, cf.~\cite[Corollary 1.14]{RapoportZink}. This implies that the natural map
\[
G_b\times_E X_S^{\mathrm{alg}}\to \underline{\Aut}(\mathcal E_b^{\mathrm{gr}})
\]
is an isomorphism.

We get the evident functor from $\bigsqcup_{b\in B(G)} [\ast/\underline{G_b(E)}]$ to this moduli problem, and it is clearly fully faithful. To see that it is surjective, take any strictly totally disconnected $S$ and an exact $\otimes$-functor $\mathcal E^{\mathrm{gr}}$ from $\Rep_E G$ to such $\mathbb Q$-graded vector bundles. For any point $s\in S$, note that $\mathbb Q$-graded vector bundles of the given form on $X_{K(s)}$ are equivalent to $\Isoc_E$, so at $s\in S$ there is an isomorphism with some $\mathcal E_b^{\mathrm{gr}}$. The type of the $\mathbb Q$-filtration is locally constant, so after replacing $S$ by an open neighborhood of $s$, we can assume that
\[
\underline{\mathrm{Isom}}(\mathcal E_b^{\mathrm{gr}},\mathcal E^{\mathrm{gr}})
\]
defines an $\underline{\Aut}(\mathcal E_b^{\mathrm{gr}})$-torsor over $X_S^{\mathrm{alg}}$, i.e.~a $G_b$-torsor over $X_S$. This defines a map $S\to \Bun_{G_b}$, taking $s$ into $\Bun_{G_b}^1$, and by Theorem~\ref{thm:trivialGtorsor} and Lemma~\ref{lemma:trivial loc sys} it follows that after replacing $S$ by an open neighborhood, we can assume that the torsor is trivial. This concludes the proof of surjectivity.
\end{proof}

\section{Non-semistable points}

\subsection{Structure of  $\underline{\Aut (\E_b)}$}

Next, we aim to describe the non-semi-stable strata $\mathrm{Bun}_G^b$. Before discussing the general case, consider the case $G=\GL_2$, with $b$ corresponding to the vector bundle $\mathcal O\oplus \mathcal O(1)$. The functor
\[
(S\in \Perf_k)\mapsto \mathrm{Aut}_{X_S}(\mathcal O_{X_S}\oplus \mathcal O_{X_S}(1))
\]
is given by the group v-sheaf
\[
\left(\begin{array}{cc} \underline{E^\times} & \BC(\O(1))\\ 0 & \underline{E^\times}\end{array}\right).
\]
This is an extension of the locally profinite group $G_b(E) = E^\times\times E^\times$ by a ``unipotent group'', namely a Banach--Colmez space $\BC(\mathcal O(1))$.

In general, fix any $b\in B(G)$ and consider the associated $G$-bundle $\mathcal E_b$ on $X_S$. For any algebraic representation $\rho: G\to \GL_n$, the corresponding vector bundle $\rho_\ast \mathcal E_b$ has its Harder--Narasimhan filtration $(\rho_\ast \mathcal E_b)^{\geq \lambda}\subset \rho_\ast \mathcal E_b$, $\lambda\in \mathbb Q$. If $G$ is quasisplit and we fix a Borel $B\subset G$, then this defines a reduction of $\mathcal E_b$ to a parabolic $P\subset G$ containing $B$.

Now inside the automorphism v-sheaf
$$\tilde{G}_b = \underline{\mathrm{Aut}}(\mathcal E_b): (S\in \Perf_k)\mapsto \mathrm{Aut}_{X_S}(\mathcal E_b|_{X_S})$$
(which necessarily preserves the Harder--Narasimhan filtration of $\rho_*\E_b$ for any $\rho \in \Rep_E(G)$) one can consider for any $\lambda>0$ the subgroup
\[
\tilde{G}_b^{\geq \lambda}\subset \tilde{G}_b
\]
of all automorphisms $\gamma: \mathcal E_b\iso \mathcal E_b$ such that
\[
(\gamma-1)(\rho_\ast \mathcal E_b)^{\geq \lambda'}\subset (\rho_\ast \mathcal E_b)^{\geq \lambda'+\lambda}
\]
for all $\lambda'$ and all representations $\rho$ of $G$. We also set $\tilde{G}_b^{>\lambda} = \bigcup_{\lambda'>\lambda} \tilde{G}_b^{\geq \lambda'}$, noting that this union is eventually constant.

As $G_b(E)$ is the automorphism group of the isocrystal corresponding to $b$, and $H^0(X_S,\O_{X_S})=\underline{E}(S)$, we have a natural injection
\[
\underline{G_b(E)}\hookrightarrow \tilde{G}_b.
\]
Now, for any automorphism $\gamma$ of $\E_b$ and any representation $\rho$, $\gamma$ induces an automorphism of the $\mathbb Q$-graded vector bundle 
$$
\bigoplus_{\lambda\in \Q} \mathrm{Gr}^\lambda ( \rho_* \E_b).
$$
Using Proposition~\ref{prop:GBunHNsplit}, we deduce that the preceding injection has a section and
$$
\tilde{G}_b= \tilde{G}_b^{>0} \rtimes \underline{G_b(E)}.
$$
For a $G$-bundle $\E$ on $X_S$ we note $\mathrm{ad}\,\E$ for its adjoint bundle deduced by pushforward by the adjoint representation $G\to \GL ( \mathrm{Lie} (G))$. This is in fact a Lie algebra bundle.

\begin{proposition}\label{prop:gradedpiecesoffancyGb} 
One has 
$$
\tilde{G}_b= \tilde{G}_b^{>0} \rtimes \underline{G_b(E)},
$$
and for any $\lambda>0$, there is a natural isomorphism
\[
\tilde{G}_b^{\geq \lambda}/\tilde{G}_b^{>\lambda}\xrightarrow{\ \sim\ } \mathcal{BC}((\mathrm{ad}\,\mathcal E_b)^{\geq \lambda}/(\mathrm{ad}\,\mathcal E_b)^{>\lambda}),
\]
the Banach--Colmez space associated to the slope $-\lambda$ isoclinic part of $(\Lie (G)\otimes_E \breve{E}, \mathrm{Ad}(b)\sigma)$.

In particular, $\tilde{G}_b$ is an extension of $\underline{G_b(E)}$ by a successive extension of positive Banach--Colmez spaces, and thus $\tilde{G}_b\to \ast$ is representable in locally spatial diamonds, of dimension $\langle 2\rho,\nu_b\rangle$ (where as usual $\rho$ is the half-sum of the positive roots).
\end{proposition}

We refer to \cite[Section IV]{SaavedraRivano} and \cite{Ziegler} for some general discussion of filtered and graded fibre functors.

\begin{proof} 
We already saw the first part. For the second part, suppose $S=\Spa (R,R^+)$ is affinoid. Let $X_R^{\mathrm{alg}}$ be the schematical curve. We use the GAGA correspondence, Proposition~\ref{prop:ampleGAGA}. Now, we apply Proposition \ref{prop:filtrations} to $X_R^{\mathrm{alg}}$ and the $G$-bundle $\E_b$ associated to $b$ on $X_R^{\mathrm{alg}}$. Let $H$ be the inner twisting of $G\times_E X_R^{\mathrm{alg}}$ by $\E_b$ as a reductive group scheme over $X_R^{\mathrm{alg}}$. It is equipped with a filtration $(H^{\geq \lambda})_{\lambda\geq 0}$ satisfying
\begin{itemize}
\item $H^{\geq 0}/H^{>0}\cong G_b\times_E X_R^{\mathrm{alg}}$,
\item for $\lambda>0$, $H^{\geq \lambda}/H^{>\lambda} = (\mathrm{ad}\, \E_b)^{\geq \lambda} / (\mathrm{ad}\, \E_b)^{>\lambda}$,
\item $\tilde{G}_b^{\geq \lambda} (S) = H^{\geq \lambda} (X_{R}^{\mathrm{alg}})$,
\end{itemize}
functorially in $(R,R^+)$; the first part uses Proposition~\ref{prop:GBunHNsplit}. Since $H^1 (X_R^{\mathrm{alg}},\O (\mu))=0$ as soon as $\mu >0$, we deduce by induction on $\mu>0$, starting with $\mu\gg 0$ and using the computation of $H^{\geq \mu}/H^{>\mu}$, that $H^1_{\et}(X_R^{\mathrm{alg}}, H^{\geq  \mu})=0$ for $\mu >0$. From this we deduce that 
$$
H^{\geq \lambda} (X_R^{\mathrm{alg}}) / H^{>\lambda} (X_R^{\mathrm{alg}}) = (H^{\geq \lambda}/ H^{>\lambda} )(X_R^{\mathrm{alg}}).
$$
Finally, it remains to compute the dimension. This is given by
\[
\sum_{\lambda>0} \lambda \cdot\dim\left((\mathrm{ad}\, \mathcal E_b)^{\geq \lambda}/(\mathrm{ad}\, \mathcal E_b)^{>\lambda}\right),
\]
which is given by $\langle 2\rho,\nu_b\rangle$.
\end{proof}

\begin{proposition}\label{prop:filtrations} Let $G$ be a reductive group over a field $K$, and let $X$ be a scheme over $K$. Let $\mathcal E$ be a $G$-bundle on $X$ with automorphism group scheme $H/X$ (an inner form of $G\times_K X$, cf. Proposition~\ref{prop:def inner twisting}). Consider a $\mathbb Q$-filtration on the fibre functor $\Rep_K (G) \to \{\text{Vector bundles on }X\}$
 associated with $\mathcal E$. Defining groups $H^{\geq \lambda}\subset H$ for $\lambda\geq 0$ as before, they are smooth group schemes, $H^{\geq 0}$ is a parabolic subgroup with unipotent radical $H^{>0}$, the Lie algebra of $H^{\geq \lambda}$ is given by $(\mathrm{ad}\, \mathcal E)^{\geq \lambda}\subset \Lie \mathrm{ad}\, \mathcal E=\Lie H$, and for $\lambda>0$ the quotient $H^{\geq \lambda}/H^{>\lambda}$ is a vector group, thus
\[
H^{\geq \lambda}/H^{>\lambda}\cong (\mathrm{ad}\, \mathcal E)^{\geq \lambda}/(\mathrm{ad}\, \mathcal E)^{>\lambda}
\]
\end{proposition}

\begin{proof}
The Lie algebra of $H$ is $\mathrm{ad}\,\E$. All statements can be checked \'etale locally on $X$. According to \cite[Theorem 1.3]{Ziegler} the $\mathbb{Q}$-filtration on the fiber functor is split locally on $X$. Moreover $\E$ is split \'etale locally on $X$. We can thus suppose that $\E$ is the trivial $G$-bundle and the filtration given by some $\nu:\mathbb{D}_{/X} \to G\times_K X$, where $\mathbb D$ is the pro-torus with character group $\mathbb Q$. Then the statement is easily checked, see \cite{SaavedraRivano}.
\end{proof}

\subsubsection{The quasi-split case} Suppose now that $G$ is moreover quasi-split. Fix $A\subset T\subset B$ with $A$ a maximal split torus and $T$ a maximal torus of $G$ inside a Borel subgroup $B$. Up to $\sigma$-conjugating $b$ one can suppose that $\nu_b:\mathbb{D}\to A$ and $\nu_b\in X_*(A)_{\mathbb Q}^+$, where $\mathbb D$ is the pro-torus over $E$ with character group $\mathbb Q$. Let $M_b$ be the centralizer of $\nu_b$ and $P_b^+$ the parabolic subgroup associated to $\nu_b$, the weights of $\nu_b$ in $\Lie (P_b^+)$ are $\geq 0$. One has $B\subset P_b^+$, $P_b^+$ is a standard parabolic subgroup with standard Levi subgroup $M_b$. Let $P_b^{-}$ be the opposite parabolic subgroup, the weights of $\nu_b$ in $\Lie (P_b^{-})$ are $\leq 0$. One has $b\in M_b (\breve{E})$ and we denote it $b_M$ as an element of $M_b (\breve{E})$.

Then, if 
\begin{align*}
Q& = \E_{b_M}\overset{M_b}{\times} P_b^-, \\
R_u Q &= \E_{b_M} \overset{M_b}{\times} R_u P_b^-
\end{align*}
as $X_R^{\mathrm{alg}}$-group-schemes, one has 
\begin{align*}
\tilde{G}_b ( R,R^+) &= Q (X_R^{\mathrm{alg}}) \\
\tilde{G}_b^{>0} (R,R^+) &= R_u Q (X_R^{\mathrm{alg}}).
\end{align*}

For $\GL_2$ and the bundle $\mathcal O\oplus \mathcal O(1)$, the group $Q$ is the upper triangular subgroup of the group scheme $\GL(\mathcal O\oplus \mathcal O(1))$ over $X_R^{\mathrm{alg}}$, and accordingly
\[
\tilde{G}_b=\left(\begin{array}{cc} \underline{E^\times} & \BC(\O(1))\\ 0 & \underline{E^\times}\end{array}\right).
\]

\subsection{Description of Harder--Narasimhan strata}

Now we can describe the structure of the stratum $\mathrm{Bun}_G^b$.

\begin{proposition}\label{prop:BunGb} Let $b\in B(G)$ be any element given by some $G$-isocrystal. The induced map $x_b: \ast\to \Bun_G^b$ is a surjective map of v-stacks, and $\ast\times_{\Bun_G^b} \ast\cong \tilde{G}_b$, so that
\[
\Bun_G^b\cong [\ast/\tilde{G}_b]
\]
is the classifying stack of $\tilde{G}_b$-torsors. In particular, the map $\tilde{G}_b\to \pi_0\tilde{G}_b \cong \underline{G_b(E)}$ induces a map
\[
\Bun_G^b\to [\ast/\underline{G_b(E)}]
\]
that admits a splitting.
\end{proposition}

\begin{proof} Let $S=\Spa (R,R^+)\in \Perf_k$ be strictly totally disconnected and let $\mathcal E$ be a $G$-bundle on $X_R^{\mathrm{alg}}$, the schematical curve, all of whose geometric fibers are isomorphic to $\mathcal E_b$. In particular, the Harder--Narasimhan polygon of $\rho_\ast \mathcal E$ is constant for all representations $\rho: G\to \GL_n$, and thus by Theorem~\ref{thm:kedlayaliu}, the vector bundle $\rho_\ast\mathcal E$ admits a relative Harder--Narasimhan filtration. This defines a $\mathbb{Q}$-filtration on the fiber functor $\Rep_{E} (G)\to \{\text{vector bundles on }X_R^{\mathrm{alg}}\}$ defined by $\E$, and exactness can be checked on geometric points where it holds by the classification of $G$-bundles. Since for any $\rho$, the Harder--Narasimhan polygon of $\rho_*\E_b$ and the one of $\rho_*\E$ are equal, the two filtered fiber functors on $\Rep_E(G)$ defined by $\E$ and $\E_b$ are of the same type. Thus, \'etale locally on $X_R^{\mathrm{alg}}$ those two filtered fiber functors are isomorphic. Let $H=\underline{\Aut} (\E_b)$ and $H^{\geq 0}= \underline{\Aut}_{\mathrm{filtered}} (\E_b)$ as group schemes over $X_R^{\mathrm{alg}}$, cf. Proposition~\ref{prop:gradedpiecesoffancyGb}. Now, look at 
$$
T=\underline{\mathrm{Isom}}_{\mathrm{filtered}} (\E_b,\E).
$$
This is an $H^{\geq 0}$-torsor over $X_R^{\mathrm{alg}}$ that is a reduction to $H^{\geq 0}$ of the $H$-torsor $\underline{\mathrm{Isom}} (\E_b,\E)$. Let us look at the image of $[T]\in H^1_{\et} (X_R^{\mathrm{alg}},H^{\geq 0})$ in $H^1_{\et} ( X_R^{\mathrm{alg}},H^{\geq 0}/H^{>0})$, that is to say the $H^{\geq 0}/H^{>0}$-torsor $T/H^{>0}$. This parametrizes isomorphisms of graded fiber functors between the two obtained by semi-simplifying the filtered fiber functors attached to $\E_b$ and $\E$. By Proposition~\ref{prop:GBunHNsplit}, this torsor is locally trivial. Now the triviality of $T$ follows from the vanishing of $H^1 (X_R^{\mathrm{alg}}, H^{>0})$. In fact, for $\lambda>0$, $H^1 (X_R^{\mathrm{alg}}, H^{\geq \lambda}/H^{>\lambda})=0$ since $H^1(X_R^{\mathrm{alg}},\O(\lambda))=0$.

It is clear that $\ast\times_{\Bun_G^b}\ast$ is given by $\tilde{G}_b$, so the rest follows formally.
\end{proof}

\begin{remark}[Followup to Remark \ref{remark:torseur representable I}]
From the vanishing of $H^1_{v} (S,\tilde{G}_b^{>0})$ for $S$ affinoid perfectoid one deduces that for such $S$, any $\tilde{G}_b$-torsor is of the form $T\times \tilde{G}_b^{>0}$ where $T\to S$ is a $\underline{G_b(E)}$-torsor. Here the action of $g_1\ltimes g_2\in \underline{G_b(E)}\ltimes \tilde{G}_b^{>0}$ 
on $T\times \tilde{G}_b^{>0}$ is given by $(x,y)\mapsto (g_1\cdot x, g_1 g_2 y g_1^{-1})$. In particular {\it any $\tilde{G}_b$-torsor is representable in locally spatial diamonds.}
\end{remark}

\chapter{Geometry of diamonds}

In this chapter, we extend various results on schemes to the setting of diamonds, showing that many advanced results in \'etale cohomology of schemes have analogues for diamonds.

In Section~\ref{sec:artinstacks}, we introduce a notion of Artin v-stacks, and discuss some basic properties; in particular, we show that $\Bun_G$ is a cohomologically smooth Artin v-stack. Moreover, we can define a notion of dimension for Artin v-stacks, which we use to determine the connected components of $\Bun_G$. In Section~\ref{sec:ULA}, we develop the theory of universally locally acyclic sheaves. In Section~\ref{sec:formalsmooth}, we introduce a notion of formal smoothness for maps of v-stacks. In Section~\ref{sec:jacobian}, we use the previous sections to prove the Jacobian criterion for cohomological smoothness, by establishing first formal smoothness, and universal local acyclicity. In Section~\ref{sec:partialcompactsupport}, we prove a result on the vanishing of certain partially compactly supported cohomology groups, ensuring that for example $\Spd k\powerseries{{x_1,\ldots,x_d}}$ behaves like a strictly local scheme for $D_\et$. In Section~\ref{sec:hyperboliclocalization}, we establish Braden's theorem on hyperbolic localization in the world of diamonds. Finally, in Section~\ref{sec:drinfeld}, we establish several version of Drinfeld's lemma in the present setup. The theme here is the idea $\pi_1((\Div^1)^I)=W_E^I$. Unfortunately, we know no definition of $\pi_1$ making this true, but for example it becomes true when considering $\Lambda$-local systems for any $\Lambda$.

\section{Artin stacks}\label{sec:artinstacks}

\subsection{Generalities}
\subsubsection{Definition and basic properties}

In this paper, we consider many small v-stacks like $\Bun_G$ as above. However, they are stacky in some controlled way, in that they are Artin v-stacks in the sense of the following definition.

\begin{definition}\label{def:artinvstack} An Artin v-stack is a small v-stack $X$ such that the diagonal $\Delta_X: X\to X\times X$ is representable in locally spatial diamonds, and there is some surjective map $f: U\to X$ from a locally spatial diamond $U$ such that $f$ is separated and cohomologically smooth.
\end{definition}

\begin{remark}\label{remark:Artin stack qs}  We are making the assumption that $f$ is separated, because only in this case we have defined cohomological smoothness. This means that we are imposing some (probably unwanted) very mild separatedness conditions on Artin v-stacks. In particular, it implies that {\it $\Delta_X$ is quasiseparated}: Let $f: U\to X$ be as in the definition, and assume without loss of generality that $U$ is a disjoint union of spatial diamonds (replacing it by an open cover if necessary), so in particular $U$ is quasiseparated. As $f$ is separated, the map $U\times_X U\to U$ is separated, and in particular $U\times_X U$ is again quasiseparated. This is the pullback of $\Delta_X: X\to X\times X$ along the surjection $U\times U\to X\times X$, so $\Delta_X$ is quasiseparated.
\end{remark}

\begin{remark}
{\it The stack $\Bun_G$ is not quasiseparated.} In fact,  $[\ast/\underline{G(E)}]$ is already not quasiseparated since the sheaf of automorphisms of the trivial $G$-bundle, $\underline{G(E)}
$, is not quasicompact. This is different from the ``classical situation'' of the stack of $G$-bundles on a proper smooth curve, this one being quasiseparated (although not separated). In the ``classical schematical case'' of Artin stacks it is a very mild assumption to suppose that Artin stacks are quasiseparated. In our situation this would be a much too strong assumption, but it is still a very mild assumption to suppose that the diagonal is quasiseparated.
\end{remark}

\begin{remark}\label{remark:diagonal representable} By Remark~\ref{remark:Artin stack qs}, for any Artin v-stack $X$, the diagonal $\Delta_X$ is quasiseparated. Conversely, let $X$ be any small v-stack, and assume that there is some surjective separated map $U\to X$ from a small v-sheaf. Then:
\begin{altenumerate}
\item[{\rm (i)}] If $U$ is quasiseparated, then $\Delta_X$ is quasiseparated, by the argument of Remark~\ref{remark:Artin stack qs}.
\item[{\rm (ii)}] If $U$ is a locally spatial diamond and $U\to X$ is representable in locally spatial diamonds, then $\Delta_X$ is quasiseparated (as we may without loss of generality assume that $U$ is quasiseparated, so that (i) applies), and to check that $\Delta_X$ is representable in locally spatial diamonds, it suffices to see that $\Delta_X$ is representable in diamonds. Indeed, \cite[Proposition 13.4 (v)]{ECoD} shows that if $\Delta_X$ is quasiseparated and representable in diamonds, then representability in locally spatial diamonds can be checked v-locally on the target. But the pullback of $\Delta_X$ along $U\times U\to X\times X$ is $U\times_X U$, which is a locally spatial diamond as we assumed that $U\to X$ is representable in locally spatial diamonds.
\item[{\rm (iii)}] Finally, in the situation of (ii), checking whether $\Delta_X$ is representable in diamonds can be done after pullback along a map $V\to X\times X$ that is surjective as a map of pro-\'etale stacks, by \cite[Proposition 13.2 (iii)]{ECoD}.
\end{altenumerate}

In particular, if there is a map $f: U\to X$ from a locally spatial diamond $U$ such that $f$ is separated, cohomologically smooth, representable in locally spatial diamonds, and surjective as a map of pro-\'etale stacks, then $X$ is an Artin v-stack. If one only has a map $f: U\to X$ from a locally spatial diamond such that $f$ is separated, cohomologically smooth, representable in locally spatial diamonds, and surjective as a map of v-stacks, then it remains to prove that $\Delta_X$ is representable in diamonds, which can be done after pullback along a map $V\to X\times X$ that is surjective as a map of pro-\'etale stacks.
\end{remark}

\begin{remark}\label{remark:coho smooth presentation surjectivity on points}
Since cohomologically smooth morphisms are open, to prove that a separated, representable in locally spatial diamonds, cohomologically smooth morphism $U\to X$ is surjective, it suffices to verify it on geometric points.
\end{remark}

\begin{remark}\label{remark:base artin stack} If $X$ is a small v-stack with a map $g: X\to S$ to some ``base'' small v-stack $S$, one might introduce a notion of an ``Artin v-stack over $S$'', asking instead that $\Delta_{X/S}: X\to X\times_S X$ is representable in locally spatial diamonds; note that the condition on the chart $f: U\to X$ will evidently remain the same as in the absolute case. We note that as long as the diagonal of $S$ is representable in locally spatial diamonds (for example, $S$ is an Artin v-stack itself), this agrees with the absolute notion. Indeed, if $\Delta_{X/S}$ and $\Delta_S$ are representable in locally spatial diamonds, then also $\Delta_X$ is representable in locally spatial diamonds, as $X\times_S X\to X\times X$ is a pullback of $\Delta_S$ and thus representable in locally spatial diamonds, so $\Delta_X$ is the composite of the two maps $X\to X\times_S X\to X\times X$ both of which are representable in locally spatial diamonds. Conversely, assume that $X$ and $S$ are such that their diagonals are representable in locally spatial diamonds. Then both $X$ and $X\times_S X$ are representable in locally spatial diamonds over $X\times X$, thus any map between them is.
\end{remark}

\begin{example}
Any locally spatial diamond is an Artin v-stack.
\end{example}

Before giving other examples let us state a few properties.

\begin{proposition}\label{prop:basic facts about Artin v stacks}\leavevmode
\begin{altenumerate}
\item Any fibre product of Artin v-stacks is an Artin v-stack.
\item Let $S\to \ast$ be a pro-\'etale surjective,  representable in locally spatial diamonds, separated and cohomologically smooth morphism of v-sheaves.  The v-stack $X$ is an Artin v-stack if and only if $X\times S$ is an Artin v-stack.
\item If $X$ is an Artin v-stack and $f: Y\to X$ is representable in locally spatial diamonds, then $Y$ is an Artin v-stack. 
\end{altenumerate}
\end{proposition}

\begin{proof}
For point (i), if $X=X_2\times_{X_1} X_3$ is such a fibre product and $f_i: U_i\to X_i$ are separated, representable in locally spatial diamonds, and cohomologically smooth surjective maps from locally spatial diamonds $U_i$, then $U=(U_1\times_{X_2} U_2)\times_{U_2} (U_2\times_{X_2} U_3)$ is itself a locally spatial diamond (using that $\Delta_{X_2}$ is representable in locally spatial diamonds), and the projection $f: U\to X$ is a separated, representable in locally spatial diamonds, and cohomologically smooth surjection. For the diagonal, since $\Delta_{X_2}$ and $\Delta_{X_3}$ are representable in locally spatial diamonds, $\Delta_{X_2}\times \Delta_{X_3}: X_2\times X_3\to (X_2\times X_3)\times (X_2\times X_3)$ is representable in locally spatial diamonds. Since $\Delta_{X_1}$ is representable in locally spatial diamonds, its pullback by $X_2\times X_3\to X_1\times X_1$, that is to say $u:X_2\times_{X_1} X_3\to X_2\times X_3$, is representable in locally spatial diamonds. Thus, $\Delta_{X_2\times_{X_1} X_3}$ is a map between stacks that are representable in locally spatial diamonds over $(X_2\times X_3)\times (X_2\times X_3)$, and thus is representable in locally spatial diamonds.
  
For point (ii), suppose $X\times S$ is an Artin v-stack. If $U$ is a locally spatial diamond and $U\to X\times S$ is separated, representable in locally spatial diamonds, cohomologically smooth, and surjective, then the composite $U\to X\times S\to X$ is too. It remains to see that $\Delta_X$ is representable in locally spatial diamonds. By Remark \ref{remark:diagonal representable} it suffices to prove that the pullback of $\Delta_X$ by $X\times X\times S\to X\times X$ is representable in locally spatial diamonds. But this pullback is the composite of $\Delta_{X\times S}$ with $X\times X\times S\times S\to X\times X\times S$, and we conclude since the projection $S\times S\to S$ is representable in locally spatial diamonds for evident reasons.

For point (iii), if $U$ is a locally spatial diamond and $U\to X$ is surjective, separated, representable in locally spatial diamonds, and cohomologically smooth, then $V=U\times_X Y$ is a locally spatial diamond, and $V\to Y$ is surjective, separated, representable in locally spatial diamonds, and cohomologically smooth. It remains to see that $\Delta_Y$ is representable in locally spatial diamonds. By Remark~\ref{remark:diagonal representable}, it suffices to see that $\Delta_Y$ is representable in diamonds. But we can write $\Delta_Y$ as the composite $Y\to Y\times_X Y\to Y\times_k Y$. The first map is $0$-truncated and injective and thus representable in diamonds by \cite[Proposition 11.10]{ECoD}, while the second map is a pullback of $\Delta_X$.
\end{proof}

We can now give more examples.

\begin{example}\label{ex:artinvstacks} \leavevmode
\begin{altenumerate}
\item According to point (ii) of Proposition \ref{prop:basic facts about Artin v stacks}, the v-stack $X$ is an Artin v-stack if and only if $X\times \Spd E$, resp. $X\times \Spa (\Fq\laurentseries{t^{1/p^\infty}})$, is an Artin v-stack. {\it To check that $X$ is an Artin v-stack we can thus replace the base point $\ast$ by $\Spd E$, resp.~$\Spa \Fq\laurentseries{t^{1/p^\infty}}$.}
\item For example, {\it any small v-sheaf $X$ such that $X\to \ast$ is representable in locally spatial diamonds is an Artin v-stack}; e.g.~$X=\ast$. 
\item Using point (iii) of Proposition \ref{prop:basic facts about Artin v stacks} and \cite[Proposition 11.20]{ECoD} we deduce that {\it any locally closed substack of an Artin v-stack is an Artin v-stack.}
\item Let $G$ be a locally profinite group that admits a closed embedding into $\GL_n(E)$ for some $n$. Then {\it the classifying stack $[\ast/\underline{G}]$ is an Artin v-stack.} For this 
it suffices to see that $[\Spd E/\underline{G}]=\Spd E\times [\ast/\underline{G}]$ is an Artin v-stack. Now let $H=\GL_{n,E}^\diamond$; then there is a closed immersion $\underline{G}\times \Spd E\hookrightarrow H$. The map $H\to \Spd E$ is representable in locally spatial diamonds, separated, and cohomologically smooth; hence so is $H/\underline{G}\to [\Spd E/\underline{G}]$ (by \cite[Proposition 13.4 (iv), Proposition 23.15]{ECoD}), and $H/\underline{G}$ is a locally spatial diamond (itself cohomologically smooth over $\Spd E$ by \cite[Proposition 24.2]{ECoD} since this becomes cohomologically smooth over the separated \'etale cover $H/\underline{K}\to H/\underline{G}$ for some compact open pro-$p$ subgroup $K$ of $G$). It is clear that the diagonal of $[\ast/\underline{G}]$ is representable in locally spatial diamonds.
\end{altenumerate}
\end{example}

\begin{remark}
If $G$ is a smooth algebraic group over the field $k$ then $\Spec (k)\to [\Spec (k)/G]$ is a smooth presentation of the Artin stack $[\Spec (k)/G]$. However, in the situation of point (4) of Example~\ref{ex:artinvstacks} the map $f:\ast \to [\ast / \underline{G}]$ is not cohomologically smooth (unless $G$ is finite) since for its pullback $\tilde{f}: \underline{G}\to \ast$, the sheaf $\tilde{f}^! \Lambda$ is the sheaf of distributions on $G$ with values in $\Lambda$.
\end{remark}

\subsubsection{Smooth morphisms of Artin v-stacks}

Notions that can be checked locally with respect to cohomologically smooth maps can be extended to Artin v-stacks (except possibly for subtleties regarding separatedness). In particular:

\begin{definition}\label{def:smoothartin} Let $f: Y\to X$ be a map of Artin v-stacks. Assume that there is some separated, representable in locally spatial diamonds, and cohomologically smooth surjection $g: V\to Y$ from a locally spatial diamond $V$ such that $f\circ g: V\to X$ is separated. Then $f$ is cohomologically smooth if for any (equivalently, one) such $g$, the map $f\circ g: V\to X$ (which is separated by assumption, and automatically representable in locally spatial diamonds) is cohomologically smooth.
\end{definition}

In the preceding definition the ``equivalently, one'' assertion is deduced from \cite[Proposition 23.13]{ECoD} that says that cohomological smoothness is ``cohomologically smooth local on the source''.  More precisely, if checked for one then for all $g:V\to X$ separated cohomologically smooth (not necessarily surjective) from a locally spatial diamond $V$, $f\circ g$ is separated cohomologically smooth.

\begin{convention}\label{conv:smoothartin} In the following, whenever we say that a map $f: Y\to X$ of Artin v-stacks is cohomologically smooth, we demand that there is some separated, representable in locally spatial diamonds, and cohomologically smooth surjection $g: V\to Y$ from a locally spatial diamond $V$ such that $f\circ g: V\to X$ is separated. Note that this condition can be tested after taking covers $U\to X$ by separated, representable in locally spatial diamonds, and cohomologically smooth surjections; i.e. after replacing $Y$ by $Y\times_X U$ and $X$ by $U$. If $X$ and $Y$ have the property that one can find a cover $U\to X$, $V\to Y$, as above with $U$ and $V$ perfectoid spaces, and $\Delta_X$ is representable in perfectoid spaces, then the condition is automatic, as all maps of perfectoid spaces are locally separated. That being said there is no reason that this is true in general since there are morphisms of spatial diamonds that are not locally separated.
\end{convention}

We will not try to give a completely general $6$-functor formalism that includes functors $Rf_!$ and $Rf^!$ for stacky maps $f$ (this would require some $\infty$-categorical setting). However, we can extend the functor $Rf^!$ to cohomologically smooth maps of Artin v-stacks. Let $\Lambda$ be a ring killed by some integer $n$ prime to $p$, or an adic ring as in \cite[Section 26]{ECoD}.

\begin{definition} Let $f: Y\to X$ be a cohomologically smooth map of Artin v-stacks. The dualizing complex $Rf^!\Lambda\in D_\et(Y,\Lambda)$ is the invertible object equipped with isomorphisms
\[
Rg^!(Rf^!\Lambda)\cong R(f\circ g)^! \Lambda
\]
for all separated, representable in locally spatial diamonds, and cohomologically smooth maps $g: V\to Y$ from a locally spatial diamond $V$, such that for all cohomologically smooth maps $h: V'\to V$ between such $g': V'\to Y$ and $g: V\to Y$, the composite isomorphism
\[
R(g')^!(Rf^!\Lambda)\cong R(f\circ g')^!\Lambda\cong R(f\circ g\circ h)^!\Lambda\cong Rh^!(R(f\circ g)^! \Lambda)\cong Rh^!(Rg^!(Rf^!\Lambda))\cong R(g')^!(Rf^!\Lambda)
\]
is the identity.
\end{definition}

As $Rf^!\Lambda$ is locally concentrated in one degree, it is easy to see that $Rf^!\Lambda$ is unique up to unique isomorphism. Let us be more precise. Let $\mathcal{C}$ be the category whose objects are separated cohomologically smooth morphisms $V\to Y$ with $V$ a locally spatial diamond, and morphisms $(V'\xrightarrow{g'} Y) \to (V\xrightarrow{g} Y)$ are couples $(h,\alpha)$ where $h:V'\to V$ is separated cohomologically smooth and $\alpha: g\circ h \Rightarrow g'$ is a $2$-morphism. Then the rule $$(V\xrightarrow{g} Y)\longmapsto R\sHom_{\Lambda}  (Rg^!\Lambda, R(f\circ g)^! \Lambda)$$ defines an element of
\[
2\text-\varprojlim_{(V\to Y)\in \mathcal{C}}\{ \text{invertible objects in } D_{\et} ( V,\Lambda)\}\cong \{ \text{invertible objects in } D_{\et} (Y,\Lambda)\}.
\]

\begin{remark}
If $g:V\to Y$ is a compactifiable  representable in locally spatial diamonds morphism of small v-stacks with $\dimtrg g<\infty$ such that $f\circ g$ satisfies the same hypothesis, it is not clear that $Rg^! (Rf^! \Lambda)\cong R (f\circ g)^!\Lambda$. This is a priori true only when $V$ is a locally spatial diamond and $g$ is separated cohomologically smooth, the only case we will need.
\end{remark}

\begin{definition} Let $f: Y\to X$ be a cohomologically smooth map of Artin v-stacks. The functor
\[
Rf^!: D_\et(X,\Lambda)\to D_\et(Y,\Lambda)
\]
is given by $Rf^! = Rf^!\Lambda\dotimes_\Lambda f^\ast$.
\end{definition}

\begin{remark} Checking after a cohomologically smooth cover, one sees that $Rf^!$ preserves all limits (and colimits) and hence admits a left adjoint $Rf_!$.
\end{remark}

\begin{definition} Let $f: Y\to X$ be a cohomologically smooth map of Artin v-stacks and let $\ell\neq p$ be a prime. Then $f$ is pure of $\ell$-dimension $d\in \tfrac 12 \mathbb Z$ if $Rf^! \mathbb F_\ell$ sits in homological degree $2d$.
\end{definition}

As $Rf^!\mathbb F_\ell$ is invertible, it is v-locally (and thus, a posteriori, \'etale locally) isomorphic to $\mathbb F_l[n]$ for some $n\in\mathbb Z$ (this can be deduced from Proposition~\ref{prop:constructible contre localement constant}~(ii)), so any cohomologically smooth map $f: Y\to X$ of Artin v-stacks decomposes uniquely into a disjoint union of $f_d: Y_d\to X$ that are pure of $\ell$-dimension $d$. A priori this decomposition may depend on $\ell$ and include half-integers $d$, but this will not happen in any examples that we study.

\subsection{The case of $\Bun_G$}
\subsubsection{Smooth charts on $\Bun_G$}

One important example is the following. We use Beauville--Laszlo uniformization to construct cohomologically smooth charts on $\Bun_G$. More refined charts will be constructed in Theorem~\ref{thm:chartbsmooth}. For $\bar{\mu} \in X_*(T)^+/\Gamma$ we note $\Gr_{G,\bar{\mu}}$ for the subsheaf of $\Gr_G$ such that 
$\Gr_{G,\bar{\mu}}\times_{\Spd (E)}\Spd (E')= \coprod_{\mu'\equiv \mu} \mathrm{Gr}_{G,\mu'}$ where $E'|E$ is a finite degree Galois extension splitting $G$. We will use the following simple proposition.

\begin{proposition}\label{prop:openschubertcell0} For any $\mu\in X_*(T)^+$, the open Schubert cell $\Gr_{G,\mu}/\Spd E'$ is cohomologically smooth of $\ell$-dimension $\langle 2\rho,\mu\rangle$.
\end{proposition}

We defer the proof to Proposition~\ref{prop:openschubertcell} as we do not want to make a digression on $\Gr_G$ here.

\begin{theorem}\label{thm:bunGartin} The stack $\Bun_G$ is a cohomologically smooth Artin v-stack of $\ell$-dimension $0$. The Beauville--Laszlo map defines a separated cohomologically smooth cover
$$
\coprod_{\bar{\mu} \in X_*(T)^+ /\Gamma} [\underline{G(E)}\backslash \mathrm{Gr}_{G,\bar{\mu}}]\lto \Bun_G.
$$
\end{theorem}

\begin{proof} We check first that $\Delta_{\Bun_G}$ is representable in locally spatial diamonds. For this, it suffices to see that for a perfectoid space $S$ with two $G$-bundles $\mathcal E_1,\mathcal E_2$ on $X_S$, the functor of isomorphisms between $\mathcal E_1$ and $\mathcal E_2$ is representable by a locally spatial diamond over $S$. By the Tannakian formalism, one can reduce to vector bundles.
For example, according to Chevalley, one can find a faithful linear representation $\rho:G\to \GL_n$, a representation $\rho':\GL_n\to \GL(W)$, and a line $D\subset W$ such that $G$ is the stabilizer of $D$ inside $\GL_n$. Then $G$-bundles on $X_S$ embed fully faithfully into rank $n$ vector bundles $\E$ together with a sub-line bundle $\mathcal{L}$ inside $\rho'_*\E$. In terms of those data, isomorphisms between $(\E_1,\mathcal{L}_1)$ and $(\E_2,\mathcal{L}_2)$ are given by a couple $(\alpha,\beta)$ where $\alpha:\E_1\iso \E_2$, and $\beta:\mathcal{L}_1\iso \mathcal{L}_2$ satisfy $(\rho'_*\alpha)_{|\mathcal{L}_1}=\beta$. Since the category of locally spatial diamonds is stable under finite projective limits we are reduced to the case of the linear group. Now the result is given by Lemma \ref{lemma:surjection ouvert}.

It remains to construct cohomologically smooth charts for $\Bun_G$. We first prove that the morphism 
$$
\pi: \coprod_{\bar{\mu} \in X_*(T)^+ /\Gamma} [\underline{G(E)}\backslash \mathrm{Gr}_{G,\bar{\mu}}]\lto \Bun_G\times_k \Spd \breve{E}
$$
is separated cohomologically smooth. Since this is surjective at the level of geometric points we deduce that it is a v-cover, cf. Remark \ref{remark:coho smooth presentation surjectivity on points}.

To verify this, note that for a perfectoid space $S$ mapping to $\Bun_G\times_k \Spd \breve{E}$ corresponding to a $G$-bundle $\mathcal E$ on $X_S$ as well as a map $S\to \Spd E$ inducing an untilt $S^\sharp/E$ and a closed immersion $i: S^\sharp\to X_S$, the fibre of $\pi$ over $S$ parametrizes modifications of $\mathcal E$ of locally constant type that are trivial at each geometric point. This is open in the space of all modifications of $\mathcal E$ of locally constant type, which is v-locally isomorphic to $\bigsqcup_{\bar\mu} \Gr_{G,\bar{\mu},E}\times_{\Spd E} S\to S$. Thus, Proposition~\ref{prop:openschubertcell0} gives the desired cohomological smoothness.

Moreover, the preceding argument shows that when restricted to $[\underline{G(E)}\backslash \mathrm{Gr}_{G,\bar{\mu}}]$, the map $\pi$ has $\ell$-dimension equal to $\langle 2\rho,\mu\rangle$. Thus, it now suffices to see that $[\underline{G(E)}\backslash \mathrm{Gr}_{G,\bar{\mu}}]$ is an $\ell$-cohomologically smooth Artin v-stack of $\ell$-dimension equal to $\langle 2\rho,\mu\rangle$. But the map
\[
[\underline{G(E)}\backslash \mathrm{Gr}_{G,\bar{\mu}}]\to [\Spd E/\underline{G(E)}]
\]
is representable in locally spatial diamonds and $\ell$-cohomologically smooth of $\ell$-dimension equal to $\langle 2\rho,\mu\rangle$, as $\mathrm{Gr}_{G,\bar\mu}\to \ast$ is by Proposition~\ref{prop:openschubertcell0}. We conclude by using that $[\ast/\underline{G(E)}]\to \ast$ is an Artin v-stack, cohomologically smooth of $\ell$-dimension $0$, by Example~\ref{ex:artinvstacks}~(4).
\end{proof}

\begin{lemma}\label{lemma:surjection ouvert}
For $\E_1, \E_2$ vector bundles on $X_S$, the sheaf $T/S\mapsto \{\text{surjections } \E_1|_{X_T}\twoheadrightarrow \E_2|_{X_T}\}$, resp. $T/S\mapsto \mathrm{Isom} (\E_1|_{X_T},\E_2|_{X_T})$, is representable by an open subdiamond of $\BC ( \E_1^\vee \otimes \E_2)$. In particular, those are locally spatial diamonds. 
\end{lemma}
\begin{proof}
The case of isomorphisms is reduced to the case of surjections since a morphism $u$ of vector bundles is an isomorphism if and only if $u$ and $u^\vee$ are surjective. For any morphism $g: \E_1\to \E_2$, the support of its cokernel is a closed subset of $|X_S|$, whose image in $|S|$ is thus closed; this implies the result.
\end{proof}

\begin{remark}
It would be tempting to study $D_{\et}(\Bun_G,\Lambda)$ using the preceding charts. But, contrary to the sheaves coming from the geometric Satake correspondence, the sheaves on $\mathrm{Gr}_G$ obtained via pullback from $\Bun_G$ are not locally constant on open Schubert strata. We will prefer other smooth charts to study $D_{\et}(\Bun_G,\Lambda)$, see Theorem~\ref{thm:chartbsmooth}.
\end{remark}

Moreover, each Harder--Narasimhan stratum $\Bun_G^b$ gives another example.

\begin{proposition}\label{prop:bunGbartin} For every $b\in B(G)$, the stratum $\Bun_G^b$ is a cohomologically smooth Artin v-stack of $\ell$-dimension $-\langle 2\rho,\nu_b\rangle$.
\end{proposition}

\begin{proof} Under the identification $\Bun_G^b\cong [\ast/\tilde{G}_b]$, note that we have a map $[\ast/\tilde{G}_b]\to [\ast/\underline{G_b(E)}]$ where the target is a cohomologically smooth Artin v-stack of dimension $0$, while the fibre admits a cohomologically smooth surjection from $\ast$ (as positive Banach--Colmez spaces are cohomologically smooth) of $\ell$-dimension $\langle 2\rho,\nu_b\rangle$. This gives the result.
\end{proof}

\subsubsection{Connected components of $\Bun_G$}

A consequence is that we can classify the connected components of $\Bun_G$.

\begin{corollary}\label{cor:pi0bunG} The Kottwitz map induces a bijection
\[
\kappa: \pi_0(\Bun_G)\to \pi_1(G)_\Gamma.
\]
\end{corollary}

\begin{proof} The Kottwitz map is well-defined and surjective. It remains to see that it is injective. To see this, recall that the basic elements of $B(G)$ biject via $\kappa$ to $\pi_1(G)_\Gamma$. Thus, it suffices to see that any nonempty open subsheaf $U$ of $\Bun_G$ contains a basic point. Note that the topological space $(X_*(T)_{\Q}^+)^\Gamma\times \pi_1(G)_\Gamma$ equipped with the product topology given by the order on  $(X_*(T)_{\Q}^+)^\Gamma$ and the discrete topology on $\pi_1(G)_\Gamma$, is (T$0$), and an increasing union of finite open subspaces; and $|\Bun_G|$ maps continuously to it. Pick some finite open $V\subset (X_*(T)_{\Q}^+)^\Gamma\times \pi_1(G)_\Gamma$ such that its preimage in $U$ is a nonempty open $U'\subset U$. Then $U'$ is a nonempty finite (T$0$) space, and thus has an open point by Lemma~\ref{lemma:finite T0 space open point}.

Thus, there is some $b\in B(G)$ such that $\Bun_G^b\subset U\subset \Bun_G$ is open. Combining Theorem~\ref{thm:bunGartin} and Proposition~\ref{prop:bunGbartin}, this forces $-\langle 2\rho,\nu_b\rangle = 0$, i.e.~$\nu_b$ is central, which means that $b$ is basic.
\end{proof}

\begin{lemma}\label{lemma:finite T0 space open point}
If $X$ is a nonempty finite spectral space, that is to say a finite (T$0$) topological space, there exists an open point $x\in X$.
\end{lemma}
\begin{proof}
Take $x$ maximal for the specialization relation, i.e. $x$ is a maximal point. Then, since $X$ is (T$0$), $X\setminus \{x\} = \cup_{y\neq x} \overline{\{y\}}$, a finite union of closed spaces thus closed.
\end{proof}

\section{Universally locally acyclic sheaves}\label{sec:ULA}
\subsection{Definition and basic properties}

In many of our results, and in particular in the (formulation and) proof of the geometric Satake equivalence, a critical role is played by the notion of universally locally acyclic (ULA) sheaves. Roughly speaking, for a morphism $f: X\to S$ of schemes, these are constructible complexes of \'etale sheaves $A$ on $X$ whose relative cohomology is constant in all fibres of $S$, even locally. Technically, one requires that for all geometric points $\overline x$ of $X$ mapping to a geometric point $\overline s$ of $S$ and a generization $\overline t$ of $\overline s$ in $S$, the natural map
\[
R\Gamma(X_{\overline x},A)\to R\Gamma(X_{\overline x}\times_{S_{\overline s}} \overline t,A)
\]
is an isomorphism, where $X_{\overline x}$ is the strict henselization of $X$ at $\overline x$ (and $S_{\overline s}$ is defined similarly). Moreover, the same property should hold universally after any base change along $S^\prime\to S$.\footnote{Recently, Gabber proved that this is automatic when $S$ is noetherian and $f$ is of finite type, cf.~\cite[Corollary 6.6]{LuZhengGeneralBase}.} By \cite[Corollary 3.5]{IllusieGeneralNearby}, universal local acyclicity is equivalent to asking that, again after any base change, the map
\[
R\Gamma(X_{\overline x},A)\to R\Gamma(X_{\overline x}\times_{S_{\overline s}} S_{\overline t},A)
\]
is an isomorphism; we prefer the latter formulation as strict henselizations admit analogues for adic spaces, while the actual fibre over a point is only a pseudo-adic space in Huber's sense \cite{Huber}.

In the world of adic spaces, there are not enough specializations to make this an interesting definition; for example, there are no specializations from $\Gr_{G,\mu}$ into $\Gr_{G,\leq\mu}\setminus \Gr_{G,\mu}$. Thus, we need to adapt the definition by adding a condition on preservation of constructibility that is automatic in the scheme case under standard finiteness hypothesis, but becomes highly nontrivial in the case of adic spaces. Here, for a diamond $X$ with a geometric point $\overline{x}$, we let $X_{\overline{x}}=\Spa(C(\overline{x}),C(\overline{x})^+)$ be the strict localization of $X$ at $\overline{x}$ (which is the initial diamond pro-\'etale over $X$ with a lift of $\overline{x}$).

\begin{definition}\label{def:ULA} Let $f: X\to S$ be a compactifiable map of locally spatial diamonds with locally $\dimtrg f<\infty$ and let $A\in D_\et(X,\Lambda)$ for some $\Lambda$ with $n\Lambda=0$ with $n$ prime to $p$.

\begin{altenumerate}
\item[{\rm (i)}] The sheaf of complexes $A$ is $f$-locally acyclic if
\begin{altenumerate}
\item[{\rm (a)}] For all geometric points $\overline x$ of $X$ with image $\overline s$ in $S$ and a generization $\overline t$ of $\overline s$, the map
\[
R\Gamma(X_{\overline x},A)\to R\Gamma(X_{\overline x}\times_{S_{\overline s}} S_{\overline t},A)
\]
is an isomorphism.
\item[{\rm (b)}] For all separated \'etale maps $j: U\to X$ such that $f\circ j$ is quasicompact, the complex $R(f\circ j)_!(A|_U)\in D_\et(S,\Lambda)$ is perfect-constructible.
\end{altenumerate}
\item[{\rm (ii)}] The sheaf of complexes $A$ is $f$-universally locally acyclic if for any map $S^\prime\to S$ of locally spatial diamonds with base change $f^\prime: X^\prime=X\times_S S^\prime\to S^\prime$ and $A^\prime\in D_\et(X^\prime,\Lambda)$ the pullback of $A$, the sheaf of complexes $A^\prime$ is $f^\prime$-locally acyclic.
\end{altenumerate}
\end{definition}

Recall that if $(K,K^+)$ is an affinoid Huber field, $S=\Spa (K,K^+)$, then $|S|=|\Spec (K^+/K^{00})|$ as a topological spectral space, that is identified with the totally ordered set of open prime ideals in  $K^+$. For any $s\in S$, $S_s\subset S$ is pro-constructible generalizing. For example, the maximal generalization is $\Spa (K,\O_K)=\cap_{a\in \O_K} \{ |a|\leq 1\}$.

\begin{remark}\label{rem:defULA} In the setup of condition (a), note that $X_{\overline x}$ is a strictly local space, i.e.~of the form $\Spa(C,C^+)$ where $C$ is algebraically closed and $C^+\subset C$ is an open and bounded valuation subring; thus, $R\Gamma(X_{\overline x},A)=A_{\overline x}$ is just the stalk of $A$. Moreover, $S_{\overline t}\subset S_{\overline s}$ is a quasicompact pro-constructible generalizing subspace, and thus $X_{\overline x}\times_{S_{\overline s}} S_{\overline t}\subset X_{\overline x}$ is itself a quasicompact pro-constructible generalizing subset that is strictly local. Its closed point $\overline y$ is the minimal generization of $\overline x$ mapping to $\overline t$, and $R\Gamma(X_{\overline x}\times_{S_{\overline s}} S_{\overline t},A) = A_{\overline y}$ is the stalk at $\overline y$. Thus, condition (a) means that $A$ is ``overconvergent'' along the horizontal lifts of generizations of $S$.
\end{remark}

\begin{remark}\label{remark:autre definition locale acyclicite}\leavevmode
\begin{altenumerate}
\item Another way to phrase the ``relative overconvergence condition'' (a), is to say that if $\bar{s}$ is a geometric point of $S$, $\bar{t}$ a generization of $\bar{s}$, $j:X\times_S S_{\bar{t}} \hookrightarrow X\times_{S} S_{\bar{s}}$, a pro-constructible generalizing immersion, and $B=A|_{X\times_S S_{\bar{s}}}$, then $B=Rj_* j^* B$ (use quasicompact base change).
\item Still another way to phrase it is to say that for any $\Spa (C,C^+)\to S$, if $B=A|_{X\times_S \Spa (C,C^+)}$, and $j:X\times_S \Spa (C,\O_C)\hookrightarrow X\times_S \Spa (C,C^+)$, then $B=Rj_*j^*B$.
\item Still another way is to say that if $\bar{x}\mapsto \bar{s}$ and $f_{\bar{x}}:X_{\bar{x}}\to S_{\bar{s}}$ then $Rf_{\bar{x}*} A|_{X_{\bar{x}}}$ is overconvergent i.e.~constant.
\end{altenumerate}
\end{remark}

In fact, asking for condition (a) universally, i.e.~after any base change, amounts to asking that $A$ is overconvergent.

\begin{proposition}\label{prop:universallyhorizontallyoverconvergent} Let $f: X\to S$ be a compactifiable map of locally spatial diamonds with locally $\dimtrg f<\infty$ and let $A\in D_\et(X,\Lambda)$ for some $\Lambda$ with $n\Lambda=0$ with $n$ prime to $p$. The condition (a) of Definition~\ref{def:ULA} holds after any base change $S'\to S$ if and only if $A$ is overconvergent, i.e.~for any specialization $\overline{y}\rightsquigarrow \overline{x}$ of geometric points of $X$, the map $A_{\overline{x}}\to A_{\overline{y}}$ is an isomorphism.
\end{proposition}

\begin{proof} The condition is clearly sufficient. For necessity, take the base change along $X_{\overline{x}}\to S$. Then $\overline{x}$ lifts to a section $\overline{x}': X_{\overline{x}}\to X\times_S X_{\overline{x}}$, and applying the relative overconvergence condition to $\overline{x}'\mapsto \overline{x}$ and the generization $\overline{y}$ of $\overline{x}$, we see that $A_{\overline{x}}\to A_{\overline{y}}$ is an isomorphism.
\end{proof}

\begin{proposition}\label{prop:locally acyclic v local}
Local acyclicity descends along v-covers of the target. More precisely, in the setup of Definition \ref{def:ULA}, if $S^\prime\to S$ is a v-cover and $A^\prime$ is $f^\prime$-locally acyclic, then automatically $A$ is $f$-locally acyclic.
\end{proposition}
\begin{proof}
Condition (a) follows by lifting geometric points, and condition (b) descends by \cite[Proposition 20.13]{ECoD}.
\end{proof}

\begin{proposition}\label{prop:constructible contre localement constant}
Let $Y$ be a spatial diamond.
\begin{altenumerate}
\item If $\mathcal{F}$ is a constructible \'etale sheaf of $\Lambda$-modules on $Y$, then $\mathcal{F}$ is locally constant if and only if $\mathcal{F}$ is overconvergent.
\item If $A\in D_{\et,pc}(Y,\Lambda)$, then $A$ is overconvergent if and only if it is locally a constant perfect complex of $\Lambda$-modules.
\end{altenumerate}
\end{proposition}

\begin{proof}
For a geometric point $\overline{y}$ of $Y$, writing $Y_{\overline{y}}=\Spa(C,C^+)=\varprojlim_{\overline{y}\to U} U$ as a limit of the \'etale neighborhoods, according to \cite[Proposition 20.7]{ECoD}, 
$$
2\text-\varinjlim_{\overline{y}\to U} \mathrm{Cons} (U,\Lambda) = \mathrm{Cons} (Y_{\overline{y}},\Lambda).
$$
An \'etale sheaf on $Y_{\overline{y}}=\Spa(C,C^+)$ is locally constant if and only if it is constant if and only if it is overconvergent. This gives point (1). Point (2) goes the same way using \cite[Proposition 20.15]{ECoD}. 
\end{proof}

\begin{remark}
The preceding argument shows that if $\mathcal{F}$ is constructible on $Y$ then $\mathcal{F}$ is locally constant in a neighborhood of any maximal point of $Y$. For example, if $Y=X^\diamond$ with $X$ a $K$-rigid space, then any constructible sheaf on $Y$ is locally constant in a neighborhood of all classical Tate points of $X$. Thus, the difference between constructible and locally constant sheaves shows up at rank $>1$ valuations.
\end{remark}

\begin{example}
Let $j:\mathbb{B}^1_K\setminus \{0\}\hookrightarrow \mathbb{B}^1_K$ be the inclusion of the punctured disk inside the  disk. Then $j_!\Lambda$ is not constructible since not locally constant around $\{0\}$. Nevertheless, if $R\in |K^\times|$ and $x$ is the coordinate on $\mathbb{B}^1_K$,  $j_R: \{ R\leq |x|\leq 1\} \hookrightarrow \mathbb{B}^1_K$, $j_{R !}\Lambda$ is constructible and $j_!\Lambda = \varinjlim_{R \to 0} j_{R!}\Lambda$. The category of \'etale sheaves of $\Lambda$-modules on a spatial diamond is the Ind-category of constructible \'etale sheaves, cf.~\cite[Proposition 20.6]{ECoD}.
\end{example}

\begin{proposition}\label{prop:ULAbase} Assume that $f: S\to S$ is the identity. Then $A\in D_\et(S,\Lambda)$ is $f$-locally acyclic if and only if it is locally constant with perfect fibres.
\end{proposition}

\begin{proof} Applying part (b) of the definition, we see that $A$ is perfect-constructible. On the other hand, part (a) says that $A$ is overconvergent. This implies that $A$ is locally constant by Proposition~\ref{prop:constructible contre localement constant}.
\end{proof}

Let us finish with a basic example of universally locally acyclic sheaves relevant to the smooth base change theorem. A more general result will be given in Proposition \ref{prop:ULAsmoothlocal}.

\begin{proposition}\label{prop:smoothULA} Assume that $f: X\to S$ is a separated map of locally spatial diamonds that is $\ell$-cohomological smooth for all divisors $\ell$ of $n$, where $n\Lambda=0$. If $A\in D_\et(X,\Lambda)$ is locally constant with perfect fibres, then $A$ is $f$-universally locally acyclic.
\end{proposition}

\begin{proof} It is enough to show that $A$ is $f$-locally acyclic, as the hypotheses are stable under base change. Condition (a) follows directly from $A$ being locally constant. Condition (b) follows from the preservation of constructible sheaves of complexes under $Rf_!$ if $f$ is quasicompact, separated and cohomologically smooth, see \cite[Proposition 23.12 (ii)]{ECoD}.
\end{proof}

\subsection{Proper push-forward, smooth pull-back}

In the ``classical algebraic case'', if $Y\xrightarrow{g} X\xrightarrow{f} S$ are morphisms of finite type between noetherian schemes, using proper and smooth  base change:
\begin{altenumerate}
\item if $g$ is proper and $A\in D_c^b (Y,\Lambda)$ is $f\circ g$ -locally acyclic then $Rg_* A$ is $f$-locally acyclic;
\item  if $g$ is smooth and $A\in D_c^b (X,\Lambda)$ is $f$-locally acyclic then $g^* A$ is $f\circ g$-locally acyclic. Moreover if $g$ is surjective then  $A\in D_c^b (X,\Lambda)$ is $f$-locally acyclic if and only if $g^* A$ is $f\circ g$-locally acyclic.
\end{altenumerate}
 We are going to see that the same phenomenon happens in our context. The fact that local acyclicity is smooth local on the source is essential to define local acyclicity for morphisms of Artin v-stacks, cf.~Definition~\ref{def:local acyclicity for Artin stacks}.

\begin{proposition}\label{prop:properULA} Let $g: Y\to X$, $f: X\to S$ be maps of locally spatial diamonds where $g$ is proper and $f$ is compactifiable and locally $\dimtrg g,\dimtrg f<\infty$. Assume that $A\in D_\et(Y,\Lambda)$ is $f\circ g$-locally acyclic (resp.~$f\circ g$-universally locally acylic). Then $Rg_\ast A\in D_\et(X,\Lambda)$ is $f$-locally acyclic (resp.~$f$-universally locally acyclic).
\end{proposition}

\begin{proof} It is enough to consider the locally acyclic case, as the hypotheses are stable under base change. For condition (a), we use Remark \ref{remark:autre definition locale acyclicite}~(2). Let $\bar{s}$ be a geometric point of $S$, with $S_{\bar{s}} = \Spa(C,C^+)$. Let us look at the cartesian diagram 
$$
\begin{tikzcd}
Y\times_{\Spa(C,C^+)} \Spa(C,\mathcal O_C) \ar[d] \ar[r,hook,"k"] & Y \ar[d,"g"] \\
X\times_{\Spa(C,C^+)} \Spa(C,\mathcal O_C) \ar[r,hook,"j"] & X
\end{tikzcd}
$$
one has by local acyclicity of $f\circ g$, $A=Rk_*k^*A$. Applying $Rg_\ast$, this gives the desired
$$
Rg_* A= Rj_* j^* ( Rg_* A ).
$$

For condition (b), take any separated \'etale map $j: U\to X$ such that $f\circ j$ is quasicompact, and set $j^\prime: V=U\times_X Y\to Y$, which is an \'etale map such that $f\circ g\circ j^\prime$ is quasicompact. Let $g^\prime: V\to U$ denote the pullback of $g$. Using proper base change and $Rg_\ast = Rg_!$, we see that
\[
R(f\circ j)_! j^\ast Rg_\ast A = R(f\circ j)_! Rg^\prime_! j^{\prime\ast} A = R(f\circ g\circ j^\prime)_! j^{\prime\ast} A,
\]
which is perfect-constructible by the assumption that $A$ is $f\circ g$-locally acyclic.
\end{proof}

In particular we have the following that generalizes the ``proper and smooth case''.

\begin{corollary}\label{coro:proper push locally acyclic}
Let $f:X\to S$ be a proper map of locally spatial diamonds with $\dimtrg f<\infty$ and $A\in D_{\et} (X, \Lambda)$ that is $f$-locally acyclic. Then $Rf_* A$ is locally a constant perfect complex of $\Lambda$-modules.
\end{corollary}

Next proposition says that {\it local acyclicity is ``cohomologically smooth local'' on the source.}

\begin{proposition}\label{prop:ULAsmoothlocal} Let $f: X\to S$ be a compactifiable map of locally spatial diamonds with locally $\dimtrg f<\infty$. For the statements in the locally acyclic case below, assume that $S$ is spatial and that the cohomological dimension of $U_\et$ for all quasicompact separated \'etale $U\to S$ is $\leq N$ for some fixed integer $N$.

Let $A\in D_\et(X,\Lambda)$ where $n\Lambda=0$ for some $n$ prime to $p$ and let $g: Y\to X$ be a separated map of locally spatial diamonds that is $\ell$-cohomologically smooth for all $\ell$ dividing $n$.
\begin{altenumerate}
\item[{\rm (i)}] If $A$ is $f$-locally acyclic (resp.~$f$-universally locally acyclic), then $g^\ast A$ is $f\circ g$-locally acyclic (resp.~$f\circ g$-universally locally acyclic).
\item[{\rm (ii)}] Conversely, if $g^\ast A$ is $f\circ g$-locally acyclic (resp.~$f\circ g$-universally locally acyclic) and $g$ is surjective, then $A$ is $f$-locally acyclic (resp.~$f$-universally locally acyclic).
\end{altenumerate}
\end{proposition}

\begin{proof} It is enough to handle the locally acyclic case with the assumption on $S$; then the universally locally acyclic case follows by testing after pullbacks to strictly totally disconnected spaces, using Proposition~\ref{prop:locally acyclic v local}.
Let us treat point (i).
 We can assume that $X$ and $Y$ are qcqs, i.e.~spatial.
In fact, this is clear for condition (a). For condition (b), if $j:V\to Y$ is separated \'etale such that $f\circ g\circ j$ is quasicompact, up to replacing $S$ by an open cover we can suppose that $S$ is spatial and thus $V$ is spatial (since $f$, $g$, and $j$ are separated, $f\circ g\circ j$ is separated quasicompact, and thus $S$ spatial implies $X$ spatial). Then one can replace $Y$, resp. $X$, by the quasicompact open subsets $j(V)$, resp. $(g\circ j)(V)$, that are separated over $S$ and thus spatial too. 
 
Condition (a) follows as pullbacks preserve stalks. For condition (b), let $j: V\to Y$ be any quasicompact separated \'etale map. Then by the projection formula for $g\circ j$, one has
\[
R(f\circ g\circ j)_! j^\ast g^\ast A = Rf_!(A\dotimes_\Lambda R(g\circ j)_! \Lambda).
\]
As $g\circ j: V\to X$ is a quasicompact separated $\ell$-cohomologically smooth map, it follows that $R(g\circ j)_!\Lambda\in D_\et(X,\Lambda)$ is perfect-constructible by \cite[Proposition 23.12 (ii)]{ECoD}. Thus, the desired result follows from Lemma \ref{lemma:ula produit tensoriel}.

In the converse direction, i.e.~for part (ii), condition (a) of $A$ being $f$-locally acyclic follows by lifting geometric points from $X$ to $Y$ and noting that stalks do not change. For condition (b), we may replace $X$ by $U$ to reduce to the assertion that $Rf_! A\in D_\et(S,\Lambda)$ is perfect-constructible. Consider the thick triangulated subcategory $\mathcal C$ of $D_\et(X,\Lambda)$ of all $B\in D_\et(X,\Lambda)$ such that $Rf_!(A\dotimes_\Lambda B)\in D_\et(S,\Lambda)$ is perfect-constructible. We have to see that $\Lambda\in \mathcal{C}$. We know that for all perfect-constructible $C\in D_\et(Y,\Lambda)$, the perfect-constructible complex $Rg_! C$ lies in $\mathcal C$. Indeed, 
using \cite[Proposition 20.17]{ECoD}, 
this reduces to the case $C=j_!\Lambda$ where $j: U\to Y$ is a quasicompact separated \'etale map, and then
\[
Rf_!(A\dotimes_\Lambda R(g\circ j)_!\Lambda) = R(f\circ g)_!(g^\ast A\dotimes_\Lambda Rj_!\Lambda),
\]
which is perfect-constructible as $g^\ast A$ is $f\circ g$-locally acyclic. Thus, it is enough to show that the set of $Rg_! C\in D_\et(X,\Lambda)$ with $C\in D_\et(Y,\Lambda)$ perfect constructible form a set of compact generators of $D_\et(X,\Lambda)$. Equivalently, for any complex $B\in D_\et(X,\Lambda)$ with $R\Hom_{D_\et(X,\Lambda)}(Rg_!C,B)=0$ for all perfect-constructible $C\in D_\et(Y,\Lambda)$, then $B=0$. The hypothesis is equivalent to $R\Hom_{D_\et(Y,\Lambda)}(C,Rg^!B)=0$ for all such $C$. By \cite[Proposition 20.17]{ECoD} and the standing assumptions on finite cohomological dimension (on $S$, $f$ and $g$), this implies that $Rg^!B=0$. As $g$ is $\ell$-cohomologically smooth, this is equivalent to $g^\ast B=0$, which implies $B=0$ as $g$ is surjective.
\end{proof}

\begin{lemma}\label{lemma:ula produit tensoriel}
Let $f:X\to S$ be a compactifiable map of locally spatial diamonds with locally $\dimtrg f<\infty$. Suppose there exists an integer $N$ such that the cohomological dimension of $U_{\et}$ is bounded by $N$ for all $U\to X$ separated \'etale. Let $A\in D_{\et} (X,\Lambda)$ be $f$-locally acyclic and $B\in D_{\et} (X,\Lambda)$ be perfect-constructible. Then $A\dotimes_{\Lambda}B$ satisfies condition (b) of Definition \ref{def:ULA}: for any $j:U\to X$ separated \'etale such that $f\circ j$ is quasicompact, $R(f\circ j)_! j^*(A\dotimes_\Lambda B)$ is perfect-constructible.
\end{lemma}
\begin{proof}
We can suppose $X$ is spatial.
According to \cite[Proposition 20.17]{ECoD}, $B$ lies in the triangulated subcategory generated by $j'_! \Lambda$ where $j':U'\to X$ is separated quasicompact  \'etale. For such a $B$, using the projection formula, $A\dotimes_\Lambda Rj'_!\Lambda=Rj'_! j'^*A$. Thus, if $V=U\times_X U'$ with projection $k:U\times_X U' \to U$, 
$$
j^* Rj'_! j'^*A= Rk_! k^*j^*A.
$$ 
We thus have 
$$
R(f\circ j)_! j^* ( A\dotimes_\Lambda Rj'_! B) = R (f\circ j\circ k)_! (j\circ k)^* (A)
$$
and we can conclude.
\end{proof}

\subsection{Local acyclicity and duality}

In this section, we prove that universal local acyclicity behaves well with respect to Verdier duality.

\subsubsection{Compatibility with base change}

We note that for $f$-ULA sheaves, the formation of the (relative) Verdier dual
\[
\mathbb D_{X/S}(A):=R\sHom_\Lambda(A,Rf^!\Lambda)
\]
commutes with base change in $S$.

\begin{proposition}\label{prop:ULAdual} Let $f: X\to S$ be a compactifiable map of locally spatial diamonds with locally $\dimtrg f<\infty$ and let $A\in D_\et(X,\Lambda)$ be $f$-universally locally acyclic. Let $g: S^\prime\to S$ be a map of locally spatial diamonds with pullback $f^\prime: X^\prime=X\times_S S^\prime\to S^\prime$, $\tilde{g}: X^\prime\to X$. Then the composite
\[
\tilde{g}^\ast \mathbb D_{X/S}(A)\to \mathbb D_{X^\prime/S^\prime}(\tilde{g}^\ast A)
\]
of the natural maps
\[
\tilde{g}^\ast R\sHom_\Lambda(A,Rf^!\Lambda)\to R\sHom_\Lambda(\tilde{g}^\ast A,\tilde{g}^\ast Rf^!\Lambda)\to R\sHom_\Lambda(\tilde{g}^\ast A,Rf^{\prime!}\Lambda)
\]
is an isomorphism.

More generally, for any $B\in D_\et(S,\Lambda)$, the map
\[
\tilde{g}^\ast R\sHom_\Lambda(A,Rf^!B)\to R\sHom_\Lambda(\tilde{g}^\ast A,Rf^{\prime!} g^\ast B)
\]
is an isomorphism.
\end{proposition}

\begin{proof} The assertion is local, so we may assume that $X$, $S$ and $S^\prime$ are spatial. By choosing a strictly totally disconnected cover $S^{\prime\prime}$ of $S^\prime$, one reduces the result for $S^\prime\to S$ to the cases of $S^{\prime\prime}\to S^\prime$ and $S^{\prime\prime}\to S$, so we may assume that $S^\prime$ is strictly totally disconnected. In that case, by \cite[Proposition 20.17]{ECoD}, whose hypothesis apply as $X'\to S'$ is of finite $\dimtrg$ and $S'$ is strictly totally disconnected, it suffices to check on global sections over all quasicompact separated \'etale maps $V^\prime\to X^\prime$. According to Lemma \ref{lemma:limite boules} we can write $S^\prime$ as a cofiltered limit of quasicompact open subsets $S_i^\prime$ of finite-dimensional balls over $S$. Then $V^\prime$ comes via pullback from a quasicompact separated \'etale map $V_i\to X\times_S S_i^\prime$ for $i$ large enough by \cite[Proposition 11.23]{ECoD}. We thus have a diagram with cartesian squares
$$
\begin{tikzcd}
V\ar[r]\ar[d] & V'\ar[d] \ar[rd,"h_i"] \\
X' \ar[r] \ar[d] & X\times_S S'_i \ar[r]\ar[d] & X\ar[d] \\
S' \ar[r] & S'_i\ar[r] & S.
\end{tikzcd}
$$
The result we want to prove is immediate when $S'\to S$ is cohomologically smooth. Up to replacing $X\to S$ by $V'\to S'_i$ we are thus reduced to prove that $R\Gamma (X',\tilde{g}^*\mathbb{D}_{X/S}(A))\iso R\Gamma(X',\mathbb{D}_{X'/S'}(\tilde{g}^*A))$.

Thus, it suffices to check the result after applying $Rf^\prime_\ast$. In that case,
\[\begin{aligned}
Rf^\prime_\ast \tilde{g}^\ast R\sHom_\Lambda(A,Rf^!B)&=g^\ast Rf_\ast R\sHom_\Lambda(A,Rf^!B)\\
&=g^\ast R\sHom_\Lambda(Rf_!A,B),
\end{aligned}\]
using \cite[Proposition 17.6, Theorem 1.8 (iv)]{ECoD}, using that $f$ satisfies $\dimtrg f<\infty$ and hence $Rf_\ast$ has finite cohomological dimension. On the other hand,
\[\begin{aligned}
Rf^\prime_\ast R\sHom_\Lambda(\tilde{g}^\ast A,Rf^{\prime!}B)&=R\sHom_\Lambda(Rf^\prime_!\tilde{g}^\ast A,B)\\
&=R\sHom_\Lambda(g^\ast Rf_!A,B)
\end{aligned}\]
using \cite[Theorem 1.8 (iv), Theorem 1.9 (ii)]{ECoD}. But by condition (b) of being $f$-locally acyclic, the complex $Rf_!A\in D_\et(S,\Lambda)$ is perfect-constructible, and thus the formation of $R\sHom_\Lambda(Rf_!A,B)$ commutes with any base change by Lemma \ref{lemma:commutation du dual au changement de base}.
\end{proof}

\begin{lemma}\label{lemma:limite boules}
Let $S$ be a spatial diamond and $X\to S$ be a morphism from an affinoid perfectoid space to $S$. Then one can write $X=\varprojlim_i U_i$ where $U_i$ is a quasicompact open subset inside a finite dimensional ball over $S$, and the projective limit is cofiltered.
\end{lemma}

\begin{proof}
If $I=\O(X)^+$, one has a closed immersion over $S$ defined by elements of $I$, 
$
X\hookrightarrow \mathbb{B}^I_S
$
where $\mathbb{B}^I_S$ is the spatial diamond over $S$ that represents the functor $T/S\mapsto (\O(T)^+)^I$ (an ``infinite dimensional perfectoid ball over $S$'' when $S$ is perfectoid). Now,
$$
\mathbb{B}^I_S =\varprojlim_{J\subset I} \mathbb{B}^J_S
$$
where $J$ goes through the set of finite subsets of $I$ and $\mathbb{B}^I_S\to \mathbb{B}^J_S$ is the corresponding projection. For each such $J$ the composite $X\hookrightarrow \mathbb{B}^I_S\twoheadrightarrow \mathbb{B}^J_S$ is a spatial morphism of spatial diamonds. Its image is a pro-constructible generalizing subset of $\mathbb{B}^J_S$ and can thus be written as $\bigcap_{\alpha\in A_J} U_\alpha$ where $U_\alpha$ is a quasicompact open subset of $\mathbb{B}^J_S$. Then one has 
\begin{align*}
X=\varprojlim_{J\subset I} \varprojlim_{\alpha\in A_J} U_\alpha. & \qedhere
\end{align*}
\end{proof}

\begin{lemma}\label{lemma:commutation du dual au changement de base}
Let $X$ be a spatial diamond and $A\in D_{\et} (X,\Lambda)$ perfect-constructible, and let $B\in D_\et(X,\Lambda)$ be arbitrary. Then the formation of $R\sHom_{\Lambda} (A,B)$ commutes with any base change.
\end{lemma}
\begin{proof}
Using \cite[Proposition 20.16 (iii)]{ECoD} this is reduced to the case when $A=j_! (\mathcal{L}|_{Z})$ where $j:U\to X$ is separated quasicompact \'etale, $Z\subset U$ closed constructible, and $\mathcal{L}\in D_{\et}(U,\Lambda)$ locally constant with perfect fibres. If $j':U\setminus Z\to X$, that is again quasicompact (since $Z$ is constructible inside $U$) separated \'etale, using the exact sequence $0\to j'_!\mathcal{L}\to j_!\mathcal{L}\to j_! (\mathcal{L}_{|Z})\to 0$, this is reduced to the case of $A$ of the form $j_!\mathcal{L}$. In this case $R\sHom_\Lambda(A,-)$ is given by $Rj_\ast(\mathcal L^\vee\dotimes_\Lambda j^\ast -)$, and this commutes with any base change by quasicompact base change, \cite[Proposition 17.6]{ECoD}.
\end{proof}

One has to be careful that, in general, {\it the naive dual of a perfect constructible complex is not constructible.} The following lemma says that in fact it is overconvergent, so never constructible unless locally constant.

\begin{lemma}\label{lemma:dual constructible surconvergent}
For $X$ a spatial diamond and $A\in D_{\et} (X,\Lambda)$ perfect constructible, $R \sHom_{\Lambda} (A,\Lambda)$ is overconvergent.
\end{lemma}

\begin{proof}
Using Lemma \ref{lemma:commutation du dual au changement de base} this is reduced to the case when $X=\Spa (C,C^+)$. Moreover, one can assume that $A=j_!\Lambda$ for some quasicompact open immersion $j: U\to X$. Then $R\sHom_\Lambda(A,\Lambda) = Rj_\ast \Lambda=\Lambda$, which is overconvergent.
\end{proof}

\subsubsection{Twisted inverse images}

Also, if $A$ is $f$-ULA, then one can relate appropriately $A$-twisted versions of $f^\ast$ and $Rf^!$.

\begin{proposition}\label{prop:ULAastshriekpullback} Let $f: X\to S$ be a compactifiable map of locally spatial diamonds with locally $\dimtrg f<\infty$ and let $A\in D_\et(X,\Lambda)$ be $f$-universally locally acyclic. Then for all $B\in D_\et(S,\Lambda)$, the natural map
\[
\mathbb D_{X/S}(A)\dotimes_\Lambda f^\ast B\to R\sHom_\Lambda(A,Rf^!B)
\]
given as the composite
\[
R\sHom_\Lambda(A,Rf^!\Lambda)\dotimes_\Lambda f^\ast B\to R\sHom_\Lambda(A,Rf^!\Lambda\dotimes_\Lambda f^\ast B)\to R\sHom_\Lambda(A,Rf^!B)
\]
is an isomorphism.
\end{proposition}

\begin{proof} First, we note that both sides commute with any base change, by Proposition~\ref{prop:ULAdual}.

It suffices to check that we get an isomorphism on stalks at all geometric points $\Spa(C,C^+)$ of $X$. For this, we may base change along the associated map $\Spa(C,C^+)\to S$ to reduce to the case that $S=\Spa(C,C^+)$ is strictly local, and we need to check that we get an isomorphism at the stalk of a section $s: S\to X$. We may also assume that $X$ is spatial, in which case $X$ is of bounded cohomological dimension, so \cite[Proposition 20.17]{ECoD} applies, and perfect-constructible complexes are the same thing as compact objects in $D_\et(X,\Lambda)=D(X_\et,\Lambda)$.

Next, we note that the functor $B\mapsto R\sHom_\Lambda(A,Rf^!B)$ is right adjoint to $A^\prime\mapsto Rf_!(A\dotimes_\Lambda A^\prime)$. The latter functor preserves perfect-constructible complexes, i.e.~compact objects, by condition (b), see Lemma \ref{lemma:ula produit tensoriel}. Thus, $B\mapsto R\sHom_\Lambda(A,Rf^!B)$ commutes with arbitrary direct sums, see Lemma \ref{lemma:adjunction compact direct sum}. Obviously, the functor $B\mapsto \mathbb D_{X/S}(A)\dotimes_\Lambda f^\ast B$ also commutes with arbitrary direct sums, so it follows that it suffices to check the assertion for $B=j_!\Lambda$ for some quasicompact open immersion $j: S^\prime=\Spa(C,C^{\prime+})\to S$ (the shifts of those compact objects generate $D_{\et} (S,\Lambda)$). If $S^\prime=S$, then $B=\Lambda$ and the result is clear. Otherwise, the stalk of $j_!\Lambda$ at the closed point is zero, and thus the stalk of the left-hand side $\mathbb D_{X/S}(A)\dotimes_\Lambda f^\ast B$ at our fixed section is zero. It remains to see that the stalk of
\[
R\sHom_\Lambda(A,Rf^! j_!\Lambda)
\]
at the (closed point of) the section $s: S\to X$ is zero. This stalk is given by the filtered colimit over all quasicompact open neighborhoods $U\subset X$ of $s(S)$ of
\[
R\Hom_{D_\et(U,\Lambda)}(A|_U,Rf^! j_!\Lambda|_U) = R\Hom_{D_\et(S,\Lambda)}(Rf_{U!} A|_U,j_!\Lambda),
\]
where $f_U: U\to S$ denotes the restriction of $f$ (the possibility to use only open embeddings in place of general \'etale maps results from the observation that the intersection of all these open subsets is the strictly local space $S$ already;  the set of open neighborhoods of $s(S)$ is cofinal among \'etale neighborhoods of $s(S)$).

Now we claim that the inverse systems of all such $U$ and of the compactifications $\overline{U}^{/S}$ are cofinal. Note that the intersection of all $\overline{U}^{/S}$ (taken inside $\overline{X}^{/S}$) is simply $s(S)$: Indeed, given any point $x\in \overline{X}^{/S}\setminus s(S)$, there are disjoint open neighborhoods $x\in V$ and $s(S)\subset U$. In fact, the maximal Hausdorff quotient $|\overline{X}^{/S}|^B$ is compact Hausdorff by \cite[Proposition 13.11]{ECoD} and its points can be identified with rank $1$ points of $\overline{X}^{/S}$, which are the same as rank $1$ points of $X$. But as $s(S)\subset \overline{X}^{/S}$ is closed, no point outside $s(S)$ admits the same rank $1$ generalization, so $x$ and $s(S)$ define distinct points of the Hausdorff spaces $|\overline{X}^{/S}|$, so that the desired disjoint open neighborhoods $x\in V$ and $s(S)\subset U$ exist. Then $x\not\in \overline{U}^{/S}$. Thus,
$$
s(S)=\bigcap_{U\supset s(S)} \overline{U}^{/X}.
$$
 Now given any open neighborhood $U$ of $s(S)$, the complement $|\overline{X}^{/S}|\setminus U$ is quasicompact, which then implies that there is some $U^\prime$ such that $\overline{U^\prime}^{/S}\subset U$. It follows that the direct systems of
\[
Rf_{U!} A|_U
\]
and
\[
Rf_{\overline{U}^{/S}\ast} Ri_{\overline{U}^{/S}}^! A
\]
are equivalent, where $i_{\overline{U}^{/S}}: \overline{U}^{/S}\to X$ and $f_{\overline{U}^{/S}}: \overline{U}^{/S}\to S$ are the evident maps. Now observe that if $j^\prime_X: X_\eta = X\times_S \Spa(C,\OO_C)\hookrightarrow X$ denotes the proconstructible generalizing immersion, then condition (a) in being $f$-locally acyclic implies that $A=Rj^\prime_{X\ast} A|_{X_\eta}$ (see Remark \ref{remark:autre definition locale acyclicite}), and then
\[\begin{aligned}
Rf_{\overline{U}^{/S}\ast} Ri_{\overline{U}^{/S}}^! A &= Rf_{\overline{U}^{/S}\ast} Ri_{\overline{U}^{/S}}^! Rj^\prime_{X\ast} A|_{X_\eta}\\
&=Rf_{\overline{U}^{/S}\ast} Rj^\prime_{\overline{U}^{/S}\ast} Ri_{\overline{U}^{/S}_\eta}^! A|_{X_\eta}\\
&=Rj^\prime_{S\ast} Rf_{\overline{U}^{/S}_\eta\ast} Ri_{\overline{U}^{/S}_\eta}^! A|_{X_\eta}
\end{aligned}\]
with hopefully evident notation; in particular, $j^\prime_S: \Spa(C,\OO_C)\to S=\Spa(C,C^+)$ denotes the pro-open immersion of the generic point on the base.

In summary, we can rewrite the stalk of $R\sHom_\Lambda(A,Rf^! j_!\Lambda)$ at $s(S)$ as the filtered colimit of
\[
R\Hom_{D_\et(S,\Lambda)}(Rj^\prime_{S\ast} Rf_{\overline{U}^{/S}_\eta\ast} Ri_{\overline{U}^{/S}_\eta}^! A,j_!\Lambda),
\]
and we need to prove that this vanishes. This follows from the observation that for all $M\in D_\et(\Spa(C,\OO_C),\Lambda)=D(\Lambda)$, one has
\[
R\Hom_{D_\et(S,\Lambda)}(Rj^\prime_{S\ast} M,j_!\Lambda)=0.
\]
For this, note that $Rj^\prime_{S\ast} M=M$ is just the constant sheaf given by the complex of $\Lambda$-modules $M$, and one has a triangle
\[
R\Hom_{D_\et(S,\Lambda)}(M,j_!\Lambda)\to R\Hom_{D_\et(S,\Lambda)}(M,\Lambda)\to R\Hom_{D_\et(S,\Lambda)}(M,i_\ast\Lambda),
\]
where $i$ denotes the complementary closed immersion. Both the second and last term are given by $R\Hom_\Lambda(M,\Lambda)$, finishing the proof.
\end{proof}

We used the following classical lemma, cf.~\cite[Theorem 5.1]{Neeman}.

\begin{lemma}\label{lemma:adjunction compact direct sum}
Let $\mathcal{C}$ and $\mathcal{D}$ be triangulated categories such that $\mathcal{C}$ is compactly generated. Let $F:\mathcal{C}\to \mathcal{D}$ and $G:\mathcal{D}\to \mathcal{C}$ be such that $G$ is right adjoint to $F$. If $F$ sends compact objects to compact objects then $G$ commutes with arbitrary direct sums.
\end{lemma}
\begin{proof}
Since $\mathcal{C}$ is compactly generated it suffices to prove that for any compact object $A$ in $\mathcal{C}$ and any collection $(B_i)_i$ of objects of $\mathcal{D}$, $\Hom (A,\oplus_i G(B_i))\iso \Hom (A, G(\oplus_i B_i))$. By compactness of $A$, $\Hom (A,\oplus_i G(B_i))=\oplus_i \Hom (A,G(B_i))$, by adjunction this is equal to $\oplus_i \Hom (F(A),B_i)$,  since $F(A)$ is compact this is equal to $\Hom (F(A),\oplus_i B_i)$, and by adjunction this is $\Hom (A,G(\oplus_i B_i))$.
\end{proof}

\begin{remark}\label{rem:locale acyclicite caracterise via twist}\leavevmode
\begin{altenumerate}
\item In fact, outside of the overconvergence condition (a) in Definition \ref{def:ULA}, {\it the property of Proposition~\ref{prop:ULAastshriekpullback} characterizes locally acyclic complexes} under the assumption that locally on $X$ the exists an integer $N$ such that for any $U\to X$ quasicompact separated \'etale the cohomological dimension of $U_{\et}$ is bounded by $N$. More precisely, if $j:U\to X$ is separated \'etale with $f\circ j$ quasicompact then $R\Hom_{\Lambda} ( R(f\circ j)_!A, B)=  R\Gamma(U,R\sHom_{\Lambda} (A,Rf^! B))$. Thus, if for all $B$ one has $\mathbb D_{X/S}(A)\dotimes_\Lambda f^\ast B\iso R\sHom_\Lambda(A,Rf^!B)$ then $R(f\circ j)_!A$ is compact since $R\Gamma(U,-)$ commutes with arbitrary direct sums.
\item Outside of the overconvergence condition (a) in Definition \ref{def:ULA}, the property of Proposition~\ref{prop:ULAastshriekpullback} universally on $S$ characterize universally locally acyclic objects. In fact, using \cite[Proposition 20.13]{ECoD}, the constructibility property is reduced to the case when the base is strictly totally disconnected, in which case we can apply point (i).
\end{altenumerate}
\end{remark}

Finally, let us note that all the previous results extend to the setting where the base $S$ is a general small v-stack, taking the following definition.

\begin{definition}\label{def:ULAgeneralbase} Let $f: X\to S$ be a map of small v-stacks that is compactifiable and representable in locally spatial diamonds with locally $\dimtrg f<\infty$. Let $A\in D_\et(X,\Lambda)$. Then $A$ is $f$-universally locally acyclic if for any map $S^\prime\to S$ from a locally spatial diamond $S^\prime$ with pullback $f^\prime: X^\prime=X\times_S S^\prime\to S^\prime$, the complex $A|_{X^\prime}\in D_\et(X^\prime,\Lambda)$ is $f^\prime$-locally acyclic.
\end{definition}

\subsubsection{Dualizability}

From the previous two propositions, we can deduce an analogue of a recent result of Lu-Zheng, \cite{LuZhengDuality}, characterizing universal local acyclicity in terms of dualizability in a certain monoidal category. We actually propose a different such characterization closer to how dualizability will appear later in the discussion of geometric Satake. In terms of applications to abstract properties of universal local acyclicity, such as its preservation by Verdier duality, it leads to the same results.

Fix a base small v-stack $S$, and a coefficient ring $\Lambda$ (killed by some integer prime to $p$). We define a $2$-category $\mathcal C_S$ as follows. The objects of $\mathcal C_S$ are maps $f: X\to S$ of small v-stacks that are compactifiable, representable in locally spatial diamonds, with locally $\dimtrg f<\infty$. For any $X,Y\in \mathcal C_S$, the category of maps $\Fun_{\mathcal C_S}(X,Y)$ is the category $D_\et(X\times_S Y,\Lambda)$. Note that any such $A\in D_\et(X\times_S Y,\Lambda)$ defines in particular a functor
\[
D_\et(X,\Lambda)\to D_\et(Y,\Lambda): B\mapsto R\pi_{2!}(A\dotimes_\Lambda \pi_1^\ast B)
\]
with kernel $A$, where $\pi_1: X\times_S Y\to X$, $\pi_2: X\times_S Y\to Y$ are the two projections. The composition in $\mathcal C_S$ is now defined to be compatible with this association. More precisely, the composition
\[
\Fun_{\mathcal C_S}(X,Y)\times \Fun_{\mathcal C_S}(Y,Z)\to \Fun_{\mathcal C_S}(X,Z)
\]
is defined to be the functor
\[
D_\et(X\times_S Y,\Lambda)\times D_\et(Y\times_S Z,\Lambda)\to D_\et(X\times_S Z,\Lambda): (A,B)\mapsto A\star B = R\pi_{13!}(\pi_{12}^\ast A\dotimes_\Lambda \pi_{23}^\ast B)
\]
where $\pi_{ij}$ denotes the various projections on $X\times_S Y\times_S Z$. It follows from the projection formula that this indeed defines a $2$-category $\mathcal C_S$. The identity morphism is given by $R\Delta_! \Lambda=R\Delta_\ast \Lambda\in D_\et(X\times_S X,\Lambda)$, where $\Delta: X\hookrightarrow X\times_S X$ is the diagonal (which is a closed immersion, as the morphism $X\to S$ is assumed to be compactifiable, in particular $0$-truncated and separated). We note that $\mathcal C_S$ is naturally equivalent to $\mathcal C_S^\op$. Indeed, $D_\et(X\times_S Y,\Lambda)$ is invariant under switching $X$ and $Y$, and the definition of composition (and coherences) is compatible with this switch.

Recall that in any $2$-category $\mathcal C$, there is a notion of adjoints. Namely, a morphism $f: X\to Y$ is a left adjoint of $g: Y\to X$ if there are maps $\alpha: \mathrm{id}_X\to gf$ and $\beta: fg\to \mathrm{id}_Y$ such that the composites
\[
f\xrightarrow{f\alpha} fgf\xrightarrow{\beta f} f,\ g\xrightarrow{\alpha g} gfg\xrightarrow{g\beta} g
\]
are the identity. If a right adjoint $g$ of $f$ exists, it is (together with the accompanying data) moreover unique up to unique isomorphism. As is clear from the definition, any functor of $2$-categories preserves adjunctions. In particular, this applies to pullback functors $\mathcal C_S\to \mathcal C_{S'}$ for maps $S'\to S$ of small v-stacks, or to the functor from $\mathcal C_S$ to triangulated categories taking $X$ to $D_\et(X,\Lambda)$ and $A\in \Fun_{\mathcal C_S}(X,Y)$ to the functor $R\pi_{2!}(A\dotimes_\Lambda \pi_1^\ast -)$ with kernel $A$.

\begin{theorem}\label{thm:ULAdualizable} Let $S$ be a small v-stack and $X\in \mathcal C_S$, and $A\in D_\et(X,\Lambda)$. The following conditions are equivalent.
\begin{altenumerate}
\item[{\rm (i)}] The complex $A$ is $f$-universally locally acyclic.
\item[{\rm (ii)}] The natural map
\[
p_1^\ast \mathbb D_{X/S}(A)\dotimes_\Lambda p_2^\ast A\to R\sHom_\Lambda(p_1^\ast A,Rp_2^! A)
\]
is an isomorphism, where $p_1,p_2:X\times_S X\to X$ are the two projections.
\item[{\rm (iii)}] The object $A\in \Fun_{\mathcal C_S}(X,S)$ is a left adjoint in $\mathcal C_S$. In that case, its right adjoint is given by
\[
\mathbb D_{X/S}(A)\in D_\et(X,\Lambda)=\Fun_{\mathcal C_S}(S,X).
\]
\end{altenumerate}
\end{theorem}

\begin{proof} That (i) implies (ii) follows from Proposition~\ref{prop:ULAdual} and Proposition~\ref{prop:ULAastshriekpullback}. For (ii) implies (iii), we claim that $A\in \Fun_{\mathcal C_S}(X,S)$ is indeed a left adjoint of $\mathbb D_{X/S}(A)\in \Fun_{\mathcal C_S}(A)$. The composites are given by
\[
A\star \mathbb D_{X/S}(A) = Rf_!(\mathbb D_{X/S}(A)\dotimes_\Lambda A)\in D_\et(S,\Lambda)=\Fun_{\mathcal C_S}(S,S)
\]
and
\[
\mathbb D_{X/S}(A)\star A = p_1^\ast A\dotimes_\Lambda p_2^\ast \mathbb D_{X/S}(A)\in D_\et(X\times_S X,\Lambda)=\Fun_{\mathcal C_S}(X,X).
\]
Then we take $\beta: A\star \mathbb D_{X/S}(A)\to \mathrm{id}_S$ to be given by the map $Rf_!(\mathbb D_{X/S}(A)\dotimes_\Lambda A)\to \Lambda$ adjoint to the map $\mathbb D_{X/S}(A)\dotimes_\Lambda A\to Rf^!\Lambda$ which is just the tautological pairing. On the other hand, for $\alpha: \mathrm{id}_X\to \mathbb D_{X/S}(A)\star A$, we have to produce a map
\[
R\Delta_! \Lambda\to p_1^\ast A\dotimes_\Lambda p_2^\ast \mathbb D_{X/S}(A).
\]
Using (ii), the right-hand side is naturally isomorphic to $R\sHom_\Lambda(p_1^\ast A,Rp_2^! A)$. Now maps from $R\Delta_!\Lambda$ are adjoint to sections of
\[
R\Delta^! R\sHom_\Lambda(p_1^\ast A,Rp_2^!A)\cong R\sHom_\Lambda(A,A)
\]
(using \cite[Theorem 1.8 (v)]{ECoD}), which has the natural identity section. It remains to prove that certain composites are the identity. This follows from a straightforward diagram chase.

Finally, it remains to prove that (iii) implies (i). We can assume that $S$ is strictly totally disconnected. It follows that the functor $Rf_!(A\dotimes_\Lambda -)$ admits a right adjoint that commutes with all colimits. This implies that condition (b) in Definition~\ref{def:ULA} is satisfied. In fact, more precisely we see that the right adjoint $R\sHom_\Lambda(A,Rf^!-)$ is given by $A'\dotimes_\Lambda f^\ast-$ for some $A'\in D_\et(X,\Lambda)$, and by using the self-duality of $\mathcal C_S^{\mathrm{op}}$, we also see that $R\sHom_\Lambda(A',Rf^!-)$ is given by $A\dotimes_\Lambda f^\ast-$. Applied to the constant sheaf, this shows in particular that $A\cong R\sHom_\Lambda(A',Rf^!\Lambda)$ is a Verdier dual. For condition (a), we can assume that $S=\Spa(C,C^+)$ and reduce to checking overconvergence along sections $s: S\to X$. In fact, using part (2) of Remark~\ref{remark:autre definition locale acyclicite}, let $j: \Spa(C,\OO_C)\to \Spa(C,C^+)$ be the pro-open immersion, with pullback $j_X: X\times_{\Spa(C,C^+)} \Spa(C,\OO_C)\to X$, and $f_\eta: X\times_{\Spa(C,C^+)}\Spa(C,\OO_C)\to \Spa(C,\OO_C)$ the restriction of $f$. To see the overconvergence, it is enough to see that $A\cong Rj_{X\ast} A_0$ for some $A_0$. But
\[\begin{aligned}
A&\cong R\sHom_\Lambda(A',Rf^!\Lambda)\cong R\sHom_\Lambda(A',Rf^!Rj_\ast\Lambda)\\
&\cong R\sHom_\Lambda(A',Rj_{X\ast}Rf_\eta^!\Lambda)\cong Rj_{X\ast}R\sHom_\Lambda(j_X^\ast A',Rf_\eta^!\Lambda),
\end{aligned}\]
giving the desired overconvergence.
\end{proof}

Before moving on, let us observe the following relative variant.

\begin{proposition}\label{prop:generalULAdualizable} Let $S$ be a small v-stack and $X,Y\in \mathcal C_S$. If $Y/S$ is proper and $A\in \Fun_{\mathcal C_S}(X,Y)=D_\et(X\times_S Y,\Lambda)$ is $p_2$-universally locally acyclic, then it is a left adjoint and the right adjoint is given by
\[
\mathbb D_{X\times_S Y/Y}(A)\in D_\et(X\times_S Y,\Lambda)\cong D_\et(Y\times_S X,\Lambda)=\Fun_{\mathcal C_S}(Y,X).
\]
\end{proposition}

The assumption that $Y/S$ is proper is important here. Already if $X=S$ and $A=\Lambda\in D_\et(Y,\Lambda)$, which is always $\mathrm{id}_Y$-universally locally acyclic, being a left adjoint in $\mathcal C_S$ implies that there is some $B\in D_\et(Y,\Lambda)$ for which $Rf_\ast\cong Rf_!(B\dotimes_\Lambda-)$.

\begin{proof} We need to produce the maps $\alpha$ and $\beta$ again. Let us give the construction of $\alpha$, which is the harder part. First, using the various projections $\pi_{ij}$ on $X\times_S Y\times_S X$, we have
\[
\mathbb D_{X\times_S Y/Y}(A)\star A\cong R\pi_{13!}(\pi_{12}^\ast \mathbb D_{X\times_S Y/Y}(A)\dotimes_\Lambda \pi_{23}^\ast A)\cong R\pi_{13\ast}R\sHom_\Lambda(\pi_{12}^\ast A,R\pi_{23}^! A)
\]
using that $A$ is $p_2$-universally locally acyclic, and properness of $\pi_{13}$ (which is a base change of $Y\to S$). Now giving a map $R\Delta_! \Lambda\to \mathbb D_{X\times_S Y/Y}(A)\star A$, for $\Delta=\Delta_{X/S}$, amounts to finding a section of
\[
R\Delta^! R\pi_{13\ast} R\sHom_\Lambda(\pi_{12}^\ast A,R\pi_{23}^! A)\cong Rp_{1\ast} R\Delta_{X\times_S Y/Y}^!R\sHom_\Lambda(\pi_{12}^\ast A,R\pi_{23}^! A)\cong Rp_{1\ast} R\sHom_\Lambda(A,A),
\]
where we can take the identity.
\end{proof}

Theorem~\ref{thm:ULAdualizable} has the following notable consequences.

\begin{corollary}\label{cor:ULAselfdual} Let $f: X\to S$ be a compactifiable map of locally spatial diamonds with locally $\dimtrg f<\infty$ and let $A\in D_\et(X,\Lambda)$ be $f$-universally locally acyclic. Then $\mathbb D_{X/S}(A)$ is again $f$-universally locally acyclic, and the biduality map
\[
A\to \mathbb D_{X/S}(\mathbb D_{X/S}(A))
\]
is an isomorphism.

If $f_i: X_i\to S$ for $i=1,2$ are compactifiable maps of small v-stacks that are representable in locally spatial diamonds with locally $\dimtrg f_i<\infty$ and $A_i\in D_\et(X_i,\Lambda)$ are $f_i$-universally locally acyclic, then also $A_1\boxtimes A_2\in D_\et(X_1\times_S X_2,\Lambda)$ is $f_1\times_S f_2$-universally locally acyclic, and the natural map
\[
\mathbb D_{X_1/S}(A_1)\boxtimes \mathbb D_{X_2/S}(A_2)\to \mathbb D_{X_1\times_S X_2/S}(A_1\boxtimes A_2)
\]
is an isomorphism.
\end{corollary}

\begin{proof} By Theorem~\ref{thm:ULAdualizable}, the object $\mathbb D_{X/S}(A)\in \Fun_{\mathcal C_S}(S,X)$ is a right adjoint of $A\in \Fun_{\mathcal C_S}(X,S)$. But $\mathcal C_S\cong \mathcal C_S^\op$; under this equivalence, this means that $\mathbb D_{X/S}(A)\in \Fun_{\mathcal C_S}(S,X)$ is a left adjoint of $A\in \Fun_{\mathcal C_S}(S,X)$. Thus, applying Theorem~\ref{thm:ULAdualizable} again, the result follows.

For the second statement, note that $A_i\in D_\et(X_i,\Lambda)$ define left adjoints, hence so does
\[
A_1\star A_2=A_1\boxtimes A_2\in D_\et(X_1\times_S X_2,\Lambda)=\Fun_{\mathcal C_S}(X_1\times_S X_2,S).
\]
Its right adjoint is the similar composition, giving the claim.
\end{proof}

The final statement admits the following generalization concerning ``compositions'' of universally locally acyclic sheaves.

\begin{proposition}\label{prop:composeULA} Let $g: Y\to X$ and $f: X\to S$ be compactifiable maps of small v-stacks representable in locally spatial diamonds with locally $\dimtrg f,\dimtrg g<\infty$, and let $A\in D_\et(X,\Lambda)$ be $f$-universally acyclic and $B\in D_\et(Y,\Lambda)$ be $g$-universally locally acyclic. Then $g^\ast A\dotimes_\Lambda B$ is $f\circ g$-universally locally acyclic, and there is natural isomorphism
\[
\mathbb D_{Y/S}(g^\ast A\dotimes_\Lambda B)\cong g^\ast \mathbb D_{X/S}(A)\dotimes_\Lambda \mathbb D_{Y/X}(B).
\]
\end{proposition}

\begin{proof} It is easy to see that condition (a) of being locally acyclic holds universally, so it suffices to identify the functor $R\sHom_\Lambda(g^\ast A\dotimes_\Lambda B,R(f\circ g)^!-)$. We compute:
\[\begin{aligned}
R\sHom_\Lambda(g^\ast A\dotimes_\Lambda B,R(f\circ g)^!-) &= R\sHom_\Lambda(B,R\sHom_\Lambda(g^\ast A,Rg^! Rf^! -))\\
&=R\sHom_\Lambda(B,Rg^!R\sHom_\Lambda(A,Rf^!-))\\
&=\mathbb D_{Y/X}(B)\dotimes_\Lambda g^\ast R\sHom(A,Rf^!-)\\
&=\mathbb D_{Y/X}(B)\dotimes_\Lambda g^\ast \mathbb D_{X/S}(A)\dotimes_\Lambda g^\ast f^\ast-,
\end{aligned}\]
implying that it commutes with colimits, hence its left adjoint preserves perfect-constructible complexes (after reduction to strictly totally disconnected $S$ and $X$ and $Y$ spatial). Moreover, evaluating this functor at $\Lambda$ gives the identification of the Verdier dual.
\end{proof}

Let us also note another corollary of Theorem~\ref{thm:ULAdualizable} concerning retracts.

\begin{corollary}\label{cor:ULAretract} Let $f: X\to S$ and $g: Y\to S$ be maps of small v-stacks that are compactifiable and representable in locally spatial diamonds with locally $\dimtrg f,\dimtrg g<\infty$. Assume that $f$ is a retract of $g$ over $S$, i.e.~there are maps $i: X\to Y$, $r: Y\to X$ over $S$ such that $ri=\mathrm{id}_X$. If $\Lambda$ is $g$-universally locally acyclic, then $\Lambda$ is $f$-universally locally acyclic.
\end{corollary}

\begin{proof} One can check this directly from the definitions, or note that the map in $\mathcal C_S$ given by $\Lambda\in D_\et(X,\Lambda)=\Fun_{\mathcal C_S}(X,S)$ is a retract of the map given by $\Lambda\in D_\et(Y,\Lambda)=\Fun_{\mathcal C_S}(Y,S)$, from which one can easily obtain adjointness.
\end{proof}

Moreover, in some cases the converse to Proposition~\ref{prop:properULA} holds.

\begin{proposition}\label{prop:properproetaleULA} Let $g: Y\to X$, $f: X\to S$ be maps of locally spatial diamonds where $g$ is proper and quasi-pro-\'etale and $f$ is compactifiable and locally $\dimtrg f<\infty$. Then $A\in D_\et(Y,\Lambda)$ is $f\circ g$-universally locally acyclic if and only if $Rg_\ast A\in D_\et(X,\Lambda)$ is $f$-universally locally acyclic.
\end{proposition}

\begin{proof} One direction is given by Proposition~\ref{prop:properULA}. For the converse, assume that $Rg_\ast A$ is $f$-universally locally acylic. To see that $A$ is $h=f\circ g$-universally locally acyclic, it suffices by Theorem~\ref{thm:ULAdualizable} that the map
\[
p_{1,Y}^\ast R\sHom(A,Rh^!\Lambda)\dotimes_\Lambda p_{2,Y}^\ast A\to R\sHom(p_{1,Y}^\ast A,Rp_{2,Y}^! A)
\]
in $D_\et(Y\times_S Y,\Lambda)$ is an isomorphism, where $p_{1,Y},p_{2,Y}: Y\times_S Y\to Y$ are the two projections. As $g\times_S g: Y\times_S Y\to X\times_S X$ is proper and quasi-pro-\'etale, pushforward along $g\times_S g$ is conservative: By testing on stalks, this follows from the observation that for a profinite set $T$, the global sections functor $R\Gamma(T,-)$ is conservative on $D(T,\Lambda)$ (as one can write any stalk as a filtered colimit of functors that are direct summands of the global sections functor). Applying $R(g\times_S g)_\ast=R(g\times_S g)_!$ to the displayed map, we get the map
\[
p_{1,X}^\ast R\sHom(Rg_\ast A,Rf^!\Lambda)\dotimes_\Lambda p_{2,X}^\ast Rg_\ast A\to R\sHom(p_{1,X}^\ast Rg_\ast A,Rp_{2,X}^! Rg_\ast A)
\]
where $p_{1,X},p_{2,X}: X\times_S X\to X$ are the two projections. This is an isomorphism precisely when $Rg_\ast A$ is $f$-universally locally acyclic.
\end{proof}

The following corollary shows that smooth base change generalizes to universally locally acyclic maps. The general version of this corollary was suggested by David Hansen.

\begin{corollary}[ULA base change]\label{cor:ULAbasechange} Consider a cartesian diagram of small v-stacks
\[\xymatrix{
X'  \ar[r]^{\tilde{g}}\ar[d]^{f^\prime} & X \ar[d]^f \\
S'\ar[r]^g & S
}\]
with $f$ representable in locally spatial diamonds, compactifiable, locally $\dimtrg f<\infty$. Assume that $\Lambda$ is $f$-universally locally acyclic. Then the base change map
\[
f^\ast Rg_\ast A\to R\tilde{g}_\ast f^{\prime\ast} A
\]
is an isomorphism. More generally, if $B\in D_\et(X,\Lambda)$ is $f$-universally locally acyclic and $A\in D_\et(S',\Lambda)$, then
\[
(f^\ast Rg_\ast A)\dotimes_\Lambda B\iso R\tilde{g}_\ast(f^{\prime\ast} A\dotimes_\Lambda \tilde{g}^\ast B).
\]
\end{corollary}

Such base change results are false without some hypothesis on $f$. For example, if $S'$ is a countable union of copies of $S$ and $X$ is a geometric point of $S$, this base change would assert that taking stalks commutes with (countable) products, which fails in general.

\begin{proof} We apply Proposition~\ref{prop:ULAastshriekpullback} to the universally locally acyclic $\mathbb D_{X/S}(B)$, so that by Corollary~\ref{cor:ULAselfdual}, we get
\[
f^\ast Rg_\ast A\dotimes_\Lambda B\cong R\sHom_\Lambda(\mathbb D_{X/S}(B),Rf^!Rg_\ast A).
\]
By \cite[Theorem 1.9 (iii)]{ECoD}, $Rf^! Rg_\ast A\cong R\tilde{g}_\ast Rf^{\prime!} A$, and then one can rewrite further as
\[
R\sHom_\Lambda(\mathbb D_{X/S}(B),R\tilde{g}_\ast Rf^{\prime!} A)\cong R\tilde{g}_\ast R\sHom_\Lambda(\tilde{g}^\ast \mathbb D_{X/S}(B),Rf^{\prime!}A).
\]
Now another application of Proposition~\ref{prop:ULAastshriekpullback} and Proposition~\ref{prop:ULAdual} gives the result.
\end{proof}

Finally let us note the following consequence of Theorem \ref{thm:ULAdualizable} and \cite{LuZhengDuality}.

\begin{proposition}
Let $K$ be a complete non-archimedean field with residue characteristic $p$, $f:X\to S$ a separated morphism of $K$-schemes locally of finite type, and $A\in D^b_c(X,\Lambda)$. 
Then $A$ is $f$-universally acyclic if and only if its analytification $A^{ad}$ is $f^{ad,\diamond}$-universally locally acyclic, where $f^{ad,\diamond}:X^{ad,\diamond}\to S^{ad,\diamond}$. 
\end{proposition}

\begin{proof} The criterion of Theorem~\ref{thm:ULAdualizable}~(ii) applies similarly in the algebraic case by \cite{LuZhengDuality}, and all operations are compatible with passing to analytic adic spaces (and diamonds).\footnote{For this compatibility, the case of $\otimes$, $f^\ast$ and $f_!$ is easy. The essential remaining case is $Rj_\ast$ for an open immersion $j$ (both internal Hom and $i^!$ for closed immersions reduce to that case, and general $f^!$ can be reduced to closed immersions and smooth maps, where the latter reduces to $f^\ast$). The case of $Rj_\ast$ is \cite[Theorem 3.8.1]{Huber}.}
\end{proof}

For example, if $S=\Spec K$ then any $A$ is $f$-universally locally acyclic and thus $A^{ad}$ is $f^{ad,\diamond}$-universally acyclic. This gives plenty of examples of ULA sheaves.

\subsection{Local acyclicity for morphisms of Artin v-stacks}

Using the descent results Remark~\ref{rem:defULA} and Proposition~\ref{prop:ULAsmoothlocal}, one can extend the previous definition and results to the case of maps of Artin v-stacks as follows.

\begin{definition}\label{def:local acyclicity for Artin stacks} Let $f: X\to S$ be a map of Artin v-stacks and assume that there is some separated, representable in locally spatial diamonds, and cohomologically smooth surjection $g: U\to X$ from a locally spatial diamond $U$ such that $f\circ g: U\to S$ is compactifiable with locally $\dimtrg(f\circ g)<\infty$. Then $A\in D_\et(X,\Lambda)$ is $f$-universally locally acyclic if $g^\ast A$ is $f\circ g$-universally locally acyclic.
\end{definition}

All previous results concerning universally locally acyclic complexes also hold in this setting (assuming that the relevant operations are defined in the case of interest -- we did not define $Rf_!$ and $Rf^!$ for general stacky maps), and follow by the reduction to the case when $S$ and $X$ are locally spatial diamonds. In particular, the characterization in terms of dualizability gives the following.

\begin{proposition}\label{prop:ULAdualizablestack} Let $f: X\to S$ be a cohomologically smooth map of Artin v-stacks, and let $A\in D_\et(X,\Lambda)$. Consider $X\times_S X$ with its two projections $p_1,p_2: X\times_S X\to X$. Then $A$ is $f$-universally locally acyclic if and only if the natural map
\[
p_1^\ast R\sHom_\Lambda(A,\Lambda)\dotimes_\Lambda p_2^\ast A\to R\sHom_\Lambda(p_1^\ast A,p_2^\ast A)
\]
is an isomorphism.
\end{proposition}

\begin{proof} Taking a chart for $S$, we can assume that $S$ is a locally spatial diamond, and then taking a presentation for $X$ we can assume that also $X$ is a locally spatial diamond, noting that the condition commutes with smooth base change. In that case, replacing some occurences of $p_2^\ast$ by $Rp_2^!$ using cohomological smoothness, the result follows from Theorem~\ref{thm:ULAdualizable}.
\end{proof}

There is a simple characterization of $\ell$-cohomological smoothness in terms of universal local acyclicity.

\begin{proposition}\label{prop:ULAcohomsmooth} Let $f: X\to S$ be a compactifiable map of v-stacks that is representable in locally spatial diamonds with locally $\dimtrg f<\infty$. Then $f$ is $\ell$-cohomologically smooth if and only if $\mathbb F_\ell$ is $f$-universally locally acyclic and its Verdier dual $Rf^!\mathbb F_\ell$ is invertible.
\end{proposition}

Note that in checking whether $\mathbb F_\ell$ is $f$-universally locally acyclic, condition (a) of Definition~\ref{def:ULA} is automatic. Also, by Theorem~\ref{thm:ULAdualizable}, the condition that $\mathbb F_\ell$ is $f$-universally locally acyclic is equivalent to the condition that the natural map
\[
p_1^\ast Rf^! \mathbb F_\ell\to Rp_2^! \mathbb F_\ell
\]
is an isomorphism, where $p_1,p_2: X\times_S X\to X$ are the two projections. Thus, {\it $f$ is $\ell$-cohomologically smooth if and only if $Rf^! \mathbb F_\ell$ is invertible and its formation commutes with any base change}.

\begin{proof} The conditions are clearly necessary. For the converse, we may assume that $S$ is strictly totally disconnected. By Proposition~\ref{prop:ULAastshriekpullback}, the natural transformation of functors
\[
Rf^!\Lambda\dotimes_\Lambda f^\ast\to Rf^!
\]
is an equivalence. As $Rf^!\Lambda$ is assumed to be invertible (it commutes with base change by Proposition~\ref{prop:ULAdual}), this shows that the condition of \cite[Definition 23.8]{ECoD} is satisfied.
\end{proof}

In particular, we can resolve a question from \cite{ECoD}.

\begin{corollary}\label{cor:Zpsmooth} The map $f: \Spd \mathcal O_E\to \Spd \mathbb F_q$ is $\ell$-cohomologically smooth for all $\ell\neq p$.
\end{corollary}

\begin{proof} This is clear if $E$ is of equal characteristic, so assume that $E$ is $p$-adic. First, we prove that $\mathbb F_\ell$ is $f$-ULA. This follows from $f^\prime: \Spd \mathcal O_{\tilde{E}}\cong \Spd \mathbb F_q\powerseries{t^{1/p^\infty}}\to \Spd \mathbb F_q$ being $\ell$-cohomologically smooth, where $\tilde{E}/E$ is some totally ramified $\mathbb Z_p$-extension, by the argument of the proof of \cite[Proposition 24.3]{ECoD} (in essence, the compactly supported pushforward for any base change of $f$ are the $\mathbb Z_p$-invariants inside the compactly supported pushforward for the corresponding base change of $f^\prime$, so constructibility of the latter implies constructibility of the former).

It remains to show that $Rf^!\mathbb F_\ell$ is invertible. If $j: \Spd E\to \Spd \mathcal O_E$ is the open immersion with complement $i: \Spd \mathbb F_q\to \Spd \mathcal O_E$, we have $Ri^! Rf^! \mathbb F_\ell=\mathbb F_\ell$ by transitivity, and $j^\ast Rf^!\mathbb F_\ell\cong \mathbb F_\ell(1)[2]$ by \cite[Proposition 24.5]{ECoD}, so we get a triangle
\[
i_\ast \mathbb F_\ell\to Rf^!\mathbb F_\ell\to Rj_\ast \mathbb F_\ell(1)[2].
\]
Using this, one computes $Rf^!\mathbb F_\ell\cong \mathbb F_\ell(1)[2]$, as desired.
\end{proof}

\section{Formal smoothness}\label{sec:formalsmooth}
\subsection{Definition}

A key step in the proof of Theorem~\ref{thm:jacobiancriterion}, the Jacobian criterion of cohomological smoothness, is the following notion of formal smoothness.

\begin{definition}\label{def:formalsmoothness} Let $f: Y\to X$ be a map of v-stacks. Then $f$ is formally smooth if for any affinoid perfectoid space $S$ of characteristic $p$ with a Zariski closed subspace $S_0\subset S$, and any commutative diagram
\[\xymatrix{
S_0\ar[r]^{g_0}\ar[d] & Y\ar[d]^f\\
S\ar[r]^h & X,
}\]
there is some \'etale map $S'\to S$ containing $S_0$ in its image and a map $g: S'\to Y$ fitting in a commutative diagram
\[\xymatrix{
S'\times_S S_0\ar[r]\ar[d] & S_0\ar[r]^{g_0}\ar[d] & Y\ar[d]^f\\
S'\ar[urr]^g\ar[r] & S\ar[r]^h & X.
}\]
\end{definition}

This kind of formal smoothness is closely related to the notion of {\it absolute neighborhood retracts} (\cite{BorsukRetractLivre}, \cite{Dold}). In fact, suppose $Y\to X$ is formally smooth with $Y$ and $X$ affinoid perfectoid. Choose a Zariski closed embedding $Y\hookrightarrow \mathbb{B}^I_X$ for some set $I$. Then there exists an \'etale neighborhood $U\to \mathbb{B}^I_X$ of $Y$ such that the closed embedding
$$
i:Y\times_{\mathbb{B}^I_X} U\hookrightarrow U
$$
admits a retraction $r:U\to Y\times_{\mathbb{B}^I_X} U$, $r\circ i =Id$. Thus, $Y/X$ is a retract of an (\'etale) neighborhood inside a ball$/X$.

\subsection{Examples and basic properties}

We will see that formally smooth morphisms share analogous properties to cohomologically smooth morphisms.
Let's begin with the following  observations:
\begin{altenumerate}
\item The composite of two formally smooth morphisms is formally smooth,
\item The formally smooth property is stable under pullback: if $Y\to X$ is formally smooth and $X'\to X$ is any map then $Y\times_X X'\to X'$ is formally smooth.
\item \'Etale maps are formally smooth.
\item For morphisms of locally spatial diamonds, formal smoothness is \'etale local on the source and the target.
\end{altenumerate}

Let us observe that, like cohomologically smooth morphisms, formally smooth morphisms are universally open.

\begin{proposition}\label{prop:formally smooth are open}
Formally smooth maps are universally open.
\end{proposition}
\begin{proof}
Let $Y\to X$ be formally smooth. We can suppose $X$ is affinoid perfectoid. Since any open subsheaf of $Y$ is formally smooth over $X$ we are reduced to prove that the image of $Y\to X$ is open. Let $S\to Y$ be a morphism with $S$ affinoid perfectoid. Choose a Zariski closed embedding $S\hookrightarrow \mathbb{B}^I_X$ for some set $I$. The formal smoothness assumption implies that there exists an \'etale neighborhood $U\to \mathbb{B}^I_X$ of $S\subset\mathbb{B}^I_X$ such that $S\times_{\mathbb B^I_X} U\to S\to Y$ extends to a map $U\to Y$; in particular, the image of $S\to Y\to X$ is contained in the image of $U\to Y\to X$, and it suffices to prove that the latter is open. We can suppose $U$ is quasicompact and separated over $\mathbb B^I_X$. Writing $\mathbb{B}^I_X=\varprojlim_J \mathbb{B}^J_X$ where $J$ goes through the set of finite subsets of $I$, there exists a some $J\subset I$ finite and $V\to \mathbb{B}^J_X$ such that $U\to \mathbb{B}^I_X$ is the pullback of $V\to \mathbb{B}^J_X$ via the projection $\mathbb{B}^I_X\to \mathbb{B}^J_X$, cf. \cite[Proposition 6.4]{ECoD}. Since $\mathbb{B}^I_X\to \mathbb{B}^J_X$ is a v-cover, 
$$
\mathrm{Im} (U\to X)=\mathrm{Im} (V\to X).
$$
Now, using that $V\to \mathbb{B}^J_X\to X$ is open, since cohomologically smooth for example, this is an open subset of $X$.
\end{proof}

Let us begin with some concrete examples. In the following, $\mathbb{B}\to \ast$ is the v-sheaf $\O^+$ on $\mathrm{Perf}_k$ and $\mathbb A^1\to \ast$ is the v-sheaf $\O$.

\begin{proposition}
The morphisms $\mathbb{B}\to \ast$, $\mathbb A^1\to \ast$ and $\Spd \mathcal O_E\to \ast$ are formally smooth.
\end{proposition}

\begin{proof}
Let $S_0=\Spa(R_0,R_0^+)\hookrightarrow S=\Spa(R,R^+)$ be a Zariski closed embedding of affinoid perfectoid spaces. Then $R\to R_0$ is surjective, which immediately shows that $\mathbb A^1\to \ast$ is formally smooth. The case of $\mathbb B\to \ast$ follows as $\mathbb B\subset \mathbb A^1$ is open. For $\Spd \mathcal O_E$, note that any untilt of $S_0$ can be given by some element $\xi\in W_{\mathcal O_E}(R_0^+)$ of the form $\xi_0=\pi+\sum_{n=0}^\infty \pi^i[r_{i,0}]$ where all $r_i\in R_0^{\circ\circ}$. But $R^{\circ\circ}\to R_0^{\circ\circ}$ is surjective (cf.~the discussion after \cite[Definition 5.7]{ECoD}), so one can lift all $r_{i,0}\in R_0^{\circ\circ}$ to $r_i\in R^{\circ\circ}$, and then $\xi=\pi+\sum_{n=0}^\infty \pi^i [r_i]$ defines an untilt of $S$ over $\mathcal O_E$ lifting the given one on $S_0$.
\end{proof}

\begin{corollary}
Is $f:Y\to X$ is a smooth morphism of analytic adic spaces over $\mathbb Z_p$ then $f^\diamond:Y^\diamond\to X^\diamond$ is formally smooth.
\end{corollary}

\begin{proof} Any smooth morphism is locally \'etale over a finite-dimensional ball.
\end{proof}

Let us remark the following.

\begin{proposition}\label{prop:formallysmoothetalesurj}
 If $f: Y\to X$ is a formally smooth and surjective map of v-stacks, then $f$ is surjective as a map of \'etale stacks. Equivalently, in case $X$ is a perfectoid space, the map $f$ splits over an \'etale cover of $X$.
\end{proposition}

\begin{proof}
We can suppose $X$ is affinoid. There exists a surjective morphism $X'\to X$ with $X'$ affinoid perfectoid and a section $s:X'\to Y$ of $Y\to X$ over $X'$. Let us choose a Zariski closed embedding $X'\hookrightarrow \mathbb{B}^I_X$. Applying the formal smoothness property we deduce there is an \'etale neighborhood $U\to \mathbb{B}^I_X$ of $X'\subset \mathbb{B}^I_X$ and a section over $U$ of $Y\to X$. It thus suffices to see that $U\to X$ admits a section over an \'etale cover of $X$. As in the proof of Proposition \ref{prop:formally smooth are open} there exists a finite subset $J\subset I$, and quasicompact \'etale map $V\to \mathbb{B}^J_X$ such that $U\to \mathbb{B}^I_X$ is the pullback of $V\to \mathbb{B}^J_X$ via the projection $\mathbb{B}^I_X\to \mathbb{B}^J_X$. This reduces us to the case $I$ is finite. We may also replace $V$ by its image in $\mathbb B^J_X$. At geometric points, the splitting follows from \cite[Lemma 9.5]{ECoD}. Approximating a section over a geometric point over an \'etale neighborhood then gives the desired splitting on an \'etale cover.
\end{proof}

According to \cite[Proposition 23.13]{ECoD} cohomological smoothness is cohomologically smooth local on the source. The same holds for formally smooth morphisms.

\begin{corollary}\label{cor:formally smooth formally smooth local}
Let $f:Y\to X$ be a morphism of v-stacks. Suppose there exists a v-surjective formally smooth morphism of v-stacks $g:Y'\to Y$ such that $f\circ g$ is formally smooth. Then $f$ is formally smooth.
\end{corollary}

\begin{proof} Given a test diagram $g_0: S_0\to Y$, $h: S\to X$ as in Definition~\ref{def:formalsmoothness}, we can first lift $S_0\to Y$ \'etale locally to $Y'$ by Proposition~\ref{prop:formallysmoothetalesurj}, and the required \'etale neighborhoods lift to $S$ by \cite[Proposition 6.4]{ECoD} applied to $S_0$ as the intersection of all open neighborhoods in $S$. Thus, the diagram can be lifted to a similar test diagram for $Y'\to X$, which admits a solution by assumption.
\end{proof}

Let us remark the following.

\begin{proposition}\label{prop:BunGformallysmooth} The stack $\Bun_G\to \ast$ is formally smooth.
\end{proposition}

\begin{proof} Let $S_0=\Spa(R_0,R_0^+)\subset S=\Spa(R,R^+)$ be a Zariski closed immersion of affinoid perfectoid spaces over $\Spd k$, and fix a pseudouniformizer $\varpi\in R$. Let $\mathcal E_0$ be a $G$-bundle on $X_{S_0}$. Pick any geometric point $\Spa(C,C^+)\to S_0$; we intend to find an \'etale neighborhood $U\to S$ of $\Spa(C,C^+)$ in $S$ such that the $G$-bundle over $U\times_S S_0$ extends to $U$.

Note that the pullback of $\mathcal E_0$ to $Y_{C,[1,q]}$ is a trivial $G$-bundle, by Theorem~\ref{thm:classifyGbun}. From \cite[Proposition 5.4.21]{GabberRamero} (applied with $R=\varinjlim_V \O^+(Y_{V,[1,q]})$, $t=\pi$ and $I=(\pi)$, where $V\to S_0$ runs over \'etale neighborhoods of $\Spa(C,C^+)$ in $S_0$; all of these lift to $S$) it follows that after passing to an \'etale neighborhood as above, we can assume that the pullback of $\mathcal E_0$ to $Y_{S_0,[1,q]}$ is a trivial $G$-bundle. In that case, $\mathcal E_0$ is given by some matrix $A\in G(B_{R_0,[1,1]})$ encoding the descent. Applying \cite[Proposition 5.4.21]{GabberRamero} again, with
\[
R=\varinjlim_{S_0\subset U\subset S} \O^+(Y_{U,[1,1]}),\ t=\pi,\ I=(\pi),
\]
then shows that we may lift $A$ into a neighborhood, as desired.
\end{proof}

The following is the analog of Proposition~\ref{prop:relativebanachcolmez1}~(iii).

\begin{proposition}\label{prop:formal smoothness of BCs}
Let $S$ be a perfectoid space and let $[\mathcal E_1\to \mathcal E_0]$ be a map of vector bundles on $X_S$ such that $\mathcal E_0$ is everywhere of positive Harder--Narasimhan slopes, and $\mathcal E_1$ is everywhere of negative Harder--Narasimhan slopes. Then $\BC([\mathcal E_1\to \mathcal E_0])\to S$ is formally smooth.
\end{proposition}

\begin{proof} Using the exact sequence
\[
0\to \BC(\mathcal E_0)\to \BC([\mathcal E_1\to \mathcal E_0])\to \BC(\mathcal E_1[1])\to 0
\]
and Proposition~\ref{prop:h1loczero}~(iii) to get \'etale local surjectivity of the second map, one reduces to the individual cases of $\BC(\mathcal E_0)$ and $\BC(\mathcal E_1[1])$. For $\mathcal E=\mathcal E_0$, we can use Corollary~\ref{cor:niceshortexseq} to choose \'etale locally on $S$ (and after replacing $\mathcal E$ by the direct sum with another bundle) a short exact sequence
\[
0\to \mathcal G\to \mathcal O_{X_S}(\tfrac 1r)^m\to \mathcal E\to 0
\]
where $\mathcal G$ is fiberwise on $S$ semistable of positive slope. Moreover, by Proposition~\ref{prop:h1loczero}~(iii), one can also ensure that $H^1(X_{S'},\mathcal G|_{X_{S'}})=0$ for all affinoid perfectoid spaces $S'\to S$. In particular, if $S_0\subset S$ is a Zariski closed immersion of affinoid perfectoid spaces, the map $\mathcal O_{X_S}(\tfrac 1r)^m(S_0)\to \mathcal E(S_0)$ is surjective, and we can replace $\mathcal E$ by $\mathcal O_{X_S}(\frac 1r)^m$. But then Proposition~\ref{prop:standardbanachcolmez}~(iv) shows that this Banach--Colmez space is representable by a perfectoid open unit disc, which is formally smooth.

For $\mathcal E=\mathcal E_1$, we can use Theorem~\ref{thm:O1ample} to find a short exact sequence
\[
0\to \mathcal E\to \mathcal O_{X_S}(d)^m\to \mathcal G\to 0
\]
for some $d,m>0$ (so $\mathcal G$ necessarily has only positive slopes), and this induces an exact sequence
\[
0\to \BC(\mathcal O_{X_S}(d)^m)\to \BC(\mathcal G)\to \BC(\mathcal E[1])\to 0
\]
where the middle term is formally smooth by the preceding, and the map $\BC(\mathcal G)\to \BC(\mathcal E[1])$ is formally smooth (as \'etale locally surjective and its fibre $\BC(\mathcal O_{X_S}(d)^m)$ is formally smooth). We conclude by Corollary~\ref{cor:formally smooth formally smooth local}.
\end{proof}

\section{A Jacobian criterion}\label{sec:jacobian}

The goal of this section is to prove that certain geometrically defined diamonds are cohomologically smooth when one expects them to be. We regard this result as the most profound in the theory of diamonds so far: While we cannot control much of the geometry of these diamonds, in particular we have no way to relate them to (perfectoid) balls in any reasonable way, we can still prove relative Poincar\'e duality for them. The spaces considered below also appear quite naturally in a variety of contexts, so we expect the result to have many applications.

The setup is the following. Let $S$ be a perfectoid space and let $Z\to X_S$ be a smooth map of sous-perfectoid adic spaces --- defining this concept of smoothness will be done in a first subsection, but it is essentially just a family of smooth rigid spaces over $X_S$, in the usual sense. One can then consider the v-sheaf $\mathcal M_Z$ of sections of $Z\to X_S$, sending any perfectoid space $S'\to S$ to the set of maps $X_{S'}\to Z$ lifting $X_{S'}\to X_S$. In general, we cannot prove that $\mathcal M_Z$ is a locally spatial diamond, but this turns out to be true when $Z$ is quasiprojective in the sense that it is a Zariski closed subspace of an open subset of (the adic space) $\mathbb P^n_{X_S}$ for some $n\geq 0$.

In general, the space $\mathcal M_Z\to S$ is not (cohomologically) smooth: If tangent spaces of $\mathcal M_Z\to S$ would exist, one would expect their fibre over $S'\to \mathcal M_Z$, given by some section $s: X_{S'}\to Z$, to be given by $H^0(X_{S'},s^\ast T_{Z/X_S})$, where $T_{Z/X_S}$ is the tangent bundle of $Z\to X_S$; and then an obstruction space would be given by $H^1(X_{S'},s^\ast T_{Z/X_S})$. Thus, one can expect smoothness to hold only when $H^1 (X_{S'},s^\ast T_{Z/X_S})$ vanishes. This holds true, locally on $S'$, if all slopes of $s^\ast T_{Z/X_S}$ are positive (by Proposition~\ref{prop:h1loczero}~(iii)), suggesting the following definition.

\begin{definition}\label{def:Mzsmooth} Let $\mathcal M_Z^{\mathrm{sm}}\subset \mathcal M_Z$ be the open subfunctor of all sections $s: X_{S'}\to Z$ such that $s^\ast T_{Z/X_S}$ has everywhere positive Harder--Narasimhan slopes.
\end{definition}

Roughly speaking, one expects $\mathcal M_Z^{\mathrm{sm}}$ to look infinitesimally like the Banach--Colmez space $\mathcal{BC}(s^\ast T_{Z/X_S})$; these indeed are cohomologically smooth when all slopes are positive, by Proposition~\ref{prop:relativebanachcolmez1}~(iii). Unfortunately, we are unable to prove a direct relation of this sort; however, we will be able to relate these spaces via a ``deformation to the normal cone''.

Our goal is to prove the following theorem.

\begin{theorem}\label{thm:jacobiancriterion} Let $S$ be a perfectoid space and let $Z\to X_S$ be a smooth map of sous-perfectoid spaces such that $Z$ admits a Zariski closed immersion into an open subset of (the adic space) $\mathbb P^n_{X_S}$ for some $n\geq 0$. Then $\mathcal M_Z$ is a locally spatial diamond, the map $\mathcal M_Z\to S$ is compactifiable, and $\mathcal M_Z^{\mathrm{sm}}\to S$ is cohomologically smooth.

Moreover, for a geometric point $x: \Spa C\to \mathcal M_Z^{\mathrm{sm}}$ given by a map $\Spa C\to S$ and a section $s: X_C\to Z$, the map $\mathcal M_Z^{\mathrm{sm}}\to S$ is at $x$ of $\ell$-dimension equal to the degree of $s^\ast T_{Z/X_S}$.
\end{theorem}

\begin{remark} The map $\mathcal M_Z^{\mathrm{sm}}\to S$ is a natural example of a map that is only locally of finite dimension, but not globally so (as there are many connected components of increasing dimension).
\end{remark}

\begin{remark}
In the ``classical context'' of algebraic curves the preceding theorem is the following (easy) result. Let $X/k$ be a proper smooth curve and $Z\to X$ be quasi-projective smooth. Consider  $\mathcal{M}_Z$ the functor on $k$-schemes that sends $S$ to morphisms $s:X\times_k S\to Z$ over $X$. This is representable by a quasi-projective scheme over $\Spec(k)$. Let $\mathcal{M}_Z^{\mathrm{sm}}$ be the open sub-scheme defined by the condition that if $s:X\times_k S\to Z$ is an $S$-point of $\mathcal{M}_Z^{\mathrm{sm}}$ then the vector bundle $s^* T_{Z/X}$ has no $H^1$ fiberwise on $S$. Then $\mathcal{M}_Z^{\mathrm{sm}}\to \Spec (k)$ is smooth.
\end{remark}

\begin{remark}
Suppose that $W$ is a smooth quasi-projective $E$-scheme. The moduli space $\mathcal{M}_Z$ with $Z=W\times_E X_S$ classifies morphisms $X_S\to W$ i.e. $\mathcal{M}_Z$ is a moduli of morphisms from families of Fargues--Fontaine curves to $W$. This is some kind of ``Gromov--Witten'' situation.
\end{remark}

\begin{remark}
We could have made the a priori weaker assumption that $Z$ admits a Zariski closed immersion inside an open subset of $\mathbb{P} (\E)$ where $\E$ is a vector bundle on $X_S$. Nevertheless, since the result is local on $S$ and we can suppose it is affinoid perfectoid, and since when $S$ is affinoid perfectoid $\O_{X_S}(1)$ ``is ample'' i.e. there is a surjection $\O_{X_S} (-N)^n\twoheadrightarrow \E$ for $N,n\gg 0$, this assumption is equivalent to the one we made i.e.~we can suppose $\E$ is free.
\end{remark}

\begin{example}[The Quot diamond] Let $\E$ be a vector bundle on $X_S$. We denote by 
$$
\mathrm{Quot}_\E\lto S
$$
the moduli space over $S$ of locally free quotients of $\E$. Fixing the rank of such a quotient, one sees that $\mathrm{Quot}_\E$ is a finite disjoint union of spaces $\mathcal{M}_Z$ with $Z\to X_S$ a Grassmannian of quotients of $\E$. This is thus representable in locally spatial diamonds, compactifiable, of locally finite $\dimtrg$.

Let $\mathrm{Quot}_\E^{\mathrm{sm}}\subset\mathrm{Quot}_\E$ be the open subset parametrizing quotients $u:\E\to \mathcal{F}$ such that fiberwise, the greatest slope of $\ker u$ is strictly less than the smallest slope of $\mathcal{F}$. According to Theorem \ref{thm:jacobiancriterion} this is cohomologically smooth over $S$. 

Fix an integer $n\geq 1$. For some $N\in \mathbb{Z}$ and $r\in\mathbb{N}_{\geq 1}$, let $\mathrm{Quot}_{\O(N)^r}^{n,\mathrm{sm},\circ}$ be the open subset of $\mathrm{Quot}_{\O(N)^r}^{\mathrm{sm}}$
where the quotient has rank $n$ and its slopes are greater than $N$. When $N$ and $r$ vary one constructs, as in the ``classical case'', cohomologically smooth charts on $\Bun_{\GL_n}$ using $\mathrm{Quot}_{\O(N)^r}^{n,\mathrm{sm},\circ}$. In fact, the morphism 
$$
\mathrm{Quot}_{\O(N)^r}^{n,\mathrm{sm},\circ} \lto \Bun_{\GL_n}
$$
given by the quotient of $\O(N)^r$ is separated cohomologically smooth. When pulled back by a morphism $S\to \Bun_{\GL_n}$ with $S$ perfectoid, this is an open subset of a positive Banach--Colmez space.
\end{example}

We will not use the Quot diamond in the following. In section \ref{sec:localcharts}, using the Jacobian criterion, we will construct charts on $\Bun_G$ for any $G$ that are better suited to our needs.

\subsection{Smooth maps of sous-perfectoid adic spaces}

We need some background about smooth morphisms of adic spaces in non-noetherian settings. We choose the setting of sous-perfectoid adic spaces as defined by Hansen-Kedlaya, \cite{HansenKedlaya}, cf.~\cite[Section 6.3]{Berkeley}. Recall that an adic space $X$ is sous-perfectoid if it is analytic and admits an open cover by $U=\Spa(R,R^+)$ where each $R$ is a sous-perfectoid Tate algebra, meaning that there is some perfectoid $R$-algebra $\tilde{R}$ such that $R\to \tilde{R}$ is a split injection in the category of topological $R$-modules.

The class of sous-perfectoid rings $R$ is stable under passage to rational localizations, finite \'etale maps, and $R\langle T_1,\ldots,T_n\rangle$. As smooth maps should be built from these basic examples, we can hope for a good theory of smooth maps of sous-perfectoid spaces.

Recall that a map $f: Y\to X$ of sous-perfectoid adic spaces is \'etale if locally on the source and target it can be written as an open immersion followed by a finite \'etale map.

\begin{definition}\label{def:smoothsous-perfectoid} Let $f: Y\to X$ be a map of sous-perfectoid adic spaces. Then $f$ is smooth if one can cover $Y$ by open subsets $V\subset Y$ such that there are \'etale maps $V\to \mathbb B^d_X$ for some integer $d\geq 0$.
\end{definition}

It can immediately be checked that analytifications of smooth schemes satisfy this condition.

\begin{proposition}\label{prop:analytificationsmooth} Let $X=\Spa(A,A^+)$ be an affinoid sous-perfectoid adic space, and let $f_0: Y_0\to \Spec A$ be a smooth map of schemes. Let $f: Y\to X$ be the analytification of $f_0: Y_0\to \Spec A$, representing the functor taking $\Spa(B,B^+)\to \Spa(A,A^+)$ to the $\Spec B$-valued points of $Y_0\to \Spec A$. Then $f: Y\to X$ is smooth.
\end{proposition}

\begin{proof} Locally, $f_0$ is the composite of an \'etale and the projection from affine space. This means that its analytification is locally \'etale over the projection from the analytification of affine space, which is a union of balls, giving the result.
\end{proof}

Let us analyze some basic properties of smooth maps of sous-perfectoid adic spaces.

\begin{proposition}\label{prop:smoothsous-perfectoid} Let $f: Y\to X$ and $g: Z\to Y$ be maps of sous-perfectoid adic spaces.
\begin{altenumerate}
\item The property of $f$ being smooth is local on $Y$.
\item If $f$ and $g$ are smooth, then so is $f\circ g: Z\to X$.
\item If $h: X'\to X$ is any map of sous-perfectoid adic spaces and $f$ is smooth, then the fibre product $Y'=Y\times_X X'$ in adic spaces exists, is sous-perfectoid, and $f': Y'\to X'$ is smooth.
\item If $f$ is smooth, then $f$ is universally open.
\item If $f$ is smooth and surjective, then there is some \'etale cover $X'\to X$ with a lift $X'\to Y$.
\end{altenumerate}
\end{proposition}

Regarding part (i), we note that we will see in Proposition~\ref{prop:smoothetalelocal} that the property of $f$ being smooth is in fact \'etale local on $Y$ (and thus smooth local on $X$ and $Y$, using (v)).

\begin{proof} Part (i) is clear from the definition. For part (ii), the composite is locally a composite of an \'etale map, a projection from a ball, an \'etale map, and another projection from a ball; but we can swap the two middle maps, and use that composites of \'etale maps are \'etale. Part (iii) is again clear, by the stability properties of sous-perfectoid rings mentioned above. For part (iv), it is now enough to see that $f$ is open, and we can assume that $f$ is a composite of an \'etale map and the projection from a ball, both of which are open. For part (v), using that $f$ is open, we can work locally on $Y$ and thus assume again that it is a composite of an \'etale map and the projection from a ball; we can then replace $Y$ by its open image in the ball. By \cite[Lemma 9.5]{ECoD}, for any geometric point $\Spa(C,C^+)\to X$ of $X$, we can find a lift to $Y$. Writing the geometric point as the limit of affinoid \'etale neighborhoods, the map to $Y\subset \mathbb B^d_X$ can be approximated at some finite stage, and then openness of $Y$ ensures that it will still lie in $Y$. This gives the desired \'etale cover of $X$ over which $f$ splits.
\end{proof}

Of course, the most important structure of a smooth morphism is its module of K\"ahler differentials. Recall that if $Y$ is sous-perfectoid, then one can define a stack (for the \'etale topology) of vector bundles on $Y$, such that for $Y=\Spa(R,R^+)$ affinoid with $R$ sous-perfectoid, the category of vector bundles is equivalent to the category of finite projective $R$-modules; see \cite{KedlayaLiu1}, \cite[Theorem 5.2.8, Proposition 6.3.4]{Berkeley}. By definition, a vector bundle on $Y$ is an $\mathcal O_Y$-module that is locally free of finite rank.

\begin{definition}\label{def:kahlerdiff} Let $f: Y\to X$ be a smooth map of sous-perfectoid adic spaces, with diagonal $\Delta_f: Y\to Y\times_X Y$. Let $\mathcal I_{Y/X}\subset \mathcal O_{Y\times_X Y}$ be the ideal sheaf. Then
\[
\Omega^1_{Y/X} := \mathcal I_{Y/X}/\mathcal I_{Y/X}^2
\]
considered as $\mathcal O_{Y\times_X Y}/\mathcal I_{Y/X} = \mathcal O_Y$-module.
\end{definition}

It follows from the definition that there is a canonical $\mathcal O_X$-linear derivation $d: \mathcal O_Y\to \Omega^1_{Y/X}$, given by $g\mapsto g\otimes 1-1\otimes g$.

\begin{proposition}\label{prop:kahlerdiff} Let $f: Y\to X$ be a smooth map of sous-perfectoid adic spaces. Then $\Omega^1_{Y/X}$ is a vector bundle on $Y$. There is a unique open and closed decomposition $Y=Y_0\sqcup Y_1\sqcup\ldots\sqcup Y_n$ such that $\Omega^1_{Y/X}|_{Y_d}$ is of rank $d$ for all $d=0,\ldots,n$. In that case, for any nonempty open subset $V\subset Y_d$ with an \'etale map $V\to \mathbb B^{d'}_X$, necessarily $d'=d$.
\end{proposition}

We will say that $f$ is smooth of dimension $d$ if $\Omega^1_{Y/X}$ is of rank $d$. By the proposition, this is equivalent to asking that $Y$ can be covered by open subsets $V$ that admit \'etale maps $V\to \mathbb B^d_X$. In particular, $f$ is smooth of dimension $0$ if and only if it is \'etale.

\begin{proof} It is enough to show that if $f$ is a composite of an \'etale map $Y\to \mathbb B^d_X$ with the projection to $X$, then $\Omega^1_{Y/X}$ is isomorphic to $\mathcal O_Y^d$. Indeed, this implies that $\Omega^1_{Y/X}$ is a vector bundle in general, of the expected rank; and the decomposition into open and closed pieces is then a general property of vector bundles.

Let $Y'=\mathbb B^d_X$. Then $Y\times_X Y\to Y'\times_X Y'$ is \'etale, and the map $Y\to Y'\times_{Y'\times_X Y'} (Y\times_X Y)$ is an open immersion (as the diagonal of the \'etale map $Y\to Y'$). It follows that $\mathcal I_{Y/X}$ is the pullback of $\mathcal I_{Y'/X'}$. But $Y'\hookrightarrow Y'\times_X Y'$ is of the form
\[
\Spa(R\langle T_1,\ldots,T_n\rangle,R^+\langle T_1,\ldots,T_n\rangle)\hookrightarrow \Spa(R\langle T_1^{(1)},\ldots,T_n^{(1)},T_1^{(2)},\ldots,T_n^{(2)}\rangle,R^+\langle T_1^{(1)},\ldots,T_n^{(2)}\rangle)
\]
if $X=\Spa(R,R^+)$, and the ideal sheaf is given by $(T_1^{(1)}-T_1^{(2)},\ldots,T_n^{(1)}-T_n^{(2)})$. This defines a regular sequence after any \'etale localization, by the lemma below. This gives the claim.
\end{proof}

\begin{lemma}\label{lem:regularsequence} Let $X=\Spa(R\langle T_1,\ldots,T_n\rangle,R^+\langle T_1,\ldots,T_n\rangle)$ where $R$ is a sous-perfectoid Tate ring, let $Y=\Spa(S,S^+)$ where $S$ is a sous-perfectoid Tate ring, and let $f: Y\to X$ be a smooth map. Then $T_1,\ldots,T_n$ define a regular sequence on $S$ and $(T_1,\ldots,T_n)S\subset S$ is a closed ideal.
\end{lemma}

\begin{proof} By induction, one can reduce to the case $n=1$. The claim can be checked locally, so we can assume that $Y$ is \'etale over $\mathbb B^d_X$ for some $d$; replacing $X$ by $\mathbb B^d_X$, we can then assume that $f$ is \'etale. Let $Y_0\subset Y$ be the base change to $X_0=\Spa(R,R^+)=V(T)\subset X$; then $Y$ and $Y_0\times_{X_0} X$ are both \'etale over $X$ and become isomorphic over $X_0\subset X$. By spreading of \'etale maps, this implies that they are isomorphic after base change to $X'=\Spa(R\langle T'\rangle,R^+\langle T'\rangle)$ where $T'=\varpi^n T$ for some $n$ (and $\varpi$ is a pseudouniformizer of $R$). This easily implies the result.
\end{proof}

Locally around a section, any smooth space is a ball:

\begin{lemma}\label{lem:smootharoundsection} Let $f: Y\to X$ be a smooth map of sous-perfectoid spaces with a section $s: \Spa(K,K^+)\to Y$ for some point $\Spa(K,K^+)\to X$. Then there are open neighborhoods $U\subset X$ of $\Spa(K,K^+)$ and $V\subset Y$ of $s(\Spa(K,K^+))$ such that $V\cong \mathbb B^d_U$.
\end{lemma}

\begin{proof} We can assume that $X$ and $Y$ are affinoid. If $f$ is \'etale, then any section extends to a small neighborhood (e.g.~by \cite[Lemma 15.6, Lemma 12.17]{ECoD}), and any section is necessarily \'etale and thus open, giving the result in that case. In general, we may work locally around the given section, so we can assume that $f$ is the composite of an \'etale map $Y\to \mathbb B^d_X$ and the projection to $X$. Using the \'etale case already handled, we can assume that $Y$ is an open subset of $\mathbb B^d_X$. Any section $\Spa(K,K^+)\to \mathbb B^d_X$ has a cofinal system of neighborhoods that are small balls over open subsets of $X$: The section is given by $d$ elements $T_1,\ldots,T_d\in K^+$, and after picking a pseudouniformizer $\varpi$ and shrinking $X$, one can find global sections $T_1',\ldots,T_d'$ of $\mathcal O_X^+(X)$ such that $T_i\equiv T_i'\mod \varpi^n$. Then $\{|T_1'|,\ldots,|T_d'|\leq |\varpi|^n\}$ is a small ball over $X$, and the intersection of all these is $\Spa(K,K^+)$. Thus, one of these neighborhoods is contained in $Y$, as desired.
\end{proof}

\begin{proposition}\label{prop:smoothmapsbetweensmooth} Let $f_i: Y_i\to X$, $i=1,2$, be smooth maps of sous-perfectoid adic spaces, and let $g: Y_1\to Y_2$ be a map over $X$.
\begin{altenumerate}
\item If $g$ is smooth, then the sequence
\[
0\to g^\ast \Omega^1_{Y_2/X}\to \Omega^1_{Y_1/X}\to \Omega^1_{Y_1/Y_2}\to 0
\]
is exact.
\item Conversely, if $g^\ast \Omega^1_{Y_2/X}\to \Omega^1_{Y_1/X}$ is a locally split injection, then $g$ is smooth.
\end{altenumerate}

In particular, if $g^\ast \Omega^1_{Y_2/X}\to \Omega^1_{Y_1/X}$ is an isomorphism, then $g$ is \'etale.
\end{proposition}

\begin{proof} Part (i) follows from a routine reduction to the case of projections from balls, where it is clear. For part (ii), we may assume that $Y_1\to \mathbb B^{d_1}_X$ and $Y_2\to \mathbb B^{d_2}_X$ are \'etale. It suffices to see that the composite $Y_1\to Y_2\to \mathbb B^{d_2}_X$ is smooth, as $g$ is the composite of its base change $Y_1\times_{\mathbb B^{d_2}_X} Y_2\to Y_2$ with the section $Y_1\to Y_1\times_{\mathbb B^{d_2}_X} Y_2$ of the \'etale map $Y_1\times_{\mathbb B^{d_2}_X} Y_2\to Y_1$; any such section is automatically itself \'etale. Thus, we may assume that $Y_2=\mathbb B^{d_2}_X$. Locally on $Y_1$, we may find a projection $g': \mathbb B^{d_1}_X\to \mathbb B^{d_1-d_2}_X$ so that
\[
(g^{\prime\ast} \Omega^1_{\mathbb B^{d_1-d_2}_X/X})|_{Y_1}
\]
is an orthogonal complement of $g^\ast \Omega^1_{Y_2/X}$. (Indeed, looking at $\mathcal O_X^{d_1}\cong \Omega^1_{Y_1/X}\to g^\ast \Omega^1_{Y_2/X}$ as a point of the Grassmannian and using the standard affine cover of the Grassmannian, one shows that one may take $g'$ to simply be a projection to a subset of $d_1-d_2$ of the coordinates.) Thus, we can assume that $d_1=d_2=:d$, and $g^\ast \Omega^1_{Y_2/X}\to \Omega^1_{Y_1/X}$ is an isomorphism.

Our aim is now to prove that $g: Y_1\to Y_2 = \mathbb B^d_X$ is \'etale. We may assume that all of $X$, $Y_1$ and $Y_2$ are affinoid. Passing to the fibre over a point $S=\Spa(K,K^+)\to X$, this follows from a result of Huber, \cite[Proposition 1.6.9 (iii)]{Huber}. The resulting \'etale map $Y_{1,S}\to Y_{2,S}$ deforms uniquely to a quasicompact separated \'etale map $Y_{1,U}'\to Y_{2,U}$ for a small enough neighborhood $U\subset X$ of $S$, by \cite[Lemma 12.17]{ECoD}. Moreover, the map $Y_{1,U}\to Y_{2,U}$ lifts uniquely to $Y_{1,U}\to Y_{1,U}'$ for $U$ small enough, by the same result. Replacing $X$ by $U$, $Y_1$ by $Y_{1,U}$ and $Y_2$ by $Y_{1,U}'$, we can now assume that $g: Y_1\to Y_2$ is a map between sous-perfectoid spaces smooth over $X$ that is an isomorphism on one fibre. It is enough to see that it is then an isomorphism in a neighborhood. To see this, we may in fact work locally on $Y_2$.

For this, we study $Y_1\subset Y_1\times_X Y_2\to Y_2$: Here $Y_1\times_X Y_2\to Y_2$ is smooth, and $Y_1\subset Y_1\times_X Y_2$ is locally the vanishing locus of $d$ functions (as $Y_2\subset Y_2\times_X Y_2$ is). Moreover, over fibres lying over the given point of $X$, the map $Y_1\to Y_2$ becomes an isomorphism, and in particular gives a section of $Y_1\times_X Y_2\to Y_2$. By Lemma~\ref{lem:smootharoundsection}, after shrinking $Y_2$, we can assume that there is an open neighborhood $V\subset Y_1\times_X Y_2$ such that $V\cong \mathbb B^d_{Y_2}$. Inside there, $Y_1$ is (locally) given by the vanishing of $d$ functions, and is only a point in one fibre. Now the result follows from the next lemma, using $Y_2$ in place of $X$.
\end{proof}

\begin{lemma} Let $X=\Spa(A,A^+)$ be a sous-perfectoid affinoid adic space with a point $X'=\Spa(K,K^+)\to X$. Let $f_1,\ldots,f_n\in A^+\langle T_1,\ldots,T_n\rangle$ be functions such that
\[
K\to K\langle T_1,\ldots,T_n\rangle/(f_1,\ldots,f_n)
\]
is an isomorphism. Then, after replacing $X$ by an open neighborhood of $X'$, the map
\[
A\to A\langle T_1,\ldots,T_n\rangle/(f_1,\ldots,f_n)
\]
is an isomorphism.
\end{lemma}

\begin{proof} For any ring $B$ with elements $g_1,\ldots,g_n\in B$, consider the homological Koszul complex
\[
\mathrm{Kos}(B,(g_i)_{i=1}^n) = [B\to B^n\to\ldots\to B^n\xrightarrow{(g_1,\ldots,g_n)} B].
\]
We claim that, after shrinking $X$, we can in fact arrange that
\[
A\to \mathrm{Kos}(A\langle T_1,\ldots,T_n\rangle,(f_i)_{i=1}^n)
\]
is a quasi-isomorphism.

Note that all terms of these complexes are free Banach-$A$-modules, and thus the formation of this complex commutes with all base changes; and one can use descent to establish the statement. In particular, we can reduce first to the case that $X$ is perfectoid, and then to the case that $X$ is strictly totally disconnected. In that case, the map $A\to K$ is automatically surjective, and so we can arrange that under the isomorphism $K\cong K\langle T_1,\ldots,T_n\rangle/(f_1,\ldots,f_n)$, all $T_i$ are mapped to $0$. Moreover, applying another change of basis, we can arrange that the image of $f_i$ in $K\langle T_1,\ldots,T_n\rangle/(T_1,\ldots,T_n)^2$ is given by $a_i T_i$ for some nonzero scalar $a_i\in K^+$. Note that we are in fact allowed to also localize on $\mathbb B^n_X$ around the origin, as away from the origin the functions $f_1,\ldots,f_n$ locally generate the unit ideal (in the fibre, but thus in a small neighborhood). Doing such a localization, we can now arrange that $f_i\equiv T_i\mod \varpi$ for some pseudouniformizer $\varpi\in A^+$. But now in fact
\[
A^+\to \mathrm{Kos}(A^+\langle T_1,\ldots,T_n\rangle,(f_i)_{i=1}^n)
\]
is a quasi-isomorphism, as can be checked modulo $\varpi$, where it is the quasi-isomorphism
\[
A^+/\varpi\to \mathrm{Kos}(A^+/\varpi[T_1,\ldots,T_n],(T_i)_{i=1}^n).\qedhere
\]
\end{proof}

Let us draw some consequences. First, we have the following form of the Jacobian criterion in this setting.

\begin{proposition}\label{prop:jacobiancriterionsous-perfectoid} Let $f: Y\to X$ be a smooth map of sous-perfectoid adic spaces, and let $f_1,\ldots,f_r\in \mathcal O_Y(Y)$ be global functions such that $df_1,\ldots,df_r\in \Omega^1_{Y/X}(Y)$ can locally be extended to a basis of $\Omega^1_{Y/X}$. Then $Z=V(f_1,\ldots,f_r)\subset Y$ is a sous-perfectoid space smooth over $X$.
\end{proposition}

\begin{proof} We can assume that all $f_i\in \mathcal O_Y^+(Y)$ by rescaling, and we can locally find $f_{r+1},\ldots,f_n\in \mathcal O_Y^+(Y)$ such that $df_1,\ldots,df_n$ is a basis of $\Omega^1_{Y/X}$. This induces an \'etale map $Y\to \mathbb B^n_X$ by Proposition~\ref{prop:smoothmapsbetweensmooth}, and then $V(f_1,\ldots,f_r)\subset Y$ is the pullback of $\mathbb B^r_X\subset \mathbb B^n_X$, giving the desired result.
\end{proof}

Moreover, we can prove that being smooth is \'etale local on the source.

\begin{proposition}\label{prop:smoothetalelocal} Let $f: Y\to X$ be a map of sous-perfectoid adic spaces. Assume that there is some \'etale cover $j: V\to Y$ such that $f\circ j$ is smooth. Then $f$ is smooth.
\end{proposition}

\begin{proof} By \'etale descent of vector bundles on sous-perfectoid adic spaces, $\Omega^1_{Y/X}:= \mathcal I_{Y/X}/\mathcal I_{Y/X}^2$ is a vector bundle, together with an $\mathcal O_X$-linear derivation $d: \mathcal O_Y\to \Omega^1_{Y/X}$. We claim that locally we can find functions $f_1,\ldots,f_n\in \mathcal O_Y$ such that $df_1,\ldots,df_n\in \Omega^1_{Y/X}$ is a basis. To do this, it suffices to find such functions over all fibres $\Spa(K,K^+)\to X$, as any approximation will then still be a basis (small perturbations of a basis are still a basis). But over fibres, the equivalence of the constructions in \cite[1.6.2]{Huber} shows that the $df$ for $f\in \mathcal O_X$ form generators of $\Omega^1_{Y/X}$.

Thus, assume that there are global sections $f_1,\ldots,f_n$ such that $df_1,\ldots,df_n\in \Omega^1_{Y/X}$ are a basis. Rescaling the $f_i$ if necessary, they define a map $g: Y\to \mathbb B^d_X$ that induces an isomorphism $g^\ast \Omega^1_{\mathbb B^d_X/X}\to \Omega^1_{Y/X}$. By Proposition~\ref{prop:smoothmapsbetweensmooth}, the map $Y\to \mathbb B^d_X$ is \'etale locally on $Y$ \'etale. We may assume that $Y$ and $X$ are affinoid; in particular, all maps are separated. Then by \cite[Lemma 15.6, Proposition 11.30]{ECoD} also $Y\to \mathbb B^d_X$ is \'etale.
\end{proof}

Finally, we note that if $Y$ and $Y'$ are both smooth over a sous-perfectoid space $X$, then the concept of Zariski closed immersions $Y\hookrightarrow Y'$ over $X$ is well-behaved.

\begin{proposition}\label{prop:zariskiclosedsmooth} Let $f: Y\to X$, $f': Y'\to X$ be smooth maps of sous-perfectoid adic spaces, and let $g: Y\to Y'$ be a map over $X$. The following conditions are equivalent.
\begin{altenumerate}
\item There is a cover of $Y'$ by open affinoid $V'=\Spa(S',S'^+)$ such that $V=Y\times_{Y'} V'=\Spa(S,S^+)$ is affinoid and $S'\to S$ is surjective, with $S^+\subset S$ the integral closure of the image of $S'^+$.
\item For any open affinoid $V'=\Spa(S',S'^+)\subset Y'$, the preimage $V=Y\times_{Y'} V'=\Spa(S,S^+)$ is affinoid and $S'\to S$ is surjective, with $S^+\subset S$ the integral closure of the image of $S'^+$.
\end{altenumerate}

Moreover, in this case the ideal sheaf $\mathcal I_{Y\subset Y'}\subset \mathcal O_{Y'}$ is pseudocoherent in the sense of \cite{KedlayaLiu2}, and locally generated by sections $f_1,\ldots,f_r\in \mathcal O_{Y'}$ such that $df_1,\ldots,df_r\in \Omega^1_{Y'/X}$ can locally be extended to a basis.
\end{proposition}

\begin{proof} We first analyze the local structure under condition (1), so assume that $Y'=\Spa(S',S'^+)$ and $Y=\Spa(S,S^+)$ are affinoid, with $S'\to S$ surjective and $S^+\subset S$ the integral closure of the image of $S'^+$. It follows that $g^\ast \Omega^1_{Y'/X}\to \Omega^1_{Y/X}$ is surjective, and letting $d'$ and $d$ be the respective dimensions of $Y'$ and $Y$ (which we may assume to be constant), we see that $r=d'-d\geq 0$ and that locally we can find $f_1,\ldots,f_r\in \mathcal I_{Y\subset Y'}$ so that $df_1,\ldots,df_r$ generate the kernel of $g^\ast \Omega^1_{Y'/X}\to \Omega^1_{Y/X}$ (as the kernel is generated by the closure of the image of $\mathcal I_{Y\subset Y'}$). By Proposition~\ref{prop:jacobiancriterionsous-perfectoid}, the vanishing locus of the $f_i$ defines a sous-perfectoid space $Z\subset Y'$ that is smooth over $X$. The induced map $Y\to Z$ induces an isomorphism on differentials, hence is \'etale by Proposition~\ref{prop:smoothmapsbetweensmooth}; but it is also a closed immersion, hence locally an isomorphism.

We see that the ideal sheaf $\mathcal I_{Y\subset Y'}$ is locally generated by sections $f_1,\ldots,f_r$ as in the statement of the proposition. By the proof of Proposition~\ref{prop:jacobiancriterionsous-perfectoid} and Lemma~\ref{lem:regularsequence}, it follows that the ideal sheaf $\mathcal I_{Y\subset Y'}$ is pseudocoherent in the sense of \cite{KedlayaLiu2}.

To finish the proof, it suffices to show that (1) implies (2). By the gluing result for pseudocoherent modules of \cite{KedlayaLiu2}, the pseudocoherent sheaf $\mathcal I_{Y\subset Y'}$ over $V'$ corresponds to a pseudocoherent module $I\subset S'$, and then necessarily $V=\Spa(S,S^+)$ where $S=S'/I$ with $S^+\subset S$ the integral closure of the image of $S'^+$.
\end{proof}

\begin{definition}\label{def:zariskiclosedsmooth} In the setup of Proposition~\ref{prop:zariskiclosedsmooth}, the map $g$ is a Zariski closed immersion if the equivalent conditions are satisfied.
\end{definition}

\subsection{Maps from $X_S$ into $\mathbb P^n$}

Our arguments make critical use of the assumption that in Theorem~\ref{thm:jacobiancriterion}, the space $Z\to X_S$ is locally closed in $\mathbb P^n_{X_S}$. For this reason, we analyze the special case of $\mathbb P^n$ in this section.

\begin{proposition}\label{prop:gromovwittenPn} Let $n\geq 0$ and consider the small v-sheaf $\mathcal M_{\mathbb P^n}$ taking any perfectoid space $S$ to the set of maps $X_S\to \mathbb P^n_{E}$. Then $\mathcal M_{\mathbb P^n}\to \ast$ is partially proper and representable in locally spatial diamonds, and admits a decomposition into open and closed subspaces
\[
\mathcal M_{\mathbb P^n} = \bigsqcup_{m\geq 0} \mathcal M_{\mathbb P^n}^m
\]
such that each $\mathcal M_{\mathbb P^n}^m\to \ast$ has finite $\dimtrg$, and the degree of the pullback of $\mathcal O_{\mathbb P^n}(1)$ to $X_{\mathcal M_{\mathbb P^n}^m}$ is $m$. In fact, there is a canonical open immersion
\[
\mathcal M_{\mathbb P^n}^m\hookrightarrow (\BC(\mathcal O(m)^{n+1})\setminus \{0\})/\underline{E^\times}.
\]
\end{proposition}

\begin{proof} The degree of the pullback $\mathcal L/X_S$ of $\mathcal O_{\mathbb P^n}(1)$ to $X_S$ defines an open and closed decomposition according to all $m\in \mathbb Z$. Fix some $m$. Then over the corresponding subspace $\mathcal M_{\mathbb P^n}^m$, we can fix a trivialization $\mathcal L\cong \mathcal O_{X_S}(m)$, which amounts to an $\underline{E^\times}$-torsor. After this trivialization, one parametrizes $n+1$ sections of $\mathcal L\cong \mathcal O_{X_S}(m)$ without common zeroes. The condition of no common zeroes is an open condition on $S$: Indeed, the common zeroes form a closed subspace of $|X_S|$, and the map $|X_S|\to |S|$ is closed (see the proof of Lemma \ref{lemma:surjection ouvert}). This implies the desired description.
\end{proof}

\begin{proposition}\label{prop:embedZ} Let $S$ be a perfectoid space and let $Z\to X_S$ be a smooth map of sous-perfectoid adic spaces such that $Z$ admits a Zariski closed embedding into an open subspace of $\mathbb P^n_{X_S}$. Then the induced functor
\[
\mathcal M_Z\to \mathcal M_{\mathbb P^n_{X_S}}
\]
is locally closed. More precisely, for any perfectoid space $T\to \mathcal M_{\mathbb P^n_{X_S}}$, the preimage of $\mathcal M_Z$ is representable by some perfectoid space $T_Z\subset T$ that is \'etale locally Zariski closed in $T$, i.e.~there is some \'etale cover of $T$ by affinoid perfectoid $T'=\Spa(R,R^+)\to T$ such that $T_Z\times_T T'=\Spa(R_Z,R_Z^+)$ is affinoid perfectoid, with $R\to R_Z$ surjective and $R_Z^+\subset R_Z$ the integral closure of the image of $R^+$.

In particular, the map $\mathcal M_Z\to S$ is representable in locally spatial diamonds and compactifiable, of locally finite $\dimtrg$.
\end{proposition}

\begin{proof} Choose an open subspace $W\subset \mathbb P^n_{X_S}$ such that $Z$ is Zariski closed in $W$. For any perfectoid space $T$ with a map $T\to \mathcal M_{\mathbb P^n_{X_S}}$ corresponding to a map $X_T\to \mathbb P^n_{X_S}$ over $X_S$, the locus $T_W\subset T$ where the section factors over $W$ is open. Indeed, this locus is the complement of the image in $|T|$ of the preimage of $|\mathbb P^n_{X_S}\setminus W|$ under $|X_T|\to |\mathbb P^n_{X_S}|$, and $|X_T|\to |T|$ is closed.

Replacing $T$ by $T_W$, we can assume that the section $X_T\to \mathbb P^n_{X_S}$ factors over $W$. We may also assume that $T=\Spa(R,R^+)$ is affinoid perfectoid and that $S=T$. Pick a pseudouniformizer $\varpi\in R$, in particular defining the cover
\[
Y_{S,[1,q]} = \{|\pi |^q\leq |[\varpi]|\leq |\pi|\}\subset \Spa W_{\mathcal O_E}(R^+)
\]
of $X_S$. The pullback of the line bundle $\mathcal O_{\mathbb P^n}(1)$ to $X_S$ along this section, and then to $Y_{S,[1,q]}$, is \'etale locally trivial, as when $S$ is a geometric point, $Y_{S,[1,q]}$ is affinoid with ring of functions a principal ideal domain by Corollary~\ref{cor:principalidealdomains}. Replacing $W$ by a small \'etale neighborhood of this section and correspondingly shrinking $S$, we can assume that the pullback of $\mathcal O_{\mathbb P^n}(1)$ to $W_{[1,q]} = W\times_{X_S} Y_{S,[1,q]}$ is trivial. In that case the pullback $Z_{[1,q]}\to Y_{S,[1,q]}$ of $Z\to X_S$ is Zariski closed in an open subset of
\[
\mathbb A^{n+1}_{Y_{S,[1,q]}}.
\]
Inside $\mathbb A^{n+1}_{Y_{S,[1,q]}}$, the image of $Y_{\Spa(K,K^+),[1,q]}$ (via the given section) for a point $\Spa(K,K^+)\to S$ is an intersection of small balls over $Y_{S',[1,q]}$ for small neighborhoods $S'\subset S$ of $\Spa(K,K^+)$. Thus, one of these balls is contained in the open subset of which $Z_{[1,q]}$ is a Zariski closed subset. Thus, after this further localization, we can assume that there is a Zariski closed immersion
\[
Z_{[1,q]}\hookrightarrow \mathbb B^{n+1}_{Y_{S,[1,q]}},
\]
and in particular $Z_{[1,q]}$ is affinoid and cut out by global functions on $\mathbb B^{n+1}_{Y_{S,[1,q]}}$ by Proposition~\ref{prop:zariskiclosedsmooth}. Pulling back these functions along the given section $Y_{S,[1,q]}\to \mathbb B^{n+1}_{Y_{S,[1,q]}}$, it suffices to see that if $S=\Spa(R,R^+)$ is an affinoid perfectoid space of characteristic $p$ with a choice of pseudouniformizer $\varpi\in R$ and
\[
f\in B_{R,[1,q]} = \mathcal O(Y_{S,[1,q]})
\]
is a function, then there is a universal perfectoid space $S'\subset S$ for which the pullback of $f$ is zero, and $S'\subset S$ is Zariski closed. This is given by Lemma \ref{lemma:ben oui quoi}.
\end{proof}

\begin{lemma}\label{lemma:ben oui quoi}
Let $S=\Spa (R,R^+)\in \Perf_{\Fq}$ be affinoid perfectoid with a fixed pseudo-uniformizer $\varpi$, $I\subset (0,\infty)$ a compact interval with rational ends, and $Z\subset |Y_{S,I}|$  a closed subset defined by the vanishing locus of an ideal $J\subset \O(Y_{S,I})$. Then, via the open projection $\upsilon:|Y_{S,I}|\to |S|$, the closed subset $|S|\setminus \upsilon (|Y_{S,I}|\setminus Z)$ is Zariski closed. The corresponding Zariski closed perfectoid subspace of $S$ is universal for perfectoid spaces $T\to S$ such that $J\mapsto 0$ via $\O (Y_{S,I})\to \O(Y_{T,I})$. 
\end{lemma}

\begin{proof}
Since $Y_{S,I}^\diamond\to S$ is cohomologically smooth, $\upsilon$ is open. We can suppose $J=(f)$ with $f\in \O(Y_{S,I})$.
For any untilt of $\mathbb F_q\laurentseries{\varpi^{1/p^\infty}}$ over $E$ such that $|\pi|^b\leq |\varpi^\sharp|\leq |\pi|^a$ if $I=[a,b]$, we get a corresponding untilt $R^\sharp$ of $R$ over $E$, with a map $B_{R,I}\to R^\sharp$. The locus where the image of $f$ in $R^\sharp$ vanishes is Zariski closed by Proposition~\ref{prop:zariskiclosedstronglyzariskiclosed}. Intersecting these Zariski closed subsets over varying such untilts gives the vanishing locus of $f$, as in any fibre, $f$ vanishes as soon at it vanishes at infinitely many untilts (e.g., by Corollary~\ref{cor:principalidealdomains}), and all rings are sous-perfectoid, in particular uniform, so vanishing at all points implies vanishing.
\end{proof}

\subsection{Formal smoothness of $\mathcal{M}_Z^{\mathrm{sm}}$}

The key result we need is the following.

\begin{proposition}\label{prop:neighborhoodretract} Let $S=\Spa(R,R^+)$ be an affinoid perfectoid space over $\Fq$ and let $Z\to X_S$ be a smooth map of sous-perfectoid adic spaces that is Zariski closed in an open subspace of $\mathbb P^n_{X_S}$. Then $\mathcal M_Z^{\mathrm{sm}}\to S$ is formally smooth.
\end{proposition}

\begin{proof} Pick a test diagram as in Definition~\ref{def:formalsmoothness}; we can and do assume that the $S$ from there is the given $S$, replacing the $S$ in this proposition if necessary. 
This means we have a diagram 
$$
\begin{tikzcd}
& Z\ar[d] \\
X_{S_0}\ar[r,hook] \ar[ru,"s_0"] & X_S \ar[u,bend right,dotted]
\end{tikzcd}
$$
and, up to replacing $S$ by an \'etale neighborhood of $S_0$ we try to extend the section $s_0$ to a section over $X_S$ (the dotted line in the diagram).
Fix a geometric point $\Spa(C,C^+)\to S_0$; we will always allow ourselves to pass to \'etale neighborhoods of this point. Fix a pseudouniformizer $\varpi\in R$ and consider the affinoid cover $Y_{S,[1,q]}\to X_S$; recall that
\[
Y_{S,[1,q]}=\{|\pi|^q\leq |[\varpi]|\leq |\pi|\}\subset \Spa W_{\mathcal O_E}(R^+)
\]
and we also consider its boundary annuli
\[
Y_{S,[1,1]} = \{|[\varpi]|=|\pi|\}, Y_{S,[q,q]} = \{|\pi|^q=|[\varpi]|\}\subset Y_{S,[1,q]}.
\]
Let $Z_{[1,q]}\to Z$ be its pullback; with pullback $Z_{[1,1]},Z_{[q,q]}\subset Z_{[1,q]}$ of $Y_{S,[1,1]},Y_{S,[q,q]}\subset Y_{S,[1,q]}$. In particular, $Z$ is obtained from $Z_{[1,q]}$ via identification of its open subsets $Z_{[1,1]}, Z_{[q,q]}$ along the isomorphism $\phi: Z_{[1,1]}\to Z_{[q,q]}$.

Arguing as in the proof of Proposition~\ref{prop:embedZ}, we can after \'etale localization on $S$ embed
\[
Z_{[1,q]}\hookrightarrow \mathbb B^{n+1}_{Y_{S,[1,q]}}
\]
as a Zariski closed subset. We thus have a diagram
$$
\begin{tikzcd}[column sep=huge]
Z_{[1,q]} \ar[r,hook,"\text{Zariski closed}"] & \mathbb{B}^{n+1}_{Y_{S,[1,q]}} \ar[d] \\
Y_{S_0,[1,q]} \ar[r,hook]\ar[u]  & Y_{S,[1,q]}.
\end{tikzcd}
$$
In particular, $Z_{[1,q]}$ is affinoid.

Next, consider the K\"ahler differentials $\Omega^1_{Z_{[1,q]}/Y_{S,[1,q]}}$. Again, as $B_{C,[1,q]}$ is a principal ideal domain, its restriction to the section $Y_{\Spa(C,C^+),[1,q]}\subset Z_{[1,q]}$ is trivial, and thus it is trivial in a small neighborhood. It follows that after a further \'etale localization we can assume that $\Omega^1_{Z_{[1,q]}/Y_{S,[1,q]}}\cong \mathcal O_{Z_{[1,q]}}^r$ is trivial. On the Zariski closed subset $Z_{0,[1,q]}\subset Z_{[1,q]}$ (defined as the pullback of $Y_{S_0}\subset Y_S$), this implies that we may find functions $f_1,\ldots,f_r\in \mathcal O(Z_{0,[1,q]})$ vanishing on the section $Y_{S_0,[1,q]}\to Z_{0,[1,q]}$ and locally generating the ideal of this closed immersion  (use Proposition \ref{prop:zariskiclosedsmooth}). In particular,
\[
df_1,\ldots,df_r\in \Omega^1_{Z_{0,[1,q]}/Y_{S_0,[1,q]}}
\]
are generators at the image of the section $Y_{S_0,[1,q]}\to Z_{0,[1,q]}$, and thus in an open neighborhood. Picking lifts of the $f_i$ to $\mathcal O(Z_{[1,q]})$ and shrinking $Z_{[1,q]}$, Proposition~\ref{prop:smoothmapsbetweensmooth} implies that they define an \'etale map
\[
Z_{[1,q]}\to \mathbb B^r_{Y_{S,[1,q]}}.
\]
Moreover, over $\{0\}_{Y_{S_0,[1,q]}}\subset \mathbb B^r_{Y_{S,[1,q]}}$, this map admits a section. Shrinking further around this section, we can thus arrange that there are open immersions
\[
(\pi^N \mathbb B)^r_{Y_{S,[1,q]}}\subset Z_{[1,q]}\subset \mathbb B^r_{Y_{S,[1,q]}},
\]
and that the section over $Y_{S_0,[1,q]}$ is given by the zero section.

The isomorphism $\phi: Z_{[1,1]}\to Z_{[q,q]}$ induces a map
\[
\phi': (\pi^N \mathbb B)^r_{Y_{S,[1,1]}}\to \mathbb B^r_{Y_{S,[q,q]}}.
\]
Recall that for any compact interval $I\subset (0,\infty)$, the space
\[
Y_{S,I} = \Spa(B_{R,I},B_{(R,R^+),I}^+)
\]
is affinoid. The map $\phi'$ is then given by a map
\[
\alpha: B_{R,[q,q]}\langle T_1,\ldots,T_r\rangle\to B_{R,[1,1]}\langle \pi^{-N} T_1,\ldots,\pi^{-N} T_r\rangle
\]
linear over the isomorphism $\phi: B_{R,[q,q]}\to B_{R,[1,1]}$. The map $\alpha$ is determined by the images of $T_1,\ldots,T_n$ which are elements
\[
\alpha_i\in B_{(R,R^+),[1,1]}^+\langle \pi^{-N} T_1,\ldots,\pi^{-N} T_r\rangle.
\]
These have the property that on the quotient $B_{R_0,[1,1]}$ they vanish at $T_1=\ldots=T_r=0$ (as over $S_0$, the zero section is $\phi$-invariant). Moreover, over the geometric point $\Spa(C,C^+)\to S_0$ fixed at the beginning, we can apply a linear change of coordinates in order to ensure that the derivative at the origin is given by a standard matrix for an isocrystal of negative slopes; i.e., there are cycles $1,\ldots,r_1$; $r_1+1,\ldots,r_2$; $\ldots$; $r_{a-1}+1,\ldots,r_a=r$ and positive integers $d_1,\ldots,d_a$ such that
\[
\alpha_i\equiv T_{i+1}\ \mathrm{in}\ B_{C,[1,1]}[T_1,\ldots,T_r]/(T_1,\ldots,T_r)^2
\]
if $i\neq r_j$ for some $j=1,\ldots,a$, and
\[
\alpha_{r_j}\equiv \pi^{-d_j} T_{r_{j-1}+1}\ \mathrm{in}\ B_{C,[1,1]}[T_1,\ldots,T_r]/(T_1,\ldots,T_r)^2.
\]
(Here, we set $r_0=0$.) Approximating this linear change of basis over an \'etale neighborhood, we respect the condition that the $\alpha_i$'s vanish at $T_1=\ldots=T_r=0$ over $S_0$, while we can for any large $M$ arrange
\[
\alpha_i\equiv T_{i+1}\ \mathrm{in}\ B_{(R,R^+),[1,1]}^+/\pi^M[\pi^{-N} T_1,\ldots,\pi^{-N}T_r]/(\pi^{-N}T_1,\ldots,\pi^{-N}T_r)^2
\]
if $i\neq r_j$ and
\[
\alpha_{r_j}\equiv \pi^{-d_j} T_{r_{j-1}+1}\ \mathrm{in}\ B_{(R,R^+),[1,1]}^+/\pi^M[\pi^{-N} T_1,\ldots,\pi^{-N}T_r]/(\pi^{-N}T_1,\ldots,\pi^{-N}T_r)^2.
\]
Moreover, rescaling all $T_i$ by powers of $\pi$, and passing to a smaller neighborhood around $S_0$, we can then even ensure that
\[
\alpha_i\in T_{i+1} + \pi^M B_{(R,R^+),[1,1]}^+\langle T_1,\ldots,T_r\rangle
\]
for $i\neq r_j$ and
\[
\alpha_{r_j}\in \pi^{-d_j} T_{r_{j-1}+1} + \pi^M B_{(R,R^+),[1,1]}^+\langle T_1,\ldots,T_r\rangle.
\]
At this point, the integers $d_1,\ldots,d_a$ are fixed, while we allow ourselves to choose $M$ later, depending only on these.

From this point on, we will no longer change $S$ and $S_0$, and instead will merely change coordinates in the balls (by automorphisms). More precisely, we study the effect of replacing $T_i$ by $T_i+\epsilon_i$ for some
\[
\epsilon_i\in \pi^d \ker(B_{(R,R^+),[1,q]}^+\to B_{R_0,[1,q]})
\]
where we take $d$ to be at least the maximum of all $d_j$. This replaces $\alpha_i$ by a new power series $\alpha_i'$, given by
\[
\alpha_i'(T_1,\ldots,T_r) = \alpha_i(T_1,\ldots,T_i + \epsilon_i,\ldots,T_r) - \phi(\epsilon_i)
\]
and the $\alpha_i'$'s still vanish at $T_1=\ldots=T_r=0$ over $S_0$. Their nonconstant coefficients will still have the same properties as for $\alpha_i$ (the linear coefficients are unchanged, while all other coefficients are divisible by $\pi^M$), and the constant coefficient satisfies
\[
\alpha_i'(0,\ldots,0)\equiv \alpha_i(0,\ldots,0) + \epsilon_{i+1} - \phi(\epsilon_i)\ \mathrm{in}\ B_{(R,R^+),[1,1]}^+/\pi^{M+d}
\]
if $i\neq r_j$ and
\[
\alpha_{r_j}'(0,\ldots,0)\equiv \alpha_{r_j}(0,\ldots,0) + \pi^{-d_j} \epsilon_{r_{j-1}+1} - \phi(\epsilon_{r_j})\ \mathrm{in}\ B_{(R,R^+),[1,1]}^+/\pi^{M+d}.
\]
Assume that by some inductive procedure we already achieved $\alpha_i(0,\ldots,0)\in \pi^{N'} B_{(R,R^+),[1,1]}^+$ for some $N'\geq M$. By Lemma~\ref{lem:quantitativeh1vanishing} below, there is some constant $c$ depending only on $d_1,\ldots,d_a$ such that we can then find $\epsilon_i\in \pi^{N'-c} B_{(R,R^+),[1,q]}^+$, vanishing over $R_0$, with
\[
\alpha_i(0,\ldots,0)= \phi(\epsilon_i) - \epsilon_{i+1}
\]
for $i\neq r_j$ and
\[
\alpha_{r_j}(0,\ldots,0) = \phi(\epsilon_{r_j}) - \pi^{-d_j} \epsilon_{r_{j-1}+1}.
\]
This means that $\alpha_i'(0,\ldots,0)\in \pi^{M+N'-c} B_{(R,R^+),[1,1]}^+$, so if we choose $M>c$ in the beginning (which we can), then this inductive procedure converges, and in the limit we get a change of basis after which the zero section defines a $\phi$-invariant section of $Z_{[1,q]}$, thus a section $s: X_S\to Z$, as desired. Note that we arranged that this section agrees with $s_0$ over $S_0$, as all coordinate changes did not affect the situation over $S_0$.
\end{proof}

We used the following quantitative version of vanishing of $H^1(X_S,\mathcal E)$ for $\mathcal E$ of positive slopes.

\begin{lemma}\label{lem:quantitativeh1vanishing} Fix a standard Dieudonn\'e module of negative slopes, given explicitly on a basis $e_1,\ldots,e_r$ by fixing cycles $1,\ldots,r_1$; $r_1+1\ldots,r_2$; $\ldots$; $r_{a-1}+1,\ldots,r_a=r$ and positive integers $d_1,\ldots,d_a>0$, via
\[
\phi(e_i)=e_{i+1}\ \mathrm{for}\ i\neq r_j,\ \phi(e_{r_j}) = \pi^{-d_j} e_{r_{j-1}+1}.
\]
Then there is an integer $c\geq 0$ with the following property.

Let $S=\Spa(R,R^+)$ be an affinoid perfectoid space over $\Fq$ with Zariski closed subspace $S_0=\Spa(R_0,R_0^+)$, and a pseudouniformizer $\varpi\in R$. Let
\[
I_{[1,q]}^+=\ker(B_{(R,R^+),[1,q]}^+\to B_{(R_0,R_0^+),[1,q]}^+),\ I_{[1,1]}^+ = \ker(B_{(R,R^+),[1,1]}^+\to B_{(R_0,R_0^+),[1,1]}^+).
\]
Then for all $f_1,\ldots,f_r\in I_{[1,1]}^+$ one can find $g_1,\ldots,g_r\in \pi^{-c} I_{[1,q]}^+$ such that
\[
f_i = \phi(g_i) - g_{i+1}\ \mathrm{for}\ i\neq r_j,\ f_{r_j} = \phi(g_{r_j}) - \pi^{-d_j} g_{r_{j-1}+1}.
\]
\end{lemma}

\begin{proof} We may evidently assume that $a=1$; set $d=d_1$. By linearity, we can assume that all but one of the $f_i$'s is equal to zero, say (by cyclic rotation) $f_1=\ldots=f_{r-1}=0$. Thus, it suffices to see that for all positive integers $r$ and $d$ there is $c\geq 0$ such that for all $f=f_r\in I_{[1,1]}^+$ one can find some $g\in \pi^{-c} I_{[1,q^r]}^+$ (for the evident definition of $I_{[1,q^r]}^+$) such that
\[
f = \phi^r(g) - \pi^{-d} g.
\]
Indeed, one then takes $g_1=g,g_2=\phi(g),\ldots,g_r=\phi^{r-1}(g)$. Replacing $E$ by its unramified extension of degree $r$, we can then assume that $r=1$. At this point, we want to reduce to the qualitative version given by Lemma~\ref{lem:qualitativeh1vanishing} below, saying that the map
\[
\phi - \pi^{-d}: I_{[1,q]}\to I_{[1,1]}
\]
is surjective. Indeed, assume a constant $c$ as desired would not exist. Then for any integer $i\geq 0$ we can find some Zariski closed immersion $S_{0,i}=\Spa(R_{0,i},R_{0,i}^+)\subset S_i=\Spa(R_i,R_i^+)$, with choices of pseudouniformizers $\varpi_i\in R_i$, as well as elements $f_i\in I_{[1,1],i}^+$ such that there is no $g\in \pi^{-2i} I_{[1,q],i}^+$ with $f_i = \phi(g_i) - \pi^{-d} g_i$. Then we can define $R^+=\prod_i R_i^+$ with $\varpi=(\varpi_i)_i\in R^+$, and $R=R^+[\tfrac 1{\varpi}]$, which defines an affinoid perfectoid space $S=\Spa(R,R^+)$, containing a Zariski closed subspace $S_0\subset S$ defined similarly. Moreover, the sequence $(\pi^i f_i)_i$ defines an element of $f\in I_{[1,1]}^+$. As $\phi-\pi^{-d}: I_{[1,q]}\to I_{[1,1]}$ is surjective by Lemma~\ref{lem:qualitativeh1vanishing}, we can find some $g\in I_{[1,q]}$ with $f=\phi(g)-\pi^{-d} g$. Then $\pi^c g\in I_{[1,q]}^+$ for some $c$, and restricting $g$ to $S_{0,i}\subset S_i$ with $i>c$ gives the desired contradiction.
\end{proof}

We reduced to the following qualitative version.

\begin{lemma}\label{lem:qualitativeh1vanishing} Let $d$ be a positive integer, let $S_0=\Spa(R_0,R_0^+)\subset S=\Spa(R,R^+)$ be a Zariski closed immersion of affinoid perfectoid spaces over $\Fq$, and let $\varpi\in R$ be a pseudouniformizer. Let
\[
I_{[1,q]} = \ker(B_{R,[1,q]}\to B_{R_0,[1,q]})\ ,\ I_{[1,1]} = \ker(B_{R,[1,1]}\to B_{R_0,[1,1]}).
\]
Then the map
\[
\phi - \pi^{-d}: I_{[1,q]}\to I_{[1,1]}
\]
is surjective.
\end{lemma}

\begin{proof} By the snake lemma and the vanishing $H^1(X_S,\mathcal O_{X_S}(d))=0$ (Proposition~\ref{prop:standardbanachcolmez}~(iii)), the lemma is equivalent to the surjectivity of
\[
H^0(X_S,\mathcal O_{X_S}(d))\to H^0(X_{S_0},\mathcal O_{X_{S_0}}(d)).
\]
For $d\leq [E:\mathbb Q_p]$ (or if $E$ is of equal characteristic), this follows directly from Proposition~\ref{prop:standardbanachcolmez}~(iv) and the surjectivity of $R^{\circ\circ}\to R_0^{\circ\circ}$. In general, we can either note that the proof of Proposition~\ref{prop:standardbanachcolmez}~(iii) also proves the lemma, or argue by induction by choosing an exact sequence
\[
0\to \mathcal O_{X_S}(d-2)\to \mathcal O_{X_S}(d-1)^2\to \mathcal O_{X_S}(d)\to 0
\]
(the Koszul complex for two linearly independent sections of $H^0(X_S,\mathcal O_{X_S}(1))$), and use the vanishing of $H^1(X_{S_0},\mathcal O_{X_{S_0}}(d-2))=0$ for $d>2$, Proposition~\ref{prop:standardbanachcolmez}~(iii). This induction gets started as long as $E\neq \mathbb Q_p$. For $E=\mathbb Q_p$, we can write $\mathcal O_{X_S}(d)$ as a direct summand of $\pi_\ast \pi^\ast \mathcal O_{X_S}(d)$ for any extension $\pi: X_{S,E}\to X_S$ with $E\neq \mathbb Q_p$.
\end{proof}

\subsection{Universal local acyclicity of $\mathcal{M}_Z^{\mathrm{sm}}\to S$}

The next step in the proof of Theorem~\ref{thm:jacobiancriterion} is to show that $\mathbb F_\ell$ is universally locally acyclic.

\begin{proposition}\label{prop:jacobiancriterionULA} Let $S$ be a perfectoid space and let $Z\to X_S$ be a smooth map of sous-perfectoid spaces such that $Z$ is Zariski closed inside an open subset of $\mathbb P^n_{X_S}$ for some $n\geq 0$. Then, for any $\ell\neq p$, the sheaf $\mathbb F_\ell$ is universally locally acyclic for the map
\[
\mathcal M_Z^{\mathrm{sm}}\to S.
\]
\end{proposition}

\begin{proof} Recall from Proposition~\ref{prop:embedZ} that $\mathcal M_Z\to \mathcal M_{\mathbb P^n_{X_S}}$ is a locally closed immersion, and the open embedding
\[
\mathcal M_{\mathbb P^n_{X_S}}\hookrightarrow \bigsqcup_{m\geq 0} (\BC(\mathcal O_{X_S}(m)^{n+1})\setminus \{0\})/\underline{E^\times}
\]
from Proposition~\ref{prop:gromovwittenPn}. In the following, we fix some $m$ and work on the preimage of
\[
(\BC(\mathcal O_{X_S}(m)^{n+1})\setminus \{0\})/\underline{E^\times}.
\]

We choose a surjection $g: T\to \BC(\mathcal O_{X_S}(m)^{n+1}\setminus\{0\})/\underline{E^\times}$ from a perfectoid space $T$ as in Lemma~\ref{lem:nicesurjection}; in particular, $g$ is separated, representable in locally spatial diamonds, cohomologically smooth, and formally smooth. Moreover, locally $T$ admits a Zariski closed immersion into the perfectoid ball $\tilde{\mathbb B}^n_S$ over $S$. Taking the pullback of $T$ to $\mathcal M_Z$, we get a surjection $T_Z\to \mathcal M_Z$ for some perfectoid space $T_Z$ such that \'etale locally $T_Z$ admits a Zariski closed immersion into a space \'etale over $\tilde{\mathbb B}^n_S$.

It follows that one can cover $\mathcal M_Z$ via maps $h_0: T_0\to \mathcal M_Z$ that are separated, representable in locally spatial diamonds, cohomologically smooth, and formally smooth, and such that $T_0$ admits a Zariski closed immersion into some space \'etale over $\tilde{\mathbb B}^n_S$. We can then also do the same for $\mathcal M_Z^{\mathrm{sm}}$. By Proposition~\ref{prop:neighborhoodretract}, we can, up to further replacement of $T_0$ by an \'etale cover, assume that the map $h_0$ extends to a map $h: T\to \mathcal M_Z^{\mathrm{sm}}$ for some perfectoid space $T$ \'etale over $\tilde{\mathbb B}^n_S$. Moreover, as $T_0\to \mathcal M_Z^{\mathrm{sm}}$ is formally smooth, we can, after a further \'etale localization, lift the map $T\to \mathcal M_Z^{\mathrm{sm}}$ to a retraction $T\to T_0$; thus, $T_0$ is a retract of a space that is \'etale over a perfectoid ball. Now the result follows from Corollary~\ref{cor:ULAretract}.
\end{proof}

We used the following presentation of certain projectivized Banach--Colmez spaces.

\begin{lemma}\label{lem:nicesurjection} Let $S$ be a perfectoid space over $\Fq$ and let $\mathcal E$ be a vector bundle on $X_S$ that is everywhere of nonnegative Harder--Narasimhan slopes. There is a perfectoid space $T\to S$ that is locally Zariski closed in a perfectoid ball $\tilde{\mathbb B}^n_S$ over $S$ and that admits a surjective map
\[
T\to (\BC(\mathcal E)\setminus \{0\})/\underline{E^\times}
\]
over $S$ that is separated, representable in locally spatial diamonds, cohomologically smooth, and formally smooth.
\end{lemma}

\begin{proof} The target parametrizes line bundles $\mathcal L$ on $X_S$ of slope zero together with a section of $\mathcal E\otimes \mathcal L$ that is nonzero fibrewise on $S$. Parametrizing in addition an injection $\mathcal L\hookrightarrow \mathcal O_{X_S}(1)$ defines a map that is separated, representable in locally spatial diamonds, cohomologically smooth, and formally smooth (by Proposition~\ref{prop:relativebanachcolmez1} and Proposition~\ref{prop:formal smoothness of BCs}). Over this cover, one has locally on $S$ an untilt $S^\sharp$ over $E$ corresponding to the support of the cokernel of $\mathcal L\to \mathcal O_{X_S}(1)$, and one parametrizes nonzero sections of $\mathcal E(1)$ that vanish at $S^\sharp\hookrightarrow X_S$. This is Zariski closed (by \cite[Theorem 7.4, Remark 7.5]{BhattScholzePrism}) inside the space of all sections of $\mathcal E(1)$. We see that it suffices to prove the similar result with $(\BC(\mathcal E)\setminus \{0\})/\underline{E^\times}$ replaced by $\BC(\mathcal E(1))\times \Spd E$, and this reduces to the individual factors. For $\BC(\mathcal E(1))$, the result follows from the argument in Proposition~\ref{prop:formal smoothness of BCs}. For $\Spd E$, there is nothing to do in equal characteristic, so assume that $E$ is $p$-adic. Then we reduce to $[\ast/\underline{\mathcal O_E^\times}]$ as the fibres of $\Spd E\to [\ast/\underline{\mathcal O_E^\times}]$ over perfectoid spaces are given by $\BC(\mathcal L)\setminus \{0\}$ for some line bundle $\mathcal L$ of slope $1$, and this in turn admits covers of the desired form. Finally, for $[\ast/\underline{\mathcal O_E^\times}]$, we can pass to the \'etale cover $[\ast/\underline{1+p^2\mathcal O_E}]\cong [\ast/\underline{\mathcal O_E}]$, or to $[\ast/\underline{E}]$. This, finally, admits a surjection from a perfectoid open unit disc $\BC(\mathcal O_{X_S}(1))$ with the desired properties by passing to Banach--Colmez spaces in an exact sequence
\[
0\to \mathcal O_{X_S}\to \mathcal O_{X_S}(\tfrac 12)\to \mathcal O_{X_S}(1)\to 0
\]
and using Proposition~\ref{prop:h1loczero}.
\end{proof}

\subsection{Deformation to the normal cone}

The final step in the proof of Theorem~\ref{thm:jacobiancriterion} is a deformation to the normal cone.

By Proposition~\ref{prop:embedZ} and Proposition~\ref{prop:jacobiancriterionULA} (and Proposition~\ref{prop:ULAcohomsmooth}), in order to prove Theorem~\ref{thm:jacobiancriterion} it only remains to prove that $Rf^! \mathbb F_\ell$ is invertible and sitting in the expected cohomological degree. Picking a v-cover $T\to \mathcal M_Z^{\mathrm{sm}}$ by some perfectoid space $T$ and using that the formation of $Rf^! \mathbb F_\ell$ commutes with any base change by Proposition~\ref{prop:ULAdual}, it suffices to prove the following result.

\begin{proposition}\label{prop:smoothatsection} Let $S$ be a perfectoid space and let $Z\to X_S$ be a smooth map of sous-perfectoid spaces such that $Z$ admits a Zariski closed immersion into an open subset of (the adic space) $\mathbb P^n_{X_S}$ for some $n\geq 0$. Let $f: \mathcal M_Z\to S$ be the moduli space of sections of $Z\to X_S$. Moreover, let $s: X_S\to Z$ be a section such that $s^\ast T_{Z/X_S}$ is everywhere of positive Harder--Narasimhan slopes, and of degree $d$.

Let $t: S\to \mathcal M_Z$ be the section of $f$ corresponding to $s$. Then $t^\ast Rf^! \mathbb F_\ell$ is \'etale locally on $S$ isomorphic to $\mathbb F_\ell[2d]$.
\end{proposition}

\begin{proof} We will prove this by deformation to the normal cone. In order to avoid a general discussion of blow-ups etc., we will instead take an approach based on the local structure of $Z$ near a section as exhibited in the proof of Proposition~\ref{prop:neighborhoodretract}.

We are free to make v-localizations on $S$ (as being cohomologically smooth can be checked after a v-cover), and replace $Z$ by an open neighborhood of $s(X_S)$. With this freedom, we can follow the proof of Proposition~\ref{prop:neighborhoodretract} and ensure that $S=\Spa(R,R^+)$ is strictly totally disconnected with pseudouniformizer $\varpi$, the pullback $Z_{[1,q]}\to Y_{S,[1,q]}$ of $Z\to X_S$ to
\[
Y_{S,[1,q]} = \{|\pi|^q\leq |[\varpi]|\leq |\pi|\}\subset \Spa W_{\mathcal O_E}(R^+)
\]
satisfies
\[
\pi^N \mathbb B^r_{Y_{S,[1,q]}}\subset Z_{[1,q]}\subset \mathbb B^r_{Y_{S,[1,q]}}
\]
and the gluing isomorphism is given by power series
\[
\alpha_i\in T_{i+1} + \pi^M B_{(R,R^+),[1,1]}^+\langle T_1,\ldots,T_r\rangle
\]
resp.
\[
\alpha_{r_j}\in \pi^{-d_j}T^{r_{j-1}+1}+\pi^M B_{(R,R^+),[1,1]}^+\langle T_1,\ldots,T_r\rangle
\]
with notation following the proof of Proposition~\ref{prop:neighborhoodretract}. Moreover, the constant coefficients of all $\alpha_i$ vanish. These in fact define a map
\[
\phi: \pi^d\mathbb B^r_{Y_{S,[1,1]}}\to \mathbb B^r_{Y_{S,[q,q]}}
\]
preserving the origin, where $d$ is the maximum of the $d_j$.

For any $n\geq N,d$, we can look at the subset
\[
Z_{[1,q]}^{(n)} = \pi^n \mathbb B^r_{Y_{S,[1,q]}}\cup \phi(\pi^n \mathbb B^r_{Y_{S,[1,1]}}) \subset Z_{[1,q]},
\]
which descends to an open subset $Z^{(n)}\subset Z$. Letting $T_i^{(n)} = \pi^{-n} T_i$, the gluing is then given by power series $\alpha_i^{(n)}$ given by
\[
\alpha_i^{(n)} = \pi^{-n}\alpha_i(\pi^n T_1,\ldots,\pi^n T_n)
\]
which satisfy the same conditions, but the nonlinear coefficients of $\alpha_i^{(n)}$ become more divisible by $\pi$. The limit
\[
\alpha_i^{(\infty)} = \varinjlim_{n\to \infty} \alpha_i^{(n)} \in B_{R,[1,1]}\langle T_1,\ldots,T_r\rangle
\]
exists, and is linear in the $T_i$.

Let $S' = S\times \underline{\mathbb N_{\geq N}\cup\{\infty\}}$, using the profinite set $\mathbb N_{\geq N}\cup\{\infty\}$. Let $Z'\to X_{S'}$ be the smooth map of sous-perfectoid spaces obtained by descending
\[
Z'_{[1,q]} = \mathbb B^r_{Y_{S',[1,q]}}\cup \phi'(\mathbb B^r_{Y_{S',[q,q]}})
\]
along the isomorphism $\phi'$ given by the power series
\[
\alpha_i' = (\alpha_i^{(N)},\alpha_i^{(N+1)},\ldots,\alpha_i^{(\infty)})\in B_{R',[1,1]}\langle T_1,\ldots,T_r\rangle.
\]
Then the fibre of $Z'\to X_{S'}$ over $S\times\{n\}$ is given by $Z^{(n)}$, while its fibre over $S\times \{\infty\}$ is given by an open subset $Z^{(\infty)}$ of the Banach--Colmez space $\BC(s^\ast T_{Z/X_S})$. Moreover, letting $S'^{(>N)}\subset S'$ be the complement of $S\times\{N\}$, there is natural isomorphism $\gamma: S'^{(>N)}\to S'$ given by the shift $S\times\{n+1\}\cong S\times\{n\}$, and this lifts to an open immersion
\[
\gamma: Z'^{(>N)} = Z'\times_{X_{S'}} X_{S'^{(>N)}}\hookrightarrow Z'.
\]

We need to check that $Z'\to X_{S'}$ still satisfies the relevant quasiprojectivity assumption.

\begin{lemma}\label{lem:deformationquasiprojective} The space $Z'\to X_{S'}$ admits a Zariski closed immersion into an open subset of $\mathbb P^m_{X_{S'}}$ for some $m\geq 0$.
\end{lemma}

\begin{proof} One may perform a parallel construction with $Z$ replaced by an open subset of $\mathbb P^m_{X_S}$, reducing us to the case that $Z$ is open in $\mathbb P^m_{X_S}$. In that case, the key observation is that the blow-up of $\mathbb P^m_{X_S}$ at the section $s: X_S\to \mathbb P^m_{X_S}$ is still projective, which is an easy consequence of $X_S$ admitting enough line bundles.
\end{proof}

Let $f': \mathcal M_{Z'}\to S'$ be the projection, with fibres $f^{(n)}$ and $f^{(\infty)}$. By Proposition~\ref{prop:jacobiancriterionULA}, both $\mathbb F_\ell$ and $Rf'^! \mathbb F_\ell$ are $f'$-universally locally acyclic. In particular, the formation of $Rf'^!\mathbb F_\ell$ commutes with base change, and we see that the restriction of $Rf'^!\mathbb F_\ell$ to the fibre over $\infty$ is \'etale locally isomorphic to $\mathbb F_\ell[2d]$, as an open subset of $\BC(s^\ast T_{Z/X_S})$. As $S$ is strictly totally disonnected, one can choose a global isomorphism with $\mathbb F_\ell[2d]$.

The map from $\mathbb F_\ell[2d]$ to the fibre of $Rf'^! \mathbb F_\ell$ over $\infty$ extends to a small neighborhood; passing to this small neighborhood, we can assume that there is a map
\[
\beta: \mathbb F_\ell[2d]\to Rf'^! \mathbb F_\ell
\]
that is an isomorphism in the fibre over $\infty$. We can assume that this map is $\gamma$-equivariant (passing to a smaller neighborhood). Let $Q$ be the cone of $\beta$. Then $Q$ is still $f'$-universally locally acyclic, as is its Verdier dual
\[
\mathbb D_{\mathcal M_{Z'}/S'}(Q) = R\sHom_{\mathcal M_{Z'}}(Q,Rf'^! \mathbb F_\ell).
\]
In particular, $Rf'_!\mathbb D_{\mathcal M_{Z'}/S'}(Q)\in D_\et(S',\mathbb F_\ell)$ is constructible, and its restriction to $S\times \{\infty\}$ is trivial. This implies (e.g.~by \cite[Proposition 20.7]{ECoD}) that its restriction to $S\times \{n,n+1,\ldots,\infty\}$ is trivial for some $n\gg 0$. Passing to this subset, we can assume that $Rf'_!\mathbb D_{\mathcal M_{Z'}/S'}(Q)=0$. Taking Verdier duals and using Corollary~\ref{cor:ULAselfdual}, this implies that $Rf'_\ast Q=0$.

In particular, for all $n\geq n_0$, we have $Rf^{(n)}_\ast Q|_{\mathcal M_{Z^{(n)}}}=0$. Using the $\gamma$-equivariance, this is equivalent to
\[
Rf^{(n)}_\ast (Q|_{\mathcal M_Z^{(n_0)}})|_{\mathcal M_{Z^{(n)}}} = 0,
\]
regarding $\mathcal M_{Z^{(n)}}\subset \mathcal M_{Z^{(n_0)}}$ as an open subset. Taking the colimit over all $n$ and using that the system $\mathcal M_{Z^{(n)}}\subset \mathcal M_{Z^{(n_0)}}$ has intersection $s(S)\subset \mathcal M_Z$ and is cofinal with a system of spatial diamonds of finite cohomological dimension (as can be checked in the case of projective space), \cite[Proposition 14.9]{ECoD} implies that
\[
s^\ast Q|_{\mathcal M_{Z^{(n_0)}}} = \varinjlim_n Rf^{(n)}_\ast (Q|_{\mathcal M_{Z^{(n_0)}}})|_{\mathcal M_{Z^{(n)}}} = 0
\]
(by applying it to the global sections on any quasicompact separated \'etale $\tilde{S}\to S$), and thus the map
\[
s^\ast \beta|_{\mathcal M_Z}: \mathbb F_\ell[2d]\to s^\ast Rf^! \mathbb F_\ell
\]
is an isomorphism, as desired. This finishes the proof of Proposition~\ref{prop:smoothatsection} and thus of Theorem~\ref{thm:jacobiancriterion}.
\end{proof}

The idea of the preceding proof is the following. Let $C\to X_S\times \mathbb A^1$ be the open subset of the deformation to the normal cone of $s:X_S\hookrightarrow Z$ (we did not develop the necessary formalism to give a precise meaning to this in the context of smooth sous-perfectoid spaces, but it could be done) whose fiber at $0\in \mathbb{A}^1$ is the normal cone of the immersion $s$ (the divisor over $0\in \mathbb{A}^1$ of the deformation to the normal cone is the union of two divisors: the projective completion of the normal cone and the blow-up of $Z$ along $X_S$, both meeting at infinity inside the projective completion).
One has a diagram 
\[\xymatrix{
X_S\times \mathbb{A}^1 \ar@{^{(}->}[r] \ar[rd] & C \ar[d] \\
& \mathbb{A}^1
}\]
where outside $t=0\in \mathbb{A}^1$ this is given by the section  $s:X_S\hookrightarrow Z$, i.e. the pullback over $\mathbb{G}_m$ of the preceding diagram gives the inclusion $X_S\times \mathbb{G}_m \hookrightarrow Z\times \mathbb{G}_m$, and at $t=0$ this is the inclusion of $X_S$ inside the normal cone of the section $s$.
Let us note moreover that $C$ is equipped with a $\mathbb{G}_m$-action compatible with the one on $\mathbb{A}^1$. 

This gives rise to an $\underline{E^\times}$-equivariant morphism with an equivariant section
\[\xymatrix{
\mathcal{M}_C \ar[d]^g \\  S\times \underline{E} \ar@/^1pc/[u]^s
}\]
whose fiber at $0\in \underline{E}$ is the zero section of $\BC (s^* T_{Z/X_S})\to S$, and is isomorphic to $\mathcal{M}_Z\times \underline{E^\times}$ equipped with the section $s$ outside of $0$. Now, the complex $s'^*Rg^!\mathbb{F}_\ell$ is $\underline{E^\times}$-equivariant on $S\times\underline{E}$. Its fiber outside $0\in\underline{E}$, i.e. its restriction to $S\times \underline{E^\times}$, is  $s^*Rf^!\mathbb{F}_\ell$, and its fiber at $0$ is $\mathbb{F}_\ell (d)[2d]$, $d=\deg (T_{Z/X_S})$ (since $g$ is universally locally acyclic the dualizing complex commutes with base change). 

Now one checks that one can replace the preceding diagram by a quasicompact $\underline{\O_E\setminus \{0\}}$-invariant open subset $U\subset \mathcal{M}_C$ together with an equivariant diagram
\[\xymatrix{
U\ar[d]^h \ar@{^{(}->}[r] & \mathcal{M_C} \ar[d]^g \\
S\times\underline{\O_E} \ar@/^1pc/[u]^t \ar@{^{(}->}[r] & S\times \underline{E} \ar@/^1pc/[u]^s.
}\]
In the preceding proof one replaces $\underline{\O_E}$ by $\underline{\pi^{\mathbb{N}\cup \{\infty\}}} \subset \underline{\O_E}$, which does not change anything for the argument. 
One concludes using that $\times \pi$ ``contracts everything to $0$'' and some constructibility argument using the fact that $U$ is spatial and some complexes are $h$-universally locally acyclic (see the argument ``$Rf'_!\mathbb D_{\mathcal M_{Z'}/S'}(Q)\in D_\et(S',\mathbb F_\ell)$ is constructible'' at the end of the proof of Proposition \ref{prop:smoothatsection}).

\section{Partial compactly supported cohomology}\label{sec:partialcompactsupport}

Let us start by recalling the following basic vanishing result. Let $C$ be a complete algebraically closed nonarchimedean field with pseudouniformizer $\varpi\in C$. Let $\Spa \mathbb Z\laurentseries{t}\times_{\mathbb Z} \Spa C = \mathbb D^\ast_C$ be the punctured open unit disc over $C$, and consider the subsets
\[
j: U = \{|t|\leq |\varpi|\}\hookrightarrow \mathbb D^\ast_C, j': U'=\{|t|\geq |\varpi|\}\hookrightarrow \mathbb D^\ast_C
\]
Note that the punctured open unit disc has two ``ends'': Towards the origin, and towards the boundary. The open subsets $U$ and $U'$ contain one ``end'' each.

\begin{lemma}\label{lem:cohomannulus} The partially compactly supported cohomology groups
\[
R\Gamma(\mathbb D^\ast_C,j_!\Lambda)=0 = R\Gamma(\mathbb D^\ast_C,j'_!\Lambda)
\]
vanish.
\end{lemma}

As usual $\Lambda$ is any coefficient ring killed by an integer $n$ prime to $p$.

\begin{proof} We treat the vanishing $R\Gamma(\mathbb D^\ast_C,j_!\Lambda)=0$, the other one being similar. Let $k:\mathbb{D}^\ast_C\hookrightarrow \mathbb{D}_C$ be the inclusion. One has an exact triangle
$$
(kj)_!\Lambda \lto Rk_* j_! \Lambda \lto i_*A\xrightarrow{\ +1 \ } 
$$ 
where $i:\{0\}\hookrightarrow \mathbb{D}_C$. One has $H^0(A)=\Lambda$, $H^1(A)=\Lambda(1)$,  $H^i(A)=0$ for $i\neq 0,1$, since 
$A=\varinjlim_n R\Gamma (U_n,\Lambda)$ with $U_n= \{|t|\leq |\varpi^n|\}\subset \mathbb{D}^\ast_C$ being a punctured disc. We thus have to prove that the preceding triangle induces an isomorphism $A\iso R\Gamma_c(U,\Lambda)[1]$.
Let $\tilde{j}:\mathbb{P}^1_C\setminus \{0,\infty\} \hookrightarrow \mathbb{P}^1_C\setminus \{0\}$.
There is a commutative diagram, obtained by applying $R\Gamma(\mathbb P^1_C,-)$ to an obvious diagram of sheaves:
$$
\begin{tikzcd}
R\Gamma_c (U,\Lambda) \ar[d,"\simeq"]\ar[r] & R\Gamma (\mathbb{D}^\ast_C,j_!\Lambda) \ar[d] \ar[r] & A \ar[r,"+1"] \ar[d,equal] & {} \\
R\Gamma_c (\mathbb{P}^1_C\setminus \{0,\infty \},\Lambda) \ar[r] & R\Gamma ( \mathbb{P}^1_C\setminus \{0\},\tilde{j}_!\Lambda)\ar[r] & A \ar[r,"+1"] & {}. 
\end{tikzcd}
$$
The left vertical map is an isomorphism by inspection (e.g., reduce to $U$ being an affinoid annulus, and the standard computation of its (compactly supported) cohomology). It thus suffices to check that $$ R\Gamma ( \mathbb{P}^1_C\setminus \{0\},\tilde{j}_!\Lambda)=0,$$ for example in the algebraic setting using comparison theorems, which is an easy exercise.
\end{proof}

Our goal now is to prove a very general version of such a result. Fix an algebraically closed field $k|\Fq$ and work on $\Perf_k$. Let $X$ be a spatial diamond such that $f: X\to \ast=\Spd k$ is partially proper with $\dimtrg f<\infty$. Then the base change $X\times_{k} S$ of $X$ to any spatial diamond $S$ is not itself quasicompact. Rather, it has two ends, and we will in this section study the cohomology with compact support towards one of the ends.

To analyze the situation, pick quasi-pro-\'etale and universally open surjections $\tilde{X}\to X$ and $\tilde{S}\to S$ from affinoid perfectoid spaces (using \cite[Proposition 11.24]{ECoD}), and pick maps $\tilde{X}\to \Spa k\laurentseries{t}$ and $\tilde{S}\to \Spa k\laurentseries{u}$ by choosing pseudouniformizers. We get a correspondence
$$
\begin{tikzcd}
& \tilde{X}\times \tilde{S}\ar[rd]\ar[ld] \\
X\times S && \Spa k\laurentseries{t}\times \Spa k\laurentseries{u} = \mathbb{D}^\ast_{k\laurentseries{u}}
\end{tikzcd}
$$
where all maps are qcqs, and the left map is (universally) open. Now $\Spa k\laurentseries{t}\times \Spa k\laurentseries{u}$ is a punctured open unit disc over $\Spa k\laurentseries{u}$, and one can write it as the increasing union of the affinoid subspaces
\[
\{|u|^b\leq |t|\leq |u|^a\}\subset \Spa k\laurentseries{t}\times \Spa k\laurentseries{u}
\]
for varying rational $0<a\leq b<\infty$. For any two choices of pseudouniformizers, a power of one divides the other, so it follows that if $\tilde{U}_{a,b}\subset \tilde{X}\times \tilde{S}$ denotes the preimage of $\{|t|^b\leq |u|\leq |t|^a\}$, then the doubly-indexed ind-system $\{\tilde{U}_{a,b}\}_{0<a\leq b<\infty}$ is independent of the choice of the maps $\tilde{X}\to \Spa k\laurentseries{t}$, $\tilde{S}\to \Spa k\laurentseries{u}$.

Let $U_{a,b}\subset X\times S$ be the image of $\tilde{U}_{a,b}$. As $\tilde{X}\times \tilde{S}\to X\times S$ is open, this is a qcqs open subset of $X\times S$. Moreover, the doubly indexed ind-system $\{U_{a,b}\}_{0<a\leq b<\infty}$ is independent of all choices made.

We let $U_a = \bigcup_{b<\infty} U_{a,b}$ and $U_b = \bigcup_{a>0} U_{a,b}$ and let
\[
j_{a,b}: U_{a,b}\to X\times S\ ,\ j_a: U_a\to X\times S\ ,\ j_b: U_b\to X\times S
\]
be the open immersions. We can now define the cohomology groups of interest, or rather the version of pushforward along $\beta: X\times S\to S$. As usual, $\Lambda$ is a coefficient ring killed by some integer $n$ prime to $p$.

\begin{definition} The functors
\[
R\beta_{!+},R\beta_{!-}: D_\et(X\times S,\Lambda)\to D_\et(S,\Lambda)
\]
are defined by
\[\begin{aligned}
R\beta_{!+}C &:= \varinjlim_a R\beta_{\ast}(j_{a!} C|_{U_a}),\\
R\beta_{!-}C &:= \varinjlim_b R\beta_{\ast}(j_{b!} C|_{U_b})\\
\end{aligned}\]
for $C\in D_\et(X\times S,\Lambda)$.
\end{definition}

The transition maps here are given by ($R\beta_\ast$ applied to) the counits of the adjunction between $j_{a,a'!}$ and $j_{a,a'}^\ast$, where $j_{a,a'}: U_a\subset U_{a'}$ is the open immersion. As the ind-systems of $U_a$ and $U_b$ are independent of all choices, these functors are canonical.

The main result is the following. Here $\alpha: X\times S\to X$ and $\beta: X\times S\to S$ are the two projections.

\begin{theorem}\label{thm:partialcompactsupportvanishing} Assume that $C=\alpha^\ast A\dotimes_\Lambda \beta^\ast B$ for $A\in D_\et(X,\Lambda)$ and $B\in D_\et(S,\Lambda)$. Then
\[
R\beta_{!+}C = 0 = R\beta_{!-} C.
\]
\end{theorem}

\begin{remark} The essential case for applications is $C=\alpha^\ast A$, i.e.~$B=\Lambda$, and $S=\Spa k\laurentseries{t}$. In other words, we take any coefficient system $A$ on $X$, pull it back to $X\times \Spa k\laurentseries{t}$, and then take the partially compactly supported cohomology (relative to $S$). However, it is sometimes useful to know the result in the relative case, i.e.~for general $S$, and then it is also natural to allow twists by $B\in D_\et(S,\Lambda)$.
\end{remark}

\begin{proof} We write the proof for $R\beta_{!+}$; the other case is exactly the same. Let $X_\bullet\to X$ be a simplicial hypercover by affinoid perfectoid spaces $X_i=\Spa(R_i,R_i^+)$ which are partially proper over $\Spa k$ (i.e., $R_i^+$ is minimal, i.e.~the integral closure of $k + R^{\circ\circ}$). As $X$ is a spatial diamond, we can arrange that the $X_i$ are the compactifications of quasi-pro-\'etale maps to $X$  (since $X$ is spatial it admits an hypercover $X_\bullet\to X$ with $X_i$ affinoid perfectoid and $X_i\to X$ quasi-pro-\'etale, since $X\to \Spd k$ is partially proper this extends to a hypercover $X_\bullet^c\to X$ where $X_i^c$ is Huber's canonical compactification over $\Spa(k)$); in particular, $g_i: X_i\to X$ satisfies $\dimtrg g_i=0<\infty$. Let $\beta_\bullet: X_\bullet\times S\to S$ be the corresponding projection. We claim that
\[
R\beta_{!+} C'
\]
is the limit of $R\beta_{\bullet,!+} (C'|_{X_\bullet\times S})$, for any $C'\in D_\et(X\times S,\Lambda)$. Writing $C'$ as a limit of its Postnikov truncations (\cite[Proposition 14.15]{ECoD}), we can assume $C'\in D^+_\et(X\times S,\Lambda)$. Now $g_i: X_i\to X$ is a qcqs map between spaces partially proper over $\ast$, so $g_i$ is proper, and hence so its base change $h_i: X_i\times S\to X\times S$. This implies that
\[
R\beta_{i,!+}(h_i^\ast C') = R\beta_{!+}(Rh_{i\ast} h_i^\ast C'),
\]
as $j_{a!}$ commutes with $Rh_{i\ast}$ by \cite[Theorem 19.2]{ECoD}. Now by \cite[Proposition 17.3]{ECoD}, one sees that $C'$ is the limit of $Rh_{i\ast} h_i^\ast C'$. But $R\beta_{!+}$ commutes with this limit, using that the filtered colimit does as everything lies in $D^+$ (with a uniform bound).

By the preceding reduction (used with $C'=C$), we may assume that $X=\Spa(R,R^+)$ is an affinoid perfectoid space. We can even assume that $X$ has no nonsplit finite \'etale covers (by taking the $X_i$ above to be compactifications of strictly totally disconnected spaces). In that case, there is a map $g: X\to Y=\Spa K$, where $K$ is the completed algebraic closure of $k\laurentseries{t}$, which is necessarily proper (as $X$ and $Y$ are partially proper over $\ast$), and as above one has
\[
R\beta_{!+}C = R\beta_{Y,!+} (Rh_\ast C)
\]
where $\beta_Y: Y\times S\to S$ is the projection and $h: X\times S\to Y\times S$ is the base change of $g$. Let $\alpha_Y: Y\times S\to Y$ be the other projection. Then the projection formula (and properness of $h$) \cite[Proposition 22.11]{ECoD} show that
\[
Rh_\ast C = Rh_\ast(\alpha^\ast A\dotimes_\Lambda \beta^\ast B)\cong Rh_\ast \alpha^\ast A\dotimes_\Lambda \beta_Y^\ast B
\]
and $Rh_\ast \alpha^\ast A\cong \alpha_Y^\ast Rg_\ast A$ by proper base change.

In other words, we can reduce to the case $X=\Spa K$; in particular $A\in D_\et(\Spa K,\Lambda) = D(\Lambda)$ is just a complex of $\Lambda$-modules. In that case, define $U_{a,b}$ and $U_a$ as above but taking $X=\tilde{X}\to \Spa k\laurentseries{t}$ the natural map. We claim that in this case for all $a>0$
\[
R\beta_\ast(j_{a!} C|_{U_a}) = 0.
\]
To prove this, it suffices to see that for all $a'>a>0$, the cone of
\[
R\beta_\ast(j_{a'!} C|_{U_{a'}})\to R\beta_\ast(j_{a!} C|_{U_a})
\]
vanishes, as $R\beta_\ast(j_{a!} C|_{U_a})$ is the limit of these cones as $a'\to \infty$. Now these cones depend on only a quasicompact part of $X\times S$, and hence their formation commutes with any base change in $S$, cf.~\cite[Proposition 17.6]{ECoD}. Therefore, we can reduce to the case $S=\Spa(L,L^+)$ for some complete algebraically closed nonarchimedean field $L$ with open and bounded valuation subring $L^+\subset L$, and check on global sections $R\Gamma(S,-)$. Moreover, the cone commutes with all direct sums in $C$, so one can assume that $A\in D_\et(\Spa K,\Lambda)=D_\et(\Lambda)$ is simply given by $A=\Lambda$.

It remains to prove the following statement: For all $B\in D_\et(\Spa(L,L^+),\Lambda)$ one has
\[
R\Gamma(\Spa K\times \Spa(L,L^+),j_{a!} B|_{U_a}) = 0.
\]
If the stalk of $B$ at the closed point vanishes, this follows from proper base change (writing $\Spa K\times \Spa(L,L^+)$ as the union of its subspaces proper over $\Spa(L,L^+)$), \cite[Theorem 19.2]{ECoD}. Thus we may assume that $B$ is concentrated at the closed point of $S$.

Recall that $U_a$ is defined using a choice of a quasi-pro-\'etale and universally open cover $\tilde{S}\to S$, together with a map $\tilde{S}\to \Spa k\laurentseries{u}$. As we only care about the closed point of $S$, we can assume that $\tilde{S}$ is the product of $S$ with a profinite set. The map $\tilde{S}\to \Spa k\laurentseries{u}$ is then given by a profinite set worth of maps $k\laurentseries{u}\to L$. The resulting subsets $\{|t|\leq |u|^a\}$ are locally constant on this profinite set, and for any two choices one is contained in the other. It follows that we may actually define $U_a$ using just one map $k\laurentseries{u}\to L$ (i.e. $\tilde{S}=S$). Thus, we can now assume that $U_a=\{|t|\leq |u|^a\}$, and we can also reduce to the case that $B$ is constant. Now using as above that the cones for $a'>a>0$ commute with any base change in $S$ and commute with direct sums in $B$, we can reduce to $B=\Lambda$ and the rank-$1$-geometric point $S=\Spa L$ where $L$ is the completed algebraic closure of $k\laurentseries{u}$.

At this point, we can further replace $\Spa K$ by $\Spa k\laurentseries{t}$: One can write $\Spa K$ as the inverse limit over finite extensions of $\Spa k\laurentseries{t}$, each of which is isomorphic to $\Spa k\laurentseries{t'}$, and although a priori $R\Gamma(\Spa K\times \Spa L,j_{a!} \Lambda)$ does not take this inverse limit to a colimit, this does happen after passing to cones for maps for $a'>a>0$, which suffices as above. Finally, we have reduced to Lemma~\ref{lem:cohomannulus}.
\end{proof}

\section{Hyperbolic localization}\label{sec:hyperboliclocalization}

In this section we extend some results of Braden, \cite{Braden}, to the world of diamonds. Our presentation is also inspired by the work of Richarz, \cite{RicharzHL}. We will use these results throughout our discussion of geometric Satake, starting in Section \ref{sec:semi infinite orbits}.

Let $S$ be a small v-stack, and let $f: X\to S$ be proper and representable in spatial diamonds with $\dimtrg f<\infty$, and assume that there is a $\mathbb G_m$-action on $X/S$, where $\mathbb G_m$ is the v-sheaf sending $\Spa(R,R^+)$ to $R^\times$. The fixed point space $X^0:=X^{\mathbb G_m}\subset X$ defines a closed subfunctor.

We make the following assumption about the $\mathbb G_m$-action. Here, $(\mathbb A^1)^+$ (resp.~$(\mathbb A^1)^-$) denotes the affine line $\Spa(R,R^+)\mapsto R$ with the natural $\mathbb G_m$-action (resp.~its inverse).

\begin{hypothesis}\label{hyp:goodGmaction} There is a decomposition of $X^0$ into open and closed subsets $X^0_1,\ldots,X^0_n$ such that for each $i=1,\ldots,n$, there are locally closed $\mathbb G_m$-stable subfunctors $X^+_i,X^-_i\subset X$ with $X^0\cap X^+_i=X^0_i$ (resp.~$X^0\cap X^-_i=X^0_i$) such that the $\mathbb G_m$-action on $X^+_i$ (resp.~$X^-_i$) extends to a $\mathbb G_m$-equivariant map $(\mathbb A^1)^+\times X^+_i\to X^+_i$ (resp.~$(\mathbb A^1)^-\times X^-_i\to X^-_i$), and such that
\[
X=\bigcup_{i=1}^n X^+_i = \bigcup_{i=1}^n X^-_i.
\]
\end{hypothesis}

We let 
\[X^+ = \bigsqcup_{i=1}^n X^+_i\ ,\ X^- = \bigsqcup_{i=1}^n X_i^-,
\]
so that there are natural maps 
\[
q^+: X^+\to X\ ,\ q^-: X^-\to X,
\]
as well as closed immersions 
\[
i^+: X^0\to X^+\ ,\ i^-: X^0\to X^-
\]
and projections
\[
p^+: X^+\to X^0\ ,\ p^-: X^-\to X^0;
\]
here $p^+$ is given by the restriction of $(\mathbb A^1)^+\times X_i^+\to X_i^+$ to $\{0\}\times X_i^+$, and $p^-$ is defined analogously.

Although the decomposition of $X^0$ into $X_i^0$ for $i=1,\ldots,n$ is a choice, ultimately the functors $X^+$ and $X^-$ are independent of any choice. Indeed, we have the following functorial description.

\begin{proposition}\label{prop:descrX+X-} Consider the functor $(X^+)^\prime$ sending any perfectoid space $T$ over $S$ to the set of $\mathbb G_m$-equivariant maps from $(\mathbb A^1)^+$ to $X$. There is a natural map $X^+\to (X^+)^\prime$, as there is a natural $\mathbb G_m$-equivariant map $(\mathbb A^1)^+\times X^+\to X^+\to X$. The map $X^+\to (X^+)^\prime$ is an isomorphism.

Analogously, $X^-$ classifies the set of $\mathbb G_m$-equivariant maps from $(\mathbb A^1)^-$ to $X$.
\end{proposition}

\begin{proof} It is enough to handle the case of $X^+$. There is a natural map $(X^+)^\prime\to X^0$ given by evaluating the $\mathbb G_m$-equivariant map on $(\mathbb A^1)^+\times (X^+)^\prime\to X$ on $\{0\}\times X^+$. Let $(X_i^+)^\prime = (X^+)^\prime\times_{X^0} X^0_i$; it is enough to prove that $X_i^+\to (X_i^+)^\prime$ is an isomorphism. For this, it is enough to prove that the map $(X_i^+)^\prime\to X$ given by evaluation at $1$ is an injection whose image is contained in the locally closed subspace $X_i^+\subset X$. This can be checked after pullback to an affinoid perfectoid base space $S=\Spa(R,R^+)$. As $X/S$ is proper (in particular, separated) and $\mathbb G_m\times (X_i^+)^\prime\subset (\mathbb A^1)^+\times (X_i^+)^\prime$ is dense, it follows that the map $(X_i^+)^\prime\to X$ is an injection. To bound its image, we can argue on geometric points. If $x\in |X|$ is any point in the image of $|(X_i^+)^\prime|$, and $\varpi\in R$ is a pseudouniformizer with induced action $\gamma$ on $X$, then the sequence $\gamma^n(x)$ converges to a point of $|X_i^0|$ for $n\to \infty$. On the other hand, if $x\not\in |X_i^+|$, then $x\in |X_j^+|$ for some $j\neq i$, which implies that $\gamma^n(x)$ converges to a point of $|X_j^0|$ for $n\to\infty$; this is a contradiction.

Thus, $(X_i^+)^\prime$ embeds into $X_i^+\subset X$, but it also contains $X_i^+$, so indeed $X_i^+=(X_i^+)^\prime$.
\end{proof}

\begin{lemma}\label{lem:plusminustransverse} The map $j: X^0\to X^+\times_X X^-$ is an open and closed immersion. More precisely, for any $i=1,\ldots,n$, the map $j_i: X_i^0\to X_i^+\times_X X_i^-$ is an isomorphism.
\end{lemma}

\begin{proof} It is enough to prove that for any $i=1,\ldots,n$, the map $j_i: X^0_i\to X^+_i\times_X X^-_i$ is an isomorphism. As it is a closed immersion (as $X^0_i\to X$ is a closed immersion and the target embeds into $X$), it is enough to prove that it is bijective on geometric rank $1$ points. Thus, we can assume $S=\Spa C$, and let $x: \Spa C=S\to X$ be a section that factors over $X^+_i\times_X X^-_i$. Then the $\mathbb G_m$-action on $x$ extends to a $\mathbb G_m$-equivariant map $g: \mathbb P^1_C\to X$. Consider the preimage of $X^+_i$ under $g$; this is a locally closed subfunctor, and it contains all geometric points. Indeed, on $(\mathbb A^1)_C^+$, the map $g$ factors over $X^+_i$ by hypothesis, and at $\infty$, it maps into $X^0_i\subset X^+_i$. This implies that the preimage of $X^+_i$ under $g$ is all of $\mathbb P^1_C$. In particular, we get a map
\[
(\mathbb A^1)^+\times \mathbb P^1_C\to (\mathbb A^1)^+\times X_i^+\to X_i^+
\]
which, when restricted to the copy of $\mathbb G_m$ embedded via $t\mapsto (t,t^{-1})$, is constant with value $x$. By continuity (and separatedness of $X_i^+$), this implies that it is also constant with value $x$ when restricted to $\mathbb A^1$ embedded via $t\mapsto (t,t^{-1})$, i.e.~the point $(0,\infty)\in (\mathbb A^1)^+\times \mathbb P^1_C$ maps to $x$. On the other hand, when restricted to $\mathbb G_m\times \{\infty\}$, the map is constant with values in $X_i^0$, and thus by continuity also on $(\mathbb A^1)^+\times\{\infty\}$. This implies that $x\in X_i^0$, as desired.
\end{proof}

In this setup, we can define two functors $D_\et(X,\Lambda)\to D_\et(X^0,\Lambda)$. We use the diagrams
$$
\begin{tikzcd}
X^{\pm}\ar[d,"p^{\pm}"] \ar[r,"q^{\pm}"] & X\\
X^0.
\end{tikzcd}
$$

\begin{definition} Define the functors
\[\begin{aligned}
L^+_{X/S} &= R(p^+)_! (q^+)^\ast: D_\et(X,\Lambda)\to D_\et(X^0,\Lambda),\\
L^-_{X/S} &= R(p^-)_\ast R(q^-)^!: D_\et(X,\Lambda)\to D_\et(X^0,\Lambda),
\end{aligned}\]
and a natural transformation $L^-_{X/S}\to L^+_{X/S}$ as follows. First, there are natural transformations
\[
R(i^+)^!=R(p^+)_! R(i^+)_! R(i^+)^!\to R(p^+)_!\ ,\ R(p^-)_\ast =  (i^-)^\ast (p^-)^\ast R(p^-)_\ast\to (i^-)^\ast,
\]
and the desired transformation $L^-_{X/S}\to L^+_{X/S}$ arises as a composite
\[
L^-_{X/S}=R(p^-)_\ast R(q^-)^!\to (i^-)^\ast R(q^-)^!\to (Ri^+)^! (q^+)^\ast\to R(p^+)_! (q^+)^\ast = L^+_{X/S},
\]
where the middle map $(i^-)^\ast R(q^-)^!\to (Ri^+)^! (q^+)^\ast$ of functors $D_\et(X,\Lambda)\to D_\et(X^0,\Lambda)$ is defined as the following composite
\[\begin{aligned}
(i^-)^\ast R(q^-)^!&\to (i^-)^\ast R(q^-)^! R(q^+)_\ast (q^+)^\ast\\
&=(i^-)^\ast R(\tilde{q}^-)_\ast R(\tilde{q}^+)^! (q^+)^\ast\\
&\to (i^-)^\ast R(\tilde{q}^-)_\ast j_\ast j^\ast R(\tilde{q}^+)^! (q^+)^\ast\\
&=(i^-)^\ast (i^-)_\ast Rj^! R(\tilde{q}^+)^! (q^+)^\ast\\
&=R(i^+)^! (q^+)^\ast,
\end{aligned}\]
using base change in the cartesian diagram
\[\xymatrix{
X^0\ar@{^(->}[r]^-j & X^+\times_X X^-\ar[d]^{\tilde{q}^+}\ar[r]^-{\tilde{q}^-} & X^+\ar[d]^{q^+}\\
& X^-\ar[r]^{q^-} & X.
}\]
Equivalently, it is enough to define for each $i=1,\ldots,n$ a natural transformation $(i_i^-)^\ast R(q_i^-)^!\to (Ri_i^+)^! (q_i^+)^\ast$ of functors $D_\et(X,\Lambda)\to D_\et(X_i^0,\Lambda)$. As
\[\xymatrix{
X_i^0\ar[r]^{i_i^+}\ar[d]^{i_i^-} & X_i^+\ar[d]^{q_i^+}\\
X_i^-\ar[r]^{q_i^-} & X
}\]
is cartesian, this arises as the composite
\[
(i_i^-)^\ast R(q_i^-)^!\to (i_i^-)^\ast R(q_i^-)^! R(q_i^+)_\ast (q_i^+)^\ast = (i_i^-)^\ast (i_i^-)_\ast R(i^+)^! (q_i^+)^\ast = R(i^+)^! (q_i^+)^\ast.
\]
\end{definition}

The following is our version of Braden's theorem, \cite{Braden}, cf.~\cite[Theorem B]{RicharzHL}.

\begin{theorem}\label{thm:braden} For any $A\in D_\et(X/\mathbb G_m,\Lambda)$ whose restriction to $X$ we continue to denote by $A$, the map
\[
L^-_{X/S} A\to L^+_{X/S} A
\]
is an isomorphism. In fact, moreover for any $A^+\in D_\et(X^+/\mathbb G_m,\Lambda)$, the map
\[
R(i^+)^! A^+\to R(p^+)_! A^+
\]
is an isomorphism, and for any $A^-\in D_\et(X^-/\mathbb G_m,\Lambda)$, the map
\[
R(p^-)_\ast A^-\to (i^-)^\ast A^-
\]
is an isomorphism, so that
\[
L^-_{X/S} A = R(p^-)_\ast R(q^-)^!A\cong (i^-)^\ast R(q^-)^! A\cong R(i^+)^! (q^+)^\ast A\cong R(p^+)_! (q^+)^\ast A = L^+_{X/S} A
\]
is a series of isomorphisms.
\end{theorem}

Before we start with the proof, we prove a certain general result about cohomology groups on spaces with ``two ends'', a flow connecting the two ends, and cohomology of sheaves, equivariant for the flow, that are compactly supported at only one end.

\begin{proposition}\label{prop:cohomflow} Let $S=\Spa(R,R^+)$ be an affinoid perfectoid space, $\varpi\in R$ a pseudouniformizer, let $f: Y\to S$ be a partially proper map of locally spatial diamonds, and assume that $Y$ is equipped with a $\mathbb G_m$-action over $S$. Assume that the quotient v-stack $Y/\mathbb G_m$ is qcqs. In that case, we can find a quasicompact open subset $V\subset Y$ such that $\mathbb G_m\times V\to Y$ is surjective and quasicompact. Write
\[
\mathbb G_{m,S} = \varinjlim_{n\geq 0} U_n\ ,\ U_n= \{x\in \mathbb G_{m,S}\mid |x|\leq |\varpi|^{-n}\},
\]
and let $j_n: V_n\subset Y$ be the open image of the cohomologically smooth map $U_n\times_S V\subset \mathbb G_{m,S}\times_S Y\to Y$. In this situation, we define for any $A\in D_\et(Y,\Lambda)$ the relative cohomology with partial supports
\[
\varinjlim_n Rf_\ast(j_{n!} A|_{V_n})\in D_\et(S,\Lambda)\ ;
\]
this functor is canonically independent of the choices made in its definition.

For any $A\in D^+_\et(Y/\mathbb G_m,\Lambda)$ (resp.~any $A\in D_\et(Y/\mathbb G_m,\Lambda)$ if $\dimtrg f<\infty$),
\[
\varinjlim_n Rf_\ast(j_{n!} A|_{V_n}) = 0.
\]
\end{proposition}

\begin{remark} Assume that $S=\Spa C$ is a geometric rank $1$ point. Then we set
\[
R\Gamma_{c+}(Y,A) = \varinjlim_n R\Gamma(Y,j_{n!} A_{V_n}),
\]
which is exactly the above functor $\varinjlim_n Rf_\ast(j_{n!} A|_{V_n})$ under the identification $D_\et(S,\Lambda) = D(\Lambda)$. Roughly speaking, the space $Y$ has two ends, one given by $\bigcup_{n<0} \gamma^n(V)$ for $V$ large enough, where $\gamma$ is the automorphism of $Y$ induced by $\varpi\in\mathbb G_m(S)$, and the other given by $\bigcup_{n>0} \gamma^n(V)$. We are considering the cohomology groups of $Y$ that have compact support in one of these directions, but not in the other. If one replaces the $\mathbb G_m$-action by its inverse, this implies a similar result for the direction of compact support interchanged.
\end{remark}

\begin{proof} One can assume that $A\in D^+_\et(Y/\mathbb G_m,\Lambda)$ by a Postnikov limit argument (in case $\dimtrg f<\infty$). Finding a v-hypercover of $Y$ by spaces with $\mathbb G_m$-action of the form $\mathbb G_m\times X_i$, where each $X_i$ is a proper spatial diamond over $S$, and using v-hyperdescent, one reduces to the case $Y=\mathbb G_m\times X_i$. Then there is the projection $\mathbb G_m\times X_i\to \mathbb G_{m,S}$, and one reduces to $Y=\mathbb G_{m,S}$. There, $A$ is a sheaf on the base $S$. One can now write $\mathbb G_m$ as an increasing union of punctured discs to reduce to Theorem~\ref{thm:partialcompactsupportvanishing}.
\end{proof}

\begin{proof}[Proof of Theorem~\ref{thm:braden}] We can assume that $A\in D^+_\et(X/\mathbb G_m,\Lambda)$ by pulling through the Postnikov limit $\varprojlim_n \tau^{\geq -n} A$, noting that $L^+_{X/S}$ commutes with limits while $(q^+)^\ast$ commutes with Postnikov limits and $R(p^+)_!$ as well by finite cohomological dimension.

By choosing a v-hypercover of $S$ by disjoint unions of strictly totally disconnected spaces $S_\bullet$, and using v-hyperdescent, we can assume that $S$ is a strictly totally disconnected space; indeed, $L^+_{X/S}$ commutes with all limits, while $(q^+)^\ast$ and $R(p^+)_!$ commute with any base change and so preserve cartesian objects, and thus also commute with the hyperdescent.

We start by proving that for any $A^+\in D_\et(X^+/\mathbb G_m,\Lambda)$, the map
\[
R(i^+)^! A^+\to R(p^+)_! A^+
\]
is an isomorphism, and similarly for any $A^-\in D_\et(X^-/\mathbb G_m,\Lambda)$, the map
\[
R(p^-)_\ast A^-\to (i^-)^\ast A^-
\]
is an isomorphism. Let $j^+: X^+\setminus X^0\hookrightarrow X^+$, $j^-: X^-\setminus X^0\hookrightarrow X^-$ denote the open embeddings. Then there are exact triangles
\[
(i^+)_\ast R(i^+)^! A^+\to A^+\to R(j^+)_\ast (j^+)^\ast A^+\ ,\ (j^-)_! (j^-)^\ast A^-\to A^-\to (i^-)_\ast (i^-)^\ast A^-.
\]
Using these triangles, we see that it is enough to see that for any $B^+\in D_\et((X^+\setminus X^0)/\mathbb G_m,\Lambda)$, $B^-\in D_\et((X^-\setminus X^0)/\mathbb G_m,\Lambda)$, one has
\[
R(p^+)_! R(j^+)_\ast B^+= 0\ ,\ R(p^-)_\ast R(j^-)_! B^-=0
\]
as objects in $D_\et(X^0,\Lambda)$. This follows from Proposition~\ref{prop:cohomflow} applied to $S=X^0$ and $Y=X^+\setminus X^0$ (resp.~$Y=X^-\setminus X^0$), and the following lemma.

\begin{lemma} The $\mathbb G_m$-action on $X^+\setminus X^0$ (resp.~$X^-\setminus X^0$) has the property that the quotient v-stack $(X^+\setminus X^0)/\mathbb G_m$ (resp.~$(X^-\setminus X^0)/\mathbb G_m$) is qcqs over $S$ (thus, over $X^0$).
\end{lemma}

\begin{proof} It is enough to do the case of $X^+\setminus X^0$, and we may restrict to $X_i^+\setminus X_i^0$. We can assume that $S=\Spa(R,R^+)$ is an affinoid perfectoid space, and fix a pseudouniformizer $\varpi\in R$. As $\mathbb G_{m,S}/\varpi^{\mathbb Z}$ is qcqs (in fact proper --- a Tate elliptic curve), it is equivalent to prove that $(X^+\setminus X^0)/\gamma^{\mathbb Z}$ is qcqs, where $\gamma$ is the automorphism given by the action of $\varpi\in \mathbb G_m(S)$.

Now we use the criterion of Lemma~\ref{lem:contractinglocspec} for the action of $\gamma$ on $|X_i^+|$. As a locally closed partially proper subspace of the proper spatial diamond $X$ over $S$, the locally spectral space $|X_i^+|$ is taut, and the condition on generizations is always fulfilled for locally spatial diamonds. The spectral closed subspace $|X_i^0|\subset |X_i^+|$ is fixed by $\gamma$, and by assumption for all $x\in |X_i^+|$, the sequence $\gamma^n(x)$ for $n\to \infty$ converges to a point of $|X_i^0|$ (as the $\mathbb G_m$-action extends to $(\mathbb A^1)^+$). It remains to see that for all $x\in |X_i^+\setminus X_i^0|$, the sequence $\gamma^n(x)$ for $n\to-\infty$ diverges in $|X_i^+|$. But $x\in |X_j^-|$ for some $j$, and $j\neq i$ by Lemma~\ref{lem:plusminustransverse}. Thus, $\gamma^n(x)$, for $n\to -\infty$, converges to a point of $|X_j^0|$, which is outside of $|X_i^+|$, so the sequence diverges in $|X_i^+|$.
\end{proof}

Now it remains to see that for any $A\in D^+_\et(X/\mathbb G_m,\Lambda)$, the map
\[
(i^-)^\ast R(q^-)^! A\to R(i^+)^! (q^+)^\ast A
\]
in $D^+_\et(X^0,\Lambda)=D^+((X^0)_\et,\Lambda)$ is an isomorphism. This can be done locally on $X^0$, so fix some $i\in \{1,\ldots,n\}$, and choose a quasicompact open neighborhood $U_0\subset X$ of $X_i^0$ that does not meet any $X_j^0$ for $j\neq i$ and such that $X_i^+\cap U_0,X_i^-\cap U_0\subset U_0$ are closed. The $\mathbb G_m$-orbit $Y=\mathbb G_m\cdot U_0\subset X$ is still open, and contains $X_i^+$ and $X_i^-$, necessarily as closed subsets.

We are now in the situation of the next proposition. To check conditions (ii) and (iii) of that proposition, note that we may find a quasicompact open subspace $V\subset Y$ such that $Y=\gamma^{\mathbb Z}\cdot V$ by averaging $U_0$ over $\{|\varpi|\leq |t|\leq 1\}\subset \mathbb G_m$. Let $W$ be the closure of $\bigcup_{n\geq 0} \gamma^n(V)\subset X$. To check (iii), it suffices (by symmetry) to see that
\[
\bigcap_{m\geq 0} \gamma^m(W)=\overline{X_i^-}
\]
in $X$. Note that $X_i^-\subset \bigcup_{n\geq 0} \gamma^n(V)$ (as for all $x\in X_i^-$, the sequence $\gamma^{-n}(x)$ converges into $X_i^0\subset V$), so $\overline{X_i^-}$ is contained in $W$, and thus in $\bigcap_{m\geq 0} \gamma^m(W)$. To prove the converse inclusion, let $W' = \bigcap_{m\geq 0} \gamma^m(W)$. If $\overline{X_i^-}\subsetneq W'$, then there is some $j\neq i$ such that $X_j^-$ contains a quasicompact open subset $A\subset W'$. Then $\tilde{A} = \bigcup_{n\geq 0} \gamma^{-n}(A)$ is a $\gamma^{-1}$-invariant open subset of $W'$ whose closure is $\gamma^{-(\mathbb N\cup\{\infty\})} \cdot \overline{A}$; in particular, replacing $A$ by $\gamma^{-n}(A)$ if necessary, we can arrange that this closure is contained in any given small neighborhood of $X_j^0$, and in particular intersects $V$ trivially. Then $\gamma^n(V)\cap A=\gamma^n(V\cap \gamma^{-n}(A))=\emptyset$ for all $n\geq 0$, and hence $A$ intersects $\bigcup_{n\geq 0} \gamma^n(V)$ trivially, and then also its closure $W$. But we assumed that $A\subset W'\subset W$, giving a contradiction.
\end{proof}

\begin{proposition} Let $S=\Spa(R,R^+)$ be a strictly totally disconnected perfectoid space, let $f: Y\to S$ be a compactifiable map of locally spatial diamonds, and assume that $Y/S$ is equipped with a $\mathbb G_m$-action, with fixed points $Y^0\subset Y$, and the following properties.
\begin{altenumerate}
\item[{\rm (i)}] There are $\mathbb G_m$-invariant closed subspaces $q^+: Y^+\subset Y$, $q^-: Y^-\subset Y$, containing $Y^0$ (via $i^+: Y^0\to Y^+$, $i^-: Y^0\to Y^-$) such that the action maps extend to maps $(\mathbb A^1)^+\times Y^+\to Y^+$ resp.~$(\mathbb A^1)^-\times Y^-\to Y^-$.
\item[{\rm (ii)}] The quotient v-stack $Y/\mathbb G_m$ is quasicompact. In particular, picking a pseudouniformizer $\varpi\in R$ with induced action $\gamma$ on $Y$, we can find some quasicompact open $V\subset Y$ such that $Y=\gamma^{\mathbb Z}\cdot V$.
\item[{\rm (iii)}]  With $V$ as in (ii), let $W^-$ be the closure of $\gamma^{\mathbb N}\cdot V$ and $W^+$ the closure of $\gamma^{-\mathbb N}\cdot V$. Then $\bigcap_{n\geq 0}\gamma^n(W^-) = Y^-$ and $\bigcap_{n\geq 0}\gamma^{-n}(W^+) = Y^+$.
\end{altenumerate}

Then $Y^0$ is a spatial diamond, the diagram
\[\xymatrix{
Y^0\ar[r]^{i^+}\ar[d]^{i^-} & Y^+\ar[d]^{q^+}\\
Y^-\ar[r]^{q^-}& Y
}\]
is cartesian, the quotient v-stacks $(Y\setminus Y^+)/\mathbb G_m$ and $(Y\setminus Y^-)/\mathbb G_m$ are qcqs, and for all $A\in D^+_\et(Y/\mathbb G_m,\Lambda)$ (resp.~all $A\in D_\et(Y/\mathbb G_m,\Lambda)$ if $\dimtrg f<\infty$) whose pullback to $Y$ we continue to denote by $A$, the map
\[
(i^-)^\ast R(q^-)^! A\to (i^-)^\ast R(q^-)^! (q^+)_\ast (q^+)^\ast A = (i^-)^\ast (i^-)_\ast R(i^+)^! (q^+)^\ast A = R(i^+)^! (q^+)^\ast A
\]
in $D_\et(Y^0,\Lambda)$ is an isomorphism.
\end{proposition}

\begin{proof} Note that $Y^0/\mathbb G_m\subset Y/\mathbb G_m$ is closed and thus $Y^0/\mathbb G_m$ is quasicompact. As the $\mathbb G_m$-action is trivial on $Y^0$, this implies that $Y^0$ is quasicompact. As $Y^0\subset Y$ is closed and $Y\to S$ is compactifiable and in particular quasiseparated, we see that $Y^0\to S$ is qcqs. That the diagram is cartesian follows from the proof of Lemma~\ref{lem:plusminustransverse}.

Next, we check that $(Y\setminus Y^+)/\mathbb G_m$ and $(Y\setminus Y^-)/\mathbb G_m$ are qcqs. By symmetry and as $\mathbb G_{m,S}/\gamma^{\mathbb Z}$ is qcqs, it suffices to see that $(Y\setminus Y^-)/\gamma^{\mathbb Z}$ is qcqs. First, we check that it is quasiseparated. Take any quasicompact open subspace $V^-\subset Y\setminus Y_i^-$; we need to see there are only finitely many $n$ with $V^-\cap \gamma^n(V^-)\neq \emptyset$. We can assume that $V^-\subset \bigcup_{n\geq 0}\gamma^n(V)$ (translating by a power of $\gamma$ if necessary), and then $V^-$ is covered by the open subsets $V^-\setminus \gamma^m(W^-)\subset V^-$ by the claim above. By quasicompacity, the intersection of $V^-$ with $\gamma^m(W^-)$ is empty for some large enough $m$, but then also the intersection of $V^-$ with $\gamma^{m'}(V^-)\subset \gamma^m(W^-)$ for $m'\geq m$ is empty.

To see that $(Y\setminus Y^-)/\gamma^{\mathbb Z}$ is quasicompact, note that $V\setminus \bigcup_{n>0} \gamma^n(V)$ is a spectral space (as it is closed in $V$) that maps bijectively to $(Y\setminus Y^-)/\gamma^{\mathbb Z}$.

Now, for the cohomological statement, we can as usual assume that $A\in D^+_\et(Y/\mathbb G_m,\Lambda)$ by a Postnikov limit argument. Then we are interested in checking that a map in $D^+_\et(Y_0,\Lambda)=D^+(Y_{0,\et},\Lambda)$ (cf.~\cite[Remark 14.14]{ECoD}) is an isomorphism, so we need to check that the sections over all quasicompact separated \'etale $Y_0'\to Y_0$ agree. Now we claim that any such quasicompact separated \'etale $Y_0'\to Y_0$ lifts to a $\mathbb G_m$-equivariant quasicompact separated \'etale map $Y'\to Y$; this will then allow us to assume $Y_0'=Y_0$ via passing to the pullback of everything to $Y'$.

To see that one may lift $Y_0'\to Y_0$ to $Y'\to Y$, consider the open subspace $V^{(n)}\subset Y$ given as the intersection of $\bigcup_{m\geq n} \gamma^m(V)$ with $\bigcup_{m\leq -n} \gamma^m(V)$. It follows from the topological situation that this is still quasicompact, and that the intersection of all $V^{(n)}$ is equal to $Y_0$ (using condition (iii)). Let $Y^{(n)}=\gamma^{\mathbb Z}\cdot V^{(n)}\subset Y$. Then $\gamma$-equivariant quasicompact separated \'etale maps to $Y^{(n)}$ are equivalent to quasicompact separated \'etale maps to $V^{(n)}$ together with isomorphisms between the two pullbacks to $V^{(n)}\cap \gamma(V^{(n)})$. The latter data extends uniquely from $Y_0$ to $V^{(n)}$ for small enough $n$ by \cite[Proposition 11.23]{ECoD}. Repeating a similar argument after taking a product with $\mathbb G_{m,S}/\gamma^{\mathbb Z}$ (which is qcqs), and observing that the $Y^{(n)}$ are cofinal with their $\mathbb G_m$-orbits, one can then attain $\mathbb G_m$-equivariance.

We have now reduced to checking the statement on global sections. Now consider the compactification $j: Y\hookrightarrow \overline{Y} = \overline{Y}^{/S}\to S$. Note that $\overline{Y}$ satisfies all the same conditions of the proposition. Restricted to $Y_0$, this gives a quasicompact open immersion $j_0: Y_0\hookrightarrow \overline{Y_0}$. By the above argument, this quasicompact open immersion spreads to a quasicompact open immersion into $\overline{Y}$, and by taking it small enough in the argument above, we can assume that it is contained in $Y$. This allows us to assume that $j$ is quasicompact. In that case the functor $Rj_\ast$ commutes with all operations in question by \cite[Proposition 17.6, Proposition 23.16 (i)]{ECoD}. Thus, we can now moreover assume that $Y$ is partially proper.

Our goal now is to prove that when $Y$ is partially proper and $A\in D^+_\et(Y/\mathbb G_m,\Lambda)$, the map
\[
(i^-)^\ast R(q^-)^! A\to  R(i^+)^! (q^+)^\ast A
\]
becomes an isomorphism after applying $Rf_{0\ast}$ where $f_0: Y_0\to S$ is the proper map. For this, we define another functor $D_\et(Y,\Lambda)\to D_\et(S,\Lambda)$, as follows. Let $j_n: V_n = \bigcup_{m\geq -n} \gamma^m(V)\hookrightarrow Y$ for $n\geq 0$. Then we consider
\[
A\mapsto F(A) = \varinjlim_n Rf_\ast(j_{n!}A|_{V_n}): D_\et(Y,\Lambda)\to D_\et(S,\Lambda).
\]

\begin{lemma} Let $j^-: Y\setminus Y^-\to Y$, $j^+: Y\setminus Y^+\to Y$ denote the open immersions.
\begin{altenumerate}
\item If $A=Rj^-_\ast A^-$ for $A^-\in D^+_\et((Y\setminus Y^-)/\mathbb G_m,\Lambda)$, then $F(A)=0$.
\item If $A=j^+_! A^+$ for $A^+\in D^+_\et((Y\setminus Y^+)/\mathbb G_m,\Lambda)$, then $F(A)=0$.
\end{altenumerate}
\end{lemma}

\begin{proof} This follows from Proposition~\ref{prop:cohomflow} and condition (iii).
\end{proof}

There are natural transformations $Rf_{0\ast} R(p^-)_\ast R(q^-)^!\to F\to Rf_{0\ast} R(p^+)_! (q^+)^\ast$, and the lemma implies that these are equivalences when evaluated on $A\in D^+_\et(Y/\mathbb G_m,\Lambda)$. Using that also $(Y^+\setminus Y^0)/\mathbb G_m$ and $(Y^-\setminus Y^0)/\mathbb G_m$ are qcqs (as closed subspaces of $(Y\setminus Y^-)/\mathbb G_m$ resp.~$(Y\setminus Y^+)/\mathbb G_m$) so that we can apply Proposition~\ref{prop:cohomflow} again as in the beginning of the proof of Theorem~\ref{thm:braden}, we get an isomorphism
\[
Rf_{0\ast} (i^-)^\ast R(q^-)^! A\cong Rf_{0\ast} R(p^-)_\ast R(q^-)^! A\cong F(A)\cong Rf_{0\ast}R(p^+)_! (q^+)^\ast A\cong Rf_{0\ast} R(i^+)^! (q^+)^\ast A.
\]
We need to see that this implies that also the map
\[
Rf_{0\ast} (i^-)^\ast R(q^-)^! A\to  Rf_{0\ast} R(i^+)^! (q^+)^\ast A
\]
defined in the statement of the proposition is an isomorphism. For this, observe that this map is an isomorphism if and only if for $A=j^+_! A^+$ with $A^+\in D_\et((Y\setminus Y^+)/\mathbb G_m,\Lambda)$, one has
\[
Rf_{0\ast} (i^-)^\ast R(q^-)^! A = 0.
\]
But this follows from the existence of some isomorphism
\[
Rf_{0\ast} (i^-)^\ast R(q^-)^! A\cong Rf_{0\ast} R(i^+)^! (q^+)^\ast A=0,
\]
using $(q^+)^\ast A = (q^+)^\ast j^+_! A^+ = 0$.
\end{proof}

Using Theorem~\ref{thm:braden}, we give the following definition.

\begin{definition}\label{def:hyperboliclocalization} Let $f: X\to S$ with $\mathbb G_m$-action be as above, satisfying Hypothesis~\ref{hyp:goodGmaction}. Let $D_\et(X,\Lambda)^{\mathbb G_m\text{-mon}}\subset D_\et(X,\Lambda)$ be the full subcategory generated under finite colimits and retracts by the image of $D_\et(X/\mathbb G_m,\Lambda)\to D_\et(X,\Lambda)$. The hyperbolic localization functor is the functor
\[
L_{X/S}: D_\et(X,\Lambda)^{\mathbb G_m\text{-mon}}\to D_\et(X^0,\Lambda)
\]
given by $L_{X/S}^-\cong L_{X/S}^+$.
\end{definition}

We observe that Theorem~\ref{thm:braden} implies the following further results.

\begin{proposition}\label{prop:hyperboliclocalizationbasechange} In the situation of Definition~\ref{def:hyperboliclocalization}, let $g: S'\to S$ be a map of small v-stacks, with pullback $f': X'=X\times_S S'\to S'$, $g_X: X'\to X$, $g^0: X^{0\prime}\to X^0$. Then there are natural equivalences
\[
g^{0\ast} L_{X/S}\cong L_{X'/S'} g_X^\ast\ ,\ L_{X/S} Rg_{X\ast}\cong Rg^0_\ast L_{X'/S'}\ ,\ L_{X/S} Rg_{X!}\cong Rg_{0!} L_{X'/S'}\ ,\ Rg^{0!} L_{X/S}\cong L_{X'/S'} Rg_X^!,
\]
the latter two in case $g$ is compactifiable and representable in locally spatial diamonds with $\dimtrg g<\infty$ (so that the relevant functors are defined).
\end{proposition}

\begin{proof} The first and third assertions are clear for $L_{X/S}^+$, while the second and fourth assertions are clear for $L_{X/S}^-$.
\end{proof}

\begin{proposition}\label{prop:hyperboliclocalizationdual} In the situation of Definition~\ref{def:hyperboliclocalization}, let $A\in D_\et(X,\Lambda)^{\mathbb G_m\text{-mon}}$ and $B\in D_\et(S,\Lambda)$. Let $L'_{X/S}$ denote the hyperbolic localization functor for the inverse $\mathbb G_m$-action. Then there is natural isomorphism
\[
R\sHom(L_{X/S}(A),Rf^{0!}B)\cong L'_{X/S} R\sHom(A,Rf^!B).
\]
In particular, taking $B=\Lambda$, hyperbolic localization commutes with Verdier duality, up to changing the $\mathbb G_m$-action.
\end{proposition}

\begin{proof} More generally, for all $A\in D_\et(X,\Lambda)$ and $B\in D_\et(S,\Lambda)$, we have a natural isomorphism $R\sHom(L_{X/S}^+(A),Rf^{0!}B)\cong L_{X/S}^{\prime-}R\sHom(A,Rf^!B)$. Indeed,
\[\begin{aligned}
R\sHom(L_{X/S}^+(A),Rf^{0!} B)&=R\sHom(R(p^+)_! (q^+)^\ast A,Rf^{0!} B)\cong R(p^+)_\ast R\sHom((q^+)^\ast A,R(p^+)^!Rf^{0!} B)\\
&\cong R(p^+)_\ast R\sHom((q^+)^\ast A,R(q^+)^! Rf^! B)\cong R(p^+)_\ast R(q^+)^!R\sHom(A,Rf^!B).
\end{aligned}\]
\end{proof}

\begin{proposition}\label{prop:hyperboliclocalizationULA} In the situation of Definition~\ref{def:hyperboliclocalization}, assume that $A\in D_\et(X,\Lambda)^{\mathbb G_m\text{-mon}}$ is $f$-universally locally acyclic. Then $L_{X/S}(A)\in D_\et(X^0,\Lambda)$ is universally locally acyclic with respect to $f^0: X^0\subset X\to S$.
\end{proposition}

\begin{proof} As the assumption is stable under base change, we may assume that $S$ is strictly totally disconnected, and it suffices to see that $L_{X/S}(A)$ is $f^0$-locally acyclic. For condition (a), we can in fact assume that $S=\Spa(C,C^+)$ is strictly local; let $j: S_0=\Spa(C,\mathcal O_C)\subset S$ be the generic open point. Then we have to see that $L_{X/S}(A) = Rj^0_\ast(L_{X/S}(A)|_{X^0\times_S S_0})$, where $j^0: X^0\times_S S_0\to X^0$ is the pullback of $j$. But this follows from Proposition~\ref{prop:hyperboliclocalizationbasechange} and the corresponding property of $A$.

For condition (b), it suffices to see that the functor $R\sHom_\Lambda(L_{X/S}(A),Rf^{0!}-)$ commutes with all direct sums, as then its left adjoint $Rf^0_!(L_{X/S}(A)\dotimes_\Lambda -)$ preserves perfect-constructible complexes. For this, we compute this functor:
\[\begin{aligned}
R\sHom_\Lambda(L_{X/S}(A),Rf^{0!} -) &\cong L_{X/S}' R\sHom_\Lambda(A,Rf^! -)\\
&\cong L_{X/S}'  (\mathbb D_{X/S}(A)\dotimes_\Lambda f^\ast -).
\end{aligned}\]
Here, we used Proposition~\ref{prop:hyperboliclocalizationdual} and Proposition~\ref{prop:ULAastshriekpullback}. The final functor clearly commutes with all direct sums, giving the desired result.
\end{proof}

\section{Drinfeld's lemma}\label{sec:drinfeld}

As a final topic of this chapter, we prove the version of Drinfeld's lemma that we will need in this paper. Contrary to the classical formulation \cite[Theorem 2.1]{DrinfeldLanglandsConjectureForGL2}, cf.~also \cite[Theorem 8.1.4]{LauThesis}, this version actually makes the Weil group of $E$, not the absolute Galois group of $E$, appear. (Also, it is worth remarking that usually, a global Galois group appears, not a local Galois group.)

In this section, we work on $\Perf_k$ where $k=\overline{\mathbb F}_q$. In that case, we can write the moduli space of degree $1$ Cartier divisors on the Fargues--Fontaine curve as $\Div^1 = \Spd \breve E/\phi^{\Z}$, where $\phi$ acts on $\Spd \breve{E} = \Spd k\times_{\Spd \mathbb F_q} \Spd E$ via the second factor. This admits a natural map
\[
\psi: \Div^1\to [\ast/\underline{W_E}]
\]
to the classifying space of the Weil group of $E$. Indeed, if $C=\hat{\overline{E}}$ is a completed algebraic closure of $E$, then there is an action of $W_E$ on $\Spd C$, with the inertia subgroup $I_E\subset W_E$ acting via its usual action, while Frobenius elements act via the composite of the usual action and the Frobenius of $\Spd C$. More precisely, $\tau\in W_E$ acts as $\tau\circ \Frob^{-\deg \tau}$ where $\deg: W_E\to \Z$ is the projection; note that this as a map over $\Spd k$ as on $\Spd k$ the two Frobenii cancel. The natural map
\[
[\Spd C/\underline{W_E}]\to [\Spd \breve E/\phi^{\Z}]
\]
is an isomorphism, thus yielding the natural map
\[
\psi: [\Spd \breve E/\phi^{\Z}]\cong [\Spd C/\underline{W_E}]\to [\ast/\underline{W_E}].
\]
One could equivalently compute
\[
\underline{W_E}\times \Spd C\cong \Spd C\times_{\Div^1} \Spd C
\]
for the natural map $\Spd C\to \Div^1$ to arrive at the result.

In particular, for any small v-stack $X$, we get a natural map
\[
\psi_X: X\times \Div^1\to X\times [\ast/\underline{W_E}].
\]
As usual, $\Lambda$ is a ring killed by some integer $n$ prime to $p$.

\begin{proposition}\label{prop:drinfeldfullyfaithful} The functor
\[
\psi_X^\ast: D_\et(X\times [\ast/\underline{W_E}],\Lambda)\to D_\et(X\times \Div^1,\Lambda)
\]
is fully faithful. If the natural pullback functor
\[
D_\et(X,\Lambda)\to D_\et(X\times \Spd C,\Lambda)
\]
is an equivalence, then $\psi_X^\ast$ is also an equivalence.
\end{proposition}

\begin{proof} We apply descent along $\ast\to [\ast/\underline{W_E}]$. This describes $D_\et(X\times [\ast/\underline{W_E}],\Lambda)$ in terms of cartesian objects in $D_\et(X\times \underline{W_E}^{\bullet},\Lambda)$, and $D_\et(X\times \Div^1,\Lambda)$ in terms of cartesian objects in $D_\et(X\times \Spd C\times \underline{W_E}^\bullet,\Lambda)$. By \cite[Theorem 1.13]{ECoD}, all functors $D_\et(X\times \underline{W_E}^{\bullet},\Lambda)\to D_\et(X\times \Spd C\times \underline{W_E}^\bullet,\Lambda)$ are fully faithful; this implies the fully faithfulness. Moreover, for essential surjectivity on cartesian objects it is enough to know essential surjectivity on the degree $0$ part of the simplicial resolution, i.e.~for $D_\et(X,\Lambda)\to D_\et(X\times \Spd C,\Lambda)$, giving the desired result.
\end{proof}

We note the following immediate corollary.

\begin{corollary}\label{cor:drinfeldfullyfaithful} For any finite set $I$, pullback along $X\times (\Div^1)^I\to X\times [\ast/\underline{W_E^I}]$ induces a fully faithful functor
\[
D_\et(X\times [\ast/\underline{W_E^I}],\Lambda)\to D_\et(X\times (\Div^1)^I,\Lambda).
\]
\end{corollary}

\begin{proof} This follows inductively from Proposition~\ref{prop:drinfeldfullyfaithful}.
\end{proof}

We need the following refinement, see Proposition \ref{prop:locsysdiv1}. For any small v-stack $Y$, let
\[
D_\lc(Y,\Lambda)\subset D_\et(Y,\Lambda)
\]
be the full subcategory of all objects that are v-locally constant with perfect fibres. (Being v-locally constant with perfect fibres is equivalent to dualizability, and on spatial diamonds such objects are actually \'etale locally constant, as follows from Proposition~\ref{prop:constructible contre localement constant}.)

\begin{proposition}\label{prop:drinfeldlemmalocallyconstant} For any finite set $I$ and any small v-stack $X$, the functor
\[
D_\lc(X\times [\ast/\underline{W_E^I}],\Lambda)\to D_\lc(X\times (\Div^1)^I,\Lambda)
\]
is an equivalence of categories.
\end{proposition}

We will mostly be using this in case $X$ is a point. The equivalence certainly fails without the local constancy condition, as there are sheaves supported on proper subsets of $|(\Div^1)^I|$, like the partial diagonals.

\begin{proof} By Corollary~\ref{cor:drinfeldfullyfaithful}, the functor is fully faithful. By induction, we can reduce to the case that $I$ has one element. By descent, we can assume that $X$ is strictly totally disconnected. Note that $X\times \Div^1$ is a spatial diamond, and using \cite[Proposition 20.15]{ECoD} we can reduce to the case that $X=\Spa(C,C^+)$ is strictly local (by writing any connected component as a cofiltered inverse limit of its open and closed neighborhoods to see that then any object is locally in the image of the functor). Moreover, the category $D_\lc$ is unchanged if we replace $\Spa(C,C^+)$ by $\Spa(C,\mathcal O_C)$, so we can assume that $X$ is even a geometric rank $1$ point.

At this point, we need to simplify the coefficient ring $\Lambda$. The algebra $\Lambda$ is a $\mathbb Z/n\mathbb Z$-algebra for some $n$ prime to $p$; we can then assume $n$ is a power of some prime $\ell\neq p$, and in fact even $n=\ell$ by an induction argument. By \cite[Proposition 20.15]{ECoD}, we can also assume that $\Lambda$ is a finitely generated $\mathbb F_\ell$-algebra. Taking a surjection from a polynomial algebra, one can then assume that $\Lambda=\mathbb F_\ell[T_1,\ldots,T_d]$. Applying \cite[Proposition 20.15]{ECoD} again, we can assume that $\Lambda$ is the localization of $\mathbb F_\ell[T_1,\ldots,T_d]$ at a closed point, or applying faithfully flat descent in the coefficients, that $\Lambda$ is the completion of $\mathbb F_\ell[T_1,\ldots,T_d]$ at a closed point, but equipped with the discrete topology. Also note that this ring is regular, so all truncations of perfect complexes are perfect, and we can assume that the complex is concentrated in degree $0$.

We are now in the following situation. We have an \'etale sheaf $A$ of $\Lambda\cong \mathbb F_{\ell^r}\powerseries{T_1,\ldots,T_d}$-modules on $S=\Spa C\times \Div^1$, such that for some finitely generated $\Lambda$-module $M$, there are \'etale local isomorphisms between $A$ and the constant $\Lambda$-module associated to $M$. Our goal is to see that after pullback along the $\underline{W_E}$-torsor
\[
\tilde{S}=\Spa C\times \Spd \hat{\overline{E}}\to S=\Spa C\times \Div^1,
\]
there is an isomorphism between $A$ and $M$. To see this, we will also need to analyze the behaviour at a carefully chosen geometric point. In fact, by Lemma~\ref{lem:genericpointcurve} we can find a point $\Spa K\to Y_C$ of the curve $Y_C$ associated with $C$ such that the induced map $\Gal(\overline{K}|K)\to I_E$ is surjective. This induces a point $y: \Spd K\to S$, and we can lift it to a geometric point $\tilde{y}: \Spd \hat{\overline{K}}\to \tilde{S}$. Define $M$ as the stalk of $A$ at $\tilde{y}$; our goal is then to prove the existence of a unique isomorphism between $A|_{\tilde{S}}$ and $M$ that is the identity at $\tilde{y}$.

To prove this, we first reduce modulo $(T_1,\ldots,T_d)^n$. Then $\Lambda_n:=\Lambda/(T_1,\ldots,T_d)^n$ is a finite ring, and the space of isomorphisms between $A/(T_1,\ldots,T_d)^n$ and $M/(T_1,\ldots,T_d)^n$ is parametrized by a space finite \'etale over $S$. By \cite[Lemma 16.3.2]{Berkeley}, all such finite \'etale covers come via pullback from finite \'etale covers of $\Div^1$, and are thus trivialized after pullback to $\tilde{S}$; this implies that there is a unique isomorphism $A/(T_1,\ldots,T_d)^n\cong M/(T_1,\ldots,T_d)^n$ reducing to the identity at $\tilde{y}$.

Taking the limit over $n$, we get an isomorphism $\hat{A}|_{\tilde{S}}\cong \hat{M}|_{\tilde{S}}$ between the pro-\'etale sheaves $\hat{A}=\varprojlim_n A/(T_1,\ldots,T_d)^n$ and $\hat{M}=\varprojlim_n M/(T_1,\ldots,T_d)^n$ after pullback to $\tilde{S}$. This gives in particular an automorphism of $\widehat{M}$ over
\[
\tilde{S}\times_S \tilde{S}\cong \underline{W_E}\times \tilde{S},
\]
and thus by connectedness of $\tilde{S}$ a continuous map $W_E\to \Aut_\Lambda(\hat{M})$ (in fact, it extends continuously to the absolute Galois group of $E$). We claim that this map is trivial on an open subgroup of $I_E$ (but not necessarily on an open subgroup of the absolute Galois group of $E$ --- here it is necessary to pass to the Weil group). Indeed, restricting the map $W_E\to \Aut_\Lambda(\hat{M})$ to $\Gal(\overline{K}|K)$ gives a map $\Gal(\overline{K}|K)\to \Aut_\Lambda(\hat{M})$ that is in fact continuous for the discrete topology on the target, as a local system of $\Lambda$-modules on $\Spd K$ is given by a continuous representation of $\Gal(\overline{K}|K)$. As $\Gal(\overline{K}|K)\to I_E$ is surjective, we get the claim.

By equivariance under an open subgroup of $I_E$, we find that the isomorphism $\hat{A}|_{\tilde{S}}\cong \hat{M}|_{\tilde{S}}$ descends, necessarily uniquely, to an isomorphism over
\[
\Spa C\times \Spd E'
\]
for some finite extension $E'|\breve E$. Now we take the pushforward of the isomorphism $\hat{A}|_{\Spa C\times \Spd E'}\cong \hat{M}|_{\Spa C\times \Spd E'}$ to the small \'etale site of $\Spa C\times \Spd E'$. As any \'etale $U\to \Spa C\times \Spd E'$ is locally connected, we have $H^0(U,M)=H^0(U,\hat{M})$ and then also $H^0(U,A)=H^0(U,\hat{A})$ (as $A$ is \'etale locally isomorphic to $M$) for all such $U$, so we get the desired isomorphism $A|_{\Spa C\times \Spd E'}\cong M|_{\Spa C\times \Spd E'}$.
\end{proof}

\chapter{$D_{\mathrm{\acute{e}t}}(\mathrm{Bun}_G)$}\label{ch:DetBunG}

In this chapter, we want to understand the basic structure of $D_\et(\Bun_G,\Lambda)$, building it up from all $D_\et(\Bun_G^b,\Lambda)$, where we continue to work in the setting where $\Lambda$ is killed by some integer $n$ prime to $p$.

Throughout this chapter, we fix an algebraically closed field $k|\Fq$ and work on $\Perf_k$. Our goal is to prove the following theorem.

\begin{theorem}[Theorem~\ref{thm:chartbsmooth}, Proposition~\ref{prop:basicMb}; Proposition~\ref{prop:DetBunstrata}, Theorem~\ref{thm:Detclassifying}; Theorem~\ref{thm:compactgeneration}; Theorem~\ref{thm:bernsteinzelevinsky}; Theorem~\ref{thm:ULAbunG}] Let $\Lambda$ be any ring killed by some integer $n$ prime to $p$.
\begin{altenumerate}
\item[{\rm (o)}] For any $b\in B(G)$, there is a map
\[
\pi_b: \mathcal M_b\to \Bun_G
\]
that is representable in locally spatial diamonds, partially proper and cohomologically smooth, where $\mathcal M_b$ parametrizes $G$-bundles $\mathcal E$ together with an increasing $\mathbb Q$-filtration whose associated graded is, at all geometric points, isomorphic to $\mathcal E_b$ with its slope grading. The v-stack $\mathcal M_b$ is representable in locally spatial diamonds, partially proper and cohomologically smooth over $[\ast/\underline{G_b(E)}]$.
\item[{\rm (i)}] Via excision triangles, there is an infinite semiorthogonal decomposition of $D_\et(\Bun_G,\Lambda)$ into the various $D_\et(\Bun_G^b,\Lambda)$ for $b\in B(G)$.
\item[{\rm (ii)}] For each $b\in B(G)$, pullback along
\[
\Bun_G^b\cong [\ast/\tilde{G}_b]\to [\ast/\underline{G_b(E)}]
\]
gives an equivalence
\[
D_\et([\ast/\underline{G_b(E)}],\Lambda)\cong D_\et(\Bun_G^b,\Lambda),
\]
and $D_\et([\ast/\underline{G_b(E)}],\Lambda)\cong D(G_b(E),\Lambda)$ is equivalent to the derived category of the category of smooth representations of $G_b(E)$ on $\Lambda$-modules.
\item[{\rm (iii)}] The category $D_\et(\Bun_G,\Lambda)$ is compactly generated, and a complex $A\in D_\et(\Bun_G,\Lambda)$ is compact if and only if for all $b\in B(G)$, the restriction
\[
i^{b\ast} A\in D_\et(\Bun_G^b,\Lambda)\cong D(G_b(E),\Lambda)
\]
is compact, and zero for almost all $b$. Here, compactness in $D(G_b(E),\Lambda)$ is equivalent to lying in the thick triangulated subcategory generated by $c\text-\mathrm{Ind}_K^{G_b(E)} \Lambda$ as $K$ runs over open pro-$p$-subgroups of $G_b(E)$.
\item[{\rm (iv)}] On the subcategory $D_\et(\Bun_G,\Lambda)^\omega\subset D_\et(\Bun_G,\Lambda)$ of compact objects, there is a Bernstein--Zelevinsky duality functor
\[
\mathbb D_{BZ}: (D_\et(\Bun_G,\Lambda)^\omega)^{\mathrm{op}}\to D_\et(\Bun_G,\Lambda)^\omega
\]
with a functorial identification
\[
R\Hom(A,B)\cong \pi_\natural(\mathbb D_{BZ}(A)\dotimes_\Lambda B)
\]
for $B\in D_\et(\Bun_G,\Lambda)$, where $\pi: \Bun_G\to \ast$ is the projection. The functor $\mathbb D_{BZ}$ is an equivalence, and $\mathbb D_{BZ}^2$ is naturally equivalent to the identity. It is compatible with usual Bernstein--Zelevinsky duality on $D(G_b(E),\Lambda)$ for basic $b\in B(G)$.
\item[{\rm (v)}] An object $A\in D_\et(\Bun_G,\Lambda)$ is universally locally acyclic (with respect to $\Bun_G\to \ast$) if and only if for all $b\in B(G)$, the restriction
\[
i^{b\ast} A\in D_\et(\Bun_G^b,\Lambda)\cong D(G_b(E),\Lambda)
\]
is admissible, i.e.~for all pro-$p$ open subgroups $K\subset G_b(E)$, the complex $(i^{b\ast} A)^K$ is perfect. Universally locally acyclic complexes are preserved by Verdier duality, and satisfy Verdier biduality.
\end{altenumerate}
\end{theorem}

\section{Classifying stacks}

First, we want to understand $D_\et([\ast/\underline{G}],\Lambda)$ for a locally pro-$p$-group $G$. Fix a coefficient ring $\Lambda$ such that $n\Lambda=0$ for some $n$ prime to $p$, and assume that $G$ is locally pro-$p$. Our aim is to prove the following theorem.

\begin{theorem}\label{thm:Detclassifying} Let $D(G,\Lambda)$ be the derived category of the category of smooth representations of $G$ on $\Lambda$-modules. There is a natural symmetric monoidal equivalence
\[
D(G,\Lambda)\simeq D_\et([\ast/\underline{G}],\Lambda)
\]
under which the functor $D(G,\Lambda)\to D(\Lambda)$ forgetting the $G$-action gets identified with the pullback functor $D_\et([\ast/\underline{G}],\Lambda)\to D_\et(\ast,\Lambda)=D(\Lambda)$ under the projection $\ast\to [\ast/\underline{G}]$.

The same result holds true for the base change $[\Spa C/\underline{G}] = [\ast/\underline{G}]\times \Spa C$ for any complete algebraically closed nonarchimedean field $C/k$; more precisely, the base change functor
\[
D_\et([\ast/\underline{G}],\Lambda)\to D_\et([\Spa C/\underline{G}],\Lambda)
\]
is an equivalence.
\end{theorem}

Note that indeed
\[
D_\et(\ast,\Lambda)=D(\Lambda).
\]
This follows from applying \cite[Theorem 1.13 (ii)]{ECoD} to the small v-stack $X=\ast$. In fact, for any complete algebraically closed field $C$, one has $D_\et(\Spa C,\Lambda) = D(\Lambda)$ and there is a sequence
$$
D(\Lambda) \lto D_{\et}(\ast,\Lambda) \xrightarrow{\ \text{fully faithful}\ } D_{\et}(\Spa C,\Lambda)=D(\Lambda)
$$
whose composite is the identity, and $D(\Lambda)\to D_{\et}(\ast, \Lambda)$ is thus an equivalence.

\begin{proof} We start by constructing a functor
\[
D(G,\Lambda)\to D_\et([\ast/\underline{G}],\Lambda)
\]
compatible with the derived tensor product and the forgetful functors. For this, one first constructs a functor from the category of smooth representations of $G$ on $\Lambda$-modules to the heart of $D_\et([\ast/\underline{G}],\Lambda)$; note that this heart is a full subcategory of the heart of $D([\ast/\underline{G}]_v,\Lambda)$, which is the category of v-sheaves on $[\ast/\underline{G}]$. Now one can send a smooth $G$-representation $V$ to the v-sheaf $\mathcal F_V$ on $[\ast/\underline{G}]$ that takes a perfectoid space $X$ with a $\underline{G}$-torsor $\tilde{X}\to X$ to the set of all continuous $G$-equivariant maps from $|\tilde{X}|$ to $V$. In fact one checks that for any perfectoid space $S$ and any locally profinite set $A$,
$$
|S\times \underline{A}|=|S|\times A,
$$
and thus $|\tilde{X}|$ has a continuous $G$-action. As v-covers induce quotient maps by \cite[Proposition 12.9]{ECoD}, this is indeed a v-sheaf. Moreover, after pullback along $\ast\to [\ast/\underline{G}]$, it is given by the functor which sends $X$ to the set of continuous $G$-equivariant maps from $|X|\times G = |X\times \underline{G}|$ to $V$. These are canonically the same (via restriction to $X\times \{1\}$) as continuous maps $|X|\to V$, so that $\mathcal F_V|_\ast = \underline{V}$ is the v-sheaf corresponding to $V$. As $V$ is discrete, this is a disjoint union of points, and in particular (after pullback to any $\Spa C$) an \'etale sheaf. According to \cite[Definition 14.13]{ECoD}, this implies that $\mathcal F_V\in D_\et([\ast/\underline{G}],\Lambda)$, as desired.

From now on, we will simply write $\underline{V}$ for $\mathcal F_V$. Given any complex of smooth $G$-representations $V^\bullet$, one can form the corresponding complex $\underline{V}^\bullet$ of v-sheaves on $[\ast/\underline{G}]$, which defines an object of $D_\et([\ast/\underline{G}],\Lambda)\subset D([\ast/\underline{G}]_v,\Lambda)$ (using \cite[Proposition 14.16]{ECoD}), giving the desired functor $D(G,\Lambda)\to D_\et([\ast/\underline{G}],\Lambda)$ compatible with the forgetful functors, using exactness of $V\mapsto \mathcal F_V$. One checks that this functor is compatible with derived tensor products by unraveling the definitions.

To check whether the functor is an equivalence, we may by \cite[Theorem 1.13 (ii)]{ECoD} replace $[\ast/\underline{G}]$ by its base change $[\Spa C/\underline{G}] = [\ast/\underline{G}]\times \Spa C$, where $C$ is some complete algebraically closed nonarchimedean field.

For the v-stack $X=[\Spa C/\underline{G}]$, we can also consider its \'etale site $X_\et\subset X_v$ consisting of all $Y\in X_v$ which are \'etale (and locally separated) over $X$. This recovers a classical site.

\begin{lemma} The \'etale site $X_\et$ is equivalent to the category $G\text{-}\Set$ of discrete $G$-sets, via sending a discrete set $S$ with continuous $G$-action to $[\underline{S}\times \Spa C/\underline{G}]$.
\end{lemma}

\begin{proof} It is clear that the functor $S\mapsto [\underline{S}\times \Spa C/\underline{G}]$ maps to $X_\et\subset X_v$ (as the pullback to $\Spa C$ is given by $\underline{S}\times \Spa C$), and is fully faithful. Conversely, if $Y\to X=[\Spa C/\underline{G}]$ is \'etale, then the pullback of $Y$ to $\Spa C$ is a discrete set, on which $\underline{G}$ acts continuously, giving the descent datum defining $Y$.
\end{proof}

\begin{lemma} There is a natural equivalence $D(G,\Lambda)\simeq D(G\text{-}\Set,\Lambda)$, such that the following diagram commutes
\[\xymatrix{
D(G,\Lambda)\ar[r]^\cong\ar[d] & D(G\text{-}\Set,\Lambda)\ar^-{\cong}[r] & D([\Spa C/\underline{G}]_\et,\Lambda)\ar[d]\\
D_\et([\ast/\underline{G}],\Lambda)\ar@{^(->}[rr] && D_\et([\Spa C/\underline{G}],\Lambda).
}\]
\end{lemma}

\begin{proof} It is enough to give an equivalence of abelian categories between smooth $G$-representations on $\Lambda$-modules, and sheaves of $\Lambda$-modules on discrete $G$-sets. The construction of the functor is as before: Send a smooth representation $V$ to the sheaf sending some discrete $G$-set $S$ to the $\Lambda$-module of continuous $G$-equivariant maps $S\to V$. This functor is clearly fully faithful. But any sheaf of $\Lambda$-modules $\mathcal F$ on $G\text{-}\Set$ comes from the smooth $G$-representation $V=\varinjlim_{H\subset G} \mathcal F(G/H)$, where $H$ runs over all open subgroups of $G$. One directly verifies that the diagram commutes.
\end{proof}

It remains to see that the natural functor
\[
D([\Spa C/\underline{G}]_\et,\Lambda)\to D_\et([\Spa C/\underline{G}],\Lambda)
\]
is an equivalence. We claim that this reduces to the case that $G$ is pro-$p$: We first reduce fully faithfulness to this case. For this, we have to see that, if $\lambda: X_v\to X_\et$ denotes the map of sites, then for any $A\in D([\Spa C/\underline{G}]_\et,\Lambda)$, the natural map
\[
A\to R\lambda_\ast \lambda^\ast A
\]
is an equivalence. This can be checked locally on $[\Spa C/\underline{G}]_\et$, meaning that we can replace $G$ by an open pro-$p$-subgroup. Similarly, for essential surjectivity, one needs to see that for all $B\in D_\et([\Spa C/\underline{G}],\Lambda)$, the map $\lambda^\ast R\lambda_\ast B\to B$ is an equivalence, which can again be checked locally.

Thus, we can assume that $G$ is pro-$p$. Note that $([\Spa C/\underline{G}]_\et,\Lambda)$ is locally of cohomological dimension $0$, as there is no continuous group cohomology of pro-$p$-groups on $\Lambda$-modules if $n\Lambda=0$ for $n$ prime to $p$. This implies (cf.~\cite[Tag 0719]{StacksProject}) that $D([\Spa C/\underline{G}]_\et,\Lambda)$ is left-complete. As $D_\et([\Spa C/\underline{G}],\Lambda)$ is also left-complete by \cite[Proposition 14.11]{ECoD}, it is enough to see that the functor
\[
D^+([\Spa C/\underline{G}]_\et,\Lambda)\to D^+_\et([\Spa C/\underline{G}],\Lambda)
\]
is an equivalence. First, we check fully faithfulness, i.e.~that the unit $\id\to R\lambda_\ast \lambda^\ast$ of the adjunction is an equivalence. For this, it is enough to see that for any \'etale sheaf of $\Lambda$-modules, i.e.~any smooth $G$-representation $V$, one has
\[
R\Gamma([\Spa C/\underline{G}]_v,\underline{V}) = V^G,
\]
i.e.~its $H^0$ is $V^G$ and there are no higher $H^i$. However, one can compute v-cohomology using the Cech nerve for the cover $\Spa C\to [\Spa C/\underline{G}]$, which produces the complex of continuous cochains, giving the desired result.

Finally, for essential surjectivity, it is now enough to check on the heart. But if $\mathcal F$ is a v-sheaf on $[\Spa C/\underline{G}]$ whose pullback to $\Spa C$ is an \'etale sheaf, then this pullback is a disjoint union of points, thus separated and \'etale, and therefore $\mathcal F$ is itself a v-stack which is \'etale over $[\Spa C/\underline{G}]$, and so defines an object in the topos $[\Spa C/\underline{G}]_\et$.
\end{proof}

\begin{corollary}\label{cor:dualrepresentation} The operation
\[
R\sHom_\Lambda(-,\Lambda): D_\et([\ast/\underline{G}],\Lambda)^\op\to D_\et([\ast/\underline{G}],\Lambda)
\]
corresponds to the derived smooth duality functor
\[
A\mapsto (A^\ast)^\sm: D(G,\Lambda)^\op\to D(G,\Lambda)
\]
induced on derived categories by the left-exact smooth duality functor
\[
V\mapsto (V^\ast)^\sm = \{f: V\to \Lambda\mid \exists H\subset G\ \mathrm{open}\ \forall h\in H, v\in V: f(hv)=f(v)\}.
\]
\end{corollary}

\begin{proof} The operation $A\mapsto (A^\ast)^\sm$ on $D(G,\Lambda)$ satisfies the adjunction
\[
\Hom_{D(G,\Lambda)}(B,(A^\ast)^\sm) = \Hom_{D(G,\Lambda)}(B\dotimes_\Lambda A,\Lambda)
\]
for all $B\in D(G,\Lambda)$. As $R\sHom_\Lambda(-,\Lambda)$ is characterized by the similar adjunction in $D_\et([\ast/\underline{G}],\Lambda)$ and the equivalence is symmetric monoidal, we get the result.
\end{proof}

\section{\'Etale sheaves on strata}

We want to describe $D_\et(\Bun_G,\Lambda)$ via its strata $\Bun_G^b$. For this, we need the following result saying roughly that connected Banach--Colmez spaces are ``contractible''.

\begin{proposition}\label{prop:banachcolmezfullyfaithful} Let $f:S^\prime\to S$ be a map of small v-stacks that is a torsor under $\BC(\mathcal E)$ resp.~$\BC(\mathcal E[1])$, where $\mathcal E$ is a vector bundle on $X_S$ that is everywhere of positive (resp.~negative) slopes. Then the pullback functor
\[
f^\ast: D_\et(S,\Lambda)\to D_\et(S^\prime,\Lambda)
\]
is fully faithful.
\end{proposition}

\begin{proof} By descent \cite[Proposition 17.3, Remark 17.4]{ECoD}, the problem is v-local on $S$, and in particular one can assume that the torsor is split. In the positive case, we can use Corollary~\ref{cor:niceshortexseq}~(iv) to find pro-\'etale locally on $S$ (and after replacing $\mathcal E$ by the direct sum with another bundle) a short exact sequence
\[
0\to \mathcal O_{X_S}(\tfrac 1{2r})^{m'}\to \mathcal O_{X_S}(\tfrac 1r)^m\to \mathcal E\to 0,
\]
inducing a similar sequence on Banach--Colmez spaces. This reduces us to the case $\mathcal E=\mathcal O_{X_S}(\tfrac 1n)$ for some $n$ (as then pullback under $\BC(\mathcal O_{X_S}(\tfrac 1r)^m)\to S$ is fully faithful, as is pullback under $\BC(\mathcal O_{X_S}(\tfrac 1r)^m)\to S'=\BC(\mathcal E)$). In that case, $\BC(\mathcal E)$ is a $1$-dimensional perfectoid open unit ball over $S$ by Proposition~\ref{prop:standardbanachcolmez}~(iv), in particular cohomologically smooth. It suffices to see that $Rf^!$ is fully faithful, for which it suffices that for all $A\in D_\et(S,\Lambda)$, the adjunction map
\[
Rf_! Rf^! A\to A
\]
is an equivalence. But note that both $Rf_!$ and $Rf^!$ commute with any base change by \cite[Proposition 22.19, Proposition 23.12]{ECoD}. Thus, we may by passage to stalks reduce to the case $S=\Spa(C,C^+)$ where $C$ is a complete algebraically closed nonarchimedean field and $C^+\subset C$ an open and bounded valuation subring, and (as we reduced to the stalk) we only need to check the statement on global sections. If the stalk of $A$ at the closed point $s\in S$ is zero, then the same holds true for $Rf_! Rf^! A$ as $Rf^! A$ agrees with $f^\ast A$ up to twist, so this follows from the commutation of $Rf_!$ with extension by zero. This allows us to reduce to the case that $A$ is constant, and then as both $Rf^!$ and $Rf_!$ commute with all direct sums, even to the case $A=\Lambda$. Thus, we are reduced to the computation of the cohomology of the perfectoid open unit disc, which is known.

The case of negative Banach--Colmez spaces follows by taking an exact sequence
\[
0\to \mathcal E\to \mathcal O_{X_S}(d)^m\to \mathcal G\to 0
\]
with $d>0$ (and hence $\mathcal G$ of everywhere positive slopes) and using the exact sequence
\[
0\to \BC(\mathcal O_{X_S}(d)^m)\to \BC(\mathcal G)\to \BC(\mathcal E[1])\to 0.\qedhere
\]
\end{proof}

Now we can formulate the desired result.

\begin{proposition}\label{prop:DetBunstrata} For any $b\in B(G)$, the map
\[
\Bun_G^b=[\ast/\tilde{G}_b]\to [\ast/\underline{G_b(E)}]
\]
induces via pullback an equivalence
\[
D_\et(G_b(E),\Lambda)\cong D_\et([\ast/\underline{G_b(E)},\Lambda])\cong D_\et([\ast/\tilde{G}_b],\Lambda).
\]
Moreover, for any complete algebraically closed nonarchimedean field $C/k$, the map
\[
D_\et([\ast/\tilde{G}_b],\Lambda)\to D_\et([\Spa C/\tilde{G}_b],\Lambda)
\]
is an equivalence.
\end{proposition}

\begin{proof} Using \cite[Theorem 1.13]{ECoD} and Theorem~\ref{thm:Detclassifying}, it is enough to prove that the functor
\[
D_\et([\Spa C/\underline{G_b(E)}],\Lambda)\to D_\et([\Spa C/\tilde G_b],\Lambda)
\]
is an equivalence. For this, it is enough to prove that pullback under the section $[\Spa C/\underline{G_b(E)}]\to [\Spa C/\tilde G_b]$ induces a fully faithful functor
\[
D_\et([\Spa C/\tilde G_b],\Lambda)\to D_\et([\Spa C/\underline{G_b(E)}],\Lambda),
\]
which follows from Proposition~\ref{prop:gradedpiecesoffancyGb} and Proposition~\ref{prop:banachcolmezfullyfaithful}.
\end{proof}

We see that $D_\et(\Bun_G,\Lambda)$ is glued from the categories $D_\et(\Bun_G^b,\Lambda)\cong D(G_b(E),\Lambda)$, which are entirely representation-theoretic.\footnote{It would be very interesting to understand the gluing of these categories in terms of pure representation theory.} In particular, this implies that the base field plays no role:

\begin{corollary}\label{cor:DetBun} For any complete algebraically closed nonarchimedean field $C$ and any locally closed substack $U\subset \Bun_G$, the functor
\[
D_\et(U,\Lambda)\to D_\et(U\times \Spa C,\Lambda)
\]
is an equivalence of categories.
\end{corollary}

Although this seems like a purely technical result, it will actually play a key role when we study Hecke operators.

\begin{proof} Fully faithfulness holds true by \cite[Theorem 1.13 (ii)]{ECoD}. To see that it is an equivalence of categories, it is enough to check on all quasicompact locally closed substacks $U\subset \Bun_G$. These are stratified into finitely many locally closed substacks $\Bun_G^b\subset \Bun_G$, and thus any object of $D_\et(U\times \Spa C,\Lambda)$ is filtered by objects $!$-extended from $D_\et(\Bun_G^b\times \Spa C,\Lambda)$. By fully faithfulness, it suffices to show that all the graded pieces lie in the essential image of $D_\et(U,\Lambda)$. Now the result follows from Proposition~\ref{prop:DetBunstrata} (and compatibility of $!$-extension with base change to $\Spa C$).
\end{proof}

\section{Local charts}\label{sec:localcharts}

For any $b\in B(G)$ we wish to construct a chart
\[
\pi_b: \mathcal M_b\to \Bun_G
\]
whose image contains $\Bun_G^b$, such that {\it $\pi_b$ is separated, representable in locally spatial diamonds and cohomologically smooth}, and whose geometry can be understood explicitly.

\begin{example}\label{ex:GL2O1O} Before we discuss the general case, let us briefly discuss the first interesting case, namely $G=\mathrm{GL}_2$ and the non-basic element $b$ corresponding to $\mathcal O(1)\oplus \mathcal O$. In that case, we let $\mathcal M_b$ be the moduli space of extensions
\[
0\to \mathcal L\to \mathcal E\to \mathcal L'\to 0
\]
where $\mathcal L$ is of degree $0$ and $\mathcal L'$ is of degree $1$. Mapping such an extension to $\mathcal E$ defines the map $\mathcal M_b\to \Bun_G$.

Note that there is a natural $\underline{E^\times}\times \underline{E^\times}$-torsor $\tilde{\mathcal M}_b\to \mathcal M_b$, parametrizing isomorphisms $\mathcal L\cong \mathcal O$ and $\mathcal L'\cong \mathcal O(1)$. On the other hand, it is clear that $\tilde{\mathcal M}_b=\BC(\mathcal O(-1)[1])$ is a negative absolute Banach--Colmez space; thus, $\mathcal M_b$ is very explicit.

Any extension $\mathcal E$ parametrized by $\mathcal M_b$ is either isomorphic to $\mathcal O\oplus \mathcal O(1)$, or to $\mathcal O(\tfrac 12)$. The fibres of $\pi_b: \mathcal M_b\to \Bun_{\mathrm{GL}_2}$ over a rank $2$ bundle $\mathcal E$ are given by an open subset of the projectivized Banach--Colmez space $(\BC(\mathcal E)\setminus \{0\})/\underline{E^\times}$. Thus, these fibres interpolate between $(\BC(\mathcal O(\tfrac 12))\setminus \{0\})/\underline{E^\times}$, which is cohomologically smooth by Proposition~\ref{prop:absolutebanachcolmez}, and an open subset of
\[
(\BC(\mathcal O\oplus \mathcal O(1))\setminus \{0\})/\underline{E^\times} = (\underline{E}\times \Spd k\powerseries{t^{1/p^\infty}}\setminus \{(0,0)\})/\underline{E^\times}.
\]
That open subspace is still cohomologically smooth, although $\underline{E}\times \Spd k\powerseries{t^{1/p^\infty}}$ is not --- the quotient by $\underline{E^\times}$ gets rid of the disconnected nature of the space. In this case, and in fact in complete generality for all $b\in B(\mathrm{GL}_n)$, one can actually check cohomological smoothness of $\pi_b$ by hand. To handle the general case, we had to prove the Jacobian criterion, Theorem~\ref{thm:jacobiancriterion}.
\end{example}

Coming back to the general case, we can in fact construct all $\mathcal M_b$ together, as follows.

\begin{definition} The v-stack $\mathcal M$ is the moduli stack taking $S\in \Perf_k$ to the groupoid of $G$-bundles $\mathcal E$ on $X_S$ together with an increasing separated and exhaustive $\mathbb Q$-filtration $(\rho_\ast \mathcal E)^{\leq \lambda}\subset \rho_\ast \mathcal E$ (ranging over algebraic representations $\rho: G\to \GL_n$, and compatible with exact sequences and tensor products) on the corresponding fibre functor such that (letting $(\rho_\ast\mathcal E)^{<\lambda}=\bigcup_{\lambda'<\lambda} (\rho_\ast \mathcal E)^{\leq \lambda'}$) the quotient
\[
(\rho_\ast \mathcal E)^\lambda = (\rho_\ast \mathcal E)^{\leq \lambda}/(\rho_\ast \mathcal E)^{<\lambda}
\]
is a semistable vector bundle of slope $\lambda$, for all $\lambda\in \mathbb Q$ and representations $\rho: G\to \GL_n$.
\end{definition}

Note that by passing to the associated graded, $\mathcal M$ maps to the moduli stack of $G$-bundles in the category of $\mathbb Q$-graded vector bundles on $X_S$ where the weight $\lambda$ piece is semistable of slope $\lambda$. By Proposition~\ref{prop:GBunHNsplit}, this is isomorphic to $\bigsqcup_{b\in B(G)} [\ast/\underline{G_b(E)}]$. In particular $\mathcal M$ decomposes naturally into a disjoint union
\[
\mathcal M = \bigsqcup_{b\in B(G)} \mathcal M_b,
\]
and for each $b\in B(G)$, we have natural maps
\[
q_b: \mathcal M_b\to [\ast/\underline{G_b(E)}].
\]

\begin{example}
When $G=\GL_n$, $\mathcal{M}$ sends $S$ to the groupoid of filtered vector bundles $0=\Fil_0\E\subsetneq \Fil_1\E\subsetneq \dots \subsetneq \Fil_r\E=\E$ for some $r$, such that $\Fil_{i+1} \E /\Fil_i \E$ is semistable and the slopes $(\mu (\Fil_{i+1} \E/\Fil_i \E))_{0\leq i<r}$ form an increasing sequence (the opposite condition to the one defining the Harder--Narasimhan filtration of a vector bundle). The maps $q_b$ send such a vector bundle to the graded vector bundle $\bigoplus_{i=0}^{r-1} \Fil_{i+1}\E/\Fil_i\E$.
\end{example}

\begin{example}
Suppose $G$ is quasisplit. Let $M_b$ be the centralizer of the slope morphism, which is a Levi subgroup. Let $P_b$ be the parabolic subgroup with Levi $M_b$ such that the weights of $\nu_b$ on $\Lie\, P_b$ are positive. There is a diagram 
$$
\begin{tikzcd}
\Bun_{P_b}\ar[d] \ar[r] & \Bun_G \\
\Bun_{M_b}
\end{tikzcd}
$$ 
induced by the inclusion $P_b\subset G$ and the quotient map $P_b \to M_b$. There is a cartesian diagram
$$
\begin{tikzcd}
\mathcal{M}_b \ar[r]\ar[d] & \Bun_{P_b} \ar[d] \\
 \Bun_{M_b}^b \ar[r,hook] & \Bun_{M_b}. 
\end{tikzcd}
$$
\end{example}

\begin{proposition}\label{prop:MbtoGbnice} For any $b\in B(G)$, the map
\[
q_b: \mathcal M_b\to [\ast/\underline{G_b(E)}].
\]
is partially proper, representable in locally spatial diamonds, and cohomologically smooth, of dimension $\langle 2\rho,\nu_b\rangle$. In fact, after pullback along $\ast\to [\ast/\underline{G_b(E)}]$, it is a successive torsor under negative Banach--Colmez spaces.

In particular, $\mathcal M_b$ is a cohomologically smooth Artin v-stack, of dimension $\langle 2\rho,\nu_b\rangle$.
\end{proposition}

\begin{proof} It suffices to check everything after pullback by the $v$-cover $\ast\to [\ast/\underline{G_b(E)}]$, inducing $\tilde{\mathcal{M}}_b\to \mathcal{M}_b$. Let $H\to X_S$ be the automorphism group of $\E_b\to X_S$, see Proposition~\ref{prop:filtrations}, the pure inner twisting of $G\times X_S$ by $\E_b$. This is equipped with a parabolic subgroup $H^{\leq 0}$, and moreover a filtration $(H^{\leq \lambda})_{\lambda \leq 0}$ with unipotent radical $H^{<0}$. This is the opposite parabolic subgroup to the one used in the proof of Proposition~\ref{prop:gradedpiecesoffancyGb}. Then $\tilde{\mathcal{M}}_b (S)$ is identified with the set of $H^{<0}$-torsors on $X_S$. The result is deduced using the description of the graded pieces of $(H^{\leq \lambda})_{\lambda<0}$ as vector bundles of negative slopes.
\end{proof}

We first prove some structural results about $\mathcal M_b$ and its universal $\underline{G_b(E)}$-torsor $\tilde{\mathcal M}_b\to \mathcal M_b$. {\it A general theme here is the subtle distinction between the absolute property of being a (locally spatial) diamond (which $\tilde{\mathcal M}_b$ is not, but it has a large open part $\tilde{\mathcal M}_b^\circ\subset \tilde{\mathcal M}_b$ that is) and the relative notion of $\tilde{\mathcal M}_b\to \ast$ being representable in (locally spatial) diamonds (which $\tilde{\mathcal M}_b$ is), and some related subtle distinctions on absolute and relative quasicompactness.}

\begin{proposition}\label{prop:basicMb} The map $\mathcal M_b\to [\ast/\underline{G_b(E)}]$ has a section $[\ast/\underline{G_b(E)}]\to \mathcal M_b$ given by the closed substack where $\mathcal E$ is (at every geometric point) isomorphic to $\mathcal E_b$, in which case $(\rho_\ast \mathcal E)^{\leq \lambda}\subset \rho_\ast \mathcal E$ is a splitting of the Harder--Narasimhan filtration of $\rho_\ast \mathcal E$ for all representations $\rho: G\to \GL_n$.

Consider the open complement $\mathcal M_b^\circ = \mathcal M_b\setminus [\ast/\underline{G_b(E)}]$, with preimage $\tilde{\mathcal M}_b^\circ = \tilde{\mathcal M}_b\setminus \{\ast\}$. Then $\tilde{\mathcal M}_b^\circ$ is a spatial diamond.

Moreover, if $U_\pi:=\nu_b^N(\pi)\in G_b(E)$ for any large enough $N$ (so that $\nu_b^N: \mathbb G_m\to G_b$ is a well-defined cocharacter), then $\tilde{\mathcal M}_b^\circ/U_\pi^{\mathbb Z}\to \ast$ is proper.
\end{proposition}

\begin{proof} To check that the substack where $\mathcal E$ is at every geometric point isomorphic to $\mathcal E_b$ is closed, note that by semicontinuity it suffices to see that everywhere on $\mathcal M_b$, the Newton point of $\mathcal E$ is bounded by $b$. By \cite[Lemma 2.2 (iv)]{RapoportRichartz}, this reduces to the case of $G=\GL_n$, where it is a simple consequence of the Harder--Narasimhan formalism (the Newton polygon of an extension is always bounded by the Newton polygon of the split extension). On this closed substack, $\mathcal E$ has two transverse filtrations, given by $(\rho_\ast \mathcal E)^{\leq \lambda}$ and the Harder--Narasimhan filtration; it follows that $\mathcal E$ upgrades to a $G$-bundle in $\mathbb Q$-graded vector bundles on the Fargues--Fontaine curve, with the weight $\lambda$-piece semistable of slope $\lambda$. We see that this gives the desired section of $\mathcal M_b\to [\ast/\underline{G_b(E)}]$, using again Proposition~\ref{prop:GBunHNsplit}.

We claim that the action of $U_\pi$ on $|\tilde{\mathcal M}_b\times_k \Spa k\laurentseries{t}|$ satisfies the hypotheses of Lemma~\ref{lem:contractinglocspec}, with fixed point locus given by the closed subspace $\ast$ considered in the previous paragraph. Writing $\tilde{\mathcal M}_b$ as a successive extension of Banach--Colmez spaces, it is clear that for all $x\in |\tilde{\mathcal M}_b\times_k \Spa k\laurentseries{t}|$ which are not in the closed substack, the sequence $U_\pi^{-n}(x)$ leaves any quasicompact open subspace for large $n$: Look at the first step in the successive extensions where $x$ does not project to the origin. Then $x$ gives an element in the fiber over the origin, which is a negative Banach--Colmez space, on which $U_\pi=\nu_b^N(\pi)$ acts via a positive power of $\pi$; thus $U_\pi^{-n}(x)$ leaves any quasicompact open subspace of this Banach--Colmez space. In particular, it follows that the fixed points locus of $U_\pi$ is precisely the origin. To apply Lemma~\ref{lem:contractinglocspec}, it remains to see that for all $x\in |\tilde{\mathcal M}_b\times_k \Spa k\laurentseries{t}|$ and quasicompact open neighborhoods $U$ of the origin, one has $\gamma^m(x)\in U$ for all sufficiently large $m$. 
This can be reduced to the case of $\GL_n$ by the Tannakian formalism, so assume $G=\GL_n$ for this argument. Now fix a map $f: \Spa(C,C^+)\to \tilde{\mathcal M}_b\times_k \Spa k\laurentseries{t}$ having $x$ in its image; it suffices to construct a map
\[
\Spa(C,C^+)\times \underline{\mathbb N\cup\{\infty\}}\to \tilde{\mathcal M}_b\times_k \Spa k\laurentseries{t}
\]
whose restriction to $\Spa(C,C^+)\times \{0\}$ is $f$, whose restriction to $\Spa(C,C^+)\times \{\infty\}$ maps to $\ast\subset \tilde{\mathcal M}_b$, and which is equivariant for the $\gamma$-action, with $\gamma$ acting on the left via shift on the profinite set $\mathbb N\cup \{\infty\}$. The map $f$ classifies some $\mathbb Q$-filtered vector bundle $\mathcal E^{\leq \lambda}\subset \mathcal E$ of rank $n$ on $X_C$ with $\bigoplus_\lambda \mathcal E^\lambda\cong \mathcal E_b$ as $\mathbb Q$-graded vector bundles. After pullback to $Y_{C,[1,q]}$, the filtration is split, so we can find an isomorphism
\[
\alpha: \mathcal E|_{Y_{C,[1,q]}}\cong \mathcal E_b|_{Y_{C,[1,q]}}
\]
of $\mathbb Q$-filtered vector bundles on $Y_{C,[1,q]}$, such that $\alpha$ reduces on graded pieces to the given identification. The descent datum is now given by some isomorphism of $\mathbb Q$-filtered vector bundles
\[
\beta: \phi^\ast (\mathcal E_b|_{Y_{C,[q,q]}})\cong \mathcal E_b|_{Y_{C,[1,1]}}
\]
that reduces to the standard Frobenius on graded pieces. In other words, $\beta$ is the standard Frobenius on $\mathcal E_b$ multiplied by some
\[
\beta': \mathcal E_b|_{Y_{C,[1,1]}}\cong \mathcal E_b|_{Y_{C,[1,1]}}
\]
and with respect to the $\mathbb Q$-grading on $\mathcal E_b$, the map $\beta'$ is the identity plus a lower triangular matrix.

The action of $\gamma$ replaces $\beta'$ by its $U_\pi$-conjugate. This multiplies all lower triangular entries by powers of $\pi$, so
\[
(\beta',\gamma(\beta'),\gamma^2(\beta'),\ldots,1)
\]
composed with the standard Frobenius defines an isomorphism
\[
\phi^\ast (\mathcal E_b|_{Y_{S,[q,q]}})\cong \mathcal E_b|_{Y_{S,[1,1]}}
\]
where $S=\Spa(C,C^+)\times \underline{\mathbb N\cup \{\infty\}}$. Using this as a descent datum defines a $\mathbb Q$-filtered vector bundle on $X_S$ defining the required map
\[
S=\Spa(C,C^+)\times \underline{\mathbb N\cup\{\infty\}}\to \tilde{\mathcal M}_b\times_k \Spa k\laurentseries{t}.
\]
This finishes the verification of the hypotheses of Lemma~\ref{lem:contractinglocspec}. 

It is clear from the definition that $\tilde{\mathcal M}_b^\circ\to \ast$ is partially proper (as the theory of vector bundles on the Fargues--Fontaine curve does not depend on $R^+$). Thus, to show that $\tilde{\mathcal M}_b^\circ/U_\pi^{\mathbb Z}\to \ast$ is proper, it suffices to see that the map is quasicompact, which can be checked after base change to $\Spa k\laurentseries{t}$; then it follows from the previous discussion and Lemma~\ref{lem:contractinglocspec}.

It remains to see that $\tilde{\mathcal M}_b^\circ$ is a spatial diamond. For this, we pick a representative $b\in G(\breve E)$ of the $\sigma$-conjugacy class that is decent in the sense of \cite[Definition 1.8]{RapoportZink}. In particular, $b\in G(E_s)$ for some unramified extension $E_s|E$ of degree $s$, and $\tilde{\mathcal M}_b^\circ$ is already defined over $\Perf_{\mathbb F_{q^s}}$. Let $\Frob_s$ be the Frobenius $x\mapsto x^{q^s}$. As $b$ is decent, the action of $U_\pi=\nu_b^N(\pi)$ on $\tilde{\mathcal M}_b^\circ$ agrees with the action of a power of $\Frob_s$ for $N$ large enough, by the decency equation for $b$. We know that
\[
\tilde{\mathcal M}_b^\circ/U_\pi^{\mathbb Z}\times_k \Spa k\laurentseries{t}
\]
is a spatial diamond (as it is proper over $\Spa k\laurentseries{t}$). Replacing $U_\pi$ by $\Frob_s$, and moving the quotient by Frobenius to the other factor (which is allowed as the absolute Frobenius acts trivially on topological spaces) one sees that also
\[
\tilde{\mathcal M}_b^\circ\times_k \Spa k\laurentseries{t}/\Frob_s^{\mathbb Z}
\]
is a spatial diamond. But $\Spa k\laurentseries{t}/\Frob_s^{\mathbb Z}\to \ast$ is proper and cohomologically smooth. Thus, Lemma~\ref{lemma:critere spatial v faisceau} shows that it is a spatial v-sheaf. By \cite[Theorem 12.18]{ECoD}, to see that $\tilde{\mathcal M}_b^\circ$ is a spatial diamond, it suffices to check on points. Writing $\tilde{\mathcal M}_b$ as a successive extension of Banach--Colmez spaces, any point in $\tilde{\mathcal M}_b^\circ$ has a minimal step where it does not map to the origin. Then its image is a nontrivial point of an absolute Banach--Colmez space, and a punctured absolute Banach--Colmez space is a diamond by Proposition~\ref{prop:absolutebanachcolmez}; the result follows.
\end{proof}

The following theorem gives the desired local charts for $\Bun_G$; its proof is based on the Jacobian criterion for (cohomological) smoothness, Theorem~\ref{thm:jacobiancriterion}.

\begin{theorem}\label{thm:chartbsmooth} The map $\pi_b: \mathcal M_b\to \Bun_G$ forgetting the filtration is partially proper, representable in locally spatial diamonds, and cohomologically smooth of $\ell$-dimension $\langle 2\rho,\nu_b\rangle$.

The image of $\pi_b$ is open, and consists exactly of the set of points of $|\Bun_G|$ specializing to $b$.
\end{theorem}

\begin{proof} First, we show that
\[
\mathcal M = \bigsqcup_b \mathcal M_b\to \Bun_G
\]
is partially proper, representable in locally spatial diamonds, and cohomologically smooth. Let $S\to \Bun_G$ be some map for an affinoid perfectoid space $S$, given by some $G$-bundle $\mathcal E$ on $X_S$. Then $\mathcal M\times_{\Bun_G} S$ parametrizes $\mathbb Q$-filtrations on $\mathcal E$ whose associated graded $\mathbb Q$-bundle has the degree $\lambda$-part semistable of slope $\lambda$. Now $\mathbb Q$-filtrations are parametrized by sections of a smooth scheme $Z=\mathcal E\times^G \mathrm{Fl}\to X_S^{\mathrm{alg}}$, for some smooth scheme $\mathrm{Fl}$ over $E$ with $G$-action (classifying such $\mathbb Q$-filtrations on the forgetful functor $\mathrm{Rep}_E G\to \mathrm{Vect}_E$). Here $\mathrm{Fl}$ is a disjoint union of projective schemes. The condition on the associated graded bundle is an open condition by openness of the semistable locus. Using the Jacobian criterion Theorem~\ref{thm:jacobiancriterion}, it remains to see that $\mathcal M\subset \mathcal M_Z$ is contained in the smooth locus $\mathcal M_Z^{\mathrm{sm}}\subset \mathcal M_Z$. The tangent bundle of $Z\to X_S^{\mathrm{alg}}$ is the $\mathcal E$-twisted form of the tangent bundle of $\mathrm{Fl}/E$. By the usual dscription of tangent bundles of flag varieties, this tangent bundle admits a $\mathbb Q_{>0}$-filtration whose $\lambda$-th graded piece is given by the graded piece of weight $\lambda$ of $\mathcal E\times^G \mathrm{Lie} G$ (with respect to its universal $\mathbb Q$-filtration). After pullback to a section in $\mathcal M$, these are semistable of slope $\lambda>0$, as desired.

Knowing cohomological smoothness of $\pi_b$, its $\ell$-dimension can be computed as the difference of those of $\mathcal M_b$ and $\Bun_G$. For the final assertion, note that by cohomological smoothness, the image must be open. Conversely, if $b'\in |\Bun_G|$ is in the image, then it follows from the explicit degeneration used in the proof of Proposition~\ref{prop:basicMb} that $\Bun_G^b$ lies in the closure of $\Bun_G^{b'}$, giving the desired result.
\end{proof}

This final assertion was proved in the case of $G=\GL_n$ by Hansen \cite{HansenGLn}, as a step towards the identification of $|\Bun_G|$; indeed, \cite{Hansenetal} determines the image of $\pi_b$ for $G=\GL_n$.

\section{Compact generation}

{\it We are now revisiting the notion of a finite type smooth representation in terms of $D_{\et} (\Bun_G,\Lambda)$.}

The goal of this section is to prove the following theorem. As above, we fix some coefficient ring $\Lambda$ such that $n\Lambda=0$ for some $n$ prime to $p$.

\begin{theorem}\label{thm:compactgeneration} For any locally closed substack $U\subset \Bun_G$, the triangulated category $D_\et(U,\Lambda)$ is compactly generated. An object $A\in D_\et(U,\Lambda)$ is compact if and only if for all $b\in B(G)$ contained in $U$, the restriction
\[
i^{b\ast} A\in D_\et(\Bun_G^b,\Lambda)\cong D(G_b(E),\Lambda)
\]
along $i^b: \Bun_G^b\subset \Bun_G$ is compact, and zero for almost all $b$. Here, compactness in $D(G_b(E),\Lambda)$ is equivalent to lying in the thick triangulated subcategory generated by $c\text-\mathrm{Ind}_K^{G_b(E)} \Lambda$ as $K$ runs over open pro-$p$-subgroups of $G_b(E)$.
\end{theorem}

To prove the theorem, we exhibit a class of compact projective generators. The key result is that $\widetilde{\mathcal M}_b$ behaves like a strictly local scheme; in some vague sense, {\it $\widetilde{\mathcal M}_b$ is the strict henselization of $\Bun_G$ at $b$}.

\begin{proposition}\label{prop:tildeMbstrictlyhenselian} Let $b\in B(G)$. For any $A\in D_\et(\widetilde{\mathcal M}_b,\Lambda)$ with stalk $A_0 = i^\ast A\in D_\et(\ast,\Lambda)\cong D(\Lambda)$ at the closed point $i: \ast\subset \widetilde{\mathcal M}_b$, the map
\[
R\Gamma(\widetilde{\mathcal M}_b,A)\to A_0
\]
is an isomorphism. In particular, $R\Gamma(\widetilde{\mathcal M}_b,-)$ commutes with all direct sums.
\end{proposition}

\begin{proof} Replacing $A$ by the cone of $A\to i_\ast A_0$, we can assume that $A=j_! A'$ for some $A'\in D_\et(\widetilde{\mathcal M}_b^\circ,\Lambda)$. We have to see that
\[
R\Gamma(\widetilde{\mathcal M}_b,j_! A')=0.
\]
But this follows from Theorem~\ref{thm:partialcompactsupportvanishing} (applied with $X=\widetilde{\mathcal M}_b^\circ$ and $S=\Spa k\laurentseries{t}$, noting that base change along $S\to \ast$ follows from smooth base change), using that the partial compactification $\widetilde{\mathcal M}_b^\circ\subset \widetilde{\mathcal M}_b$ is precisely a compactification towards one of the two ends of $\widetilde{\mathcal M}_b^\circ$, as follows from the behaviour of the Frobenius action exhibited in the proof of Proposition~\ref{prop:basicMb}.
\end{proof}

\begin{remark}\leavevmode
\begin{altenumerate}
\item
Consider the $v$-sheaf $X=\Spd (k\powerseries{x_1,\dots,x_d})$ and the quasicompact open subset
\[
U= \Spd (k\powerseries{x_1,\dots,x_d}) \setminus V(x_1,\dots,x_d)
\]
that is representable by a perfectoid space. When base changed to $S=\Spa (k\laurentseries{t})$, $X$ becomes isomorphic to an open unit disk, $U$ becomes the punctured unit disk that has two ends: the origin and the exterior of the disk. The picture is thus analogous to the preceding one with $\tilde{\mathcal{M}}_b$, and for any $A\in D_{\et}(X,\Lambda)$ one has $R\Gamma (X,A)\cong i^*A$ where $i:\Spd (k)\hookrightarrow X$ is $V(x_1,\dots,x_d)$.
\item
This applies more generally to the $v$-sheaf associated to any $W(k)$-affine formal scheme $\Spf (R)$ with $R$ an $I$-adic ring for a finitely generated ideal $I$. Then, for any $A\in D_{\et} (\Spd (R),\Lambda)$ one has $R\Gamma (\Spd (R),A) = R\Gamma (\Spd (R/I), i^*A)$ with $i:\Spd (R/I)\hookrightarrow \Spd (R)$. In particular, for $A\in D_{\et} (\Spd (R)\setminus V(I), \Lambda)$, where here $\Spd (R)\setminus V(I)$ is representable by a spatial diamond and even a perfectoid space if $R$ is a $k$-algebra,  one has $$R\Gamma ( \Spd (R)\setminus V(I), A)\cong R\Gamma ( \Spd (R/I), i^*Rj_* A).$$ Thus,  Proposition \ref{prop:tildeMbstrictlyhenselian} can be seen as a result about ``nearby cycles on the strict henselization of $\Bun_G$ at $b$''.
\end{altenumerate}
\end{remark}

\begin{corollary}\label{cor:quotientoftildeMbcompact} Let $b\in B(G)$ and let $K\subset G_b(E)$ be an open pro-$p$-subgroup. Then for any $A\in D_\et(\widetilde{\mathcal M}_b/\underline{K},\Lambda)$ with pullback $A_0=i^\ast A\in D_\et([\ast/\underline{K}],\Lambda)\cong D(K,\Lambda)$ corresponding to a complex $V$ of smooth $K$-representations, the map
\[
R\Gamma(\widetilde{\mathcal M}_b/\underline{K},A)\to R\Gamma([\ast/\underline{K}],A_0)\cong V^K
\]
is an isomorphism. In particular, $R\Gamma(\widetilde{\mathcal M}_b/\underline{K},-)$ commutes with all direct sums.
\end{corollary}

\begin{proof} This follows formally from Proposition~\ref{prop:tildeMbstrictlyhenselian} by descent along $\psi: \widetilde{\mathcal M}_b\to \widetilde{\mathcal M}_b/\underline{K}$; more precisely, by writing any $A$ as a direct summand of $\psi_\ast \psi^\ast A$.
\end{proof}

\begin{remark} Let
\[
i: [\ast/\underline{G_b(E)}]\hookrightarrow \mathcal M_b\hookleftarrow \mathcal M_b^\circ: j
\]
be the usual diagram. From the corollary, one deduces that if one regards
\[
i^\ast Rj_\ast A\in D_\et([\ast/\underline{G_b(E)}],\Lambda)\cong D(G_b(E),\Lambda)
\]
as a complex of $G_b(E)$-representations, then this is given by
\[
R\Gamma(\tilde{\mathcal M}_b^\circ,A).
\]
{\it As $\tilde{\mathcal M}_b^\circ$ is qcqs (and of finite cohomological dimension) by Proposition~\ref{prop:basicMb}, this commutes with all direct sums in $A$.}

For example, in the case of $\mathrm{GL}_2$ and $\mathcal E_b=\mathcal O(1)\oplus \mathcal O$, one has
\[
\tilde{\mathcal M}_b^\circ=\BC(\mathcal O(-1)[1])\setminus \{0\}=\Spa k\laurentseries{t^{1/p^\infty}}/\underline{SL_1(D)}
\]
by Example~\ref{ex:BCO-1} and Example~\ref{ex:GL2O1O}. Thus, in this case one can compute
\[
i^\ast Rj_\ast A = R\Gamma(\Spa k\laurentseries{t^{1/p^\infty}}/\underline{SL_1(D)},A)
\]
which is a very explicit formula. If one would use the presentation
\[
\BC(\mathcal O(-1)[1])\times \Spa C=(\mathbb A^1_{C^\sharp})^\diamond/\underline{E}
\]
instead, it would be considerably more difficult to compute the answer. In particular, we critically used quasicompacity of the absolute $\tilde{\mathcal M}_b^\circ$, its base change $\tilde{\mathcal M}_b^\circ\times_k \Spa C$ is no longer quasicompact. {\it This highlights the importance of working with absolute objects, and of using the right local charts.}

In fact, Theorem~\ref{thm:chartbsmooth}, smooth base change, and this formula for $i^\ast Rj_\ast$ show that the gluing of the representation-theoretic strata $D_\et(\Bun_G^b,\Lambda)\cong D(G_b(E),\Lambda)$ into $D_\et(\Bun_G,\Lambda)$ is encoded in the spaces $\tilde{\mathcal M}_b^\circ$, showing that the local charts $\mathcal M_b$ are of fundamental and not just technical importance.
\end{remark}

Now we can prove Theorem~\ref{thm:compactgeneration}.

\begin{proof}[Proof of Theorem~\ref{thm:compactgeneration}] As for closed immersions $i$, the functor $i^\ast$ preserves compact objects, it is enough to handle the case that $U$ is an open substack. Let $b\in |U|\subset \Bun_G\cong B(G)$ be any point of $U$, and let $K\subset G_b(E)$ be an open pro-$p$-subgroup, giving rise to the map
\[
f_K: \widetilde{\mathcal M}_b/\underline{K}\to \Bun_G.
\]
By Corollary~\ref{cor:quotientoftildeMbcompact} and Theorem~\ref{thm:chartbsmooth}, we see that $A^b_K:=Rf_{K!} Rf_K^!\Lambda\in D_\et(\Bun_G,\Lambda)$ is compact; in fact,
\[\begin{aligned}
R\Hom(A^b_K,B)&\cong R\Hom(Rf_K^! \Lambda,Rf_K^! B)\cong R\Hom(Rf_K^! \Lambda,f_K^\ast B\dotimes_\Lambda Rf_K^! \Lambda)\\
&\cong R\Gamma(\widetilde{\mathcal M}_b/\underline{K},f_K^\ast B)\cong (i^{b\ast} B)^K.
\end{aligned}\]
From this computation, we see that the collection of objects $A^b_K$ for varying $b\in |U|$ and $K\subset G_b(E)$ open pro-$p$ form a class of compact generators: Indeed, if $B$ is nonzero, then $(i^{b\ast} B)^K$ must be nonzero for some $b$ and $K$.

To prove the characterization of compact objects, we argue by induction on the number of points of $|U|$, noting that any compact object must be concentrated on a quasicompact substack, and thus on finitely many points. So assume that $|U|$ is finite, $b\in |U|$ is a closed point and $j: V=U\setminus \{b\}\subset U$ is the open complement, so we know the result for $V$. It suffices to prove that $j^\ast$ preserves compact objects. Indeed, then $A\in D_\et(U,\Lambda)$ is compact only if $j^\ast A$ and $i^{b\ast} A$ are compact, implying by induction that all stalks are compact. For the converse, one has to see that $!$-extension from strata preserves compact objects. By induction, this is true for all strata except $\Bun_G^b$; and to check it for this one, reduce to the sheaf corresponding to the representation $c\text-\mathrm{Ind}_K^{G_b(E)} \Lambda$ of $G_b(E)$, where $K$ is pro-$p$; in that case, the $!$-extension is the cone of $j_! j^\ast A_K^b\to A_K^b$, which is compact. Indeed, as $f_K$ is cohomologically smooth, $f_{K!} f_K^!$ commutes with any base change and hence the restriction of $A_K^b$ to $\Bun_G^b$ is given by $g_! g^! \Lambda$ for the map
\[
g: [\ast/\underline{K}]\xrightarrow{g_1} [\ast/\underline{G_b(E)}]\xrightarrow{g_2} \Bun_G^b.
\]
Here, $g_{2!}g_2^!$ is the identity by Proposition~\ref{prop:DetBunstrata}, while $g_{1!} g_1^! \Lambda = g_{1!} \Lambda$ produces the sheaf corresponding to the representation $c\text-\mathrm{Ind}_K^{G_b(E)}\Lambda$. Thus, the restriction of $A_K^b$ to $\Bun_G^b$ gives the representation $c\text-\mathrm{Ind}_K^{G_b(E)}\Lambda$, as desired.

To see that $j^\ast$ preserves compact objects, we can check on the given generators. For generators $A^{b'}_K$ corresponding to $b'\neq b$ we get $j^\ast A^{b'}_K = A^{b'}_K$, so there is nothing to prove. On the other hand, $j^\ast A^b_K = Rf_{K!}^\circ f_K^{\circ!}\Lambda$ for
\[
f_K^\circ: \widetilde{\mathcal M}_b^\circ/\underline{K}\to V\subset \Bun_G.
\]
The compactness of $j^\ast A^b_K$ then follows from $R\Gamma(\widetilde{\mathcal M}_b^\circ/\underline{K},-)$ commuting with all direct sums. But this is true as $\widetilde{\mathcal M}_b^\circ$ is a spatial diamond of finite $\dimtrg$, by Proposition~\ref{prop:basicMb} and Proposition~\ref{prop:MbtoGbnice}, and taking cohomology under $K$ is exact.
\end{proof}

\section{Bernstein--Zelevinsky duality}

We note that one can define a Bernstein--Zelevinsky involution on (the compact objects of) $D_\et(\Bun_G,\Lambda)$.\footnote{The ``classical'' Bernstein--Zelevinsky duality on smooth representations was first defined on the level of irreducibles of $\GL_n(E)$ by Zelevinsky \cite{ZelevinskyDuality}, then generalized by Aubert \cite{AubertDuality} to all groups but still on the level of Grothendieck groups, independently discovered by Bernstein and Schneider--Stuhler \cite{SchneiderStuhler}. Its categorical formulation as (derived) Hom into the Hecke algebra is discussed in \cite{FarguesDuality}.} More precisely, we have the following result. Here, in anticipation of some functor introduced later, we write
\[
\pi_\natural: D_\et(\Bun_G,\Lambda)\to D(\ast,\Lambda)=D(\Lambda): A\mapsto R\pi_!(A\dotimes_\Lambda R\pi^!\Lambda)
\]
for the left adjoint of $\pi^\ast$, where $\pi:\Bun_G\to \ast$ is the (cohomologically smooth!) projection.\footnote{The dualizing sheaf $R\pi^!\Lambda$ is, in fact, isomorphic to $\Lambda[0]$, by fixing a Haar measure on $G_b(E)$ for each basic $b\in B(G)$, so $\pi_\natural$ is isomorphic to $R\pi_!$. Here, we use that local systems concentrated in degree $0$ on $\mathrm{Bun}_G$ embed fully faithfully into local systems on the semistable locus, which follows easily from purity and the smoothness of all Harder--Narasimhan strata.}

\begin{theorem}\label{thm:bernsteinzelevinsky} For any compact object $A\in D_\et(\Bun_G,\Lambda)$, there is a unique compact object $\mathbb D_{\mathrm{BZ}}(A)\in D_\et(\Bun_G,\Lambda)$ with a functorial identification
\[
R\Hom(\mathbb D_{\mathrm{BZ}}(A),B)\cong \pi_\natural(A\dotimes_\Lambda B)
\]
for $B\in D_\et(\Bun_G,\Lambda)$. Moreover, the functor $\mathbb D_{\mathrm{BZ}}$ is a contravariant autoequivalence of $D_\et(\Bun_G,\Lambda)^\omega$, and $\mathbb D_{\mathrm{BZ}}^2$ is naturally isomorphic to the identity.

If $U\subset \Bun_G$ is an open substack and $A$ is concentrated on $U$, then so is $\mathbb D_{\mathrm{BZ}}(A)$. In particular, $\mathbb D_{\mathrm{BZ}}$ restricts to an autoequivalence of the compact objects in $D_\et(\Bun_G^b,\Lambda)\cong D(G_b(E),\Lambda)$ for $b\in B(G)$ basic, and in that setting it is the usual Bernstein--Zelevinsky involution.
\end{theorem}

\begin{proof} By the Yoneda lemma, the uniqueness of $\mathbb D_{\mathrm{BZ}}(A)$ is clear. For simplicity, choose Haar measures on $G_b(E)$ for all basic $b\in B(G)$, leading to an isomorphism $R\pi^!\Lambda\cong \Lambda$ (at generic points, but then by spreading out everywhere, noting that both are invertible), and hence $\pi_\natural\cong R\pi_!$.

For the existence, it suffices to check on a system of generators. For any $b\in B(G)$ and $K\subset G_b(E)$ pro-$p$, consider the map
\[
g_K: [\ast/\underline{K}]\to \Bun_G
\]
factoring over $[\ast/\underline{G_b(E)}]\to \Bun_G^b\subset \Bun_G$, and consider $A=Rg_{K!} \Lambda$. These are, up to shift and twist, the sheaves supported on $\Bun_G^b$ corresponding to the representation $c\text-\mathrm{Ind}_K^{G_b(E)} \Lambda$, and are compact generators of $D_\et(\Bun_G,\Lambda)$. Then the functor
\[
B\mapsto \pi_\natural(A\dotimes_\Lambda B)\cong R\pi_!(A\dotimes_\Lambda B) = R\pi_!(Rg_{K!} \Lambda\dotimes_\Lambda B) = R(\pi\circ g_K)_! g_K^\ast B
\]
is given by $B\mapsto (i^{b\ast} B)^K$ when $i^{b\ast} B$ is regarded as $G_b(E)$-representation. Here, we note that for pro-$p$-groups $K$, the lower $!$-functor along $[\ast/\underline{K}]\to \ast$ maps isomorphically to the lower $\ast$-functor, which is cohomology. By Corollary~\ref{cor:quotientoftildeMbcompact}, this agrees with $R\Hom(A_K^b,B)$, so $\mathbb D_{\mathrm{BZ}}(A) = A_K^b$. This also shows that if $A$ is concentrated on an open substack $U\subset \Bun_G$, then so is $\mathbb D_{\mathrm{BZ}}(A)$.

Now note that
\[
R\Hom(\mathbb D_{\mathrm{BZ}}(A),B)\cong \pi_\natural(A\dotimes_\Lambda B)\cong \pi_\natural(B\dotimes_\Lambda A)\cong R\Hom(\mathbb D_{\mathrm{BZ}}(B),A).
\]
In particular, taking $B=\mathbb D_{\mathrm{BZ}}(A)$, we see that there is a natural functorial map $\mathbb D_{\mathrm{BZ}}^2(A)\to A$. We claim that this is an equivalence. It suffices to check on generators. We have seen that the Bernstein--Zelevinsky dual of $A=Rg_{K!} \Lambda$ is $A^b_K$. Its restriction to $\Bun_G^b$ corresponds again to the representation $i^b_! c\text-\mathrm{Ind}_K^{G_b(E)} \Lambda$, so one easily checks that the map $\mathbb D_{\mathrm{BZ}}^2(A)\to A$ is an isomorphism over $\Bun_G^b$. To see that it is an isomorphism everywhere, one needs to see that if $B=Rj_\ast B'$, $B'\in D_\et(U,\Lambda)$ for some open substack $j: U\subset \Bun_G$ not containing $\Bun_G^b$, then
\[
\pi_\natural(A^b_K\dotimes_\Lambda B)=0.
\]
Twisting a few things away and using the definition of $A^b_K=Rf_{K!} Rf_K^! \Lambda$, this follows from the assertion that for all $A'\in D_\et(\widetilde{\mathcal M}_b^\circ/\underline{K},\Lambda)$, with $j_K: \widetilde{\mathcal M}_b^\circ/\underline{K}\hookrightarrow \widetilde{\mathcal M}_b/\underline{K}$ the open immersion, one has
\[
R\Gamma_c(\widetilde{\mathcal M}_b/\underline{K},Rj_{K\ast} A')=0.
\]
Using the trace map for $\widetilde{\mathcal M}_b\to \widetilde{\mathcal M}_b/\underline{K}$, this follows from Theorem~\ref{thm:partialcompactsupportvanishing}, applied as before with $X=\widetilde{\mathcal M}_b^\circ$ and $S=\Spa k\laurentseries{t}$, noting that base change along $S\to \ast$ holds by smooth base change and is conservative.

The comparison to Bernstein--Zelevinsky duality follows formally by taking $B$ corresponding to the regular representation of $G_b(E)$, in which case $\pi_\natural(A\dotimes_\Lambda B)$ is isomorphic to the underlying chain complex of $A$. Moreover, as the regular representation has two commuting $G_b(E)$-actions, there is a residual $G_b(E)$-action, which is the usual action on $A$. This gives the usual definition of the Bernstein--Zelevinsky involution as $R\Hom$ into the regular representation.
\end{proof}

\section{Verdier duality}

It turns out that one can also understand how Verdier duality acts on $D_\et(\Bun_G,\Lambda)$. The key result is the following.

\begin{theorem}\label{thm:localbidual} Let $j: V\hookrightarrow U$ be an open immersion of open substacks of $\Bun_G$. For any $A\in D_\et(V,\Lambda)$, the natural map
\[
j_! R\sHom_\Lambda(A,\Lambda)\to R\sHom_\Lambda(Rj_\ast A,\Lambda)
\]
is an isomorphism in $D_\et(U,\Lambda)$.
\end{theorem}

Note that one always has
\[
Rj_\ast R\sHom_\Lambda(A,\Lambda) = R\sHom_\Lambda(j_! A,\Lambda);
\]
the theorem asserts that this is also true with $j_!$ and $Rj_\ast$ exchanged, which is related to a local biduality statement: If $A\in D_\et(V,\Lambda)$ is reflexive, the theorem implies formally that $j_! A\in D_\et(U,\Lambda)$ is reflexive.

\begin{proof} We can assume that $U$ and $V$ are quasicompact, and then by induction, we can assume that $V=U\setminus \{b\}$ for some closed $b\in |U|$. The map is clearly an isomorphism over $V$, so it suffices to see that for the compact objects
\[
A^b_K=Rf_{K!}f_K^!\Lambda,\ f_K: \widetilde{\mathcal M}_b/\underline{K}\to U\subset \Bun_G,
\]
one gets an isomorphism after applying $R\Hom(A^b_K,-)$. As $R\Hom(A^b_K,B) = (i^{b\ast} B)^K$, we see that the left-hand side
\[
R\Hom(A^b_K,j_! R\sHom_\Lambda(A,\Lambda))=0
\]
vanishes. On the other hand, using the left adjoint $\pi_\natural: D_\et(\Bun_G,\Lambda)\to D(\ast,\Lambda)\cong D(\Lambda)$ of pullback, we have
\[
R\Hom(A^b_K,R\sHom_\Lambda(Rj_\ast A,\Lambda))\cong R\Hom(A^b_K\dotimes_\Lambda Rj_\ast A,\Lambda)\cong R\Hom(\pi_\natural(A^b_K\dotimes_\Lambda Rj_\ast A),\Lambda).
\]
But by Theorem~\ref{thm:bernsteinzelevinsky}, and the identification of the Bernstein--Zelevinsky dual of $A^b_K$ as $i^b_! c\text-\mathrm{Ind}_K^{G_b(E)}\Lambda$, one has
\[
\pi_\natural(A^b_K\dotimes_\Lambda Rj_\ast A)\cong R\Hom(i^b_! c\text-\mathrm{Ind}_K^{G_b(E)}\Lambda,Rj_\ast A)\cong R\Hom(c\text-\mathrm{Ind}_K^{G_b(E)}\Lambda,Ri^{b!} Rj_\ast A)=0.\qedhere
\]
\end{proof}

Recall that as a cohomologically smooth Artin stack of dimension $0$, $\Bun_G$ admits a dualizing complex $D_{\Bun_G}\in D_\et(\Bun_G,\Lambda)$ that is locally isomorphic to $\Lambda[0]$.

\begin{theorem}\label{thm:charreflexive} For any open substack $U\subset \Bun_G$, an object $A\in D_\et(U,\Lambda)$ is reflexive, i.e.~the natural map
\[
A\to R\sHom_\Lambda(R\sHom_\Lambda(A,D_U),D_U)
\]
is an equivalence, if and only if for all $b\in B(G)$ lying in $U$ with corresponding locally closed stratum $i^b: \Bun_G^b\to U$, the restriction
\[
i^{b\ast} A\in D_\et(\Bun_G^b,\Lambda) = D_\et([\ast/\tilde G_b],\Lambda)\cong D_\et([\ast/\underline{G_b(E)}],\Lambda) = D(G_b(E),\Lambda)
\]
is reflexive as a complex of admissible $G_b(E)$-representations; this means that the complex of $K$-invariants is reflexive in $D(\Lambda)$ for all open pro-$p$-subgroups $K\subset G_b(E)$.
\end{theorem}

In the definition of reflexivity, we can replace $D_U$ by $\Lambda$ (as this changes the dual by a twist, and then the bidual stays the same). The theorem follows immediately from the following result.

\begin{lemma} Let $U\subset \Bun_G$ be an open substack and $A\in D_\et(U,\Lambda)$. For any $b\in B(G)$ lying in $U$, there is a natural isomorphism
\[
i^{b\ast} R\sHom_\Lambda(R\sHom_\Lambda(A,\Lambda),\Lambda)\cong R\sHom_\Lambda(R\sHom_\Lambda(i^{b\ast} A,\Lambda),\Lambda).
\]
\end{lemma}

\begin{proof} We may assume that $U\subset \Bun_G$ is the set of generalizations of $b$, and let $j: V=U\setminus \{b\}\hookrightarrow U$. Let $B=j^\ast A$. Using the exact triangle
\[
j_! B\to A\to i^b_\ast i^{b\ast} A\to,
\]
and the invertibility of $i^{b!} \Lambda$ (as $\Bun_G^b$ is also cohomologically smooth), it is enough to prove that
\[
i^{b\ast} R\sHom_\Lambda(R\sHom_\Lambda(j_! B,\Lambda),\Lambda) = 0.
\]
But $R\sHom_\Lambda(j_! B,\Lambda) = Rj_\ast R\sHom_\Lambda(B,\Lambda)$, and by Theorem~\ref{thm:localbidual}, 
\[
R\sHom_\Lambda(Rj_\ast R\sHom_\Lambda(B,\Lambda),\Lambda) = j_! R\sHom_\Lambda(R\sHom_\Lambda(B,\Lambda),\Lambda).\qedhere
\]
\end{proof}

\section{ULA sheaves}

Finally, we want to classify the objects $A\in D_\et(\Bun_G,\Lambda)$ that are universally locally acyclic with respect to $\Bun_G\to \ast$. Our goal is to prove the following theorem.  {\it This gives a geometric interpretation of the classical notion of admissible representation in terms of $D_{\et} (\Bun_G,\Lambda)$.}

\begin{theorem}\label{thm:ULAbunG} Let $A\in D_\et(\Bun_G,\Lambda)$. Then $A$ is universally locally acyclic with respect to $\Bun_G\to \ast$ if and only if for all $b\in B(G)$, the pullback $i^{b\ast} A$ to $i^b: \Bun_G^b\hookrightarrow \Bun_G$ corresponds under $D_\et(\Bun_G^b,\Lambda)\cong D(G_b(E),\Lambda)$ to a complex $M_b$ of smooth $G_b(E)$-representations for which $M^K_b$ is a perfect complex of $\Lambda$-modules for all open pro-$p$ subgroups $K\subset G_b(E)$.
\end{theorem}

We want to use Proposition~\ref{prop:ULAdualizablestack}. As preparation, we need to understand $D_\et(\Bun_G\times \Bun_G,\Lambda)$. More generally, we have the following result.

\begin{proposition}\label{prop:BunGproduct} Let $G_1$ and $G_2$ be two reductive groups over $E$, and let $G=G_1\times G_2$. Then $\Bun_G\cong \Bun_{G_1}\times \Bun_{G_2}$, giving rise to an exterior tensor product
\[
-\boxtimes-: D_\et(\Bun_{G_1},\Lambda)\times D_\et(\Bun_{G_2},\Lambda)\to D_\et(\Bun_G,\Lambda).
\]

For all compact objects $A_i\in D_\et(\Bun_{G_i},\Lambda)$, $i=1,2$, the exterior tensor product $A_1\boxtimes A_2\in D_\et(\Bun_G,\Lambda)$ is compact, these objects form a class of compact generators, and for all further objects $B_i\in D_\et(\Bun_{G_i},\Lambda)$, $i=1,2$, the natural map
\[
R\Hom(A_1,B_1)\dotimes_\Lambda R\Hom(A_2,B_2)\to R\Hom(A_1\boxtimes A_2,B_1\boxtimes B_2)
\]
is an isomorphism.
\end{proposition}

\begin{remark} The proposition says that as $\Lambda$-linear presentable stable $\infty$-categories, the exterior tensor product functor
\[
\mathcal D_\et(\Bun_{G_1},\Lambda)\otimes_{\mathcal D(\Lambda)} \mathcal D_\et(\Bun_{G_2},\Lambda)\to \mathcal D_\et(\Bun_G,\Lambda)
\]
is an equivalence. Here, we use Lurie's tensor product \cite[Section 4.8]{LurieHA}.
\end{remark}

\begin{proof} We use the compact generators $A_i=A^{b_i}_{K_i}$ for certain $b_i\in B(G_i)$, $K_i\subset G_{i,b_i}(E)$ open pro-$p$. These give rise to $b=(b_1,b_2)\in B(G)$ and $K=K_1\times K_2\subset G_b(E)=G_{1,b_1}(E)\times G_{2,b_2}(E)$, and using
\[
\mathcal M_b\cong \mathcal M_{b_1}\times \mathcal M_{b_2}
\]
and the K\"unneth formula, one concludes that $A_1\boxtimes A_2\cong A^b_K$, which is again compact. As $B(G)=B(G_1)\times B(G_2)$ and open pro-$p$ subgroups of the form $K_1\times K_2\subset G_b(E)$ are cofinal, we see that these objects form a set of compact generators.

Moreover, as $R\Hom(A^b_K,B)=(i^{b\ast} B)^K$ for all $B\in D_\et(\Bun_G,\Lambda)$ and similarly for $A^{b_i}_{K_i}$, we also see that the map
\[
R\Hom(A^{b_1}_{K_1},B_1)\dotimes_\Lambda R\Hom(A^{b_2}_{K_2},B_2)\to R\Hom(A^b_K,B_1\boxtimes B_2)
\]
is an isomorphism. As these objects generate, the same follows for all compact $A_1,A_2$.
\end{proof}

Now we can prove Theorem~\ref{thm:ULAbunG}.

\begin{proof}[Proof of Theorem~\ref{thm:ULAbunG}] By Proposition~\ref{prop:ULAdualizablestack}, we see that $A$ being universally locally acyclic is equivalent to the map
\[
p_1^\ast R\sHom(A,\Lambda)\dotimes_\Lambda p_2^\ast A\to R\sHom(p_1^\ast A,p_2^\ast A)
\]
in $D_\et(\Bun_G\times \Bun_G,\Lambda)\cong D_\et(\Bun_{G\times G},\Lambda)$ being an isomorphism.

By Proposition~\ref{prop:BunGproduct}, this is equivalent to being an isomorphism after applying $R\Hom(A_1\boxtimes A_2,-)$ for varying compact $A_i\in D_\et(\Bun_{G_i},\Lambda)$. Using Proposition~\ref{prop:BunGproduct}, the left-hand side is given by
\[
R\Hom(A_1,R\sHom(A,\Lambda))\dotimes_\Lambda R\Hom(A_2,A)\cong R\Hom(\pi_\natural(A_1\dotimes_\Lambda A),\Lambda)\dotimes_\Lambda R\Hom(A_2,A).
\]
The right-hand side is, using $\pi: \Bun_G\to \ast$ for the projection,
\[\begin{aligned}
R\Hom(A_1\boxtimes A_2,R\sHom(p_1^\ast A,p_2^\ast A))&\cong R\Hom((A_1\dotimes_\Lambda A)\boxtimes A_2,p_2^\ast A)\\
&\cong R\Hom(p_1^\ast(A_1\dotimes_\Lambda A),p_2^\ast R\sHom_\Lambda(A_2,A))\\
&\cong R\Hom(A_1\dotimes_\Lambda A,Rp_{1\ast} p_2^\ast R\sHom_\Lambda(A_2,A))\\
&\cong R\Hom(A_1\dotimes_\Lambda A,\pi^\ast R\Hom(A_2,A))\\
&\cong R\Hom(\pi_\natural(A_1\dotimes_\Lambda A),R\Hom(A_2,A)),
\end{aligned}\]
using usual adjunctions and smooth base change for $p_2$ and $\pi$ several times. Under these isomorphisms, the map
\[
R\Hom(\pi_\natural(A_1\dotimes_\Lambda A),\Lambda)\dotimes_\Lambda R\Hom(A_2,A)\to R\Hom(\pi_\natural(A_1\dotimes_\Lambda A),R\Hom(A_2,A))
\]
is the natural map. This is an isomorphism as soon as $\pi_\natural(A_1\dotimes_\Lambda A)\in D(\Lambda)$ is perfect for all compact $A_1$. In fact, the converse is also true: If one takes $A_2=\mathbb D_{\mathrm{BZ}}(A_1)$, then $R\Hom(A_2,A)=\pi_\natural(A_1\dotimes_\Lambda A)$, and hence it follows that for $M=\pi_\natural(A_1\dotimes_\Lambda A)\in D(\Lambda)$, the map
\[
R\Hom(M,\Lambda)\dotimes_\Lambda M\to R\Hom(M,M)
\]
is an isomorphism, which means that $M$ is dualizable in $D(\Lambda)$, i.e.~perfect.

Now we use the system of compact generators given by $i^b_! c\text-\mathrm{Ind}_K^{G_b(E)} \Lambda$ for varying $b\in B(G)$, with locally closed immersion $i^b: \Bun_G^b\to \Bun_G$, and $K\subset G_b(E)$ open pro-$p$. This translates the condition on perfectness of $R\pi_!(A_1\dotimes_\Lambda R\pi^!\Lambda\dotimes_\Lambda A)$ into the desired condition on stalks.
\end{proof}

\chapter{Geometric Satake}

As before, we fix a nonarchimedean local field $E$ with residue field $\Fq$ of characteristic $p$ and a uniformizer $\pi\in E$. We also fix a reductive group $G$ over $E$, and a coefficient ring $\Lambda$ killed by some integer $n$ prime to $p$.

Recall that for any perfectoid space $S$ over $\Fq$, we defined the ``curve'' $\mathcal Y_S$ over $\mathcal O_E$, as well as $Y_S=\mathcal Y_S\setminus V(\pi)$ and the quotient $X_S=Y_S/\phi^\Z$. In this chapter, we are interested in studying modifications of $G$-torsors on these spaces, and perverse sheaves on such. Our discussion will mirror this three-step procedure of the construction of $X_S$: If one has understood the basic theory over $\mathcal Y_S$, the basic results carry over easily to $Y_S$ and then to $X_S$. While as in previous chapters our main focus is on $X_S$, in this chapter we will actually make critical use of $\mathcal Y_S$ in order to degenerate to the Witt vector affine Grassmannian, and hence to apply some results from the setting of schemes (notably the decomposition theorem). As the discussion here is very much of a local sort, one can usually reduce easily to the case that $G$ is split, and hence admits a (split) reductive model over $\mathcal O_E$, and we will often fix such a split model of $G$.

For any $d\geq 0$, we consider the moduli space $\Div^d_{\mathcal Y}$ parametrizing degree $d$ Cartier divisors $D\subset \mathcal Y_S$. For affinoid $S$, one can form the completion $B^+$ of $\mathcal O_{X_S}$ along $D$. Inverting $D$ defines a localization $B$ of $B^+$. One can then define a positive loop group $L^+_{\Div^d_{\mathcal Y}} G$ and loop group $L_{\Div^d_{\mathcal Y}} G$, with values given by $G(B^+)$ resp.~$G(B)$; for brevity, we will simply write $L^+G$ and $LG$ here. One can then define the local Hecke stack
\[
\Hloc_{G,\Div^d_{\mathcal Y}} = [L^+ G\backslash LG / L^+ G]\to \Div^d_{\mathcal Y}
\]
We will often break symmetry, and first take the quotient on the right to define the Beilinson--Drinfeld Grassmannian
\[
\Gr_{G,\Div^d_{\mathcal Y}} = LG / L^+ G\to \Div^d_{\mathcal Y}
\]
so that
\[
\Hloc_{G,\Div^d_{\mathcal Y}} = L^+ G\backslash \Gr_{G,\Div^d_{\mathcal Y}}.
\]

The Beilinson--Drinfeld Grassmannian $\Gr_{G,\Div^d_{\mathcal Y}}\to \Div^d_{\mathcal Y}$ is a small v-sheaf that can be written as an increasing union of closed subsheaves that are proper and representable in spatial diamonds, by bounding the relative position; this is one main result of \cite{Berkeley}. On the other hand, $L^+ G$ can be written as an inverse limit of truncated positive loop groups, which are representable in locally spatial diamonds and cohomologically smooth; moreover, on each bounded subset, it acts through such a finite-dimensional quotient. This essentially reduces the study of all bounded subsets of $\Hloc_{G,\Div^d_{\mathcal Y}}$ to Artin stacks.

For any small v-stack $S\to \Div^d_{\mathcal Y}$, we let
\[
\Hloc_{G,S/\Div^d_{\mathcal Y}}=\Hloc_{G,\Div^d_{\mathcal Y}}\times_{\Div^d_{\mathcal Y}} S
\]
be the pullback. Let
\[
D_\et(\Hloc_{G,S/\Div^d_{\mathcal Y}},\Lambda)^{\mathrm{bd}}\subset D_\et(\Hloc_{G,S/\Div^d_{\mathcal Y}},\Lambda)
\]
be the full subcategory of all objects with quasicompact support over $\Div^d_{\mathcal Y}$. This is a monoidal category under convolution $\star$. Here, we use the convolution diagram
\[
\Hloc_{G,S/\Div^d_{\mathcal Y}}\times_S \Hloc_{G,S/\Div^d_{\mathcal Y}}\xleftarrow{(p_1,p_2)} L^+ G\backslash L G\times^{L^+ G} L G / L^+ G\xrightarrow{m} \Hloc_{G,S/\Div^d_{\mathcal Y}}
\]
and define
\[
A\star B = Rm_\ast(p_1^\ast A\dotimes_\Lambda p_2^\ast B).
\]

On $D_\et(\Hloc_{G,S/\Div^d_{\mathcal Y}},\Lambda)^{\mathrm{bd}}$, one can define a relative perverse $t$-structure (where an object is perverse if and only if it is perverse over any geometric fibre of $S$), see Section~\ref{sec:perverse}. In particular, this $t$-structure is compatible with any base change in $S$. For this $t$-structure, the convolution $\star$ is left $t$-exact (and $t$-exactness only fails for issues related to non-flatness over $\Lambda$). To prove this, we reinterpret convolution as fusion, and use some results on hyperbolic localization.

Moreover, one can restrict to the complexes $A\in D_\et(\Hloc_{G,S/\Div^d_{\mathcal Y}},\Lambda)^{\mathrm{bd}}$ that are universally locally acyclic over $S$. This condition is also preserved under convolution. For $d=1$, or in general when $S$ maps to the locus of distinct untilts $(\Div^d_{\mathcal Y})_{\neq}\subset \Div^d_{\mathcal Y}$, one can describe the category of universally locally acyclic by the condition that the restriction to any Schubert cell is locally constant with perfect fibres. To prove that all such sheaves are universally locally acyclic, we also introduce (for $d=1$) the affine flag variety, in Section~\ref{sec:affineflag}, and use their Demazure resolutions.

\begin{definition} The Satake category
\[
\Sat(\Hloc_{G,S/\Div^d_{\mathcal Y}},\Lambda)\subset D_\et(\Hloc_{G,S/\Div^d_{\mathcal Y}},\Lambda)^{\mathrm{bd}}
\]
is the category of all $A\in D_\et(\Hloc_{G,S/\Div^d_{\mathcal Y}},\Lambda)^{\mathrm{bd}}$ that are perverse, flat over $\Lambda$ (i.e., for all $\Lambda$-modules $M$, also $A\dotimes_\Lambda M$ is perverse), and universally locally acyclic over $(\Div^1)^I$.
\end{definition}

Intuitively, $\Sat(\Hloc_{G,S/\Div^d_{\mathcal Y}},\Lambda)$ are the ``flat families of perverse sheaves on $\Hloc_{G,S/\Div^d_{\mathcal Y}}\to S$'', where flatness refers both to the geometric aspect of flatness over $S$ (encoded in universal local acyclicity) and the algebraic aspect of flatness in the coefficients $\Lambda$. The Satake category $\Sat(\Hloc_{G,S/\Div^d_{\mathcal Y}},\Lambda)$ is a monoidal category under convolution. The forgetful functor
\[
\Sat(\Hloc_{G,S/\Div^d_{\mathcal Y}},\Lambda)\to D_\et(\Gr_{G,S/\Div^d_{\mathcal Y}},\Lambda)
\]
is fully faithful. If $d=1$ and $S=\Spd \overline{\mathbb F}_q$, then one can compare it to the category considered by Zhu \cite{ZhuWitt} and Yu \cite{YuWittSatake}, defined in terms of the Witt vector affine Grassmannian. Moreover, the categories for $S=\Spd \mathcal O_C$ and $S=\Spd C$ are naturally equivalent to the category for $S=\Spd \overline{\mathbb F}_q$, via the base change functors; here $C=\hat{\overline{E}}$. Thus, the Satake category is, after picking a reductive model of $G$, naturally the same for the Witt vector affine Grassmannian and the $B_{\mathrm{dR}}^+$-affine Grassmannian. At this point, we could in principle use Zhu's results \cite{ZhuWitt} (refined integrally by Yu \cite{YuWittSatake}) to identify the Satake category with the category of representations of $\hat{G}$, at least when $G$ is unramified. However, for the applications we actually need finer knowledge of the functoriality of the Satake equivalence including the case for $d>1$; we thus prove everything we need directly.

More precisely, we now switch to $(\Div^1_X)^d$ in place of $\Div^d_{\mathcal Y}$, replacing also the use of $\mathcal Y_S$ by $X_S$; the two local Hecke stacks are locally isomorphic, so this poses no problems. For any finite set $I$, let
\[
\Hloc_G^I=\Hloc_{G,(\Div^1_X)^I}
\]
and consider the monoidal category
\[
\Sat_G^I(\Lambda) = \Sat(\Hloc_G^I,\Lambda).
\]
In fact, the monoidal structure naturally upgrades to a symmetric monoidal structure. This relies on the fusion product, for which it is critical to allow general finite sets $I$. Namely, given finite sets $I_1,\ldots,I_k$, letting $I=I_1\sqcup\ldots\sqcup I_k$, one has an isomorphism
\[
\Hloc_G^I\times_{(\Div^1)^I} (\Div^1)^{I;I_1,\ldots,I_k}\cong \prod_{j=1}^k \Hloc_G^{I_j}\times_{(\Div^1)^I} (\Div^1)^{I;I_1,\ldots,I_k}
\]
where $(\Div^1)^{I;I_1,\ldots,I_k}\subset (\Div^1)^I$ is the open subset where $x_i\neq x_{i'}$ whenever $i,i'\in I$ lie in different $I_j$'s. The exterior tensor product then defines a functor
\[
\boxtimes_{j=1}^k: \prod_{j=1}^k \Sat_G^{I_j}(\Lambda)\to \Sat_G^{I;I_1,\ldots,I_k}(\Lambda)
\]
where $\Sat_G^{I;I_1,\ldots,I_k}(\Lambda)$ is the variant of $\Sat_G^I(\Lambda)$ for $\Hloc_G^I\times_{(\Div^1)^I} (\Div^1)^{I;I_1,\ldots,I_k}$. However, the restriction functor
\[
\Sat_G^I(\Lambda)\to \Sat_G^{I;I_1,\ldots,I_k}(\Lambda)
\]
is fully faithful, and the essential image of the exterior product lands in its essential image. Thus, we get a natural functor
\[
\bigast_{j=1}^k: \prod_{j=1}^k \Sat_G^{I_j}(\Lambda)\to \Sat_G^I(\Lambda),
\]
independent of the ordering of the $I_j$. In particular, for any $I$, we get a functor
\[
\Sat_G^I(\Lambda)\times \Sat_G^I(\Lambda)\to \Sat_G^{I\sqcup I}(\Lambda)\to \Sat_G^I(\Lambda),
\]
using functoriality of $\Sat_G^J(\Lambda)$ in $J$, which defines a symmetric monoidal structure $\ast$ on $\Sat_G^I(\Lambda)$, commuting with $\star$. This is called the fusion product. In general, for any symmetric monoidal category $(\mathcal C,\ast)$ with a commuting monoidal structure $\star$, the monoidal structure $\star$ necessarily agrees with $\ast$; thus, the fusion product refines the convolution product. (As usual in geometric Satake, we actually need to change $\ast$ slightly by introducing certain signs into the commutativity constraint, depending on the parity of the support of the perverse sheaves.)

Moreover, restricting $A\in \Sat_G^I(\Lambda)$ to $\Gr_G^I$ and taking the pushforward to $(\Div^1)^I$, all cohomology sheaves are local systems of $\Lambda$-modules on $(\Div^1)^I$. By a version of Drinfeld's lemma, these are equivalent to representations of $W_E^I$ on $\Lambda$-modules. This defines a symmetric monoidal fibre functor
\[
F^I: \Sat_G^I(\Lambda)\to \Rep_{W_E^I}(\Lambda),
\]
where $\Rep_{W_E^I}(\Lambda)$ is the category of continuous representations of $W_E^I$ on finite projective $\Lambda$-modules. Using a version of Tannaka duality, one can then build a Hopf algebra in the $\Ind$-category of $\Rep_{W_E^I}(\Lambda)$ so that $\Sat_G^I(\Lambda)$ is given by its category of representations (internal in $\Rep_{W_E^I}(\Lambda)$). For any finite set $I$, this is given by the tensor product of $I$ copies of the corresponding Hopf algebra for $I=\{\ast\}$, which in turn is given by some affine group scheme $\widecheck{G}$ over $\Lambda$ with $W_E$-action.

\begin{theorem}[Theorem~\ref{thm:geometricsatake}] There is a canonical isomorphism $\widecheck{G}\cong \hat{G}$ with the Langlands dual group, under which the action of $W_E$ on $\widecheck{G}$ agrees with the usual action of $W_E$ on $\hat{G}$ up to an explicit cyclotomic twist. If $\sqrt{q}\in \Lambda$, the cyclotomic twist can be trivialized, and $\Sat_G^I(\Lambda)$ is naturally equivalent to the category of $(\hat{G}\rtimes W_E)^I$-representations on finite projective $\Lambda$-modules.
\end{theorem}

For the proof, one can restrict to $\Lambda=\mathbb Z/\ell^n\mathbb Z$; passing to a limit over $n$, one can actually build a group scheme over $\mathbb Z_\ell$. Its generic fibre is reductive, as the Satake category with $\mathbb Q_\ell$-coefficients is (geometrically) semisimple: For this, we use the degeneration to the Witt vector affine Grassmannian and the decomposition theorem for schemes. To identify the reductive group, we argue first for tori, and then for rank $1$ groups, where everything reduces to $G=\mathrm{PGL}_2$ which is easy to analyze by using the minuscule Schubert cell. Here, the pinning includes a cyclotomic twist as of course the cohomology of the minuscule Schubert variety $\mathbb P^1$ of $\Gr_{\mathrm{\PGL}_2}$ contains a cyclotomic twist. Afterwards, we apply hyperbolic localization in order to construct symmetric monoidal functors $\Sat_G\to \Sat_M$ for any Levi $M$ of $G$, inducing dually maps $\widecheck{M}\to \widecheck{G}$. This produces many Levi subgroups of $\widecheck{G}_{\mathbb Q_\ell}$ from which it is easy to get the isomorphism with $\hat{G}_{\mathbb Q_\ell}$, including a pinning. As these maps $\widecheck{M}\to \widecheck{G}$ are even defined integrally, and $\hat{G}(\mathbb Z_\ell)\subset \hat{G}(\mathbb Q_\ell)$ is a maximal compact open subgroup by Bruhat--Tits theory, generated by the rank $1$ Levi subgroups, one can then deduce that $\widecheck{G}\cong \hat{G}$ integrally, again with an explicit (cyclotomic) pinning.

We will also need the following addendum regarding a natural involution. Namely, the local Hecke stack $\Hloc_G^I$ has a natural involution $\mathrm{sw}$ given by reversing the roles of the two $G$-torsors; in the presentation in terms of $LG$, this is induced by the inversion on $LG$. Then $\mathrm{sw}^\ast$ induces naturally an involution of $\Sat_G^I(\Lambda)$, and thus involution can be upgraded to a symmetric monoidal functor commuting with the fibre functor, thus realizing a $W_E$-equivariant automorphism of $\check{G}\cong \hat{G}$.

\begin{proposition}[Proposition~\ref{prop:chevalleyinvolution}] The action of $\mathrm{sw}^\ast$ on $\mathrm{Sat}_G^I$ induces the automorphism of $\hat{G}$ that is the Chevalley involution of the split group $\hat{G}$, conjugated by $\hat{\rho}(-1)$.
\end{proposition}

\section{The Beilinson--Drinfeld Grassmannian}

First, we define the base space of the Beilinson--Drinfeld Grassmannian for any $d\geq 0$.

\begin{definition}\label{def:divdY} For any $d\geq 0$, consider the small v-sheaves on $\Perf_{\Fq}$ given by
\[
\Div^d_{\mathcal Y} = (\Spd \mathcal O_E)^d/\Sigma_d,\ \Div^d_Y = (\Spd E)^d/\Sigma_d,\ \Div^d_X = \Div^d = (\Spd E/\phi^\Z)^d/\Sigma_d,
\]
where $\Sigma_d$ is the symmetric group.
\end{definition}

As always, quotients of small v-sheaves are taken inside v-sheaves, and are still small (and in particular exist).

\begin{proposition}\label{prop:divdY} For any $d\geq 0$, there is a functorial injection
\begin{altenumerate}
\item[{\rm (i)}] from $\Div^d_{\mathcal Y}(S)$ into the set of closed Cartier divisors on $\mathcal Y_S$,
\item[{\rm (ii)}] from $\Div^d_Y(S)$ into the set of closed Cartier divisors on $Y_S$,
\item[{\rm (iii)}] from $\Div^d_X(S)$ into the set of closed Cartier divisors on $X_S$.
\end{altenumerate}
Moreover, in case (i) and (ii), if $S=\Spa(R,R^+)$ is affinoid perfectoid, then for any closed Cartier divisor $D\subset \mathcal Y_S$ resp.~$D\subset Y_S$ in the image of this embedding, the adic space $D=\Spa(Q,Q^+)$ is affinoid. In case (iii), the same happens locally in the analytic topology on $S$.
\end{proposition}

\begin{proof} We handle case (i) first. Over $(\Spd \mathcal O_E)^d$ and for $S=\Spa(R,R^+)$ affinoid, we get $d$ untilts $R_i^\sharp$, $i=1,\ldots,d$ of $R$, and there are elements $\xi_i\in W_{\mathcal O_E}(R^+)$ generating the kernel of $\theta_i: W_{\mathcal O_E}(R^+)\to R_i^{\sharp +}$. Each of the $\xi_i$ defines a closed Cartier divisor by Proposition~\ref{prop:degreeonecartier}. Then $\xi=\prod_i \xi_i$ defines another closed Cartier divisor, given by $\Spa(A,A^+)$ for $A=W_{\mathcal O_E}(R^+)[\tfrac 1{[\varpi]}]/\xi$, and $A^+$ the integral closure of $W_{\mathcal O_E}(R^+)/\xi$, where $\varpi\in R$ is a pseudouniformizer.

Now the ideal sheaf of this closed Cartier divisor is a line bundle, and by \cite[Proposition 19.5.3]{Berkeley}, line bundles on $\mathcal Y_S$ satisfy v-descent. Thus, even if we are only given a map $S\to \Div^d_{\mathcal Y}=(\Spd \mathcal O_E)^d/\Sigma_d$, we can still define a line bundle $\mathcal I\subset \mathcal O_{\mathcal Y_S}$, and it still defines a closed Cartier divisor as this can be checked v-locally. Also, $V(\mathcal I)\subset \mathcal Y_S$ is quasicompact over $S$, as this can again be checked v-locally. This implies that it is contained in some affinoid $\mathcal Y_{S,[0,n]}$, and hence $D=\Spa(A,A^+)$ is affinoid in general.

The case (ii) now follows formally by passing to an open subset, and case (iii) by passing to the quotient under Frobenius. 
\end{proof}

\begin{remark}
As in \cite{FarguesClassFieldTheory} one checks that $\Div^d_{\mathcal{Y}}(S)$ is the set of ``relative Cartier divisors'' of degree $d$, that is to say Cartier divisors that give degree $d$ Cartier divisors when pulled back via any geometric point $\Spa (C,C^+)\to S$. The same holds for $\Div^d_Y$ and $\Div^d_X$.
\end{remark}

In the following we will consider a perfectoid space $S$ equipped with a map $f: S\to \Div^d_{\mathcal Y}$ (resp.~$f: S\to \Div^d_Y$, resp.~$f: S\to \Div^d_X$). We denote by $D_S\subset \mathcal Y_S$ (resp.~$D_S\subset Y_S$, $D_S\subset X_S$) the corresponding closed Cartier divisor. Let $\mathcal I_S\subset \mathcal O_{\mathcal Y_S}$ (resp.~$\mathcal I_S\subset \mathcal O_{Y_S}$, $\mathcal I_S\subset \mathcal O_{X_S}$) be the corresponding invertible ideal sheaf.

Let us note the following descent result for vector bundles.

\begin{proposition}\label{prop:vdescentCartier} Sending $S$ as above to the category of vector bundles on $D_S$ defines a v-stack.
\end{proposition}

\begin{proof} Any vector bundle on $D_S$ defines a v-sheaf on $\Perf_{S}$: This reduces formally to the structure sheaf of $D_S$, which then further reduces to the structure of $\mathcal O_{\mathcal Y_S}$ (resp.~$\mathcal O_{Y_S}$, resp.~$\mathcal O_{X_S}$). It remains to prove that v-descent of vector bundles is effective. The case of $X_S$ reduces to $Y_S$ as locally on $S$, the relevant $D_S$ is isomorphic; and clearly $Y_S$ reduces to $\mathcal Y_S$.

Now assume first that $S=\Spa(C,C^+)$ for some complete algebraically closed $C$. Then $D_S$ is given by a finite sum of degree $1$ Cartier divisors on $\mathcal Y_S$, and one can reduce by induction to the case of degree $1$ Cartier divisors, where the result is \cite[Lemma 17.1.8]{Berkeley} applied to the corresponding untilt of $S$.

On the other hand, assume that $T\to S$ is an \'etale cover with a vector bundle $\mathcal E_T$ on $D_T$ equipped with a descent datum to $D_S$; we want to prove descent to $D_S$. By the argument of de Jong--van der Put \cite[Proposition 3.2.2]{deJongvanderPut}, cf.~\cite[Proposition 8.2.20]{KedlayaLiu1}, one can reduce to the case that $T\to S$ is a finite \'etale cover. Then $D_T\to D_S$ is also finite \'etale (as $\mathcal Y_T\to \mathcal Y_S$ is), and the result follows from usual finite \'etale descent.

Now let $S$ be general and $T\to S$ a v-cover with a vector bundle $\mathcal E_T$ on $D_T$ equipped with a descent datum to $D_S$. For any geometric point $\Spa(C,C^+)\to S$, one can descent $\mathcal E_{T\times_S \Spa(C,C^+)}$ to a vector bundle $\mathcal E_{\Spa(C,C^+)}$ on $D_{\Spa(C,C^+)}$. Now we follow the proof of \cite[Lemma 17.1.8]{Berkeley} to see that one can in fact descend $\mathcal E_T$ in an \'etale neighborhood of $\Spa(C,C^+)$, which is enough by the previous paragraph. We can assume that $S$ and $T$ are affinoid. As $\mathcal E_{\Spa(C,C^+)}$ is necessarily free, also $\mathcal E_{T\times_S \Spa(C,C^+)}$ is free, and by an \'etale localization we can assume that $\mathcal E_T$ is free. Then the descent datum is given by some matrix with coefficients in $\mathcal O_{D_{T\times_S T}}$. Moreover, by approximating the basis coming via pullback from $\mathcal E_{\Spa(C,C^+)}$, we can ensure that this matrix has coefficients in $1+[\varpi]\mathcal O_{\mathcal Y_{T\times_S T}}^+(\mathcal Y_{T\times_S T,[0,n]})/\xi$ for some $n$ so that $D_S\subset \mathcal Y_{S,[0,n]}$; here $\xi$ is a generator of $\mathcal I_S$. Now one uses that the v-cohomology group $H^1_v(S,\mathcal O(\mathcal Y_{[0,n]})^+)$ is almost zero, as follows from almost vanishing in the perfectoid case, and writing it as a direct summand of the positive structure sheaf of the base change to $\mathcal O_E[\pi^{1/p^\infty}]^\wedge$. Then the argument from \cite[Lemma 17.1.8]{Berkeley} applies, showing that one can successively improve the basis to produce a basis invariant under the descent datum in the limit.
\end{proof}

Assuming that $D_S$ is affinoid, as is the case locally on $S$, we let
\[
B^+_{\Div^d_{\mathcal Y}}(S)\ (\mathrm{resp.}\ B^+_{\Div^d_Y}(S),\ \mathrm{resp.}\ B^+_{\Div^d_X}(S))
\]
be (the global sections of) the completion of $\mathcal O_{\mathcal Y_S}$ along $\mathcal I_S$ (resp.~of $\mathcal O_{Y_S}$ along $\mathcal I_S$, resp.~of $\mathcal O_{X_S}$ along $\mathcal I_S$), and
\[
B_{\Div^d_{(-)}}(S) = B^+_{\Div^d_{(-)}}(S)[\tfrac 1{\mathcal I_S}].
\]
This defines v-sheaves $B^+_{\Div^d_{(-)}}\subset B_{\Div^d_{(-)}}$ over $\Div^d_{(-)}$ in all three cases. In the case of $d=1$, those rings are the ones that are usually denoted $B^+_{\mathrm{dR}}$, resp.~$B_{\mathrm{dR}}$.

\begin{definition}\label{def:loopfunctors} Let $Z$ be an affine scheme over $\mathcal O_E$. The positive loop space $L^+_{\Div^d_{\mathcal Y}}Z$ (resp.~loop space $L_{\Div^d_{\mathcal Y}}Z$) of $Z$ is the v-sheaf over $\mathrm{Div}^d_{\mathcal Y}$ given by
\[
S\mapsto L^+_{\Div^d_{\mathcal Y}} Z(S) = Z(B^+_{\Div^d_{\mathcal Y}}(S))\ (\mathrm{resp.}\ S\mapsto L_{\Div^d_{\mathcal Y}} Z(S) = Z(B_{\Div^d_{\mathcal Y}}(S))).
\]
Similarly, if $Z$ is an affine scheme over $E$, one defines the positive loop space $L^+_{\Div^d_Y} Z$ and $L^+_{\Div^d_X} Z$ (resp.~loop space $L_{\Div^d_Y} Z$ and $L_{\Div^d_X} Z$).
\end{definition}

We note that we use affinity of $Z$ to see that these are actually v-sheaves --- this makes it possible to reduce to the v-sheaf property of the structure sheaf. (It is likely that they define v-sheaves for general schemes $Z$, using Bhatt's Tannaka result \cite{BhattTannaka} together with descent results for perfect complexes, but we do not pursue this here.) Now we can define the local Hecke stacks.

\begin{definition}\label{def:Hkloc} Let $G$ be a reductive group over $\mathcal O_E$ (resp.~over $E$, resp.~over $E$). The local Hecke stack $\Hloc_{G,\Div^d_{\mathcal Y}}$ (resp.~$\Hloc_{G,\Div^d_Y}$, resp.~$\Hloc_{G,\Div^d_X}$) is the functor sending an affinoid perfectoid $S\to \Div^d_{\mathcal Y}$ (resp.~$S\to \Div^d_Y$, resp.~$S\to \Div^d_X$, assuming that $D_S$ is affinoid) to the groupoid of pairs of $G$-bundles $\mathcal E_1,\mathcal E_2$ over $B^+_{\Div^d_{\mathcal Y}}(S)$ (resp.~over $B^+_{\Div^d_Y}(S)$, resp.~over $B^+_{\Div^d_X}(S)$) together with an isomorphism $\mathcal E_1\cong \mathcal E_2$ over $B_{\Div^d_{\mathcal Y}}(S)$ (resp.~over $B_{\Div^d_Y}(S)$, resp.~over $B_{\Div^d_X}(S)$).
\end{definition}

The $G$-bundles here are taken in the algebraic sense, as living on the spectrum of the respective rings. As in Chapter III, we will generally think of $G$-bundles in Tannakian terms, as tensor functors from $\mathrm{Rep}_E G$ to vector bundles; and over an affine or affinoid base, vector bundles will always correspond to finite projective modules, so the notion of $G$-bundle is rather insensitive to the underlying geometric formalism.

Also note that in the case of $\Div^d_X$, the local Hecke functor is only defined on a certain full subcategory of affinoid perfectoid $S\to \Div^d_X$, namely those where $D_S$ is affinoid; but any $S\to \Div^d_X$ admits an open cover by such by Proposition~\ref{prop:divdY}, so we have still defined the functor on a basis.

\begin{proposition}\label{prop:Hkloc} The local Hecke stack $\Hloc_{G,\Div^d_{\mathcal Y}}$ (resp.~$\Hloc_{G,\Div^d_Y}$, resp.~$\Hloc_{G,\Div^d_X}$) is a small v-stack. There is a natural isomorphism of \'etale stacks over $\Div^d_{\mathcal Y}$ (resp.~over $\Div^d_Y$, resp.~over $\Div^d_X$)
\[
\Hloc_{G,\Div^d_{\mathcal Y}}\cong (L^+_{\Div^d_{\mathcal Y}} G)\backslash (L_{\Div^d_{\mathcal Y}} G) / (L^+_{\Div^d_{\mathcal Y}} G)
\]
(resp.
\[
\Hloc_{G,\Div^d_Y}\cong (L^+_{\Div^d_Y} G)\backslash (L_{\Div^d_Y} G) / (L^+_{\Div^d_Y} G),
\]
resp.
\[
\Hloc_{G,\Div^d_X}\cong (L^+_{\Div^d_X} G)\backslash (L_{\Div^d_X} G) / (L^+_{\Div^d_X} G).)
\]
\end{proposition}

\begin{proof} The category of vector bundles over $B^+_{\Div^d_{\mathcal Y}}$ (resp.~$B^+_{\Div^d_Y}$, $B^+_{\Div^d_X}$) satisfies v-descent: It is enough to check this modulo powers of the ideal $\mathcal I_S$, where the result is Proposition~\ref{prop:vdescentCartier}. By the Tannakian formalism, it follows that the category of $G$-bundles also satisfies v-descent, so one can descend $\mathcal E_1,\mathcal E_2$. The isomorphism between $\mathcal E_1$ and $\mathcal E_2$ over $B_{\Div^d_{\mathcal Y}}(S)$ (resp.~over $B_{\Div^d_Y}(S)$, resp.~over $B_{\Div^d_X}(S)$) is then given by a section of an affine scheme over the respective ring, which again satisfies v-descent. Smallness follows from the argument of Proposition~\ref{prop:BunGsmall}.

Any $G$-bundle over $B^+_{\Div^d_{\mathcal Y}}(S)$ is \'etale locally on $S$ trivial. Indeed, if $S$ is a geometric point then $B^+_{\Div^d_{\mathcal Y}}(S)$ is a product of complete discrete valuation rings with algebraically closed residue field, so that all $G$-torsors are trivial. In general, note that triviality of the $G$-torsor over $B^+_{\Div^d_{\mathcal Y}}(S)$ is implied by triviality modulo $\mathcal I_S$ (as one can always lift sections over nilpotent thickenings). Then the result follows from \cite[Proposition 5.4.21]{GabberRamero}. Trivializing $\mathcal E_1$ and $\mathcal E_2$ \'etale locally then directly produces the given presentations.
\end{proof}

Similarly, one can define the Beilinson--Drinfeld Grassmannians.

\begin{definition}\label{def:BDaffgrass} Let $G$ be a reductive group over $\mathcal O_E$ (resp.~over $E$, resp.~over $E$). The Beilinson--Drinfeld Grassmannian $\Gr_{G,\Div^d_{\mathcal Y}}$ (resp.~$\Gr_{G,\Div^d_Y}$, resp.~$\Gr_{G,\Div^d_X}$) is the functor sending an affinoid perfectoid $S\to \Div^d_{\mathcal Y}$ (resp.~$S\to \Div^d_Y$, resp.~$S\to \Div^d_X$, assuming again that $D_S$ is affinoid) to the groupoid of $G$-bundles $\mathcal E$ over $B^+_{\Div^d_{\mathcal Y}}(S)$ (resp.~over $B^+_{\Div^d_Y}(S)$, resp.~over $B^+_{\Div^d_X}(S)$) together with a trivialization of $\mathcal E$ over $B_{\Div^d_{\mathcal Y}}(S)$ (resp.~over $B_{\Div^d_Y}(S)$, resp.~over $B_{\Div^d_X}(S)$).
\end{definition}

\begin{proposition}\label{prop:BDaffgrass} The Beilinson--Drinfeld Grassmannian $\Gr_{G,\Div^d_{\mathcal Y}}$ (resp.~$\Gr_{G,\Div^d_Y}$, resp.~$\Gr_{G,\Div^d_X}$) is a small v-sheaf. There is a natural isomorphism of \'etale sheaves over $\Div^d_{\mathcal Y}$ (resp.~over $\Div^d_Y$, resp.~over $\Div^d_X$)
\[
\Gr_{G,\Div^d_{\mathcal Y}}\cong (L_{\Div^d_{\mathcal Y}} G) / (L^+_{\Div^d_{\mathcal Y}} G)
\]
(resp.
\[
\Gr_{G,\Div^d_Y}\cong (L_{\Div^d_Y} G) / (L^+_{\Div^d_Y} G),
\]
resp.
\[
\Gr_{G,\Div^d_X}\cong (L_{\Div^d_X} G) / (L^+_{\Div^d_X} G).)
\]
\end{proposition}

\begin{proof} The proof is identical to the proof of Proposition~\ref{prop:Hkloc}.
\end{proof}

The positive loop group $L^+_{\Div^d_{\mathcal Y}} G$ admits the natural filtration by closed subgroups
\[
(L^+_{\Div^d_{\mathcal Y}} G)^{\geq m}\subset L^+_{\Div^d_{\mathcal Y}} G
\]
defined for all $m\geq 1$ by the kernel of
\[
G(B^+_{\Div^d_{\mathcal Y}})\to G(B^+_{\Div^d_{\mathcal Y}}/\mathcal I_S^m);
\]
we refer to these as the principal congruence subgroups of $L^+_{\Div^d_{\mathcal Y}} G$. Similar definitions of course apply also over $\Div^d_Y$ and $\Div^d_X$. For $d=1$, one can easily describe the graded pieces of this filtration (and again the result also holds for $\Div^1_Y$ and $\Div^1_X$).

\begin{proposition}\label{prop:gradedcongruencesubgroups} There are natural isomorphisms
\[
L^+_{\Div^1_{\mathcal Y}} G/(L^+_{\Div^1_{\mathcal Y}} G)^{\geq 1}\cong G^\diamond
\]
and
\[
(L^+_{\Div^1_{\mathcal Y}} G)^{\geq m}/(L^+_{\Div^1_{\mathcal Y}} G)^{\geq m+1}\cong (\Lie G)^\diamond\{m\}
\]
where $\{m\}$ signifies a ``Breuil-Kisin twist'' by $\mathcal I_S^m/\mathcal I_S^{m+1}$.
\end{proposition}

Here $G^\diamond$ is defined as in \cite[Section 27, before Proposition 27.5]{ECoD}, and sends $S=\Spa(R,R^+)$ to a choice $R^\sharp$ of untilt of $R$, and an element of $G(R^\sharp)$.

\begin{proof} The first equality follows directly from the definitions, while the second comes from the exponential.
\end{proof}

For general $d$, we still have the following result.

\begin{proposition}\label{prop:gradedcongruencesubgroups2} For any $d$ and $m\geq 1$, the quotient
\[
(L^+_{\Div^d_{\mathcal Y}} G)^{\geq m}/(L^+_{\Div^d_{\mathcal Y}} G)^{\geq m+1}
\]
sends a perfectoid space $S\to \Div^d_{\mathcal Y}$ with corresponding Cartier divisor $D_S\subset Y_S$ with ideal sheaf $\mathcal I_S$ to
\[
(\Lie G\otimes_{\mathcal O_E} \mathcal I_S^m/\mathcal I_S^{m+1})(S)
\]
where $\mathcal I_S^m/\mathcal I_S^{m+1}$ is a line bundle on $D_S$. This is representable in locally spatial diamonds, partially proper, and cohomologically smooth of $\ell$-dimension equal to $d$ times the dimension of $G$.

Moreover, one can filter
\[
(L^+_{\Div^d_{\mathcal Y}} G)^{\geq m}/(L^+_{\Div^d_{\mathcal Y}} G)^{\geq m+1}\times_{\Div^d_{\mathcal Y}} (\Div^1_{\mathcal Y})^d
\]
with subquotients given by twists of
\[
(\Lie G)^\diamond\times_{\Div^1_{\mathcal Y},\pi_i} (\Div^1_{\mathcal Y})^d
\]
where $\pi_i: (\Div^1_{\mathcal Y})^d\to \Div^1_{\mathcal Y}$ is the projection to the $i$-th factor.
\end{proposition}

\begin{proof} The description of the subquotient follows from the exponential sequence again. To check that it is representable in locally spatial diamonds, partially proper, and cohomologically smooth, we can pull back to $(\Div^1_{\mathcal Y})^d$, and then these properties follow from the existence of the given filtration. For this in turn, note that over $(\Div^1_{\mathcal Y})^d$, we have $d$ ideal sheaves $\mathcal I_1,\ldots,\mathcal I_d$, and one can filter $\mathcal O_{D_S}$ by $\mathcal O_{D_S}/\mathcal I_1$, $\mathcal I_1/\mathcal I_1\mathcal I_2$, $\ldots$, $\mathcal I_1\cdots\mathcal I_{d-1}/\mathcal I_1\cdots \mathcal I_d$, each of which is, after pullback to an affinoid perfectoid space $S$, isomorphic to $\mathcal O_{S_i^\sharp}$.
\end{proof}

One can also show that the first quotient is cohomologically smooth, but this is slightly more subtle.

\begin{proposition}\label{prop:firstquotientsmooth} For any $d$, the quotient
\[
L^+_{\Div^d_{\mathcal Y}} G/(L^+_{\Div^d_{\mathcal Y}} G)^{\geq 1}\to \Div^d_{\mathcal Y}
\]
parametrizes over a perfectoid space $S\to \Div^d_{\mathcal Y}$ maps $D_S\to G$. For any quasiprojective smooth scheme $Z$ over $\mathcal O_E$, the sheaf
\[
T_Z\to \Div^d_{\mathcal Y}
\]
taking a perfectoid $S$ over $\Div^d_{\mathcal Y}$ to maps $D_S\to Z$ (of locally ringed spaces) is representable in locally spatial diamonds, partially proper, and cohomologically smooth over $\Div^d_{\mathcal Y}$ of $\ell$-dimension equal to $d$ times the dimension of $Z$; in particular, this applies to this quotient group.
\end{proposition}

\begin{proof} The description of the quotient group is clear. To analyze $T_Z$, we first note that if $Z$ is an affine space, then the result holds true, as was proved in the previous proposition. In fact, after pullback via the quasi-pro-\'etale surjective morphism $(\Div^1_{\mathcal{Y}})^d\to \Div^d_{\mathcal{Y}}$, there is a sequence of morphisms
$$
T_{\mathbb{A}^n} =W_1 \lto W_2\lto \cdots \lto W_{d+1}=(\Div^1_{\mathcal{Y}})^d
$$
where, for $S$ affinoid perfectoid with $S\to  W_{i+1}$ giving rise to the untilts $(S_1^\sharp,\dots, S_d^\sharp)\in (\Div^1_{\mathcal{Y}})^d (S)$, 
$W_i\times_{W_{i+1}} S \to S$ is locally on $S$ isomorphic to $\mathbb{A}^{n,\diamond}_{S_i^\sharp}$.

If $Z'\to Z$ is any separated \'etale map between schemes over $\mathcal O_E$, we claim that $T_{Z'}\to T_Z$ is also separated \'etale. For this, we analyze the pullback along any $S\to T_Z$ given by some perfectoid space $S$ and a map $D_S\to Z$ (of locally ringed spaces). Then $D'=D_S\times_Z Z'\to D_S$ is separated \'etale (here, the fibre product is taken as in \cite[Proposition 3.8]{HuberAGeneralization}), and the fibre product $T_{Z'}\times_{T_Z} S$ parametrizes $S'\to S$ with a lift $D_{S'}\to D'$ over $D_S$. By Lemma~\ref{lem:etalecartierdivisor} below, this is representable by a perfectoid space separated \'etale over $S$. In case $Z'\to Z$ is an open immersion, it follows that $T_{Z'}\to T_Z$ is injective and \'etale, thus an open immersion.

Now note that for any geometric point of $T_Z$, the corresponding map $D_S\to Z$ has finite image, and is thus contained in some open affine subscheme. It follows that $T_Z$ admits an open cover by $T_{Z'}$ for affine $Z'$. If $Z$ is affine, then one sees directly that $T_Z$ is representable in locally spatial diamonds and partially proper by taking a closed immersion into $\mathbb A^n_{\mathcal O_E}$ for some $n$. For cohomological smoothness, we observe that we can in fact choose these affines so that they admit \'etale maps to $\mathbb A^d_{\mathcal O_E}$, as again we only need to arrange this at finitely many points at a time. Now the result follows from the discussion of $\mathbb A^d_{\mathcal O_E}$ and of separated \'etale maps.
\end{proof}

\begin{lemma}\label{lem:etalecartierdivisor} Let $S$ be a perfectoid space with a map $S\to \Div^d_{\mathcal Y}$ giving rise to the Cartier divisor $D_S\subset \mathcal Y_S$. Let $D'\to D_S$ be a separated \'etale map. Then there is a separated \'etale map $S'\to S$ such that for $T\to S$, maps $T\to S'$ over $S$ are equivalent to lifts $D_T\to D'$ over $D_S$.
\end{lemma}

\begin{proof} By descent of separated \'etale maps \cite[Proposition 9.7]{ECoD}, we can assume that $S$ is strictly totally disconnected. Exhausting $D'$ by a rising union of quasicompact subspaces, we can assume that $D'$ is quasicompact. In any geometric fibre, $D'$ is then a disjoint union of open subsets (as any geometric fibre is, up to nilpotents, a finite disjoint union of untilts $\Spa(C_i^\sharp,C_i^{\sharp+})$), and this description spreads into a small neighborhood by \cite[Proposition 11.23, Lemma 15.6]{ECoD}. We can thus reduce to the case that $D'\to D_S$ is an open immersion. Now the lemma follows from the observation that the map $|D_S|\to |S|$ is closed. 
\end{proof}

\section{Schubert varieties}

Now we recall the Schubert varieties. Assume in this section that $G$ is a split reductive group over $\mathcal O_E$ (or over $E$, but in that case we can always choose a model over $\mathcal O_E$). Fix a split torus and Borel $T\subset B\subset G$. Note that we can always pass to the situation of split $G$ by making a finite \'etale extension of $\mathcal O_E$ resp.~$E$; this way, the results of this section are useful in the general case. Similarly, the cases of $X_S$ and $Y_S$ reduce easily to the case of $\mathcal Y_S$, so we only do the latter case explicitly here.

Assume first that $d=1$. In that case, for every geometric point $S=\Spa(C,C^+)\to \Div^1_{\mathcal Y}=\Spd \mathcal O_E$ given by an untilt $S^\sharp=\Spa(C^\sharp,C^{\sharp +})$ of $S$, one has $B^+_{\Div^1_{\mathcal Y}}=B^+_{\mathrm{dR}}(C^\sharp)$ and $B_{\Div^1_{\mathcal Y}}=B_{\mathrm{dR}}(C^\sharp)$ for the usual definition of $B_{\mathrm{dR}}^+$ and $B_{\mathrm{dR}}$ (relative to $\mathcal O_E$). Recall that $B^+_{\mathrm{dR}}(C^\sharp)$ is a complete discrete valuation ring with residue field $C^\sharp$, fraction field $B_{\mathrm{dR}}(C^\sharp)$, and uniformizer $\xi$. It follows that by the Cartan decomposition
\[
G(B_{\mathrm{dR}}(C^\sharp)) = \bigsqcup_{\mu\in X_\ast(T)^+} G(B_{\mathrm{dR}}^+(C^\sharp))\mu(\xi)G(B_{\mathrm{dR}}^+(C^\sharp)),
\]
so as a set
\[
\Hloc_{G,\Div^1_{\mathcal Y}}(S)/{\cong}\ = X_\ast(T)^+,
\]
the dominant cocharacters of $T$. Recall that on $X_\ast(T)^+$, we have the dominance order, where $\mu\geq \mu'$ if $\mu-\mu'$ is a sum of positive coroots with $\mathbb Z_{\geq 0}$-coefficients.

\begin{remark} Since we work over $\mathcal{Y}$ and do not restrict ourselves to $Y$, we include the case of the Cartier divisor $\pi=0$. For this divisor, $C^\sharp=C$ and $B^+_{\mathrm{dR}}(C^\sharp)= W_{\O_E}(C)$.
\end{remark}

\begin{definition}\label{def:schubertvarieties} For any $\mu\in X_\ast(T)^+$, let
\[
\Hloc_{G,\Div^1_{\mathcal Y},\leq \mu}\subset \Hloc_{G,\Div^1_{\mathcal Y}}
\]
be the subfunctor of all those maps $S\to \Hloc_{G,\Div^1_{\mathcal Y}}$ such that at all geometric points $S'=\Spa(C,C^+)\to S$, the corresponding $S'$-valued point is given by some $\mu'\in X_\ast(T)^+$ with $\mu'\leq \mu$. Moreover,
\[
\Gr_{G,\Div^1_{\mathcal Y},\leq \mu}\subset \Gr_{G,\Div^1_{\mathcal Y}}
\]
is the preimage of $\Hloc_{G,\Div^1_{\mathcal Y},\leq \mu}\subset \Hloc_{G,\Div^1_{\mathcal Y}}$.
\end{definition}

Recall the following result.

\begin{proposition}[{\cite[Proposition 20.3.6]{Berkeley}}]\label{prop:schubertvarietyspatialdiamond} The inclusion
\[
\Hloc_{G,\Div^1_{\mathcal Y},\leq \mu}\subset \Hloc_{G,\Div^1_{\mathcal Y}}
\]
is a closed subfunctor and
\[
\Hloc_{G,\Div^1_{\mathcal Y}} = \varinjlim_\mu \Hloc_{G,\Div^1_{\mathcal Y},\leq \mu};
\]
thus, similar properties hold for $\Gr_{G,\Div^1_{\mathcal Y}}$. Here, the index category is the partially ordered set of $\mu$'s under the dominance order, which is a disjoint union (over $\pi_1(G)$) of filtered partially ordered sets.

The map $\Gr_{G,\Div^1_{\mathcal Y},\leq\mu}\to \Div^1_{\mathcal Y}$ is proper and representable in spatial diamonds.
\end{proposition}

\begin{proof} It is enough to prove the assertions over $\Gr_{G,\Div^1_{\mathcal Y}}$ as this is a v-cover of $\Hloc_{G,\Div^1_{\mathcal Y}}$. Then \cite[Proposition 20.3.6]{Berkeley} gives the results, except for the assertion that
\[
\Gr_{G,\Div^1_{\mathcal Y}} = \varinjlim_\mu \Gr_{G,\Div^1_{\mathcal Y},\leq \mu}.
\]
For this, note that the map from right to left is clearly an injection. For surjectivity, note that for any quasicompact $S$ with a map $S\to \Gr_{G,\Div^1_{\mathcal Y}}$, only finitely many strata can be met, as the meromorphic isomorphism of $G$-bundles necessarily has bounded poles. This, coupled with the fact $\Gr_{G,\Div^1_{\mathcal Y}}\to \Div^1_{\mathcal Y}$ is separated while $\Gr_{G,\Div^1_{\mathcal Y},\leq\mu}\to \Div^1_{\mathcal Y}$ is proper, implies that the map $\bigsqcup_\mu \Gr_{G,\Div^1_{\mathcal Y},\leq \mu}\to \Gr_{G,\Div^1_{\mathcal Y}}$ is a v-cover, whence we get the desired surjectivity.
\end{proof}

In particular,
\[
\Hloc_{G,\Div^1_{\mathcal Y},\mu} = \Hloc_{G,\Div^1_{\mathcal Y},\leq\mu}\setminus \bigcup_{\mu'<\mu} \Hloc_{G,\Div^1_{\mathcal Y},\leq\mu'}\subset \Hloc_{G,\Div^1_{\mathcal Y},\leq \mu}
\]
is an open subfunctor, and similarly its preimage $\Gr_{G,\Div^1_{\mathcal Y},\mu}\subset \Gr_{G,\Div^1_{\mathcal Y},\leq\mu}$. By the Cartan decomposition, the space $\Hloc_{G,\Div^1_{\mathcal Y},\mu}$ has only one point in every geometric fibre over $\Div^1_{\mathcal Y}$. This point can in fact be defined as a global section
\[
[\mu]: \Div^1_{\mathcal Y}\to \Gr_{G,\Div^1_{\mathcal Y},\mu}
\]
given by $\mu(\xi)\in (L_{\Div^1_{\mathcal Y}} G)(S)$ whenever $\xi$ is a local generator of $\mathcal I_S$; up to the action of $L^+_{\Div^1_{\mathcal Y}} G$, this is independent of the choice of $\xi$.

\begin{proposition}\label{prop:openschubertcell} The map
\[
[\mu]: \Div^1_{\mathcal Y}\to \Hloc_{G,\Div^1_{\mathcal Y},\mu}
\]
given by $\mu$ is a v-cover. This gives an isomorphism
\[
\Hloc_{G,\Div^1_{\mathcal Y},\mu}\cong [\Div^1_{\mathcal Y}/(L^+_{\Div^1_{\mathcal Y}} G)_\mu]
\]
where $(L^+_{\Div^1_{\mathcal Y}} G)_\mu\subset L^+_{\Div^1_{\mathcal Y}} G$ is the closed subgroup stabilizing the section $[\mu]$ of $\Gr_{G,\Div^1_{\mathcal Y}}/\Div^1_{\mathcal Y}$. Recalling the principal congruence subgroups 
\[
(L^+_{\Div^1_{\mathcal Y}} G)^{\geq m}\subset L^+_{\Div^1_{\mathcal Y}} G,
\]
we let
\[
(L^+_{\Div^1_{\mathcal Y}} G)_\mu^{\geq m} = (L^+_{\Div^1_{\mathcal Y}} G)_\mu\cap (L^+_{\Div^1_{\mathcal Y}} G)^{\geq m}\subset (L^+_{\Div^1_{\mathcal Y}} G)_\mu.
\]
Then
\[
(L^+_{\Div^1_{\mathcal Y}} G)_\mu/(L^+_{\Div^1_{\mathcal Y}} G)_\mu^{\geq 1}\cong (P_\mu^-)^\diamond\subset L^+_{\Div^1_{\mathcal Y}} G/(L^+_{\Div^1_{\mathcal Y}} G)^{\geq 1}\cong G^\diamond
\]
and
\[
(L^+_{\Div^1_{\mathcal Y}} G)_\mu^{\geq m}/(L^+_{\Div^1_{\mathcal Y}} G)_\mu^{\geq m+1}\cong (\Lie G)_{\mu\leq m}^\diamond\{m\}\subset (L^+_{\Div^1_{\mathcal Y}} G)^{\geq m}/(L^+_{\Div^1_{\mathcal Y}} G)^{\geq m+1}\cong (\Lie G)^\diamond\{m\},
\]
where $P_\mu^-\subset G$ is the parabolic with Lie algebra $(\Lie G)_{\mu\leq 0}$, and $(\Lie G)_{\mu\leq m}\subset \Lie G$ is the subspace on which $\mu$ acts via weights $\leq m$ via the adjoint action.

In particular,
\[
\Gr_{G,\Div^1_{\mathcal Y},\mu}\cong L^+_{\Div^1_{\mathcal Y}} G/(L^+_{\Div^1_{\mathcal Y}} G)_\mu
\]
is cohomologically smooth of $\ell$-dimension $\langle 2\rho,\mu\rangle$ over $\Div^1_{\mathcal Y}$.
\end{proposition}

\begin{proof} We first handle the case $G=\GL_n$ with its standard upper-triangular Borel and diagonal torus. In that case, $\mu$ is given by some sequence $k_1\geq\ldots\geq k_n$ of integers, and $\Gr_{G,\Div^1_{\mathcal Y},\mu}$ parametrizes $B^+_{\Div^1_{\mathcal Y}}$-lattices
\[
\Xi\subset B_{\Div^1_{\mathcal Y}}^n
\]
that are of relative position $\mu$ at all points. Let $S=\Spa(R,R^+)$ be an affinoid perfectoid space with a map $S\to \Div^1_{\mathcal Y}=\Spd \mathcal O_E$ given by an untilt $S^\sharp=\Spa(R^\sharp,R^{\sharp+})$ over $\mathcal O_E$ of $S$. By the proof of \cite[Proposition 19.4.2]{Berkeley}, the $R^\sharp$-modules
\[
\mathrm{Fil}^i_\Xi (R^\sharp)^n = (\xi^i \Xi\cap B_{\mathrm{dR}}^+(R^\sharp)^n)/(\xi^i \Xi\cap \xi B_{\mathrm{dR}}^+(R^\sharp)^n)
\]
are finite projective of rank equal to the number of occurrences of $-i$ among $k_1,\ldots,k_n$. Localizing, we may assume that they are finite free. We may then pick a basis $e_1,\ldots,e_n$ of $(R^\sharp)^n$ so that any $\mathrm{Fil}^i_\Xi (R^\sharp)^n$ is freely generated by a subset $e_1,\ldots,e_{n_i}$ of $e_1,\ldots,e_n$. Lifting $e_{n_{i-1}+1},\ldots,e_{n_i}$ to elements of $f_{n_{i-1}+1},\ldots,f_{n_i}\in \xi^i\Xi\cap B_{\mathrm{dR}}^+(R^\sharp)^n$, and setting $g_{n_{i-1}+1}=\xi^{-i} f_{n_{i-1}+1}$, $\ldots$, $g_{n_i} = \xi^{-i} f_{n_i}$, or equivalently $g_j = \xi^{k_j} f_j$ for $j=1,\ldots,n$, one sees that $f_1,\ldots,f_n$ form a $B_{\mathrm{dR}}^+(R^\sharp)$-basis of $B_{\mathrm{dR}}^+(R^\sharp)^n$, and $g_1,\ldots,g_n$ will form a $B_{\mathrm{dR}}^+(R^\sharp)$-basis of $\Xi$. Thus, changing basis to the $f_i$'s, one has moved $\Xi$ to the lattice
\[
\xi^{k_1} B_{\mathrm{dR}}^+(R^\sharp)\oplus\ldots\oplus \xi^{k_n} B_{\mathrm{dR}}^+(R^\sharp).
\]
This is the lattice corresponding to $[\mu]\in \Gr_{\GL_n,\Div^1_{\mathcal Y}}$, showing that the map
\[
\Div^1_{\mathcal Y}\to \Hloc_{\GL_n,\Div^1_{\mathcal Y},\mu}
\]
is indeed surjective.

Moreover, the stabilizer $(L^+_{\Div^1_{\mathcal Y}}\GL_n)_\mu$ of $\xi^{k_1} B_{\mathrm{dR}}^+(R^\sharp)\oplus\ldots\oplus \xi^{k_n} B_{\mathrm{dR}}^+(R^\sharp)$ in $L^+_{\Div^1_{\mathcal Y}} \GL_n$ is the set of all matrices $A=(A_{ij})\in \GL_n(B^+_{\Div^1_{\mathcal Y}})$ such that for $i<j$, $A_{ij}\in \xi^{k_i-k_j} B^+_{\Div^1_{\mathcal Y}}$. This easily implies the description of
\[
(L^+_{\Div^1_{\mathcal Y}}\GL_n)_\mu/(L^+_{\Div^1_{\mathcal Y}}\GL_n)_\mu^{\geq 1}\cong (P_\mu^-)^\diamond\subset \GL_n^\diamond
\]
and
\[
(L^+_{\Div^1_{\mathcal Y}}\GL_n)_\mu^{\geq m}/(L^+_{\Div^1_{\mathcal Y}}\GL_n)_\mu^{\geq m+1}\cong (\Lie \GL_n)_{\mu\leq m}^\diamond\{m\}\subset (\Lie \GL_n)^\diamond\{m\}.
\]
The description also implies that $(L^+_{\Div^1_{\mathcal Y}}\GL_n)_\mu$ contains $L^+_{\Div^1_{\mathcal Y}} P_\mu^-$ and $(L^+_{\Div^1_{\mathcal Y}} U_a)^{\geq \mu(a)}$ for any positive root $a$.

In general, picking a closed immersion of $G$ into $\GL_n$ (compatible with the torus and the Borel), one sees that
\[
(L^+_{\Div^1_{\mathcal Y}} G)_\mu/(L^+_{\Div^1_{\mathcal Y}} G)_\mu^{\geq 1}\subset (P_\mu^-)^\diamond\subset G^\diamond
\]
and
\[
(L^+_{\Div^1_{\mathcal Y}} G)_\mu^{\geq m}/(L^+_{\Div^1_{\mathcal Y}} G)_\mu^{\geq m+1}\subset (\Lie G)_{\mu\leq m}^\diamond\{m\}\subset (\Lie G)^\diamond\{m\}
\]
as these subquotients embed into the similar subquotient for $\GL_n$. Moreover, one sees that $(L^+_{\Div^1_{\mathcal Y}} G)_\mu$ contains $L^+_{\Div^1_{\mathcal Y}} P_\mu^-$ and $(L^+_{\Div^1_{\mathcal Y}} U_a)^{\geq \mu(a)}$ for any positive root $a$. These imply that the two displayed inclusions are actually equalities.

A consequence of these considerations is that the map
\[
L^+_{\Div^1_{\mathcal Y}} G/(L^+_{\Div^1_{\mathcal Y}} G)_\mu\to L^+_{\Div^1_{\mathcal Y}}\GL_n/(L^+_{\Div^1_{\mathcal Y}} \GL_n)_\mu
\]
is a closed immersion (as this happens on all subquotients for the principal congruence filtration). The target is isomorphic to $\Gr_{\GL_n,\Div^1_{\mathcal Y},\mu}$, which contains $\Gr_{G,\Div^1_{\mathcal Y},\mu}$ as a closed subspace (by \cite[Proposition 20.3.7]{Berkeley}). We see that we get an inclusion
\[
L^+_{\Div^1_{\mathcal Y}} G/(L^+_{\Div^1_{\mathcal Y}} G)_\mu\hookrightarrow \Gr_{G,\Div^1_{\mathcal Y},\mu}\subset \Gr_{\GL_n,\Div^1_{\mathcal Y},\mu}
\]
of closed subspaces, with the same geometric points: This implies that it is an isomorphism (e.g., as the map is then necessarily a closed immersion, thus qcqs, so one can apply \cite[Lemma 12.5]{ECoD}). From here, all statements follow.
\end{proof}

\begin{remark}
The map
\[
\Gr_{G,\Div^1_{\mathcal Y},\mu} = L^+_{\Div^1_{\mathcal{Y}}} G / (L^+_{\Div^1_{\mathcal{Y}}} G)_\mu \to L^+_{\Div^1_{\mathcal{Y}}} G  / 
(L^+_{\Div^1_{\mathcal{Y}}} G)^{\geq 1}_\mu \cong (G/P_{-\mu})^\diamond
\]
is the Bia{\l}ynicki-Birula map, see \cite{CaraianiScholze}.
\end{remark}

Passing to general $d$, we first note that any geometric fibre of
\[
\Hloc_{G,\Div^d_{\mathcal Y}}\to \Div^d_{\mathcal Y}
\]
is isomorphic to a product of geometric fibres of $\Hloc_{G,\Div^1_{\mathcal Y}}\to \Div^1_{\mathcal Y}$. More precisely, if $f: \Spa(C,C^+)\to \Div^d_{\mathcal Y}$ is a geometric point, it is given by an unordered tuple $\Spa(C_i^\sharp,C_i^{\sharp +})$, $i\in I$ with $|I|=d$, of untilts over $\mathcal O_E$. Some of these may be equal, so one can partition $I$ into sets $I_1,\ldots,I_r$ of equal untilts. Then we really have $r$ untilts, given by maps $f_1,\ldots,f_r: \Spa(C,C^+)\to \Div^1_{\mathcal Y}$, and one has an isomorphism
\[
\Hloc_{G,\Div^d_{\mathcal Y}}\times_{\Div^d_{\mathcal Y},f} \Spa(C,C^+)\cong \prod_{i=1}^r \Hloc_{G,\Div^1_{\mathcal Y}}\times_{\Div^1_{\mathcal Y},f_i} \Spa(C,C^+),
\]
and similarly
\[\begin{aligned}
L^+_{\Div^d_{\mathcal Y}} G\times_{\Div^d_{\mathcal Y},f} \Spa(C,C^+)&\cong \prod_{i=1}^r L^+_{\Div^1_{\mathcal Y}} G\times_{\Div^1_{\mathcal Y},f_i} \Spa(C,C^+),\\
L_{\Div^d_{\mathcal Y}} G\times_{\Div^d_{\mathcal Y},f} \Spa(C,C^+)&\cong \prod_{i=1}^r L_{\Div^1_{\mathcal Y}} G\times_{\Div^1_{\mathcal Y},f_i} \Spa(C,C^+),\\
\Gr_{G,\Div^d_{\mathcal Y}}\times_{\Div^d_{\mathcal Y},f} \Spa(C,C^+)&\cong \prod_{i=1}^r \Gr_{G,\Div^1_{\mathcal Y}}\times_{\Div^1_{\mathcal Y},f_i} \Spa(C,C^+).
\end{aligned}\]
Indeed, it suffices to prove this on the level of the positive loop and loop group, where in turn it follows from a similar decomposition of $B^+_{\Div^d_{\mathcal Y}}$ after pullback, which is clear.

In particular, we can define the following version of Schubert varieties.

\begin{definition}\label{def:generalschubertvariety} For any unordered collection $\mu_\bullet=(\mu_j)_{j\in J}$ of elements $\mu_j\in X_\ast(T)^+$ with $|J|=d$, let
\[
\Hloc_{G,\Div^d_{\mathcal Y},\leq \mu_\bullet}\subset \Hloc_{G,\Div^d_{\mathcal Y}}
\]
be the subfunctor of all those $S\to \Hloc_{G,\Div^d_{\mathcal Y}}$ such that at all geometric points $\Spa(C,C^+)\to S$, then equipped with an (unordered) tuple of $d$ untilts $\Spa(C_i^\sharp,C_i^{\sharp +})$, $i\in I$ with $|I|=d$, there is some bijection between $\psi: I\cong J$ such that the relative position of $\mathcal E_1$ and $\mathcal E_2$ at $\Spa(C_i^\sharp,C_i^{\sharp +})$ is bounded by
\[
\sum_{j\in J,C_{\psi(j)}^\sharp\cong C_i^\sharp} \mu_j.
\]

Let
\[
\Gr_{G,\Div^d_{\mathcal Y},\leq \mu_\bullet}\subset \Gr_{G,\Div^d_{\mathcal Y}}
\]
be the preimage of $\Hloc_{G,\Div^d_{\mathcal Y},\leq \mu_\bullet}\subset \Hloc_{G,\Div^d_{\mathcal Y}}$.
\end{definition}

\begin{proposition}\label{prop:generalschubertvariety} The inclusion
\[
\Hloc_{G,\Div^d_{\mathcal Y},\leq \mu_\bullet}\subset \Hloc_{G,\Div^d_{\mathcal Y}}
\]
is a closed subfunctor. The map $\Gr_{G,\Div^d_{\mathcal Y},\leq \mu_\bullet}\to \Div^d_{\mathcal Y}$ is proper, representable in spatial diamonds, and of finite $\dimtrg$.
\end{proposition}

\begin{proof} This can be checked after pullback to $(\Div^1_{\mathcal Y})^d$. Then it follows from \cite[Proposition 20.5.4]{Berkeley}.
\end{proof}

Moreover, we have the following result. Here, we let
\[
(L^+_{\Div^d_{\mathcal Y}} G)^{<m} = L^+_{\Div^d_{\mathcal Y}} G/(L^+_{\Div^d_{\mathcal Y}} G)^{\geq m}
\]
be the quotient by the principal congruence subgroup.

\begin{proposition}\label{prop:finitedimaction} For any $\mu_\bullet=(\mu_j)_{j\in J}$ as above, the action of $L^+_{\Div^d_{\mathcal Y}} G$ on $\Gr_{G,\Div^d_{\mathcal Y},\leq \mu_\bullet}$ factors over $(L^+_{\Div^d_{\mathcal Y}} G)^{<m}$ where $m$ is chosen so that for $\mu=\sum_{j\in J} \mu_j$, all weights of $\mu$ on $\Lie G$ are $\leq m$.
\end{proposition}

\begin{proof} We need to see that the action of $(L^+_{\Div^d_{\mathcal Y}} G)^{\geq m}$ is trivial. As everything is separated, this can be checked on geometric points, where one reduces to $d=1$ by a decomposition into products. Then it follows from Proposition~\ref{prop:openschubertcell}.
\end{proof}

\section{Semi-infinite orbits}\label{sec:semi infinite orbits}

For this section, we continue to assume that $G$ is split, and again we only spell out the case of $\Div^d_{\mathcal Y}$; analogous results hold for $\Div^d_Y$ and $\Div^d_X$, and follow easily from the case presented. 

Previously, we stratified the affine Grassmanian using the Cartan decomposition, the strata being affine Schubert cells. We now use the Iwasawa decomposition to obtain another stratification by semi-infinite orbits.

Fix a cocharacter $\lambda: \mathbb G_m\to T\subset G$, inducing a Levi $M_\lambda$ with Lie algebra $(\Lie G)_{\lambda=0}$, a parabolic $P_\lambda=P_\lambda^+$ with Lie algebra $(\Lie G)_{\lambda\geq 0}$ and its unipotent radical $U_\lambda$ with Lie algebra $(\Lie G)_{\lambda>0}$. We get an action of the v-sheaf $\mathbb G_m$ (taking an affinoid perfectoid space $S=\Spa(R,R^+)$ of characteristic $p$ to $R^\times$) on $\mathrm{Gr}_{G,\Div^d_{\mathcal Y}}$ via the composition of the Teichm\"uller map
\[
[\cdot]: \mathbb G_m\to L^+_{\Div^d_{\mathcal Y}} \mathbb G_m,
\]
the map
\[
L^+_{\Div^d_{\mathcal Y}}\lambda: L^+_{\Div^d_{\mathcal Y}}\mathbb G_m\to L^+_{\Div^d_{\mathcal Y}} G
\]
and the action of $L^+_{\Div^d_{\mathcal Y}} G$ on $\mathrm{Gr}_{G,\Div^d_{\mathcal Y}}$. We wish to apply Braden's theorem in this setup. For this purpose, we need to verify Hypothesis~\ref{hyp:goodGmaction}. To construct the required stratification, we use the affine Grassmannian
\[
\Gr_{P_\lambda,\Div^d_{\mathcal Y}}
\]
associated to the parabolic $P_\lambda$. Note that this admits a map
\[
\Gr_{P_\lambda,\Div^d_{\mathcal Y}}\to \Gr_{M_\lambda,\Div^d_{\mathcal Y}}\to \Gr_{\overline{M}_\lambda,\Div^d_{\mathcal Y}}
\]
where $M_\lambda$ is the Levi quotient of $P_\lambda$ and $\overline{M}_\lambda$ is the maximal torus quotient of $M_\lambda$ (the cocenter). Then $\Gr_{\overline{M}_\lambda,\Div^d_{\mathcal Y}}$ admits a surjection from a disjoint union of copies of $(\Div^1_{\mathcal Y})^d$ parametrized by $X_\ast(\overline{M}_\lambda)^d$. While there are many identifications between these copies, the sum $\overline{\mu}:=\sum_{i=1}^d \overline{\mu}_i\in X_\ast(\overline{M}_\lambda)$ defines a well-defined locally constant function
\begin{equation}\label{eq:loc cst Gr M lambda}
\Gr_{\overline{M}_\lambda,\Div^d_{\mathcal Y}}\to X_\ast(\overline{M}_\lambda).
\end{equation}
More precisely, for $\Spa (C,C^+)\to \Gr_{\overline{M}_\lambda,\Div^d_{\mathcal Y}}$ a geometric point, let $C_1^\sharp,\dots,C_r^\sharp$ be the corresponding distinct untilts with $1\leq r\leq d$. Then 
$$
\Gr_{\overline{M}_\lambda,\Div^d_{\mathcal Y}}\times_{\Div^d_{\mathcal{Y}}} \Spa (C,C^+) \cong \prod_{i=1}^r 
\mathrm{Gr}_{\overline{M}_\lambda,\Div^1_{\mathcal{Y}}} \times_{\Div^1_{\mathcal{Y}}} \Spa (C,C^+)
$$
where the morphism $\Spa (C,C^+)\to \Div^1_{\mathcal{Y}}$ is given by $C_i^\sharp$ on the $i$-th component of the product. This is identified with
$$
\prod_{i=1}^r X_*(\overline{M}_\lambda) \times\Spa (C,C^+)
$$
and the weighted sum morphism $X_*(\overline{M}_\lambda)^r\to X_*(\overline{M}_\lambda)$ (weighing each term with the multiplicity of $C_i^\sharp$ as an untilt of $C$ in the morphism $\Spa(C,C^+)\to \Div^d_{\mathcal Y}$) defines thus a function
$$
\Gr_{\overline{M}_\lambda,\Div^d_{\mathcal Y}}\times_{\Div^d_{\mathcal{Y}}} \Spa (C,C^+)  \lto X_*(\overline{M}_\lambda).
$$
This defines the locally constant function of \eqref{eq:loc cst Gr M lambda}.

For $\nu\in X_\ast(\overline{M}_\lambda)$ let
\[
\Gr_{P_\lambda,\Div^d_{\mathcal Y}}^{\nu}\subset \Gr_{P_\lambda,\Div^d_{\mathcal Y}}
\]
be the corresponding open and closed subset obtained as the preimage.

\begin{proposition}\label{prop:hypbraden} The map
\[
\Gr_{P_\lambda,\Div^d_{\mathcal Y}}=\bigsqcup_{\nu} \Gr_{P_\lambda,\Div^d_{\mathcal Y}}^{\nu} \to \Gr_{G,\Div^d_{\mathcal Y}}
\]
is bijective on geometric points, and it is a locally closed immersion on each $\Gr_{P_\lambda,\Div^d_{\mathcal Y}}^{\nu}$. The union $\bigcup_{\nu'\leq \nu} \Gr_{P_\lambda,\Div^d_{\mathcal Y}}^{\nu'}$ has closed image in $\Gr_{G,\Div^d_{\mathcal Y}}$. The action of $\mathbb G_m$ via $L^+\lambda$ on $\Gr_{P_\lambda,\Div^d_{\mathcal Y}}$ extends to an action of the monoid $\mathbb A^1$, and the $\mathbb G_m$-fixed points agree with $\Gr_{M_\lambda,\Div^d_{\mathcal Y}}$.
\end{proposition}

Applying this proposition also in the case of the inverse $\mathbb G_m$-action, and pulling back to a relative Schubert variety, verifies Hypothesis~\ref{hyp:goodGmaction} in this situation.

\begin{proof} The action of $\mathbb G_m$ on $P_\lambda$ via conjugation extends to an action of the monoid $\mathbb A^1$. Applying loop spaces to this observation and the observation that the map $L_{\Div^d_{\mathcal Y}}P_\lambda\to \Gr_{P_\lambda,\Div^d_{\mathcal Y}}$ is equivariant for the action of $L^+_{\Div^d_{\mathcal Y}}\mathbb G_m$ on the source via conjugation and on the target via the given action (as we quotient by the right action of $L^+_{\Div^d_{\mathcal Y}}P_\lambda$) gives the action of the monoid $L^+_{\Div^d_{\mathcal Y}}\mathbb A^1$, and thus of $\mathbb A^1$ via restricting to Teichm\"uller elements. As everything is separated, this also shows that $\mathbb G_m$-fixed points necessarily lie in the image of $L_{\Div^d_{\mathcal Y}}M_\lambda$, thus the $\mathbb G_m$-fixed points agree with $\Gr_{M_\lambda,\Div^d_{\mathcal Y}}$.

Bijectivity of the map
\[
\mathrm{Gr}_{P_\lambda,\Div^d_{\mathcal Y}}\to \mathrm{Gr}_{G,\Div^d_{\mathcal Y}}
\]
on geometric points follows from the Iwahori decomposition. It remains to prove that the map is a locally closed immersion on each $\Gr_{P_\lambda,\Div^d_{\mathcal Y}}^{\nu}$, and the union over $\nu'\leq \nu$ is closed. Picking a closed embedding into $\mathrm{GL}_n$, this reduces to the case $G=\mathrm{GL}_n$, and by writing any standard parabolic as an intersection of maximal parabolics, we can assume that $P_\lambda\subset \mathrm{GL}_n$ is a maximal parabolic. Passing to a higher exterior power, we can even assume that $P_\lambda\subset \mathrm{GL}_n$ is the mirabolic, fixing a one-dimensional quotient of the standard representation. In that case, $\mathrm{Gr}_{\mathrm{GL}_n,\Div^d_{\mathcal Y}}$ parametrizes finite projective $B^+_{\Div^d_{\mathcal Y}}$-modules $M$ with an identification $M\otimes_{B^+_{\Div^d_{\mathcal Y}}} B_{\Div^d_{\mathcal Y}}\cong B_{\Div^d_{\mathcal Y}}^n$, and $\mathrm{Gr}_{P_\lambda,\Div^d_{\mathcal Y}}$ parametrizes such $M$ for which the image $L\subset B_{\Div^d_{\mathcal Y}}$ of $M$ in the quotient $B^n_{\Div^d_{\mathcal Y}}\to B_{\Div^d_{\mathcal Y}}: (x_1,\ldots,x_n)\mapsto x_n$ is a line bundle over $B^+_{\Div^d_{\mathcal Y}}$. It also suffices to prove the result after pullback along $(\Div^1_{\mathcal Y})^d\to \Div^d_{\mathcal Y}$. Now the result follows from the next lemma.
\end{proof}

In the following, the ``relative position'' of a $B^+$-lattice $L\subset B$ to the standard lattice $B^+\subset B$ refers to the image under the map
\[
\Gr_{\mathbb G_m,\Div^d_{\mathcal Y}}\to X_\ast(\mathbb G_m)=\mathbb Z
\]
defined above.

\begin{lemma} Let $S=\Spa(R,R^+)$ be an affinoid perfectoid space over $\Fq$ with untilts $S_i^\sharp = \Spa(R_i^\sharp,R_i^{\sharp +})$ over $\mathcal O_E$ for $i=1,\ldots,n$. Let $\xi_i\in W_{\mathcal O_E}(R^+)$ generate the kernel of $\theta_i: W_{\mathcal O_E}(R^+)\to R_i^{\sharp +}$ and let $\xi=\xi_1\cdots\xi_n$. Let $B^+$ be the $\xi$-adic completion of $W_{\mathcal O_E}(R^+)[\tfrac 1{[\varpi]}]$ where $\varpi\in R$ is a pseudouniformizer, and let $B=B^+[\frac 1\xi]$. Finally, let
\[
L\subset B
\]
be a finitely generated $B^+$-module that is open and bounded, i.e.~there is some integer $N$ such that $\xi^N B^+\subset L\subset \xi^{-N} B^+\subset B$.

For any $m\in \mathbb Z$, let $S_m\subset S$ be the subset of those points at which the relative position of $L$ to $B^+\subset B$ is given by $m$. Then $\bigcup_{m'\geq m} S_{m'}$ is closed, and if $S_m=S$ then the $B^+$-module $L$ is a line bundle.
\end{lemma}

\begin{proof} We can assume that $L\subset B^+$ via multiplying by a power of $\xi$. Let $s\in S$ be any point, corresponding to a map $\Spa(K(s),K(s)^+)\to S_m$. Let $B^+_s$ be the version of $B^+$ constructed from $(K(s),K(s)^+)$. Then $B^+_s$ is a finite product of discrete valuation rings, and the image $L_s$ of $L\otimes_{B^+} B^+_s$ in $B_s^+$ is necessarily free of rank $1$. Then $s\in S_m$ if and only if the length of $B_s^+/L_s$ as $B^+_s$-module is given by $m$. Localizing on $S$ if necessary, we can find an element $l\in L\subset B$ whose image in $L_s$ is a generator. In a neighborhood of $s$, the element $l$ generates a submodule $L'=B^+\cdot l\subset L$ whose relative position to $B^+$ is bounded above by $m$ at all points by the next lemma, and then the relative position of $L\subset B$ is also bounded above by $m$ at all points as $L'\subset L$. This gives the desired semicontinuity of the stratification (noting that as $L$ is open and bounded, only finitely many values of $m$ can appear). If $S_m=S$, then the containment $L'\subset L$ has to be an equality, and hence $L=L'$ is generated by $l$, so $L$ is a line bundle.
\end{proof}

\begin{lemma} In the situation of the previous lemma, let $f\in B^+$ be any element, and consider the map
\[
|S|\to \mathbb Z_{\geq 0}\cup \{\infty\}
\]
sending any point $s$ of $S$ to the length of $B^+_s/f$ as $B^+_s$-module. This map is semicontinuous in the sense that for any $m\geq 0$, the locus where it is $\leq m$ is open.
\end{lemma}

\begin{proof} For any $i=1,\ldots,n$, one can look at the closed subspace $S_i\subset S$ where the image of $f$ in $R_i^\sharp$ vanishes. On the open complement of all $S_i$, the function is identically $0$. By induction, we can thus pass to a closed subspace $S_i\subset S$, where we can consider the function $f_i=\frac f{\xi_i}$; the length function for $f$ is then the length function for $f_i$ plus one. This gives the result.
\end{proof}

\begin{example}\label{ex:semi infinite orbit}
Suppose $\lambda\in X_*(T)$ is regular dominant. Then $P_\lambda=B$. We then obtain the stratification by semi-infinite orbits 
$$
S_\nu = \mathrm{Gr}_{B,\Div^d_\mathcal{Y}}^{\nu}
$$
for $\nu \in X_*(T)$. One has $S_\nu\hookrightarrow \mathrm{Gr}_{G,\Div^d_\mathcal{Y}}$,  a locally closed immersion, and 
$$
\mathrm{Gr}_{G,\Div^d_{\mathcal{Y}}} = \bigcup_{\nu\in X_*(T)} \mathrm{Gr}_{B,\Div^d_\mathcal{Y}}^{\nu}
$$
(disjoint union) at the level of points. 
\end{example}

We can now apply Theorem~\ref{thm:braden}. Here, we also use the opposite parabolic $P_\lambda^-\subset G$. If $S\to \Div^d_{\mathcal Y}$ is a small v-stack, we denote
\[
\Gr_{G,S/\Div^d_{\mathcal Y}} = \Gr_{G,\Div^d_{\mathcal Y}} \times_{\Div^d_{\mathcal Y}} S
\]
and similarly
\[
\Hloc_{G,S/\Div^d_{\mathcal Y}} = \Hloc_{G,\Div^d_{\mathcal Y}} \times_{\Div^d_{\mathcal Y}} S.
\]
For any $A\in D_\et(\Gr_{G,S/\Div^d_{\mathcal Y}},\Lambda)$, we call $A$ {\it bounded} if it arises via pushforward from some finite union $\Gr_{G,S/\Div^d_{\mathcal Y},\leq\mu_\bullet}\subset \Gr_{G,S/\Div^d_{\mathcal Y}}$. We let
\[
D_\et(\Gr_{G,S/\Div^d_{\mathcal Y}},\Lambda)^{\mathrm{bd}}\subset D_\et(\Gr_{G,S/\Div^d_{\mathcal Y}},\Lambda)
\]
be the corresponding full subcategory.

\begin{corollary}\label{cor:hyperboliclocalizationGr} Let $S\to \Div^d_{\mathcal Y}$ be any small v-stack. Consider the diagram
\[\xymatrix{
&\Gr_{P_\lambda,\Div^d_{\mathcal Y}}\ar[dl]_{q^+}\ar[dr]^{p^+}&\\
\Gr_{G,\Div^d_{\mathcal Y}}&&\Gr_{M_\lambda,\Div^d_{\mathcal Y}}\\
&\Gr_{P_\lambda^-,\Div^d_{\mathcal Y}}\ar[ul]_{q^-}\ar[ur]^{p^-} &
}\]
and denote by $q^+_S$ etc.~the base change along $S\to \Div^d_{\mathcal Y}$. Consider the full subcategory
\[
D_\et(\Gr_{G,S/\Div^d_{\mathcal Y}},\Lambda)^{\mathbb G_m\text-\mathrm{mon},\mathrm{bd}}\subset D_\et(\Gr_{G,S/\Div^d_{\mathcal Y}},\Lambda)^{\mathrm{bd}}
\]
of all $A\in D_\et(\Gr_{G,S/\Div^d_{\mathcal Y}},\Lambda)$ that are bounded and $\mathbb G_m$-monodromic in the sense of Definition~\ref{def:hyperboliclocalization}.

On $D_\et(\Gr_{G,S/\Div^d_{\mathcal Y}},\Lambda)^{\mathbb G_m\text-\mathrm{mon},\mathrm{bd}}$, the natural map
\[
R(p^-_S)_\ast R(q^-_S)^!\to R(p^+_S)_!(q^+_S)^\ast
\]
is an equivalence, inducing a ``constant term'' functor
\[
\mathrm{CT}_{P_\lambda}: D_\et(\Gr_{G,S/\Div^d_{\mathcal Y}},\Lambda)^{\mathbb G_m\text-\mathrm{mon},\mathrm{bd}}\to D_\et(\Gr_{M_\lambda,S/\Div^d_{\mathcal Y}},\Lambda)^{\mathrm{bd}}.
\]
This functor commutes with any base change in $S$ and preserves the condition of being universally locally acyclic over $S$ (which is well-defined for bounded $A$).
\end{corollary}

\begin{proof} This follows from Proposition~\ref{prop:hypbraden} and Theorem~\ref{thm:braden}, Proposition~\ref{prop:hyperboliclocalizationbasechange} and Proposition~\ref{prop:hyperboliclocalizationULA}.
\end{proof}

\begin{example}[Follow-up to Example \ref{ex:semi infinite orbit}]\label{ex:suite semi infinite orbit ex}
In the context of Example \ref{ex:semi infinite orbit}, suppose $d=1$. Then, $\mathrm{Gr}_{T,S/\Div^1_{\mathcal Y}} = X_*(T)\times S$, and the corresponding semi-infinite orbits are denoted by $S_\nu\subset \Gr_{G,\Div^1_{\mathcal Y}}$ for $\nu\in X_*(T)$. Thus, 
$$
\mathrm{CT}_B(A)  = \bigoplus_{\nu\in X_*(T)} R(p_\nu)_! (A|_{S_\nu})
$$ 
with $p_{\nu}: S_\nu\times_{\Div^1_\mathcal{Y}} S\to S\subset \Gr_{T,S/\Div^1_{\mathcal Y}}$ the embedding indexed by $\nu$.
\end{example}

As a final topic here, let us analyze more closely the semi-infinite orbits in the special fibre, i.e.~for the Witt vector affine Grassmannian $\Gr_G^{\mathrm{Witt}}$ (so that $(\Gr_G^{\mathrm{Witt}})^\diamond\cong \Gr_{G,\Spd \Fq/\Div^1_{\mathcal Y}}$), which is an increasing union of perfections of projective varieties over $\Fq$ by \cite{BhattScholzeWitt}, cf.~also \cite{ZhuWitt}. For any $\lambda\in X_\ast(T)$ as above, we have the semi-infinite orbit
\[
S_\lambda = LU\cdot [\lambda]\subset \Gr_G^{\mathrm{Witt}}.
\]

\begin{proposition}\label{prop:semiinfiniteorbitaffine} For any $\mu\in X_\ast(T)^+$, the intersection $S_\lambda\cap \Gr_{G,\leq \mu}^{\mathrm{Witt}}$ is representable by an affine scheme.
\end{proposition}

\begin{proof} Picking a closed immersion $G\hookrightarrow \GL_n$, one can reduce to $G=\GL_n$. In that case, there is an ample line bundle $\mathcal L$ on $\Gr_G^{\mathrm{Witt}}$ constructed in \cite{BhattScholzeWitt}. We first claim that the pullback of $\mathcal L$ to $\Gr_B^{\mathrm{Witt}}$ is trivial. Indeed, recall that if $\Spec R\to \Gr_G^{\mathrm{Witt}}$ corresponds to a lattice $\Xi\subset W_{\mathcal O_E}(R)[\tfrac 1\pi]^n$, then $\mathcal L$ is given by $\det(\pi^{-m} W_{\mathcal O_E}(R)/\Xi)$ for any large enough $m$, using the determinant
\[
\det: \Perf(W_{\mathcal O_E}(R)\ \mathrm{on}\ R)\to \Pic(R),
\]
which is multiplicative in exact triangles. On $\Gr_B^{\mathrm{Witt}}$, one has a universal filtration of $\Xi$ compatible with the standard filtration on the standard lattice, which induces a similar filtration on $\Xi/\pi^m W_{\mathcal O_E}(R)$, where all the graded quotients are locally constant (and constant on $S_\lambda$). This means that the line bundle is naturally trivialized over each connected component $S_\lambda$ of $\Gr_B^{\mathrm{Witt}}$.

We claim that this section over $S_\lambda$ extends uniquely to a section over the closed subset $\bigcup_{\lambda'\leq \lambda} S_{\lambda'}$ that vanishes over the complement of $S_\lambda$, showing that the intersection of $S_\lambda$ with each $\Gr_{G,\leq \mu}^{\mathrm{Witt}}$ must be affine. To see this, by the v-descent results of \cite{BhattScholzeWitt}, it suffices to check that for any rank $1$ valuation ring $V$ with a map $\Spec V\to \Gr_G^{\mathrm{Witt}}$ whose generic point $\Spec K$ maps into $S_\lambda$, the section of $\mathcal L$ over $\Spec K$ extends to $\Spec V$ and is nonzero in the special fibre precisely when all of $\Spec V$ maps into $S_\lambda$. Now the filtration 
\[
0=\Xi_{K,0}\subset \Xi_{K,1}\subset \ldots \subset \Xi_{K,n}=\Xi_K
\]
with
\[
\Xi_{K,i} = \Xi_K\cap W_{\mathcal O_E}(K)[\tfrac{1}{\pi}]^i\subset W_{\mathcal O_E}(K)[\tfrac{1}{\pi}]^i
\]
has the property that
\[
\Xi_{K,i}/\Xi_{K,{i-1}} = \pi^{\lambda_i} W_{\mathcal O_E}(K)
\]
for the cocharacter $\lambda=(\lambda_1,\ldots,\lambda_n)$. Moreover, the filtration by the $\Xi_{K,i}$ extends integrally to the filtration
\[
0=\Xi_0\subset \Xi_1\subset \ldots \subset \Xi_n=\Xi
\]
with
\[
\Xi_i=\Xi\cap W_{\mathcal O_E}(V)[\tfrac{1}{\pi}]^i\subset W_{\mathcal O_E}(V)[\tfrac{1}{\pi}]^i,
\]
which is still a filtration by finite projective $W_{\mathcal O_E}(V)$-modules by \cite[Lemma 14.2.3]{Berkeley}. The injection of $\Xi_i/\Xi_{i-1}$ into $W_{\mathcal O_E}(V)[\tfrac{1}{\pi}]$ (projecting to the $i$-th coordinate) has image contained in
\[
W_{\mathcal O_E}(V)[\tfrac{1}{\pi}]\cap \pi^{\lambda_i} W_{\mathcal O_E}(K) = \pi^{\lambda_i} W_{\mathcal O_E}(V),
\]
so we get natural injections $\Xi_i/\Xi_{i-1}\hookrightarrow \pi^{\lambda_i} W_{\mathcal O_E}(V)$, that are isomorphisms after inverting $\pi$ or $[a]$ for a pseudouniformizer $a\in V$. Now the relevant line bundle can be written as the tensor product of the line bundles given by the determinants of the complexes
\[
\pi^{\lambda_i} W_{\mathcal O_E}(V)/(\Xi_i/\Xi_{i-1})\in \Perf(W_{\mathcal O_E}(R)\ \mathrm{on}\ R).
\]
These line bundles are indeed naturally trivial over $K$ as the perfect complex is acyclic there. Now this complex is concentrated in degree $0$, and is torsion, so admits a filtration by complexes of the form $V/a\cong [aV\hookrightarrow V]$ for pseudouniformizers $a\in V$. The associated line bundle on $V$ is then given by the alternating tensor product $V\otimes_V (aV)^{-1} = a^{-1} V$, and the natural section by $1\in a^{-1} V$. We see that the section is indeed integral, and that it is nonzero in the special fibre if and only if all the above complexes are acyclic, equivalently if $\Xi_i/\Xi_{i-1}\to \pi^{\lambda_i} W_{\mathcal O_E}(V)$ is an isomorphism. But this is precisely the condition that all of $\Spec V$ maps into $S_\lambda$.
\end{proof}

The union $\bigcup_{\lambda,\langle 2\rho,\lambda\rangle\leq d} S_\lambda\subset \Gr_G^{\mathrm{Witt}}$ is closed, thus so is
\[
\bigcup_{\lambda,\langle 2\rho,\lambda\rangle\leq d} S_\lambda\cap \Gr_{G,\leq\mu}^{\mathrm{Witt}}\subset \Gr_{G,\leq\mu}^{\mathrm{Witt}}.
\]
For $d=\langle 2\rho,\mu\rangle$, this is all of $\Gr_{G,\leq\mu}^{\mathrm{Witt}}$, while for $d=-\langle 2\rho,\mu\rangle$ it contains only a point, corresponding to $[\lambda]$ for $\lambda$ the antidominant representative of the Weyl orbit of $\mu$. Also, only $d$ of the same parity as $\langle 2\rho,\mu\rangle$ are relevant. By Proposition~\ref{prop:semiinfiniteorbitaffine}, the successive complements
\[
\bigcup_{\lambda,\langle 2\rho,\lambda\rangle\leq d} S_\lambda\cap \Gr_{G,\leq\mu}^{\mathrm{Witt}}\setminus \bigcup_{\lambda,\langle 2\rho,\lambda\rangle\leq d-2} S_\lambda\cap \Gr_{G,\leq\mu}^{\mathrm{Witt}} = \bigsqcup_{\lambda,\langle 2\rho,\lambda\rangle=d} S_\lambda\cap \Gr_{G,\leq\mu}^{\mathrm{Witt}}
\]
are affine. This means that at each step, the dimension can drop by at most $1$. However, in $\langle 2\rho,\mu\rangle$ steps, it drops by $\langle 2\rho,\mu\rangle$. We get the following corollary on Mirkovi\'c--Vilonen cycles, cf.~\cite[Theorem 3.2]{MirkovicVilonen}, and \cite{GHKR}, \cite[Corollary 2.8]{ZhuWitt} for a different proof based on point counting, the classical Satake isomorphism, and the Kato-Lusztig formula \cite{KatoKatoLusztig}, \cite{LusztigKatoLusztig}.

\begin{corollary}\label{cor:dimensionmirkovicvilonen} The scheme $S_\lambda\cap \Gr_{G,\leq\mu}^{\mathrm{Witt}}$ is equidimensional of dimension $\langle \rho,\mu+\lambda\rangle$.
\end{corollary}

\section{Equivariant sheaves}

Now we go back to the setting of general reductive groups $G$ over $\mathcal O_E$ (resp.~over $E$ if we work over $\Div^d_Y$ or $\Div^d_X$). As usual, let $\Lambda$ be some coefficient ring killed by some integer $n$ prime to $p$. We want to study $D_\et(-,\Lambda)$ for the local Hecke stack
\[
\Hloc_{G,\Div^d_{\mathcal Y}} = L^+_{\Div^d_{\mathcal Y}}G\backslash \Gr_{G,\Div^d_{\mathcal Y}}
\]
or its versions for $\Div^d_Y$ and $\Div^d_X$. Neither this nor its bounded versions $\Hloc_{G,\Div^d_{\mathcal Y},\leq \mu_\bullet}$ (say, when $G$ is split) is an Artin stack as $L^+_{\Div^d_{\mathcal Y}}G$ is not finite-dimensional. However, Proposition~\ref{prop:finitedimaction} shows that on the bounded version, the action factors over a finite-dimensional quotient.

First, we observe that on the level of $D_\et(-,\Lambda)$, one can then forget about the rest of the action.

\begin{proposition}\label{prop:finitedimquotientenough} Let $H$ be a group small v-sheaf over a small v-sheaf $S$ that admits a filtration $H^{\geq m}\subset H$ by closed subgroups such that, v-locally on $S$, for each $m\geq 1$ each quotient $H^{\geq m}/H^{\geq m+1}$ admits a further finite filtration with graded pieces given by $(\mathbb A^1_{S^\sharp})^\diamond$ for some untilt $S^\sharp$ of $S$ (that may depend on the graded piece). Let $X$ be some small v-sheaf over $S$ with an action of $H$ that factors over $H^{<m}=H/H^{\geq m}$ for some $m>0$. Then the pullback functor
\[
D_\et(H^{<m}\backslash X,\Lambda)\to D_\et(H\backslash X,\Lambda)
\]
is an equivalence.
\end{proposition}

\begin{proof} Both stacks live over the classifying stack $[H^{<m}\backslash S]$ of $H^{<m}$ over $S$. Applying descent along $S\to [H^{<m}\backslash S]$, one reduces to the case that $H^{<m}$ is trivial. In that case, the map $X/H\to X/H^{<m}=X$ admits a section $s: X\to X/H$, and it is enough to prove that $s^\ast$ is fully faithful. Doing descent once more, it is enough to prove that for any affinoid perfectoid space $S'$ over $S$ over which a filtration by $\mathbb A^1$'s exists, pullback
\[
D_\et(S',\Lambda)\to D_\et(S'\times H,\Lambda)
\]
is fully faithful. Replace $S$ by $S'$ and let $f: H\to S$ be the projection. We need to see that for all $A\in D_\et(S,\Lambda)$, the map
\[
A\to Rf_\ast f^\ast A
\]
is an isomorphism; doing this for all $S$, it is enough to check it on global sections, i.e.
\[
R\Gamma(S,A)\to R\Gamma(S\times H,f^\ast A)
\]
is an isomorphism. Using Postnikov towers, we can assume that $A\in D^+_\et(S,\Lambda)$. We can write $H$ as a filtered colimit of subgroups $H_j\subset H$ such that each $H_j$ is a successive extension as before, but now the quotients are balls inside each $\mathbb A^1_{S^\sharp}$. In particular, each $H_j$ is a spatial diamond, and it is enough to prove that
\[
R\Gamma(S,A)\to R\Gamma(S\times H_j, f^\ast A)
\]
is an isomorphism for all $j$. Now each $H_j=\varprojlim_m H_j^{<m}$ is an inverse limit of spatial diamonds, so by \cite[Proposition 14.9]{ECoD} it is enough to prove that
\[
R\Gamma(S,A)\to R\Gamma(S\times H_j^{<m},f_n^\ast A)
\]
is an isomorphism for all $m$ and $j$. But this follows easily from each $H_j^{<m}$ being a successive extension of balls inside $\mathbb A^1_{S^\sharp}$.
\end{proof}

Using hyperbolic localization, we can prove the following important conservativity result.

\begin{proposition}\label{prop:conservativity} Assume that $B\subset G$ is a Borel. Let $S\to \Div^d_{\mathcal Y}$ be any small v-sheaf. Let $A\in D_\et(\Hloc_{G,S/\Div^d_{\mathcal Y}},\Lambda)$ with support quasicompact over $S$. Assume that the hyperbolic localization $\mathrm{CT}_B(A)=0$ of the pullback of $A$ to $\Gr_{G,S/\Div^d_{\mathcal Y}}$ vanishes. Then $A=0$.

The similar assertion holds with $\Div^d_Y$ and $\Div^d_X$ in place of $\Div^d_{\mathcal Y}$.
\end{proposition}

\begin{proof} Note that the formation of $\mathrm{CT}_B$ commutes with any base change in $S$, by Corollary~\ref{cor:hyperboliclocalizationGr}. We can thus assume that $S=\Spa(C,C^+)$ is strictly local. Up to replacing $d$ by a smaller integer, removing double points, we can assume that the map $S\to \Div^d_{\mathcal Y}$ is given by $d$ distinct untilts $S_i^\sharp$ over $\mathcal O_E$, $i=1,\ldots,d$. Let $E'|E$ be an extension splitting $G$, assumed unramified in our situation where we work over $\mathcal Y$. We can then lift all $S_i^\sharp$ to $\mathcal O_{E'}$, and thereby reduce to the case of split $G$. The corresponding geometric fibre
\[
\Hloc_{G,S/\Div^d_{\mathcal Y}}
\]
has a stratification enumerated by $\mu_1,\ldots,\mu_d\in X_\ast(T)^+$, with strata
\[
[S/(\prod_{i=1}^d (L^+_{\Div^1_{\mathcal Y}} G)_{\mu_i}\times_{\Div^1_{\mathcal Y}} S)].
\]
If $A$ is nonzero, we can find a maximal such stratum on which $A$ is nonzero. Now we apply Corollary~\ref{cor:hyperboliclocalizationGr}, see Example \ref{ex:suite semi infinite orbit ex}. One has an isomorphism
\[
S\times X_\ast(T)^d\cong \Gr_{T,S/\Div^d_{\mathcal Y}}.
\]
Over the copy of $S$ enumerated by the antidominant representatives of (the Weyl group orbits of) $\mu_1,\ldots,\mu_n$ the functor $\mathrm{CT}_B$ is the pullback of $A$ to a section of the stratum corresponding to $\mu_1,\ldots,\mu_n\in X_\ast(T)^+$ (which, as we recall, correspond to a maximal stratum where $A$ is nonzero). This shows that the restriction of $A$ to a section over this maximal stratum is zero. This gives the desired contradiction, so $A=0$.
\end{proof}

\section{Affine flag variety}\label{sec:affineflag}

At a few isolated spots, it will be useful to use the affine flag variety, the main point being that the Schubert varieties in the affine flag variety admit explicit resolutions of singularities, given by Demazure resolutions (also known as Bott--Samelson resolutions). It will be enough to appeal to these in the setting of a split reductive group $G$, with a reductive model over $\mathcal O_E$ and Borel $B\subset G$ defined over $\mathcal O_E$, for $d=1$, and for a small v-stack $S\to \Div^1_{\mathcal Y}$ factoring over $\Spd \mathcal O_C$ where $C=\widehat{\overline{E}}$, so we restrict attention to this setting.

Consider the base change $G_A$ of $G$ to $A=W_{\mathcal O_E}(\mathcal O_{C^\flat})$. We have Fontaine's map $\theta: A\to \mathcal O_C$, and we can define an ``Iwahori'' group scheme $\mathcal I\to G_A$, flat over $A$, whose points in a $\ker\theta$-torsionfree $A$-algebra $R$ are given those elements $g\in G(R)$ such that $\theta(g)\in G(R\otimes_A \mathcal O_C)$ lies in $B(R\otimes_A \mathcal O_C)$. Similarly, for any parabolic $P\subset G$ containing $B$, we get a ``parahoric'' group scheme $\mathcal P\to G_A$, flat over $A$, whose points in a $\ker\theta$-torsionfree $A$-algebra $R$ are those $g\in G(R)$ such that $\theta(g)\in P(R\otimes_A \mathcal O_C)$. In particular, this applies to the parabolics $P_i$ corresponding to the simple reflections $s_i$; let $\mathcal P_i$ be the corresponding parahorics. Still more generally, for any affine simple reflection $s_i$, one can define a parabolic $\mathcal P_i\to G_A$ flat over $A$, and such that $\mathcal I\to G_A$ factors over $\mathcal P_i$. (The construction of these parahoric group schemes over $A$ can be reduced to the case of $W_{\mathcal O_E}(k)[[u]]$ via a faithfully flat embedding $W_{\mathcal O_E}(k)[[u]]\hookrightarrow A$ along which everything arises via base change, and then one can appeal to the work of Bruhat--Tits \cite[Section 3.9.4]{BruhatTits}.)

\begin{definition}\label{def:affineflag} In the situation above, including a small v-stack $S$ over $\Spd \mathcal O_C$, mapping to $\Div^1_{\mathcal Y}$, let
\[
\Fl_{G,S}\to S
\]
be the \'etale quotient $LG/L^+\mathcal I$, where $L^+\mathcal I(R,R^+)=\mathcal I(B_{\mathrm{dR}}^+(R^\sharp))$.
\end{definition}

Note here that as $S$ lives over $\Spd \mathcal O_C$, any $\Spa(R,R^+)$ over $S$ comes with an untilt $R^\sharp$ over $\mathcal O_C$, in which case $B_{\mathrm{dR}}^+(R^\sharp)$ is an $A=W_{\mathcal O_E}(\mathcal O_C)$-algebra, so that $\mathcal I(B_{\mathrm{dR}}^+(R^\sharp))$ is well-defined.

\begin{proposition}\label{prop:affineflagtoaffgr} There is a natural projection map
\[
\Fl_{G,S}\to \Gr_{G,S/\Div^1_{\mathcal Y}}
\]
that is v-locally isomorphic to a product with $(G/B)^\diamond$. In particular, it is proper, representable in spatial diamonds, and cohomologically smooth.
\end{proposition}

\begin{proof} This follows from the identification $L^+G/L^+\mathcal I\cong (G/B)^\diamond$, which follows from the definition, and the similar properties of $(G/B)^\diamond\to \Spd \mathcal O_E$.
\end{proof}

We analyze the stratification of $\Fl_{G,S}$ into $L^+\mathcal I$-orbits. Let $N(T)\subset T$ be the normalizer of $T$, and
\[
\tilde{W} = N(T)(B_\dR(C'))/T(B_\dR^+(C'))
\]
be the affine Weyl group, for any complete algebraically closed field $C'$ over $\mathcal O_E$ with a map $C^\flat\to C'^\flat$; this is naturally independent of the choice of $C'$. As $T(B_\dR(C'))/T(B_\dR^+(C'))\cong X_\ast(T)$, there is a short exact sequence
\[
0\to X_\ast(T)\to \tilde{W}\to W\to 0,
\]
where $W$ is the usual Weyl group of $G$.

\begin{proposition}\label{prop:bruhatdecomposition} The decomposition of $\Fl_{G,S}(C')$ into $L^+ \mathcal I(C')$-orbits is given by
\[
\Fl_{G,S}(C') = \bigsqcup_{w\in \tilde{W}} L^+\mathcal I(C)\cdot w.
\]
\end{proposition}

\begin{proof} If $C'$ lives over $E$, we can choose an isomorphism $B_\dR(C')\cong C'((\xi))$ and the result follows from the classical result. If $C'$ lives over the residue field $\Fq$ of $E$, this reduces to the assertion for the Witt vector affine flag variety, for which we refer to \cite{ZhuWitt}.
\end{proof}

Recall that $\tilde{W}$ acts on $X_\ast(T)$. Fixing the alcove $\mathfrak a$ corresponding to the Iwahori group $\mathcal I$, one gets a set of affine simple reflections $s_i$ as the reflections along the faces of the alcove; these generate a normal subgroup $W_{\mathrm{aff}}\subset \tilde{W}$. Letting $\Omega\subset \tilde{W}$ denote the stabilizer of the alcove, there is a split short exact sequence
\[
1\to W_{\mathrm{aff}}\to \tilde{W}\to \Omega\to 1.
\]
One gets the Bruhat order on $\tilde{W}$: If $w_i=w_{i,0}\omega_i\in \tilde{W} = W_{\mathrm{aff}}\rtimes \Omega$ for $i=1,2$ are two elements, then $w_1\leq w_2$ if $\omega_1=\omega_2$ and in one (hence every) presentation of $w_2$ as a product of affine simple reflections, $w_1$ is obtained by removing some factors.

\begin{definition}\label{def:affineflagschubert} For $w\in \tilde{W}$, the affine Schubert cell is the subfunctor $\Fl_{G,w,S}\subset \Fl_{G,S}$ of all maps $\Spa(R,R^+)\to \Fl_{G,S}$ that on all geometric points lie in the $L^+ \mathcal I$-orbit of $w$. The affine Schubert variety is the subfunctor $\Fl_{G,\leq w,S}$ of all maps $\Spa(R,R^+)\to \Fl_{G,S}$ that on all geometric points lie in the $L^+\mathcal I$-orbit of $w^\prime$ for some $w^\prime\leq w$.
\end{definition}

\begin{theorem}\label{thm:affineflagschubert} For each $w\in \tilde{W}$, the subfunctor $\Fl_{G,\leq w,S}\subset \Fl_{G,S}$ is closed, and $\Fl_{G,\leq w,S}\to \Spd \mathcal O_C$ is proper and representable in spatial diamonds, of finite $\dimtrg$. The subfunctor $\Fl_{G,w,S}\subset \Fl_{G,\leq w,S}$ is open and dense.
\end{theorem}

\begin{proof} We will prove the theorem by constructing the Demazure resolution of $\Fl_{G,\leq w,S}$. Write $w=w_0\omega\in \tilde{W}=W_{\mathrm{aff}}\rtimes\Omega$, and fix a decomposition $w_0=\prod_{j=1}^l s_{i_j}$ as a product of affine simple reflections of minimal length, so $l(w)=l(w_0)=l$. We write $\dot{w}$ for the element $w$ with such a choice of decomposition.

For each affine simple reflection $s_i$, we have a corresponding parahoric group $\mathcal P_i\to G_A$ corresponding to the face of $\mathfrak a$; one has $L^+ \mathcal P_i/L^+ \mathcal I\cong (\mathbb P^1)^\diamond$.

\begin{definition}\label{def:demazure} The Demazure variety corresponding to $\dot{w}$ is the \'etale sheaf
\[
\mathrm{Dem}_{\dot{w},S} = L^+ \mathcal P_{s_{i_1}}\times^{L^+\mathcal I} L^+\mathcal P_{s_{i_2}}\times^{L^+\mathcal I} \ldots \times^{L^+\mathcal I} L^+\mathcal P_{s_{i_l}}/L^+\mathcal I\to S,
\]
equipped with the left $L^+\mathcal I$-action and the $L^+\mathcal I$-equivariant map
\[
\mathrm{Dem}_{\dot{w},S}\to \Fl_{G,S}
\]
given by $(p_1,\ldots,p_l)\mapsto p_1\cdots p_l\cdot \omega$.
\end{definition}

It is clear from the definition that $\mathrm{Dem}_{\dot{w}}\to S$ is a successive $(\mathbb P^1)^\diamond$-fibration over $S$, and in particular is a spatial diamond, proper over $S$ of finite $\dimtrg$. As $\Fl_{G,S}\to S$ is partially proper, it follows that the image of $\mathrm{Dem}_{\dot{w},S}\to \Fl_{G,S}$ is proper. Moreover, the image can be identified on geometric points, and we see that $\mathrm{Dem}_{\dot{w},S}\to \Fl_{G,\leq w,S}$ is surjective, $\Fl_{G,\leq w,S}\subset \Fl_{G,S}$ is closed, and $\Fl_{G,\leq w,S}$ is proper over $S$. In particular, $\Fl_{G,w,S}\subset \Fl_{G,\leq w,S}$ is open, as the complement is a finite union of closed subfunctors. As $\Fl_{G,S}\to \Gr_{G,S}$ is locally a product with $(G/B)^\diamond$, it follows from \cite[Theorem 19.2.4]{Berkeley} that $\Fl_{G,S}\times_{\Gr_{G,S}} \Gr_{G,\leq \mu,S}$ is a spatial diamond, and thus so is $\Fl_{G,\leq w,S}$, as it is a closed subspace for $\mu$ large enough.

Also, by checking on geometric points and reducing to the classical case, the map $\mathrm{Dem}_{\dot{w},S}\to \Fl_{G,\leq w,S}$ is an isomorphism over $\Fl_{G,w,S}$ whose preimage is given by
\[
(L^+ \mathcal P_{s_{i_1}}\setminus L^+\mathcal I)\times^{L^+\mathcal I} (L^+\mathcal P_{s_{i_2}}\setminus L^+\mathcal I)\times^{L^+\mathcal I} \ldots \times^{L^+\mathcal I} (L^+\mathcal P_{s_{i_l}}\setminus L^+\mathcal I)/L^+\mathcal I.
\]
This implies that $\Fl_{G,w,S}\subset \Fl_{G,\leq w,S}$ is dense, as desired. As usual, a consequence of this discussion is that the Bruhat order is independent of the choice of $\dot{w}$.
\end{proof}

Using Demazure resolutions, one can prove that the standard sheaves on the affine flag variety are universally locally acyclic.

\begin{proposition}\label{prop:ULAaffineflag} For any $w\in \tilde{W}$, let $j_w: \Fl_{G,w,S}\hookrightarrow \Fl_{G,\leq w,S}$ be the open embedding. Then $j_{w!}\Lambda\in D_\et(\Fl_{G,\leq w,S},\Lambda)$ is universally locally acyclic over $S$.
\end{proposition}

\begin{proof} Using Proposition~\ref{prop:properULA}, it suffices to prove the same for $\tilde{j}_w: \Fl_{G,w,S}\hookrightarrow \mathrm{Dem}_{\dot{w},S}$ and $\tilde{j}_{w!}\Lambda$. Then $\tilde{j}_{w!}\Lambda$ can be resolved in terms of $\Lambda$ and all $i_{\dot{w}',\dot{w},\ast}\Lambda$ for
\[
i_{\dot{w}',\dot{w}}: \mathrm{Dem}_{\dot{w}',S}\to \mathrm{Dem}_{\dot{w},S}
\]
the closed immersion from another Demazure variety, corresponding to a subword of $\dot{w}'$ of $\dot{w}$; note that combinatorially, we are dealing with the situation of a normal crossing divisor at the boundary. By cohomological smoothness of all $\mathrm{Dem}_{\dot{w}',S}\to S$ and Proposition~\ref{prop:properULA}, the result follows.
\end{proof}

\section{ULA sheaves}

We will be interested in universally locally acyclic sheaves on the local Hecke stack.

\begin{definition}\label{def:localHeckeULA} Let $S\to \Div^d_{\mathcal Y}$ be any small v-stack. An object
\[
A\in D_\et(\Hloc_{G,S/\Div^d_{\mathcal Y}},\Lambda)
\]
is universally locally acyclic over $S$ if it is bounded, and its pullback to
\[
\Gr_{G,S/\Div^d_{\mathcal Y}}
\]
is universally locally acyclic over $S$.

Let
\[
D^\ULA_\et(\Hloc_{G,S/\Div^d_{\mathcal Y}},\Lambda)\subset D_\et(\Hloc_{G,S/\Div^d_{\mathcal Y}},\Lambda)
\]
be the corresponding full subcategory.
\end{definition}

This definition is a priori not symmetric in the two bundles $\mathcal E_1$, $\mathcal E_2$ parametrized by the local Hecke stack. However, we can check that it actually is.

\begin{proposition}\label{prop:ULAsymmetry} Consider the automorphism
\[
\mathrm{sw}: \Hloc_{G,S/\Div^d_{\mathcal Y}}\cong \Hloc_{G,S/\Div^d_{\mathcal Y}}
\]
switching $\mathcal E_1$ and $\mathcal E_2$. Then $A\in D_\et(\Hloc_{G,S/\Div^d_{\mathcal Y}},\Lambda)$ is universally locally acyclic over $S$ if and only if $\mathrm{sw}^\ast A$ is universally locally acyclic over $S$.
\end{proposition}

\begin{proof} Fix any large enough substack $U\subset \Hloc_{G,S/\Div^d_{\mathcal Y}}$ quasicompact over $S$ containing the support of $A$. Let $(L_{\Div^d_{\mathcal Y}} G)_U\subset L_{\Div^d_{\mathcal Y}} G$ be the preimage of $U$. Universal local acyclicity after pullback to $\Gr_{G,\Div^d_{\mathcal Y}}$ is equivalent to universal local acyclicity after pullback to
\[
(L_{\Div^d_{\mathcal Y}} G)_U/(L^+_{\Div^d_{\mathcal Y}} G)^{\geq m}
\]
for any $m>0$, by Proposition~\ref{prop:gradedcongruencesubgroups2} and Proposition~\ref{prop:firstquotientsmooth}. We need to see that this is equivalent to universal local acyclicity after pullback to
\[
(L^+_{\Div^d_{\mathcal Y}} G)^{\geq m}\backslash (L_{\Div^d_{\mathcal Y}} G)_U
\]
for any $m>0$. For this, we note that these two pro-systems in $m$ are pro-isomorphic. By the next lemma, the transition maps back and forth are also cohomologically smooth, which implies the desired equivalence.
\end{proof}

In the following lemma, we call a map $f$ universally locally acyclic if $\Lambda$ is $f$-universally locally acyclic.

\begin{lemma}\label{lem:twooutofsixsmooth} Let
\[
X_4\xrightarrow{f_3} X_3\xrightarrow{f_2} X_2\xrightarrow{f_1} X_1\xrightarrow{f_0} X_0
\]
be surjective maps of locally spatial diamonds that are compactifiable and of locally finite $\dimtrg$. Assume that $f_0\circ f_1$ and $f_1\circ f_2$ are cohomologically smooth. Then $f_0$ is universally locally acyclic. If $f_1$ is universally locally acyclic and $f_0\circ f_1$ is cohomologically smooth, then $f_1$ is cohomologically smooth. Thus, if $f_0\circ f_1$, $f_1\circ f_2$ and $f_2\circ f_3$ are cohomologically smooth, then $f_0$ and $f_1$ are cohomologically smooth.
\end{lemma}

We would expect that $f_3: X_4\to X_3$ should be unnecessary in order for $f_0$ to be cohomologically smooth.

\begin{proof} We claim that for any map $g_0: Y_0\to X_0$, with pullbacks $g_i: Y_i\to X_i$ and $\tilde{f}_i: Y_{i+1}\to Y_i$, the natural transformation
\[
\tilde{f}_0^\ast Rg_0^!\to Rg_1^! f_0^\ast
\]
is an isomorphism. Indeed, we have natural maps
\[
\tilde{f}_2^\ast \tilde{f}_1^\ast \tilde{f}_0^\ast Rg_0^!\to \tilde{f}_2^\ast \tilde{f}_1^\ast Rg_1^! f_0^\ast\to \tilde{f}_2^\ast Rg_2^! f_1^\ast f_0^\ast\to Rg_3^! f_2^\ast f_1^\ast f_0^\ast
\]
and the composite of any two maps is an isomorphism. By the two-out-of-six-lemma, this implies that all maps are isomorphisms. By surjectivity of $f_1$ and $f_2$, this implies that $\tilde{f}_0^\ast Rg_0^!\to Rg_1^! f_0^\ast$ is an isomorphism. Applying this with $Y_0=X_1$ and to the constant sheaf $\Lambda$ then shows, by the criterion of Theorem~\ref{thm:ULAdualizable}, that $\Lambda$ is $f_0$-universally locally acyclic.

Now assume that $f_1$ is universally locally acyclic and $f_0\circ f_1$ is cohomologically smooth, then
\[
R(f_0\circ f_1)^!\Lambda\cong f_1^\ast Rf_0^!\Lambda\dotimes_\Lambda Rf_1^!\Lambda
\]
is invertible. This implies that both tensor factors are invertible, and in particular $Rf_1^!\Lambda$ is invertible, so $f_1$ is cohomologically smooth. For the final statement, we now know that the hypotheses imply that $f_0$ and $f_1$ are universally locally acyclic, so the displayed equation implies that $f_0$ and $f_1$ are cohomologically smooth.
\end{proof}

Using the conversativity result Proposition~\ref{prop:conservativity}, we can characterize universally locally acyclic sheaves in terms of their hyperbolic localization. Note that we can always reduce to the case of quasisplit $G$ by \'etale localization on $S$.

\begin{proposition}\label{prop:ULAGrcohomology} Let $B\subset G$ be a Borel with torus quotient $T$. Let $S$ be a small v-stack with a map $S\to \Div^d_{\mathcal Y}$, and let
\[
A\in D_\et(\Hloc_{G,S/\Div^d_{\mathcal Y}},\Lambda)^{\mathrm{bd}}.
\]
Then $A$ is universally locally acyclic over $S$ if and only if the hyperbolic localization
\[
\mathrm{CT}_B(A)\in D_\et(\Gr_{T,S/\Div^d_{\mathcal Y}},\Lambda)^{\mathrm{bd}}
\]
is universally locally acyclic over $S$. This, in turn, is equivalent to the property that
\[
R\pi_{T,S,\ast} \mathrm{CT}_B(A)\in D_\et(S,\Lambda)
\]
is locally constant with perfect fibres.
\end{proposition}

Here
\[
\pi_{T,S}: \Gr_{T,S/\Div^d_{\mathcal Y}}\to S
\]
is the projection.

\begin{proof} The forward direction follows from Corollary~\ref{cor:hyperboliclocalizationGr} and the ind-properness of $\pi_{T,S}$ and Corollary \ref{coro:proper push locally acyclic}. For the converse direction, we may assume that $S$ is strictly totally disconnected and $G$ is split. Note that to prove universal local acyclicity of $A$, it is enough to prove that the map
\[
p_1^\ast R\sHom(A,R\pi_{G,S}^!\Lambda)\dotimes_\Lambda p_2^\ast A\to R\sHom(p_1^\ast A,Rp_2^! A)
\]
is an isomorphism (by Theorem~\ref{thm:ULAdualizable}). (Implicitly, we pass here to a bounded part of $\Hloc_{G,\Div^d_{\mathcal Y}}$ and replace the quotient by $L^+_{\Div^d_{\mathcal Y}} G$ by a finite-dimensional quotient in order to be in the setting of Artin stacks.) By Proposition~\ref{prop:conservativity} applied to $G\times G$, it is enough to prove this after applying $\mathrm{CT}_{B^-\times B}$, where $B^-$ is the opposite Borel. Using that hyperbolic localization commutes with exterior tensor products, and Proposition~\ref{prop:hyperboliclocalizationdual}, this translates exactly into the similar isomorphism characterizing universal local acyclicity of $\mathrm{CT}_B(A)$. The final statement follows from Proposition~\ref{prop:properproetaleULA}.
\end{proof}

In the case of one leg, one can completely characterize universally locally acyclic sheaves.

\begin{proposition}\label{prop:ULAoneleg} Assume that $G$ is split. Let $S\to \Div^1_{\mathcal Y}$ be any small v-stack. Consider
\[
A\in D_\et(\Hloc_{G,S/\Div^1_{\mathcal Y}},\Lambda)^{\mathrm{bd}}.
\]
Then $A$ is universally locally acyclic over $S$ if and only if for all $\mu\in X_\ast(T)^+$, the restriction of $A$ to the section $[\mu]: S\to \Hloc_{G,S/\Div^1_{\mathcal Y}}$ is locally constant with perfect fibres in $D_\et(S,\Lambda)$.
\end{proposition}

If $G$ is not split, a similar characterization holds, by applying the result \'etale locally to reduce to the case of split $G$. Again, there is also the obvious version for $\Div^1_Y$ and $\Div^1_X$.

\begin{proof} First, we prove that if all fibres are locally constant with perfect fibres, then $A$ is universally locally acyclic. This easily reduces to the case of $j_{\mu!}\Lambda$ where
\[
j_\mu: \Hloc_{G,\Div^1_{\mathcal Y},\mu}\hookrightarrow \Hloc_{G,\Div^1_{\mathcal Y}}
\]
is the inclusion of an open Schubert cell, and $S=\Div^1_{\mathcal Y}$. We can also argue v-locally on $\Div^1_{\mathcal Y}$ and so base change to the case $S=\Spd \mathcal O_C$. In that case, Proposition~\ref{prop:affineflagtoaffgr} and Proposition~\ref{prop:ULAsmoothlocal} show that it suffices to prove the similar assertion for the affine flag variety, where it follows from Proposition~\ref{prop:ULAaffineflag}.

Now for the converse, we argue by induction on the support of $A$. On a maximal Schubert cell $\Gr_{G,S/\Div^1_{\mathcal Y},\mu}$ where $A$ is nonzero, its restriction is universally locally acyclic, and as on the Hecke stack this stratum is the classifying space of a (pro-)cohomologically smooth group, it follows that the restriction of $A$ along the section $[\mu]: S\to \Hloc_{G,S/\Div^1_{\mathcal Y}}$ is locally constant with perfect fibres. Replacing $A$ by the cone of
\[
j_{\mu!} A|_{\Hloc_{G,S/\Div^1_{\mathcal Y},\mu}}\to A,
\]
the claim follows.
\end{proof}

In the following corollaries, we no longer assume that $G$ is split.

\begin{corollary}\label{cor:sixoperationsULAGr} Let $S\to \Div^1_{\mathcal Y}$ be any small v-stack. Then
\[
D_\et^\ULA(\Hloc_{G,S/\Div^1_{\mathcal Y}},\Lambda)
\]
is stable under Verdier duality and $-\dotimes_\Lambda-$, $R\sHom_\Lambda(-,-)$ as well as $j_! j^\ast$, $Rj_\ast j^\ast$, $j_! Rj^!$, $Rj_\ast Rj^!$ where $j$ is the locally closed immersion of a Schubert cell. Moreover, all of these operations commute with all pullbacks in $S$.
\end{corollary}

\begin{proof} Stability under Verdier duality and compatibility with base change in $S$ follow from Corollary~\ref{cor:ULAselfdual}. For the other assertions, one can reduce to the case that $G$ is split by working locally on $S$, where it follows from the previous proposition, and juggling with the six functors, and the stalkwise characterization of the previous proposition.
\end{proof}

\begin{corollary}\label{cor:degenerateSatake} For a complete algebraically closed extension $C$ of $E$ with residue field $k$, taking $S=\Spd \mathcal O_C$, $S=\Spd C$ and $S=\Spd k$, the functors
\[
D_\et^\ULA(\Hloc_{G,\Spd C/\Div^1_{\mathcal Y}},\Lambda)\leftarrow D_\et^\ULA(\Hloc_{G,\Spd \mathcal O_C/\Div^1_{\mathcal Y}},\Lambda)\to D_\et^\ULA(\Hloc_{G,\Spd k/\Div^1_{\mathcal Y}},\Lambda)
\]
are equivalences.
\end{corollary}

\begin{proof} Use that the formation of $R\sHom$ commutes with any base change in $S$, and that the category of locally constant sheaves with perfect fibres on any such $S$ is equivalent to the category of perfect $\Lambda$-modules.
\end{proof}

In fact, the previous results extend to the case of general $d$ as long as $S\to \Div^d_{\mathcal Y}$ has image in the open subset $(\Div^d_{\mathcal Y})_{\neq}\subset \Div^d_{\mathcal Y}$ where all untilts are distinct. After passing to a finite \'etale cover of $S$, we can then in fact assume that $S$ maps to $(\Div^1_{\mathcal Y})^d_{\neq}$.

\begin{proposition}\label{prop:ULAdistinctlegs} Assume that $G$ is split. Let $S\to (\Div^1_{\mathcal Y})^d_{\neq}\to \Div^d_{\mathcal Y}$ be a small v-stack. Consider
\[
A\in D_\et(\Hloc_{G,S/\Div^d_{\mathcal Y}},\Lambda)^{\mathrm{bd}}.
\]
Then $A$ is universally locally acyclic over $S$ if and only if for all $\mu_1,\ldots,\mu_d\in X_\ast(T)^+$, the restriction of $A$ to the section $[\mu_\bullet]: S\to \Hloc_{G,S/\Div^d_{\mathcal Y}}$ is locally constant with perfect fibres in $D_\et(S,\Lambda)$.

The category
\[
D_\et^\ULA(\Hloc_{G,S/\Div^d_{\mathcal Y}},\Lambda)
\]
is stable under Verdier duality and $-\dotimes_\Lambda-$, $R\sHom_\Lambda(-,-)$ as well as $j_! j^\ast$, $Rj_\ast j^\ast$, $j_! Rj^!$, $Rj_\ast Rj^!$ where $j$ is the locally closed immersion of a Schubert cell. Moreover, all of these operations commute with all pullbacks in $S$.
\end{proposition}

\begin{proof} We have the decomposition
\[
\Hloc_{G,S/\Div^d_{\mathcal Y}}\cong \prod_{i=1}^d \Hloc_{G,S/_{\pi_i}\Div^1_{\mathcal Y}}
\]
where $\pi_1,\ldots,\pi_d: S\to \Div^1_{\mathcal Y}$ are the $d$ projections, and the product on the right is taken over $S$. One can then stratify according to Schubert cells parametrized by tuples $\mu_\bullet=(\mu_1,\ldots,\mu_d)$ and the above arguments imply the result. Here, in the beginning, to see that $j_{\mu_\bullet !}\Lambda$ is universally locally acyclic, one uses that exterior tensor products preserve universal local acyclicity, see Corollary~\ref{cor:ULAselfdual}, to reduce to the case of one leg.
\end{proof}

\section{Perverse Sheaves}\label{sec:perverse}

For any small v-stack $S\to \Div^d_{\mathcal Y}$, we define a (relative) perverse $t$-structure on
\[
D_\et(\Hloc_{G,S/\Div^d_{\mathcal Y}},\Lambda)^{\mathrm{bd}}.
\]

\begin{defprop}\label{def:perversetstructure} Let $S\to \Div^d_{\mathcal Y}$ be a small v-stack. There is a unique $t$-structure $({}^p D^{\leq 0}, {}^p D^{\geq 0})$ on $D_\et(\Hloc_{G,S/\Div^d_{\mathcal Y}},\Lambda)^{\mathrm{bd}}$ such that
\[
A\in {}^p D_\et^{\leq 0}(\Hloc_{G,S/\Div^d_{\mathcal Y}},\Lambda)^{\mathrm{bd}}
\]
if and only if for all geometric points $\Spa(C,C^+)\to S$ and open Schubert cells of
\[
\Hloc_{G,\Spa(C,C^+)/\Div^d_{\mathcal Y}},
\]
parametrized by some $\mu_1,\ldots,\mu_r\in X_\ast(T)^+$ (where $r$ is the number of distinct untilts at $\Spa(C,C^+)\to \Div^d_{\mathcal Y}$), the pullback of $A$ to this open Schubert cell sits in cohomological degrees $\leq - \sum_{i=1}^r \langle 2\rho,\mu_i\rangle$.
\end{defprop}

\begin{proof} We note that on any bounded closed subset of $Z\subset \Hloc_{G,\Div^d_{\mathcal Y}}$ there is a presentable stable $\infty$-category $\mathcal D_\et(Z\times_{\Div^d_{\mathcal Y}} S,\Lambda)$ refining the derived category, and the given class of objects is stable under all colimits and extensions (and is generated by a set of objects). Thus, the existence and uniqueness of the $t$-structure follow from \cite[Proposition 1.4.4.11]{LurieHA}. Moreover, one easily checks that when enlarging $Z$, the inclusion functors are $t$-exact, so these glue to a $t$-structure in the direct limit.
\end{proof}

Let
\[
\Perv(\Hloc_{G,S/\Div^d_{\mathcal Y}},\Lambda)\subset D_\et(\Hloc_{G,S/\Div^d_{\mathcal Y}},\Lambda)^{\mathrm{bd}}
\]
be the heart of the perverse $t$-structure. On it, pullback to the affine Grassmannian is fully faithful.

\begin{proposition}\label{prop:pullbackperversetoaffGrass} The pullback functor
\[
\Perv(\Hloc_{G,S/\Div^d_{\mathcal Y}},\Lambda)\to D_\et(\Gr_{G,S/\Div^d_{\mathcal Y}},\Lambda)^{\mathrm{bd}}
\]
is fully faithful.

Moreover, if
\[
A\in {}^p D_\et^{\leq 0}(\Hloc_{G,S/\Div^d_{\mathcal Y}},\Lambda)^{\mathrm{bd}}\ \mathrm{and}\ B\in {}^p D_\et^{\geq 0}(\Hloc_{G,S/\Div^d_{\mathcal Y}},\Lambda)^{\mathrm{bd}},
\]
then $R\sHom_\Lambda(A,B)\in D^{\geq 0}_\et(\Hloc_{G,S/\Div^d_{\mathcal Y}},\Lambda)^{\mathrm{bd}}$.
\end{proposition}

\begin{proof} For the final statement, we need to see that if $C\in D^{\leq -1}_\et(\Hloc_{G,S/\Div^d_{\mathcal Y}},\Lambda)^{\mathrm{bd}}$, then there are no nonzero maps $C\to R\sHom_\Lambda(A,B)$; equivalently, there are no nonzero maps $C\dotimes_\Lambda A\to B$. But this follows from the simple observation that $C\dotimes_\Lambda A$ lies in ${}^p D^{\leq -1}$.

Now using this property of $R\sHom_\Lambda(A,B)$, descent implies that it is enough to see that if $A,B\in \Perv(\Hloc_{G,S/\Div^d_{\mathcal Y}},\Lambda)$, then any map between their pullbacks to $\Gr_{G,S/\Div^d_{\mathcal Y}}$ is automatically invariant under the action of $L^+_{\Div^d_{\mathcal Y}} G$. This follows from Lemma~\ref{lem:connectedgroupsacttrivially} applied to a finite-dimensional approximation of
\[
\Gr_{G,S/\Div^d_{\mathcal Y}}\times_{\Div^d_{\mathcal Y}} L^+_{\Div^d_{\mathcal Y}} G \to \Gr_{G,S/\Div^d_{\mathcal Y}}.\qedhere
\]
\end{proof}

We used the following lemma on actions of connected groups on \'etale sheaves.

\begin{lemma}\label{lem:connectedgroupsacttrivially} Let $f: Y\to X$ be a compactifiable cohomologically smooth map of locally spatial diamonds with a section $s: X\to Y$. Assume that all geometric fibres of $f$ are connected. Then for all $A\in D_\et^{\geq 0}(X,\Lambda)$, the map
\[
\mathcal H^0(Rf_\ast f^\ast A)\to \mathcal H^0(A)
\]
given by evaluation at the section $s$ is an isomorphism.
\end{lemma}

\begin{proof} Note that
\[
Rf_\ast f^\ast A\cong R\sHom(Rf_! Rf^!\Lambda,A)
\]
where $Rf_! Rf^! \Lambda$ sits in cohomological degrees $\leq 0$ with $\mathcal H^0\cong \Lambda$. Indeed, this reduces easily to the case of discrete $\Lambda$, and then to $\Lambda=\mathbb F_\ell$, and can be checked on geometric stalks. But if $S=\Spa(C,C^+)$ and $i: \{s\}\hookrightarrow S$ is the closed point, then $i^\ast Rf_! Rf^!\mathbb F_\ell\in D(\mathbb F_\ell)$ with dual
\[
R\Hom(i^\ast Rf_! Rf^! \mathbb F_\ell,\mathbb F_\ell)\cong R\Gamma(Y,f^\ast i_\ast \mathbb F_\ell)
\]
which sits in degrees $\geq 0$ and is equal to $\mathbb F_\ell$ in degree $0$, as the geometric fibres are connected. Using the section, we get $Rf_! Rf^!\Lambda\cong \Lambda\oplus B$ for some $B$ that sits in cohomological degrees $\leq -1$, and the lemma follows.
\end{proof}

{\it Unfortunately, it is a priori not easy to describe the category ${}^p D^{\geq 0}$.} It is however possible to describe it via hyperbolic localization. This also implies that pullbacks in $S$ are $t$-exact.

\begin{proposition}\label{prop:perversehyperboliclocalization} For any $S'\to S\to \Div^d_{\mathcal Y}$, pullback along
\[
\Hloc_{G,S'/\Div^d_{\mathcal Y}}\to \Hloc_{G,S/\Div^d_{\mathcal Y}}
\]
is $t$-exact for the perverse $t$-structure. Moreover, if $G$ is split, then
\[
\mathrm{CT}_B: D_\et(\Hloc_{G,S/\Div^d_{\mathcal Y}},\Lambda)^{\mathrm{bd}}\to D_\et(\Gr_{T,S/\Div^d_{\mathcal Y}},\Lambda)
\]
satisfies the following exactness property. There is the natural locally constant map $\Gr_{T,\Div^d_{\mathcal Y}}\to X_\ast(T)$ measuring the sum of relative positions, and by pairing with $2\rho$, we get a locally constant map $\mathrm{deg}: \Gr_{T,\Div^d_{\mathcal Y}}\to \mathbb Z$. Then $\mathrm{CT}_B[\mathrm{deg}]$ is $t$-exact for the perverse $t$-structure on the source, and the standard $t$-structure on the right. As $\mathrm{CT}_B[\mathrm{deg}]$ is conservative, this implies in particular that
\[
A\in {}^p D_\et^{\leq 0}(\Hloc_{G,S/\Div^d_{\mathcal Y}},\Lambda)^{\mathrm{bd}}\ (\mathrm{resp.}\ A\in {}^p D_\et^{\geq 0}(\Hloc_{G,S/\Div^d_{\mathcal Y}},\Lambda)^{\mathrm{bd}})
\]
if and only if
\[
\mathrm{CT}_B(A)[\mathrm{deg}]\in D_\et^{\leq 0}(\Gr_{T,S/\Div^d_{\mathcal Y}},\Lambda)\ (\mathrm{resp.}\ \mathrm{CT}_B(A)[\mathrm{deg}]\in D_\et^{\geq 0}(\Gr_{T,S/\Div^d_{\mathcal Y}},\Lambda)).
\]
\end{proposition}

\begin{proof} To prove $t$-exactness of pullbacks, we need to see that pullback commutes with $t$-truncations. By descent, it is enough to check that this holds v-locally on $S$; this allows us to reduce to the case that $G$ is split. It is then enough to prove $t$-exactness of $\mathrm{CT}_B[\mathrm{deg}]$, as by conservativity of $\mathrm{CT}_B[\mathrm{deg}]$ (see Proposition~\ref{prop:conservativity}) this gives the characterization in terms of the $t$-structure in terms of $\mathrm{CT}_B[\mathrm{deg}]$, and the latter characterization is clearly preserved under pullback (as hyperbolic localization commutes with pullback, see Corollary~\ref{cor:hyperboliclocalizationGr}).

We can assume that $S\to \Div^d_{\mathcal Y}$ lifts to $S\to (\Div^1_{\mathcal Y})^d$. We have a stratification of $S$ into finitely many strata $i_a: S_a\hookrightarrow S$, obtained by pulling back the partial diagonals of $(\Div^1_{\mathcal Y})^d$. Accordingly, we get triangles expressing $A$ as a successive extension of $i_{a!} i_a^\ast A$ resp.~$Ri_{a\ast} Ri_a^! A$; if $A\in {}^p D^{\leq 0}$ (resp.~$A\in {}^p D^{\geq 0}$), then also all $i_{a!} i_a^\ast A\in {}^p D^{\leq 0}$ (resp.~$Ri_{a\ast} Ri_a^! A\in {}^p D^{\geq 0})$. This allows us to reduce to the cases of $i_{a!} i_a^\ast A$ and $Ri_{a\ast} Ri_a^! A$. As hyperbolic localization commutes with all functors by Proposition~\ref{prop:hyperboliclocalizationbasechange}, we can then reduce to the case that $S=S_a$ for some $a$. Reducing $d$ if necessary, we can then assume that $S$ maps into the locus of distinct untilts $(\Div^1_{\mathcal Y})^d_{\neq}\subset (\Div^1_{\mathcal Y})^d$. There is then a stratification in terms of open Schubert cells
\[
j_{\mu_\bullet}: \Hloc_{G,S/\Div^d_{\mathcal Y},\mu_\bullet}\hookrightarrow \Hloc_{G,S/\Div^d_{\mathcal Y}}
\]
parametrized by $\mu_\bullet=(\mu_1,\ldots,\mu_d)$, $\mu_i\in X_\ast(T)^+$. Now $A\in {}^p D^{\leq 0}$ if and only if all
\[
j_{\mu_\bullet}^\ast A\in D^{\leq -d_{\mu_\bullet}}
\]
for $d_{\mu_\bullet} = \sum_{i=1}^d \langle 2\rho,\mu_i\rangle$, and dually $A\in {}^p D^{\geq 0}$ if and only if all
\[
Rj_{\mu_\bullet}^! A\in D^{\geq -d_{\mu_\bullet}}.
\]
Using excision triangles, we can then assume that
\[
A=j_{\mu_\bullet !} A_{\mu_\bullet},\ A_{\mu_\bullet}\in D^{\leq -d_{\mu_\bullet}}(\Hloc_{G,S/\Div^d_{\mathcal Y},\mu_\bullet},\Lambda)
\]
resp.
\[
A=Rj_{\mu_\bullet \ast} A_{\mu_\bullet},\ A_{\mu_\bullet}\in D^{\geq -d_{\mu_\bullet}}(\Hloc_{G,S/\Div^d_{\mathcal Y},\mu_\bullet},\Lambda).
\]
Moreover, filtering by cohomology sheaves, we can actually assume that $A_{\mu_\bullet}$ is concentrated in degree $-d_{\mu_\bullet}$. Recall that
\[
[\mu_\bullet]: S\to \Hloc_{G,S/\Div^d_{\mathcal Y},\mu_\bullet}
\]
is a v-cover, and the automorphism group of the stratum is an inverse limit of smooth and connected groups (as follows from Proposition~\ref{prop:openschubertcell} and the K\"unneth formula); this implies that for complexes concentrated in one degree, pullback under $[\mu_\bullet]^\ast$ is fully faithful, cf.~Lemma~\ref{lem:connectedgroupsacttrivially}. We can thus assume that $A_{\mu_\bullet}$ comes via pullback from some $B\in D(S,\Lambda)$ concentrated in cohomological degree $-d_{\mu_\bullet}$. Note that at this point, the desired statement (that $\mathrm{CT}_B(A)[\mathrm{deg}]$ sits in the correct degrees) can be checked after pullback along $\Spa(C,C^+)\to S$, so we can assume $S=\Spa(C,C^+)$ is strictly local, and it is enough to check that $\mathrm{CT}_B(A)[\mathrm{deg}]$ sits in the correct degrees in the fibre over the closed point of $S$. This fibre in turn depends only on the restriction of $A$ to the fibre over the closed point of $S$, by Proposition~\ref{prop:hyperboliclocalizationbasechange}. We can thus assume that $B$ is in fact constant. We can assume $\Lambda=\mathbb Z/n\mathbb Z$ for some $n$ prime to $p$, and then by d\'evissage that $\Lambda=\mathbb F_\ell$ for some $\ell\neq p$. One can then further reduce to the case $B=\mathbb F_\ell[d_{\mu_\bullet}]$. Also, by the K\"unneth formula, we can then reduce to the case $d=1$. Thus, finally
\[
A=j_{\mu!} \mathbb F_\ell[\langle 2\rho,\mu\rangle]
\]
resp.
\[
A=Rj_{\mu\ast} \mathbb F_\ell[\langle 2\rho,\mu\rangle],
\]
and we want to see that $\mathrm{CT}_B(A)[\mathrm{deg}]\in D^{\leq 0}$ (resp.~$\mathrm{CT}_B(A)[\mathrm{deg}]\in D^{\geq 0}$). By Proposition~\ref{prop:hyperboliclocalizationdual}, it suffices to handle the first case. Note that $A$ is now universally locally acyclic, and the claim can be checked in the universal case $S=\Div^1_{\mathcal Y}$. As $\Gr_{T,\Div^1_{\mathcal Y}}\to \Div^1_{\mathcal Y}$ is a disjoint union of $X_\ast(T)$ many copies of $\Div^1_{\mathcal Y}$, and the image is universally locally acyclic, thus locally constant, it is in fact enough to check the result after pullback to the special fibre $S=\Spd \Fq\to \Div^1_{\mathcal Y}$, where
\[
\Gr_{G,\Spd \Fq/\Div^1_{\mathcal Y}}\cong (\Gr_G^{\mathrm{Witt}})^\diamond.
\]

Using \cite[Section 27]{ECoD}, we can now translate all computations to the setting of schemes. Let $\lambda\in X_\ast(T)$ be any element, giving rise to the semi-infinite orbit $S_\lambda\subset \Gr_G^{\mathrm{Witt}}$, i.e.
\[
S_\lambda = \Gr_B^{\mathrm{Witt}}\times_{\Gr_T^{\mathrm{Witt}}} [\lambda].
\]
By Corollary~\ref{cor:dimensionmirkovicvilonen}, the dimension of $S_\lambda\cap \Gr_{G,\mu}^{\mathrm{Witt}}$ is bounded by $\langle \rho,\lambda+\mu\rangle$. The restriction of
\[
\mathrm{CT}_B(j_{\mu!} \mathbb F_\ell[\langle 2\rho,\mu\rangle])
\]
to $[\lambda]\in \Gr_T^{\mathrm{Witt}}$ is given by
\[
R\Gamma_c((S_\lambda\cap \Gr_{G,\mu}^{\mathrm{Witt}})_{\overline{\mathbb F}_q},\mathbb F_\ell)[\langle 2\rho,\mu\rangle]
\]
and thus sits in degrees $\leq 2\langle \rho,\lambda+\mu\rangle - \langle 2\rho,\mu\rangle = \langle 2\rho,\lambda\rangle$, giving the desired bound.
\end{proof}

We note that if $d=1$, $G$ is split, and $S=\Spd k\to \Div^1_{\mathcal Y}$ for $k=\overline{\mathbb F}_q$, then under the full inclusion
\[
\Perv(\Hloc_{G,\Spd k/\Div^1_{\mathcal Y}},\Lambda)\subset D_\et(\Gr_{G,\Spd k/\Div^1_{\mathcal Y}},\Lambda)^{\mathrm{bd}},
\]
the identification $\Gr_{G,\Spd k/\Div^1_{\mathcal Y}}\cong (\Gr_{G,k}^{\mathrm{Witt}})^\diamond$ and the full embedding
\[
D_\et(\Gr_{G,k}^{\mathrm{Witt}},\Lambda)^{\mathrm{bd}}\hookrightarrow  D_\et(\Gr_{G,\Spd k/\Div^1_{\mathcal Y}},\Lambda)^{\mathrm{bd}}
\]
from \cite[Proposition 27.2]{ECoD}, the category
\[
\Perv(\Hloc_{G,\Spd k/\Div^1_{\mathcal Y}},\Lambda)
\]
identifies with the full subcategory
\[
\Perv_{L^+G}(\Gr_{G,k}^{\mathrm{Witt}},\Lambda)\subset D_\et(\Gr_{G,k}^{\mathrm{Witt}},\Lambda)^{\mathrm{bd}}
\]
of $L^+G$-equivariant perverse sheaves on $\Gr_{G,k}^{\mathrm{Witt}}$; this was considered by Zhu \cite{ZhuWitt} and Yu \cite{YuWittSatake}. In particular, this discussion implies the following result that we will need later.

\begin{proposition}\label{prop:satakesemisimple} Assume that $G$ is split, so that for any $\mu\in X_\ast(T)^+$ we have the inclusion
\[
j_\mu: \Hloc_{G,\Div^1_{\mathcal Y},\mu}\hookrightarrow \Hloc_{G,\Div^1_{\mathcal Y}}
\]
of the open Schubert cell, of dimension $d_\mu=\langle 2\rho,\mu\rangle$. Then
\[
{}^p j_{\mu!} \Lambda[d_\mu]={}^p \mathcal H^0(j_{\mu!}\Lambda[d_\mu]),\ {}^p Rj_{\mu\ast} \Lambda[d_\mu]={}^p\mathcal H^0(Rj_{\mu\ast}\Lambda[d_\mu])
\]
are universally locally acyclic, and their image under $\mathrm{CT}_B[\mathrm{deg}]$ is locally finite free over $\Lambda$. Their formation commutes with any base change in $\Lambda$. The natural map
\[
{}^p Rj_{\mu\ast} \Lambda[d_\mu](d_\mu)\to \mathbb D({}^p j_{\mu!} \Lambda[d_\mu])
\]
is an isomorphism.

Moreover, if $\Lambda$ is a $\mathbb Z_\ell$-algebra, then there is some integer $a=a(\mu)$ (independent of $\Lambda$) such that the kernel and cokernel of the map
\[
{}^p \mathcal H^0(j_{\mu!}\Lambda[d_\mu])\to {}^p \mathcal H^0(Rj_{\mu\ast} \Lambda[d_\mu])
\]
are killed by $\ell^a$.
\end{proposition}

We remark that the final statement ultimately makes use of the decomposition theorem (and thus requires the degeneration to the Witt vector affine Grassmannian).

\begin{proof} Consider $A=j_{\mu!} \Lambda[d_\mu]\in {}^p D^{\leq 0}$, which is universally locally acyclic. Then $\mathrm{CT}_B(A)[\mathrm{deg}]$ sits in degrees $\leq 0$, and is universally locally acyclic. Moreover, its degree $0$ part is locally finite free over $\Lambda$. Indeed, this can be computed in terms of the top compactly supported cohomology group of the Mirkovi\'c--Vilonen cycles $S_\lambda\cap \Gr_{G,\mu}^{\mathrm{Witt}}$, which (as for any separated variety) is finite free over $\Lambda$. As $\mathrm{CT}_B(A)[\mathrm{deg}]$ is $t$-exact, this implies that $A'={}^p \mathcal H^0(A)$ has the property that $\mathrm{CT}_B(A')[\mathrm{deg}]$ is locally finite free over $\Lambda$. Applying Verdier duality and using Proposition~\ref{prop:hyperboliclocalizationdual}, we see that $\mathrm{CT}_B(\mathbb D(A))[\mathrm{deg}]\cong\mathbb D(w_0^\ast \mathrm{CT}_B(A)[\mathrm{deg}])$ (where $w_0$ is the longest Weyl group element), which then sits in cohomological degrees $\geq 0$, and is finite free in degree $0$, with degree $0$ parts also under Verdier duality. This shows that the natural map
\[
{}^p Rj_{\mu\ast} \Lambda[d_\mu](d_\mu)\to \mathbb D({}^p j_{\mu!} \Lambda[d_\mu])
\]
is an isomorphism. The proof also shows that the formation commutes with any base change in $\Lambda$.

For the final statement, we can first of all reduce by universal local acyclicity and Corollary~\ref{cor:degenerateSatake} to the same statement on $\Gr_{G,k}^{\mathrm{Witt}}$. By base change, we can assume that $\Lambda=\mathbb Z/\ell^N\mathbb Z$, and we can even formally pass to the inverse limit over $N$, and then invert $\ell$; it is thus enough to show that on the perfectly projective scheme $\Gr_{G,k,\leq\mu}^{\mathrm{Witt}}$, the map
\[
{}^p j_{\mu!} \mathbb Q_\ell[d_\mu]\to {}^p Rj_{\mu\ast} \mathbb Q_\ell[d_\mu]
\]
is an isomorphism. This follows from \cite[Lemma 2.1]{ZhuWitt}, cf.~also \cite[Proposition 1]{GaitsgoryCentral}, \cite{LusztigKatoLusztig}. Let us recall the argument. It is enough to prove injectivity, as then surjectivity follows by Poincar\'e duality (as the two sheaves are Verdier dual, as we have already proved), using that we are working with field coefficients now. Let $j_{\mu!\ast} \mathbb Q_\ell[d_\mu]$ be the image of the displayed map (i.e., the intersection complex of $\Gr_{G,k,\leq\mu}^{\mathrm{Witt}}$). It is enough to see that for $i: \Gr_{G,k,<\mu}^{\mathrm{Witt}}\hookrightarrow \Gr_{G,k,\leq\mu}^{\mathrm{Witt}}$ the complementary closed, that $i^\ast j_{\mu!\ast} \mathbb Q_\ell[d_\mu]$ lies in ${}^p D^{\leq -2}$. Indeed, we have a short exact sequence
\[
0\to i_\ast K\to {}^p j_{\mu!} \mathbb Q_\ell[d_\mu]\to j_{\mu!\ast} \mathbb Q_\ell[d_\mu]\to 0
\]
for some perverse sheaf $K$ on $\Gr_{G,k,<\mu}^{\mathrm{Witt}}$; but this gives a map $i^\ast j_{\mu!\ast} \mathbb Q_\ell[d_\mu]\to K[1]$, so if $i^\ast j_{\mu!\ast} \mathbb Q_\ell[d_\mu]\in {}^p D^{\leq -2}$, then necessarily $K=0$.

To prove that $i^\ast j_{\mu!\ast} \mathbb Q_\ell[d_\mu]\in {}^p D^{\leq -2}$, it suffices to prove that all geometric fibres of $j_{\mu!\ast} \mathbb Q_\ell[d_\mu]$ are concentrated in degrees of the same parity as $d_\mu$; indeed, any other stratum in $\Gr_{G,k,\leq\mu}$ has dimension of the same parity as $d_\mu$, so the trivial bound $i^\ast j_{\mu!\ast} \mathbb Q_\ell[d_\mu]\in {}^p D^{\leq -1}$ gets amplified by one on each stratum for parity reasons. This parity claim about the intersection complex can be checked smooth locally. We have the smooth map $\Fl_{G,k}^{\mathrm{Witt}}\to \Gr_{G,k}^{\mathrm{Witt}}$ from the Witt vector affine flag variety, and choosing the element $w$ in the Iwahori-Weyl group corresponding to the generic stratum on the preimage of the Schubert cell, we get the smooth map $\Fl_{G,k,\leq w}^{\mathrm{Witt}}\to \Gr_{G,k,\leq \mu}^{\mathrm{Witt}}$. It is thus enough to prove the similar claim about the intersection complex of $\Fl_{G,k,\leq w}^{\mathrm{Witt}}$. Choosing a reduced expression $\dot{w}=s_1\cdots s_r\cdot \omega$ as above, we get the Demazure-Bott-Samuelson resolution
\[
\pi_{\dot{w}}: \mathrm{Dem}_{\dot{w}}^{\mathrm{Witt}}\to \Fl_{G,k,\leq w}^{\mathrm{Witt}}.
\]
This has the property that all geometric fibres admit stratifications into affine spaces, cf.~\cite[Section 1.4.2]{ZhuWitt}. In particular, all geometric fibres of $R\pi_{\dot{w}\ast} \mathbb Q_\ell$ sit only in even degrees. On the other hand, by the decomposition theorem, the intersection complex is a direct summand of $R\pi_{\dot{w}\ast}\mathbb Q_\ell[d_\mu]$, giving the claim.
\end{proof}

As a consequence of Proposition~\ref{prop:perversehyperboliclocalization}, we see, perhaps surprisingly, that containment in ${}^p D^{\geq 0}$ can be checked in geometric fibres over $S$. (Note however that we are using a relative perverse $t$-structure.) This gives a complete justification for calling it a relative perverse $t$-structure.

\begin{corollary}\label{cor:perversecoconnective} Let $S\to \Div^d_{\mathcal Y}$ be any small v-stack and let
\[
A\in D_\et(\Hloc_{G,S/\Div^d_{\mathcal Y}},\Lambda)^{\mathrm{bd}}.
\]
Then $A\in {}^p D_\et^{\geq 0}$ if and only if this holds true after pullback to all strictly local $\Spa(C,C^+)\to S$. In particular,
\[
A\in \Perv(\Hloc_{G,S/\Div^d_{\mathcal Y}},\Lambda)
\]
if and only if for all strictly local $\Spa(C,C^+)\to S$, the pullback of $A$ to $\Hloc_{G,\Spa(C,C^+)/\Div^d_{\mathcal Y}}$ is perverse.
\end{corollary}

Also note that over geometric points, we are simply considering the usual perverse $t$-structure corresponding to the stratification in terms of open Schubert cells, and then ${}^p D^{\geq 0}$ admits its usual characterization in terms of $!$-restriction to the open Schubert cells.

\begin{proof} It suffices to check after a cover, as pullback is $t$-exact. This allows us to reduce to the case that $G$ is split. But then it follows from the condition in terms of the hyperbolic localization.
\end{proof}

\subsection{The Satake category} We also get the following characterization. The condition asked here is stronger than perversity.

\begin{proposition}\label{prop:stronglyperverse} Let $S$ be any small v-stack over $\Div^d_{\mathcal Y}$ and assume that $G$ is split. Then
\[
A\in D_\et^\ULA(\Hloc_{G,S/\Div^d_{\mathcal Y}},\Lambda)
\]
has the property that $A$ is a flat perverse sheaf (in the sense that $A\dotimes_\Lambda M$ is perverse for all $\Lambda$-modules $M$) if and only if
\[
R\pi_{T\ast} \mathrm{CT}_B(A)[\mathrm{deg}]\in D_\et(S,\Lambda)
\]
is \'etale locally on $S$ isomorphic to a finite projective $\Lambda$-module in degree $0$.
\end{proposition}

\begin{proof} The functor $R\pi_{T\ast} \mathrm{CT}_B(A)[\mathrm{deg}]$ preserves universally locally acyclic sheaves and hence takes values in sheaves that are locally constant with perfect fibres. By Proposition~\ref{prop:perversehyperboliclocalization} and as any bounded part of $\Gr_{T,\Div^d_{\mathcal Y}}\to \Div^d_{\mathcal Y}$ is finite over the base, the condition $A\in {}^p D^{\leq 0}$ is equivalent to $R\pi_{T\ast} \mathrm{CT}_B(A)[\mathrm{deg}]\in D^{\leq 0}$. The flatness then ensures that this is locally isomorphic to a perfect complex of Tor-amplitude $[0,0]$, i.e.~a finite projective $\Lambda$-module in degree $0$.
\end{proof}

In the following definition, $S\to \Div^d_{\mathcal Y}$ is any small v-stack, and $G$ is general.

\begin{definition}\label{def:satakegeneral} Let
\[
\Sat(\Hloc_{G,S/\Div^d_{\mathcal Y}},\Lambda)\subset D_\et(\Hloc_{G,S/\Div^d_{\mathcal Y}},\Lambda)
\]
be the full subcategory of all objects that are universally locally acyclic and flat perverse.
\end{definition}

This definition has the virtue that it is invariant under switching $\mathrm{sw}^\ast$. Let us give some examples of objects in the Satake category, when $d=1$. Assume for simplicity that $G$ is split. For any $\mu\in X_\ast(T)^+$, we get the open Schubert cell
\[
j_\mu: \Hloc_{G,\Div^1_{\mathcal Y},\mu}\hookrightarrow \Hloc_{G,\Div^1_{\mathcal Y}}
\]
of dimension $d_\mu=\langle 2\rho,\mu\rangle$. The following proposition gives the analogue of highest weight modules in the Satake category.

\begin{proposition}\label{prop:highestweightsatake} The perverse sheaves
\[
{}^p j_{\mu!} \Lambda[d_\mu] = {}^p \mathcal H^0(j_{\mu!} \Lambda[d_\mu]),\ {}^p Rj_{\mu\ast} \Lambda[d_\mu]={}^p \mathcal H^0(Rj_{\mu\ast} \Lambda[d_\mu])
\]
lie in the Satake category $\Sat(\Hloc_{G,\Div^1_{\mathcal Y}},\Lambda)$.
\end{proposition}

\begin{proof} This follows from Proposition~\ref{prop:stronglyperverse} and Proposition~\ref{prop:satakesemisimple}.
\end{proof}

\begin{defprop}\label{prop:fibrefunctorsatake} The functor
\[
R\pi_{G,S\ast}: \Sat(\Hloc_{G,S/\Div^d_{\mathcal Y}},\Lambda)\to D_\et(S,\Lambda)
\]
of pullback to $\Gr_{G,S/\Div^d_{\mathcal Y}}$ and pushforward along $\pi_{G,S}: \Gr_{G,S/\Div^d_{\mathcal Y}}\to S$ takes values in complexes $C\in D_\et(S,\Lambda)$ such that all $\mathcal H^i(C)$ are local systems of finite projective $\Lambda$-modules, and each functor
\[
\mathcal H^i(R\pi_{G,S\ast}): \Sat(\Hloc_{G,S/\Div^d_{\mathcal Y}},\Lambda)\to \LocSys(S,\Lambda).
\]
is exact.

Let
\[
F_{G,S}=\bigoplus_{i\in \mathbb Z} \mathcal H^i(R\pi_{G,S\ast}): \Sat(\Hloc_{G,S/\Div^d_{\mathcal Y}},\Lambda)\to \LocSys(S,\Lambda).
\]
The functor $F_{G,S}$ is exact, faithful, and conservative. Moreover, if $f: A\to B$ is a map in
\[
\Sat(\Hloc_{G,S/\Div^d_{\mathcal Y}},\Lambda)
\]
such that $\ker F_{G,S}(f)$ is a direct summand of $F_{G,S}(A)$, then $f$ admits a kernel in $\Sat(\Hloc_{G,S/\Div^d_{\mathcal Y}},\Lambda)$; similarly for cokernels.
\end{defprop}

The final statement in particular ensures the condition of ``existence of coequalizers of $F_{G,S}$-split parallel pairs'' appearing in the Barr--Beck theorem.

\begin{remark} It is not clear whether there are natural isomorphisms $F_{G,S}(\mathrm{sw}^\ast A)\cong F_{G,S}(A)$, so this fibre functor is (at least a priori) destroying part of the symmetry. What makes this question slightly delicate is that it is asking for extra structure, and ideally one would like to produce this structure in a clean way; it is conceivable that one can reduce to geometric points and then use the affirmative answers we give later under stronger assumptions on $S$.
\end{remark}

\begin{proof} Localize on $S$ to reduce to the case that $G$ is split, and fix a Borel $B\subset G$ with torus $T$. Using the stratification of $\Gr_G$ into the strata $S_\nu$, we get a filtration on $R\pi_{G,S\ast} A$ whose associated graded is given by $\bigoplus_\nu Rp_{\nu!} A|_{S_\nu}$. Restricting to connected components of $\Gr_G$ and for $A$ in the Satake category, these are concentrated in degrees of the same parity, so the corresponding spectral sequence necessarily degenerates. Thus, most of this follows from Proposition~\ref{prop:stronglyperverse} and Proposition~\ref{prop:conservativity}. Faithfulness of $F_{G,S}$ reduces to conservativity and the Barr--Beck type assertion, so it remains to prove the Barr--Beck type assertion. For this, consider the kernel of $f$ in the category of all perverse sheaves on $\Hloc_{G,S/\Div^d_{\mathcal Y}}$. We need to see that this is still universally locally acyclic, and flat perverse. These properties can be checked after applying hyperbolic localization, shift by $\mathrm{deg}$, and pushforward to $S$ (using various $t$-exactness properties), where they follow from the assumption of being a direct summand.
\end{proof}

The Satake category also carries a Verdier duality functor. Again, it is not clear that this functor commutes naturally with $\mathrm{sw}^\ast$ (we will settle it later under stronger assumptions on $S$).

\begin{proposition}\label{prop:satakeverdierdual} The image of the fully faithful functor
\[
\Sat(\Hloc_{G,S/\Div^d_{\mathcal Y}},\Lambda)\to D_\et(\Gr_{G,S/\Div^d_{\mathcal Y}},\Lambda)
\]
is stable under Verdier duality $\mathbb D_{\Gr_{G,S/\Div^d_{\mathcal Y}}/S}$. The induced functor
\[
\mathbb D: \Sat(\Hloc_{G,S/\Div^d_{\mathcal Y}},\Lambda)^{\mathrm{op}}\to \Sat(\Hloc_{G,S/\Div^d_{\mathcal Y}},\Lambda)
\]
is an equivalence, with $\mathbb D^2=\mathrm{id}$. Moreover, it makes the diagram
\[\xymatrix{
\Sat(\Hloc_{G,S/\Div^d_{\mathcal Y}},\Lambda)^{\mathrm{op}}\ar[r]^{\mathbb D}\ar[d]^{F_{G,S}} & \Sat(\Hloc_{G,S/\Div^d_{\mathcal Y}},\Lambda)\ar[d]^{F_{G,S}}\\
\LocSys(S,\Lambda)^{\mathrm{op}}\ar[r]^-{V\mapsto V^\ast} & \LocSys(S,\Lambda)
}\]
commute naturally.
\end{proposition}

\begin{proof} The Verdier dual $\mathbb D(A)\in D_\et(\Gr_{G,S/\Div^d_{\mathcal Y}},\Lambda)$ can actually be defined already in $D_\et(\Hloc_{G,S/\Div^d_{\mathcal Y}},\Lambda)^{\mathrm{bd}}$ by using Verdier duality along bounded subsets of $\Hloc_G\to [\ast/L^+G]$. It follows from Verdier duality that it commutes with the passage to cohomology, i.e.~the functor $F_{G,S}$, and from this one can deduce that it is flat perverse and hence lies in the Satake category. Biduality follows from Corollary~\ref{cor:ULAselfdual}. 
\end{proof}

Moreover, the formation of the Satake category is compatible with constant term functors. We define a locally constant function
\[
\mathrm{deg}_P: \Gr_{M,\Div^d_{\mathcal Y}}\to \mathbb Z
\]
as the composite of the projection to $X_\ast(\overline{M})$ considered before and the map $X_\ast(\overline{M})\to \mathbb Z$ given by pairing with $2\rho_G-2\rho_M$.

\begin{proposition}\label{prop:satakeconstantterm} Let $P\subset G$ be a parabolic with Levi $M$. Let $S\to \Div^d_{\mathcal Y}$ be any small v-stack. Consider the diagram
\[
\Gr_{G,S/\Div^d_{\mathcal Y}}\xleftarrow{q_S} \Gr_{P,S/\Div^d_{\mathcal Y}}\xrightarrow{p_S}\Gr_{M,S/\Div^d_{\mathcal Y}}.
\]
Then the functor $Rp_{S!} q_S^\ast[\mathrm{deg}_P]$ defines a functor
\[
\mathrm{CT}_{P,S}[\mathrm{deg}_P]: \Sat(\Hloc_{G,S/\Div^d_{\mathcal Y}},\Lambda)\to \Sat(\Hloc_{M,S/\Div^d_{\mathcal Y}},\Lambda).
\]
These functors are compatible with composition, i.e.~if $P'\subset P$ is a further parabolic with image $Q\subset M$ and Levi $M'$, then there is a natural equivalence
\[
\mathrm{CT}_{P',S}[\mathrm{deg}_{P'}]\cong \mathrm{CT}_{Q,S}[\mathrm{deg}_{Q}]\circ \mathrm{CT}_{P,S}[\mathrm{deg}_P]: \Sat(\Hloc_{G,S/\Div^d_{\mathcal Y}},\Lambda)\to \Sat(\Hloc_{M',S/\Div^d_{\mathcal Y}},\Lambda)
\]
(and for triple compositions, the obvious diagram commutes).
\end{proposition}

\begin{proof} Let $\lambda: \mathbb G_m\to G$ be a cocharacter such that $P=P_\lambda$. This induces in particular a Levi splitting $M\hookrightarrow P$ as the centralizer of $\lambda$. We can then divide the diagram
\[
\Gr_{G,S/\Div^d_{\mathcal Y}}\xleftarrow{q_S} \Gr_{P,S/\Div^d_{\mathcal Y}}\xrightarrow{p_S}\Gr_{M,S/\Div^d_{\mathcal Y}}.
\]
by $L^+_{\Div^d_{\mathcal Y}} M$ to see that one can refine $Rp_{S!} q_S^\ast$ into a functor
\[
D_\et(\Hloc_{G,S/\Div^d_{\mathcal Y}},\Lambda)\to D_\et(\Hloc_{M,S/\Div^d_{\mathcal Y}},\Lambda)
\]
via first pulling back to $L^+_{\Div^d_{\mathcal Y}} M\backslash \Gr_{G,S/\Div^d_{\mathcal Y}}$. It is clear that these functors are compatible with composition.

We want to see that the image is contained in the Satake category. First, by Proposition~\ref{prop:hyperboliclocalizationULA}, we see that the image is universally locally acyclic. Now the claim follows from Proposition~\ref{prop:stronglyperverse} and the compatibility with composition (used for the Borel $B\subset P$), after passing to an \'etale cover to assume that $G$ is split.
\end{proof}

\section{Convolution}

For any $d$ and small v-stack $S\to \Div^d_{\mathcal Y}$, the category
\[
D_\et(\Hloc_{G,S/\Div^d_{\mathcal Y}},\Lambda)
\]
is naturally a monoidal category. Indeed, with all loop groups taken over $\Div^d_{\mathcal Y}$, there is a convolution morphism
\[
\Hloc_{G,\Div^d_{\mathcal Y}}\times_{\Div^d_{\mathcal Y}} \Hloc_{G,\Div^d_{\mathcal Y}}\xleftarrow{a} L^+ G\backslash L G\times^{L^+ G} L G/L^+ G\xrightarrow{b} L^+ G\backslash L G/L^+ G = \Hloc_{G,\Div^d_{\mathcal Y}}
\]
where the morphism $a$ is an $L^+ G$-torsor, and the right morphism is ind-proper (its fibres are the fibres of $\Gr_{G,\Div^d_{\mathcal Y}}\to \Div^d_{\mathcal Y}$). If one denotes by $a_S$ and $b_S$ the pullbacks along $S\to \Div^d_{\mathcal Y}$, one can then define the convolution product $\star$ on
\[
D_\et(\Hloc_{S/\Div^d_{\mathcal Y}},\Lambda)^{\mathrm{bd}}
\]
via $A_1\star A_2 = Rb_{S\ast} a_S^\ast(A_1\boxtimes A_2)$ for $A_1,A_2\in D_\et(\Hloc_{S/\Div^d_{\mathcal Y}},\Lambda)^{\mathrm{bd}}$. It is easy to see that this is associative by writing out the corresponding convolution diagrams with multiple factors.

In fact, modulo the problem that $[\Div^d_{\mathcal Y}/L^+_{\Div^d_{\mathcal Y}} G]\to [\Div^d_{\mathcal Y}/L_{\Div^d_{\mathcal Y}} G]$ is not representable in locally spatial diamonds (only ind-representable), the category $D_\et(\Hloc_{S/\Div^d_{\mathcal Y}},\Lambda)^{\mathrm{bd}}$ is precisely the category of endomorphisms of $[\Div^d_{\mathcal Y}/L^+_{\Div^d_{\mathcal Y}} G]\times_{\Div^d_{\mathcal Y}} S$ in the $2$-category $\mathcal C_T$ defined in Subsection IV.2.3.3, for $T=[\Div^d_{\mathcal Y}/L_{\Div^d_{\mathcal Y}} G]\times_{\Div^d_{\mathcal Y}} S$. This problem is corrected by passing to bounded sheaves -- one can extend the formalism to the case of maps that are ind-representable in locally spatial diamonds, with closed immersions in the ind-system, using categories of bounded sheaves as morphisms.

Convolution interacts nicely with the classes of sheaves we have previously singled out. In particular, it preserves (flat) perverse sheaves; this observation goes back to Lusztig \cite{LusztigKatoLusztig}.

\begin{proposition}\label{prop:convolutionnice} Let $A_1,A_2\in D_\et(\Hloc_{S/\Div^d_{\mathcal Y}},\Lambda)^{\mathrm{bd}}$.
\begin{altenumerate}
\item[{\rm (i)}] If $A_1$ and $A_2$ are universally locally acyclic, then $A_1\star A_2$ is universally locally acyclic.
\item[{\rm (ii)}] If $A_1$ and $A_2$ lie in ${}^p D^{\leq 0}$, then $A_1\star A_2\in {}^p D^{\leq 0}$.
\item[{\rm (iii)}] If $A_1,A_2\in \Sat(\Hloc_{G,S/\Div^d_{\mathcal Y}},\Lambda)$, then also $A_1\star A_2\in \Sat(\Hloc_{G,S/\Div^d_{\mathcal Y}},\Lambda)$.
\end{altenumerate}
\end{proposition}

\begin{proof} Part (i) follows from Proposition~\ref{prop:properULA} and Proposition~\ref{prop:composeULA}. For part (ii), we first make some reductions. Namely, the claim can be checked if $S=\Spa(C,C^+)$ is strictly local and $G$ split. Moreover, by a d\'evissage one can assume that $A_1$ and $A_2$ are the $!$-extensions of the constant sheaves on open Schubert cells; in particular, these are universally locally acyclic. By the K\"unneth formula one can then reduce to the case $d=1$. In that case, we can pass to the universal case $S=\Div^1_{\mathcal Y}$. Over $(\Div^1_{\mathcal Y})^2$, we can consider the moduli space $\tilde{\Hloc}_{G,(\Div^1_{\mathcal Y})^2}$ of $G$-bundles $\mathcal E_0$, $\mathcal E_1$, $\mathcal E_2$ over $B^+_{(\Div^1_{\mathcal Y})^2}$ together with isomorphisms between $\mathcal E_0$ and $\mathcal E_1$ over $B^+_{(\Div^1_{\mathcal Y})^2}[\tfrac 1{\mathcal I_1}]$ and between $\mathcal E_1$ and $\mathcal E_2$ over $B^+_{(\Div^1_{\mathcal Y})^2}[\tfrac 1{\mathcal I_2}]$, where $\mathcal I_1,\mathcal I_2\subset \mathcal O_{\mathcal Y_S}$ are the ideal sheaves parametrizing the two Cartier divisors. Away from the diagonal, this is isomorphic to $\Hloc_{G,(\Div^1_{\mathcal Y})^2_{\neq}/\Div^2_{\mathcal Y}}$, while over the diagonal it is isomorphic to
\[
L^+_{\Div^1_{\mathcal Y}} G\backslash L_{\Div^1_{\mathcal Y}} G\times^{L^+_{\Div^1_{\mathcal Y}} G} L_{\Div^1_{\mathcal Y}} G/L^+_{\Div^1_{\mathcal Y}} G.
\]
There are two natural projections
\[
p_1,p_2: \tilde{\Hloc}_{G,(\Div^1_{\mathcal Y})^2}\to \Hloc_{G,\Div^1_{\mathcal Y}}
\]
keeping track of $\mathcal E_0$ and $\mathcal E_1$ resp.~$\mathcal E_1$ and $\mathcal E_2$, and a projection
\[
m: \tilde{\Hloc}_{G,(\Div^1_{\mathcal Y})^2}\to \Hloc_{G,(\Div^1_{\mathcal Y})^2/\Div^2_{\mathcal Y}}
\]
keeping track of $\mathcal E_0$ and $\mathcal E_2$. One can the form $B=Rm_\ast(p_1^\ast A_1\dotimes_\Lambda p_2^\ast A_2)$. Recall that we reduced to the case that $A_1$ and $A_2$ are moreover universally locally acyclic. By Proposition~\ref{prop:properULA} and Proposition~\ref{prop:composeULA}, one sees that also
\[
B\in D_\et^\ULA(\Hloc_{G,(\Div^1_{\mathcal Y})^2/\Div^2_{\mathcal Y}},\Lambda).
\]
Away from the diagonal, this is simply the exterior tensor product of $A_1$ and $A_2$ and in particular lies in ${}^p D^{\leq 0}$. Looking at $\mathrm{CT}_B(B)[\mathrm{deg}]$, we get a universally locally acyclic sheaf on (a bounded subset of) $\Gr_{T,(\Div^1_{\mathcal Y})^2/\Div^2_{\mathcal Y}}$ whose restriction away from the diagonal lies in degrees $\leq 0$. This implies that the whole sheaf lies in degrees $\leq 0$: As any bounded subset of $\Gr_{T,(\Div^1_{\mathcal Y})^2/\Div^2_{\mathcal Y}}$ is finite over $(\Div^1_{\mathcal Y})^2$, it suffices to check this for the pushforward to $(\Div^1_{\mathcal Y})^2$. But this pushforward is locally constant with perfect fibres, and the complement of the diagonal is dense.

Thus, using Proposition~\ref{prop:perversehyperboliclocalization}, the restriction of $B$ to the diagonal lies in ${}^p D^{\leq 0}$. But this restriction is precisely $A_1\star A_2$, giving the desired result. Finally, part (iii) easily follows from (i), (ii), and the observation that convolution commutes with Verdier duality.
\end{proof}

\subsection{Dualizability} Next, we observe that all objects are dualizable.

\begin{proposition}\label{prop:satakedualizable} All objects of the monoidal category $\Sat(\Hloc_{G,S/\Div^d_{\mathcal Y}},\Lambda)$ are (left and right) dualizable. The right dual of $A\in \Sat(\Hloc_{G,S/\Div^d_{\mathcal Y}},\Lambda)$ is given by $\mathrm{sw}^\ast \mathbb D(A)$ where $\mathrm{sw}: \Hloc_{G,\Div^d_{\mathcal Y}}\cong \Hloc_{G,\Div^d_{\mathcal Y}}$ is the switching isomorphism (induced by inversion on $L_{\Div^d_{\mathcal Y}} G$).
\end{proposition}

\begin{remark} In the classical setting, this is asserted without indication of proof in \cite[end of Section 11]{MirkovicVilonen}.
\end{remark}

\begin{proof} All objects of $D^\ULA(\Hloc_{G,S/\Div^d_{\mathcal Y}},\Lambda)$ are left dualizable, with right dual given by $\mathrm{sw}^\ast \mathbb D(A)$: This follows from Proposition~\ref{prop:generalULAdualizable} (modulo the technical nuisance that everything is only ind-representable here; everything adapts to that setting). Here, $\mathrm{sw}$ simply arises by swapping source and target. Also note that the condition of being universally locally acyclic is invariant under $\mathrm{sw}^\ast$ by Proposition~\ref{prop:ULAsymmetry}, so also using Proposition~\ref{prop:satakeverdierdual}, the functor $\mathrm{sw}^\ast \mathbb D(A)$ preserves the Satake category.
\end{proof}

\begin{remark} Again, we stress that all results above also hold if $G$ is reductive over $E$, and we replace $\Div^d_{\mathcal Y}$ with $\Div^d_Y$ or $\Div^d_X$. Indeed, the case of $\Div^d_Y$ follows from the case of $\Div^d_{\mathcal Y}$ as it is an open subset, at least if $G$ admits a reductive model over $\mathcal O_E$. In general, this happens \'etale locally, making it possible to reduce to this case. Then the case of $\Div^d_X$ follows from the case of $\Div^d_Y$ as any map $S\to \Div^d_X$ can locally be lifted to a map $S\to \Div^d_Y$ in such a way that the corresponding pullbacks of the local Hecke stacks are isomorphic.
\end{remark}

\section{Fusion}

Now let $G$ be a reductive group over $E$. From now on, we fix the base field $k=\overline{\mathbb F}_q$ and work on the category $\Perf_k$. For brevity, we define for any finite set $I$ with $d=|I|$ the local Hecke stack
\[
\Hloc_G^I = \Hloc_{G,\Div^d_X}\times_{\Div^d_X} (\Div^1_X)^I
\]
and correspondingly
\[
\Gr_G^I = \Gr_{G,\Div^d_X}\times_{\Div^d_X} (\Div^1_X)^I.
\]

\begin{definition}\label{def:satakecategory} For any finite set $I$, the Satake category
\[
\Sat_G^I(\Lambda)
\]
is the category $\Sat(\Hloc_G^I,\Lambda)$ of all
\[
A\in D_\et(\Hloc_G^I,\Lambda)
\]
that are universally locally acylic and flat perverse over $(\Div^1_X)^I$.
\end{definition}

By Proposition~\ref{prop:fibrefunctorsatake}, we get a (not yet monoidal) fibre functor
\[
F^I: \Sat_G^I(\Lambda)\to \LocSys((\Div^1_X)^I,\Lambda).
\]
The target category $\LocSys((\Div^1_X)^I,\Lambda)$ is in fact very explicit.

\begin{proposition}\label{prop:locsysdiv1} The category $\LocSys((\Div^1_X)^I,\Lambda)$ is naturally equivalent to the category $\Rep_{W_E^I}(\Lambda)$ of continuous representations of $W_E^I$ on finite projective $\Lambda$-modules.
\end{proposition}

\begin{proof} This is a consequence of Proposition~\ref{prop:drinfeldlemmalocallyconstant}.
\end{proof}

For any map $f: I\to J$ of finite sets, there is a natural monoidal functor $\Sat_G^I(\Lambda)\to \Sat_G^J(\Lambda)$. Indeed, there is a natural closed immersion $\Gr_G^I\times_{(\Div^1_X)^I} (\Div^1_X)^J\hookrightarrow \Gr_G^J$, so pull-push along
\[
\Gr_G^I\leftarrow \Gr_G^I\times_{(\Div^1_X)^I} (\Div^1_X)^J\to \Gr_G^J
\]
defines the desired functor (noting that this preserves the required equivariance condition).\footnote{We thank Tony Feng for pointing out that on Hecke stacks the map is not a closed immersion.} It is easy to see that this functor is compatible with composition of maps of finite sets. Moreover, the functors $\Sat_G^I(\Lambda)\to \Sat_G^J(\Lambda)$ make the diagram
\[\xymatrix{
\Sat_G^I(\Lambda)\ar[r]\ar[d]^{F^I} & \Sat_G^J(\Lambda)\ar[d]^{F^J}\\
\Rep_{W_E^I}(\Lambda)\ar[r] & \Rep_{W_E^J}(\Lambda)
}\]
commute naturally, where the lower functor is pullback under $W_E^J\to W_E^I$.

Actually, the functor $I\mapsto \Sat_G^I(\Lambda)$ has further functoriality, given by the fusion product. Namely, for finite sets $I_1,\ldots,I_k$ with disjoint union $I=I_1\sqcup\ldots\sqcup I_k$, there is a natural monoidal functor
\[
\Sat_G^{I_1}(\Lambda)\times \ldots\times \Sat_G^{I_k}(\Lambda)\to \Sat_G^I(\Lambda),
\]
functorial in $I_1,\ldots,I_k$ and compatible with composition. To construct this, let
\[
j: (\Div^1_X)^{I;I_1,\ldots,I_k}\subset (\Div^1_X)^I
\]
be the open subset where $x_i\neq x_{i'}$ whenever $i,i'\in I=I_1\sqcup\ldots\sqcup I_k$ lie in different $I_j$'s, and let
\[
\Sat_G^{I;I_1,\ldots,I_k}(\Lambda)\subset D_\et^\ULA(\Hloc_G^I\times_{(\Div^1_X)^I} (\Div^1_X)^{I;I_1,\ldots,I_k},\Lambda)
\]
be defined similarly as $\Sat_G^I(\Lambda)$.

\begin{proposition}\label{prop:restrictfullyfaithful} The restriction functor
\[
j^\ast: \Sat_G^I(\Lambda)\to \Sat_G^{I;I_1,\ldots,I_k}(\Lambda)
\]
is fully faithful. Similarly,
\[
j^\ast: \LocSys((\Div^1_X)^I,\Lambda)\to \LocSys((\Div^1_X)^{I;I_1,\ldots,I_k},\Lambda)
\]
is fully faithful.
\end{proposition}

\begin{proof} For the first part, it suffices to prove that for all $A\in \Sat_G^I(\Lambda)$, the natural map
\[
A\to {}^p\mathcal H^0(Rj_\ast j^\ast A)
\]
is an isomorphism. Let $i: Z\hookrightarrow (\Div^1_X)^I$ be the complementary closed. It suffices to see that $i_\ast i^! A\in {}^p D^{\geq 2}$. Working locally to reduce to the case $G$ split, and applying the $t$-exact hyperbolic localization functor $R\pi_{T\ast}\mathrm{CT}_B[\mathrm{deg}]$, taking values in local systems of finite projective $\Lambda$-modules on $\Sat_G^I(\Lambda)$, this follows from the observation that $i_\ast i^! \Lambda\in D^{\geq 2}$, which follows from the observation that $Z$ admits a stratification (by partial diagonals) with smooth strata of $\ell$-codimension $\geq 1$ inside the smooth $(\Div^1_X)^I$.

This final argument in fact proves directly the second part.
\end{proof}

On the other hand, over $(\Div^1_X)^{I;I_1,\ldots,I_k}$, one has
\[
\Hloc_G^I\times_{(\Div^1_X)^I} (\Div^1_X)^{I;I_1,\ldots,I_k}\cong \prod_{j=1}^k \Hloc_G^{I_j}\times_{\prod_j (\Div^1_X)^{I_j}} (\Div^1_X)^{I;I_1,\ldots,I_k},
\]
so there is a natural monoidal functor
\[
\mathrm{Sat}_G^{I_1}(\Lambda)\times \ldots\times \mathrm{Sat}_G^{I_k}(\Lambda)\to \Sat_G^{I;I_1,\ldots,I_k}(\Lambda)
\]
given by exterior product. Actually, recall that when forming symmetric monoidal tensor products, there are implicit sign rules when commuting factors. We change these here by hand. Namely, note that each
\[
\Hloc_G^I = (\Hloc_G^I)^{\mathrm{even}}\sqcup (\Hloc_G^I)^{\mathrm{odd}}
\]
decomposes into open and closed subsets given by the even and the odd part; the even part contains those Schubert varieties for which $d_{\mu_\bullet}=\sum_i \langle 2\rho,\mu_i\rangle$ is even, while the odd part contains those for which $d_{\mu_\bullet}$ is odd. Note that the dominance order can only nontrivially compare elements with the same parity, so these are really open and closed subsets. Also note that it follows from Proposition~\ref{prop:perversehyperboliclocalization} that for sheaves concentrated on the even (resp.~odd) part, the functor $F^I$ is concentrated in even (resp.~odd) degrees. Now we impose that when forming the above exterior product, we introduce a minus sign whenever we commute two sheaves concentrated on the odd parts. A different way to say it is that there is a natural commutative diagram
\[\xymatrix{
\mathrm{Sat}_G^{I_1}(\Lambda)\times \ldots\times \mathrm{Sat}_G^{I_k}(\Lambda)\ar[r]\ar[d]^{(F^{I_1},\ldots,F^{I_k})} & \Sat_G^{I;I_1,\ldots,I_k}(\Lambda)\ar[d]^{F^{I;I_1,\ldots,I_k}}\\
\LocSys((\Div^1_X)^{I_1},\Lambda)\times\ldots\times \LocSys((\Div^1_X)^{I_k},\Lambda)\ar[r]^-{\boxtimes} & \LocSys((\Div^1_X)^{I;I_1,\ldots,I_k},\Lambda)
}\]
functorial in $I_1,\ldots,I_k$, and under permutations of the sets $I_1,\ldots,I_k$. Indeed, note that the functors $F^I$ invoke a shift by $\mathrm{deg}$, which exactly introduces this sign rule. This in fact pins down this choice of signs by faithfulness of the functors.

\begin{defprop}\label{prop:fusionproduct} The image of
\[
\mathrm{Sat}_G^{I_1}(\Lambda)\times \ldots\times \mathrm{Sat}_G^{I_k}(\Lambda)\to \Sat_G^{I;I_1,\ldots,I_k}(\Lambda)
\]
lands in $\Sat_G^I(\Lambda)\subset \Sat_G^{I;I_1,\ldots,I_k}(\Lambda)$, defining the fusion product
\[
\ast: \mathrm{Sat}_G^{I_1}(\Lambda)\times \ldots\times \mathrm{Sat}_G^{I_k}(\Lambda)\to\Sat_G^I(\Lambda),
\]
a functor of monoidal categories, functorial in $I_1,\ldots,I_k$. It makes the diagram
\[\xymatrix{
\mathrm{Sat}_G^{I_1}(\Lambda)\times \ldots\times \mathrm{Sat}_G^{I_k}(\Lambda)\ar[r]^\ast\ar[d]^{(F^{I_1},\ldots,F^{I_k})} & \Sat_G^I(\Lambda)\ar[d]^{F^I}\\
\LocSys((\Div^1_X)^{I_1},\Lambda)\times\ldots\times \LocSys((\Div^1_X)^{I_k},\Lambda)\ar[r]^-{\boxtimes} & \LocSys((\Div^1_X)^I,\Lambda)
}\]
commute functorially in $I_1,\ldots,I_k$ and permutations of $I_1,\ldots,I_k$.
\end{defprop}

\begin{proof} We can define a convolution local Hecke stack
\[
\Hloc_G^{I;I_1,\ldots,I_k}\to (\Div^1_X)^I
\]
as follows. It parametrizes $G$-bundles $\mathcal E_0,\ldots,\mathcal E_k$ over $B^+_{(\Div^1_X)^I}$ together with isomorphisms of $\mathcal E_{j-1}$ and $\mathcal E_j$ after inverting $\mathcal I_i$ for all $i\in I_j$, for $j=1,\ldots,k$. Here $\mathcal I_i\subset \mathcal O_{X_S}$ is the ideal defining the $i$-th Cartier divisor. There are natural projections
\[
p_j: \Hloc_G^{I;I_1,\ldots,I_k}\to \Hloc_G^{I_j},\ j=1,\ldots,k
\]
remembering $\mathcal E_{j-1}$ and $\mathcal E_j$, and
\[
m: \Hloc_G^{I;I_1,\ldots,I_k}\to \Hloc_G^I
\]
remembering $\mathcal E_0$ and $\mathcal E_k$. Given $A_j\in \mathrm{Sat}_G^{I_j}(\Lambda)$, one can then define
\[
B=Rm_\ast(p_1^\ast A_1\dotimes_\Lambda\ldots\dotimes_\Lambda p_k^\ast A_k)\in D_\et(\Hloc_G^I,\Lambda).
\]
This is still universally locally acyclic, by Proposition~\ref{prop:properULA} and Proposition~\ref{prop:composeULA}. After pullback to $(\Div^1_X)^{I;I_1,\ldots,I_k}$, the map $m$ is an isomorphism, and we simply get the exterior product of all $A_i$. Moreover, working locally to reduce to the case $G$ is split, we see that $R\pi_{T\ast} \mathrm{CT}_B(B)[\mathrm{deg}]$ is a local system of finite projective $\Lambda$-modules, as it is locally constant with perfect fibres, and over the dense open subset $(\Div^1_X)^{I;I_1,\ldots,I_k}$, the perfect complex is a finite projective module in degree $0$. This means that $B\in \Sat_G^I(\Lambda)$, as desired.
\end{proof}

In particular, for any finite set $I$, this structure makes $\Sat_G^I(\Lambda)$ into an $E_\infty$-monoid object in monoidal categories, functorially in $I$, by using the composite
\[
\Sat_G^I(\Lambda)\times\ldots\times \Sat_G^I(\Lambda)\to \Sat_G^{I\sqcup\ldots\sqcup I}(\Lambda)\to \Sat_G^I(\Lambda),
\]
using the functor corresponding to the natural map $I\sqcup \ldots\sqcup I\to I$. Recall that $E_\infty$-monoid structures on monoidal categories are the same as symmetric monoidal category structures refining the given monoidal category structure. Thus, each $\Sat_G^I(\Lambda)$ has become naturally a symmetric monoidal category with the fusion product, refining the monoidal convolution product; and everything is functorial in $I$. Moreover, by the final part of Definition/Proposition~\ref{prop:fusionproduct}, the functor
\[
F^I: \Sat_G^I(\Lambda)\to \LocSys((\Div^1_X)^I,\Lambda)\cong \Rep_{W_E^I}(\Lambda)
\]
is a symmetric monoidal functor, functorially in $I$.

A consequence of these symmetric monoidal structures are the following natural isomorphisms.

\begin{corollary}\label{cor:niceisomorphisms} For $A\in \Sat_G^I(\Lambda)$, there are natural isomorphisms
\[
F^I(\mathrm{sw}^\ast A)\cong F^I(A),\ \mathbb D(\mathrm{sw}^\ast A)\cong \mathrm{sw}^\ast \mathbb D(A).
\]
Moreover, $\mathbb D$ is naturally a symmetric monoidal functor, and $\mathbb D\circ F^I\cong (F^I)^\ast$ as symmetric monoidal functors.
\end{corollary}

\begin{proof} By Proposition~\ref{prop:satakedualizable}, all $A\in \Sat_G^I(\Lambda)$ are dualizable, with dual $\mathrm{sw}^\ast \mathbb D(A)$. In a symmetric monoidal category, this means that there are natural isomorphisms
\[
\mathrm{sw}^\ast \mathbb D \mathrm{sw}^\ast \mathbb D(A)\cong A.
\]
As both $\mathrm{sw}^\ast$ and $\mathbb D$ are self-inverse, this amounts to the commutation of $\mathbb D$ and $\mathrm{sw}^\ast$.

Also, as $F^I$ is symmetric monoidal, it follows that $F^I(\mathrm{sw}^\ast \mathbb D(A))$ and $F^I(A)$ are naturally dual. But by Proposition~\ref{prop:satakeverdierdual}, the dual of $F^I(A)$ is also $F^I(\mathbb D(A))$. Replacing $\mathbb D(A)$ by $A$, we find a natural isomorphism $F^I(\mathrm{sw}^\ast A)\cong F^I(A)$.

Finally, it is easy to see that the whole construction of the fusion product is compatible with Verdier duality, making Verdier duality a symmetric monoidal functor, compatibly with the fibre functor.
\end{proof}

Moreover, the constant term functors are compatible with the fusion product. More precisely, given a parabolic $P\subset G$ with Levi $M$, we have the constant term functors
\[
\mathrm{CT}_P^I[\mathrm{deg}_P]: \Sat_G^I(\Lambda)\to \Sat_M^I(\Lambda).
\]

\begin{proposition}\label{prop:constanttermfusion} For any finite set $I$ decomposed into finite sets $I=I_1\sqcup\ldots\sqcup I_k$, the diagram
\[\xymatrix{
\mathrm{Sat}_G^{I_1}(\Lambda)\times \ldots\times \mathrm{Sat}_G^{I_k}(\Lambda)\ar[rr]\ar[d]^{(\mathrm{CT}_P^{I_1}[\mathrm{deg}_P],\ldots,\mathrm{CT}_P^{I_k}[\mathrm{deg}_P])} && \Sat_G^I(\Lambda)\ar[d]^{\mathrm{CT}_P^I[\mathrm{deg}_P]}\\
\mathrm{Sat}_M^{I_1}(\Lambda)\times \ldots\times \mathrm{Sat}_M^{I_k}(\Lambda)\ar[rr] && \Sat_M^I(\Lambda)
}\]
commutes functorially in $I$ and permutations of $I_1,\ldots,I_k$.
\end{proposition}

\begin{proof} After passing to the open subset $(\Div^1_X)^{I;I_1,\ldots,I_k}$, this follows from the K\"unneth formula, so Proposition~\ref{prop:restrictfullyfaithful} gives the result.
\end{proof}

In particular, the functor
\[
\mathrm{CT}_P^I[\mathrm{deg}_P]: \Sat_G^I(\Lambda)\to \Sat_M^I(\Lambda)
\]
is naturally symmetric monoidal with respect to the fusion product. Moreover, everything is compatible with composition, for another parabolic $P'\subset P$.

\section{Tannakian reconstruction}

Our next goal is to construct a group scheme whose category of representations recovers the symmetric monoidal category $\Sat_G^I(\Lambda)$. More precisely, we want to use some relative Tannaka duality over $\Rep_{W_E^I}(\Lambda)$. To achieve this, we need the following proposition. Given any finite and Galois-stable subsets $W_i\subset X_\ast(T)^+$, $i\in I$, closed under the dominance order, we have a quasicompact closed substack
\[
\Hloc_{G,(W_i)_i}^I\subset \Hloc_G^I
\]
and we get a corresponding full subcategory
\[
\Sat_{G,(W_i)_i}^I(\Lambda)\subset \Sat_G^I(\Lambda).
\]

\begin{proposition}\label{prop:leftadjointsatake} The functor
\[
F^I: \Sat_{G,(W_i)_i}^I(\Lambda)\to \Rep_{W_E^I}(\Lambda)
\]
admits a left adjoint $L^I_{(W_i)_i}$, satisfying the following properties.
\begin{altenumerate}
\item[{\rm (i)}] There is a natural isomorphism
\[
L^I_{(W_i)_i}(V)\cong L^I_{(W_i)_i}(1)\otimes V,\ V\in \Rep_{W_E^I}(\Lambda),
\]
where $1\in \Rep_{W_E^I}(\Lambda)$ is the tensor unit, and we use that $\Sat_G^I(\Lambda)$ is tensored over $\Rep_{W_E^I}(\Lambda)$.
\item[{\rm (ii)}] There is a natural isomorphism
\[
L^I_{(W_i)_i}(1)\cong \bigast_{i\in I} L_{W_i}^{\{i\}}(1)
\]
as the fusion of $L_{W_i}^{\{i\}}(1)\in \Sat_{G,W_i}^{\{i\}}(\Lambda)$.
\item[{\rm (iii)}] If $I=\{i\}$ has one element and $W=W_i$, then the left adjoint is the restriction of the left adjoint to
\[
F' = \bigoplus_m \mathcal H^m(R\pi_{G\ast}): \Perv(\Hloc_{G,W}^{\{i\}},\Lambda)\to \Mod_{W_E}(\Lambda).
\]
\end{altenumerate}
\end{proposition}

\begin{proof} It is enough to find the value $L^I_{(W_i)_i}(1)$ satisfying (ii) and (iii). Indeed, then the formula in (i) defines the left adjoint in general. Assume now that (iii) holds, and let us denote $P_{W_i}=L_{W_i}^{\{i\}}(1)$. Then for part (ii) we first observe that
\[
F^{\prime I}= \bigoplus_m \mathcal H^m(R\pi_{G\ast}): \Perv(\Hloc_{G,(W_i)_i}^I,\Lambda)\to \mathrm{Shv}_\et((\Div^1_X)^I,\Lambda)
\]
admits a left adjoint, and this left adjoint evaluated on the unit, $P_{(W_i)_i}$, is generically on $(\Div^1_X)^I$ given by an exterior tensor product of the corresponding left adjoints for $I$ being a singleton. Indeed, note that there is a natural map
\[
P_{(W_i)_i}\to \bigast_{i\in I} P_{W_i}
\]
adjoint to the section of
\[
F^{\prime I}(\bigast_{i\in I} P_{W_i})\cong \bigboxtimes_{i\in I} F'(P_{W_i})
\]
given by the exterior tensor product of the classes given by (iii). To check that this is an isomorphism generically, we can by \'etale descent reduce to the case that $G$ is split. In that case, one can make the left adjoint explicit in terms of the left adjoint to hyperbolic localization. Writing hyperbolic localization as a composite of $!$-pullback and $\ast$-pushforward, this left adjoint is then given in terms of $\ast$-pullback and $!$-pushforward, and the perverse ${}^p \mathcal H^0$. As generically, everything decomposes geometrically into a product, it follows from the K\"unneth formula that the left adjoint commutes with exterior products. But as any $B\in \Sat_{G,(W_i)_i}^I(\Lambda)$ is equal to ${}^p \mathcal H^0(Rj_\ast j^\ast B)$ as in the discussion of the fusion product, we see that
\[\begin{aligned}
F^{\prime I}(B)&\cong \Hom(P_{(W_i)_i},B)\cong \Hom(P_{(W_i)_i},{}^p \mathcal H^0(Rj_\ast j^\ast B))\\
&\cong \Hom(P_{(W_i)_i},Rj_\ast j^\ast B)\\
&\cong \Hom(j^\ast P_{(W_i)_i},j^\ast B)\\
&\cong \Hom(j^\ast \bigast_{i\in I} P_{W_i},j^\ast B)\\
&\cong \Hom(\bigast_{i\in I} P_{W_i},{}^p\mathcal H^0(Rj_\ast j^\ast B))\cong \Hom(\bigast_{i\in I} P_{W_i},B).
\end{aligned}\]

It remains to prove part (iii). We can assume that $\Lambda=\mathbb Z/\ell^c\mathbb Z$, using base change. Note first that
\[
F' = \bigoplus_m \mathcal H^m(R\pi_{G\ast}): \Perv(\Hloc_{G,W}^{\{i\}},\Lambda)\to \Mod_{W_E}(\Lambda)
\]
admits a left adjoint $L_W'$, by the adjoint functor theorem. We need to see that when evaluated at the unit, $P_W:=L_W'(1)$ is universally locally acyclic, and flat perverse. By the characterization of these properties, it suffices to show that $F'(P_W)\in \Mod_{W_E}(\Lambda)$ is a representation on a finite projective $\Lambda$-module. This does not depend on the $W_E$-action, so we can check these things after pullback along $\Spd C\to \Div^1_X$, where $C$ is a completed algebraic closure of $E$. In particular, we can assume that $G$ is split. For any open Schubert cell
\[
j_\mu: \Hloc_{G,\Spd C,\mu}\hookrightarrow \Hloc_{G,\Spd C,W}
\]
for $\mu\in W$, of dimension $d_\mu=\langle 2\rho,\mu\rangle$, we can compute
\[
\Hom(P_W,{}^p Rj_{\mu\ast} \Lambda[d_\mu])=F'({}^p Rj_{\mu\ast}\Lambda[d_\mu]).
\]
By Proposition~\ref{prop:highestweightsatake}, this is a finite free $\Lambda$-module. Using adjunction, we thus see that
\[
\Hom({}^p j_\mu^\ast P_W,\Lambda[d_\mu])
\]
is a finite free $\Lambda$-module for all $\mu\in W$. Now ${}^p j_\mu^\ast P_W$ is concentrated on an open Schubert cell $\Hloc_{G,\Spd C,\mu}$, which is covered by $\Spd C$, and concentrated in degree $-d_\mu$. It is thus given by the constant sheaf $M[d_\mu]$ for some $\Lambda$-module $M$, and we know that $\Hom(M,\Lambda)$ is finite free over $\Lambda$. As we reduced to $\Lambda=\mathbb Z/\ell^m\mathbb Z$, this implies that $M$ is free.

Now argue by induction on $W$, and take a maximal element $\mu\in W$; let $\overline{W} = W\setminus \{\mu\}$. We get an exact sequence
\[
0\to K\to {}^p j_{\mu!} j_\mu^\ast P_W\to P_W\to Q\to 0
\]
in $\Perv(\Hloc_{G,\Spd C},\Lambda)$ supported on $W$. In fact, we necessarily have $Q=P_{\overline{W}}$ (as they represent the same functor), for which we know by induction that $F'(Q)$ is a finite free $\Lambda$-module. We claim that $K=0$. As $K$ lies in the kernel of ${}^p j_{\mu!} j_\mu^\ast P_W\to {}^p Rj_{\mu\ast} j_\mu^\ast P_W$, it follows from Proposition~\ref{prop:satakesemisimple} that $\ell^a K=0$ for some $a$ independent of $\Lambda$. Using functoriality of the construction for $\Lambda'=\mathbb Z/\ell^{a+c}\mathbb Z\to \Lambda=\mathbb Z/\ell^c \mathbb Z$ and that ${}^p j_{\mu!} j_\mu^\ast P_W$ lies in the Satake category (so in particular it is flat over $\Lambda$), we see the image of $K'$ in $K$ is equal to $0$. On the other hand, as all constructions are compatible with base change, the map $K'\to K$ had to be surjective. It follows that $K=0$, as desired. (Alternatively, we could have reduced to $\mathbb Z_\ell$-coefficients, in which case ${}^p j_{\mu!} j_\mu^\ast P_W\to {}^p Rj_{\mu\ast} j_\mu^\ast P_W$ is injective (as the kernel is both $\ell$-torsion free and killed by $\ell^a$), implying $K=0$ directly.)
\end{proof}

Now we use the following general Tannakian reconstruction result. This is essentially an axiomatization of \cite[Proposition 11.1]{MirkovicVilonen}. Recall that a symmetric monoidal category is rigid if all of its objects are dualizable.

\begin{proposition}\label{prop:tannakareconstruction} Let $\mathcal A$ be a rigid symmetric monoidal category, and let $\mathcal C$ be a symmetric monoidal category with a tensor action of $\mathcal A$. Moreover, let
\[
F: \mathcal C\to \mathcal A
\]
be a symmetric monoidal $\mathcal A$-linear conservative functor, such that $\mathcal C$ admits and $F$ reflects coequalizers of $F$-split parallel pairs. Assume that $\mathcal C$ can be written as a filtered union of full subcategories $\mathcal C_i$, stable under coequalizers of $F$-split parallel pairs and the $\mathcal A$-action, such that $F|_{\mathcal C_i}$ is representable by some $X_i\in \mathcal C$.

Then
\[
\mathcal H=\varinjlim_i F(X_i)^\vee\in \Ind(\mathcal A)
\]
admits a natural structure as a bialgebra (with commutative multiplication and associative comultiplication), and $\mathcal C$ is naturally equivalent to the symmetric monoidal category of representations of $\mathcal H$ in $\mathcal A$. If $\mathcal C$ is rigid, then $\mathcal H$ admits an inverse, i.e.~is a Hopf algebra.
\end{proposition}

Here, the symmetric monoidal category of representations of $\mathcal H$ is the category of comodules over $\mathcal H$ as a coalgebra, endowed with the symmetric monoidal structure coming from the commutative multiplication on $\mathcal H$.

\begin{proof} Consider $F_i=F|_{\mathcal C_i}: \mathcal C_i\to \mathcal A$. This admits the left adjoint $A\mapsto A\otimes X_i$, as
\[
\Hom_{\mathcal C_i}(A\otimes X_i,Y)\cong \Hom_{\mathcal C_i}(X_i,A^\vee\otimes Y)\cong F(A^\vee\otimes Y)\cong A^\vee\otimes F(Y)\cong \Hom(A,F(Y)).
\]
By the Barr--Beck monadicity theorem, it follows that $\mathcal C_i$ is equivalent to the category of modules over the monad
\[
A\mapsto F(A\otimes X_i)\cong A\otimes F(X_i).
\]
Note that the monad structure here is equivalently turning $F(X_i)$ into an associative algebra $A_i\in \mathcal A$, and its category of modules is the category of modules over $A_i$. Passing to duals, we note that $F(X_i)^\vee$ is a coalgebra, and its category of comodules is equivalent to the category of modules over $A_i$, i.e.~to $\mathcal C_i$. Now we can take a colimit over $i$ and see that
\[
\mathcal H=\varinjlim_i F(X_i)^\vee
\]
is naturally a coalgebra whose category of comodules in $\mathcal A$ is equivalent to $\mathcal C$. The functor is the following: Any $X\in \mathcal C$ defines the object $F(X)\in \mathcal A$ and for any $i$ large enough so that $X\in \mathcal C_i$ a map $F(X)\otimes X_i\to X$ (by adjunction), thus a map
\[
F(X)\otimes F(X_i)\cong F(F(X)\otimes X_i)\to F(X),
\]
and hence dually we get the map
\[
F(X)\to F(X)\otimes F(X_i)^\vee\to F(X)\otimes \mathcal H.
\]

Moreover, for any $i$, $j$ there is some $k$ such that
\[
\mathcal C_i\otimes \mathcal C_j\subset \mathcal C_k:
\]
indeed, $\mathcal C_i$ (resp.~$\mathcal C_j$) is generated by $X_i$ (resp.~$X_j$) under tensors with $\mathcal A$ and coequalizers of $F$-split parallel pairs, so $\mathcal C_i\otimes \mathcal C_j$ is generated by $X_i\otimes X_j$ under these operations. Thus, for any $k$ such that $X_i\otimes X_j\in \mathcal C_k$, we actually have $\mathcal C_i\otimes \mathcal C_j\subset \mathcal C_k$. Let $X_k\in \mathcal C_k$ represent $F|_{\mathcal C_k}$; then we have a natural map
\[
X_k\to X_i\otimes X_j.
\]
Indeed, this is adjoint to a map $1\to F(X_i\otimes X_j)=F(X_i)\otimes F(X_j)$, for which we use the tensor product of the unit maps $1\to F(X_i)$, $1\to F(X_j)$. This means that there is a natural map
\[
\mathcal H\otimes \mathcal H=\varinjlim_{i,j} F(X_i)^\vee\otimes F(X_j)^\vee\cong \varinjlim_{i,j} F(X_i\otimes X_j)^\vee\to \varinjlim_k F(X_k)^\vee=\mathcal H,
\]
which turns $\mathcal H$ into a commutative algebra, where the unit is induced by the maps $X_i\to 1$ adjoint to $1=F(1)\in \mathcal A$ (inducing maps $1=F(1)\to F(X_i)^\vee$).

It is a matter of unraveling definitions that this makes $\mathcal H$ into a Hopf algebra whose symmetric monoidal category of representations in $\mathcal A$ is exactly $\mathcal C$. If $\mathcal C$ is rigid, one also sees that $\mathcal A$ admits an inverse. Indeed, one can write
\[
F(X_i)^\vee\cong F(X_i^\vee)\cong \sHom(X_i,X_i^\vee)\cong \sHom(X_i\otimes X_i,1)
\]
and the switching of the two factors defines the desired involution on $\mathcal A$. Here $\sHom\in \mathcal A$ denotes the internal Hom over $\mathcal A$.
\end{proof}

We can apply Proposition~\ref{prop:tannakareconstruction} to $\Sat_G^I(\Lambda)$ to get Hopf algebras
\[
\mathcal H_G^I(\Lambda)\in \Ind(\Rep_{W_E^I}(\Lambda)).
\]

\begin{proposition}\label{prop:hopfalgebraexteriortensor} The exterior tensor product
\[
\bigboxtimes_{i\in I}: \prod_{i\in I} \Sat_G^{\{i\}}(\Lambda)\to \Sat_G^I(\Lambda)
\]
induces an isomorphism
\[
\bigotimes_{i\in I} \mathcal H_G^{\{i\}}(\Lambda)\cong \mathcal H_G^I(\Lambda).
\]
\end{proposition}

\begin{proof} This is a consequence of the construction of the Hopf algebras together with Proposition~\ref{prop:leftadjointsatake}~(ii), noting that
\[
F^I(\bigast_{i\in I} A_i)\cong \bigotimes_{i\in I} F(A_i).\qedhere
\]
\end{proof}

We see that all information about the categories $\Sat_G^I(\Lambda)$ is in the Hopf algebra
\[
\mathcal H_G(\Lambda)=\mathcal H_G^{\{\ast\}}(\Lambda)\in \Ind(\Rep_{W_E}(\Lambda)).
\]
Note also that the construction of $\mathcal H_G(\Lambda)$ is compatible with base change in $\Lambda$, so it is enough to consider the case $\Lambda=\mathbb Z/n\mathbb Z$ with $n$ prime to $p$. In fact, note that we can formally take the inverse limit over $n$ to define
\[
\Sat_G(\hat{\mathbb Z}^p)=\varprojlim_n \Sat_G(\mathbb Z/n\mathbb Z)
\]
with a fibre functor into
\[
\Rep_{W_E}^\cont(\hat{\mathbb Z}^p)=\varprojlim_n \Rep_{W_E}(\mathbb Z/n\mathbb Z),
\]
the category of continuous representations of $W_E$ on finite free $\hat{\mathbb Z}^p=\varprojlim_n \mathbb Z/n\mathbb Z$-modules, yielding a Hopf algebra
\[
\mathcal H_G\in \Ind(\Rep_{W_E}^\cont(\hat{\mathbb Z}^p)).
\]
This can equivalently be thought of as an affine group scheme $\widecheck{G}$ over $\hat{\mathbb Z}^p$, with an action of $W_E$ that is in a suitable sense continuous.

\section{Identification of the dual group}

Our goal is to identify $\widecheck{G}$ with the Langlands dual group of $G$. Recall that the universal Cartan of $G$ defines a cocharacter group $X_\ast$ as an \'etale sheaf on $\Spec(E)$, i.e.~equivalently a finite free abelian group $X_\ast$ together with an action of the absolute Galois group of $E$, and in particular of $W_E$. It comes with the $W_E$-stable set of coroots $\Phi^\vee\subset X_\ast$, and the subset of positive roots $\Phi^\vee_+$. Dually, we have the cocharacters $X^\ast$ and dominant Weyl chamber $(X^\ast)^+\subset X^\ast$, and the roots $\Phi\subset X^\ast$, containing the positive roots $\Phi_+\subset \Phi$. These data give rise to a pinned Chevalley group scheme $\hat{G}$ over $\hat{\mathbb Z}^p$ (or already over $\mathbb Z$, but we will only consider it over $\hat{\mathbb Z}^p$) corresponding to the dual root data $(X^\ast,\Phi,X_\ast,\Phi^\vee)$. Being pinned, there are distinguished torus and Borel $\hat{T}\subset \hat{B}\subset\hat{G}$, isomorphisms $X^\ast(\hat{T})\cong X_\ast$ under which the positive coroots $\Phi^\vee_+$ correspond to the weights of $\hat{T}$ on $\Lie \hat{B}/\Lie\hat{T}$, so
\[
\Lie\hat{B}/\Lie\hat{T} = \bigoplus_{a\in \Phi^\vee_+} \Lie \hat{U}_a
\]
for root subgroups $\hat{U}_a\subset \hat{B}$. Moreover, one has fixed pinnings
\[
\psi_a: \Lie \hat{U}_a\cong \hat{\mathbb Z}^p
\]
for all simple roots $a$. We want to endow $\hat{G}$ with a $W_E$-action. We already have the $W_E$-action on $(X^\ast,\Phi,X_\ast,\Phi^\vee)$, but we need to twist the action on the pinning. More precisely, let us write the pinning instead with a Tate twist as
\[
\psi_a: \Lie \hat{U}_a\cong \hat{\mathbb Z}^p(1).
\]
Then $W_E$ acts naturally on the pinning as well, and thereby induces an action of $W_E$ on $\hat{G}$.

We aim to prove the following theorem. Recall that we write $\widecheck{G}$ for the Tannaka group arising from the Satake category. Generally, we will denote by $\widecheck{-}$ various objects defined via the Satake category, while by $\hat{-}$ we will denote objects formally defined as Langlands duals.

\begin{theorem}\label{thm:geometricsatake} There is a canonical $W_E$-equivariant isomorphism $\widecheck{G}\cong \hat{G}$.
\end{theorem}

We note that the formulation of this theorem is slightly more precise than the formulation in \cite{MirkovicVilonen}, where no canonical isomorphism is given. Also, we handle the case of non-split groups. Note that in particular, $\widecheck{G}$ only depends on $G$ up to inner automorphisms; this is not clear.

To prove the theorem, we can work over $\mathbb Z_\ell$ for some $\ell\neq p$: Indeed, the statement of the theorem is equivalent to having isomorphisms over $\mathbb Z/n\mathbb Z$ for all $n$ prime to $p$ (by the Tannakian perspective), so the reduction follows from the Chinese remainder theorem.

We will now first prove the theorem when the group $G$ is split; more precisely, if we have fixed a split torus and Borel $T\subset B\subset G$ and trivializations of all simple root groups $U_a\subset B$. Afterwards, we will verify that the isomorphism does not depend on this pinning (essentially, as pinnings vary algebraically, while automorphisms of $\hat{G}/\mathbb Z_\ell$ form an $\ell$-adic group), and finally use Galois descent to deduce the result in general.

Note first that if $G=T$ is a torus, then $\Gr_{T,\Div^1_X}\cong X_\ast(T)\times \Div^1_X$, and it is clear that $\Sat_T$ is just the category of $X_\ast(T)$-graded objects in $\Rep_{W_E}^\cont(\mathbb Z_\ell)$. This implies that $\widecheck{T}\cong \hat{T}$ is the dual torus with $X^\ast(\hat{T})=X_\ast(T)$.

We have the symmetric monoidal constant term functor
\[
\mathrm{CT}_B[\mathrm{deg}]: \Sat_G\to \Sat_T,
\]
and it commutes with the fibre functors by the identity $\bigoplus_i \mathcal H^i(R\pi_{G\ast})\cong \mathcal H^0(R\pi_{T\ast} \mathrm{CT}_B[\mathrm{deg}])$. This gives rise to a $W_E$-equivariant map $\hat{T}=\widecheck{T}\to \widecheck{G}$. Using the objects $A_\mu={}^p j_{\mu!} \mathbb Z_\ell[d_\mu]$, whose $\mu$-weight space is $1$-dimensional, we see that the map $\widecheck{T}\to \widecheck{G}$ must be a closed immersion.

We have the following information about the generic fibre $\widecheck{G}_{\mathbb Q_\ell}$, following \cite[Section 7]{MirkovicVilonen}. First, it follows from Proposition~\ref{prop:satakesemisimple} that its category of representations $\Sat_G(\mathbb Q_\ell)$ is given by
\[
\Sat_G(\mathbb Q_\ell)\cong \bigoplus_\mu \Rep_{W_E}^\cont(\mathbb Q_\ell)\otimes A_\mu.
\]
(Here
\[
\Sat_G(\mathbb Q_\ell) = \Sat_G(\mathbb Z_\ell)[\tfrac 1\ell],
\]
where $\Sat_G(\mathbb Z_\ell)=\varprojlim_m \Sat_G(\mathbb Z/\ell^m\mathbb Z)$.) The category of representations of $\widecheck{G}_{\mathbb Q_\ell}$ as an abstract group scheme is then given by
\[
\Sat_G(\mathbb Q_\ell)\otimes_{\Rep_{W_E}^\cont(\mathbb Q_\ell)} \mathrm{Vect}(\mathbb Q_\ell)\cong \bigoplus_\mu \mathrm{Vect}(\mathbb Q_\ell)\otimes A_\mu,
\]
and in particular is semisimple. As $A_\mu\ast A_{\mu'}$ contains $A_{\mu+\mu'}$ as a direct summand and $X_\ast^+$ is finitely generated as a monoid, we see that $\Sat_G(\mathbb Q_\ell)$ has a finite number of tensor generators. This implies that $\widecheck{G}_{\mathbb Q_\ell}$ is of finite type by \cite[Proposition 2.20]{DeligneMilne}. Moreover, it is connected as $\Sat_G(\mathbb Q_\ell)$ does not have nontrivial finite tensor subcategories (as for any $A_\mu$ with $\mu\neq 0$, the tensor category generated by $A_\mu$ contains all $A_{n\mu}$), cf.~\cite[Corollary 2.22]{DeligneMilne}. As $\Sat_G(\mathbb Q_\ell)$ is semisimple, we even know that $\widecheck{G}_{\mathbb Q_\ell}$ is reductive by \cite[Proposition 2.23]{DeligneMilne}. For any simple object $A_\mu$, the weights of $A_\mu$ on $\hat{T}_{\mathbb Q_\ell}\to \widecheck{G}_{\mathbb Q_\ell}$ are contained in the set of all $\lambda\in X_\ast=X^\ast(\hat{T})$ such that the dominant representative of $\lambda$ is bounded by $\mu$ in the dominance order, and contains $\mu$ (with weight $1$). This implies that $\widecheck{T}_{\mathbb Q_\ell}\to \widecheck{G}_{\mathbb Q_\ell}$ is a maximal torus of $\widecheck{G}_{\mathbb Q_\ell}$. We can also define a subgroup $\widecheck{B}\subset \widecheck{G}$ as the stabilizer of the filtration associated to the cohomological grading of $F$ (stabilizing the filtration $\bigoplus_{m\leq i} R^m\pi_{G\ast}$ on the fibre functor $F=\bigoplus_m R^m\pi_{G\ast}$); then $\widecheck{B}_{\mathbb Q_\ell}\subset \widecheck{G}_{\mathbb Q_\ell}$ is a Borel.

Now we analyze the case $G=\PGL_2$. In that case, we have the minuscule cocharacter $\mu: \mathbb G_m\to G$ giving rise to the minuscule Schubert cell $\Gr_{G,\Div^1_X,\mu}\cong \mathbb P^1_{\Div^1_X}$. Then
\[
F(A_\mu)=H^0(\mathbb P^1)\oplus H^2(\mathbb P^1)=\mathbb Z_\ell\oplus \mathbb Z_\ell(-1)
\]
as $W_E$-representation. This is a representation of $\widecheck{G}$, giving a natural map $\widecheck{G}\to \GL(\mathbb Z_\ell\oplus \mathbb Z_\ell(-1))$. We claim that this is a closed immersion, with image given by $\SL(\mathbb Z_\ell\oplus \mathbb Z_\ell(-1))$. Note that $\widecheck{T}$ acts on $\mathbb Z_\ell\oplus \mathbb Z_\ell(-1)$ with weight $\pm 1$, and in particular lands inside $\SL(\mathbb Z_\ell\oplus \mathbb Z_\ell(-1))$. As $\widecheck{G}_{\mathbb Q_\ell}$ is reductive of rank $1$, it necessarily follows that
\[
\widecheck{G}_{\mathbb Q_\ell}\to \SL(\mathbb Q_\ell\oplus \mathbb Q_\ell(-1))
\]
is an isomorphism, and integrally we get a map $\widecheck{G}\to \SL(\mathbb Z_\ell\oplus \mathbb Z_\ell(-1))$. This gives a map $\widecheck{G}_{\mathbb F_\ell}\to \SL(\mathbb F_\ell\oplus \mathbb F_\ell(-1))$. Let $H\subset \SL(\mathbb F_\ell\oplus \mathbb F_\ell(-1))$ be the closed subgroup that is the image of $\widecheck{G}_{\mathbb F_\ell}$. Note that the irreducible representations of $\widecheck{G}_{\mathbb F_\ell}$ are in bijection with dominant cocharacters, corresponding to the simple objects $B_\mu=j_{\mu!\ast} \mathbb F_\ell$ on $\Gr_{G,\Spd C}$; each $B_\mu$ has a highest weight vector given by weight $\mu$. It follows that $H$ satisfies the hypothesis of the next lemma.

\begin{lemma}\label{lem:subgroupsofSL2} Let $H$ be a closed subgroup of $\SL_2/\mathbb F_\ell$ containing the diagonal torus such that its set of irreducible representations injects into $\mathbb Z_{\geq 0}$ via consideration of highest weight vectors. Then $H=\SL_2$.
\end{lemma}

\begin{proof} Using a power of Frobenius, one can assume that $H$ is reduced, and thus smooth. By \cite[Corollary 2.22]{DeligneMilne} and consideration of highest weight vectors, one also sees that $H$ must be connected. Then $H$ is either the torus, or a Borel, or $\SL_2$. The first cases lead to too many irreducible representations.
\end{proof}

Thus, the map $\widecheck{G}_{\mathbb F_\ell}\to \SL(\mathbb F_\ell\oplus \mathbb F_\ell(-1))$ is surjective. Together with the isomorphism on the generic fibre, this implies formally that $\widecheck{G}\to \SL(\mathbb Z_\ell\oplus \mathbb Z_\ell(-1))$ is an isomorphism by the following lemma (used on the level of the corresponding Hopf algebras).

\begin{lemma}\label{lem:genericisospecialinj} Let $f: M\to N$ be a map of flat $\mathbb Z_\ell$-modules such that $M/\ell\to N/\ell$ is injective and $M[\tfrac 1\ell]\to N[\tfrac 1\ell]$ is an isomorphism. Then $f$ is an isomorphism.
\end{lemma}

\begin{proof} As $M$ is flat, $f: M\to N$ is injective; moreover, for any $x\in N$ there is some minimal $k$ such that $\ell^k n=f(m)$ lies in the image of $M$. But if $k>0$, then $m$ lies in the kernel of $M/\ell\to N/\ell$, a contradiction.
\end{proof}

The subgroup $\widecheck{B}\subset \widecheck{G}$ is then given by the Borel stabilizing the line $\mathbb Z_\ell\subset \mathbb Z_\ell\oplus \mathbb Z_\ell(-1)$. Its unipotent radical is the space of maps $\mathbb Z_\ell(-1)\to \mathbb Z_\ell$, which is canonically isomorphic to $\mathbb Z_\ell(1)$. This finishes the proof of the theorem for $G=\PGL_2$.

If now $G$ is of rank $1$, we get the map $G\to G_{\mathrm{ad}}\cong \PGL_2$, where the isomorphism $G_{\mathrm{ad}}\cong \PGL_2$ is uniquely determined by our choice of pinning. The map
\[
\Gr_{G,\Div^1_X}\to \Gr_{G_{\mathrm{ad}},\Div^1_X}
\]
is an isomorphism when restricted to each connected component, inducing an isomorphism
\[
\Gr_{G,\Div^1_X}\cong \pi_1(G)\times_{\pi_1(G_{\mathrm{ad}})} \Gr_{G_{\mathrm{ad}},\Div^1_X}.
\]
Here of course $\pi_1(G_{\mathrm{ad}})\cong \mathbb Z/2\mathbb Z$. This implies that $\Sat_G$ can be equivalently described as the category of $A\in \Sat_{G_{\mathrm{ad}}}$ together with a refinement of the $\mathbb Z/2\mathbb Z$-grading to a $\pi_1(G)$-grading. This implies that $\widecheck{G}=\widecheck{G_{\mathrm{ad}}}\times^{\mu_2} \widecheck{Z}$ where $\widecheck{Z}$ is the split torus with character group $\pi_1(G)$. Thus, one gets an isomorphism $\widecheck{G}\cong \hat{G}$ also in this case, including the isomorphism $\psi_a$ on the root group.

Coming back to a general split group $G$, let $a$ be any simple coroot. We now look at the corresponding minimal Levi subgroups $M_a\subset G$ properly containing $T$, with parabolic $P_a\subset B$. We have the symmetric monoidal constant term functor
\[
\mathrm{CT}_{P_a}[\mathrm{deg}_{P_a}]: \Sat_G\to \Sat_{M_a},
\]
commuting with the functors to $\Sat_T$. This induces a map $\widecheck{M_a}\to \widecheck{G}$, commuting with the inclusion of $\widecheck{T}$ into both. In particular, passing to Lie algebras, we see that $a\in X_\ast= X^\ast(\widecheck{T})$ is a root of $\widecheck{G}$, and $a^\vee\in X^\ast=X_\ast(\widecheck{T})$ is a coroot of $\widecheck{G}$. Moreover, if $s_a\in W$ is the corresponding simple reflection for $G$, we also see that $s_a\in \widecheck{W}$, the Weyl group of the reductive group $\widecheck{G}_{\mathbb Q_\ell}$. Using this information for all $a$, we see that $W\subset \widecheck{W}$, and that under $X_\ast=X^\ast(\widecheck{T})$ resp.~$X^\ast=X_\ast(\widecheck{T})$, we have
\[
\Phi^\vee\subset \Phi(\widecheck{G}_{\mathbb Q_\ell}),\ \Phi\subset \Phi^\vee(\widecheck{G}_{\mathbb Q_\ell}).
\]
Moreover, for any irreducible object $A_\mu\in \Sat_G(\mathbb Q_\ell)$, the weights of $A_\mu$ are contained in the convex hull of the $W$-orbit of $\mu$. This implies that these inclusions must be isomorphisms --- indeed, the directions of the edges emanating from $\mu$, for $\mu$ regular, correspond to $\Phi(\widecheck{G}_{\mathbb Q_\ell})$. Together with the isomorphisms on simple root groups, we get a unique isomorphism
\[
\widecheck{G}_{\mathbb Q_\ell}\cong \hat{G}_{\mathbb Q_\ell}.
\]
Under this isomorphism, the map $\hat{M_a}\cong \widecheck{M_a}\to \widecheck{G}$ is compatible with the map $\hat{M_a}\to \hat{G}$ induced by Langlands duality. It follows that
\[
\widecheck{G}(\breve{\mathbb Z}_\ell)\subset \widecheck{G}(\breve{\mathbb Q}_\ell)\cong \hat{G}(\breve{\mathbb Q}_\ell)
\]
is a subgroup containing all $\hat{M_a}(\breve{\mathbb Z}_\ell)$. But these generate $\hat{G}(\breve{\mathbb Z}_\ell)$, so $\hat{G}(\breve{\mathbb Z}_\ell)\subset \widecheck{G}(\breve{\mathbb Z}_\ell)$. Now pick a representation $\widecheck{G}\to \GL_N$ (given by some object of $\Sat_G$) that is a closed immersion over $\mathbb Q_\ell$. By the inclusion $\hat{G}(\breve{\mathbb Z}_\ell)\subset \widecheck{G}(\breve{\mathbb Z}_\ell)$, we see that the map $\hat{G}_{\mathbb Q_\ell}\cong \widecheck{G}_{\mathbb Q_\ell}\to \GL_N$ extends to a map $\hat{G}\to \GL_N$. By Lemma~\ref{lem:reductiveclosed}, this is necessarily a closed immersion, at least if $\ell\neq 2$ or $\hat{G}$ is simply connected. We can always reduce to the case that $\hat{G}$ is simply connected by arguing with the adjoint group $G_{\mathrm{ad}}$ (whose dual group $\hat{G_{\mathrm{ad}}}$ is simply connected) first, as in the discussion of rank-$1$-groups above. It then follows that $\widecheck{G}\to \GL_N$ factors over $\hat{G}\hookrightarrow \GL_N$, giving a map $\widecheck{G}\to \hat{G}$ that is an isomorphism in the generic fibre, and surjective in the special fibre (as any $\overline{\mathbb F}_\ell$-point of $\hat{G}$ lifts to $\breve{\mathbb Z}_\ell$, and then to $\widecheck{G}(\breve{\mathbb Z}_\ell)$), and hence an isomorphism by Lemma~\ref{lem:genericisospecialinj}.

\begin{lemma}[{\cite[Corollary 5.2]{PrasadYuQuasiReductive}}]\label{lem:reductiveclosed} Let $H$ be a reductive group over $\mathbb Z_\ell$, $H'$ some affine group scheme of finite type over $\mathbb Z_\ell$, and let $\rho: H\to H'$ be a homomorphism that is a closed immersion in the generic fibre. Assume that $\ell\neq 2$, or that no almost simple factor of the derived group of $H_{\overline{\mathbb Q}_\ell}$ is isomorphic to $\mathrm{SO}_{2n+1}$ (e.g., the derived group of $H$ is simply connected). Then $\rho$ is a closed immersion.
\end{lemma}

This finishes the proof of Theorem~\ref{thm:geometricsatake} when $G$ is split, and endowed with a splitting. Now we prove independence of the choice of splitting. For this, we note that in fact the cohomological grading on $F$ alone determines $\widecheck{T}\subset \widecheck{G}$ as its stabilizer, and $\widecheck{B}\subset \widecheck{G}$ as the stabilizer of the associated filtration. It remains to check that the isomorphisms
\[
\psi_a: \Lie \hat{U}_a\cong \mathbb Z_\ell(1)
\]
are independent of the choices. For this, consider the flag variety $\Fl$ over $E$, parametrizing Borels $B\subset G$. Each such Borel comes with its torus $T$, which is the universal Cartan and thus descends to $E$. Equivalently, note that tori over $\Fl$ are equivalent to \'etale $\mathbb Z$-local systems, and as $\Fl$ is simply connected all of them come via pullback from $E$; this then gives the so-called universal Cartan $T$ over $E$, which is split as $G$ is split. Let $a$ be a simple coroot of $G$. At each point of $\Fl$, we get the corresponding parabolic $P_a\supset B$, with Levi $M_a$. Let $\tilde{\Fl}_a\to \Fl$ parametrize pinnings of $M_a$, i.e.~isomorphisms of $U_a$ with the additive group; this is a $\mathbb G_m$-torsor. Over $\tilde{\Fl}_a$, the universal group $M_a$ is constant, with adjoint group $M_{a,\mathrm{ad}}\cong \PGL_2$. Consider
\[
S=\tilde{\Fl}^\diamond/\phi^\Z\to \Spd E/\phi^\Z=\Div^1_X.
\]
Applying the constant term functor for $P_a$ over $S$ gives a symmetric monoidal functor
\[
\Sat(\Hloc_{G,\Div^1_X}\times_{\Div^1_X} S,\mathbb Z_\ell)\to \Sat(\Hloc_{M_a,\Div^1_X}\times_{\Div^1_X} S,\mathbb Z_\ell);
\]
here, being symmetric monoidal is verified by repeating the construction of the fusion product after the smooth pullback $S\to \Div^1_X$. Both sides admit fibre functors to $\LocSys(S,\mathbb Z_\ell)$; this contains $\LocSys(\Div^1_X,\mathbb Z_\ell)=\Rep_{W_E}^\cont(\mathbb Z_\ell)$ fully faithfully, and we can consider the symmetric monoidal full subcategories on which the fibre functors land in this subcategory. As the constant term functor is compatible with fibre functors, it induces a symmetric monoidal functor on these full subcategories, which are then easily seen to be equivalent to $\Sat_G$ and $\Sat_M$ (reconstructing both starting from Schubert cells). This shows that the constant term functor $\Sat_G\to \Sat_{M_a}$ is naturally independent of the choice of Borel, reducing us to the rank $1$ case. In the rank $1$ case, we can then further reduce to $\PGL_2$, and we have the minuscule Schubert variety, which is the flag variety $\Fl\cong \mathbb P^1$ of $G\cong \PGL_2$. There are canonical isomorphisms
\[
H^0(\Fl)=\mathbb Z_\ell,\ H^2(\Fl)=\mathbb Z_\ell(-1),
\]
and $\hat{U}_a$ is canonically isomorphic to $\Hom(H^2(\Fl),H^0(\Fl))\cong \mathbb Z_\ell(1)$.

Thus, we have shown that if $G$ is split, the isomorphism $\widecheck{G}\cong \hat{G}$ is canonical. Finally, the general case follows by Galois descent from a finite Galois extension $E'|E$ splitting $G$.

\section{Chevalley involution}

Any Chevalley group scheme $\hat{G}$ comes with the Chevalley involution, induced by the map on root data which on $X_\ast$ is given by $\mu\mapsto -w_0(\mu)$ where $w_0$ is the longest Weyl group element. Under the geometric Satake equivalence, this has a geometric interpretation: Namely, it essentially corresponds to the switching equivalence $\mathrm{sw}^\ast$. Note that one can upgrade
\[
\mathrm{sw}^\ast: \Sat_G^I(\Lambda)\to \Sat_G^I(\Lambda)
\]
to a symmetric monoidal functor by writing it as the composition of Verdier duality and the duality functor $\mathrm{sw}^\ast \mathbb D$ in $\Sat_G^I(\Lambda)$; moreover, this symmetric monoidal functor commutes with the fibre functor $F^I$ (as symmetric monoidal functors), cf.~Corollary~\ref{cor:niceisomorphisms}. Thus, $\mathrm{sw}^\ast$ induces an automorphism of the Tannaka group $\widecheck{G}$, commuting with the $W_E$-action.

\begin{proposition}\label{prop:chevalleyinvolution} Under the isomorphism $\widecheck{G}\cong \hat{G}$ with the dual group, the isomorphism $\mathrm{sw}^\ast$ is given by the Chevalley involution, up to conjugation by $\hat{\rho}(-1)\in \hat{G}_{\mathrm{ad}}(\mathbb Z_\ell)$.
\end{proposition}

\begin{remark} There is a different construction of the commutativity constraint on $\Sat_G$, not employing the fusion product, that relies on the Chevalley involution --- this is essentially a categorical version of the classical Gelfand trick to prove commutativity of the Satake algebra. For the Satake category, this construction was first proposed by Ginzburg \cite{GinzburgSatake}, who however overlooked the sign $\hat{\rho}(-1)$. Zhu's proof \cite{ZhuWitt} of the geometric Satake equivalence for $\Gr_G^{\mathrm{Witt}}$ used this approach, taking careful control of the signs; these are related to the work of Lusztig--Yun \cite{LusztigYun}. We remark that Zhu gives a different construction of the commutation of $\mathrm{sw}^\ast$ with the fibre functor, using instead that the two actions (on left and right) of $H^\ast([\ast/L^+G],\mathbb Q_\ell)$ on $H^\ast(\Hloc_G,A)$ agree for $A\in \Sat_G(\mathbb Q_\ell)$.
\end{remark}

\begin{proof} We note that this is really a proposition: The statement only asks about the commutation of a certain diagram, not some extra structure. For the statement, we can also forget about the $W_E$-action. In particular, enlarging $E$, we can assume that $G$ is split. As in the proof of Theorem~\ref{thm:geometricsatake}, one can reduce to the case that $G$ is adjoint, so $\hat{G}$ is semisimple and simply connected. We also fix a pinning of $G$.

Now, being pinned, $G$ has its own Chevalley involution $\theta: G\to G$, and by the functoriality of all constructions under isomorphisms, the induced automorphism of $\Sat_G$ corresponds to the Chevalley involution of $\hat{G}$. In other words, we need to see that the automorphism $\theta^\ast \mathrm{sw}^\ast: \Sat_G\to \Sat_G$ (which is symmetric monoidal, and commutes with the fibre functors) induces conjugation by $\widecheck{\rho}(-1)$ on $\widecheck{G}$.

We claim that the natural cohomological grading on the fibre functor $F: \Sat_G(\Lambda)\to \Rep_{W_E}(\Lambda)$ is compatible with $\mathrm{sw}^\ast$. In other words, we need to see that in Corollary~\ref{cor:niceisomorphisms}, the isomorphism $F(A)\cong F(\mathrm{sw}^\ast A)$ is compatible with the grading, which follows from its construction. In particular, it follows that the automorphism of $\widecheck{G}$ restricts to the identity on the corresponding cocharacter $2\widecheck{\rho}: \mathbb G_m\subset \widecheck{G}$. This implies already that it preserves $\widecheck{T}$ and the Borel $\widecheck{B}$  (as the centralizer and dynamical parabolic). Any such automorphism of $\widecheck{G}$ is given by conjugation by some element $s\in \widecheck{T}_{\mathrm{ad}}\subset \widecheck{G}_{\mathrm{ad}}$. We need to see that $s=\widecheck{\rho}(-1)$. Equivalently, the automorphism acts by negation on any simple root space $\widecheck{U}_a$ of $\widecheck{G}$.

We claim that the symmetric monoidal automorphism $\theta^\ast \mathrm{sw}^\ast: \Sat_G\to \Sat_G$ (commuting with the fibre functor) is compatible with the constant term functors $\mathrm{CT}_P$, for any standard parabolic $P\supset B$, and the similar functor on Levi subgroups. Let $\theta'$ be the composition of $\theta$ with conjugation by $w_0$. We know, by the proof of Theorem~\ref{thm:geometricsatake}, that any inner automorphism of $G$ induces the identity on $\widecheck{G}$. Thus, it suffices to prove the similar claim for $\theta^{\prime\ast}\mathrm{sw}^\ast: \Sat_G\to \Sat_G$. Let $P^-\subset G$ be the opposite parabolic of $P$; then $\theta'(P^-)=P$, and the induced automorphism of the Levi $M$ is given by the corresponding automorphism $\theta'_M$ defined similarly as $\theta'$. Now Proposition~\ref{prop:hyperboliclocalizationdual} and the fusion definition of the symmetric monoidal structure (along with the definition of $\mathrm{sw}^\ast$ as the composite of Verdier duality and internal duality) give the claim.

These observations reduce us to the case $G=\PGL_2$. We note that in this case the Chevalley involution is the identity, so we can ignore $\theta$. We have the minuscule Schubert variety $i_\mu: \Gr_{G,\mu}=\mathbb P^1\subset \Gr_G$ and the sheaf $A = i_{\mu\ast} \mathbb Z_\ell[1](\tfrac 12)\in \Sat_G$ (assuming without loss of generality $\sqrt{q}\in \Lambda$ to introduce a half-Tate twist), and we know $\widecheck{G} = \SL(F(A))$, where
\[
F(A)=F(A)_1\oplus F(A)_{-1} = H^0(\mathbb P^1)(\tfrac 12)\oplus H^2(\mathbb P^1)(\tfrac 12).
\]
The image of $A$ under $\mathrm{sw}^\ast$ is isomorphic to $A$ itself; fix an isomorphism. Then on the one hand
\[
F(A)\cong F(\mathrm{sw}^\ast A)
\]
as the functor $\mathrm{sw}^\ast: \Sat_G\to \Sat_G$ commutes with the fibre functor $F$, while on the other hand
\[
F(\mathrm{sw}^\ast A)\cong F(A)
\]
as the two objects are isomorphic. We need to see that the composite isomorphism is given by the diagonal action of $(u,-u)$ for some $u\in \mathbb Z_\ell^\times$ (this claim is independent of the chosen isomorphism between $A$ and $\mathrm{sw}^\ast A$). We already know that the isomorphism is graded, so it is given by diagonal multiplication by $(u_1,u_2)$ for some units $u_1,u_2\in \mathbb Z_\ell^\times$.

Recall that the first isomorphism is constructed as the composite of Verdier duality and internal duality in $\Sat_G$. Now the Verdier dual of $A$ is $A$ itself (because of the half-Tate twist), and the Verdier duality pairing
\[
F(A)\otimes F(A)\to \mathbb Z_\ell
\]
is the tautological pairing; in particular, restricted to $F(A)_{-1}\otimes F(A)_1$ and $F(A)_1\otimes F(A)_{-1}$ it is the same map, up to the natural commutativity constraint on $\mathbb Z_\ell$-modules. It follows that the internal dual $A^\vee$ of $A$ is also isomorphic to $A$, and picking such an identification we need to understand the induced pairing
\[
F(A)\otimes F(A)\to F(1)=\mathbb Z_\ell,
\]
and show that when restricted to $F(A)_{-1}\otimes F(A)_1$ and $F(A)_1\otimes F(A)_{-1}$, the two induced maps differ by a sign (up to the natural commutativity constraint); this claim is again independent of the chosen isomorphism between $A^\vee$ and $A$. But this is a question purely internal to the symmetric monoidal category $\Sat_G\cong \Rep(\SL_2)$ with its fibre functor. In there, we have the tautological representation $V=\mathbb Z_\ell^2 = \mathbb Z_\ell e_1\oplus \mathbb Z_\ell e_2$, and it has the determinant pairing $V\otimes V\to \mathbb Z_\ell$ as $\SL_2$-representation, realizing the internal duality. The determinant pairing is alternating, so takes opposite signs on $(e_1,e_2)$ and $(e_2,e_1)$, as desired.
\end{proof}

\chapter{$D_{\blacksquare}(X)$}\label{ch:solid}

In order to deal with smooth representations of $G(E)$ on $\overline{\mathbb Q}_\ell$-vector spaces (not Banach spaces), we extend (a modified form of) the $6$-functor formalism from \cite{ECoD} to the larger class of solid pro-\'etale sheaves. The results in this chapter were obtained in discussions with Clausen, and Mann has obtained analogues of some of these results in the case of schemes. (Strangely enough, in some respects the formalism actually works better for diamonds than for schemes.)

More precisely, we want to find a ``good'' category of $\mathbb{Q}_\ell$-sheaves on $[\ast/\underline{G(E)}]$ that corresponds to smooth representations of $G(E)$ with values in $\mathbb{Q}_\ell$-vector spaces, and extends to a category of $\mathbb{Q}_\ell$-sheaves on $\Bun_G$ with a good formalism of six operations that allows us to extend the preceding results for \'etale torsion coefficients. 
The first idea would be to take pro-systems of \'etale torsion sheaves as $\mathbb{Z}_\ell$-coefficients and invert formally $\ell$; this formalism is easy to construct, see \cite[Section 26]{ECoD}. This would give rise to continuous representations of $G(E)$ in $\mathbb{Q}_\ell$-Banach spaces, and we do not want that:
\begin{altitemize}
\item  supercuspidal representations of $G(E)$ in $\overline{\mathbb{Q}}_\ell$-vector spaces are defined over a finite degree extension of $\mathbb{Q}_\ell$, and after twist admit an invariant lattice that allows us to complete them $\ell$-adically, but we do not want to make such a choice.
\item we want to construct semi-simple Langlands parameters using the Bernstein center and not some $\ell$-adic completion of it.
\item in usual discussions of the cohomology of the Lubin--Tate tower, or more general Rapoport--Zink spaces, it {\it is} possible to use $\overline{\mathbb Q}_\ell$-coefficients while talking about usual smooth representations. We want to be able to achieve the same on the level of $\Bun_G$.
\end{altitemize}

We could take $\mathbb{Q}_\ell$-pro-\'etale sheaves. This would give rise to representations of $G(E)$ (seen as a condensed group) with values in condensed $\mathbb{Q}_\ell$-vector spaces. This category is too big; there is no hope to obtain a formalism of six operations in this too general context. We need to ask for some ``completeness'' of the sheaves, for which we take inspiration from the theory of solid abelian groups developed in \cite{Condensed}.

The idea is the following. We define a category of solid pro-\'etale $\mathbb{Q}_\ell$-sheaves on $\Bun_{G}$ with a good formalism of (a modified form of) six operations. More precisely, for any small v-stack, we define a full subcategory
\[
D_\solid(X,\mathbb Z_\ell)\subset D(X_v,\mathbb Z_\ell),
\]
compatible with pullback, and equipped with a symmetric monoidal tensor product (for which pullback is symmetric monoidal). A complex is solid if and only if each cohomology sheaf is solid, and this can be checked v-locally. The subcategory $D_\solid(X,\mathbb Z_\ell)$ is stable under all (derived) limits and colimits, and the inclusion into $D(X_v,\mathbb Z_\ell)$ admits a left adjoint. If $X$ is a diamond, then $D_\solid(X,\mathbb Z_\ell)$ is also a full subcategory of $D(X_\qproet,\mathbb Z_\ell)$. If $X$ is a spatial diamond, then on the abelian level, the category of solid $\mathbb Z_\ell$-sheaves is the $\Ind$-category of the $\Pro$-category of constructible \'etale sheaves killed by some power of $\ell$. In this way, one can bootstrap many results from the usual \'etale case.

For any map $f: Y\to X$ of small v-stacks, the pullback functor $f^\ast$ admits a right adjoint $Rf_\ast: D_\solid(Y,\mathbb Z_\ell)\to D_\solid(X,\mathbb Z_\ell)$ that in fact commutes with {\it any} base change, see Proposition~\ref{prop:solidbasechange}. Similarly, the formation of $R\sHom$ commutes with any base change. Both of these operations can a priori be taken in all v-sheaves, but turn out to preserve solid sheaves. This already gives us four operations.

Unfortunately, $Rf_!$ does not have the same good properties as usual. In particular, if $f$ is proper (and finite-dimensional), $Rf_\ast$ does not in general satisfy a projection formula. As a remedy, it turns out that for {\it all} $f$, the functor $f^\ast$ admits a left adjoint
\[
f_\natural: D_\solid(Y,\mathbb Z_\ell)\to D_\solid(X,\mathbb Z_\ell),
\]
given by ``relative homology''. This is a completely novel feature, and already for closed immersions this takes usual \'etale sheaves to complicated solid sheaves. Again $f_\natural$ commutes with any base change, and also satisfies the projection formula (which is just a condition here, as there is automatically a natural map).

When $f$ is ``proper and smooth'', one can moreover relate relative $f_\natural$ (``homology'') and $Rf_\ast$ (``cohomology'') in the expected way. One also gets a formula for the dualizing complex of $f$ in terms of such functors. These results even extend to universally locally acyclic complexes. This solid $5$-functor formalism thus has some excellent formal properties. We are somewhat confused about exactly how expressive it is, and whether it is preferable over the standard $6$-functor formalism.\footnote{One can also treat this $5$-functor formalism as a $6$-functor formalism in which $f^!=f^\ast$ for all maps $f$, i.e.~``all maps are \'etale''.} One advantage is certainly that $f_\natural$ is more canonical, and even defined much more generally, than $Rf_!$ (whose construction for stacky maps requires the resolution of subtle homotopy coherence issues, and also can only be defined for {\it certain} (finite-dimensional) maps). The main problem with the solid formalism is that a stratification of a stack does not lead to a semi-orthogonal decomposition on the level of $D_\solid$.

On the other hand, for our concerns here, $D_\solid(\Bun_G,\mathbb Z_\ell)$ is much too large. On $[\ast/\underline{G(E)}]$ this gives rise to representations of $G(E)$ (as a condensed group) with values in solid $\mathbb{Z}_{\ell}$-modules. The category of discrete $\mathbb{Q}_\ell$-vector spaces injects into the category of solid $\mathbb{Z}_{\ell}$-modules. In fact, $\mathbb{Q}_\ell$-vector spaces are the same as Ind-finite dimensional vector spaces. Any finite dimensional $\mathbb{Q}_\ell$-vector space is complete and thus ``solid''. This means that if $V$ is a $\mathbb{Q}_\ell$-vector space then it gives rise to the solid condensed sheaf $V\otimes_{\mathbb{Q}_\ell^{disc}} \mathbb{Q}_{\ell}$ 
 whose value on the profinite set $S$ is 
$$
\varinjlim_{\substack{W\subset V\\ \text{finite dim.}}} \mathrm{Cont}(S,W).
$$

The category of smooth representations of $G(E)$ with values in discrete $\mathbb{Q}_\ell$-vector spaces injects in the category of solid $\mathbb{Q}_{\ell}$-pro-\'etale shaves on $[\ast/\underline{G}(E)]$. In fact, since $\ell\neq p$ and $G(E)$ is locally pro-$p$, representations of the condensed group $G(E)$ on the condensed $\mathbb{Q}_\ell$-vector space $V\otimes_{\mathbb{Q}_\ell^{disc}} \mathbb{Q}_\ell$ are the same as smooth representations of $G(E)$ on $V$.

We then cut out a subcategory $D_{\lis} (\Bun_G, \mathbb{Q}_\ell)$ of $D_{\solid} (\Bun_G,\mathbb{Z}_\ell)$ that gives back the category  of smooth representations of $G(E)$ in $\mathbb{Q}_\ell$-vector spaces when we look at $[\ast/\underline{G(E)}]$. (Of course, we can also stick with $\mathbb Z_\ell$-coefficients.)

\section{Solid sheaves}

In the following, $\hat{\mathbb Z}$ always denotes the pro-\'etale sheaf $\underline{\hat{\mathbb Z}} = \varprojlim_n \mathbb Z/n\mathbb Z$ where $n$ runs over nonzero integers. We will quickly restrict attention to $\hat{\mathbb Z}^p=\varprojlim_{(n,p)=1} \mathbb Z/n\mathbb Z$, allowing only $n$ prime to $p$.

Let $X$ be a spatial diamond. We have the associated quasi-pro-\'etale site $X_\qproet$. Not all of its objects are cofiltered limits of \'etale maps, but one has a fully faithful functor
\[
\mathrm{Pro}(X_\et^{\mathrm{qcqs}})\to X_\qproet^{\mathrm{qcqs}}
\]
on the level of quasicompact and quasiseparated objects; this follows from \cite[Proposition 11.23 (ii)]{ECoD}. Moreover, \cite[Proposition 11.24]{ECoD} ensures that this full subcategory is a basis for the quasi-pro-\'etale topology. In the following constructions we will often work in this subcategory.

For any quasi-pro-\'etale $j: U\to X$ that can be written as a cofiltered inverse limit of qcqs \'etale $j_i: U_i\to X$, we let
\[
j_\natural \hat{\mathbb Z} = \varprojlim_i j_{i!} \hat{\mathbb Z};
\]
as the pro-system of the $U_i$ is unique, this is well-defined. Note that there is a tautological section of $j_\natural\hat{\mathbb Z}$ over $U$. Equivalently, if one writes $\hat{\mathbb Z}[U]$ for the free pro-\'etale sheaf of $\hat{\mathbb Z}$-modules generated by $U$ (noting that $j^\ast$ admits a left adjoint on pro-\'etale sheaves, being a slice), there is a natural map $\hat{\mathbb Z}[U]\to \varprojlim_i \widehat{\mathbb{Z}}[U_i] = j_\natural \hat{\mathbb Z}$.

When $X=\Spa (C)$, $j_{\natural}\widehat{\mathbb{Z}}$ is the sheaf denoted $\widehat{\mathbb{Z}}[U]^\solid$ in \cite{Condensed}, and the same notation will be appropriate here in general.

\begin{definition}\label{def:solid} Let $\mathcal F$ be a pro-\'etale sheaf of $\hat{\mathbb Z}$-modules on $X$. Then $\mathcal F$ is solid if for all $j: U\to X$ as above, the map
\[
\Hom(j_\natural\hat{\mathbb Z},\mathcal F)\to \mathcal F(U)
\]
is an isomorphism.
\end{definition}

Let us begin with the following basic example. We note $\nu:X_{\qproet}\to X_{\et}$ the projection to the \'etale site.

\begin{proposition}\label{prop:etale implique solide}
For any \'etale sheaf $\mathcal{F}$ of $\mathbb Z/n\mathbb Z$-modules on $X_{\et}$, $\nu^*\mathcal{F}$ is solid.
\end{proposition}

\begin{proof}
This is a consequence of \cite[Proposition 14.9]{ECoD}. 
\end{proof}

The notion of solid sheaf is well-behaved:

\begin{theorem}\label{thm:solid} The category of solid $\hat{\mathbb Z}$-sheaves on $X$ is an abelian subcategory of all pro-\'etale $\hat{\mathbb Z}$-sheaves on $X$, stable under all limits, colimits, and extensions. It is generated by the finitely presented objects $j_\natural \hat{\mathbb Z}$ for quasi-pro-\'etale $j: U\to X$ as above, and the inclusion admits a left adjoint $\mathcal F\mapsto \mathcal F^\solid$ that commutes with all colimits.

Let $\mathcal F$ be a pro-\'etale $\hat{\mathbb Z}$-sheaf on $X$. The following conditions are equivalent.
\begin{altenumerate}
\item[(1)] The $\hat{\mathbb Z}$-sheaf $\mathcal F$ is finitely presented in the category of all pro-\'etale $\hat{\mathbb Z}$-sheaves, and is solid.
\item[(2)] The $\hat{\mathbb Z}$-sheaf $\mathcal F$ is solid, and finitely presented in the category of solid $\hat{\mathbb Z}$-sheaves.
\item[(3)] The $\hat{\mathbb Z}$-sheaf $\mathcal F$ can be written as a cofiltered inverse limit of torsion constructible \'etale sheaves.
\end{altenumerate}
For any such $\mathcal F$, the underlying pro-\'etale sheaf is representable by a spatial diamond. The category of $\mathcal F$ satisfying (1) -- (3) is stable under kernels, cokernels, and extensions, in particular an abelian category, and is equivalent to the $\Pro$-category of torsion constructible \'etale sheaves.

Moreover, the category of all solid $\hat{\mathbb Z}$-sheaves on $X$ is equivalent to the $\Ind$-category of the $\Pro$-category of torsion constructible \'etale sheaves.
\end{theorem}

\begin{question} If $\mathcal F$ is a pro-\'etale $\hat{\mathbb Z}$-sheaf on $X$ whose underlying pro-\'etale sheaf is representable by a spatial diamond, or even is just qcqs, is $\mathcal F$ necessarily solid? If so, these conditions would also be equivalent to (1) -- (3).
\end{question}

Let us remark the following lemma.

\begin{lemma}\label{lem:constructible implique diamond}
Any torsion constructible \'etale sheaf on the spatial diamond $X$ is represented by a spatial diamond.
\end{lemma}
\begin{proof}
Let $\mathcal{F}\to X$ be such a sheaf. We can find a surjection of \'etale sheaves $\mathcal F'=\bigoplus_i j_{i!} \mathbb Z/n_i\mathbb Z\to \mathcal F$ for some quasicompact separated \'etale maps $j_i: U_i\to X$ and nonzero integers $n_i$ (where the direct sum is finite). Then $\mathcal F'$ is quasicompact separated \'etale over $X$, and thus a spatial diamond; and the surjective map $\mathcal F'\to \mathcal F$ is also quasicompact separated \'etale, in particular universally open, and so also $\mathcal F$ is a spatial diamond.
\end{proof}

Before starting the proof, we record a key proposition. Its proof is a rare instance that requires w-contractible objects --- in most proofs, strictly totally disconnected objects suffice.

\begin{proposition}\label{prop:limitofconstructible} Let $X$ be a spatial diamond and let $\mathcal F_i$, $i\in I$, be a cofiltered system of torsion constructible \'etale sheaves. Then for all $j>0$ the higher inverse limit
\[
R^j \varprojlim_i \mathcal F_i=0,
\]
taken in the category of pro-\'etale sheaves on $X$, vanishes.
\end{proposition}

\begin{proof} The pro-\'etale site of $X$ has a basis given by the w-contractible $Y\to X$, that is strictly totally disconnected perfectoid spaces $Y$ such that the closed points in $|Y|$ are a closed subset and $\pi_0 Y$ is an extremally disconnected profinite set; equivalently, any pro-\'etale cover $\tilde{Y}\to Y$ splits, cf.~\cite[Section 2.4]{BhattScholzeProetale} for a discussion of w-contractibility. Thus, it suffices to check sections on $Y$. In other words, we may assume that $X$ is w-contractible, and prove that
\[
R^j\varprojlim_i \mathcal F_i(X)=0
\]
for all $j>0$. As $X$ is strictly totally disconnected, sheaves on $X_\et$ are equivalent to sheaves on $|X|$. Moreover, we have the closed immersion $f: \pi_0 X\to |X|$ given by the closed points, and pullback along this map induces isomorphisms $\mathcal F_i(X)\cong (f^\ast \mathcal F_i)(\pi_0 X)$. Let $S=\pi_0 X$ be the extremally disconnected profinite set.

Then any torsion constructible sheaf on $S$ is locally constant with finite fibres. In particular, $f^\ast \mathcal F_i$ maps isomorphically to $\sHom(\sHom(f^\ast \mathcal F_i,\underline{S^1}),\underline{S^1})$ where $\underline{S^1} = \underline{\mathbb R/\mathbb Z}$ is the sheaf of continuous maps to the circle. (We could for the moment also use $\mathbb Q/\mathbb Z$, but it will become critical that $S^1$ is compact.) It follows that
\[
R^j\varprojlim_i \mathcal F_i(X) = \Ext^j_S(\varinjlim_i \sHom(f^\ast \mathcal F_i,\underline{S^1}),\underline{S^1}).
\]
Thus, the result follows from the injectivity of $\underline{S^1}$ as stated in the next lemma.
\end{proof}

\begin{lemma}\label{lem:S1injective} Let $S$ be an extremally disconnected profinite set. Then the abelian sheaf $\underline{\mathbb R/\mathbb Z}$ on $S$ is injective.
\end{lemma}

\begin{proof} First, we note that $\underline{\mathbb R/\mathbb Z}$ is flasque in a strong sense. Namely, if $U\subset S$ is any open subset with closure $\overline{U}\subset S$, then $\overline{U}$ is the Stone-\v{C}ech compactification of $U$ (as $\overline{U}\subset S$ is open by one definition of extremally disconnected spaces, and then $\beta U\sqcup (S\setminus \overline{U})\to S$ admits a splitting, which in particular gives a splitting of the surjection $\beta U\to \overline{U}$ that is the identity on $U$, thus implying that $\beta U\cong \overline{U}$) and hence any section of $\underline{\mathbb R/\mathbb Z}$ over $U$ extends uniquely to $\overline{U}$ as $\mathbb R/\mathbb Z$ is compact Hausdorff. Also, all sections over $\overline{U}$ extend to $S$, as $\overline{U}$ is open and closed in $S$.

Let $\mathcal F\hookrightarrow \mathcal G$ be an injection of sheaves on $S$ with a map $\mathcal F\to \underline{\mathbb R/\mathbb Z}$. Using Zorn's lemma, choose a maximal subsheaf of $\mathcal G$ containing $\mathcal F$ with an extension of the map to $\underline{\mathbb R/\mathbb Z}$. Replacing $\mathcal F$ by this maximal subsheaf, we can assume that $\mathcal F$ is maximal already. If $\mathcal F\to \mathcal G$ is not an isomorphism, then it is not an isomorphism on global sections (any local section not in the image can be extended by zero to form a global section not in the image), so we can find a map $\mathbb Z\to \mathcal G$ such that $\mathcal F'=\mathcal F\times_{\mathcal G} \mathbb Z\subset \mathbb Z$ is a proper subsheaf, and we can replace $\mathcal G$ by $\mathbb Z$ and assume that $\mathcal F$ is a proper subsheaf of $\mathbb Z$.

For each integer $n$, we can look at the open subset $j_n: U_n\subset S$ where $n\in \mathbb Z$ lies in $\mathcal F$. On this open subset, we have a map $n\mathbb Z\to \underline{\mathbb R/\mathbb Z}$, and by the above this extends uniquely to $\overline{j}_n: \overline{U_n}\hookrightarrow S$. The extension $\overline{j}_{n!} \mathbb Z\to \underline{\mathbb R/\mathbb Z}$ necessarily agrees with the restriction of the given map $\mathcal F\to \underline{\mathbb R/\mathbb Z}$ on the intersection $\overline{j}_{n!}\mathbb Z\cap \mathcal F\subset \overline{j}_{n!} \mathbb Z$, as this contains the dense subset $j_{n!} \mathbb Z$ and $\mathbb R/\mathbb Z$ is separated.

Thus, by maximality of $\mathcal F$, we see that necessarily all $U_n$ are open and closed, hence so are all $V_n=U_n\setminus \bigcup_{m<n} U_m$. Thus $V_1,V_2,\ldots\subset S$ are pairwise disjoint open and closed subsets such that $\mathcal F=\bigoplus_n n\mathbb Z|_{V_n}$. But one can then extend to $\bigoplus_n \mathbb Z|_{V_n}$ as the continuous maps from $V_n$ to $\mathbb R/\mathbb Z$ form a divisible group. By maximality of $\mathcal F$, this means that $V_n$ is empty for all $n>1$, and hence $\mathcal F = \mathbb Z|_{V_1}$ is a direct summand of $\mathbb Z$, in which case the possibility of extension is clear.
\end{proof}

Now we can give the proof of Theorem~\ref{thm:solid}.

\begin{proof}[Proof of Theorem~\ref{thm:solid}] The $\Pro$-category of torsion constructible \'etale sheaves is an abelian category, and by Proposition~\ref{prop:limitofconstructible} the functor to pro-\'etale $\hat{\mathbb Z}$-sheaves is exact. It is also fully faithful: For this, it suffices to see that if $\mathcal F_i$, $i\in I$, is a cofiltered inverse system of torsion constructible \'etale sheaves and $\mathcal G$ is any \'etale sheaf, then
\[
\varinjlim_i \Hom(\mathcal F_i,\mathcal G)\to \Hom(\varprojlim_i \mathcal F_i,\mathcal G)
\]
is an isomorphism. But the underlying pro-\'etale sheaf of each $\mathcal F_i$ is a spatial diamond over $X$, so by \cite[Proposition 14.9]{ECoD} (applied with $j=0$) we see that the similar result holds true when then taking homomorphisms of pro-\'etale sheaves (without the abelian group structure). Enforcing compatibility with addition amounts to a similar diagram for $\varprojlim_i \mathcal F_i\times \varprojlim_i \mathcal F_i$ to which the same argument applies.

We see that the $\Pro$-category of torsion constructible \'etale sheaves is a full subcategory $\mathcal C$ of all pro-\'etale $\hat{\mathbb Z}$-modules on $X$, stable under the formation of kernels and cokernels. All of these sheaves are solid: As the condition of being solid is stable under all limits, it suffices to see that any \'etale sheaf is solid; this is Proposition \ref{prop:etale implique solide}.

Also, by \cite[Lemma 11.22]{ECoD} and Lemma \ref{lem:constructible implique diamond} all objects of $\mathcal C$ have as underlying pro-\'etale sheaf a spatial diamond. Using Proposition~\ref{prop:derivedsolid} below (whose proof is direct), one also checks that $\mathcal C$ is stable under extensions.

Next, we prove that (3) implies (1), so let $\mathcal F$ be in $\mathcal C$; in particular, the underlying pro-\'etale sheaf is a spatial diamond. Then for any pro-\'etale $\hat{\mathbb Z}$-module $\mathcal G$, one can describe $\Hom(\mathcal F,\mathcal G)$ as the maps $\mathcal F\to \mathcal G$ of pro-\'etale sheaves satisfying additivity and $\hat{\mathbb Z}$-linearity, i.e.~certain maps $\mathcal F\times \mathcal F\to \mathcal G$ resp.~$\mathcal F\times \hat{\mathbb Z}\to \mathcal G$ agree. This description commutes with filtered colimits (as for spatial diamonds $Y$, the functor $\mathcal G\mapsto \mathcal G(Y)$ commutes with filtered colimits).

Now we can describe the full category of solid $\hat{\mathbb Z}$-sheaves. Indeed, using that $j_\natural \hat{\mathbb Z}$ is finitely presented in all pro-\'etale $\hat{\mathbb Z}$-modules by the previous paragraph, we see from the definition that the category of solid $\hat{\mathbb Z}$-sheaves is stable under all filtered colimits. In particular, we get an exact functor from the $\Ind$-category of $\mathcal C$ to solid $\hat{\mathbb Z}$-sheaves. This is also fully faithful, as all objects of $\mathcal C$ are finitely presented. Moreover, $\Ind(\mathcal C)$ is an abelian category for formal reasons. We see that $\Ind(\mathcal C)$ is a full subcategory of the category of pro-\'etale $\hat{\mathbb Z}$-sheaves stable under kernels and cokernels, and all of its objects are solid. Conversely, any solid $\hat{\mathbb Z}$-sheaf admits a surjection from a direct sum of objects of the form $j_\natural \hat{\mathbb Z}\in \mathcal C$, and the kernel of any such surjection is still solid, so we may write any solid $\hat{\mathbb Z}$-sheaf as the cokernel of a map in $\Ind(\mathcal C)$. As $\Ind(\mathcal C)$ is stable under cokernels, we see that $\Ind(\mathcal C)$ is exactly the category of solid $\hat{\mathbb Z}$-sheaves.

As filtered colimits of solid sheaves stay solid, it is now formal that (1) implies (2), and (2) implies (3) as $\mathcal C\subset \Ind(\mathcal C)$ are the finitely presented objects (as $\mathcal C$ is idempotent-complete). This finishes the proof of the equivalences.

The identification with $\Ind(\mathcal C)$ shows that the category of solid $\hat{\mathbb Z}$-sheaves is stable under kernels, cokernels, and filtered colimits. The latter two imply stability under all colimits, and stability under all limits is clear from the definition. One also easily checks stability under extensions by reduction to $\mathcal C$ (again, cf.~proof of Proposition~\ref{prop:derivedsolid}). For the existence of the left adjoint, note that it exists on the free pro-\'etale $\hat{\mathbb Z}$-modules generated by $U$, with value $j_\natural \hat{\mathbb Z}$, i.e.
$$
\widehat{\mathbb{Z}}[U]^\solid = j_{\natural} \widehat{\mathbb{Z}}.
$$
As these generate all pro-\'etale $\hat{\mathbb Z}$-modules, one finds that the left adjoint $\mathcal F\mapsto \mathcal F^\solid$ exists in general: one can write any $\mathcal{F}$ as a colimit $\varinjlim_\alpha \widehat{\mathbb{Z}}[U_\alpha]$ and $(\varinjlim_\alpha \widehat{\mathbb{Z}}[U_\alpha])^\solid= \varinjlim_\alpha \widehat{\mathbb{Z}}[U_\alpha]^\solid$.
\end{proof}

We have the following proposition on the functorial behaviour of the notion of solid $\hat{\mathbb Z}$-sheaves.

\begin{proposition}\label{prop:solidfunctor} Let $f: Y\to X$ be a map of spatial diamonds. Then pullbacks of solid $\hat{\mathbb Z}$-sheaves are solid, and the functor $f^\ast$ commutes with solidification. Moreover, if $f$ is surjective, and $\mathcal F$ is a pro-\'etale $\hat{\mathbb Z}$-sheaf on $X$ such that $f^\ast \mathcal F$ is solid, then $\mathcal F$ is solid.
\end{proposition}

\begin{proof} Recall that $f^\ast$ commutes with all limits (and of course colimits) by \cite[Lemma 14.4]{ECoD}. To check that $f^\ast$ commutes with solidification, it suffices to check on the pro-\'etale $\hat{\mathbb Z}$-modules $\hat{\mathbb Z}[U]$ generated by some quasi-pro-\'etale $j: U\to X$ that is a cofiltered inverse limit of qcqs \'etale maps $U_i\to X$, and in that case the claim follows from $f^\ast$ commuting with all limits,
\[
f^*(\widehat{\mathbb{Z}}[U]^\solid)= f^*(\varprojlim_i \widehat{\mathbb{Z}}[U_i])=\varprojlim_i f^*\widehat{\mathbb{Z}}[U_i]=\varprojlim_i \widehat{\mathbb{Z}}[f^*U_i]=\widehat{\mathbb{Z}}[f^*U]^\solid.
\]
In particular, applied to solid $\hat{\mathbb Z}$-sheaves on $X$, this implies that their pullback to $Y$ is already solid.

Now assume that $f$ is surjective and $\mathcal F$ is a pro-\'etale $\hat{\mathbb Z}$-sheaf on $X$ such that $f^\ast \mathcal F$ is solid. Let $j: U=\varprojlim_i U_i\to X$ be a quasi-pro-\'etale map. We first check that
\[
\Hom(j_\natural \hat{\mathbb Z},\mathcal F)\to \mathcal F(U)
\]
is injective. Indeed, assume that $g: j_\natural \hat{\mathbb Z}\to \mathcal F$ lies in the kernel. Then after pullback to $Y$, this map vanishes since $f^*\widehat{\mathbb{Z}}[U]^\solid=\widehat{\mathbb{Z}}[f^*U]^\solid$. But for any quasi-pro-\'etale $V\to X$, the map $\mathcal F(V)\to (f^\ast \mathcal F)(V\times_X Y)$ is injective (using \cite[Proposition 14.7]{ECoD} one has $f^*\mathcal{F}(V\times_X Y)=(\lambda_X^*\mathcal{F})(V\times_X Y)$ and we conclude since $V\times_X Y\to V$ is a v-cover), so it follows that $f=0$.

We see that an element of $\mathcal F(U)$ determines at most one map $j_\natural \hat{\mathbb Z}\to \mathcal F$, and this assertion stays true after any pullback. By \cite[Proposition 14.7]{ECoD}, it suffices to construct the map v-locally; but it exists after pullback to $Y\to X$, thus proving existence.
\end{proof}

In particular, it makes sense to make the following definition.

\begin{definition}\label{def:solidvstack} Let $Y$ be a small v-stack and let $\mathcal F$ be a v-sheaf of $\hat{\mathbb Z}$-modules on $Y$. Then $\mathcal F$ is solid if for all maps $f: X\to Y$ from a spatial diamond $X$, the pullback $f^\ast \mathcal F$ comes via pullback from a solid $\hat{\mathbb Z}$-sheaf on $X_\qproet$.
\end{definition}

Regarding passage to the derived category, we make the following definition.

\begin{definition}\label{def:derivedsolid} Let $X$ be a small v-stack. Let $D_\solid(X,\hat{\mathbb Z})\subset D(X_v,\hat{\mathbb Z})$ be the full subcategory of all $A$ such that each cohomology sheaf $\mathcal H^i(A)$ is solid.
\end{definition}

As being solid is stable under kernels, cokernels, and extensions, this defines a triangulated subcategory.

If $X$ is a diamond, one could alternatively define a full subcategory of $D(X_\qproet,\hat{\mathbb Z})$ by the same condition, and pullback from the quasi-pro-\'etale to the v-site defines a functor. This functor is an equivalence, by repleteness (to handle Postnikov towers, cf.~\cite[Section 3]{BhattScholzeProetale}) and the following proposition that is an amelioration of \cite[Proposition 14.7]{ECoD} for solid sheaves.

\begin{proposition}\label{prop:solidproetv} Let $X$ be a diamond and let $\mathcal F$ be a sheaf of $\hat{\mathbb Z}$-modules on $X_\qproet$ that is solid. Let $\lambda: X_v\to X_\qproet$ be the map of sites. Then $\mathcal F\to R\lambda_\ast \lambda^\ast \mathcal F$ is an isomorphism.
\end{proposition}

\begin{proof} We may assume that $X$ is spatial (or strictly totally disconnected). Then $\mathcal F$ is a filtered colimit of finitely presented solid $\hat{\mathbb Z}$-sheaves, and the functor $R\lambda_\ast$ commutes with filtered colimits in $D^{\geq 0}$. We may thus assume that $\mathcal F$ is finitely presented; in that case $\mathcal F$ is a cofiltered limit of torsion constructible \'etale sheaves, and $\lambda^\ast$ commutes with all limits by \cite[Lemma 14.4]{ECoD}. Thus, we can assume that $\mathcal F$ is an \'etale sheaf, where the claim is \cite[Proposition 14.7]{ECoD}.
\end{proof}

Moreover, solid objects in the derived category satisfy a derived and internal version of Definition~\ref{def:solid}.

\begin{proposition}\label{prop:derivedsolid} Let $X$ be a spatial diamond. For all $A\in D_\solid(X,\hat{\mathbb Z})$, the map
\[
R\sHom(j_\natural \hat{\mathbb Z},A)\to Rj_\ast A|_U
\]
is an isomorphism for all quasi-pro-\'etale $j: U\to X$.
\end{proposition}

\begin{proof} By taking a Postnikov limit, we can assume that $A\in D^+_\solid(X,\hat{\mathbb Z})$, and then one reduces to the case that $A=\mathcal F[0]$ is concentrated in degree $0$. Now by a resolution of Breen \cite[Section 3]{BreenAdditif} (appropriately sheafified), there is a resolution of any $\hat{\mathbb Z}$-sheaf $\mathcal G$ where all terms are finite direct sums of sheaves of the form $\hat{\mathbb Z}[\mathcal G^i\times \hat{\mathbb Z}^j]$. If $\mathcal G$ is a spatial diamond, then all $\mathcal G^i\times \hat{\mathbb Z}^j$ are spatial diamonds, hence
\[
R\sHom(\hat{\mathbb Z}[\mathcal G^i\times \hat{\mathbb Z}^j],\mathcal F)
\]
commutes with all filtered colimits. Applied to $\mathcal G=j_\natural \hat{\mathbb Z}$, Breen's resolution then implies that $R\sHom(\mathcal G,-)$ commutes with all filtered colimits.

We may thus assume that $\mathcal F$ is finitely presented. But then Theorem~\ref{thm:solid} implies that $\mathcal F$ is a limit of constructible \'etale sheaves, so one can reduce to the case that $\mathcal F$ is an \'etale sheaf. But then Breen's resolution shows that
\[
R\sHom(j_\natural \hat{\mathbb Z},\mathcal F) = \varinjlim_i R\sHom(j_{i!} \hat{\mathbb Z},\mathcal F) = \varinjlim_i Rj_{i\ast} \mathcal F|_{U_i}
\]
and \cite[Proposition 14.9]{ECoD} shows that this identifies with $Rj_\ast \mathcal F|_U$.
\end{proof}

\begin{proposition}\label{prop:leftadjointsolid} Let $X$ be a spatial diamond. The inclusion
\[
D_\solid(X,\hat{\mathbb Z})\subset D(X_\qproet,\hat{\mathbb Z})
\]
admits a left adjoint
\[
A\mapsto A^\solid: D(X_\qproet,\hat{\mathbb Z})\to D_\solid(X,\hat{\mathbb Z}).
\]
Moreover, $D_\solid(X,\hat{\mathbb Z})$ identifies with the derived category of solid $\hat{\mathbb Z}$-sheaves on $X$, and $A\mapsto A^\solid$ with the left derived functor of $\mathcal F\mapsto \mathcal F^\solid$. The formation of $A\mapsto A^\solid$, for $A\in D(X_\qproet,\hat{\mathbb Z})$, commutes with any base change $X'\to X$ of spatial diamonds.
\end{proposition}

\begin{proof} This follows easily from Proposition~\ref{prop:derivedsolid}.
\end{proof}

\begin{proposition}\label{prop:tensorproductsolid} Let $X$ be a spatial diamond. The kernel of $A\mapsto A^\solid$ is a tensor ideal. In particular, there is a unique symmetric monoidal structure $-\soliddotimes -$ on $D_\solid(X,\hat{\mathbb Z})$ making $A\mapsto A^\solid$ symmetric monoidal. It is the left derived functor of the induced symmetric monoidal structure on solid $\hat{\mathbb Z}$-sheaves. This symmetric monoidal structure commutes with all colimits (in each variable) and any pullback.
\end{proposition}

\begin{proof} To check that the kernel is a tensor ideal, take any quasi-pro-\'etale $j: U\to X$ written as a cofiltered inverse limit of separated \'etale $j_i: U_i\to X$, and any further quasi-pro-\'etale $j': U'\to X$. Then for any solid $A\in D_\solid(X,\hat{\mathbb Z})$, we know by Proposition~\ref{prop:derivedsolid} that the map
\[
R\sHom(j_\natural \hat{\mathbb Z},A)\to Rj_\ast A|_U
\]
is an isomorphism. Taking sections over $U'\to X$, this translates into
\[
R\Hom(j_\natural \hat{\mathbb Z}\dotimes_{\hat{\mathbb Z}} \hat{\mathbb Z}[U'],A)\to R\Hom(\hat{\mathbb Z}[U\times_X U'],A)
\]
being an isomorphism. In other words, taking the tensor product of $\hat{\mathbb Z}[U]\to j_\natural \hat{\mathbb Z}$ with $\hat{\mathbb Z}[U']$ still lies in the kernel, but these generate the tensor ideal generated by the kernel.

It is now formal that there is a unique symmetric monoidal structure $-\soliddotimes-$ on $D_\solid(X,\hat{\mathbb Z})$ making $A\mapsto A^\solid$ symmetric monoidal (given by the solidification of the tensor product in all solid pro-\'etale sheaves). As solidification commutes with all colimits, so does this tensor product. On generators $j: U\to X$, $j': U'\to X$ as above, it is given by $j_\natural \hat{\mathbb Z}\soliddotimes j'_\natural \hat{\mathbb Z} = (j\times_X j')_\natural \hat{\mathbb Z}$, which still sits in degree $0$; this implies that the functor is a left derived functor. Moreover, this description commutes with any base change.
\end{proof}

Moreover, the inclusion into all v-sheaves also admits a left adjoint, if $X$ is a diamond. We will later improve on this proposition when working with $\hat{\mathbb Z}^p$-coefficients.

\begin{proposition}\label{prop:vsolidification} For any diamond $X$, the fully faithful embedding
\[
D_\solid(X,\hat{\mathbb Z})\subset D(X_v,\hat{\mathbb Z})
\]
admits a left adjoint $A\mapsto A^\solid$. The formation of $A^\solid$ commutes with quasi-pro-\'etale base change $X'\to X$.
\end{proposition}

\begin{proof} Assume first that $X$ is strictly totally disconnected. It suffices to construct the left adjoint on a set of generators, such as the pro-\'etale sheaves of $\hat{\mathbb Z}$-modules generated by some strictly totally disconnected $Y\to X$. By \cite[Lemma 14.5]{ECoD}, there is a strictly totally disconnected affinoid pro-\'etale $j: Y'\to X$ such that $Y\to X$ factors over a map $Y\to Y'$ that is surjective and induces a bijection of connected components. Then for any $B\in D_\solid(Y',\hat{\mathbb Z})$, the map
\[
R\Gamma(Y',B)\to R\Gamma(Y,B)
\]
is an isomorphism. Indeed, by Postnikov limits this easily reduces to $B=\mathcal F[0]$ for a solid sheaf of $\hat{\mathbb Z}$-modules, and then to a finitely presented solid sheaf, and finally to a constructible \'etale sheaf, for which the result is proved at the end of the proof of \cite[Lemma 14.4]{ECoD}. This means that the left adjoint $A\mapsto A^\solid$ when evaluated on $\hat{\mathbb Z}[Y]$ exists and is given by $j_\natural \hat{\mathbb Z}$.

The formation of $Y'\to X$ from $Y\to X$ commutes with any quasi-pro-\'etale base change of strictly totally disconnected $X'\to X$. This implies that $A\mapsto A^\solid$ commutes with such base changes. By descent, this implies the existence of the left adjoint in general, and its commutation with quasi-pro-\'etale base change.
\end{proof}

As usual, we also want to have a theory with coefficients in a ring $\Lambda$. As before, we assume that $\Lambda$ is constant in the sense that it comes via pullback from the point. In our case, this means that it comes via pullback from the pro-\'etale site of a point, i.e.~is a condensed ring \cite{Condensed}, and we need to assume that it is solid over $\hat{\mathbb Z}$; in other words, we allow as coefficients any solid $\hat{\mathbb Z}$-algebra $\Lambda$. Via pullback, this gives rise to a v-sheaf of $\hat{\mathbb Z}$-algebras on any small v-stack $X$, and we can consider $D(X_v,\Lambda)$.

\begin{example} We may consider $\Lambda = \overline{\mathbb{Z}}_\ell$ as the solid condensed ring $\varinjlim_{L|\mathbb{Q}_\ell\text{ finite}} \O_L$. 
\end{example}

\begin{definition}\label{def:solidcoefficients} Let $D_\solid(X,\Lambda)\subset D(X_v,\Lambda)$ be the full subcategory of all $A\in D(X_v,\Lambda)$ such that the image of $A$ in $D(X_v,\hat{\mathbb Z})$ is solid.
\end{definition}

On the level of $\infty$-categorical enrichments, we thus see that $\mathcal D_\solid(X,\Lambda)$ is the category of $\underline{\Lambda}$-modules in $\mathcal D_\solid(X,\hat{\mathbb Z})$. It is then formal that the inclusion $D_\solid(X,\Lambda)\subset D(X_v,\Lambda)$ admits a symmetric monoidal left adjoint $A\mapsto A^\solid$, compatible with forgetting the $\Lambda$-structure.

\begin{remark} Let us briefly compare the present theory with the one developed in \cite{Condensed}. Over a geometric point $X=\Spa C$, $D_\solid(X,\hat{\mathbb Z})$ is the derived category of solid $\hat{\mathbb Z}$-modules in the sense of \cite{Condensed}. For general $\Lambda$, we are now simply considering $\Lambda$-modules in $\mathcal D_\solid(X,\hat{\mathbb Z})$. This is in general different from the theory of $\Lambda_\solid$-modules, which would ask for a stronger completeness notion relative to $\Lambda$. Our present theory corresponds to the analytic ring structure on $\Lambda$ induced from $\hat{\mathbb Z}_\solid$.

One might wonder whether for any analytic ring $\mathcal A$ in the sense of \cite{Condensed} one can define a category $D(X,\mathcal A)$ of ``$\mathcal A$-complete'' pro-\'etale sheaves on any spatial diamond $X$. This does not seem to be the case; it is certainly not formal. In fact, already for $\mathcal A=\mathbb Z_\solid$, problems occur and there is certainly no abelian category; it is still possible to define a nice derived category, though. For general $\mathcal A$, defining $D(X,\mathcal A)$ also seems to require extra data beyond the analytic ring structure on $\mathcal A$.
\end{remark}

\section{Four functors}

Now we discuss some functors on solid sheaves. For this, {\it we assume from now on that we work with coefficients $\Lambda$ given by a solid $\hat{\mathbb Z}^p$-algebra} (so we stay away from $p$-adic coefficients). For any map $f: Y\to X$ of small v-stacks, we have the pullback functor $f^\ast: D_\solid(X,\Lambda)\subset D_\solid(Y,\Lambda)$. This admits a right adjoint $Rf_\ast$; in fact, one can simply import $Rf_\ast$ from the full $D(Y,\Lambda)$:

\begin{proposition}\label{prop:pushforwardsolid} Let $f: Y\to X$ be a map of small v-stacks and let $A\in D_\solid(Y,\Lambda)\subset D(Y,\Lambda)$. Then $Rf_{v\ast} A\in D(X_v,\Lambda)$ lies in $D_\solid(X,\Lambda)$. In particular, $Rf_{v\ast}: D(Y_v,\Lambda)\to D(X_v,\Lambda)$ restricts to a functor $Rf_\ast: D_\solid(Y,\Lambda)\to D_\solid(X,\Lambda)$ that is right adjoint to $f^\ast$.
\end{proposition}

\begin{proof} We can formally reduce to the case $\Lambda=\hat{\mathbb Z}^p$. The formation of $Rf_{v\ast}$ commutes with any pullback (as everything is a slice in the v-site), so using Proposition \ref{prop:solidfunctor} we can assume that $X$ is a spatial diamond. Moreover, taking a simplicial resolution of $Y$ by disjoint unions of spatial diamonds, and using that $D_\solid(X,\hat{\mathbb Z}^p)\subset D(X_v,\hat{\mathbb Z}^p)$ is stable under all derived limits (as it is stable under all products), we can also assume that $Y$ is a spatial diamond.

We may assume $A\in D^+_\solid(Y,\hat{\mathbb Z}^p)$ by a Postnikov limit, then that $A=\mathcal F[0]$ is concentrated in degree $0$, then that $\mathcal F$ is finitely presented by writing it as a filtered colimit, and finally that $\mathcal F$ is a constructible \'etale sheaf by writing it as a cofiltered limit. Now the result follows from \cite[Proposition 17.6]{ECoD}.
\end{proof}

\begin{proposition}\label{prop:leftadjointvsolid} For any small v-stack $X$, the inclusion
\[
D_\solid(X,\Lambda)\subset D(X_v,\Lambda)
\]
admits a left adjoint
\[
A\mapsto A^\solid: D(X_v,\Lambda)\to D_\solid(X,\Lambda).
\]
The functor $A\mapsto A^\solid$ commutes with any base change.

The kernel of $A\mapsto A^\solid$ is a tensor ideal. In particular, there is a unique symmetric monoidal structure $-\soliddotimesLambda -$ on $D_\solid(X,\Lambda)$ making $A\mapsto A^\solid$ a symmetric monoidal functor. The functor $-\soliddotimesLambda -$ commutes with all colimits (in each variable) and with all pullbacks $f: Y\to X$.
\end{proposition}

We note that in the case of overlap with previous definitions of $A\mapsto A^\solid$ and $-\soliddotimesLambda-$, the definitions agree, by uniqueness of the previous definitions.

\begin{proof} Again, one can formally reduce to the case $\Lambda=\hat{\mathbb Z}^p$. By descent, we can reduce to the case that $X$ is strictly totally disconnected. (Note that $Y\mapsto \mathcal D_\solid(Y,\hat{\mathbb Z}^p)$ is a v-sheaf of $\infty$-categories --- this is clear for $\mathcal D(Y_v,\hat{\mathbb Z}^p)$, and follows for $\mathcal D_\solid$ as being solid can be checked v-locally by Proposition~\ref{prop:solidfunctor}. Thus, if the left adjoints exist v-locally and commute with base change, they assemble into the desired left adjoint, cf.~\cite[Proposition 4.7.4.19]{LurieHA}.) In this case, we already know existence of the left adjoint $A\mapsto A^\solid$ by Proposition~\ref{prop:vsolidification}.

We check that the left adjoint $A\mapsto A^\solid$ commutes with any base change $f: Y\to X$. We already know that pullbacks of solid objects stay solid, so we have to see that if $A\in D(X_v,\hat{\mathbb Z}^p)$ satisfies $A^\solid=0$, then also $(f_v^\ast A)^\solid=0$. But this statement is adjoint to the statement that $Rf_{v\ast}$ preserves $D_\solid$, i.e.~Proposition~\ref{prop:pushforwardsolid}.

We need to see that the class of all $A\in D(X_v,\hat{\mathbb Z}^p)$ with $A^\solid=0$ is a $\otimes$-ideal. But we have seen that for all $f: Y\to X$, also $f_v^\ast A$ lies in the corresponding class for $Y$, and then so does $f_{v\natural} f_v^\ast A$ (as pullback preserves $D_\solid$), where we write $f_{v\natural}$ for the left adjoint of $f_v^\ast$ (which exists as it is a slice). But $f_{v\natural} f_v^\ast A = A\dotimes_{\hat{\mathbb Z}^p} f_{v\natural} \hat{\mathbb Z}^p$ by the projection formula for slices, so this gives the desired claim.
\end{proof}

It turns out that for $\Lambda=\hat{\mathbb Z}^p$, the functor $-\soliddotimes -$ is actually almost exact. If one would work with $\Lambda=\mathbb F_\ell$-coefficients, it would even be exact.

\begin{proposition}\label{prop:solidtensorexact} Let $X$ be a small v-stack and $A,B\in D_\solid(X,\hat{\mathbb Z}^p)$ be concentrated in degree $0$. Then $A\soliddotimes B$ sits in cohomological degrees $-1$ and $0$.

If $X$ is a spatial diamond and $\mathcal F=\varprojlim_i \mathcal F_i$ and $\mathcal G=\varprojlim_j \mathcal G_j$ are finitely presented solid $\hat{\mathbb Z}^p$-sheaves written as cofiltered limits of constructible \'etale sheaves killed by some integer prime to $p$, then the natural map
\[
\mathcal F\soliddotimes \mathcal G\to R\varprojlim_{i,j} \mathcal F_i\dotimes \mathcal G_j
\]
is an isomorphism.
\end{proposition}

\begin{proof} It suffices to prove the final assertion, as the statement on $A\soliddotimes B$ can be checked after pullback to spatial diamonds, and then $A$ and $B$ can be written as filtered colimits of finitely presented solid $\hat{\mathbb Z}^p$-sheaves (and $\mathcal F_i\dotimes \mathcal G_j$ sits in degrees $-1$ and $0$ as $\hat{\mathbb Z}^p$ has global dimension $1$). Resolving $\mathcal F$ and $\mathcal G$, we can reduce to the case $\mathcal F=j_\natural \hat{\mathbb Z}^p$, $\mathcal G=j'_\natural \hat{\mathbb Z}^p$. But their solid tensor product is indeed given by $(j\times_X j')_\natural \hat{\mathbb Z}^p$.
\end{proof}

At this point, we have defined $D_\solid(X,\Lambda)\subset D(X_v,\Lambda)$ for any small v-stack $X$, and this subcategory is preserved by pullback and pushforward, and in particular this gives such functors for $D_\solid(X,\Lambda)$. Moreover, $D_\solid(X,\Lambda)$ has a natural symmetric monoidal structure $-\soliddotimesLambda-$, commuting with colimits in both variables, and with pullbacks. Moreover, we have a functor
\[
R\sHom_\Lambda(-,-): D_\solid(X,\Lambda)^{\mathrm{op}}\times D_\solid(X,\Lambda)\to D_\solid(X,\Lambda),
\]
a partial right adjoint to $-\soliddotimesLambda -$ as usual. Again, it can be obtained from the corresponding functor on $D(X_v,\Lambda)$ via restriction. In fact, for all $A\in D(X_v,\Lambda)$ and $B\in D_\solid(X,\Lambda)$, one has $R\sHom_\Lambda(A,B)\in D_\solid(X,\Lambda)$. This can be reduced to $\Lambda=\hat{\mathbb Z}^p$ and the case $A=f_\natural \hat{\mathbb Z}^p$ for some $f: Y\to X$, and then it amounts to $Rf_{v\ast} f_v^\ast B\in D_\solid(X,\hat{\mathbb Z}^p)$, which follows from Proposition~\ref{prop:pushforwardsolid}.

There is the following general base change result. We stress the absence of \emph{any} conditions.

\begin{proposition}\label{prop:solidbasechange} Let
\[\xymatrix{
Y'\ar[r]^{g'}\ar[d]^{f'} & Y\ar[d]^f\\
X'\ar[r]^g & X
}\]
be a cartesian diagram of small v-stacks. For all $A\in D_\solid(Y,\Lambda)$, the base change map
\[
g^\ast Rf_\ast A\to Rf'_\ast g^{\prime\ast} A
\]
in $D_\solid(X',\Lambda)$ is an isomorphism.

Similarly, for any map $f: Y\to X$ of small v-stacks and all $A,B\in D_\solid(X,\Lambda)$, the map
\[
f^\ast R\sHom(A,B)\to R\sHom(f^\ast A,f^\ast B)
\]
in $D_\solid(Y,\Lambda)$ is an isomorphism.
\end{proposition}

\begin{proof} The base change is a direct consequence of Proposition~\ref{prop:pushforwardsolid}, noting that in the v-site, everything is a slice (and hence satisfies base change). The statement about $R\sHom$ follows similarly from the compatibility with the $R\sHom$ as formed on the v-site, as was noted above.
\end{proof}

The projection formula, however, fails to hold.

\begin{warning}\label{war:solidprojformula} If $f: Y\to X$ is a proper map of small v-stacks that is representable in spatial diamonds with $\dimtrg f<\infty$, the map
\[
A\soliddotimes  Rf_\ast B\to Rf_\ast (f^\ast A\soliddotimes B)
\]
may fail to be an isomorphism for $A\in D_\solid(X,\hat{\mathbb Z}^p)$ and $B\in D_\solid(Y,\hat{\mathbb Z}^p)$. In fact, already if $X=\mathbb B_C$ is a perfectoid ball and $f=j: Y=\Spa C\to X$ is the inclusion of a point (which is quasi-pro-\'etale), then this fails for $A=j_\natural \hat{\mathbb Z}^p$ and $B=\hat{\mathbb Z}^p$. In fact, the map becomes $j_\natural \hat{\mathbb Z}^p\to Rj_\ast \hat{\mathbb Z}^p$, which is far from an isomorphism: For example, on global sections the left-hand side becomes $\hat{\mathbb Z}^p[-2]$, while the right-hand side becomes $\hat{\mathbb Z}^p$.
\end{warning}

There is the following result on change of algebraically closed base field, an analogue of \cite[Theorem 19.5]{ECoD}.

\begin{proposition}\label{prop:solidchangeofbasefield} Let $X$ be a small v-stack.
\begin{altenumerate}
\item[{\rm (i)}] Assume that $X$ lives over $k$, where $k$ is a discrete algebraically closed field of characteristic $p$, and $k'/k$ is an extension of discrete algebraically closed base fields, $X'=X\times_k k'$. Then the pullback functor
\[
D_\solid(X,\Lambda)\to D_\solid(X',\Lambda)
\]
is fully faithful.
\item[{\rm (ii)}] Assume that $X$ lives over $k$, where $k$ is a discrete algebraically closed field of characteristic $p$. Let $C/k$ be an algebraically closed complete nonarchimedean field, and $X'=X\times_k \Spa(C,C^+)$ for some open and bounded valuation subring $C^+\subset C$ containing $k$. Then the pullback functor
\[
D_\solid(X,\Lambda)\to D_\solid(X',\Lambda)
\]
is fully faithful.
\item[{\rm (iii)}] Assume that $X$ lives over $\Spa(C,C^+)$, where $C$ is an algebraically closed complete nonarchimedean field with an open and bounded valuation subring $C^+\subset C$, $C'/C$ is an extension of algebraically closed complete nonarchimedean fields, and $C^{\prime+}\subset C'$ an open and bounded valuation subring containing $C^+$, such that $\Spa(C',C^{\prime+})\to \Spa(C,C^+)$ is surjective. Then for $X'=X\times_{\Spa(C,C^+)} \Spa(C',C^{\prime+})$, the pullback functor
\[
D_\solid(X,\Lambda)\to D_\solid(X',\Lambda)
\]
is fully faithful.
\end{altenumerate}
\end{proposition}

\begin{proof} We can assume $\Lambda=\hat{\mathbb Z}^p$. As in \cite[Theorem 19.5]{ECoD}, it suffices to prove (iii) and the restricted case of (ii) where $C$ is the completed algebraic closure of $k\laurentseries{t}$ (and hence $C^+=\mathcal O_C$).

Let $f: X'\to X$ be the map. We have to see that for all $A\in D_\solid(X,\hat{\mathbb Z}^p)$, the map
\[
A\to Rf_\ast f^\ast A
\]
is an equivalence. This can be checked locally in the v-topology, so we can assume that $X=\Spa(R,R^+)$ is an affinoid perfectoid space. By Postnikov limits, we can also assume that $A\in D^+_\solid(X,\hat{\mathbb Z}^p)$, and then that $A$ is concentrated in degree $0$. In case (iii), we can now conclude by writing $A$ as a filtered colimit of finitely presented solid $\hat{\mathbb Z}^p$-modules, and these as cofiltered limits of constructible \'etale sheaves, noting that both operations commute with $Rf_\ast$ and $f^\ast$ (as $f$ is qcqs in case (iii)), and hence reducing us to \cite[Theorem 19.5]{ECoD}.

It remains to handle case (ii) when $C$ is the completed algebraic closure of $k\laurentseries{t}$. In that case $X'$ lives over a punctured open unit disc $\mathbb D^\ast_X$ over $X$, and fixing a pseudouniformizer $\varpi\in R$, this can be written as the increasing union of quasicompact open subspaces $X'_n = \{|t|^n\leq |\varpi|\leq |t|^{1/n}\}\subset X'$, with maps $f_n: X'_n\to X$. It suffices to prove that for all $n$, the map
\[
A\to Rf_{n\ast} f_n^\ast A
\]
is an isomorphism. These functors commute again with filtered colimits of sheaves, and hence the previous reductions apply and reduce the assertion to the \'etale case, which was handled in the proof of \cite[Theorem 19.5]{ECoD}.
\end{proof}

As an application, let us record the following versions of Proposition~\ref{prop:drinfeldfullyfaithful} and Corollary~\ref{cor:drinfeldfullyfaithful}, where we fix an algebraically closed field $k|\Fq$ and work on $\Perf_k$.

\begin{corollary}\label{cor:drinfeldfullyfaithfulsolid0} For any small v-stack $X$, the functor
\[
\psi_X^\ast: D_\solid(X\times [\ast/\underline{W_E}],\Lambda)\to D_\solid(X\times \Div^1,\Lambda)
\]
is fully faithful. If the natural pullback functor
\[
D_\solid(X,\Lambda)\to D_\solid(X\times \Spd \widehat{\overline{E}},\Lambda)
\]
is an equivalence, then $\psi_X^\ast$ is also an equivalence.
\end{corollary}

\begin{proof} By descent along $X\to X\times [\ast/\underline{W_E}]$ this reduces to Proposition~\ref{prop:solidchangeofbasefield}.
\end{proof}

\begin{corollary}\label{cor:drinfeldfullyfaithfulsolid} For any small v-stack $X$ and finite set $I$, pullback along $X\times (\Div^1)^I\to X\times [\ast/\underline{W_E^I}]$ induces a fully faithful functor
\[
D_\solid(X\times [\ast/\underline{W_E^I}],\Lambda)\to D_\solid(X\times (\Div^1)^I,\Lambda).
\]
\end{corollary}

\begin{proof} This follows inductively from Corollary~\ref{cor:drinfeldfullyfaithfulsolid0}.
\end{proof}

We also need a solid analogue of Theorem~\ref{thm:partialcompactsupportvanishing}; we only prove a restricted variant, however. As there, work over $\Perf_k$, and let $X$ be a spatial diamond such that $X\to \ast$ is proper, of finite $\dimtrg$, and take any spatial diamond $S$. As before, one can introduce the doubly-indexed ind-system $\{U_{a,b}\}_{(a,b)}\subset X\times S$, well-defined up to ind-isomorphism; and then $U_a=\bigcup_{b<\infty} U_{a,b}$ and $U_b=\bigcup_{a>0} U_{a,b}$.

\begin{definition} The functors
\[
R\beta_{!+},R\beta_{!-}: D_\solid(X\times S,\Lambda)\to D_\solid(S,\Lambda)
\]
are defined by
\[\begin{aligned}
R\beta_{!+}C &:= \varinjlim_a R\beta_{\ast}(j_{a!} C|_{U_a}),\\
R\beta_{!-}C &:= \varinjlim_b R\beta_{\ast}(j_{b!} C|_{U_b})\\
\end{aligned}\]
for $C\in D_\solid(X\times S,\Lambda)$.
\end{definition}

Here $j_{a!}$ and $j_{b!}$ denote the left adjoints to $j_a^\ast$ and $j_b^\ast$. Let $\alpha: X\times S\to X$ be the projection.

\begin{theorem}\label{thm:partialcompactsupportvanishingsolid} Assume that $C=\alpha^\ast A$ for $A\in D_\solid(X,\Lambda)$, and assume that either $A\in D_\solid^+(X,\Lambda)$, or that $X\to \ast$ is cohomologically smooth. Then
\[
R\beta_{!+}C = 0 = R\beta_{!-} C.
\]
\end{theorem}

\begin{proof} We can assume $\Lambda=\hat{\mathbb Z}^p$. All operations commute with any base change; we can thus assume that $S=\Spa K$ where $K$ is the complete algebraic closure of $k\laurentseries{t}$. We observe that if $X\to \ast$ is cohomologically smooth, then $R\beta_\ast: D_\solid(X\times S,\Lambda)\to D_\solid(S,\Lambda)$ has finite cohomological dimension; this is a statement about sheaves concentrated in degree $0$. Any such $B$ can be written as the countable limit of $Rj_{a,b,\ast} j_{a,b}^\ast B$ for the open immersions $j_{a,b}: U_{a,b}\subset X\times S$; it is thus enough to show that pushforward along $U_{a,b}\to S$ has finite cohomological dimension on solid sheaves. As $U_{a,b}\to S$ is qcqs, we can reduce to finitely presented sheaves; these are cofiltered limits of constructible sheaves. For constructible sheaves, the cohomological dimension is bounded, and each cohomology group (recall that $S=\Spa K$ is a geometric point) is finite by \cite[Theorem 25.1]{ECoD}. Thus, the cofiltered limit stays in the same range of degrees.

It follows that we can assume that $A\in D_\solid^+(X,\hat{\mathbb Z}^p)$. Arguing as in the proof of Theorem~\ref{thm:partialcompactsupportvanishing}, we can then reduce to the case that $X=\Spa(R,R^+)$ is an affinoid perfectoid space with no nonsplit finite \'etale covers, and then to $X=\Spa K$ where $K$ is still the completed algebraic closure of $k\laurentseries{t}$. In that case, as in the proof of Theorem~\ref{thm:partialcompactsupportvanishing}, one can make a more precise assertion on actual annuli; this statement is compatible with passage to filtered colimits, reducing us to the case that $A$ is a finitely presented solid sheaf. For $A\in D_\solid(\Spa K,\hat{\mathbb Z}^p)$, this means that $A$ is a cofiltered limit of finite abelian groups killed by integers prime to $p$. We can also pull this cofiltered limit through, reducing us to Theorem~\ref{thm:partialcompactsupportvanishing}.
\end{proof}

\section{Relative homology}

A unique feature of the formalism of solid sheaves is the existence of a general left adjoint to pullback, with excellent properties. We continue to work with coefficients in a solid $\hat{\mathbb Z}^p$-algebra $\Lambda$.

\begin{proposition}\label{prop:fnatural} Let $f: Y\to X$ be any map of small v-stacks.
\begin{altenumerate}
\item[{\rm (i)}] The functor $f^\ast: D_\solid(X,\Lambda)\to D_\solid(Y,\Lambda)$ admits a left adjoint
\[
f_\natural: D_\solid(Y,\Lambda)\to D_\solid(X,\Lambda).
\]
The natural map
\[
f_\natural(A\soliddotimesLambda f^\ast B)\to f_\natural A\soliddotimesLambda B
\]
is an isomorphism for all $A\in D_\solid(Y,\Lambda)$ and $B\in D_\solid(X,\Lambda)$. Similarly, the map
\[
R\sHom(f_\natural A,B)\to Rf_\ast R\sHom(A,f^\ast B)
\]
is an isomorphism.

\item[{\rm (ii)}] The formation of $f_\natural$ commutes with restriction of coefficients along a map $\Lambda'\to \Lambda$.

\item[{\rm (iii)}] For any cartesian diagram
\[\xymatrix{
Y'\ar[r]^{g'}\ar[d]^{f'} & Y\ar[d]^f\\
X'\ar[r]^g & X
}\]
of small v-stacks, the natural map
\[
f'_\natural g^{\prime\ast} A\to g^\ast f_\natural A
\]
is an isomorphism for all $A\in D_\solid(Y,\Lambda)$.
\end{altenumerate}
\end{proposition}

\begin{proof} As $f$ is a slice in the v-site, it is tautological that $f_v^\ast: D(X_v,\Lambda)\to D(Y_v,\Lambda)$ admits a left adjoint $f_{v\natural}: D(Y_v,\Lambda)\to D(X_v,\Lambda)$. One can then define $f_\natural$ as the solidification of $f_{v\natural}$. By general properties of slices, the map
\[
f_{v\natural}(A\dotimes_\Lambda f^\ast B)\to f_{v\natural} A\dotimes_\Lambda B
\]
is an isomorphism. Passing to solidifications, using that that this is symmetric monoidal, then gives that
\[
f_\natural(A\soliddotimesLambda f^\ast B)\to f_\natural A\soliddotimesLambda B
\]
is an isomorphism. The isomorphism
\[
R\sHom(f_\natural A,B)\cong Rf_\ast R\sHom(A,f^\ast B)
\]
then follows by adjointness.

For part (ii), we can assume $\Lambda'=\hat{\mathbb Z}^p$, and check on generators. These are given by $j_\natural \hat{\mathbb Z}^p\dotimes_{\hat{\mathbb Z}^p} \Lambda$ for $j: Y'\to Y$. The claim then follows from the projection formula.

Part (iii) is obtained by passing to left adjoints in Proposition~\ref{prop:solidbasechange}.
\end{proof}

Now if $f$ is ``proper and smooth'', we want to relate the left adjoint $f_\natural$ (``homology'') and the right adjoint $Rf_\ast$ (``cohomology''). Thus, assume that $f: Y\to X$ is a proper map of small v-stacks that is representable in spatial diamonds with $\dimtrg f<\infty$, and cohomologically smooth, i.e.~$\ell$-cohomologically smooth for all $\ell\neq p$ (or just all $\ell$ relevant for $\Lambda$). In this case, we want to express $Rf_\ast$ in terms of $f_\natural$. As a first step, we show that $Rf_\ast$ has bounded cohomological dimension.

\begin{proposition}\label{prop:ULAboundedcohomdim} Let $f: Y\to X$ be a proper map of small v-stacks that is representable in spatial diamonds with $\dimtrg f<\infty$, and cohomologically smooth. Then $Rf_\ast: D_\solid(Y,\Lambda)\to D_\solid(X,\Lambda)$ has bounded cohomological dimension and commutes with arbitrary direct sums. If $X$ is a spatial diamond (thus $Y$ is) and $\mathcal F$ is a finitely presented solid $\hat{\mathbb Z}^p$-sheaf on $Y$, then $Rf_\ast \mathcal F$ is a bounded complex all of whose cohomology sheaves are finitely presented solid $\hat{\mathbb Z}^p$-sheaves on $X$.
\end{proposition}

\begin{proof} We can assume $\Lambda=\hat{\mathbb Z}^p$. The commutation with arbitrary direct sums follows from bounded cohomological dimension, as one can then reduce to the case of complexes concentrated in degree $0$, where $Rf_\ast$ commutes with all direct sums as $f$ is qcqs. For the claim about bounded cohomological dimension, we can argue v-locally, and hence assume that $X$ is a spatial diamond. It suffices to prove that for all solid $\hat{\mathbb Z}^p$-sheaves $\mathcal F$ on $Y$, the complex $Rf_\ast \mathcal F$ is bounded; this reduces to the case of finitely presented solid $\hat{\mathbb Z}^p$-sheaves as $Rf_\ast$ commutes with filtered colimits of sheaves. Now if $\mathcal F$ is finitely presented, it is a cofiltered limit of constructible \'etale sheaves killed by some integer prime to $p$. As $Rf_\ast$ commutes with this limit, it is now enough to see that $Rf_\ast$ preserves constructible complexes and has bounded amplitude. But this follows from cohomological smoothness, cf.~\cite[Proposition 23.12 (ii)]{ECoD}.
\end{proof}

Next, we prove a projection formula for $Rf_\ast$.

\begin{proposition}\label{prop:ULAprojformula} Let $f: Y\to X$ be a proper map of small v-stacks that is representable in spatial diamonds with $\dimtrg f<\infty$, and cohomologically smooth. Then for all $A\in D_\solid(Y,\Lambda)$ and $B\in D_\solid(X,\Lambda)$, the projection map
\[
Rf_\ast A\soliddotimes B\to Rf_\ast(A\soliddotimes f^\ast B)
\]
is an isomorphism.
\end{proposition}

\begin{proof} We can assume $\Lambda=\hat{\mathbb Z}^p$. We note that $Rf_\ast$ and $\soliddotimes$ both have bounded cohomological dimension, so one easily reduces to the case that $A$ and $B$ are concentrated in degree $0$. We can also assume that $X$ is a spatial diamond (thus $Y$ is, too). Then we can write $A$ and $B$ as filtered colimits of finitely presented solid $\hat{\mathbb Z}^p$-sheaves, and reduce to the case that $A$ and $B$ are cofiltered limits of constructible \'etale sheaves killed by some integer prime to $p$. In that case, it follows from Proposition~\ref{prop:solidtensorexact} and Proposition~\ref{prop:ULAboundedcohomdim} that all operations commute with these cofiltered limits, and one reduces to the case that $A$ and $B$ are constructible \'etale sheaves killed by some integer prime to $p$. Now it follows from \cite[Proposition 22.11]{ECoD}.
\end{proof}

Moreover, the functor $Rf_\ast$ interacts well with $g_\natural$ for maps $g: X'\to X$.

\begin{proposition}\label{prop:ULApushforwardhomology} Let
\[\xymatrix{
Y'\ar[r]^{g'}\ar[d]^{f'} & Y\ar[d]^f\\
X'\ar[r]^g & X
}\]
be a cartesian diagram of small v-stacks, where $f: Y\to X$ is proper, representable in spatial diamonds, with $\dimtrg f<\infty$ and cohomologically smooth. Then the natural transformation
\[
g_\natural Rf'_\ast A\to Rf_\ast g'_\natural A
\]
is an isomorphism for all $A\in D_\solid(Y',\Lambda)$.
\end{proposition}

\begin{proof} We can assume $\Lambda=\hat{\mathbb Z}^p$. By Proposition~\ref{prop:ULAboundedcohomdim} both sides commute with Postnikov limits, so we can assume $A\in D^+$, and then reduce to the case that $A$ is concentrated in degree $0$. We may assume that $X$ is a spatial diamond, and one can also reduce to the case $X'$ is a spatial diamond, by writing $A$ as the geometric realization of $h'_{\bullet\natural} h_\bullet^{\prime\ast} A$ for some simplicial hypercover $h_\bullet: X'_\bullet\to X'$ by disjoint unions of spatial diamonds, and its pullback $h'_\bullet: Y'_\bullet\to Y'$ (and using Proposition~\ref{prop:ULAboundedcohomdim} to commute the geometric realization with pushforward). Under these circumstances, one can write $A$ as a filtered colimit of finitely presented solid $\hat{\mathbb Z}^p$-modules, and hence reduce to the case that $A$ is a cofiltered limit of constructible \'etale sheaves killed by some integer prime to $p$. By Proposition~\ref{prop:ULAboundedcohomdim} the complex $Rf'_\ast A$ is then bounded with all cohomology sheaves finitely presented solid $\hat{\mathbb Z}^p$-modules. As $g_\natural$ preserves pseudocoherent objects, it follows that the map $g_\natural Rf'_\ast A\to Rf_\ast g'_\natural A$ is a map of bounded to the right complexes in $D_\solid(X,\hat{\mathbb Z}^p)$ all of whose cohomology sheaves are finitely presented solid $\hat{\mathbb Z}^p$-modules. If the cone of this map is nonzero, then by looking at its first nonzero cohomology sheaf, we find some nonzero map to a constructible \'etale sheaf $B$ on $X$, killed by some integer prime to $p$. Note that, using the usual \'etale $Rf^!$ functor, there is a natural adjunction
\[
R\Hom(Rf_\ast g'_\natural A,B)\cong R\Hom(g'_\natural A,Rf^! B):
\]
Indeed, it suffices to check this when $g'_\natural A$ is replaced by a finitely presented solid $\hat{\mathbb Z}^p$-module, by a Postnikov tower (and as all cohomology sheaves of $g'_\natural A$ are of this form). Writing this as a cofiltered limit of constructible \'etale sheaves killed by some integer prime to $p$, both sides turn this cofiltered limit into a filtered colimit, so the claim reduces to the usual \'etale adjunction.

Applying $R\Hom(-,B)$ to the map $g_\natural Rf'_\ast A\to Rf_\ast g'_\natural A$ will thus produce $R\Hom(A,-)$ applied to the base change map
\[
Rf'^! g^\ast B\leftarrow g^{\prime\ast} Rf^! B,
\]
which is an isomorphism by \cite[Proposition 23.12 (iii)]{ECoD}.
\end{proof}

Now we can describe the functor $Rf_\ast$. Indeed, consider the diagram
\[\xymatrix{
Y\ar[r]^-{\Delta_f} & Y\times_X Y\ar[r]^-{\pi_2}\ar[d]^{\pi_1}&Y\ar[d]^f\\
& Y\ar[r]^f & X.
}\]
Then, under our assumption that $f: Y\to X$ is a proper map of small v-stacks that is representable in spatial diamonds with $\dimtrg f<\infty$ and cohomologically smooth, we have
\[\begin{aligned}
Rf_\ast A &\cong Rf_\ast \pi_{2\natural} \Delta_{f\natural} A\\
&\cong f_\natural R\pi_{1\ast} \Delta_{f\natural} A\\
&\cong f_\natural R\pi_{1\ast} \Delta_{f\natural} \Delta_f^\ast \pi_1^\ast A\\
&\cong f_\natural R\pi_{1\ast}(\Delta_{f\natural} \Lambda\soliddotimesLambda \pi_1^\ast A)\\
&\cong f_\natural (R\pi_{1\ast} \Delta_{f\natural} \Lambda\soliddotimesLambda A).
\end{aligned}\]
We combine this with the following observation.

\begin{proposition}\label{prop:dualizingcomplexsolid} Let $f: Y\to X$ be a proper map of small v-stacks that is representable in spatial diamonds with $\dimtrg f<\infty$ and cohomologically smooth. Then
\[
R\pi_{1\ast} \Delta_{f\natural}\Lambda\in D_\solid(Y,\Lambda)
\]
is invertible, and its inverse is canonically isomorphic to
\[
Rf^! \Lambda := \varprojlim_n Rf^! \mathbb Z/n\mathbb Z\dotimes_{\hat{\mathbb Z}^p} \Lambda.
\]

Thus, there is a canonical isomorphism
\[
f_\natural A\cong Rf_\ast (A\soliddotimesLambda Rf^! \Lambda): D_\solid(Y,\Lambda)\to D_\solid(X,\Lambda).
\]
\end{proposition}

Thus, we get a somewhat unusual formula for the dualizing complex. We remark that the fibres of $R\pi_{1\ast} \Delta_{f\natural} \hat{\mathbb Z}^p$ are given by the limit of $R\Gamma_c(U,\hat{\mathbb Z}^p)$ over all \'etale neighborhoods $U$ of the given geometric point.

\begin{remark} We see here that an important instance of $Rf^!$ admits an alternative description in terms of $g_\natural$ functors. We are a bit confused about exactly how expressive the present $5$-functor formalism is. So far, we were always able to translate any argument in terms of a $6$-functor formalism into this $5$-functor formalism, although it is often a nontrivial matter and there seems to be no completely general recipe.
\end{remark}

\begin{proof} We can assume $\Lambda=\hat{\mathbb Z}^p$. By the isomorphism $Rf_\ast\cong f_\natural (R\pi_{1\ast} \Delta_{f\natural} \hat{\mathbb Z}^p\soliddotimes A)$, it follows that $Rf_\ast: D_\solid(Y,\hat{\mathbb Z}^p)\to D_\solid(X,\hat{\mathbb Z}^p)$ admits a right adjoint, given by
\[
A\mapsto R\sHom(R\pi_{1\ast} \Delta_{f\natural} \hat{\mathbb Z}^p,f^\ast A).
\]
We claim that this right adjoint maps $D_\et(X,\mathbb Z/n\mathbb Z)$ into $D_\et(Y,\mathbb Z/n\mathbb Z)$ for any $n$ prime to $p$, and thus agrees with the right adjoint $Rf^!$ in that setting. Here, we use the embedding $D_\et\subset D_\solid$, cf.~Section VII.4.1 below. This claim can be checked v-locally, so we can assume that $X$ is a spatial diamond. Then $R\pi_{1\ast}\Delta_{f\natural}\hat{\mathbb Z}^p\in D_\solid(Y,\hat{\mathbb Z}^p)$ is a bounded complex all of whose cohomology sheaves are finitely presented solid, by Proposition~\ref{prop:ULAboundedcohomdim} and as $\Delta_f$ is quasi-pro-\'etale (so $\Delta_{f\natural}\hat{\mathbb Z}^p$ is finitely presented solid). This implies that $R\sHom(R\pi_{1\ast} \Delta_{f\natural} \hat{\mathbb Z}^p,-)$ preserves $D_\et(Y,\mathbb Z/n\mathbb Z)$.

Thus, for any $A\in D_\et(X,\mathbb Z/n\mathbb Z)$, there is a natural isomorphism
\[
Rf^! A\cong R\sHom(R\pi_{1\ast} \Delta_{f\natural} \hat{\mathbb Z}^p,f^\ast A).
\]
Applied with $A=\mathbb Z/n\mathbb Z$, this gives isomorphisms
\[
Rf^! \mathbb Z/n\mathbb Z\cong R\sHom(R\pi_{1\ast} \Delta_{f\natural} \hat{\mathbb Z}^p,\mathbb Z/n\mathbb Z).
\]

It remains to see that $R\pi_{1\ast} \Delta_{f\natural} \hat{\mathbb Z}^p$ is invertible; more precisely, we already get a natural map
\[
R\pi_{1\ast} \Delta_{f\natural} \hat{\mathbb Z}^p\to (Rf^! \hat{\mathbb Z}^p)^{-1}
\]
that we want to prove is an isomorphism. This can again be checked v-locally, so we can assume that $X$ is a spatial diamond. Then $R\pi_{1\ast}\Delta_{f\natural} \hat{\mathbb Z}^p$ is a bounded complex all of whose cohomology sheaves are finitely presented solid; so as in the proof of Proposition~\ref{prop:ULApushforwardhomology}, it is enough to check that one gets isomorphisms after applying $R\sHom(-,B)$ for any $B\in D_\et(Y,\mathbb Z/n\mathbb Z)$. But
\[\begin{aligned}
R\sHom(R\pi_{1\ast} \Delta_{f\natural} \hat{\mathbb Z}^p,B)&\cong R\pi_{1\ast} R\sHom(\Delta_{f\natural} \hat{\mathbb Z}^p,R\pi_1^! B)\\
&\cong \Delta_f^\ast R\pi_1^! B\\
&\cong B\dotimes_{\mathbb Z/n\mathbb Z} Rf^! \mathbb Z/n\mathbb Z,
\end{aligned}\]
giving the result.

The final statement follows formally from the identification of $R\pi_{1\ast} \Delta_{f\natural} \Lambda$ and the discussion leading up to the proposition.
\end{proof}

\section{Relation to $D_\et$}

Assume now that $\Lambda$ is discrete. In particular, also being a $\hat{\mathbb Z}^p$-algebra, we have $n\Lambda=0$ for some $n$ prime to $p$. We wish to understand the relation between $D_\et(X,\Lambda)$ and $D_\solid(X,\Lambda)$, and the functors defined on them.

\subsection{Naive embedding} For any small v-stack $X$, we have a fully faithful embedding
\[
D_\et(X,\Lambda)\hookrightarrow D_\solid(X,\Lambda)
\]
as full subcategories of $D(X_v,\Lambda)$. As usual, the adjoint functor theorem implies that this admits a right adjoint $R_{X\et}: D_\solid(X,\Lambda)\to D_\et(X,\Lambda)$. The full inclusion $D_\et(X,\Lambda)\subset D_\solid(X,\Lambda)$ is symmetric monoidal, and compatible with pullback. Moreover, by \cite[Proposition 17.6]{ECoD}, it also commutes with $Rf_\ast$ if $f: Y\to X$ is qcqs and one restricts to $D^+$; or in general $f$ is qcqs and of finite cohomological dimension. Moreover, one always has
\[
R_{X\et} Rf_\ast\cong Rf_\ast R_{Y\et}.
\]
Similarly, passing to right adjoints in the commutation with tensor products, we also have
\[
R_{X\et} R\sHom_{D_\solid(X,\Lambda)}(A,B)\cong R\sHom_{D_\et(X,\Lambda)}(A,R_{X\et}B)
\]
if $A\in D_\et(X,\Lambda)$ and $B\in D_\solid(X,\Lambda)$. If $A$ is perfect-constructible, then for all $B\in D_\et(X,\Lambda)$, one actually has
\[
R\sHom_{D_\et(X,\Lambda)}(A,B)\cong R\sHom_{D_\solid(X,\Lambda)}(A,B):
\]
by descent, it suffices to check this when $X$ is spatial diamond, and then one reduces to $A=j_! \Lambda$ for some quasicompact separated \'etale map $j: U\to X$. In that case, it follows from $Rj_\ast$ commuting with the embedding $D_\et(X,\Lambda)\to D_\solid(X,\Lambda)$, as it is qcqs and has cohomological dimension $0$.

\subsection{Dual embedding} For a small v-stack $X$, let $D_\et^\dagger(X,\Lambda)\subset D_\et(X,\Lambda)$ be the full subcategory of overconvergent objects. Recall that $A\in D_\et(X,\Lambda)$ is overconvergent if for any strictly local $\Spa(C,C^+)\to X$, the map
\[
R\Gamma(\Spa(C,C^+),A)\to R\Gamma(\Spa(C,\mathcal O_C),A)
\]
is an isomorphism.

\begin{proposition}\label{prop:naivebiduality} Assume that $\Lambda=\mathbb Z/n\mathbb Z$ with $n$ prime to $p$. For any overconvergent $A\in D_\et^\dagger(X,\Lambda)$, let
\[
A^\vee = R\sHom_{D_\solid(X,\Lambda)}(A,\Lambda)\in D_\solid(X,\Lambda).
\]
Then the functor
\[
D_\et^\dagger(X,\Lambda)^\op\to D_\solid(X,\Lambda): A\mapsto A^\vee
\]
is fully faithful, $t$-exact (for the standard $t$-structure), compatible with pullback, and the map
\[
A\to R\sHom_{D_\solid(X,\Lambda)}(A^\vee,\Lambda)
\]
is an isomorphism.
\end{proposition}

\begin{proof} As the formation of $R\sHom$ in the solid context commutes with any base change, all assertions can be proved by v-descent, so we can assume that $X$ is strictly totally disconnected. Then $D_\et^\dagger(X,\Lambda)\cong D(\pi_0 X,\Lambda)$. The heart of the standard $t$-structure is then an abelian category with compact projective generators $i_\ast \Lambda$ for open and closed subsets $i: S\subset \pi_0 X$, and the whole category is the $\Ind$-category of the constructible complexes of $\Lambda$-modules on $\pi_0 X$ (which are locally constant with finite fibres). Passage to the naive dual is an autoequivalence on constructible complexes (as $\Lambda$ is selfinjective), and thus embeds the whole $\Ind$-category fully faithfully into the $\Pro$-category of constructible complexes of $\Lambda$-modules on $\pi_0 X$, which sits fully faithfully inside the category of finitely presented solid sheaves on $X$. This already establishes that the functor is fully faithful and $t$-exact, and we already observed at the beginning that it commutes with any pullback.

It remains to prove that
\[
A\to R\sHom_{D_\solid(X,\Lambda)}(A^\vee,\Lambda)
\]
is an isomorphism. Again, we can assume $\Lambda=\mathbb F_\ell$ so that all operations are $t$-exact. Again, the statement is clear if $A$ is constructible, and in general it follows from Breen's resolution that the $\Pro$-structure on $A^\vee$ dualizes to a filtered colimit on applying $R\sHom_{D_\solid(X,\Lambda)}(-,B)$.
\end{proof}

The functor $A\mapsto A^\vee$ is also close to being symmetric monoidal. Note that it is lax-symmetric monoidal, i.e.~there is a natural functorial map
\[
A^\vee\soliddotimesLambda B^\vee\to (A\dotimes_\Lambda B)^\vee.
\]

\begin{proposition}\label{prop:dualembeddingsymmmon} Assume that $A\in D_\et^\dagger(X,\Lambda)$ has finite Tor-amplitude over $\Lambda=\mathbb Z/n\mathbb Z$, i.e.~for all quotients $\Lambda\to \mathbb F_\ell$, the complex $A\dotimes_\Lambda \mathbb F_\ell\in D_\et^\dagger(X,\mathbb F_\ell)$ is bounded. Then for all $B\in D_\et^\dagger(X,\Lambda)$, the maps
\[
A^\vee\soliddotimesLambda B^\vee\to (A\dotimes_\Lambda B)^\vee\ ,\ A\dotimes_\Lambda B\to R\sHom_{D_\solid(X,\Lambda)}(A^\vee,B)
\]
are isomorphisms.
\end{proposition}

\begin{proof} The second follows from the first: Using Proposition~\ref{prop:naivebiduality},
\[
R\sHom_{D_\solid(X,\Lambda)}(A^\vee,B)\cong R\sHom_{D_\solid(X,\Lambda)}(A^\vee,R\sHom_{D_\solid(X,\Lambda)}(B^\vee,\Lambda))\cong R\sHom_{D_\solid(X,\Lambda)}(A^\vee\soliddotimesLambda B^\vee,\Lambda),
\]
which one can further rewrite to $A\dotimes_\Lambda B$ assuming the first isomorphism.

We can assume $\Lambda=\mathbb F_\ell$, and that $A$ is concentrated in degree $0$. Now as functors of $B$, all operations are $t$-exact, so we can reduce to the case that also $B$ is concentrated in degree $0$. We can assume that $X$ is strictly totally disconnected, and then $D_\et^\dagger(X,\mathbb F_\ell)\cong D(\pi_0 X,\mathbb F_\ell)$. Then $A$ and $B$ are filtered colimits of constructible sheaves on $\pi_0 X$, and $R\sHom(-,\mathbb F_\ell)$ is a contravariant autoequivalence on constructible $\mathbb F_\ell$-sheaves on $\pi_0 X$. Then the result follows by observing that $A\mapsto A^\vee$ simply exchanges the $\Ind$-category of constructible $\mathbb F_\ell$-sheaves on $\pi_0 X$ with its $\Pro$-category.
\end{proof}

As noted above, the functor $A\mapsto A^\vee$ is compatible with pullback. Regarding pushforward, we have the following result.

\begin{proposition}\label{prop:dualembeddingpushforward} Let $f: Y\to X$ be a proper map of small v-stacks that is representable in spatial diamonds with $\dimtrg f<\infty$. Let $A\in D_\et^\dagger(Y,\Lambda)$ with dual $A^\vee\in D_\solid(Y,\Lambda)$. Then there is a natural isomorphism
\[
(Rf_\ast A)^\vee\cong f_\natural A^\vee.
\]
\end{proposition}

Note that $Rf_\ast A$ is again overconvergent, by proper base change.

\begin{proof} One has
\[
R\sHom_{D_\solid(X,\Lambda)}(f_\natural A^\vee,\Lambda)\cong Rf_\ast R\sHom_{D_\solid(Y,\Lambda)}(A^\vee,\Lambda)\cong Rf_\ast A,
\]
so by biduality one gets a natural map
\[
f_\natural A^\vee\to (Rf_\ast A)^\vee;
\]
we claim that this is an isomorphism. This can be checked v-locally on $X$, so we can assume that $X$ is w-contractible. One can assume $A$ is bounded above (i.e.~$A\in D^-$) as both functors take very coconnective objects to very connective objects; by shifting, we can assume $A\in D^{\leq 0}$. Now using a Postnikov limit and the assumption $\dimtrg f<\infty$, we can also assume that $A\in D^+$, and hence reduce to $A$ sitting in degree $0$. Now we can choose a hypercover of $Y$ by perfectoid spaces $Y_i$ that are the canonical compactifications (relative to $X$) of w-contractible spaces. One can then replace $Y$ by one of the $Y_i$, so assume that $Y$ is the canonical compactification of a w-contractible space. In particular, $D_\et^\dagger(Y,\Lambda)\cong D(\pi_0 Y,\Lambda)$, and all operations can be computed on the level of $\pi_0 f: \pi_0 Y\to \pi_0 X$ instead. Here, the result amounts again to the duality between $\Ind$- and $\Pro$-objects in the category of constructible sheaves on profinite sets.
\end{proof}

\section{Dualizability}\label{sec:soliddualizability}

It turns out that most of the results above on Poincar\'e duality hold verbatim if the assumption that $f$ is cohomologically smooth is relaxed to the assumption that $\mathbb F_\ell$ is $f$-universally locally acyclic for all $\ell\neq p$. In fact, even more generally, one can obtain certain results comparing twisted forms of $f_\natural$ and $Rf_\ast$ for any $f$-universally locally acyclic complex $A$.

Assume that $\Lambda$ is a quotient of $\hat{\mathbb Z}^p$ of the form $\varprojlim_n \mathbb Z/n\mathbb Z$ where $n$ now runs only over some integers prime to $p$. If $f: X\to S$ is a compactifiable map of small v-stacks that is representable in locally spatial diamonds with locally $\dimtrg f<\infty$, we want to define the category $D^\ULA(X/S,\Lambda)$ of $f$-universally locally acyclic complexes with coefficients $\Lambda$ as the limit of the full subcategories
\[
D^\ULA(X/S,\mathbb Z/n\mathbb Z)\subset D_\et(X,\mathbb Z/n\mathbb Z)
\]
of $f$-universally locally acyclic objects in $D_\et(X,\mathbb Z/n\mathbb Z)$, for $n$ running over the same set of integers prime to $p$. As usual, the limit should be taken on the level of $\infty$-categorical enhancements. One way to phrase it without enhancements is to say that $D^\ULA(X/S,\Lambda)$ is the category of all $A\in D_\solid(X,\Lambda)$ such that $A_n=A\dotimes_\Lambda \mathbb Z/n\mathbb Z$ lies in $D_\et(X,\mathbb Z/n\mathbb Z)$ for all such $n$, is $f$-universally locally acyclic, and $A$ is the derived limit of the $A_n$.

Given such an $A$, in particular all $A_n$ are overconvergent, and the functor
\[
A\mapsto A^\vee = R\sHom_{D_\solid(X,\Lambda)}(A,\Lambda) = R\varprojlim_n A_n^\vee\in D_\solid(X,\Lambda)
\]
defines another fully faithful (contravariant) embedding
\[
D^\ULA(X/S,\Lambda)^\op\hookrightarrow D_\solid(X,\Lambda)
\]
of $f$-universally locally acyclic complexes into $D_\solid(X,\Lambda)$. We can also precompose with Verdier duality $\mathbb D_{X/S}$ to obtain a covariant fully faithful embedding
\[
D^\ULA(X/S,\Lambda)\hookrightarrow D_\solid(X,\Lambda): A\mapsto \mathbb D_{X/S}(A)^\vee.
\]

\begin{example} Assume that $S=\Spa C$ is a geometric point, and $X$ is the analytification of an algebraic variety $X^{\mathrm{alg}}/\Spec C$. Then any constructible complex on $X^{\mathrm{alg}}$ is universally locally acyclic over $S$, yielding a fully faithful embedding
\[
D^b_c(X^{\mathrm{alg}},\mathbb Z_\ell)\hookrightarrow D^\ULA(X/S,\mathbb Z_\ell)\hookrightarrow D_\solid(X,\mathbb Z_\ell),
\]
embedding the usual bounded derived category of constructible $\mathbb Z_\ell$-sheaves on $X^{\mathrm{alg}}$ into $D_\solid(X,\mathbb Z_\ell)$. The image lands in bounded complexes with finitely presented solid cohomology sheaves; in fact, in compact objects. Thus, this fully faithful embedding extends to a fully faithful embedding
\[
\Ind D^b_c(X^{\mathrm{alg}},\mathbb Z_\ell)\hookrightarrow D_\solid(X,\mathbb Z_\ell).
\]
The category on the left is the one customarily associated to $X^{\mathrm{alg}}$. This functor takes the sheaf $i_\ast \mathbb Z_\ell$, for a point $i: \Spec C\to X^{\mathrm{alg}}$, to the solid sheaf $i_\natural \mathbb Z_\ell$.

In many papers in geometric Langlands and related fields, one often finds the following construction. If $Y$ is a stack on the category of schemes over $\Spec C$, let
\[
\mathcal D(Y,\mathbb Z_\ell) := \varprojlim_{X^{\mathrm{alg}}\to Y} \Ind \mathcal D^b_c(X^{\mathrm{alg}},\mathbb Z_\ell)
\]
where $X^{\mathrm{alg}}$ runs over schemes of finite type over $\Spec C$, and the transition functors are given by $Rf^!$. This, in fact, embeds naturally into $D_\solid(Y^\diamond,\mathbb Z_\ell)$ via the previous embedding, noting that it intertwines $Rf^!$ with the usual pullback $f^\ast$ on solid sheaves. In fact,
\[
\mathbb D_{X'/S}(Rf^! A)^\vee\cong (f^\ast \mathbb D_{X/S}(A))^\vee\cong f^\ast \mathbb D_{X/S}(A)^\vee
\]
for a map $f: X^{\prime \mathrm{alg}}\to X^{\mathrm{alg}}$ of algebraic varieties over $\Spec C$.
\end{example}

Now for $A\in D^\ULA(X/S,\Lambda)$, we analyze the functor
\[
f_\natural(\mathbb D_{X/S}(A)^\vee\soliddotimesLambda -): D_\solid(X,\Lambda)\to D_\solid(S,\Lambda).
\]
We note that from the definition one sees that this functor commutes with all colimits, the formation of this functor commutes with any base change, and it satisfies the projection formula. In fact, this functor extends the functor $Rf_!(A\dotimes_\Lambda -)$.

\begin{proposition}\label{prop:generalizeRfshriek} Assume that $A\in D^\ULA(X/S,\Lambda)$ has bounded Tor-amplitude. Let $\mathbb Z/n\mathbb Z$ be a discrete quotient of $\Lambda$. For $B\in D_\et(X,\mathbb Z/n\mathbb Z)$, there is a natural equivalence
\[
f_\natural(\mathbb D_{X/S}(A)^\vee\soliddotimesLambda B)\cong Rf_!(A_n\dotimes_{\mathbb Z/n\mathbb Z} B)\in D_\et(S,\mathbb Z/n\mathbb Z).
\]
\end{proposition}

\begin{proof} We can assume $\Lambda=\mathbb Z/n\mathbb Z$. Note that for any $C\in D_\et(S,\Lambda)$, one has
\[\begin{aligned}
R\Hom_\Lambda(f_\natural(\mathbb D_{X/S}(A)^\vee\soliddotimesLambda B),C)&\cong R\Hom_\Lambda(\mathbb D_{X/S}(A)^\vee\soliddotimesLambda B,f^\ast C)\\
&\cong R\Hom_\Lambda(B,R\sHom_{D_\solid(X,\Lambda)}(\mathbb D_{X/S}(A)^\vee,f^\ast C))\\
&\cong R\Hom_\Lambda(B,\mathbb D_{X/S}(A)\dotimes_\Lambda f^\ast C)\\
&\cong R\Hom_\Lambda(B,R\sHom_{D_\et(X,\Lambda)}(A,Rf^! C))\\
&\cong R\Hom_\Lambda(A\dotimes_\Lambda B,Rf^! C)\\
&\cong R\Hom_\Lambda(Rf_!(A\dotimes_\Lambda B),C).
\end{aligned}\]
Here, we use Proposition~\ref{prop:naivebiduality} and Proposition~\ref{prop:ULAastshriekpullback}. In particular, there is a natural map
\[
f_\natural(\mathbb D_{X/S}(A)^\vee\soliddotimesLambda B)\to Rf_!(A\dotimes_\Lambda B).
\]
We claim that this is an isomorphism. This can be checked v-locally, so we can assume that $S$ is strictly totally disconnected. We can assume that $X$ is a spatial diamond by localization. As the functor commutes with all colimits in $B$, we can also assume that $B=j_!\Lambda$ for some quasicompact separated \'etale $j: V\to X$. Replacing $X$ by $V$, we can then even assume $B=\Lambda$.

Now $\mathbb D_{X/S}(A)$ lies in $D^+$ and then again $\mathbb D_{X/S}(A)^\vee$ in $D^-$. It follows that $\mathbb D_{X/S}(A)^\vee$ is a complex that is bounded above, and finitely presented solid in each degree. Thus $f_\natural \mathbb D_{X/S}(A)^\vee$ is of the same form, and so is the cone $Q$ of $f_\natural \mathbb D_{X/S}(A)^\vee\to Rf_! A$. If $Q$ is nonzero, we can look at the largest $i$ such that $\mathcal H^i(Q)$ is nonzero. This is finitely presented solid, so a cofiltered limit of constructible \'etale sheaves. But $R\Hom(Q,C)=0$ for all $C\in D_\et(S,\Lambda)$, so it follows that indeed $Q=0$.
\end{proof}

If $f$ is moreover proper, one can also prove the following version of $A$-twisted Poincar\'e duality.

\begin{proposition}\label{prop:Atwistedsolidpoincareduality} Assume that $f: X\to S$ is a proper map of small v-stacks that is representable in spatial diamonds with $\dimtrg f<\infty$. Let $A\in D^\ULA(X/S,\Lambda)$ with bounded Tor-amplitude. Then there is a natural equivalence
\[
f_\natural(\mathbb D_{X/S}(A)^\vee\soliddotimesLambda -)\cong Rf_\ast R\sHom_{D_\solid(X,\Lambda)}(A^\vee,-)
\]
of functors $D_\solid(X,\Lambda)\to D_\solid(S,\Lambda)$.
\end{proposition}

\begin{proof} First, we construct the natural transformation. Let $\pi_1,\pi_2: X\times_S X\to X$ be the two projections. Giving a map
\[
f_\natural(\mathbb D_{X/S}(A)^\vee\soliddotimesLambda B)\to Rf_\ast R\sHom_{D_\solid(X,\Lambda)}(A^\vee,B)
\]
is equivalent to giving a map
\[
f^\ast f_\natural (\mathbb D_{X/S}(A)^\vee\soliddotimesLambda B)\soliddotimesLambda A^\vee\to B.
\]
But
\[\begin{aligned}
f^\ast f_\natural (\mathbb D_{X/S}(A)^\vee\soliddotimesLambda B)\soliddotimesLambda A^\vee&\cong \pi_{1\natural} \pi_2^\ast(\mathbb D_{X/S}(A)^\vee\soliddotimesLambda B)\soliddotimesLambda A^\vee\\
&\cong \pi_{1\natural}(\pi_2^\ast\mathbb D_{X/S}(A)^\vee\soliddotimesLambda \pi_2^\ast B\soliddotimesLambda \pi_1^\ast A^\vee).\\
\end{aligned}\]
Thus, it suffices to construct a functorial map
\[
\pi_2^\ast\mathbb D_{X/S}(A)^\vee\soliddotimesLambda \pi_1^\ast A^\vee\soliddotimesLambda \pi_2^\ast B\to \pi_1^\ast B.
\]
For this in turn it suffices to construct a natural map
\[
\pi_2^\ast\mathbb D_{X/S}(A)^\vee\soliddotimesLambda \pi_1^\ast A^\vee\to \Delta_\natural \Lambda
\]
where $\Delta: X\hookrightarrow X\times_S X$ is the diagonal. Here $\Delta_\natural \Lambda\cong (\Delta_\ast \Lambda)^\vee$ by Proposition~\ref{prop:dualembeddingpushforward} and $\pi_2^\ast\mathbb D_{X/S}(A)^\vee\soliddotimesLambda \pi_1^\ast A^\vee\cong (\pi_2^\ast \mathbb D_{X/S}(A)\dotimes_\Lambda \pi_1^\ast A)^\vee$ by Proposition~\ref{prop:dualembeddingsymmmon}. Thus, we have to find a map
\[
\Delta_\ast \Lambda\to \pi_2^\ast \mathbb D_{X/S}(A)\dotimes_\Lambda \pi_1^\ast A
\]
or equivalently a section of $R\Delta^!(\pi_2^\ast \mathbb D_{X/S}(A)\dotimes_\Lambda \pi_1^\ast A)$. But $\pi_2^\ast \mathbb D_{X/S}(A)\dotimes_\Lambda \pi_1^\ast A\cong R\sHom(\pi_2^\ast A,R\pi_1^! A)$ as $A$ is $f$-universally locally acyclic, and then
\[
R\Delta^!(\pi_2^\ast \mathbb D_{X/S}(A)\dotimes_\Lambda \pi_1^\ast A)\cong R\Delta^!R\sHom(\pi_2^\ast A,R\pi_1^! A)\cong R\sHom(A,A),
\]
where we find the identity section.

To show that the map is an isomorphism, we can now localize on $S$, and in particular assume that $S$ is strictly totally disconnected. By the bounded assumption on $A$ and finite cohomological dimension of $f$, the functor $Rf_\ast R\sHom_\Lambda(A^\vee,-)$ commutes with all direct sums, and hence we can assume that $B$ is finitely presented solid (concentrated in degree $0$). Then we can write $B$ as a cofiltered limit of constructible \'etale sheaves, and the left-hand side commutes with such limits; so we can reduce to $B$ being a constructible \'etale sheaf, where the result follows from Proposition~\ref{prop:generalizeRfshriek} and Proposition~\ref{prop:naivebiduality}.
\end{proof}

From the perspective of using sheaves as kernels of induced functors, we have the following picture. We can introduce a variant of the category $\mathcal C_S$ introduced above. Namely, for any small v-stack $S$, let us consider the $2$-category $\mathcal C_{S,\solid}$ whose objects are relatively $0$-truncated small v-stacks $X$ over $S$, and whose categories of morphisms
\[
\Fun_{\mathcal C_{S,\solid}}(X,Y) = D_\solid(X\times_S Y,\Lambda)
\]
are given by solid complexes. Again, to any $X\in \mathcal C_{S,\solid}$, we can associate the triangulated category $D_\solid(X,\Lambda)$ and to any $A\in D_\solid(X\times_S Y,\Lambda)$ the functor
\[
p_{2\natural}(A\soliddotimesLambda p_1^\ast): D_\solid(X,\Lambda)\to D_\solid(Y,\Lambda)
\]
with kernel $A$. The composition in $\mathcal C_{S,\solid}$ is defined by the convolution
\[
D_\solid(X\times_S Y,\Lambda)\times D_\solid(Y\times_S Z,\Lambda)\to D_\solid(X\times_S Z,\Lambda): (A,B)\mapsto A\star B = p_{13\natural}(p_{12}^\ast A\soliddotimesLambda p_{23}^\ast B).
\]

We wish to compare $\mathcal C_S$ and $\mathcal C_{S,\solid}$. Note that the naive embedding $D_\et(X\times_S Y,\Lambda)\hookrightarrow D_\solid(X\times_S Y,\Lambda)$ is not compatible with the convolution (as one employs $R\pi_{13!}$ while the other employs $\pi_{13\natural}$). On the other hand, we can restrict to the sub-$2$-category $\mathcal C_S^\dagger\subset \mathcal C_S$ whose objects are only the proper $X/S$ representable in spatial diamonds of finite $\dimtrg$, and with
\[
\Fun_{\mathcal C_S^\dagger}(X,Y)=D_{\et,\mathrm{ftor}}^\dagger(X\times_S Y,\Lambda),
\]
where the subscript $\mathrm{ftor}$ stands for finite Tor-dimension over $\Lambda$. Then there is a fully faithful embedding $\mathcal C_S^\dagger\hookrightarrow \mathcal C_{S,\solid}^\mathrm{co}$ where the superscript $\mathrm{co}$ means that we change the direction of the arrows within each $\Fun_{\mathcal C_{S,\solid}}(X,Y)$. Indeed, for any $X,Y\in \mathcal C_S^\dagger$, the functor $A\mapsto A^\vee$ defines a fully faithful embedding
\[
\Fun_{\mathcal C_S^\dagger}(X,Y)=D_{\et,\mathrm{ftor}}^\dagger(X\times_S Y,\Lambda)\hookrightarrow D_\solid(X\times_S Y,\Lambda)^\op=\Fun_{\mathcal C_{S,\solid}^\mathrm{co}}(X,Y).
\]
This is compatible with composition by Proposition~\ref{prop:dualembeddingsymmmon} and Proposition~\ref{prop:dualembeddingpushforward}. This discussion leads to another proof of Proposition~\ref{prop:Atwistedsolidpoincareduality}:

\begin{corollary}\label{cor:dualizablesolid} Let $f: X\to S$ be a proper map of small v-stacks that is representable in spatial diamonds with $\dimtrg f<\infty$. Let $A\in D_\et(X,\Lambda)$ be $f$-universally locally acyclic and of finite Tor-dimension over $\Lambda$. Then $A^\vee\in D_\solid(X,\Lambda)=\Fun_{\mathcal C_{S,\solid}}(X,S)$ is right adjoint to
\[
\mathbb D_{X/S}(A)^\vee\in D_\solid(X,\Lambda)=\Fun_{\mathcal C_{S,\solid}}(S,X).
\]
In particular, the functor
\[
f_\natural(A^\vee\soliddotimesLambda -): D_\solid(X,\Lambda)\to D_\solid(S,\Lambda)
\]
is right adjoint to the functor
\[
\mathbb D_{X/S}(A)^\vee\soliddotimesLambda f^\ast -: D_\solid(S,\Lambda)\to D_\solid(X,\Lambda),
\]
so
\[
f_\natural(A^\vee\soliddotimesLambda -)\cong Rf_\ast R\sHom_\Lambda(\mathbb D_{X/S}(A)^\vee,-): D_\solid(X,\Lambda)\to D_\solid(S,\Lambda).
\]
\end{corollary}

Moreover, when applied to the Satake category, we get a fully faithful embedding
\[
(\Sat_G^I)^\op\hookrightarrow D_\solid(\Hloc_G^I,\hat{\mathbb Z}^p): A\mapsto A^\vee
\]
compatible with the monoidal structure (and functorially in $I$), where the right-hand side is given by
\[
\Fun_{\mathcal C_{S,\solid}}([(\Div^1)^I_X/L^+_{(\Div^1_X)^I} G],[(\Div^1)^I_X/L^+_{(\Div^1_X)^I} G])
\]
for $S=[(\Div^1_X)^I/L_{(\Div^1_X)^I} G]$. Precomposing with Verdier duality, we get a covariant fully faithful embedding
\[
\Sat_G^I\hookrightarrow D_\solid(\Hloc_G^I,\hat{\mathbb Z}^p): A\mapsto \mathbb D_{X/S}(A)^\vee.
\]
By Proposition~\ref{prop:generalizeRfshriek}, when one uses objects in the Satake category as kernels to define Hecke operators, this fully faithful embedding makes it possible to extend Hecke operators from $D_\et$ to $D_\solid$.

\section{Lisse-\'etale sheaves}

The category $D_\solid(X,\Lambda)$ is huge: Already if $X$ is a point and $\Lambda=\mathbb F_\ell$, it is the derived category of solid $\mathbb F_\ell$-vector spaces, which is much larger than the category of usual discrete $\mathbb F_\ell$-vector spaces. When applied to $\Bun_G$, we would however really like to study smooth representations on discrete $\Lambda$-modules.

As coefficients, we will from now on choose a discrete $\mathbb Z_\ell$-algebra $\Lambda$ for some $\ell\neq p$, or rather the corresponding condensed ring $\Lambda:=\mathbb Z_\ell\otimes_{\mathbb Z_{\ell,\mathrm{disc}}} \Lambda_{\mathrm{disc}}$. (For a technical reason, we have to restrict attention to a particular prime $\ell$.)

It turns out that when $X$ is an Artin v-stack, one can define a full subcategory $D_\lis(X,\Lambda)\subset D_\solid(X,\Lambda)$ that when specialized to $X=\Bun_G$ has the desired properties. Here the subscript ``lis'' is an abbreviation of ``lisse'' (french smooth), and is not meant to evoke lisse sheaves in the sense of locally constant sheaves, but lisse-\'etale sheaves in the sense of Artin stacks \cite{LaumonMoretBailly}.

\begin{definition}\label{def:Dlisse} Let $X$ be an Artin v-stack. The full subcategory $D_\lis(X,\Lambda)\subset D_\solid(X,\Lambda)$ is the smallest triangulated subcategory stable under all direct sums that contains $f_\natural \Lambda$ for all maps $f: Y\to X$ that are separated, representable in locally spatial diamonds, and $\ell$-cohomologically smooth.
\end{definition}

In principle, one could give this definition even when $X$ is any small v-stack, but in that case there might be very few objects.

\begin{proposition}\label{prop:Dlissetensorpullback} Let $X$ be an Artin v-stack. The full subcategory $D_\lis(X,\Lambda)\subset D_\solid(X,\Lambda)$ is stable under $-\soliddotimesLambda-$. Moreover, if $f: Y\to X$ is a map of Artin v-stacks, then $f^\ast$ maps $D_\lis(X,\Lambda)\subset D_\solid(X,\Lambda)$ into $D_\lis(Y,\Lambda)\subset D_\solid(Y,\Lambda)$.
\end{proposition}

\begin{proof} As tensor products and pullbacks commute with all direct sums, it suffices to check the claim on the generators $g_\natural \Lambda$ for maps $g: Z\to X$ that are separated, representable in locally spatial diamonds, and $\ell$-cohomologically smooth. Now the result follows as pullbacks and products of such maps are of the same form.
\end{proof}

\begin{proposition}\label{prop:lissification} Let $X$ be an Artin v-stack. The inclusion $D_\lis(X,\Lambda)\subset D_\solid(X,\Lambda)$ admits a right adjoint
\[
A\mapsto A^\lis: D_\solid(X,\Lambda)\to D_\lis(X,\Lambda).
\]
The kernel of $A\mapsto A^\lis$ is the class of all $A\in D_\solid(X,\Lambda)$ such that $A(Y)=0$ for all $f: Y\to X$ that are separated, representable in locally spatial diamonds, and $\ell$-cohomologically smooth.
\end{proposition}

\begin{proof} The existence of the right adjoint is formal. We note that the $\infty$-category $\mathcal D_\solid(X,\Lambda)$ is not itself presentable, but rather is the large filtered colimit of presentable $\infty$-categories $\mathcal D_\solid(X_\kappa,\Lambda)$ for uncountable strong limit cardinals $\kappa$ (restricting the v-site to $\kappa$-small perfectoid spaces). Also note that $D_\lis(X,\Lambda)$ is contained in $D_\solid(X_\kappa,\Lambda)$ for some $\kappa$: This can be checked when $X$ is a spatial diamond and for the generators $f_\natural \Lambda\cong f_\natural \mathbb Z_\ell\dotimes_{\mathbb Z_\ell} \Lambda$ when $f: Y\to X$ is in addition quasicompact, in which case $f_\natural \mathbb Z_\ell$ is the limit of $f_\natural \mathbb Z/\ell^m\mathbb Z$ all of which lie in $D_\et(X,\mathbb Z/\ell^m\mathbb Z)$, so we conclude by \cite[Remark 17.4]{ECoD}. It follows that the right adjoints to $D_\lis(X,\Lambda)\to D_\solid(X_\kappa,\Lambda)$ for all large enough $\kappa$ glue to the desired right adjoint.

The description of the kernel is formal.
\end{proof}

Using Proposition~\ref{prop:lissification}, we can then also define $R\sHom_\lis(A,B)\in D_\lis(X,\Lambda)$ for $A,B\in D_\lis(X,\Lambda)$ and $Rf_{\lis\ast}: D_\lis(Y,\Lambda)\to D_\lis(X,\Lambda)$ for a map $f: Y\to X$ of Artin v-stacks, satisfying the usual adjunction to the tensor product and pullback.

The goal of passing to $D_\lis$ is to make sheaves ``discrete'' again. Recall the following result.

\begin{proposition} For any condensed ring $A$ with underlying ring $A(\ast)$, the functor $M\mapsto M\otimes_{A(\ast)} A$ induces a fully faithful functor
\[
D(A(\ast))\hookrightarrow D(A)
\]
from the derived category of usual $A(\ast)$-modules to the derived category of condensed modules over the condensed ring $A$.
\end{proposition}

\begin{proof} We need to see that for any $M,N\in D(A(\ast))$, the map
\[
R\Hom_{A(\ast)}(M,N)\to R\Hom_{A}(M\otimes_{A(\ast)} A,N\otimes_{A(\ast)} A)
\]
is an isomorphism. The class of all $M$ for which this happens is triangulated and stable under all direct sums, so it suffices to consider $M=A(\ast)$. Then it amounts to
\[
N(\ast)\to (N\otimes_{A(\ast)}A)(\ast)
\]
being an isomorphism, which follows from evaluation at $\ast$ being symmetric monoidal.
\end{proof}

In particular, we have the following result for a geometric point.

\begin{proposition}\label{prop:Dlispoint} Let $X=\Spa C$ for some complete algebraically closed nonarchimedean field $C$. Then $D_\lis(X,\Lambda)\cong D(\Lambda)$, the derived category of (relatively) discrete $\Lambda$-modules.
\end{proposition}

\begin{proof} We need to see that for all separated $\ell$-cohomologically smooth maps $f: Y\to X$ of spatial diamonds, one has $f_\natural \Lambda\in D(\Lambda)$. This reduces to $\Lambda=\mathbb Z_\ell$. In that case, $f_\natural \mathbb Z_\ell = \varprojlim_m f_\natural \mathbb Z/\ell^m\mathbb Z$, where by Proposition~\ref{prop:generalizeRfshriek} each
\[
f_\natural \mathbb Z/\ell^m\mathbb Z\cong Rf_! Rf^! \mathbb Z/\ell^m\mathbb Z,
\]
which is a perfect complex of $\mathbb Z/\ell^m\mathbb Z$-modules, in particular discrete. Taking the limit over $m$, we get a perfect complex of $\mathbb Z_\ell$-modules, which is in particular (relatively) discrete over $\mathbb Z_\ell$.
\end{proof}

When working with torsion coefficients, one recovers $D_\et$.

\begin{proposition}\label{prop:DlisDet} Let $X$ be an Artin v-stack, and assume that $\Lambda$ is killed by a power of $\ell$. Then $D_\lis(X,\Lambda)\subset D_\solid(X,\Lambda)$ is contained in the image of the naive embedding $D_\et(X,\Lambda)\hookrightarrow D_\solid(X,\Lambda)$. If there is a separated $\ell$-cohomologically smooth surjection $U\to X$ from a locally spatial diamond $U$, such that $U_\et$ has a basis with bounded $\ell$-cohomological dimension, then it induces an equivalence $D_\lis(X,\Lambda)\cong D_\et(X,\Lambda)$.
\end{proposition}

\begin{proof} If $f: Y\to X$ is separated, representable in locally spatial diamonds, and $\ell$-cohomologically smooth, then $f_\natural \Lambda=Rf_!Rf^!\Lambda$ lies in $D_\et(X,\Lambda)$, hence $D_\lis(X,\Lambda)\subset D_\et(X,\Lambda)$. To check equality, we can work on an atlas, so by the assumption we can reduce to the case that $X$ is a locally spatial diamond for which $X_\et$ has a basis with bounded $\ell$-cohomological dimension. In that case $D_\et(X,\Lambda)\cong D(X_\et,\Lambda)$ by \cite[Proposition 20.17]{ECoD} (the proof only needs a basis with bounded cohomological dimension), which is generated by $j_\natural\Lambda$ for $j: U\to X$ quasicompact separated \'etale, which is thus also contained in $D_\lis(X,\Lambda)$.
\end{proof}

The most severe problem with the general formalism of solid sheaves is that stratifications of a space do not lead to corresponding decompositions of sheaves into pieces on the individual strata. This problem is somewhat salvaged by $D_\lis(X,\Lambda)$: We expect that it holds true if $X$ and its stratification are sufficiently nice. Here is a simple instance that will be sufficient for our purposes.

\begin{proposition}\label{prop:stratificationpoint} Let $X$ be a locally spatial diamond with a closed point $x\in X$, giving a corresponding closed subdiamond $i: Z\subset X$ with complement $j: U\subset X$. Assume that $Z=\Spa C$ is representable, with $C$ an algebraically closed nonarchimedean field. Moreover, assume that $Z$ can be written as a cofiltered intersection of qcqs open neighborhoods $V\subset X$ such that $R\Gamma(V,\mathbb F_\ell)\cong \mathbb F_\ell$.

Then one has a semi-orthogonal decomposition of $D_\lis(X,\Lambda)$ into $D_\lis(U,\Lambda)$ and $D_\lis(Z,\Lambda)\cong D(\Lambda)$.
\end{proposition}

\begin{proof} We may assume that $X$ is spatial. We analyze the quotient of $D_\lis(X,\Lambda)$ by $j_\natural D_\lis(U,\Lambda)$. This is equivalently the subcategory of all $A\in D_\lis(X,\Lambda)$ with $j^\ast A=0$, as this is the right orthogonal of $j_\natural D_\lis(U,\Lambda)$; the composite
\[
D_\lis(X,\Lambda)\to D_\lis(X,\Lambda)/j_\natural D_\lis(U,\Lambda)\to D_\lis(X,\Lambda)
\]
is given by $A\mapsto \mathrm{cone}(j_\natural j^\ast A\to A)$. This quotient is generated by the images of $f_\natural \Lambda$ for $f: Y\to X$ cohomologically smooth separated map of spatial diamonds; under the embedding of the quotient category back into $D_\lis(X,\Lambda)$, this corresponds to the cone of $j_\natural j^\ast f_\natural \Lambda\to f_\natural \Lambda$. Let $M = i^\ast f_\natural \Lambda\in D_\lis(Z,\Lambda)\cong D(\Lambda)$, which in fact is a perfect complex of $\Lambda$-modules (by the proof of Proposition~\ref{prop:Dlispoint}). Then we claim that there is an isomorphism
\[
\mathrm{cone}(j_\natural j^\ast f_\natural \Lambda\to f_\natural \Lambda)\cong \mathrm{cone}(j_\natural M\to M).
\]
To see this, it suffices to prove that there is some open neighborhood $V$ of $Z$ such that $f_\natural \Lambda|_V\cong M$, the constant sheaf associated with $M$. We can reduce to $\Lambda=\mathbb Z_\ell$. As $f_\natural \mathbb F_\ell$ is constructible, we can find some such $V$ for which $f_\natural \mathbb F_\ell|_V\cong M/\ell$. Picking such an isomorphism reducing to the identity at $x$, and choosing $V$ with the property $R\Gamma(V,\mathbb F_\ell)\cong \mathbb F_\ell$, we see that in fact the isomorphism lifts uniquely to $\mathbb Z/\ell^m \mathbb Z$ for each $m$, and thus by taking the limit over $m$ to the desired isomorphism $f_\natural \Lambda|_V\cong M$.

Thus, the quotient of $D_\lis(X,\Lambda)$ by $j_\natural D_\lis(U,\Lambda)$ is generated by the constant sheaf $\Lambda$. Moreover, the endomorphisms of $\Lambda$ in the quotient category are given by the cone of
\[
R\Gamma(X,j_\natural \Lambda)\to R\Gamma(X,\Lambda).
\]
This is equivalently the filtered colimit of $R\Gamma(V,\Lambda)$ over all qcqs open neighborhoods $V$ of $Z$; we can restrict to those for which $R\Gamma(V,\mathbb F_\ell)\cong \mathbb F_\ell$. This implies formally that $R\Gamma(V,\mathbb Z_\ell)\cong \mathbb Z_\ell$ by passing to limits and then $R\Gamma(V,\Lambda)\cong \Lambda$ by passing to filtered colimits. Thus, we get the desired semi-orthogonal decomposition.
\end{proof}

\section{$D_{\mathrm{lis}}(\mathrm{Bun}_G)$}

Our goal now is to extend the results of Chapter~\ref{ch:DetBunG} to the case of $D_\lis(\Bun_G,\Lambda)$. This will notably include the case $\Lambda=\overline{\mathbb Q}_\ell$.

Thus, let again be $E$ any nonarchimedean local field with residue field $\Fq$ and $G$ a reductive group over $E$. We work with $\Perf_k$ where $k=\overline{\mathbb F}_q$, and fix a complete algebraically closed nonarchimedean field $C/k$.

\begin{proposition}\label{prop:DlisBunGb} Let $b\in B(G)$. The pullback functors
\[\begin{aligned}
D_\lis(\Bun_G^b,\Lambda)&\to D_\lis([\ast/\underline{G_b(E)}],\Lambda)\to D_\lis([\Spa C/\underline{G_b(E)}],\Lambda),\\
\ D_\lis(\Bun_G^b,\Lambda)&\to D_\lis(\Bun_G^b\times \Spa C,\Lambda)\to D_\lis([\Spa C/\underline{G_b(E)}],\Lambda)
\end{aligned}\]
are equivalences, and all categories are naturally equivalent (as symmetric monoidal categories) to the derived category $D(G_b(E),\Lambda)$ of smooth representations of $G_b(E)$ on discrete $\Lambda$-modules.
\end{proposition}

\begin{proof} Recall that the map $s: [\ast/\underline{G_b(E)}]\to \Bun_G^b$ is cohomologically smooth and surjective; in fact, its fibres are successive extensions of positive Banach--Colmez spaces. This implies that $s_\natural \Lambda\cong \Lambda$. This, in turn, implies by the projection formula for $s_\natural$ that $s_\natural s^\ast A\cong A$ for all $A\in D_\solid(\Bun_G^b,\Lambda)$, thus giving fully faithfulness. The same applies after base change to $\Spa C$. Moreover, using pullback under the projection $\Bun_G^b\to [\ast/\underline{G_b(E)}]$, we see that $s^\ast$ is also necessarily essentially surjective.

It remains to show that the pullback $D_\lis([\ast/\underline{G_b(E)}],\Lambda)\to D_\lis([\Spa C/\underline{G_b(E)}],\Lambda)$ is an equivalence, and identify this symmetric monoidal category with $D(G_b(E),\Lambda)$. By Proposition~\ref{prop:solidchangeofbasefield}, the functor $D_\lis([\ast/\underline{G_b(E)}],\Lambda)\to D_\lis([\Spa C/\underline{G_b(E)}],\Lambda)$ is fully faithful. One can easily build a functor $D(G_b(E),\Lambda)\to D_\lis([\ast/\underline{G_b(E)}],\Lambda)$, and it is enough to see that the composite functor
\[
D(G_b(E),\Lambda)\to D_\lis([\ast/\underline{G_b(E)}],\Lambda)\to D_\lis([\Spa C/\underline{G_b(E)}],\Lambda)
\]
is an equivalence. Using that $D(G_b(E),\Lambda)$ is generated by $c\text-\mathrm{Ind}_K^{G_b(E)} \Lambda$ for $K\subset G_b(E)$ open pro-$p$, one easily sees that the functor is fully faithful, so it remains to prove essential surjectivity. Using descent along $\Spa C\to [\Spa C/\underline{G_b(E)}]$ and the equivalence $D_\lis(\Spa C,\Lambda)\cong D(\Lambda)$, any cohomology sheaf of an object in $D_\lis([\Spec C/\underline{G_b(E)}],\Lambda)$ gives a representation of the condensed group $G_b(E)$ on a condensed $\Lambda$-module of the form $M\otimes_{\mathbb Z_{\ell,\mathrm{disc}}} \mathbb Z_\ell$ for some (abstract) $\Lambda$-module $M$. (Indeed, any cohomology sheaf of an object of $D_\solid([\Spa C/\underline{G_b(E)}],\Lambda)$ is a representation of the condensed group $G_b(E)$ on a solid $\Lambda$-module, but here we must get objects whose underlying solid module lies in $D_\lis(\Spa C,\Lambda)\cong D(\Lambda)$.) As $G_b(E)$ is locally pro-$p$, any such action in fact comes from a smooth action on $M$: For $K\subset G_b(E)$ pro-$p$, the $K$-orbit of any $m\in M$ lies in some compact submodule, thus in $M\otimes_{\mathbb Z_{\ell,\mathrm{disc}}} \mathbb Z_\ell$ for some finitely generated $\mathbb Z_\ell$-submodule $M'\subset M$. The action of $K$ on $m$ then gives a continuous map $K\to \GL(M')$. As the target is locally pro-$\ell$, this map has finite image, so that the action of $K$ on $m$ is locally constant.
\end{proof}

Recall that for any $b\in B(G)$, we have the cohomologically smooth chart $\pi_b: \mathcal M_b\to \Bun_G$ near $\Bun_G^b$. This comes with a projection $q_b: \mathcal M_b\to [\ast/\underline{G_b(E)}]$ which has a natural section, given by the preimage of $\Bun_G^b\subset \Bun_G$ in $\mathcal M_b$. Over $\mathcal M_b$, we have the $\underline{G_b(E)}$-torsor $\tilde{\mathcal M}_b\to \mathcal M_b$, and for any complete algebraically closed field $C$ over $k=\overline{\mathbb F}_q$, the base change
\[
\tilde{\mathcal M}_{b,C} = \tilde{\mathcal M}_b\times_{\Spd k} \Spa C
\]
is representable by a locally spatial diamond, endowed with a distinguished point $i: \Spa C\hookrightarrow \tilde{\mathcal M}_{b,C}$.  Recall that $\tilde{\mathcal M}_{b,C}$ is a successive extension of negative Banach--Colmez spaces. Iteratively restricting to small quasicompact balls inside these negative Banach--Colmez spaces, we see that the closed subset $i: \Spa C\hookrightarrow \tilde{\mathcal M}_{b,C}$ can be written as cofiltered intersection of quasicompact open subsets $V$ for which $R\Gamma(V,\mathbb F_\ell)\cong \mathbb F_\ell$. (To see that one can choose small balls in negative Banach--Colmez spaces with this property, one can present $\BC_C(\mathcal O(-n)[1])$ as a quotient of a product of $n$ different $(\mathbb A^1_{C_i^\sharp})^\diamond$ by $\underline{E}$ (via taking a generic embedding of $\mathcal O_{X_C}(-n)$ into $\mathcal O_{X_C}$, whose cokernel is a sum of $n$ skyscraper sheaves at different untilts $C_i^\sharp/E$ of $C$), and take balls in $\BC_C(\mathcal O(-n)[1])$ that are similar products of balls, quotiented by a lattice in $\underline{E}$.)

\begin{proposition}\label{prop:compactobjectsDlis} For any $b\in B(G)$ with locally closed immersion $i^b: \Bun_G^b\to \Bun_G$, the functor
\[
i^{b\ast}: D_\lis(\Bun_G,\Lambda)\to D_\lis(\Bun_G^b,\Lambda)\cong D_\lis([\ast/\underline{G_b(E)}],\Lambda)
\]
admits a left adjoint, given by
\[
\pi_{b\natural} q_b^\ast: D_\lis([\ast/\underline{G_b(E)}],\Lambda)\to D_\lis(\Bun_G,\Lambda).
\]
The unit of the adjunction is given by the equivalence $\mathrm{id}\cong i^{b\ast}\pi_{b\natural} q_b^\ast$ arising from base change, and the identification of the pullback of $i^b$ along $\pi_b$ with $[\ast/\underline{G_b(E)}]\subset \mathcal M_b$.
\end{proposition}

\begin{proof} As $D([\ast/\underline{G_b(E)}],\Lambda)\cong D(G_b(E),\Lambda)$ is generated by $c\text-\mathrm{Ind}_K^{G_b(E)} \Lambda$ for open pro-$p$ subgroups $K\subset G_b(E)$, and as we already determined the unit of the adjunction, it suffices to verify the adjunction on these objects. Let $\mathcal M_{b,K}=\tilde{\mathcal M}_b/\underline{K}\to \mathcal M_b$. This comes with a closed immersion $i_K: [\ast/\underline{K}]\to \mathcal M_{b,K}$. It suffices to see that for all $A\in D_\lis(\mathcal M_{b,K},\Lambda)$, the map
\[
R\Gamma(\mathcal M_{b,K},A)\to R\Gamma([\ast/\underline{K}],A)
\]
is an isomorphism, where we continue to denote by $A$ any of its pullbacks. Assume first that $A=j_\natural A_0$ for some $A_0\in D_\lis(\mathcal M_{b,K}^\circ,\Lambda)$. Then the result follows from Theorem~\ref{thm:partialcompactsupportvanishingsolid}. In general we can then replace $A$ by the cone of $j_\natural A\to A$ in the displayed formula. For this statement, we can even base change to $\Spa C$ for some complete algebraically closed nonarchimedean field $C|k$, and allow more generally any $A\in D_\lis(\mathcal M_{b,K,C},\Lambda)$. We can then assume that $A=f_\natural \mathbb Z_\ell$ for some $\ell$-cohomologically smooth separated qcqs map $f: Y\to \mathcal M_{b,K,C}$. Then as in the proof of Proposition~\ref{prop:stratificationpoint}, $A$ is constant in a neighborhood of $[\Spa C/\underline{K}]$, which implies the result (as $\Spa C\subset \tilde{\mathcal M}_{b,C}$ is a cofiltered intersection of quasicompact open $V$'s with trivial cohomology).
\end{proof}

\begin{proposition}\label{prop:DlisBunG} For any quasicompact open substack $U\subset \Bun_G$, the Harder--Narasimhan stratification induces a semi-orthogonal decomposition of $D_\lis(U,\Lambda)$ into the categories $D_\lis(\Bun_G^b,\Lambda)\cong D(G_b(E),\Lambda)$ for $b\in |U|\subset B(G)$. Moreover, for any not necessarily quasicompact $U$, the functor
\[
D_\lis(U,\Lambda)\to D_\lis(U\times_{\Spd k} \Spa C,\Lambda)
\]
is an equivalence.
\end{proposition}

\begin{proof} We argue by induction on $|U|$, so take some closed element $b\in |U|\subset B(G)$ and let $i:\Bun_G^b\to U$ and $j: V\to U$ be the closed and complementary open substacks. We know that $D_\lis(U,\Lambda)\to D_\lis(U\times_{\Spd k} \Spa C,\Lambda)$ is fully faithful by Proposition~\ref{prop:solidchangeofbasefield}, and by induction $D_\lis(V,\Lambda)\to D_\lis(V\times_{\Spd k} \Spa C,\Lambda)$ is an equivalence.

Now by the previous proposition, $i^{b\ast}$ admits the left adjoint
\[
\pi_{b\natural} q_b^\ast: D_\lis([\ast/\underline{G_b(E)}],\Lambda)\to D_\lis(U,\Lambda),
\]
and in fact the proof of that proposition shows (using our standing induction assumption) that, composed with the embedding into $D_\lis(U\times_{\Spd k} \Spa C,\Lambda)$, it continues to be a left adjoint to $i^{b\ast}: D_\lis(U\times_{\Spd k} \Spa C,\Lambda)\to D([\Spa C/\underline{G_b(E)}],\Lambda)\cong D([\ast/\underline{G_b(E)}],\Lambda)$.

The unit $\mathrm{id}\to i^{b\ast}\pi_{b\natural} q_b^\ast$ of the adjunction is an equivalence. We see that $D_\lis(U\times_{\Spd k} \Spa C,\Lambda)$ has full subcategories given by $j_\natural D_\lis(V,\Lambda)$ and the essential image of $\pi_{b\natural} q_b^\ast$ (both of which lie in $D_\lis(U,\Lambda)$). To see that one has a semi-orthogonal decomposition, it suffices to see that if $A\in D_\lis(U\times_{\Spd k} \Spa C,\Lambda)$ with $i^\ast A = j^\ast A = 0$, then $A=0$. This can be checked after pullback to $\tilde{\mathcal M}_{b,C}$, where it follows from Proposition~\ref{prop:stratificationpoint}. This also shows that $D_\lis(U,\Lambda)\to D_\lis(U\times_{\Spd k} \Spa C,\Lambda)$ is an equivalence.
\end{proof}

Now we also want to analyze the compact objects as well as the universally locally acyclic objects, and various dualities. We start with the compact objects.

\begin{proposition}\label{prop:compactlisseBunG} The category $D_\lis(\Bun_G,\Lambda)$ is compactly generated. An object $A\in D_\lis(\Bun_G,\Lambda)$ is compact if and only if it has finite support and $i^{b\ast} A\in D_\lis(\Bun_G^b,\Lambda)\cong D(G_b(E),\Lambda)$ is compact for all $b\in B(G)$, i.e.~lies in the thick triangulated subcategory generated by $c\text-\mathrm{Ind}_K^{G_b(E)} \Lambda$ for open pro-$p$ subgroups $K\subset G_b(E)$.

Moreover, for each $b$ and $K\subset G_b(E)$ pro-$p$, letting
\[
f_K: \mathcal M_{b,K}\to \Bun_G
\]
be the natural map, the object $A_K^b = f_{K\natural} \Lambda\in D_\lis(\Bun_G,\Lambda)$ is compact, and these generate $D_\lis(\Bun_G,\Lambda)$.
\end{proposition}

\begin{proof} By Proposition~\ref{prop:DlisBunG}, the left adjoints $\pi_{b\natural} q_b^\ast$ to $i^{b\ast}$ generate $D_\lis(\Bun_G,\Lambda)$; as $i^{b\ast}$ commutes with colimits, these left adjoints also preserve compact objects. As each $D(G_b(E),\Lambda)$ is compactly generated, it follows that $D_\lis(\Bun_G,\Lambda)$ is compactly generated, with compact generators $A^b_K$.

To see that the given property characterizes compact objects, we argue by induction over quasicompact open substacks $U\subset \Bun_G$. Pick any closed $b\in |U|\subset B(G)$, and assume the result for the complementary open $j: V\subset U$. We first show that all of the given compact generators (coming from $b'\in |U|\subset B(G)$) have the property that all of their stalks are compact. This is clear by induction if $b'\in |V|$, so we can assume $b'=b$. Then we need to see that $j^\ast \pi_{b\natural} q_b^\ast$ preserves compact objects. But this follows from Lemma~\ref{lem:smoothcompact} below. Using the semi-orthogonal decomposition structure, it now follows that conversely, all $A$ with compact stalks are compact.
\end{proof}

\begin{lemma}\label{lem:smoothcompact} For $K\subset G_b(E)$ an open pro-$p$ subgroup, the functor
\[
R\Gamma(\mathcal M_{b,K}^\circ,-): D_\solid(\mathcal M_{b,K}^\circ,\Lambda)\to D(\Lambda)
\]
has finite cohomological dimension and commutes with all direct sums.
\end{lemma}

\begin{proof} As $\mathcal M_{b,K}^\circ=\tilde{\mathcal M}_b^\circ/\underline{K}$ where $\tilde{\mathcal M}_b^\circ$ is a spatial diamond, it suffices to prove that the functor has finite cohomological dimension. It suffices to prove this for $\tilde{\mathcal M}_b^\circ$ (as taking $K$-invariants is exact). One can formally reduce to $\Lambda=\mathbb Z_\ell$ and then to finitely presented solid $\mathbb Z_\ell$-sheaves $\mathcal F$ on $\tilde{\mathcal M}_b^\circ$. Now these can be written as cofiltered inverse limits of constructible $\mathcal F_i$. The $R\Gamma(\tilde{\mathcal M}_b^\circ,\mathcal F_i)$ are uniformly bounded; to see that their derived limit is also bounded, it is then sufficient to see that each $H^j(\tilde{\mathcal M}_b^\circ,\mathcal F_i)$ is finite. By Theorem~\ref{thm:partialcompactsupportvanishing}, this is isomorphic to $H_c^{j+1}(\tilde{\mathcal M}_b^\circ,\mathcal F_i)$. But $R\Gamma_c(\tilde{\mathcal M}_b^\circ,-)$ preserves compact objects as its right adjoint commutes with all colimits (as $\tilde{\mathcal M}_b^\circ$ is cohomologically smooth over $\Spd k$, being open in a successive extension of negative Banach--Colmez spaces).
\end{proof}

Next, we study Bernstein--Zelevinsky duality. Denoting $\pi: \Bun_G\to \ast$ the projection, the pullback $\pi^\ast$ has a left adjoint
\[
\pi_\natural: D_\lis(\Bun_G,\Lambda)\to D_\lis(\ast,\Lambda)\cong D(\Lambda).
\]
This induces a pairing
\[
D_\lis(\Bun_G,\Lambda)\times D_\lis(\Bun_G,\Lambda)\to D(\Lambda): (A,B)\mapsto \pi_\natural(A\soliddotimesLambda B).
\]

\begin{proposition}\label{prop:bernsteinzelevinskylisse} For any compact object $A\in D_\lis(\Bun_G,\Lambda)$, there is a unique compact object $\mathbb D_{\mathrm{BZ}}(A)\in D_\lis(\Bun_G,\Lambda)$ with a functorial identification
\[
R\Hom(\mathbb D_{\mathrm{BZ}}(A),B)\cong \pi_\natural(A\soliddotimesLambda B)
\]
for $B\in D_\lis(\Bun_G,\Lambda)$. Moreover, the functor $\mathbb D_{\mathrm{BZ}}$ is a contravariant autoequivalence of $D_\lis(\Bun_G,\Lambda)^\omega$, and $\mathbb D_{\mathrm{BZ}}^2$ is naturally isomorphic to the identity.

If $U\subset \Bun_G$ is an open substack and $A$ is concentrated on $U$, then so is $\mathbb D_{\mathrm{BZ}}(A)$. In particular, $\mathbb D_{\mathrm{BZ}}$ restricts to an autoequivalence of the compact objects in $D_\lis(\Bun_G^b,\Lambda)\cong D(G_b(E),\Lambda)$ for $b\in B(G)$ basic, and in that setting it is the usual Bernstein--Zelevinsky involution.
\end{proposition}

\begin{proof} The existence of $\mathbb D_{\mathrm{BZ}}$ follows as in Theorem~\ref{thm:bernsteinzelevinsky}, using the left adjoint given by Proposition~\ref{prop:compactobjectsDlis}; this construction also shows that $\mathbb D_{\mathrm{BZ}}$ preserves $D_\lis(U,\Lambda)$, and for basic $b$ it recovers the usual Bernstein--Zelevinsky involution by the same argument as in Theorem~\ref{thm:bernsteinzelevinsky}.

We also formally get a morphism $\mathbb D_{\mathrm{BZ}}^2(A)\to A$ by adjunctions. We need to see that this is an isomorphism. It suffices to check on generators, such as the Bernstein--Zelevinsky dual of $A_K^b$ (which is up to twist and shift $i^b_! c\text-\mathrm{Ind}_K^{G_b(E)} \Lambda$). As in the proof of Theorem~\ref{thm:bernsteinzelevinsky}, one easily checks that the map $\mathbb D_{\mathrm{BZ}}^2(A)\to A$ is an isomorphism over $\Bun_G^b$. To see that it is an isomorphism everywhere, one needs to see that if $B=Rj_\ast B'$, $B'\in D_\lis(U,\Lambda)$ for some open substack $j: U\subset \Bun_G$ not containing $\Bun_G^b$, then
\[
\pi_\natural(A^b_K\soliddotimesLambda B)=0.
\]
Twisting a few things away and using the definition of $A^b_K=f_{K\natural} \Lambda$, this follows from the assertion that for all $A'\in D_\lis(\mathcal M_{b,K}^\circ,\Lambda)$, with $j_K: \mathcal M_{b,K}^\circ\hookrightarrow\mathcal M_{b,K}$ the open immersion, one has
\[
R\Gamma_c(\mathcal M_{b,K},Rj_{K\ast} A')=0.
\]
Using the trace map for $\widetilde{\mathcal M}_b\to \mathcal M_{b,K}$, this follows from Theorem~\ref{thm:partialcompactsupportvanishingsolid}.
\end{proof}

As in Theorem~\ref{thm:localbidual}, this has the following consequence for Verdier duality.

\begin{proposition}\label{proposition:localbiduallisse} Let $j: V\hookrightarrow U$ be an open immersion of open substacks of $\Bun_G$. For any $A\in D_\lis(V,\Lambda)$, the natural map
\[
j_\natural R\sHom_\lis(A,\Lambda)\to R\sHom_\lis(Rj_{\lis\ast} A,\Lambda)
\]
is an isomorphism in $D_\lis(U,\Lambda)$.
\end{proposition}

\begin{proof} The proof is identical to the proof of Theorem~\ref{thm:localbidual}.
\end{proof}

Using this, one can characterize the reflexive objects as in Theorem~\ref{thm:charreflexive}; we omit it here.

Finally, one can also characterize the universally locally acyclic $A\in D_\lis(\Bun_G,\Lambda)$. Note that we have not defined a notion of universal local acyclicity for lisse-\'etale sheaves, but in our present situation we can simply import the characterization from Proposition~\ref{prop:ULAdualizablestack} and make the following definition.

\begin{definition}\label{def:lisseULABunG} A complex $A\in D_\lis(\Bun_G,\Lambda)$ is universally locally acyclic (with respect to $\Bun_G\to \ast$) if the natural map
\[
p_1^\ast R\sHom_\lis(A,\Lambda)\soliddotimesLambda p_2^\ast A\to R\sHom_\lis(p_1^\ast A,p_2^\ast A)
\]
is an isomorphism, where $p_1,p_2: \Bun_G\times \Bun_G\to \Bun_G$ are the two projections.
\end{definition}

We get the following version of Theorem~\ref{thm:ULAbunG}.

\begin{proposition}\label{prop:lisseULAbunG} Let $A\in D_\lis(\Bun_G,\Lambda)$. Then $A$ is universally locally acyclic if and only if for all $b\in B(G)$, the pullback $i^{b\ast} A$ to $i^b: \Bun_G^b\hookrightarrow \Bun_G$ corresponds under $D_\lis(\Bun_G^b,\Lambda)\cong D(G_b(E),\Lambda)$ to a complex $M_b$ of smooth $G_b(E)$-representations for which $M^K$ is a perfect complex of $\Lambda$-modules for all open pro-$p$ subgroups $K\subset G_b(E)$.
\end{proposition}

The proof is identical to the proof of Theorem~\ref{thm:ULAbunG}, and proceeds by proving first the following proposition.

\begin{proposition}\label{prop:lisseBunGproduct} Let $G_1$ and $G_2$ be two reductive groups over $E$, and let $G=G_1\times G_2$. Consider the exterior tensor product
\[
-\boxtimes-: D_\lis(\Bun_{G_1},\Lambda)\times D_\lis(\Bun_{G_2},\Lambda)\to D_\lis(\Bun_G,\Lambda).
\]

For all compact objects $A_i\in D_\lis(\Bun_{G_i},\Lambda)$, $i=1,2$, the exterior tensor product $A_1\boxtimes A_2\in D_\et(\Bun_G,\Lambda)$ is compact, these objects form a class of compact generators, and for all further objects $B_i\in D_\lis(\Bun_{G_i},\Lambda)$, $i=1,2$, the natural map
\[
R\Hom(A_1,B_1)\dotimes_\Lambda R\Hom(A_2,B_2)\to R\Hom(A_1\boxtimes A_2,B_1\boxtimes B_2)
\]
is an isomorphism.
\end{proposition}

\begin{proof} The proof is identical to the proof of Proposition~\ref{prop:BunGproduct}.
\end{proof}

\chapter{$L$-parameter}

It is time to understand the other side of the correspondence: In this chapter, we define, and study basic properties of, the stack of $L$-parameters. These results have recently been obtained by Dat--Helm--Kurinczuk--Moss \cite{DatHelmKurinczukMoss}, and also Zhu \cite{ZhuConjectures}; previous work in a related direction includes \cite{Helm}, \cite{HartlHellmann}, \cite{BellovinGee}, \cite{BooherPatrikis}, \cite[Appendix E]{LiuTianXiaoZhangZhu}.

In this chapter, we fix again a nonarchimedean local field $E$ with residue field $\Fq$ of characteristic $p$, and a reductive group $G$ over $E$, as well as a prime $\ell\neq p$. We get the dual group $\hat{G}/\mathbb Z_\ell$, which we endow with its usual ``algebraic'' action by $W_E$; the action thus factors over a finite quotient $Q$ of $W_E$, and we fix such a quotient $Q$ of $W_E$. (The difference to the cyclotomically twisted $W_E$-action disappears after base change to $\mathbb Z_\ell[\sqrt{q}]$, and we could thus obtain analogues of all results below for this other action by a simple descent along $\mathbb Z_\ell[\sqrt{q}]/\mathbb Z_\ell$.) We define a scheme whose $\Lambda$-valued points, for a $\mathbb Z_\ell$-algebra $\Lambda$, are the condensed $1$-cocycles
\[
\varphi: W_E\to \hat{G}(\Lambda),
\]
where $\Lambda=\Lambda_{\mathrm{disc}}\otimes_{\mathbb Z_{\ell,\mathrm{disc}}} \mathbb Z_\ell$ is regarded as a relatively discrete condensed $\mathbb Z_\ell$-module.

\begin{theorem}[Theorem~\ref{thm:Lparameterrepresentable}] There is a scheme $Z^1(W_E,\hat{G})$ over $\mathbb Z_\ell$ whose $\Lambda$-valued points, for a $\mathbb Z_\ell$-algebra $\Lambda$, are the condensed $1$-cocycles
\[
\varphi: W_E\to \hat{G}(\Lambda).
\]
The scheme $Z^1(W_E,\hat{G})$ is a union of open and closed affine subschemes $Z^1(W_E/P,\hat{G})$ as $P$ runs through open subgroups of the wild inertia subgroup of $W_E$, and each $Z^1(W_E/P,\hat{G})$ is a flat local complete intersection over $\mathbb Z_\ell$ of dimension $\dim G$.
\end{theorem}

To prove the theorem, following \cite{DatHelmKurinczukMoss} and \cite{ZhuConjectures} we define discrete dense subgroups $W\subset W_E/P$ by discretizing the tame inertia, and the restriction $Z^1(W_E/P,\hat{G})\to Z^1(W,\hat{G})$ is an isomorphism, where the latter is clearly an affine scheme.

We can also prove further results about the $\hat{G}$-action on $Z^1(W_E,\hat{G})$, or more precisely each $Z^1(W_E/P,\hat{G})$.

\begin{theorem}[Theorem~\ref{thm:precisecoliminIndPerf}] Assume that $\ell$ does not divide the order of $\pi_1(\hat{G})_{\mathrm{tor}}$. Then $H^i(\hat{G},\mathcal O(Z^1(W_E/P,\hat{G})))=0$ for $i>0$ and the formation of the invariants $\mathcal O(Z^1(W_E/P,\hat{G}))^{\hat{G}}$ commutes with any base change. The algebra $\mathcal O(Z^1(W_E/P,\hat{G}))^{\hat{G}}$ admits an explicit presentation in terms of excursion operators,
\[
\mathcal O(Z^1(W_E/P,\hat{G}))^{\hat{G}} = \colim_{(n,F_n\to W)} \mathcal O(Z^1(F_n,\hat{G}))^{\hat{G}}
\]
where the colimit runs over all maps from a free group $F_n$ to $W\subset W_E/P$, and $Z^1(F_n,\hat{G})\cong \hat{G}^n$ with the simultaneous twisted $\hat{G}$-conjugation.

Moreover, the $\infty$-category $\Perf(Z^1(W_E/P,\hat{G})/\hat{G})$ is generated under cones, shifts and retracts by the image of $\Rep(\hat{G})\to \Perf(Z^1(W_E/P,\hat{G})/\hat{G})$, and $\Ind\Perf(Z^1(W_E/P,\hat{G}))$ is equivalent to the $\infty$-category of modules over $\mathcal O(Z^1(W_E/P,\hat{G}))$ in $\Ind\Perf(B\hat{G})$.

All of these results also hold with $\mathbb Q_\ell$-coefficients, without the assumption on $\ell$.
\end{theorem}

With $\mathbb Q_\ell$-coefficients, these results are simple, as the representation theory of $\hat{G}$ is semisimple. However, with $\mathbb Z_\ell$-coefficients, these results are quite subtle, and we need to dive into modular representation theory of reductive groups. We prove in particular the following result. The last part of this generalizes results of Brundan \cite{Brundan} and van der Kallen \cite{vanderKallenBrundan} that treat the case $P=\mathbb Z/2\mathbb Z$ of involutions. While their argument is case-by-case, we are able to give a conceptual argument.

\begin{theorem}[Section~\ref{sec:fixedpoints}] Let $G$ be a reductive group over an algebraically closed field $L$ of characteristic $\ell$. Let $P$ be a finite group of order prime to $\ell$ acting on $G$. Then $H=G^P$ is a smooth linear algebraic group whose connected component $H^\circ$ is reductive, and with $\pi_0 H$ of order prime to $\ell$. If $P$ is solvable, the image of $\Perf(\ast/G)\to \Perf(\ast/H)$ generates the whole category under cones and retracts. Moreover, still under the assumption that $P$ is solvable, $H^\circ\subset G$ is a Donkin subgroup, i.e.~for any representation $V$ of $G$ that admits a good $G$-filtration, also $V|_{H^\circ}$ has a good $H^\circ$-filtration.
\end{theorem}

\section{The stack of $L$-parameters}
\subsection{Definition and representability}

Recall that for a reductive group $G$ over a nonarchimedean local field $E$, we have the (pinned) dual group $\hat{G}$ over $\mathbb Z_\ell$, equipped with an action of the Weil group $W_E$. In this chapter, we use the standard action (compatible with the pinning).

Now let $\Lambda$ be any $\mathbb Z_\ell$-algebra. As in the last chapter, we regard it as a condensed $\mathbb Z_\ell$-algebra, as $\Lambda_{\mathrm{disc}}\otimes_{\mathbb Z_{\ell,\mathrm{disc}}} \mathbb Z_\ell$. Its value on a profinite set $S$ is the ring of maps $S\to \Lambda$ that take values in a sub-$\mathbb{Z}_\ell$-module of finite type and are continuous. For example, if $\Lambda= \overline{\mathbb{Q}}_\ell$ then $\Lambda(S)= \varinjlim_{L\subset \overline{\mathbb{Q}}_\ell} \mathrm{Cont}(S,L)$ with $L|\mathbb{Q}_\ell$ finite.

\begin{definition}\label{def:Lparameter} An $L$-parameter for $G$, with coefficients in $\Lambda$, is a section
\[
\varphi: W_E\to \hat{G}(\Lambda)\rtimes W_E
\]
of the natural map of condensed groups
\[
\hat{G}(\Lambda)\rtimes W_E\to W_E.
\]
Equivalently, an $L$-parameter for $G$ with coefficients in $\Lambda$ is a (condensed) $1$-cocycle
\[
\varphi: W_E\to \hat{G}(\Lambda)
\]
for the given $W_E$-action on $\hat{G}$.
\end{definition}

More concretely, an $L$-parameter with values in $\Lambda$ is a $1$-cocycle $\phi:W_E\to \hat{G}(\Lambda)$ such that if $\widehat{G}\hookrightarrow \GL_N$,  the associated map $W_E\to \GL_N (\Lambda)$ is continuous. The preceding means  the matrix coefficients of its restriction to $I_E$ are maps $I_E\to \Lambda$ that take values in finite type $\mathbb{Z}_\ell$-modules and are continuous. 

\begin{remark} The standard action of $W_E$ factors over a finite quotient $Q$. This means that $L$-parameters are also equivalent to maps $W_E\to \hat{G}(\Lambda)\rtimes Q$ lifting $W_E\to Q$.
\end{remark}

The first main result is the following.

\begin{theorem}\label{thm:Lparameterrepresentable} There is a scheme $Z^1(W_E,\hat{G})$ over $\mathbb Z_\ell$ parametrizing $L$-parameters for $G$, which is a disjoint union of affine schemes of finite type over $\mathbb Z_\ell$. It is flat and a relative complete intersection of dimension $\dim G=\dim \hat{G}$.
\end{theorem}

\begin{proof} Any condensed $1$-cocycle $\varphi: W_E\to \hat{G}(\Lambda)$ is trivial on an open subgroup of the wild inertia subgroup $P_E$; note also that $P_E$ acts on $\hat{G}$ through a finite quotient. Moreover, for any $\gamma\in P_E$ acting trivially on $\hat{G}$, the locus where $\varphi(\gamma)=1$ is open and closed: Taking a closed embedding $\hat{G}\hookrightarrow \GL_N$, this follows from $A=1$ being a connected component of the locus of all $A\in \GL_N$ such that $A^{p^r}=1$, as can be checked by observing that the tangent space at $A=1$ is trivial. It follows that the moduli space of $L$-parameters decomposes as a disjoint union of open and closed subspaces according to the kernel of $\varphi$ on $P_E$.

Thus, fix now some quotient $W_E\to W_E'$ by an open subgroup of $P_E$ such that the action of $W_E$ on $\hat{G}$ factors over $W_E'$. We are interested in the moduli space of condensed $1$-cocycles $W_E'\to \hat{G}(\Lambda)$. Inside $W_E'$, we look at the discrete dense subgroup $W\subset W_E'$ generated by the image of $P_E$, a choice of generator of the tame inertia $\tau$, and a choice of Frobenius $\sigma$. Thus, $W$ sits in an exact sequence
\[
0\to I\to W\to \sigma^{\Z}\to 0
\]
where $I$ in turn sits in an exact sequence
\[
0\to P\to I\to \tau^{\Z[\tfrac 1p]}\to 0
\]
where $P$ is a finite $p$-group. Moreover, in $W/P$, the elements $\tau$ and $\sigma$ satisfy the commutation $\sigma^{-1}\tau\sigma = \tau^q$.

Now observe that any condensed $1$-cocycle $W_E'\to \hat{G}(\Lambda)$ is already determined by its restriction to the discrete group $W$, as $\hat{G}(\Lambda)$ is quasiseparated and $W\subset W_E'$ is dense. Conversely, we claim that any $1$-cocycle $W\to \hat{G}(\Lambda)$ extends uniquely to a condensed $1$-cocycle $W_E'\to \hat{G}(\Lambda)$. To check this, we may replace $E$ by a finite extension; we can thus pass to a setting where the action of $W_E'$ on $\hat{G}$ is trivial, and where $P=1$. Taking a closed immersion $\hat{G}\hookrightarrow \GL_N$, it then suffices to see that any representation of $\tau^{\Z[\tfrac 1p]}\rtimes \sigma^{\Z}$ on a finite free $\Lambda$-module extends uniquely to a representation of the condensed group $\hat{\mathbb Z}^p\rtimes \sigma^{\Z}$. For this, in turn, it suffices to see that for any $A\in \GL_N(\Lambda)$ such that $A$ is conjugate to $A^q$, the map
\[
\Z\to \GL_N(\Lambda): n\mapsto A^n
\]
extends uniquely to $\hat{\mathbb Z}^p$. The assumption on $A$ implies that all eigenvalues of $A$ at all geometric points of $\Spec \Lambda$ are roots of unity of order prime to $p$; replacing $A$ by a prime-to-$p$-power (as we may) we can thus reduce to the case that $A$ is unipotent, i.e.~$A-1$ is nilpotent. But then $n\mapsto A^n$ extends to a continuous map
\[
n\mapsto A^n=(1+(A-1))^n = \sum_{i\geq 0} \binom{n}{i}(A-1)^i,
\]
defining a map $\mathbb Z_\ell\to \GL_N(\Lambda)$ (and hence $\hat{\mathbb Z}^p\to \mathbb Z_\ell\to \GL_N(\Lambda)$).

Thus, we need to see that the space $X=Z^1(W,\hat{G})$ of all $1$-cocycles $\varphi: W\to \hat{G}(\Lambda)$ is an affine scheme of finite type over $\mathbb Z_\ell$ that is flat and a relative complete intersection of dimension $\dim \hat{G}$. It is clear that it is an affine scheme of finite type over $\mathbb Z_\ell$ as $W$ is discrete and finitely generated.

To prove the geometric properties, we find it slightly more convenient to argue with the Artin stack $[X/\hat{G}]$, which we aim to prove is flat and a relative complete intersection of dimension $0$ over $\mathbb Z_\ell$.

We can understand the deformation theory of $[X/\hat{G}]$: If $\Lambda$ is a field, then the obstruction group is $H^2(W,\hat{\mathfrak{g}}\otimes_{\mathbb Z_\ell} \Lambda)$ (where $\hat{\mathfrak{g}}$ is the Lie algebra of $\hat{G}$), the tangent space is $H^1(W,\hat{\mathfrak{g}}\otimes_{\mathbb Z_\ell} \Lambda)$, and the infinitesimal automorphisms are $H^0(W,\hat{\mathfrak{g}}\otimes_{\mathbb Z_\ell} \Lambda)$, where in all cases the action of $W$ is twisted by the local $1$-cocycle $\varphi$. Now note that by direct computation the prime-to-$p$ cohomological dimension of $W$ is $2$, and the Euler characteristic of any representation is equal to $0$. Thus, this analysis shows that we only have to prove that all fibres of $[X/\hat{G}]\to \Spec \mathbb Z_\ell$ are of dimension at most $0$.

Note that $X$ is actually naturally defined over $\mathbb Z[\tfrac 1p]$ (as $\hat{G}$ is, and the discretization $W$ of $W_E'$ is independent of $\ell$). It follows that it suffices to bound the dimension of the fibre over $\mathbb F_\ell$ (as if we can do this for all closed points of $\Spec \mathbb Z[\tfrac 1p]$, it follows over the generic fibre by constructibility of the dimension of fibers). To do this, we switch back to the picture of condensed $1$-cocycles on $W_E$. From now on, we work over $\overline{\mathbb F}_\ell$.

The stack $[Z^1(W_E,\hat{G})_{\overline{\mathbb F}_\ell}/\hat{G}]$ maps to the similar stack parametrizing $1$-cocycles $\phi^{I^{\ell}}: I^{\ell}\to \hat{G}_{\overline{\mathbb F}_\ell}$ of the prime-to-$\ell$ inertia subgroup $I^{\ell}$, up to conjugation. By deformation theory, that stack is smooth and each connected component is a quotient of $\Spec \overline{\mathbb F}_\ell$ by the centralizer group $C_{\phi^{I^\ell}}\subset \hat{G}_{\overline{\mathbb F}_\ell}$, which is a smooth group, whose identity component is reductive by \cite[Theorem 2.1]{PrasadYu}. We may thus fix $\phi^{I^\ell}: I^\ell\to \hat{G}(\overline{\mathbb F}_\ell)$ and consider the closed subscheme $X_{\phi^{I^\ell}}\subset Z^1(W_E,\hat{G})_{\overline{\mathbb F}_\ell}$ of all $1$-cocycles $\phi: W_E\to \hat{G}(\Lambda)$ whose restriction to $I^\ell$ is equal to $\phi^{I^\ell}$. Our goal is to show that $X_{\phi^{I^\ell}}$ is of dimension at most $\dim C_{\phi^{I^\ell}}$.

Consider the normalizer $\tilde{C}$ of $\phi^{I^\ell}(I^\ell)$ inside $\hat{G}_{\overline{\mathbb F}_\ell}\rtimes Q$. Then $X_{\phi^{I^\ell}}$ maps with finite fibres to the space of maps
\[
f: W_E/I^\ell\cong \mathbb Z_\ell\rtimes \sigma^{\mathbb Z}\longrightarrow \tilde{C}/\phi^{I^\ell}(I^\ell).
\]
Note that, by representability of $X_{\phi^{I^\ell}}$, the universal map $f$ factors over a quotient of the form $\mathbb Z/\ell^m\mathbb Z\rtimes \sigma^{\mathbb Z}$. Finally, we have reduced to Lemma~\ref{lem:lusztigunipotence} below.
\end{proof}

\begin{lemma}\label{lem:lusztigunipotence} Let $H$ be a smooth group scheme over $\overline{\mathbb F}_\ell$ whose identity component is reductive. Then the affine scheme parametrizing maps of groups
\[
\mathbb Z/\ell^m\mathbb Z\rtimes \sigma^{\mathbb Z}\to H,
\]
where $\sigma$ acts on $\mathbb Z/\ell^m \mathbb Z$ via multiplication by $q$, is of dimension at most $\dim H$.
\end{lemma}

\begin{proof} The image of the generator of $\mathbb Z/\ell^m \mathbb Z$ is a unipotent element of $H$. By finiteness of the number of unipotent conjugacy classes, cf.~\cite{LusztigUnipotent}, \cite[Corollary 2.6]{FulmanGuralnick}, we can stratify the scheme according to the conjugacy class of the image of $\tau$. But for each fixed conjugacy class, one has to choose the image of $\sigma$ so as to conjugate $\tau$ into $\tau^q$: This bounds the dimension of each stratum by the dimension of the conjugacy class of $\tau$ (giving the choices for $\tau$) plus the codimension of the conjugacy class of $\tau$ (giving the choices for $\sigma$, for any given $\tau$), which is the dimension of $H$.
\end{proof}

\section{The singularities of the moduli space}

The following proposition was already implicitly noted in the proof of Theorem~\ref{thm:Lparameterrepresentable}.

\begin{proposition}\label{prop:calcul complexe cotangent LocSys}
For any parameter $\varphi:W_E\to \hat{G}(\Lambda)\rtimes Q$ corresponding to $x:\Spec ( \Lambda) \to [Z^1(W_E,\hat{G})/\hat{G}]$, 
$$
x^* \mathbb{L}_{Z^1(W_E,\hat{G})/\hat{G}}^\vee = R\Gamma(W_E, (\hat{\mathfrak{g}} \otimes_{\mathbb{Z}_\ell} \Lambda)_\varphi ) [1]
$$
where $(\hat{\mathfrak{g}}\otimes_{\mathbb{Z}_\ell} \Lambda)_\varphi$ is $\hat{\mathfrak{g}}\otimes_{\mathbb{Z}_\ell} \Lambda$ equipped with the twisted action of $W_E$ deduced from $\varphi$.
\end{proposition}

\begin{proof} This would be clear if we defined the moduli problem on all animated $\mathbb Z_\ell$-algebras, by deformation theory. Then the cohomological dimension of $W_E$ would imply that this moduli problem is a derived local complete intersection, of expected dimension $0$. However, we proved that $Z^1(W_E,\hat{G})/\hat{G}$ is a local complete intersection Artin stack of dimension $0$, hence it represents the correct moduli problem even on all animated $\mathbb Z_\ell$-algebras, thus giving the result.
\end{proof}

\begin{proposition}\label{prop:dual complexe cotangent}
Let $M$ be a free $\Lambda$-module of finite rank equipped with a condensed action of $W_E$. Then $R\Gamma (W_E,M)$ is a perfect complex of $\Lambda$-modules and there is a canonical isomorphism
$$
R\Gamma (W_E,M)^* \cong R\Gamma (W_E,M^*(1))[2].
$$
\end{proposition}

\begin{proof} This follows from Poincar\'e duality applied to $\Div^1\to \ast$, using Proposition~\ref{prop:dualizingcomplexsolid} and the discussion before. It can also be proved by hand, by comparing the $W_E$-cohomology with the $W$-cohomology, for a discretization $W$ of $W_E/P$ as before.
\end{proof}

\begin{corollary}\label{coro:description complexe cotangent par dualite}
For any parameter $\varphi:W_E\to \hat{G}(\Lambda)\rtimes W_E$ corresponding to $x:\Spec ( \Lambda) \to [Z^1(W_E,\hat{G})/\hat{G}]$, 
$$
x^* \mathbb{L}_{Z^1(W_E,\hat{G})/\hat{G}}= R\Gamma (W_E, (\hat{\mathfrak{g}}^* \otimes_{\mathbb{Z}_\ell}\Lambda)_\varphi (1) ) [1]
$$
where $(\hat{\mathfrak{g}}^\ast\otimes_{\mathbb{Z}_\ell} \Lambda)_\varphi$ is $\hat{\mathfrak{g}}^\ast\otimes_{\mathbb{Z}_\ell} \Lambda$ equipped with the twisted action of $W_E$ deduced from $\varphi$.
\end{corollary}

\subsection{The characteristic zero case}\label{sec:the char zero case}

Fix an isomorphism $I_E/P_E \cong \widehat{\mathbb{Z}}^p$. There is a $\hat{G}_{\mathbb{Q}_\ell}$-equivariant ``unipotent monodromy'' morphism
$$
\mathcal{M}:Z^1(W_E,\hat{G})_{\mathbb{Q}_\ell} \lto \mathcal{N}_{\hat{G}_{\mathbb{Q}_\ell}}
$$
where $\mathcal{N}_{\hat{G}_{\mathbb{Q}_\ell}}$ is the nilpotent cone inside $\hat{\mathfrak{g}}\otimes \mathbb{Q}_\ell$.

In fact, one can lift the inclusion $\mathbb{Z}_\ell \hookrightarrow \widehat{\mathbb{Z}}^p\cong I_E/P_E$ to a morphism $\mathbb{Z}_\ell \to I_E$. Now, if $\varphi : W_E\to \hat{G}(\Lambda)$, with $\Lambda$ a $\mathbb {Q}_\ell$-algebra, is a parameter, then $\varphi|_{\mathbb{Z}_\ell} : \mathbb{Z}_\ell \to \hat{G}(\Lambda)$ is such that for $n\gg 0$, $\varphi|_{\ell^n \mathbb{Z}_\ell}$ is a morphism of condensed groups satisfying 
$$\varphi ( \sigma^m)  \varphi|_{\ell^n \mathbb{Z}_\ell} \varphi(\sigma^m)^{-1} =   \varphi_{|\ell^n \mathbb{Z}_\ell}^{q^m} $$ for $m\gg 0$. One deduces, using an embedding of $\hat{G}$ in $\GL_N$, that there is a unique $N\in \mathcal{N}_{\hat{G}}(\Lambda)$ such that for $n\gg 0$ and $x\in \ell^n\mathbb{Z}_\ell$,
$$
\varphi (x) = \exp (x N).
$$

Using these observations, we get a comparison to Weil--Deligne $L$-parameters.

\begin{definition}
For $\Lambda$ a $\mathbb{Q}_\ell$-algebra one defines $\mathrm{Par}^{\mathrm{WD}}_{\hat{G}} (\Lambda)$ to be the set of pairs $(\varphi_0,N)$ where 
\begin{altenumerate}
\item $\varphi_0 : W_E\to \hat{G}(\Lambda)$ is a $1$-cocycle that is continuous for the discrete topology on $\hat{G}(\Lambda)$ (i.e., trivial on an open subgroup of $I_E$),
\item $N\in \mathcal{N}_{\hat{G}}\otimes \Lambda$ satisfies $\mathrm{Ad}(\varphi_0 (\sigma)).\sigma N = q^{|\sigma|} N$ for all $\sigma\in W_E$.
\end{altenumerate}
\end{definition}

Then we have the following result, which is essentially Grothendieck's quasi-unipotence theorem.

\begin{proposition}[{\cite[Lemma~3.1.8]{ZhuConjectures}}]
There is a $\hat{G}$-equivariant isomorphism $$Z^1(W_E,\hat{G})\otimes \mathbb{Q}_\ell \iso \mathrm{Par}_{\hat{G}}^{\mathrm{WD}}.$$
\end{proposition}

However, we warn the reader that this isomorphism depends on some auxiliary choices, such as that of a Frobenius element.

\subsection{The singular support}\label{sec:singularsupport}
\subsubsection{General construction}

Recall the following construction, see for example \cite{ArinkinGaitsgory}. Let $A\to B$ be a flat map of commutative rings. One has the Hochschild cohomology
\[
H\! H^\bullet(B/A) = \Ext^\bullet_{B\otimes_A B}(B,B).
\]
Note that any $M\in D(B\otimes_A B)$ induces a functor $D(B)\to D(B)$, via $N\mapsto M\dotimes_B N$ (with the ``left'' $B$-module structure). Here, $M=B\in D(B\otimes_A B)$, via the multiplication $B\otimes_A B\to B$, induces the identity functor. It follows that there is a natural map
\[
H\! H^i(B/A) = \Ext^i_{B\otimes_A B}(B,B)\to \Ext^i_B(N,N)
\]
for any $N\in D(B)$. Moreover, Hochschild cohomology is naturally a graded algebra, and this map is a map of algebras
\[
H\! H^\bullet(B/A)\to \Ext^\bullet_B(N,N).
\]

There is an identification (\cite[Theorem X.3.1]{MacLaneHomology})
$$
H\! H^2 (B/A) = \Ext^1_B ( \mathbb{L}_{B/A}, B)
$$
which itself is nothing else than $\mathcal{E}\! xalcom_A (B,B)$ (\cite[Chap.0, Sec. 18.4]{EGAIV.1}). We thus have an identification 
$$
H\! H^2 (B/A)= H^1 ( \mathbb{L}_{B/A}^\vee).
$$
{\it Suppose now that $A\to B$ is syntomic, i.e.~flat and a local complete intersection.} Let $X=\Spec B\to S=\Spec A$ be the associated map of affine schemes.

\begin{definition}
 The scheme
$$
\mathrm{Sing}_{X/S}\lto X
$$
represents the functor $T/X\mapsto H^{-1} ( \mathbb{L}_{X/S}\dotimes_{\mathcal O_X} \mathcal O_T)$.
\end{definition}

 In fact, locally on $X$, $\mathbb{L}_{X/S}$ is isomorphic to a complex of vector bundles $[\E^{-1}\to \E^0]$ and then $\mathrm{Sing}_{X/S}$ is the kernel of $\mathbb{V} (\E^{-1})\to \mathbb{V} ( \E^0)$. Explicitly, $\mathrm{Sing}_{X/S}$ is the affine scheme with
\[
\mathcal O(\mathrm{Sing}_{X/S}) = \mathrm{Sym}^\bullet_B H^1(\mathbb L_{B/A}^\vee).
\]
This is an $X$-group scheme equipped with an action of $\mathbb{G}_m$. {\it The image of $\mathrm{Sing}_{X/S}\setminus \{
0\}\to X$ is the closed subset  complementary of the smooth locus of $X\to S$.}
\\

Consider now any $$N\in D^b_{coh} (X),$$ and the graded $B$-algebra $\Ext^\bullet_B(N,N)$. Using the map
\[
H^1(\mathbb L_{B/A}^\vee)=H\! H^2(B/A)\to \Ext^2_B(N,N),
\]
this is in fact naturally a (graded) $\mathcal O(\mathrm{Sing}_{X/S})$-algebra. This defines a $\mathbb{G}_m$-equivariant quasi-coherent sheaf $$\mu\End(N)$$ on $\mathrm{Sing}_{X/S}$.

{\it Suppose now moreover that $S$ is regular.}

\begin{theorem}[{\cite[Theorem 3.1]{Gulliksen},\cite[Appendix D]{ArinkinGaitsgory}}]
For $N\in D^b_{coh}(X)$, the quasi-coherent sheaf $\mu\End(N)$ on $\mathrm{Sing}_{X/S}$ is coherent.
\end{theorem}
 
\begin{definition}\label{def:singsupp}
The singular support of $N$, $\mathrm{SingSupp}(N)$, is the support of $\mu\End(N)$ as a closed conical subset of $\mathrm{Sing}_{X/S}$. 
\end{definition}

Of course, the image of $\mathrm{SingSupp}(N)\to X$ is contained in $\mathrm{Supp}(N)$.

\begin{theorem}[{\cite[Theorem 4.2.6]{ArinkinGaitsgory}}]
The following are equivalent:
\begin{altenumerate}
\item $N$ is a perfect complex,
\item $\mathrm{SingSupp} (N)$ is contained in the zero section of $\mathrm{Sing}_{X/S}$.
\end{altenumerate}
\end{theorem}

\begin{proof}
We have to prove that if $\Ext^i_B(N,N)=0$ for $i\gg 0$ then $N$ is a perfect complex. This is for example a consequence of \cite{Jorgensen}. Since $S$ is regular $X$ is Gorenstein. According to \cite{Jorgensen}, if $N$ is a $B$-module of finite type that satisfies $\Ext^i_B(N,N) =0$ for $i>n$, then $\mathrm{pd}_B N \leq n$. In general, up to taking a shift of $N$, we can find a map $N\to N'$, where $N'$ is a finitely generated $B$-module concentrated in degree $0$, such that the cone $C$ of $N\to N'$ is perfect. Suppose that $\Ext^i_B(N,N)=0$ for $i\gg 0$. In the long exact sequence 
$$
\cdots\lto \Ext^i_B(C,N)\lto \Ext^i_B(N',N)\lto \Ext^i_B(N,N)\lto \cdots
$$
one has $\Ext^i_B(C,N)=0$ for $i\gg 0$ since $C$ is perfect and $\Ext^i_B(N,N)=0$ for $i\gg 0$ by hypothesis. We deduce that $\Ext^i_B(N',N)=0$ for $i\gg 0$. In the long exact sequence 
$$
\cdots \lto \Ext^i_B(N',N)\lto \Ext^i_B(N',N')\lto \Ext^i_B(N',C)\lto \cdots
$$
we have $\Ext^i_B(N',C)=0$ for $i\gg 0$ since $C$ is perfect and $B$ has finite injective dimension over itself since it is Gorenstein. Thus, for $i\gg 0$, $\Ext^i_B(N',N')\iso \Ext^i_B(N',N)$ and this vanishes. We can thus apply Jorgensen's theorem to $N'$ to conclude that $N'$, and hence $N$, is perfect.
\end{proof}

Let us note the following corollary. 

\begin{corollary}
The image of $\mathrm{SingSupp} (N)\setminus \{0\} \to X$ is the complementary of the biggest open subset of $X$ on which $N$ is a perfect complex. 
\end{corollary}

\subsubsection{The case of $Z^1(W_E,\hat{G})$}

Now we apply the preceding theory in the case $A=\mathbb Z_\ell$ and $X=Z^1(W_E,\hat{G})$ (which is only a union of affine schemes, but this is not a problem). We can also pass to the quotient stack $Z^1(W_E,\hat{G})/\hat{G}$ as the formation of $\mathrm{Sing}$ commutes with smooth maps. According to Corollary~\ref{coro:description complexe cotangent par dualite}, there is an embedding
$$
\begin{tikzcd}[row sep=1.2cm]
\mathrm{Sing}_{[Z^1(W_E,\hat{G})/\hat{G}]/\mathbb Z_\ell} \ar[r,hook]\ar[d]& \ar[ld] [ \hat{\mathfrak{g}}^* / \hat{G}]\times_{\ast/\hat{G}} [Z^1(W_E,\hat{G})/\hat{G}] \\
\! [Z^1(W_E,\hat{G})/\hat{G}]
\end{tikzcd}
$$
where $\hat{\mathfrak{g}}=\Lie \hat{G}$ and $[\hat{\mathfrak{g}}^* /\hat{G}]$ is seen here as a vector bundle
on $\ast/\hat{G}=[\Spec \mathbb{Z}_\ell/\hat{G}]$. Let $\mathcal{N}_{\hat{G}}^* \subset \hat{\mathfrak{g}}^*$ be the nilpotent cone; by this we mean the closed subset of all $\hat{G}$-orbits whose closure contains the origin. (If there is a $\hat{G}$-equivariant isomorphism between $\hat{\mathfrak{g}}^*$ and $\hat{\mathfrak{g}}$, this would identify with the usual nilpotent cone.) Since this is stable under the adjoint action this defines a Zariski closed substack 
$$
\begin{tikzcd}[row sep=1.2cm]
\! [\mathcal{N}_{\hat{G}}^* / \hat{G} ] \times_{\ast/\hat{G}} [Z^1(W_E,\hat{G})/\hat{G}] \ar[r,hook]\ar[d]& \ar[ld] [ \hat{\mathfrak{g}}^* / \hat{G} ]\times_{\ast/\hat{G}} [Z^1(W_E,\hat{G})/\hat{G}]  \\
\! [Z^1(W_E,\hat{G})/\hat{G}] .
\end{tikzcd}
$$

\begin{proposition}
For  a $\mathbb{Z}_\ell$-field $L$ and a point $x:\Spec (L)\to [Z^1(W_E,\hat{G})/\hat{G}]$ we have 
$$
x^*\mathrm{Sing}_{[Z^1(W_E,\hat{G})/\hat{G}] /\mathbb{Z}_\ell}  \subset \mathcal{N}_{\hat{G}}^*\otimes_{\mathbb{Z}_\ell} L
$$
in the following two cases:
\begin{altenumerate}
\item $L|\mathbb{Q}_\ell$, 
\item If $n=f_{E'/E}$ with $W_{E'}=\ker (W_E\to \mathrm{Out} (\hat{G}))$, then $q^{en}-1$ is not divisible by $\ell$ for any exponent $e$ of $\hat{G}$.
\end{altenumerate}
\end{proposition}

\begin{proof} 
Assumption (ii) implies that $\ell$ is a very good prime for $\hat{G}$ and in particular the Chevalley isomorphism $\hat{\mathfrak{g}}\sslash \hat{G} = \hat{\mathfrak{t}}\sslash W$ holds, and there is an isomorphism $\hat{\mathfrak{g}}^\ast\cong \hat{\mathfrak{g}}$.

If $x$ corresponds to the parameter $\varphi$ then $x^*\mathrm{Sing}_{[Z^1(W_E,\hat{G})/\hat{G}] /\mathbb{Z}_\ell}= H^0(W_E,\hat{\mathfrak{g}}^*\otimes_{\mathbb{Z}_\ell} L(1))$ where the $W_E$ action on $\hat{\mathfrak{g}}^*\otimes_{\mathbb{Z}_\ell} L(1)$ is twisted by $\varphi$. For an element $v\in \hat{\mathfrak{g}}\cong \hat{\mathfrak{g}}^*$ in this subspace we thus have that $\sigma.v$ and $qv$ are in the same orbits under the adjoint action of $\hat{G}(L)$ (here $\sigma.v$ is given by the action of $W_E$ on $\hat{\mathfrak{g}}^*$ defining the $L$-group). We thus obtain that $v$ is conjugated under the adjoint action to $q^n v$. There is a morphism
$$
\hat{\mathfrak{g}} \lto \hat{\mathfrak{g}} \sslash \hat{G} = \hat{\mathfrak{t}} \sslash W \cong \mathbb{A}^m_{\mathbb{Z}_\ell}
$$ 
given by $m$ homogeneous polynomials of degrees the exponents of the root system. This implies that the image of $v$ in $\mathbb{A}^m (L)$ is zero and thus $v$ lies in the nilpotent cone.
\end{proof}

The supremum of the exponents of $\hat{G}$ is the Coxeter number $h$ of $G$. The preceding condition is satisfied if for example $\ell > q^{hn}-1$. We refer to \cite[Section 5.3]{DatHelmKurinczukMoss} for finer definitions and results about $\hat{G}$-banal primes; we have not tried to optimize the condition above, and it is likely that with their results one can obtain a much better condition on $\ell$.

\begin{remark}
In the non-banal case things become more complicated and the Arinkin--Gaitsgory condition of nilpotent singular support becomes important. This is also the case when interesting congruences between smooth irreducible representations of $G(E)$ occur, cf.~\cite[Section 1.5]{DatHelmKurinczukMoss}.
\end{remark}

\begin{remark} The appearance of $\hat{\mathfrak{g}}^*$ here is another indication that the assumption that $\ell$ does not divide the order of $\pi_1(\widehat{G})_{\mathrm{tor}}$ may be important: Indeed, this assumption determines the isomorphism class of $\hat{\mathfrak{g}}^*$ as a representation of the adjoint group, within the isogeny class of $G$. However, when $\ell$ is a bad prime, then the nilpotent cone is not well-behaved (for example, there may be infinitely many nilpotent orbits), and we are not sure whether the resulting notion of nilpotent singular support is in fact the correct notion.
\end{remark}

\section{The coarse moduli space}

Let us now describe the corresponding coarse moduli space, i.e.~we consider the quotient
\[
Z^1(W_E,\hat{G})\sslash\hat{G}
\]
taken in the category of schemes. Concretely, for every connected component $\Spec A\subset Z^1(W_E,\hat{G})$, we get a corresponding connected component $\Spec A^{\hat{G}}\subset Z^1(W_E,\hat{G})\sslash\hat{G}$.

\subsection{Geometric points}

For any algebraically closed field $L$ over $\mathbb Z_\ell$, the $L$-valued points of $Z^1(W_E,\hat{G})_L\sslash \hat{G}$ are in bijection with the closed $\hat{G}$-orbits in $Z^1(W_E,\hat{G})_L$.

We want to describe $L$-valued points with closed $\hat{G}$-orbit as the ``semisimple'' parameters. For this, recall (cf.~\cite{Borel79}) that parabolic subgroups of $\hat{G}_L\rtimes W_E$ surjecting onto $W_E$ are up to $\hat{G}(L)$-conjugation given by $\hat{P}_L\rtimes W_E$ for a standard parabolic $P\subset G^\ast$ of the quasisplit inner form $G^\ast$ of $G$. A Levi subgroup is given by $\hat{M}_L\rtimes W_E$ where $M\subset P$ is the standard Levi. {\it We now call them the parabolic subgroups of $\hat{G}\rtimes W_E$ i.e. we always suppose they surject to $W_E$.}
If $\hat{\Delta}$ are the simple roots of $\hat{G}$ then the standard parabolic subgroups are in bijection with the finite $W_E$-stable subsets of $\hat{\Delta}$.

\begin{definition}\label{def:semisimpleLparameter} Let $L$ be an algebraically closed field over $\mathbb Z_\ell$. An $L$-parameter $\varphi: W_E\to \hat{G}(L)\rtimes W_E$ is semisimple if whenever the image of $\varphi$ is contained in a parabolic subgroup of $\hat{G}\rtimes W_E$ then it is contained in a Levi subgroup of this parabolic subgroup.
\end{definition}

In terms of the standard parabolic subgroups this means that if some $\hat{G}(L)$-conjugate $\varphi'$ of $\varphi$ factorizes through $\hat{P}(L)\rtimes W_E$, then there exists $g\in \hat{P}(L)$ such that
\[
g \varphi' g^{-1}= i_{{}^L M}\mathrm{pr}_{{}^L M} \circ \varphi',
\]
where $\mathrm{pr}_{{}^L M}:\hat{P}(L)\rtimes W_E \to \hat{M}(L)\rtimes W_E$ is the projection onto the standard Levi subgroup, and $i_{{}^L M}: \hat{M}(L)\rtimes W_E\to \hat{G}(L)\rtimes W_E$ the inclusion.

\begin{proposition}[{\cite[Proposition 4.13]{DatHelmKurinczukMoss}}]\label{prop:coarsemodulisemisimple} Let $L$ be an algebraically closed field over $\mathbb Z_\ell$ and $\varphi: W_E\to \hat{G}(L)\rtimes W_E$ a parameter. The following are equivalent:
\begin{altenumerate}
\item  The $\hat{G}$-orbit of $\varphi$ in $Z^1(W_E,\hat{G})_L$ is closed.
\item  For any conjugate $\varphi'$ of $\varphi$ such that $\varphi':W_E\to \hat{P}(L)\rtimes W_E$ factors over a standard parabolic subgroup, $\varphi$ is $\hat{G}(L)$-conjugate to $i_{{}^L M} \mathrm{pr}_{{}^L M} \circ \varphi'$.
\item $\varphi$ is semi-simple.
\end{altenumerate}
\end{proposition}

\begin{proof} We use the Hilbert--Mumford--Kempf theorem, cf.~\cite[Corollary 3.5]{Kempf78}. Recall that this criterion says that an orbit is closed if and only if any degeneration along a $1$-parameter family induced by a $\mathbb G_m$ has limit inside the same orbit. Thus take any $\lambda: \mathbb G_m\to \hat{G}_L$. Up to conjugation (which one can move into a conjugation of the parameter) we can assume $\lambda \in X_*(\hat{T})^+$. For each $\tau \in W_E$ there is a morphism $\mathrm{ev}_\tau: Z^1(W_E,\hat{G})_L \to \hat{G}$ given by evaluating a parameter on $\tau$. Thus, if $\lim_{t\to 0} \lambda (t)\cdot \varphi$ exists, i.e. the associated morphism $\mathbb{G}_{m,L} \to Z^1(W_E,\hat{G})_L$ extends to $\mathbb{A}^1_L$, for each $\tau \in W_E$ one has $\lambda^\tau =\lambda$ and $\varphi (\tau) \in Q_\lambda(L)\rtimes \tau$, cf.~Lemma~\ref{lem:limite dynamique parametre}. One thus has $Q_\lambda = \hat{P}$ for $P$ a standard parabolic subgroup of $G^\ast$, and $\varphi:W_E\to \hat{P}\rtimes W_E$. 

For $g\in \hat{P}$, $\lim_{t\to 0} \lambda (t)g\lambda (t)^{-1}$ is the projection onto the standard Levi subgroup $\hat{M}$. Thus, using the evaluation morphism $\mathrm{ev}_\tau$ for each $\tau$ we deduce that $\lim_{t\to 0} \lambda(t)\cdot \varphi$, if it exists, is given by the composite $W_E\xrightarrow{\varphi} \hat{P}(L)\rtimes W_E \xrightarrow{\mathrm{proj}} \hat{M}(L)\rtimes W_E$. Reciprocally, since the morphism $\mathbb{G}_m\times \hat{P}\to \hat{P}$, given by $(t,g)\mapsto \lambda (t)g\lambda (t)^{-1}$ extends to $\mathbb{A}^1\times \hat{P}$ with fiber over $0\in \mathbb{A}^1$ given by the projection to $\hat{M}$, for any $\varphi : W_E\to  \hat{P}\rtimes W_E$, $\lim_{t\to 0} \lambda (t)\cdot \varphi$ exists.

From this analysis we deduce the equivalence between (i) and (ii). It is clear that (iii) implies (ii). For the proof of (ii) implies (iii) we use the results of \cite{BateMartinRoehrle} and \cite{Richardson88}. For this we see parameters as morphisms $W\to \hat{G}(\Lambda)\rtimes Q$ where $W$ is discrete finitely generated as in the proof of Theorem~\ref{thm:Lparameterrepresentable}, and $Q$ is a finite quotient of $W$. Let $\varphi:W\to \hat{G}(L)\rtimes Q$ satisfying (ii). Let $H\subset \hat{G}_L\rtimes Q$ be the Zariski closure of the image of $\varphi$.
 Then if $(x_1,\dots,x_n)\in (\hat{G}(L)\times Q)^n$ are the images of a set of generators of $W$, applying the Hilbert--Mumford--Kempf criterion we see that the $\hat{G}_L$-orbit of $(x_1,\dots,x_n)$ via the diagonal action is closed, cf.~the proof of \cite[Lemma 2.17]{BateMartinRoehrle}. We can then apply \cite{Richardson88}, cf. \cite[Proposition 2.16]{BateMartinRoehrle}, to deduce that $H$ is strongly reductive in $\hat{G}_L\rtimes Q$ and thus $\hat{G}_L$-completely reducible. Strictly speaking, since we are working in a non-connected situation, we use in fact \cite[Section 6]{BateMartinRoehrle}.
\end{proof}

\begin{lemma}\label{lem:limite dynamique parametre}
For $\lambda\in X_\ast(T)^+$ and $g\rtimes \tau \in \hat{G}(L)\rtimes W_E$, the limit $\lim_{t\to 0} \lambda(t) g \lambda(t)^{-\tau}$ exists if and only if $g\in Q_{\lambda} (L)$, the parabolic subgroup attached to $\lambda$, and $\lambda^\tau=\lambda$.
\end{lemma}

\begin{proof} The subgroups $Q_{\lambda}$ and $Q_{\lambda^\tau}$ are standard parabolic subgroups of $\hat{G}$. Let us write $g=g' \overset{.}{w} g''$ with $g'\in Q_{\lambda} (L)$, $g''\in Q_{\lambda^\tau}(L)$ and $w\in \hat{W}$. Then, writing
\[
\lambda(t) g \lambda(t)^{-\tau} = ( \lambda(t) g' \lambda (t)^{-1}) (\lambda (t) \overset{.}{w} \lambda (t)^{-\tau}  ) ( \lambda (t)^\tau g'' \lambda (t)^{-\tau}),
\]
one deduces that $\lim_{t\to 0} \lambda (t) \overset{.}{w} \lambda (t)^{-\tau}$ exists. Thus, $\lim_{t\to 0} (\lambda (\lambda^{-\tau} )^w )(t)$ exists and thus $\lambda =  (\lambda^{\tau})^w$. Since $\lambda^\tau \in X_*(\hat{T})^+$ we deduce $\lambda =\lambda^{\tau}$ and $\lambda^w=\lambda$.
\end{proof}

The proof shows that up to replacing $\hat{G}(L)\rtimes W_E$ by $\hat{G}(L)\rtimes Q$ for some finite quotient $Q$ of $W_E$ (as we can), semisimplicity of $\varphi$ is equivalent to the Zariski closure of the image of $\varphi$ being completely reducible in the terminology of \cite[Section 6]{BateMartinRoehrle}.

\subsection{A presentation of $\O(Z^1(W_E,\hat{G}))$}

It will be useful to have a presentation of the algebra $\O(Z^1(W_E,\hat{G}))$, or rather of the finite type $\mathbb Z_\ell$-algebras $\O(Z^1(W_E/P,\hat{G}))$ for open subgroups $P$ of the wild inertia (with the property that the action of $W_E$ on $\hat{G}$ factors over $W_E/P$). Pick a discrete dense subgroup $W\subset W_E/P$ as above, so that $Z^1(W_E/P,\hat{G})=Z^1(W,\hat{G})$. For any $n\geq 0$ with a map $F_n\to W$, we get a $\hat{G}$-equivariant map
\[
\mathcal O(Z^1(F_n,\hat{G}))\to \mathcal O(Z^1(W,\hat{G})),
\]
where the source is isomorphic to $\mathcal O(\hat{G}^n)$ with appropriately twisted diagonal $\hat{G}$-conjugation. Consider the category $\{(n,F_n\to W)\}$ consisting of maps from finite free groups to $W$, with maps given by commutative diagrams $F_n\to F_m\to W$; this is a sifted index category (as it admits coproducts). The map
\[
\colim_{(n,F_n\to W)} \mathcal O(Z^1(F_n,\hat{G}))\to \mathcal O(Z^1(W,\hat{G}))
\]
is an isomorphism of algebras with $\hat{G}$-action (as $1$-cocycles from $W$ to $\hat{G}$ are uniquely specified by compatible collections of $1$-cocycles $F_n\to \hat{G}$ for all $F_n\to W$). By Haboush's theorem on geometric reductivity \cite{Haboush} it follows that the map
\[
\colim_{(n,F_n\to W)} \mathcal O(Z^1(F_n,\hat{G}))^{\hat{G}}\to \mathcal O(Z^1(W,\hat{G}))^{\hat{G}}
\]
on $\hat{G}$-invariants is a universal homeomorphism of finite type $\mathbb Z_\ell$-algebras, and an isomorphism after inverting $\ell$.

\begin{definition}\label{def:algebraofexcursionoperators} The algebra of excursion operators (for $Z^1(W,\hat{G})$) is
\[
\mathrm{Exc}(W,\hat{G}) = \colim_{(n,F_n\to W)} \mathcal O(Z^1(F_n,\hat{G}))^{\hat{G}}.
\]
\end{definition}

We see in particular that the geometric points of $\mathrm{Exc}(W,\hat{G})$ and $Z^1(W,\hat{G})$ agree.

Actually, the following higher-categorical variant is true.

\begin{proposition}\label{prop:colimunderlying} Working in the derived $\infty$-category $\mathcal D(\mathbb Z_\ell)$, the map
\[
\colim_{(n,F_n\to W)} \mathcal O(Z^1(F_n,\hat{G}))\to \mathcal O(Z^1(W,\hat{G}))
\]
is an isomorphism in $\mathcal D(\mathbb Z_\ell)$.
\end{proposition}

In fact, both sides naturally admit the structure of animated $\mathbb Z_\ell$-algebras, and the map is a morphism of such. The proposition then implies that it is in fact an isomorphism of animated $\mathbb Z_\ell$-algebras.

\begin{proof} The left-hand side defines an animated $\mathbb Z_\ell$-algebra, in fact the universal animated $\mathbb Z_\ell$-algebra $A$ with a $1$-cocycle $W\to \hat{G}(A)$, and the right-hand side is given by $\pi_0 A$. Now the deformation-theoretic arguments from the proof of Theorem~\ref{thm:Lparameterrepresentable} show that $A$ is a derived complete intersection, but as $\pi_0 A$ has the correct dimension, we get $A=\pi_0 A$.
\end{proof}

We will later prove an even finer version, incorporating the $\hat{G}$-action; we defer the proof to Section~\ref{sec:integralcoefficients}.

\begin{theorem}\label{thm:coliminIndPerf} Assume that $\ell$ does not divide the order of $\pi_1(\hat{G})_{\mathrm{tor}}$. Then the map
\[
\colim_{(n,F_n\to W)} \mathcal O(Z^1(F_n,\hat{G}))\to \mathcal O(Z^1(W,\hat{G}))
\]
is an isomorphism in the presentable stable $\infty$-category $\Ind \Perf(\ast/\hat{G})$.\footnote{This map is also a map (and hence isomorphism) of $E_\infty$-algebras in $\Ind\Perf(\ast/\hat{G})$, but the question whether it is an isomorphism does not depend on the algebra structure.}

In particular, the map
\[
\mathrm{Exc}(W,\hat{G})=\colim_{(n,F_n\to W)} \mathcal O(Z^1(F_n,\hat{G}))^{\hat{G}}\to \mathcal O(Z^1(W,\hat{G}))^{\hat{G}}
\]
is an isomorphism.
\end{theorem}

In particular, we see that the algebra of excursion operators gives a presentation of $\mathcal O(Z^1(W_E/P,\hat{G}))^{\hat{G}}$. In the next subsection, we analyze it more explicitly.

\subsection{The algebra of excursion operators}

Fix a finite quotient $Q$ of $W_E$ over which the $W_E$-action on $\hat{G}$ factors. Let $(\hat{G}\rtimes Q)^n\sslash\hat{G}$ be the quotient of $(\hat{G}\rtimes Q)^n$ under simultaneous conjugation by $\hat{G}$.

\begin{proposition}\label{prop:descrexcursionoperators} The algebra of excursion operators $\mathrm{Exc}(W,\hat{G})$ is the universal $\mathbb Z_\ell$-algebra $A$ equipped with maps
\[
\Theta_n: \mathcal O((\hat{G}\rtimes Q)^n\sslash\hat{G})\to \Map(W^n,A)
\]
for $n\geq 1$, linear over $\mathcal O(Q^n)\to \Map(W^n,A)$, subject to the following relations. If $g: \{1,\ldots,m\}\to \{1,\ldots,n\}$ is any map, the induced diagram
\[\xymatrix{
\mathcal O((\hat{G}\rtimes Q)^m\sslash\hat{G})\ar[r]\ar[d] & \Map(W^m,A)\ar[d]\\
\mathcal O((\hat{G}\rtimes Q)^n\sslash\hat{G})\ar[r] & \Map(W^n,A)
}\]
commutes, where both vertical maps are the natural pullback maps. On the other hand, $g$ also induces a map $(\hat{G}\rtimes Q)^m\to (\hat{G}\rtimes Q)^n$, multiplying in every fibre over $i=1,\ldots,n$ the terms in $g^{-1}(i)$ (ordered by virtue of their ordering as a subset of $\{1,\ldots,m\}$). This map is equivariant under diagonal $\hat{G}$-conjugation, and hence descends to the quotient. Similarly, $g$ induces a map $W^m\to W^n$. Then also the induced diagram
\[\xymatrix{
\mathcal O((\hat{G}\rtimes Q)^n\sslash\hat{G})\ar[r]\ar[d] & \Map(W^n,A)\ar[d]\\
\mathcal O((\hat{G}\rtimes Q)^m\sslash\hat{G})\ar[r] & \Map(W^m,A)
}\]
commutes.

The $\ell$-torsion free quotient of $\mathrm{Exc}(W,\hat{G})$ is also the universal flat $\mathbb Z_\ell$-algebra $A'$ equipped with maps
\[
\Theta_n': \mathcal O((\hat{G}\rtimes Q)^n\sslash\hat{G})\to \Map((W_E/P)^n,A')
\]
for $n\geq 1$, linear over $\mathcal O(Q^n)\to \Map((W_E/P)^n,A')$, satisfying the same relations as in Proposition~\ref{prop:descrexcursionoperators}, where the right-hand side $\Map((W_E/P)^n,A')$ denotes the maps of condensed sets (where as usual $A'$ is considered as relatively discrete over $\mathbb Z_\ell$). In particular, the $\ell$-torsion free quotient of $\mathrm{Exc}(W,\hat{G})$ is independent of the discretization $W$ of $W_E/P$.
\end{proposition}

We do not know whether it is necessary to pass to the $\ell$-torsion free quotient for the final assertion. Note that if $\ell$ does not divide the order of $\pi_1(\hat{G})_{\mathrm{tor}}$, then $\mathrm{Exc}(W,\hat{G})\cong \O(Z^1(W_E/P,\hat{G}))^{\hat{G}}$ is flat over $\mathbb Z_\ell$. Moreover note that the $\ell$-torsion in $\mathrm{Exc}(W,\hat{G})$ is always nilpotent, so passing to this quotient is a universal homeomorphism.

\begin{proof} The datum of the $\Theta_n$ is equivalent the datum of a map of algebras
\[
\mathcal O(Z^1(F_n,\hat{G}))^{\hat{G}}\to A
\]
for each map $F_n\to W$. The relations encode the relations arising in the diagram category $(n,F_n\to W)$ corresponding to maps $F_m\to F_n$ (over $W$) either sending generators to generators, or multiplying subsets of generators. If one would also allow the inversion of elements, then this would generate all required relations. We leave it as an exercise to see that this relation, corresponding to $F_n\to F_n$ which is the identity on the first $n-1$ generators and inverts the $n$-th generator, is in fact enforced by the others. (Hint: Look at the part of $\Theta_{n+1}$ corresponding to $(\gamma_1,\ldots,\gamma_n,\gamma_n^{-1})$ and use that under multiplication of the last two variables, this maps to $(\gamma_1,\ldots,\gamma_{n-1},1)$, which arises from $(\gamma_1,\ldots,\gamma_{n-1})$.)

The second description is a priori stronger as $\Map((W_E/P)^n,A')$ injects into $\Map(W^n,A')$ as $W\subset W_E/P$ is dense. The $\ell$-torsion free quotient of $\mathrm{Exc}(W,\hat{G})$ injects into $\O(Z^1(W_E/P,\hat{G}))^{\hat{G}}$ (as we have an isomorphism after inverting $\ell$), and by density of $W\subset W_E/P$ the elements of $\Map(W^n,\mathrm{Exc}(W,\hat{G}))$ map to elements of $\Map((W_E/P)^n,\O(Z^1(W_E/P,\hat{G}))^{\hat{G}})$ (we only need to check the integrality). Thus, this already happens on the $\ell$-torsion free quotient of $\mathrm{Exc}(W,\hat{G})$, which thus has the desired universal property.
\end{proof}

Regarding the passage to $W_E$ in place of $W_E/P$, where there is no natural (finite type) algebra anymore, we still have the following result.

\begin{proposition}\label{prop:excursionoperatorsgiveLparameters} Let $L$ be an algebraically closed field over $\mathbb Z_\ell$. Then the following are in canonical bijection.
\begin{altenumerate}
\item[{\rm (i)}] Semisimple $L$-parameters $\varphi: W_E\to \hat{G}(L)\rtimes W_E$, up to $\hat{G}(L)$-conjugation.
\item[{\rm (ii)}] $L$-valued points of $Z^1(W_E,\hat{G})\sslash\hat{G}$.
\item[{\rm (iii)}] Collections of maps of $\mathbb Z_\ell$-algebras
\[
\Theta_n: \mathcal O((\hat{G}\rtimes Q)^n\sslash\hat{G})\to \Map(W_E^n,L)
\]
for $n\geq 1$, linear over $\mathcal O(Q^n)\to \Map(W_E^n,L)$, such that for any map $g:\{1,\ldots,m\}\to \{1,\ldots,n\}$, the diagrams
\[\xymatrix{
\mathcal O((\hat{G}\rtimes Q)^m\sslash\hat{G})\ar[r]^-{\Theta_m}\ar[d] & \Map(W_E^m,L)\ar[d]\\
\mathcal O((\hat{G}\rtimes Q)^n\sslash\hat{G})\ar[r]^-{\Theta_n} & \Map(W_E^n,L)
}\]
induced by pullback, and
\[\xymatrix{
\mathcal O((\hat{G}\rtimes Q)^n\sslash\hat{G})\ar[r]^-{\Theta_n}\ar[d] & \Map(W_E^n,L)\ar[d]\\
\mathcal O((\hat{G}\rtimes Q)^m\sslash\hat{G})\ar[r]^-{\Theta_m} & \Map(W_E^m,L)
}\]
induced by multiplication, commute.
\end{altenumerate}
\end{proposition}

\begin{proof} We already know that (i) and (ii) are in natural bijection. The recipe above gives a canonical map from (ii) to (iii). Now take data as in (iii). Forgetting the continuity of all maps, we see that data as in (iii) gives rise to a semisimple $1$-cocycle $\varphi: W_E\to \hat{G}(L)$ (of discrete groups), up to conjugation. We need to see that if the data in (iii) are maps of condensed sets, then $\varphi$ is also a map of condensed sets (this condition does not depend on the representative of its conjugacy class). This follows from the proof of \cite[Proposition 11.7]{VLafforgue1}, in particular the choice of finitely many elements of $\gamma_1,\ldots,\gamma_n\in W_E$ such that $\varphi(\gamma)$ is determined by the closed $\hat{G}$-orbit in $(\hat{G}\rtimes Q)^{n+1}$ determined by $(\gamma_1,\ldots,\gamma_n,\gamma)$ via $\Theta_{n+1}$, cf.~\cite[Lemma 11.10]{VLafforgue1}.
\end{proof}

\section{Excursion operators}\label{sec:excursionoperators}

One can use Proposition~\ref{prop:excursionoperatorsgiveLparameters} to construct $L$-parameters in the following general categorical situation. In order to avoid topological problems, we work in the setting of the discrete subgroup $W\subset W_E/P$; in fact, we can take here any discrete group $W$. Let $\Lambda$ be a discrete $\mathbb Z_\ell$-algebra and let $\mathcal C$ be a $\mathbb Z_\ell$-linear category. Assume that functorially in finite sets $I$, we are given a monoidal $\Rep_{\mathbb Z_\ell}(Q^I)$-linear functor
\[
\Rep_{\mathbb Z_\ell}(\hat{G}\rtimes Q)^I\to \mathrm{End}(\mathcal C)^{BW^I}: V\mapsto T_V
\]
where $\mathrm{End}(\mathcal C)$ is the category of endomorphisms of $\mathcal C$, and $\mathrm{End}(\mathcal C)^{BW_E^I}$ is the category of $F\in \mathrm{End}(\mathcal C)$ equipped with a map of groups $W^I\to \mathrm{Aut}(F)$.

The goal of this section is to prove the following theorem; this is essentially due to V.~Lafforgue \cite{VLafforgue1}.

\begin{theorem}\label{thm:constructionexcursionoperators} Given the above categorical data, there is a natural map of algebras
\[
\mathrm{Exc}(W,\hat{G})=\colim_{(n,F_n\to W)} \O(Z^1(F_n,\hat{G}))^{\hat{G}}\to \mathrm{End}(\mathrm{id}_{\mathcal C})
\]
to the Bernstein center of $\mathcal C$ (i.e., the algebra of endomorphisms of the identity of $\mathcal C$).
\end{theorem}

To prove Theorem~\ref{thm:constructionexcursionoperators}, we construct explicit ``excursion operators''. These are associated to the following data.

\begin{definition}\label{def:excursiondata} An excursion datum is a tuple $\mathcal D=(I,V,\alpha,\beta,(\gamma_i)_{i\in I})$ consisting of a finite set $I$, an object $V\in \Rep_{\mathbb Z_\ell}((\hat{G}\rtimes Q)^I)$ with maps $\alpha: 1\to V|_{\Rep_{\mathbb Z_\ell}(\hat{G})}$, $\beta: V|_{\Rep_{\mathbb Z_\ell}(\hat{G})}\to 1$ and elements $\gamma_i\in W$, $i\in I$.
\end{definition}

Here, the restriction $\Rep_{\mathbb Z_\ell}((\hat{G}\rtimes Q)^I)\to \Rep_{\mathbb Z_\ell}(\hat{G})$ is the restriction to the diagonal copy of $\hat{G}\subset \hat{G}^I\subset (\hat{G}\rtimes Q)^I$.

Now consider excursion data $\mathcal D=(I,V,\alpha,\beta,(\gamma_i)_{i\in I})$. These give rise to an endomorphism of the identity functor of $\mathcal C$, as follows.
\[
S_{\mathcal D}: \mathrm{id}=T_1\xrightarrow{T_\alpha} T_V\xrightarrow{(\gamma_i)_{i\in I}} T_V\xrightarrow{T_\beta} T_1=\mathrm{id}.
\]
Varying the $\gamma_i$, this gives a map
\[
W^I\to \mathrm{End}(\mathrm{id}_{\mathcal C})
\]
to the endomorphisms of the identity functor on $\mathcal C$.

We note that if we have two excursion data $\mathcal D=(I,V,\alpha,\beta,(\gamma_i)_{i\in I})$ and $\mathcal D'=(I,V',\alpha',\beta',(\gamma_i)_{i\in I})$ with same finite set $I$ and elements $\gamma_i\in W$, and a map $g: V\to V'$ taking $\alpha$ to $\alpha'$ and $\beta'$ to $\beta$ (by post- and pre-composition), then $S_{\mathcal D}=S_{\mathcal D'}$. Indeed, the diagram
\[\xymatrix{
T_1\ar[r]^{T_\alpha}\ar@{=}[d] & T_V\ar[r]^{(\gamma_i)_{i\in I}}\ar[d]^{T_g} & T_V\ar[r]^{T_\beta}\ar[d]^{T_g} & T_1\ar@{=}[d]\\
T_1\ar[r]^{T_{\alpha'}} & T_{V'}\ar[r]^{(\gamma_i)_{i\in I}} & T_{V'}\ar[r]^{T_{\beta'}} & T_1
}\]
commutes. Now note that $(V,\alpha,\beta)$ give rise to an element
\[
f=f(V,\alpha,\beta)\in \mathcal O(\hat{G}\backslash (\hat{G}\rtimes Q)^I/\hat{G}),
\]
the quotient under diagonal left and right multiplication. Indeed, given any $g_i\in \hat{G}\rtimes Q$, $i\in I$, one can form the composite
\[
1\xrightarrow{\alpha} V\xrightarrow{(g_i)_{i\in I}} V\xrightarrow{\beta} 1,
\]
giving an element of the base ring; as $\alpha$ and $\beta$ are equivariant for the diagonal $\hat{G}$-action, this indeed gives an element
\[
f=f(V,\alpha,\beta)\in \mathcal O(\hat{G}\backslash (\hat{G}\rtimes Q)^I/\hat{G}).
\]
Conversely, given $f$ we can look at the $(\hat{G}\rtimes Q)^I$-representation $V=V_f\subset \mathcal O((\hat{G}\rtimes Q)^I/\hat{G})$ generated by $f$. This comes with a map $\alpha_f: 1\to V_f|_{\Rep_{\mathbb Z_\ell}(\hat{G})}$ induced by the element $f$, and a map $\beta_f: V_f|_{\Rep_{\mathbb Z_\ell}(\hat{G})}\to 1$ given by evaluation at $1\in (\hat{G}\rtimes Q)^I$. If we replace $V$ by the subrepresentation generated by $\alpha$, then there is a natural map $V\to V_f$ taking $\alpha$ to $\alpha_f$ and $\beta_f$ to $\beta$. The above commutative diagrams then imply that $S_\mathcal D$ depends on $(V,\alpha,\beta)$ only through $f$, and we get a map (a priori, of $\mathbb Z_\ell$-modules)
\[
\Theta^I: \mathcal O(\hat{G}\backslash (\hat{G}\rtimes Q)^I/\hat{G})\to \Map(W^I,\mathrm{End}(\mathrm{id}_{\mathcal C})).
\]
Restricted to $\mathcal O(Q^I)$, this is given by the natural map $\mathcal O(Q^I)\to \Map(W^I,\Lambda)$ (and $\Lambda\to \mathrm{End}(\mathrm{id}_{\mathcal C})$). Also, it follows from the definitions that for any map $g: I\to J$, the diagram
\[\xymatrix{
\mathcal O(\hat{G}\backslash (\hat{G}\rtimes Q)^I/\hat{G})\ar[r]^{\Theta^I}\ar[d] & \Map(W^I,\mathrm{End}(\mathrm{id}_{\mathcal C}))\ar[d]\\
\mathcal O(\hat{G}\backslash (\hat{G}\rtimes Q)^J/\hat{G})\ar[r]^{\Theta^J} & \Map(W^J,\mathrm{End}(\mathrm{id}_{\mathcal C})),
}\]
induced by pullback along $g$, is cartesian.

We want to check that $\Theta^I$ is a map of algebras. For this, we use a version of ``convolution product = fusion product'' in this situation. Namely, given $f_1,f_2\in \mathcal O(\hat{G}\backslash (\hat{G}\rtimes Q)^I/\hat{G})$, we can build their exterior product $f_1\boxtimes f_2\in \mathcal O(\hat{G}\backslash (\hat{G}\rtimes Q)^{I\sqcup I}/\hat{G})$. Then one easily checks
\[
\Theta^{I\sqcup I}(f_1\boxtimes f_2)((\gamma_i,\gamma_i')_{i\in I}) = \Theta^I(f_1)((\gamma_i)_{i\in I})\Theta^I(f_2)((\gamma_i')_{i\in I}).
\]
Applying now functoriality for pullback under $I\sqcup I\to I$, it follows that indeed $\Theta^I(f_1f_2) = \Theta^I(f_1)\Theta^I(f_2)$.

For any $n\geq 0$, we can identify
\[
\mathcal O(\hat{G}\backslash (\hat{G}\rtimes Q)^{\{0,\ldots,n\}}/\hat{G})\otimes_{\mathcal O(Q^{\{0,\ldots,n\}})} \mathcal O(Q^{\{1,\ldots,n\}})\cong \mathcal O((\hat{G}\rtimes Q)^n\sslash\hat{G})
\]
via pullback under $(g_1,\ldots,g_n)\mapsto (1,g_1,\ldots,g_n)$. This translates $\Theta^{\{0,\ldots,n\}}$ into maps of $\mathbb Z_\ell$-algebras
\[
\Theta_n: \mathcal O((\hat{G}\rtimes Q)^n\sslash\hat{G})\to \Map(W^n,\mathrm{End}(\mathrm{id}_{\mathcal C}))
\]
over $\mathcal O(Q^n)\to \Map(W^n,\Lambda)$, still satisfying compatibility with pullback under maps $g:\{1,\ldots,m\}\to \{1,\ldots,n\}$.

Arguing also as in \cite[Lemma 10.1, equation (10.5)]{VLafforgue1} and the resulting \cite[Proposition 10.8 (iii), Definition-Proposition 11.3 (d)]{VLafforgue1}, one sees that the maps $\Theta_n$ are also compatible with the multiplication maps induced by such maps $g$, thus finishing the proof of Theorem~\ref{thm:constructionexcursionoperators}.

In particular, using the description of geometric points, Theorem~\ref{thm:constructionexcursionoperators} implies the following proposition.

\begin{corollary}\label{cor:definitionLparameter} Assume that $\Lambda=L$ is an algebraically closed field and $X\in \mathcal C$ is an object with $\mathrm{End}(X)=L$. Then there is, up to $\hat{G}(L)$-conjugation, a unique semisimple $L$-parameter
\[
\varphi_X: W\to \hat{G}(L)\rtimes W
\]
such that for all excursion data $\mathcal D=(I,V,\alpha,\beta,(\gamma_i)_{i\in I})$, the endomorphism $S_{\mathcal D}(X)\in \mathrm{End}(X)=L$,
\[
X=T_1(X)\xrightarrow{\alpha} T_V(X)\xrightarrow{(\gamma_i)_{i\in I}} T_V(X)\xrightarrow{\beta} T_1(X)=X,
\]
is given by the composite
\[
L\xrightarrow{\alpha} V\xrightarrow{(\varphi_X(\gamma_i))_{i\in I}} V\xrightarrow{\beta} L.
\]
\end{corollary}

\section{Modular representation theory}\label{sec:integralcoefficients}

The goal of this section is to give a proof of Theorem~\ref{thm:coliminIndPerf}. In fact, we prove a slight refinement of it, concerning perfect complexes, that will be useful in the construction of the spectral action.

\begin{theorem}\label{thm:precisecoliminIndPerf} Assume that $\ell$ does not divide the order of $\pi_1(\hat{G})_{\mathrm{tor}}$. Then the map
\[
\colim_{(n,F_n\to W)} \mathcal O(Z^1(F_n,\hat{G}))\to \mathcal O(Z^1(W,\hat{G}))
\]
is an isomorphism in the presentable stable $\infty$-category $\Ind \Perf(\ast/\hat{G})$. Moreover, the $\infty$-category $\Perf(Z^1(W,\hat{G})/\hat{G})$ is generated under cones and retracts by $\Perf(\ast/\hat{G})$, and $\Ind\Perf(Z^1(W,\hat{G})/\hat{G})$ identifies with the $\infty$-category of modules over $\mathcal O(Z^1(W,\hat{G}))$ in $\Ind\Perf(\ast/\hat{G})$.
\end{theorem}

The difficulties in this theorem all arise on the special fibre. Indeed, we will show below that we can reduce to the following version in characteristic $\ell$.

\begin{theorem}\label{thm:precisecoliminIndPerfreducechar} Assume that $\ell$ does not divide the order of $\pi_1(\hat{G})_{\mathrm{tor}}$, and let $L=\overline{\mathbb F}_\ell$. Then the map
\[
\colim_{(n,F_n\to W)} \mathcal O(Z^1(F_n,\hat{G})_L)\to \mathcal O(Z^1(W,\hat{G})_L)
\]
is an isomorphism in the presentable stable $\infty$-category $\Ind \Perf(\ast/\hat{G})$. Moreover, the $\infty$-category $\Perf(Z^1(W,\hat{G})_L/\hat{G})$ is generated under cones and retracts by $\Perf(\ast/\hat{G})$.
\end{theorem}

Then we have the following reduction:

\begin{proof}[Theorem~\ref{thm:precisecoliminIndPerfreducechar} implies Theorem~\ref{thm:precisecoliminIndPerf}] For the colimit claim, we need to see that for all representations $V$ of $\hat{G}$, the map
\[
\colim_{(n,F_n\to W)} R\Gamma(\hat{G},\mathcal O(Z^1(F_n,\hat{G}))\otimes V)\to R\Gamma(\hat{G},\mathcal O(Z^1(W,\hat{G}))\otimes V)
\]
in $\mathcal D(\mathbb Z_\ell)$ is an isomorphism. It is an isomorphism after inverting $\ell$, as then the representation theory of $\hat{G}$ is semisimple, and it is true on underlying complexes by Proposition~\ref{prop:colimunderlying}. Thus, it suffices to show that it is an isomorphism after reduction modulo $\ell$, or even after base change to $L$, which follows from Theorem~\ref{thm:precisecoliminIndPerfreducechar}.

For the other half, note first that if $\Perf(B\hat{G})$ generates $\Perf(Z^1(W,\hat{G})/\hat{G})$, then it follows by Barr--Beck--Lurie \cite[Theorem 4.7.4.5]{LurieHA} that $\Ind\Perf(Z^1(W,\hat{G})/\hat{G})$ is the $\infty$-category of modules over $\mathcal O(Z^1(W,\hat{G}))$ in $\Ind\Perf(B\hat{G})$. Now take any $V\in \Perf(Z^1(W,\hat{G})/\hat{G})$. As its lowest cohomology group is finitely generated, we can find some surjection from an induced vector bundle onto it, and by passing to cones reduce the perfect amplitude until $V$ is a $\hat{G}$-equivariant vector bundle on $Z^1(W,\hat{G})$. We may then again find a representation $V'$ of $\hat{G}$ and a surjection $V'\otimes \mathcal O(Z^1(W,\hat{G}))\to V$. This map splits after inverting $\ell$, showing that $V$ is a retract of an induced vector bundle up to a power of $\ell$. Thus, it suffices to show that $V/\ell$ lies in this subcategory, and this follows from Theorem~\ref{thm:precisecoliminIndPerfreducechar}.
\end{proof}

Thus, we concentrate now on Theorem~\ref{thm:precisecoliminIndPerfreducechar}, which takes place over an algebraically closed base field $L$ of characteristic $\ell$. For the proof, we need many preparations on the modular representation theory of reductive groups, for $\hat{G}$ and many of its subgroups. As everything here happens on the dual side but we do not want to clutter notation, we will change notation, only for this section, and write $G$ for reductive groups over $L$.

\subsection{Good filtrations} First, we need to recall the notion of good filtrations. Let $G$ be a reductive group over $L$; recall that ``reductive'' always means connected for us. Let $T\subset B\subset G$ be a torus and Borel for $G$. For any dominant cocharacter $\lambda$ of $T$, we have the induced representation
\[
\nabla_\lambda=H^0(G/B,\mathcal O(\lambda)).
\]
A representation $V$ of $G$ has a good filtration if it admits an exhaustive filtration
\[
0=V_{-1}\subset V_0\subset V_1\subset \ldots \subset V
\]
such that each $V_i/V_{i-1}$ is isomorphic to a direct sum of $\nabla_\lambda$'s. If one picks a total ordering $0=\lambda_0,\lambda_1,\ldots$ of the dominant cocharacters, compatible with their dominance order, one can choose $V_i\subset V$ to be the maximal subrepresentation admitting only weights $\lambda_j$ with $j\leq i$. In that case, $(V_i/V_{i-1})^\ast$ is generated by its highest weight space $W_i^\ast$, and by adjunction there is a map
\[
V_i/V_{i-1}\to W_i\otimes \nabla_{\lambda_i};
\]
then $V$ admits a good filtration if and only if all of these maps are isomorphisms. For this, it is in fact enough that $V_i\to W_i\otimes \nabla_{\lambda_i}$ is surjective: The kernel is necessarily given by $V_{i-1}$, as it has only smaller weights.

A key result is that if $V$ and $W$ admit good $G$-filtrations, then so does $V\otimes W$; this is a theorem of Donkin \cite{DonkinGoodTensor} and Mathieu \cite{MathieuGoodTensor} in general. Moreover, if $V$ admits a good filtration, then $H^i(G,V)=0$ for $i>0$: This clearly reduces to the case of $V=\nabla_\lambda$, in which case it follows from Kempf's vanishing theorem \cite{KempfVanishing}. These results imply the following standard characterization of modules admitting a good filtration.

\begin{proposition}[{\cite{DonkinFiltration81}}]\label{prop:characterizegoodfiltration} A $G$-representation $V$ admits a good filtration if and only if for all $\lambda$, one has
\[
H^i(G,V\otimes \nabla_\lambda)=0
\]
for $i>0$.
\end{proposition}

Using this, one can define a well-behaved notion of a ``good filtration dimension'' of $V$, referring to the minimal $i$ such that $H^j(G,V\otimes\nabla_\lambda)=0$ for all $\lambda$ and $j>i$. Equivalently, there is a resolution of length $i$ by representations with a good filtration. Any finite-dimensional representation has finite good filtration dimension. The notion of ``good filtration dimension'' lets us put an interesting $t$-structure on the stable $\infty$-category $\Ind\Perf(\ast/G)$ freely generated under filtered colimits by perfect complexes of $G$-representations.

\begin{definition}\label{def:goodtstructure} Consider the stable $\infty$-category $\Ind\Perf(\ast/G)$ compactly generated by bounded complexes of $G$-representations. The good filtration $t$-structure is defined as follows.
\begin{enumerate}
\item[{\rm (i)}] An object $M\in \Ind\Perf(\ast/G)$ lies in the connective part of the $t$-structure if for all $\lambda$ one has $H^i(G,M\otimes \nabla_\lambda)=0$ for $i>0$; equivalently, if $M$ has good filtration dimension $\leq 0$.
\item[{\rm (ii)}] An object $M\in \Ind\Perf(\ast/G)$ lies in the coconnective part of the $t$-structure if for all $\lambda$ one has $H^i(G,M\otimes \Delta_\lambda)=0$ for $i<0$.
\end{enumerate}
\end{definition}

Equivalently, the connective part is generated under colimits by the $\nabla_\lambda$. The existence follows easily from \cite[Proposition 1.4.4.11]{LurieHA}. We warn the reader that the subcategory $\Perf(\ast/G)\subset \Ind\Perf(\ast/G)$ is not stable under the truncation operations. We will be mostly making statements that certain objects are connective in the good filtration $t$-structure.

An important observation on this $t$-structure is the following.

\begin{proposition}\label{prop:goodtstructureleftcomplete} The good filtration $t$-structure on $\Ind\Perf(\ast/G)$ is separated.
\end{proposition}

We note that the standard $t$-structure on $\Ind\Perf(\ast/G)$ is far from separated, due to issues of infinite cohomological dimension (which is the main issue we have to address).

\begin{proof} Let $M\in \Ind\Perf(\ast/G)$. To see that $M=0$ it suffices to see that for all $\lambda$, one has $R\Gamma(G,M\otimes\nabla_\lambda)=0$. But if $M$ is $\infty$-connective for the good filtration $t$-structure, then also $R\Gamma(G,M\otimes\nabla_\lambda)$ is $\infty$-connective, and hence zero.
\end{proof}

Another key result we need is the following.

\begin{theorem}[{\cite{Koppinen}, \cite{DonkinSkew}}]\label{thm:KoppinenDonkin} The $G\times G$-representation $\mathcal O(G)$ (via left and right multiplication) admits a good filtration.
\end{theorem}

In particular, we have the following corollary. For any $n\geq 0$, let $F_n$ be the free group on $n$ letters.

\begin{corollary}\label{cor:twistedconjugation} For any map $F_n\to \Aut(G)$, the $G$-representation $\mathcal O(Z^1(F_n,G))$ admits a good filtration.
\end{corollary}

\begin{proof} Note that $Z^1(F_n,G)=G^n$, where the $G$-action is that of simultaneous twisted conjugation (by the $n$ given automorphisms of $G$). But $\mathcal O(G^n)$ admits a good filtration as representation of $G^{2n}$, and restricting to $G\subset G^{2n}$ it remains good by stability under tensor products (and as the induced representations of $G^{2n}$ are tensor products of induced representations of each factor).
\end{proof}

\subsection{Equivariant vector bundles}

Assume that $X=\mathrm{Spec}(A)$ is an affine scheme of finite type over $L$, equipped with an action of a linear-algebraic group $G$ (not assumed reductive yet). We will be interested in the question whether all $G$-equivariant vector bundles on $X$ can be resolved, up to retracts, by those that are pulled back from representations of $G$ via $X/G\to \ast/G$. It is convenient to frame this question in terms of the stable $\infty$-category $\Perf(X/G)$ of $G$-equivariant perfect complexes on $X$ (i.e.~the full subcategory of dualizable objects of the quasicoherent derived $\infty$-category $D(X/G)$). We warn the reader that, being in positive characteristic, these objects are usually not compact in $D(X/G)$, even when $G$ is reductive. Let $\Perf^{\mathrm{ind}}(X/G)\subset \Perf(X/G)$ be the full subcategory generated under cones and retracts by the image of $\Perf(\ast/G)$.

\begin{proposition}\label{prop:criterionmodulegenerated} Let $M\in \Perf(X/G)$, with dual $M^\ast$ and internal endomorphisms $M\otimes_A M^\ast\in \Perf(X/G)$. Then $M\in \Perf^{\mathrm{ind}}(X/G)$ if and only if the natural map
\[
\colim [\ldots\to M\otimes_L A\otimes_L M^\ast\to M\otimes_L M^\ast]\to M\otimes_A M^\ast
\]
in $\Ind\Perf(\ast/G)$ is an isomorphism.
\end{proposition}

The left-hand side computes the tensor product $M\otimes_A M^\ast$ when all three objects are considered in $\Ind\Perf(\ast/G)$.

\begin{proof} If $M\in \Perf^{\mathrm{ind}}(X/G)$, we need to see that it is an isomorphism. In fact, it will be an isomorphism for all $N\in \Perf(X/G)$ in place of $M^\ast$. This can be reduced to $M=M_0\otimes_L A$ for some representation $M_0$ of $G$; and then replacing $N$ by $N\otimes M_0^\ast$, we can even reduce to $M_0=L$, so $M=A$. In that case, the augmented simplicial object underlying the displayed natural map has an extra degeneracy, yielding the isomorphism.

In the other direction, let $M'$ be the image of $M$ in $\mathrm{Mod}_A(\Ind\Perf(\ast/G))$, and $N'$ the image of $N=M^\ast$ in there. Then there is a natural map $M'\otimes_A N'\to A$ in $\mathrm{Mod}_A(\Ind\Perf(\ast/G))$ as the forgetful functor $\Ind\Perf(X/G)\to \mathrm{Mod}_A(\Ind\Perf(\ast/G))$ is lax symmetric monoidal (being right adjoint to the symmetric monoidal pullback). On the other hand, using in general $-'$ to denote this forgetful functor, we also get a map $A'\to M'\otimes_{A'} N'$ by our assumption, as that tensor product agrees with $(M\otimes_A N)'$. This way, we see that $M'$ is a dualizable object of the symmetric monoidal presentable stable $\infty$-category $\mathrm{Mod}_A(\Ind\Perf(\ast/G))$ with compact unit, and hence $M'$ is a compact object, and therefore a retract of a finite complex of induced vector bundles (as those are compact generators of $\mathrm{Mod}_A(\Ind\Perf(\ast/G))$ essentially by definition). As on bounded complexes, the forgetful functor $\Ind\Perf(X/G)\to \mathrm{Mod}_A(\Ind\Perf(\ast/G))$ is fully faithful, this implies the same for $M$, as desired.
\end{proof}

A curious consequence of Proposition~\ref{prop:criterionmodulegenerated} is that when $G$ is reductive with Borel $B\subset G$, then to check whether $M$ can be resolved $G$-equivariantly by induced $G$-vector bundles, it is enough to resolve $B$-equivariantly by induced $B$-vector bundles. Slightly more generally:

\begin{corollary}\label{cor:modulegeneratedBorel} Assume that $G^\circ$ is reductive with $\pi_0 G$ of order prime to $\ell$, and let $B\subset G^\circ$ be a Borel subgroup of $G^\circ$. Let $M\in \Perf(X/G)$ and assume that the corresponding object $M|_B\in \Perf(X/B)$ lies in the subcategory generated under cones and retracts by $\Perf(\ast/B)$. Then $M\in \Perf^{\mathrm{ind}}(X/G)$.
\end{corollary}

\begin{proof} By Proposition~\ref{prop:criterionmodulegenerated}, we have to see that the natural map
\[
\mathrm{colim} [\ldots\to M\otimes_L A\otimes_L M^\ast\to M\otimes_L M^\ast]\to M\otimes_A M^\ast
\]
in $\Ind\Perf(\ast/G)$ is an isomorphism. But $\Ind\Perf(\ast/G)\to \Ind\Perf(\ast/G^\circ)$ is conservative (as $\pi_0 G$ is of order prime to $\ell$ so that its representation theory is semisimple), and $\Perf(\ast/G^\circ)\to \Perf(\ast/B)$ is fully faithful by Kempf vanishing (and of course symmetric monoidal). Thus, it suffices to prove that the same map is an isomorphism in $\Ind\Perf(\ast/B)$. But this follows from Proposition~\ref{prop:criterionmodulegenerated} in the other direction.
\end{proof}

A variant of Proposition~\ref{prop:criterionmodulegenerated} is the following, which shows that the question of generating perfect complexes by induced vector bundles has direct relations to the theory of good filtrations. We will actually only use the easy direction of this proposition, and only in order to show that the assumption that $\ell$ does not divide the order of the fundamental group is necessary.

\begin{proposition}\label{prop:criterionmodulegeneratedgoodfiltration} Assume that $G^\circ$ is reductive and $\pi_0 G$ of order prime to $\ell$, and $G$ acts on $X=\mathrm{Spec}(A)$ such that $A$ admits a good $G^\circ$-filtration. Let $M\in \Perf(X/G)$ and assume that $M$, without its $A$-action, has good $G^\circ$-filtration dimension $\leq 0$, i.e.~lies in the connective part of the good $G^\circ$-filtration $t$-structure on $\Ind\Perf(\ast/G)$. Then $M\in \Perf^{\mathrm{ind}}(X/G)$ if and only if for all $N\in \Perf(X/G)$ that have good $G^\circ$-filtration dimension $\leq 0$, also $M\otimes_A N$ has good $G^\circ$-filtration dimension $\leq 0$.
\end{proposition}

\begin{proof} If $M\in \Perf^{\mathrm{ind}}(X/G)$, then the natural map
\[
\colim [\ldots\to M\otimes_L A\otimes_L N\to M\otimes_L N]\to M\otimes_A N
\]
in $\Ind\Perf(\ast/G)$ is an isomorphism, as was proved in the beginning of the proof of Proposition~\ref{prop:criterionmodulegenerated}. But if $N$ lies in the connective part of the good $G^\circ$-filtration $t$-structure, then all terms on the left-hand side lie in this connective part, and hence so does the colimit. It follows that also $M\otimes_A N$ lies in the connective part of the good filtration $t$-structure.

For the converse, we have to see that
\[
\colim [\ldots\to M\otimes_L A\otimes_L M^\ast\to M\otimes_L M^\ast]\to M\otimes_A M^\ast
\]
is an isomorphism in $\Ind\Perf(\ast/G)$. But this map is gotten by starting with the map
\[
\colim [\ldots\to A\otimes_L A\otimes_L M^\ast\to A\otimes_L M^\ast]\to M^\ast
\]
in $\Ind\Perf(X/G)$ and applying $M\otimes_A -$. Note that this colimit becomes an isomorphism in $\Ind\Perf(\ast/G)$, as then the augmented simplicial object has extra degeneracies. By \cite[Corollary 1.5]{TouzevanderKallen}, $M^\ast$ has finite good $G^\circ$-filtration dimension. This implies that the preceding colimit can be written as a sequential colimit of finite colimits, where the finite colimits become increasingly connective in the good $G^\circ$-filtration $t$-structure. Now we can apply $M\otimes_A -$ which by assumption preserves connectivity in the good $G^\circ$-filtration $t$-structure, and conclude by Proposition~\ref{prop:goodtstructureleftcomplete}.
\end{proof}

As an application of the preceding results, we have the following key result.

\begin{proposition}\label{prop:generateskyscraper} Assume that $G^\circ$ is reductive and $\pi_0 G$ is of order prime to $\ell$. Let $G$ act on itself via conjugation. Let $i: \ast\to G$ be the inclusion of the unit element. The following are equivalent.
\begin{enumerate}
\item[{\rm (i)}] The order of $\pi_1(G^\circ)_{\mathrm{tor}}$ is not divisible by $\ell$.
\item[{\rm (ii)}] The object $i_\ast L\in \Perf(G/G)$ lies in $\Perf^{\mathrm{ind}}(G/G)$.
\item[{\rm (iii)}] The inclusion $\Perf^{\mathrm{ind}}(G/G)\subset \Perf(G/G)$ is an equality.
\end{enumerate}
\end{proposition}

We are of course mainly interested in the implication from (i) to (iii), but the backwards direction tells us that (i) is really required.

\begin{proof} It is clear that (iii) implies (ii). Let us first show that (ii) implies (i). We can assume that $G$ is connected. Let $f: \tilde{G}\to G$ be a central extension such that $\tilde{G}$ has simply connected derived group, with kernel $Z\subset \tilde{G}$ of order divisible by $\ell$. Then both $f_\ast \mathcal O(\tilde{G})$ and $i_\ast L$ are in $\Perf(G/G)$ and have good filtrations. In fact, also $I=\mathrm{ker}(\mathcal O(G)\to i_\ast L)$ has a good filtration. Proposition~\ref{prop:criterionmodulegeneratedgoodfiltration} shows that if $i_\ast L\in \Perf^{\mathrm{ind}}(G/G)$, then $I\otimes_{\mathcal O(G)} \mathcal O(\tilde{G})$ must have a good filtration. This implies that the map
\[
\mathcal O(\tilde{G})^G\to \mathcal O(Z)^G = \mathcal O(Z)
\]
must be surjective. But it is known that the $\ell$-primary part of the center is contained in the unipotent locus, hence all functions in $\mathcal O(\tilde{G})^G$ are constant on them.

Now we show that (i) implies (ii). Using the criterion of Proposition~\ref{prop:criterionmodulegenerated}, we see that we can assume that $G$ is connected. We can also assume that $G$ has simply connected derived group, via a central extension (noting that finite free maps of degree prime to $\ell$ admit a canonical splitting on structure sheaves given by the trace map). We use Corollary~\ref{cor:modulegeneratedBorel}, so it suffices to show that $i_\ast L\in \Perf(G/_{\mathrm{ad}}B)$ is generated by vector bundles induced from $B$-representations. (Here, $B$ acts on $G$ via conjugation.) Now choosing a generic dominant cocharacter, so that pairing it with the roots induces a total order on the roots of $G$, we can filter $G$ by root spaces, and (using any auxiliary filtration of the torus part) find a $B$-equivariant flag of smooth subvarieties
\[
X_0=\{1\}\subset X_1\subset \ldots\subset X_{\mathrm{dim} B} = B\subset \ldots\subset X_{\mathrm{dim} G}=G.
\]
Each $X_{i-1}\subset X_i$ is a Cartier divisor, whose corresponding $B$-equivariant line bundle is induced from a character of $B$. Indeed, for $i\leq\mathrm{dim} B$ one has a map from $X_i$ to the corresponding $\mathbb G_m$ (with trivial $B$-action) or root space $\mathbb G_a$ (with $B$ acting via the root). For $i>\mathrm{dim} B$ the situation arises via pullback from a similar filtration on the flag variety $G/B$, with each term being a closed Bruhat stratum. When the derived group of $G$ is simply connected, all relevant line bundles are induced from $B$-representations.

Thus, by descending induction on $i$ we can show $\mathcal O(X_i)\in \Perf^{\mathrm{ind}}(G/_{\mathrm{ad}} B)$, yielding the desired result for $i=0$.

Finally, we prove that (ii) implies (iii). For any $K\in \Perf((G\times G)/G)$ (where $G$ acts on both factors via conjugation), we get an endofunctor of $\Ind\Perf(G/G)$ via $p_{2\ast}(p_1^\ast - \otimes K)$. Denoting
\[
\Delta: G/G\to (G\times G)/G
\]
the diagonal, the object $K=\Delta_\ast \mathcal O_{G/G}$ induces the identity endofunctor of $\Ind\Perf(G/G)$. On the other hand, if $K=p_1^\ast K_0$ for some $K_0\in \Perf(G/G)$, then the induced functor is given by
\[
p_{2\ast} p_1^\ast(-\otimes K_0) = \pi^\ast \pi_\ast(-\otimes K_0)
\]
where $\pi: G/G\to \ast/G$ is the projection. In particular, this functor has image in $\Ind\Perf^{\mathrm{ind}}(G/G)\subset \Ind\Perf(G/G)$. Thus, if $\Delta_\ast \mathcal O_{G/G}$ lies in the subcategory of $\Perf((G\times G)/G)$ generated under cones and retracts by the image of $p_1^\ast \Perf(G/G)$, then the identity endofunctor of $\Ind\Perf(G/G)$ factors over $\Ind\Perf^{\mathrm{ind}}(G/G)$, giving the desired result. But the map
\[
q: (G\times G)/G\to G/G: (g_1,g_2)\mapsto g_1g_2^{-1}
\]
sits in a cartesian diagram
\[\xymatrix{
(G\times G)/G\ar[r]^q\ar[d]^{p_1} & G/G\ar[d]^\pi\\
G/G\ar[r]^\pi & \ast/G
}\]
and $\Delta_\ast \mathcal O_{G/G}\in \Perf((G\times G)/G)$ arises via pullback from $i_\ast L\in \Perf(G/G)$. Thus, part (ii) gives the desired claim by pullback.
\end{proof}

In fact, the proof for (iii) applies more generally, for example to the following result.

\begin{proposition}\label{prop:generateperfectcomplexespowerofG} Assume that $G^\circ$ is reductive and the orders of $\pi_0 G$ and $\pi_1(G^\circ)_{\mathrm{tor}}$ are not divisible by $\ell$. Let $\Theta_1,\ldots,\Theta_n$ be automorphisms of $G$, and let $G$ act on $G^n$ via
\[
g\cdot (g_1,\ldots,g_n) = (gg_1\Theta_1(g)^{-1},\ldots,gg_n\Theta_n(g)^{-1}).
\]
Then the inclusion $\Perf^{\mathrm{ind}}(G^n/_\Theta G)\subset \Perf(G^n/_\Theta G)$ is an equality.
\end{proposition}

\begin{proof} First, if all $\Theta_i=\mathrm{id}$ and $i:\ast\to G^n$ is the inclusion of the origin, the object $i_\ast L\in \Perf(G^n/G)$ lies in $\Perf^{\mathrm{ind}}(G^n/G)$. Indeed, this follows from Proposition~\ref{prop:generateskyscraper} applied to the group $G^n$ (upon pullback from $G^n/G^n$ to $G^n/G$). Now arguing as in the proof of (iii) in Proposition~\ref{prop:generateskyscraper}, it suffices to generate $\Delta_\ast \mathcal O_{G^n/_\Theta G}$ in $\Perf((G^n\times G^n)/_\Theta G)$. But we have the map
\[
q: (G^n\times G^n)/_\Theta G\to G^n/G: (g_1,\ldots,g_n,g_1',\ldots,g_n')\mapsto (g_1 g_1^{\prime -1},\ldots,g_n g_n^{\prime -1})
\]
where on the target $G$ acts by simultaneous (untwisted) conjugation, and $\Delta_\ast \mathcal O_{G^n/_\Theta G}$ is the pullback of $i_\ast L\in \Perf(G^n/G)$. This can be resolved by induced vector bundles, and the resulting resolution of $\Delta_\ast \mathcal O_{G^n/_\Theta G}$ shows that the identity endofunctor of $\Ind\Perf(G^n/_\Theta G)$ factors over $\Ind\Perf^{\mathrm{ind}}(G^n/_\Theta G)$, giving the result.
\end{proof}

Another situation of interest will be the following. Consider a derived fibre product
\[\xymatrix{
X\ar[r]^i\ar[d]^f & \tilde{X}\ar[d]^{\tilde{f}}\\
\ast\ar[r]^{i_0} & G
}\]
where $G$ is endowed with the adjoint action of $G$, $i_0:\ast\to G$ is the inclusion of the unit, and $\tilde{X}$ is a (possibly derived) affine scheme of finite type over $L$ equipped with a compatible action of $G$. Write
\[
X=\mathrm{Spec}(A),\ \tilde{X}=\mathrm{Spec}(\tilde{A}).
\]

\begin{proposition}\label{prop:generationclosedsubscheme} Consider a derived fibre product of affine schemes of finite type with $G$-action as above. Assume that $G^\circ$ is reductive, and $\pi_0 G$ and $\pi_1(G^\circ)_{\mathrm{tor}}$ are of order prime to $\ell$.
\begin{enumerate}
\item[{\rm (i)}] The map $L\otimes_{\mathcal O(G)} \tilde{A}\to A$ is an isomorphism in $\Ind\Perf(\ast/G)$, where the tensor product on the left is formed in $\Ind\Perf(\ast/G)$.
\item[{\rm (ii)}] Assume in addition that $\Perf(\tilde{X}/G)=\Perf^{\mathrm{ind}}(\tilde{X}/G)$ and that $\tilde{A}$ lies in the connective part of the good filtration $t$-structure. Then also $A$ lies in the connective part of the good filtration $t$-structure, and $\Perf(X/G) = \Perf^{\mathrm{ind}}(X/G)$.
\end{enumerate}
\end{proposition}

\begin{proof} Part (i) follows from Proposition~\ref{prop:generateskyscraper}, as this allows us to write the derived tensor product $L\otimes_{\mathcal O(G)} -$ in a finitary way, implying in particular that it preserves bounded complexes of $G$-representations. On such, forgetting the $G$-action is conservative, giving (i).

For part (ii), the claim about $A$ is immediate from the formula in (i). Then it suffices to see that $\Perf(X/G)$ is generated by the image of $\Perf(\tilde{X}/G)$. Arguing as in Proposition~\ref{prop:criterionmodulegenerated}, we have to see that for all $M\in \Perf(X/G)$ with dual $M^\ast$, the natural map
\[
\colim [\ldots\to M\otimes_{\tilde{A}} A\otimes_{\tilde{A}} M^\ast\to M\otimes_{\tilde{A}} M^\ast]\to M\otimes_A M^\ast
\]
in $\Ind\Perf(\ast/G)$ is an isomorphism. Here, the tensor product on the right is computed in $\Perf(X/G)$, while the tensor products on the left are computed in $\Perf(\tilde{X}/G)$; more concretely, these tensor products all evaluate to bounded complexes in $\Ind\Perf(\ast/G)$, so are the ``naive'' tensor products. We note that each term here is of the form
\[
M\otimes_A g^\ast K\otimes_A M^\ast
\]
for certain $K\in \Coh(\ast\times_G \ast/G)$, where $g: X\times_{\tilde{X}} X/G\to \ast\times_G \ast/G$ is the projection. (Here, $\ast\times_G \ast$ and $X\times_{\tilde{X}} X/G$ denote the derived intersections.) Indeed, $M\otimes_A M^\ast$ corresponds to $K=\Delta_\ast L$ where $\Delta:\ast\to \ast\times_G \ast$ is the diagonal, while $M\otimes_{\tilde{A}} M^\ast$ corresponds to $L\otimes_{\mathcal O(G)} L$, i.e.~the structure sheaf of $\ast\times_G \ast$ (where the fibre product, just like all tensor products, is derived). Moreover, all transition maps come from maps between $K$'s, so the cone of the above displayed map can be written as a sequential colimit of
\[
M\otimes_A g^\ast K_n\otimes_A M^\ast
\]
for various $K_n\in \Coh(\ast\times_G \ast/G)$. Again, each term $M\otimes_A g^\ast K_n\otimes_A M^\ast$ here is a bounded complex in $\Ind\Perf(\ast/G)$.

Endow $\Ind\Coh(\ast\times_G \ast/G)$ with a $t$-structure whose connective part is generated under colimits by the image of the connective part of the good filtration $t$-structure on $\Ind\Perf(\ast/G)$ under
\[
\Delta_\ast: \Ind\Perf(\ast/G)\to \Ind\Coh(\ast\times_G \ast/G),
\]
and let $\widehat{\Ind\Coh}(\ast\times_G \ast/G)$ be the corresponding separated quotient. The functor $K\mapsto M\otimes_A g^\ast K\otimes_A M^\ast$ defines a colimit-preserving functor
\[
\Ind\Coh(\ast\times_G \ast/G)\to \Ind\Perf(\ast/G)
\]
that has bounded $t$-amplitude. Indeed, it is bounded by the good filtration dimension of $M\otimes_A M^\ast$, and the latter is finite as $M\otimes_A M^\ast$ can also be considered as an object of $\Perf(\tilde{X}/G)=\Perf^{\mathrm{ind}}(\tilde{X}/G)$, and all objects of $\Perf^{\mathrm{ind}}(\tilde{X}/G)$ have bounded good filtration dimension.\footnote{One could also cite \cite[Corollary 1.5]{TouzevanderKallen}.} As $\Ind\Perf(\ast/G)$ is separated for the good filtration $t$-structure, it induces a functor
\[
\widehat{\Ind\Coh}(\ast\times_G \ast/G)\to \Ind\Perf(\ast/G).
\]
Thus, it suffices to see that the colimit of the $K_n$ in $\widehat{\Ind\Coh}(\ast\times_G \ast/G)$ vanishes.

There is also the forgetful functor $\Coh(\ast\times_G \ast/G)\to \Perf(\ast/G)$, which also induces a ``forgetful'' functor
\[
\widehat{\Ind\Coh}(\ast\times_G \ast/G)\to \Ind\Perf(\ast/G).
\]
The colimit of the $K_n$ vanishes in $\Ind\Perf(\ast/G)$ (as then the simplicial diagram acquires extra degeneracies). Thus, it suffices to see that this forgetful functor is conservative. Equivalently, we have to see that $\widehat{\Ind\Coh}(\ast\times_G \ast/G)$ is generated by the image of $\Ind\Perf(\ast/G)$ (under pullback).

Using naive truncations, it is clear that $\Coh(\ast\times_G \ast/G)$ is generated by $\Delta_\ast V$ for varying representations $V$ of $G$. Thus, it suffices to see that $\Delta_\ast L\in \widehat{\Ind\Coh}(\ast\times_G \ast/G)$ is in the subcategory generated by the image of $\Ind\Perf(\ast/G)$. This can be proved like Proposition~\ref{prop:generateskyscraper}, reducing to the case that $G=G^\circ$ has simply connected derived group, fixing a Borel $B\subset G$ and the complete $B$-invariant flag of smooth closed subspaces $X_i\subset G$ considered there. Argue inductively about the image of the structure sheaf of $\ast\times_{X_i}\ast/B$ in $\Coh(\ast\times_G \ast/B)$ (at each step, one can get an infinite resolution by copies of the previous object twisted by powers of a line bundle), and then project to $\ast\times_G\ast/G$.
\end{proof}

These results are already sufficient to handle the case of the space of $L$-parameters of a compact Riemann surface. Indeed, this is a certain fibre product
\[\xymatrix{
X\ar[r]\ar[d] & G^{2g}\ar[d]\\
\ast\ar[r] & G
}\]
and the preceding propositions apply to show $\Perf(X/G)=\Perf^{\mathrm{ind}}(X/G)$, and identify $A$ with the excursion algebra. In fact, the same argument applies for tame $L$-parameters of local fields. It remains to deal with the wild part.

\subsection{Fixed point subgroups}\label{sec:fixedpoints}

We will need to know some properties of the fixed points $H=G^P$ of reductive groups $G$ under a (finite) group $P$ of automorphisms of $G$ of order prime to $\ell$. (Our choice of notation $P$ is motivated by the later application to the wild inertia group.) For technical reasons, we will allow $G$ to be disconnected, but always with $G^\circ$ reductive and $\pi_0 G$ of order prime to $\ell$. First, we have the following structural result.

\begin{proposition}\label{prop:fixedpointsreductive} Let $L$ be an algebraically closed field of characteristic $\ell>0$, and let $G$ be a linear algebraic group over $L$ such that $G^\circ$ is reductive and $\pi_0 G$ is of order prime to $\ell$. Assume that $P$ is a finite group of order prime to $\ell$ acting on $G$ and let $H=G^P$ be the fixed points. Then $H$ is a smooth linear algebraic group, $H^\circ$ is reductive, and $\pi_0 H$ is of order prime to $\ell$.
\end{proposition}

We note that our proof that $\pi_0 H$ is of order prime to $\ell$ probably uses unnecessarily heavy machinery. Under the assumption that $P$ is solvable (the only case relevant to us), this can be deduced much more directly from Steinberg's theorem \cite[Theorem 8.1]{SteinbergEndomorphisms} by reducing to cyclic $P$ and simply connected $G$.

\begin{proof} We can assume $G=G^\circ$. It is a standard fact that the fixed points of a smooth affine scheme under a finite group of order prime to the characteristic is still affine and smooth. Moreover, by \cite[Theorem 2.1]{PrasadYu}, $H^\circ$ is reductive.

For the final statement, we consider the action of $G$ on $\prod_{P\setminus\{1\}} G$, where it acts on the factor enumerated by $\Theta\in P$ through $\Theta$-twisted conjugation. Note that $\mathcal O(\prod_{P\setminus\{1\}} G)$ has a good $G$-filtration. By \cite[Corollary 1.5]{TouzevanderKallen}, for any $G$-equivariant finitely generated $\mathcal O(\prod_{P\setminus\{1\}} G)$-module $M$, the good filtration dimension of $M$ is finite, and in particular $H^i(G,M)=0$ for all large enough $i$.

Now $G/H$ is a closed orbit of $G$ acting on $\prod_{P\setminus\{1\}} G$ (the orbit of the identity element). Moreover, if $\pi_0 H$ has an element of order $\ell$, then we get a subgroup $H'\subset H$ with $\pi_0 H'\cong \mathbb Z/\ell \mathbb Z$. In that case,
\[
R\Gamma(G,\mathcal O(G/H'))\cong R\Gamma(H',L)\cong R\Gamma(\mathbb Z/\ell\mathbb Z,L)
\]
has cohomology in all positive degrees, while $\mathcal O(G/H')$ correponds to a $G$-equivariant coherent sheaf on $\prod_{P\setminus\{1\}} G$ (equipped with its $\Theta$-conjugation), so this contradicts the previous paragraph.
\end{proof}

Assume from now on that $P$ is solvable. Our goal now is to prove the following result. This simultaneously generalizes the classical case of Levi subgroups, and the case of involutions known as Brundan's conjecture \cite{Brundan}, \cite{vanderKallenBrundan}. The remaining cases are for exceptional groups, and are discussed under some restrictions on $\ell$ in \cite{HagueMcNinch}.

\begin{theorem}\label{thm:fixedpointsdonkin} The subgroup $H^\circ\subset G^\circ$ is a Donkin subgroup. In other words, for any representation $V$ of $G^\circ$ that admits a good $G^\circ$-filtration, also $V|_{H^\circ}$ admits a good $H^\circ$-filtration. Equivalently, for any representation $W$ of $H^\circ$ that admits a good $H^\circ$-filtration, also $\mathrm{Ind}_{H^\circ}^{G^\circ} W$ admits a good $G^\circ$-filtration.
\end{theorem}

A notable consequence is that the well-known assertion that Levi subgroups are Donkin subgroups can be generalized to the statement that centralizers of regular semisimple elements define Donkin subgroups. Thus the proof below gives, in particular, a new proof that Levi subgroups are Donkin subgroups.

\begin{proof} First, Proposition~\ref{prop:characterizegoodfiltration} shows that the formulations are equivalent. Indeed, Proposition~\ref{prop:characterizegoodfiltration} applied to the group $H^\circ$ shows that the first statement is equivalent to the assertion that for all representations $V$ of $G^\circ$ with a good $G^\circ$-filtration and all representations $W$ of $H^\circ$ with a good $H^\circ$-filtration, one has $H^i(H^\circ,V|_{H^\circ}\otimes W)=0$. But $H^i(H^\circ,V|_{H^\circ}\otimes W)=H^i(G^\circ,V\otimes \mathrm{Ind}_{H^\circ}^{G^\circ} W)$, so Proposition~\ref{prop:characterizegoodfiltration} applied to the group $G^\circ$ translates this into the second statement.

By induction, we can assume that $P$ is a cyclic group of prime order $p\neq \ell$, so $P=\mathbb Z/p\mathbb Z$. Let $\Theta$ denote the automorphism of $G$ corresponding to $1\in \mathbb Z/p\mathbb Z=P$. We can do also evidently assume that $G$ is connected, and reduce to the case that $G$ is simple and simply connected: As the property of admitting a good filtration is detected after restriction to the derived group, and also after passing to finite covers, we can assume that $G$ is simply connected. In fact, decomposing $G$ into simple factors, we can assume that $P$ permutes the simple factors of $G$ transitively. If $G$ is not simple, then $G = \prod_P H$ with $P$-action permuting the factors, and $H\subset G=\prod_P H$ is the diagonal embedding, so the result follows from the stability of good filtrations under tensor products. Thus, we can assume that $G$ is simple (and simply connected). In particular, $H$ is connected by \cite[Theorem 8.1]{SteinbergEndomorphisms}.

We will first handle the case $W=1$; or more precisely, the assertion that $\mathcal O(G/H) = \mathrm{Ind}_H^G 1$ has a good $G$-filtration. This argument works directly for arbitrary $G$ as in the statement of the theorem, so does not need the previous reduction. Let
\[
X=\{(g_0,\ldots,g_{p-1})\in G^p\mid g_0\Theta(g_1)\cdots \Theta^{p-1}(g_{p-1})=1\}\subset G^p
\]
endowed with the simultaneous $\Theta$-conjugation by $G$; i.e.~$g\in G$ acts on $(g_i)_i=(g_0,\ldots,g_{p-1})\in X$ via
\[
g(g_i)_i=(gg_i\Theta(g)^{-1})_i.
\]
Moreover, endow $X$ with the $G$-equivariant $P=\mathbb Z/p\mathbb Z$-action taking $(g_0,g_1\ldots,g_{p-1})$ to $(g_1,\ldots,g_{p_1},g_0)$. Let $Y=X^P$. As $X$ is smooth affine and $P$ is of order prime to $\ell$, also $Y$ is a smooth affine scheme, equipped with a remaining $G$-action. Concretely,
\[
Y=\{g_0\in G\mid g_0\Theta(g_0)\cdots \Theta^{p-1}(g_0)=1\}\subset G
\]
equipped with $\Theta$-conjugation. A simple calculation on tangent spaces shows that this is a finite disjoint union of $G$-orbits. In particular, the orbit of $1\in Y\subset G$ is given by $G/H$ where $H=G^P=G^\Theta$.

By Lemma~\ref{lem:fixedpointsIndPerf} below, the map
\[
\mathrm{colim}_{[n]\in \Delta^{\mathrm{op}}} \mathcal O(\prod_{P^n} X)\to \mathcal O(X^P) = \mathcal O(Y)
\]
is an isomorphism in $\mathrm{IndPerf}(\ast/G)$. But all terms on the left-hand side admit a good $G$-filtration, hence also $\mathcal O(Y)$ lies in the connective part of the good filtration $t$-structure, and so admits a good $G$-filtration. In particular, the retract $\mathcal O(G/H)$ has a good $G$-filtration.

At this point, we make use of the reduction to the case that $G$ is simple and simply connected. Assume first that $\Theta$ is an inner automorphism. Thus, $H$ is the centralizer of some regular semisimple element $g\in G$ (with $\overline{g}\in G_{\mathrm{ad}}$ of order $p$). We can find some maximal torus $T\subset G$ containing $g$, which is then also a maximal torus of $H$. Choose an enumeration of the dominant weights $0=\lambda_0,\lambda_1,\ldots$ of $G$ such that all weights of $\nabla_{\lambda_i}$ are contained in the Weyl group orbit of $\lambda_0,\ldots,\lambda_i$. We argue by induction on $i$ that $\nabla_{\lambda_i}|_H$ has a good $H$-filtration. More precisely, let $W_i\subset W$ be the subset of elements $w\in W$ such that $w\lambda_i$ is $H$-dominant. Then we claim, by induction on $i$, that there is a surjective map
\[
\nabla_{\lambda_i}|_H\to \bigoplus_{w\in W_i} \nabla^H_{w\lambda_i}
\]
whose kernel has a good $H$-filtration (where all weights that occur in the kernel are in the Weyl group orbit of $\lambda_0,\ldots,\lambda_{i-1}$).

Note that there is indeed such a natural map, as $\nabla_{\lambda_i}|_H$ has a $w\lambda_i$-weight space of dimension $1$ for all $w\in W_i$. We need to see that the homotopy fibre $X$ of this map is connective in the good $H$-filtration $t$-structure. This homotopy fibre has only weights in the Weyl group orbit of $\lambda_0,\ldots,\lambda_{i-1}$. By the characterization of the connective part of the good filtration $t$-structure it suffices to see that
\[
R\Gamma(H,X\otimes \nabla^H_{w'\lambda_j})
\]
is connective for all $j<i$ and $w'\in W_j$. But $R\Gamma(H,\nabla^H_{w\lambda_i}\otimes \nabla^H_{w'\lambda_j})=0$, so using the definition of $X$ this can be rewritten more easily as
\[
H^m(H,\nabla_{\lambda_i}|_H\otimes \nabla^H_{w'\lambda_j})=0
\]
for $m>0$. Equivalently, for all $j<i$,
\[
H^m(H,\nabla_{\lambda_i}|_H\otimes \bigoplus_{w'\in W_j} \nabla^H_{w'\lambda_j})=0
\]
for $m>0$. But by induction on $i$ we have the surjective map
\[
\nabla_{\lambda_j}|_H\to \bigoplus_{w'\in W_j} \nabla^H_{w'\lambda_j}
\]
whose kernel has a good $H$-filtration (with smaller weights). Thus, doing also an induction on $j$, we see that it suffices to see that
\[
H^m(H,\nabla_{\lambda_i}|_H\otimes \nabla_{\lambda_j}|_H)=0
\]
for $m>0$. But
\[
H^m(H,\nabla_{\lambda_i}|_H\otimes \nabla_{\lambda_j}|_H) = H^m(G,\mathrm{Ind}_H^G 1\otimes \nabla_{\lambda_i}\otimes \nabla_{\lambda_j})=0
\]
as all three tensor factors admit a good $G$-filtration, hence so does the tensor product.

It remains to handle the case that $G$ is simple and simply connected and $\Theta$ is an outer automorphism, so necessarily $p=2$ or $p=3$. We could finish this off by appealing for $p=2$ to the case of involutions handled by Brundan \cite{Brundan} and van der Kallen \cite{vanderKallenBrundan} and for $p=3$ by noting that this only occurs when $G=\mathrm{Spin}_8$, where $H$ is either $G_2\subset \mathrm{Spin}_8$ (for the diagram automorphism) which is handled in \cite{Brundan}, or the fixed points $\mathrm{SL}_3\subset G_2\subset \mathrm{Spin}_8$ of an inner automorphism of order $3$ of $G_2$ (handled either by the above, or by another reference to \cite{Brundan}). On the other hand, we can also repeat the arguments above. Namely, we note that we can in general lift any highest weight (i.e.~dominant cocharacter) $\mu$ of $H$ to a highest weight $\lambda$ of $G$ and then a similar inductive argument shows that there is some finite set $X_\mu$ of dominant cocharacters of $H$ (with $\mu\in X_\mu$) such that there is a surjective map
\[
\nabla_\lambda|_H\to \bigoplus_{\mu'\in X_\mu} \nabla^H_{\mu'}
\]
for which the kernel has ``smaller'' weights and a good $H$-filtration. Now if $\lambda$ is general, we want to see that $\nabla_\lambda|_H$ has a good $H$-filtration. We have to see that for all dominant cocharacters $\mu$ of $H$, the $H$-representation $\nabla_\lambda|_H\otimes \nabla^H_\mu$ has no higher cohomology. Pick some dominant cocharacter $\lambda_\mu$ of $G$ lifting $\lambda$. By induction on $\mu$ and the preceding claim, it suffices to see that $\nabla_\lambda|_H\otimes \nabla_{\lambda_\mu}|_H$ has no higher $H$-cohomology. But this agrees with the $G$-cohomology of $\mathrm{Ind}_H^G 1\otimes \nabla_\lambda\otimes \nabla_{\lambda_\mu}$ which has no higher $G$-cohomology as all three tensor factors have good $G$-filtrations.
\end{proof}

We used the following key lemma in the proof. Its full force will be required later to prove Theorem~\ref{thm:precisecoliminIndPerfreducechar}.

\begin{lemma}\label{lem:fixedpointsIndPerf} Let $G$ be a linear algebraic group over an algebraically closed field $L$ of characteristic $\ell$ such that $G^\circ$ is reductive and $\pi_0 G$ is of order prime to $\ell$. Let $\Theta$ be an automorphism of $G$ of prime order $p\neq \ell$. Let
\[
X=\{(g_0,\ldots,g_{p-1})\in G^p\mid g_0\Theta(g_1)\cdots \Theta^{p-1}(g_{p-1})=1\}
\]
equipped with the $G$-action of simultaneous $\Theta$-conjugation, and the $P=\mathbb Z/p\mathbb Z$-action of cyclic permutation. Consider the corresponding augmented cosimplicial $G$-space
\[
X^P\to X\rightrightarrows \prod_P X \ldots.
\]
Then the map
\[
\mathrm{colim}_{[n]\in \Delta^{\mathrm{op}}} \mathcal O(\prod_{P^n} X)\to \mathcal O(X^P)
\]
is an isomorphism in $\mathrm{IndPerf}(\ast/G)$.
\end{lemma}

\begin{proof} The idea is to introduce the formal completion $\hat{X}$ of $X$ along $X^P$ and then observe that on the one hand, replacing $X$ by $\hat{X}$ yields an isomorphic colimit; and on the other hand, that there is a canonical $(G,P)$-equivariant retraction from $\hat{X}$ onto $X^P$. More precisely, after replacing all terms $\prod_{P^n} X$ by formal completions along $X^P$, the augmented cosimplicial $G$-space acquires an extra degeneracy. 

To execute this strategy, we have to enlarge $\mathrm{IndPerf}(\ast/G)$ in order to allow power series algebras. The following discussion is inspired by the theory of solid modules \cite{Condensed}. Note that $\mathrm{IndPerf}(\ast/G)$ is freely generated, as a presentable stable $\infty$-category, by the exact category $\mathrm{Rep}(G)$. (So $\mathrm{IndPerf}(\ast/G)$ is the $\infty$-category of contravariant functors from $\mathrm{Rep}(G)^{\mathrm{op}}$ to spectra that take exact sequences to fibre sequences.) Let $\mathrm{Pro}(\mathrm{Rep}(G))$ be the Pro-category of $\mathrm{Rep}(G)$; it is again an exact symmetric monoidal category (using that Pro-vector spaces are well-behaved). We will actually only need countable Pro-systems, so the reader is invited to restrict to this subcategory. Let us denote by $\mathrm{IndProPerf}(\ast/G)$ the corresponding category freely generated as a presentable stable $\infty$-category by the exact category $\mathrm{Pro}(\mathrm{Rep}(G))$. (It is a slight misnomer as the compact objects are not all of $\mathrm{Pro}(\mathrm{Perf}(\ast/G))$ but only the part of bounded amplitude.)

Let $\hat{\mathcal O(G)}$ be the completion of $\mathcal O(G)$ at the unit element (equipped with the action of $G$ by usual conjugation); this is an object of $\mathrm{Pro}(\mathrm{Rep}(G))$. Formally, let $\hat{G}\subset G$ be the corresponding geometric object. Critically, $\hat{\mathcal O(G)}$ is an idempotent $\mathcal O(G)$-algebra in $\mathrm{IndProPerf}(\ast/G)$. Indeed, by the usual trick using the group structure to write the diagonal as a pullback of the unit section, it suffices to see that $L\otimes_{\mathcal O(G)} \hat{\mathcal O(G)} = L$. If $\ell$ does not divide the order of $\pi_1(G^\circ)_{\mathrm{tor}}$, then this follows from Proposition~\ref{prop:generateskyscraper}, as then the left-hand side is a bounded complex and the isomorphism can be checked after forgetting the $G$-action (where it is clear). In general, we can find an embedding $G\hookrightarrow \GL_n$, and it suffices to see that
\[
\hat{\mathcal O(\GL_n)} \otimes_{\mathcal O(\GL_n)} \mathcal O(G) = \hat{\mathcal O(G)}
\]
in $\mathrm{IndProPerf}(\ast/G)$. For this, it suffices to see that the left-hand side is a bounded complex, for which it suffices to show that $\mathcal O(G)\in \Perf(\GL_n/G)$ lies in $\Perf^{\mathrm{ind}}(\GL_n/G)$. But in fact $\Perf(\GL_n/G) = \Perf^{\mathrm{ind}}(\GL_n/G)$ by the argument of Proposition~\ref{prop:generateskyscraper}~(iii).

Note that $X$ embeds into $X' = G^p$ with $X'^P = G$. The corresponding augmented cosimplicial $G$-space
\[
X'^P\to X'\rightrightarrows \prod_P X' \ldots
\]
admits an extra degeneracy and hence the map
\[
\mathrm{colim}_{[n]\in \Delta^{\mathrm{op}}} \mathcal O(\prod_{P^n} X')\to \mathcal O(X'^P)=\mathcal O(G)
\]
is an isomorphism in $\mathrm{IndPerf}(\ast/G)$. The same applies to $\widehat{\prod_{P^n} X'}$, the completion of $\prod_{P^n} X'$ along the diagonal inclusion $X'^P=G$, so also the map
\[
\mathrm{colim}_{[n]\in \Delta^{\mathrm{op}}} \mathcal O(\widehat{\prod_{P^n} X'})\to \mathcal O(X'^P)=\mathcal O(G)
\]
is an isomorphism in $\mathrm{IndProPerf}(\ast/G)$. In particular, the natural map
\[
\mathrm{colim}_{[n]\in \Delta^{\mathrm{op}}} \mathcal O(\prod_{P^n} X')\to \mathrm{colim}_{[n]\in \Delta^{\mathrm{op}}} \mathcal O(\widehat{\prod_{P^n} X'})
\]
is an isomorphism in $\mathrm{IndProPerf}(\ast/G)$. As geometric realizations of bisimplicial objects can be computed after diagonal restriction, it now suffices to show
\[
\mathrm{colim}_{[n]\in \Delta^{\mathrm{op}}} \mathcal O(\prod_{P^n} X)\otimes_{\mathcal O(\prod_{P^n} X')} \mathcal O(\widehat{\prod_{P^n} X'})\to \mathcal O(X^P)
\]
is an isomorphism in $\mathrm{IndProPerf}(\ast/G)$. We note that each term identifies with $\mathcal O(\widehat{\prod_{P^n} X})$ where $\widehat{\prod_{P^n} X}$ is the completion of $\prod_{P^n} X$ along the diagonal embedding of $X^P$. (All completions here are a priori derived but actually concentrated in degree $0$ as all varieties are noetherian (even smooth).) Indeed, this is clear without the $G$-action. To see it $G$-equivariantly (i.e., in $\mathrm{IndProPerf}(\ast/G)$), it suffices to see that these tensor product are computed by bounded complexes. If $\ell$ does not divide the order of $\pi_1(G^\circ)_{\mathrm{tor}}$, this follows from Proposition~\ref{prop:generateperfectcomplexespowerofG} which ensures
\[
\mathcal O(\prod_{P^n} X)\in \Perf(\prod_{P^n} X'/G) = \Perf^{\mathrm{ind}}(\prod_{P^n} X'/G).
\]
In general, we can argue as above by using an embedding of $G$ into $\GL_n$.

Thus, we have a cosimplicial augmented $G$-space
\[
X^P\to \hat{X}\rightrightarrows \widehat{\prod_P X}\rightrightarrows \ldots
\]
with $\widehat{\prod_{P^n} X}$ the completion of $\prod_{P^n} X$ along the diagonal embedding of $X^P$. It remains to see that this has a $G$-equivariant extra degeneracy. For this, it suffices to construct a $(G,P)$-equivariant map
\[
\widehat{\prod_P X}\to \hat{X}
\]
whose restriction along the diagonal embedding $\hat{X}\to \widehat{\prod_P X}$ is the identity of $\hat{X}$. Indeed, this defines an ``averaging'' map that can be used to construct the extra degeneracies in a standard way.

Such a map $\widehat{\prod_P X}\to \hat{X}$ needs to take a $p$-tuple $(g_{0,i},\ldots,g_{p-1,i})_{i\in P}$ of $p$-tuples $(g_{0,i},\ldots,g_{p-1,i})$ satisfying
\[
g_{0,i}\Theta(g_{1,i})\cdots \Theta^{p-1}(g_{p-1,i})=1
\]
and produce a new $p$-tuple $(h_0,\ldots,h_{p-1})$ such that
\[
h_0\Theta(h_1)\cdots \Theta^{p-1}(h_{p-1})=1.
\]
The $G$-equivariance means that this construction must be invariant under simultaneous $\Theta$-conjugation. The $P$-equivariance means that the recipe for $h_0$ must determine the recipes for $h_1,\ldots,h_{p-1}$ through suitable conjugation. And the final condition is that if one applies this in the case where all $g_{a,i}=g_a$ are independent of $i$, then $h_a=g_a$ for $a=0,\ldots,p-1$.

Now we define this retraction by
\[
h_0 = (g_{0,0}\Theta(g_{1,1})\cdots \Theta^{p-1}(g_{p-1,p-1}))^{-1/p} g_{0,0}.
\]
Here, the element of which we are taking a $p$-th root is an element of $\hat{G}$ (the completion of $G$ at the unit) and here one can uniquely take $p$-th roots. Moreover, the formation of $p$-th roots is invariant under $G$-conjugation; thus this formula has the required $G$-equivariance. If all $g_{a,i}=g_a$ are independent of $i$, then this element is actually equal to $1$ and hence $h_0 = g_0$.

We are then forced to take
\[
h_i = (g_{i,i}\Theta(g_{i+1,i+1})\cdots \Theta^{p-1}(g_{i-1,i-1}))^{-1/p} g_{i,i}
\]
for the other $i$ in order to get $P$-equivariance. Thus, we get $(G,P)$-equivariance, and the property that on equal elements $g_{a,i}=g_a$, one has $h_a=g_a$. It remains to see that we actually defined a map to $\hat{X}$, for which we have to see
\[
h_0\Theta(h_1)\cdots \Theta^{p-1}(h_{p-1})=1.
\]
But, using that formation of $p$-th roots is invariant under conjugation, one can start rewriting the first two factors
\[\begin{aligned}
h_0\Theta(h_1) &= (g_{0,0}\Theta(g_{1,1})\cdots \Theta^{p-1}(g_{p-1,p-1}))^{-1/p} g_{0,0} (\Theta(g_{1,1})\Theta^2(g_{2,2})\cdots g_{0,0})^{-1/p} \Theta(g_{1,1}) \\
&= g_{0,0} (\Theta(g_{1,1})\cdots \Theta^{p-1}(g_{p-1,p-1})g_{0,0})^{-1/p}(\Theta(g_{1,1})\Theta^2(g_{2,2})\cdots g_{0,0})^{-1/p} \Theta(g_{1,1})\\
&= g_{0,0} (\Theta(g_{1,1})\Theta^2(g_{2,2})\cdots g_{0,0})^{-2/p} \Theta(g_{1,1})\\
&= g_{0,0} \Theta(g_{1,1}) (\Theta^2(g_{2,2})\cdots g_{0,0}\Theta(g_{1,1}))^{-2/p}.
\end{aligned}\]
This procedure can be continued, shifting the $p$-th root past each factor $\Theta^i(g_{i,i})$ using conjugation invariance of forming $p$-th roots; and then it matches the next $p$-th root, so they can be combined. In the end, one gets
\[
h_0\Theta(h_1)\cdots \Theta^{p-1}(h_{p-1}) = g_{0,0} \Theta(g_{1,1})\cdots \Theta^{p-1}(g_{p-1,p-1}) (g_{0,0}\Theta(g_{1,1})\cdots \Theta^{p-1}(g_{p-1,p-1}))^{-p/p} = 1.\qedhere
\]
\end{proof}

Recall that Theorem~\ref{thm:fixedpointsdonkin} affirms that for a representation $W$ of $H^\circ$, if $W$ admits a good $H^\circ$-filtration then $\mathrm{Ind}_{H^\circ}^{G^\circ} W$ admits a good $G^\circ$-filtration. In fact, the converse is true as well. More precisely, we have the following assertions.

\begin{proposition}\label{prop:strongdonkinsubgroup} In the situation of Theorem~\ref{thm:fixedpointsdonkin}, the following results hold true.
\begin{enumerate}
\item[{\rm (i)}] Let $W$ be a representation of $H^\circ$. Then $W$ admits a good $H^\circ$-filtration if and only if $\mathrm{Ind}_{H^\circ}^{G^\circ} W$ admits a good $G^\circ$-filtration.
\item[{\rm (ii)}] Let $W$ be a representation of $H$. Then $W$ admits a good $H^\circ$-filtration if and only if $\mathrm{Ind}_H^G W$ admits a good $G^\circ$-filtration.
\item[{\rm (iii)}] Let $W$ be a representation of $H^\circ$ that admits a good $H^\circ$-filtration. Then the kernel of
\[
(\mathrm{Ind}_{H^\circ}^{G^\circ} W)|_{H^\circ}\to W
\]
admits a good $H^\circ$-filtration.
\item[{\rm (iv)}] Let $W$ be a representation of $H$ that admits a good $H^\circ$-filtration. Then the kernel of
\[
(\mathrm{Ind}_H^G W)|_H\to W
\]
admits a good $H^\circ$-filtration.
\item[{\rm (v)}] The image of the restriction $\Perf(\ast/G^\circ)\to \Perf(\ast/H^\circ)$ generates under cones and retracts.
\item[{\rm (vi)}] The image of the restriction $\Perf(\ast/G)\to \Perf(\ast/H)$ generates under cones and retracts.
\end{enumerate}
\end{proposition}

In a previous version of this manuscript, assertion (vi) was proved by an exhaustive analysis of all possible cases. One key issue in the proof is that the individual assertions do not generally allow simple reductions to the case that $G$ is simply connected and the precise form of the center matters, as the following example shows. We are able to give a better argument now as Theorem~\ref{thm:fixedpointsdonkin} gives us very good control about representations with good filtrations, making a d\'evissage to the simply connected case possible.

\begin{remark} The following example shows that the hypothesis that $P$ is of order prime to $\ell$ is important, and cannot be weakened to ``quasi-semisimple'' automorphisms (preserving a Borel and a torus) or even automorphisms fixing a pinning; also, the example shows that the precise form of the center is critical. If $G=(\mathrm{SL}_2\times\mathrm{SL}_2)/\mu_2$ with the automorphism switching the two factors, then $H=\mathrm{PGL}_2\times (\mu_2\times \mu_2)/\mu_2$. If we had $\ell=2$, then one can show that for all objects $A\in \Perf(\ast/H)$ in the image of $\Perf(\ast/G)$, the summand $A_1$ of $A$ with nontrivial central character has the property that the (homotopy) invariants of the $\mathbb Z/2\mathbb Z\subset \mathrm{PGL}_2$-action on $A_1$ are a perfect complex.\footnote{Indeed, in characteristic $2$ the irreducible representations of $\mathrm{SL}_2$ are tensor products of Frobenius pullbacks of the standard representation. To get something on which the center of $H$ acts nontrivially, the standard representation has to occur as a tensor factor on both $\mathrm{SL}_2$-factors of $G$. But the standard representation of $\mathrm{SL}_2$ restricts to the regular representation of $\mathbb Z/2\mathbb Z\subset \mathrm{PGL}_2$, so tensored with anything other finite-dimensional representation, the homotopy invariants are a perfect complex.} This implies that the nontrivial character of $H$ is not generated by $\Perf(\ast/G)$ under cones and retracts.
\end{remark}

\begin{proof} We will first show that assertions (i), (iii) and (v) are equivalent; and similarly, that assertions (ii), (iv) and (vi) are equivalent.

Let us start by showing that (v) implies (iii) and (i). Note that (v) is equivalent to the assertion that for all representations $W$ of $H^\circ$, the complex
\[
\ldots\to \mathrm{Ind}_{H^\circ}^{G^\circ} \mathrm{Ind}_{H^\circ}^{G^\circ} W\to \mathrm{Ind}_{H^\circ}^{G^\circ} W\to W\to 0
\]
is a resolution in $\Ind\Perf(\ast/H^\circ)$. Indeed, if it is a resolution, then in particular $W$ is resolved by objects in the image of $\Ind\Perf(\ast/G^\circ)\to \Ind\Perf(\ast/H^\circ)$, so the image generates. Conversely, if $\Perf(\ast/G^\circ)\to \Perf(\ast/H^\circ)$ generates under cones and retracts, then $\mathrm{Ind}_{H^\circ}^{G^\circ}: \Ind\Perf(\ast/H^\circ)\to \Ind\Perf(\ast/G^\circ)$ is conservative; but the above complex becomes split exact after this operation.

Now if $W$ is as in (i), such that $\mathrm{Ind}_{H^\circ}^{G^\circ} W$ admits a good $G^\circ$-filtration, then this resolution (and Theorem~\ref{thm:fixedpointsdonkin}) shows that $W$ admits a resolution in $\Ind\Perf(\ast/H^\circ)$ where all terms admits a good $H^\circ$-filtration, i.e.~are connective in the good $H^\circ$-filtration $t$-structure. It follows that $W$ is connective in the good $H^\circ$-filtration $t$-structure, i.e.~admits a good $H^\circ$-filtration. Similarly, if $W$ as in (iii), then the same argument applied to the first truncation of the last displayed resolution shows that this kernel admits a good $H^\circ$-filtration.

For the converse, we show that (i) implies (iii), and (iii) implies (v). Assume (i), and take any $W$ as in (iii). Then to show that this kernel admits a good $H^\circ$-filtration, it suffices to see that after applying $\mathrm{Ind}_{H^\circ}^{G^\circ}$, it admits a good $G^\circ$-filtration. But after applying $\mathrm{Ind}_{H^\circ}^{G^\circ}$, the map becomes split, and $\mathrm{Ind}_{H^\circ}^{G^\circ} \mathrm{Ind}_{H^\circ}^{G^\circ} W$ admits a good $G^\circ$-filtration by Theorem~\ref{thm:fixedpointsdonkin}. Now we show that (iii) implies (v). To show (v), it suffices to see that for any $H^\circ$-representation $W$ with a good $H^\circ$-filtration, there is some $G^\circ$-representation $V$ with a good $G^\circ$-filtration and a surjection $V|_{H^\circ}\to W$ whose kernel admits a good $H^\circ$-filtration. Indeed, one can then inductively build a resolution of $W$, using that $\Ind\Perf(\ast/H^\circ)$ is separated for the good $H^\circ$-filtration $t$-structure. But by (iii), we can simply take $V=\mathrm{Ind}_{H^\circ}^{G^\circ} W$.

A similar analysis shows that (ii), (iv) and (vi) are equivalent. Note first that in (ii), $W$ is a direct summand of $\mathrm{Ind}_{H^\circ}^H W$ (as $\pi_0 H$ is of order prime to $\ell$ by Proposition~\ref{prop:fixedpointsreductive}), and hence $\mathrm{Ind}_H^G W$ is a direct summand of $\mathrm{Ind}_{H^\circ}^G W$; using Theorem~\ref{thm:fixedpointsdonkin} this shows that if $W$ admits a good $H^\circ$-filtration then also $\mathrm{Ind}_H^G W$ admits a good $G^\circ$-filtration. So again (ii) is really about the converse direction.

Now again (vi) is equivalent to the complex
\[
\ldots\to \mathrm{Ind}_H^G \mathrm{Ind}_H^G W\to \mathrm{Ind}_H^G W\to W\to 0
\]
being a resolution in $\Ind\Perf(\ast/H)$, for all representations $W$ of $H$. Assume that this holds. Then if $W$ is as in (ii) so that $\mathrm{Ind}_H^G W$ admits a good $G^\circ$-filtration, then all terms in this resolution have a good $H^\circ$-filtration (by Theorem~\ref{thm:fixedpointsdonkin}), hence so does $W$. Similarly, if $W$ is as in (iv), then the truncation of this sequence shows that $\mathrm{ker}(\mathrm{Ind}_H^G W\to W)$ has a good $H^\circ$-filtration. Thus, (vi) implies (ii) and (iv).

Going from (ii) to (iv) is the same argument as going from (i) to (iii). The argument that (iv) implies (vi) also adapts verbatim from the argument that (iii) implies (v), noting that $\Ind\Perf(\ast/H)$ is separated for the good $H^\circ$-filtration $t$-structure (as $\pi_0 H$ is of order prime to $\ell$).

Thus, all items with an odd number are equivalent, as are all items with an even number. But it is clear that (vi) implies (v) as the image of $\Perf(\ast/H)\to \Perf(\ast/H^\circ)$ generates under retracts, as for all representations $W$ of $H^\circ$, one can split $W$ off $\mathrm{Ind}_{H^\circ}^H W$ as $\pi_0 H$ is of order prime to $\ell$. On the other hand, if $G$ is connected, we claim that (iii) implies (iv). Take any $H$-representation $W$ that admits a good $H^\circ$-filtration. Then
\[
\mathrm{Ind}_{H^\circ}^G(W|_{H^\circ}) = \mathrm{Ind}_H^G(\mathrm{Ind}_{H^\circ}^H(W|_{H^\circ}))
\]
admits $\mathrm{Ind}_H^G W$ as a direct summand, in a way compatibly with the map back to $W$. This implies that the kernel in (iv) is a direct summand of the kernel in (iii) for the representation $W|_{H^\circ}$, and hence (iii) implies (iv).

At this point, we have proved that for connected $G$, all six assertions are equivalent, and in general (vi) implies all other assertions. Moreover, assertion (vi) is amenable to induction on $P$, so we can assume that $P=\mathbb Z/p\mathbb Z$ is cyclic of prime order $p\neq \ell$. If $G$ is connected and simply connected, then the arguments with highest weights in the proof of Theorem~\ref{thm:fixedpointsdonkin} show that (vi) (which agrees with (v) in this case) holds. It remains to reduce to the case that $G$ is connected and simply connected.

Assume that $G$ is connected and semisimple and let $G'\to G$ be the simply connected cover of $G$, with (connected) fixed points $H'$. Then $H'\to H^\circ$ is a central isogeny. We argue that (iii) for $G'$ implies (iii) for $G$. Indeed, for a $H^\circ$-representation $W$, one can write $\mathrm{Ind}_{H^\circ}^G W$ as the part of $\mathrm{Ind}_{H'}^{G'} W$ on which the kernel of $G'\to G$ acts trivially. This means that the kernel in (iii) for $G$ is a direct summand of the kernel in (iii) for $G'$ (and $W$ considered as a representation of $H'$). Thus, we have handled the case that $G$ is connected and semisimple.

Now take a general $G$. Let $G'\subset G$ be the derived subgroup of $G^\circ$, with fixed points $H'\subset H$. Let $D=G/G'$ which is linearly reductive (more precisely, $D^\circ$ is a torus and $\pi_0 D$ is of order prime to $\ell$). We get an exact sequence
\[
1\to H'\to H\to \overline{H}\to 1
\]
where $\overline{H}\subset D$ is linearly reductive (again, $\overline{H}^\circ$ is a torus and $\pi_0 \overline{H}$ is of order prime to $\ell$). We argue that (vi) for $G'$ implies (vi) for $G$, thereby finishing the proof. Take any representation $W$ of $H$. Then $W$ splits off $\mathrm{Ind}_{H'}^H(W|_{H'}) = W\otimes \mathcal O(\overline{H})$ as $\overline{H}$ has semisimple representation theory and hence $1$ is a direct summand of $\mathcal O(\overline{H})$. Thus, it suffices to see that for all representations $W'$ of $H'$, the representation $\mathrm{Ind}_{H'}^H W'\in \Ind\Perf(\ast/H)$ is in the subcategory generated under colimits by the image of $\Perf(\ast/G)$. Using (vi) for $H'$, and using that any representation $V'$ of $G'$ is a retract of $(\mathrm{Ind}_{G'}^G V)|_{G'}$, it suffices to see that for any representation $V$ of $G$, the representation $\mathrm{Ind}_{H'}^H(V|_{H'})\in \Ind\Perf(\ast/H)$ is in the subcategory generated under colimits by the image of $\Perf(\ast/G)$. But this is given by $V|_H\otimes \mathcal O(\overline{H})$. Thus, it suffices to see that $\mathcal O(\overline{H})\in \Ind\Perf(\ast/\overline{H})\subset \Ind\Perf(\ast/H)$ is in the subcategory generated under colimits by the image of $\Perf(\ast/G)$. This reduces us to $\overline{H}\subset D$, where the whole representation theory is semisimple, so that $\mathcal O(\overline{H})$ is a retract of $\mathcal O(D)|_{\overline{H}}$.
\end{proof}

To combine the results of this section with the preceding section, we also need the following observation. Note that in this section, we did not need the assumption that $\ell$ does not divide the order of $\pi_1(G^\circ)_{\mathrm{tor}}$; but also the operations of this section do not interfere with this condition:

\begin{proposition}\label{prop:fixedpointsfundamentalgroup} Assume that $G^\circ$ is reductive and that $P$ is a finite solvable group of order prime to $\ell$ acting on $G$ with fixed points $H=G^P$. If $\pi_1(G^\circ)_{\mathrm{tor}}$ is of order prime to $\ell$, then also $\pi_1(H^\circ)_{\mathrm{tor}}$ is of order prime to $\ell$.
\end{proposition}

\begin{proof} By induction, we can assume that $P=\mathbb Z/p\mathbb Z$ is cyclic of prime order $p\neq \ell$. We may evidently assume that $G$ is connected, and we can also pass to the derived group and its simply connected cover. As before, we can also assume that $G$ is simple. We need to see that $\pi_1(H)_{\mathrm{tor}}$ is of order prime to $\ell$. There is probably both a good reference and a good argument, but lacking both, we quickly run through the list of possibilities. In type $A$, all inner automorphisms give Levi subgroups, which always have simply connected derived subgroup. Outer automorphisms are involutions so $p=2$ and $\ell\neq 2$, and give either symplectic or orthogonal groups, whose fundamental group is a $2$-group. In the triality case, the fixed points of an order $3$ automorphism are either $G_2$ or $\mathrm{SL}_3$, in particular simply connected. In types $B$, $C$ and $D$ outside triality, either $p>2$ and the automorphism is inner and the centralizer a Levi subgroup (having simply connected derived subgroup) or $p=2$ and hence $\ell\neq 2$ while the fixed point subgroup is a classical group, with fundamental group a $2$-group. There remain the exceptional groups. Looking at the case of involutions (cf.~\cite{vanderKallenBrundan} or its reference \cite{SpringerInvolutions}) the only possibly dangerous case is the inner involution of $E_6$ which has fixed points of type $A_5\times A_1$ which might contribute a $3$-group to the fundamental group. But $E_6$ has center of order $3$ which survives to the fixed points. There remain the inner involutions of exceptional groups of prime order $p>2$. The possible cases that are not Levi subgroups are enumerated in \cite[Theorem 4.3.3]{HagueMcNinch}, but the list there includes also centralizers of elements of non-prime order. Restricting to prime orders, there is one case for $p=5$ which is $E_8$ with subgroup of type $A_4\times A_4$, so the fundamental group is necessarily a $5$-group. In all remaining cases $p=3$. For $G_2$, this gives a group of $A_2$; for $F_4$, of type $A_2\times A_2$; for $E_6$, of type $A_2\times A_2\times A_2$; these are all fine as the fundamental group is necessarily a $3$-group. For $E_7$, it is of type $A_2\times A_5$. This might a priori contribute a $2$-group, but again the center of $E_7$ shows that this does not happen. Finally, for $E_8$, one gets a group of type $A_2\times E_6$ or $A_8$; in both cases, the fundamental group is necessarily a $3$-group.
\end{proof}

\subsection{End of proof}

Finally, we can prove Theorem~\ref{thm:precisecoliminIndPerfreducechar}. The subtle part is to give a clean account of the reduction to the tame case. We take an approach that is inspired by our construction of the spectral action in Chapter~\ref{ch:spectral} below.

As preparation, consider any discrete group $W$ and a gerbe $\mathcal G\to \ast/W$ (on the fpqc site of $\mathrm{Spec} L$) banded by some linear-algebraic group $G$ over $L$, such that $G^\circ$ is reductive and $\pi_0 G$ of order prime to $\ell$. (Picking a point of $\mathcal G$, this is given by the classifying space of $\mathcal E_{\mathcal G}$ for some extension
\[
1\to G\to \mathcal E_{\mathcal G}\xrightarrow{\pi} W\to 1.
\]
This extension will later be given by the $L$-group.) In fact, slightly more generally, we will allow that $\mathcal G$ is a finite disjoint union of such gerbes, or even slightly more generally, that $\mathcal G$ is a stack whose pullback along $\ast\to \ast/W$ is a finite disjoint union of classifying stacks of such groups $G$. We will be mostly interested in sections of $\mathcal G\to \ast/W$ and we note that sections of $\mathcal G\to \ast/W$ are necessarily contained in one connected component, so that in practice we can reduce all questions to the case of gerbes.

For any anima $S$ mapping to $\ast/W$, we can look at the derived fpqc stack $\mathrm{Map}_{\ast/W}(S,\mathcal G)$ of maps $S\to \mathcal G$ over $\ast/W$. If $S$ is a finite set equipped with a (necessarily trivial) $W$-torsor, and $\mathcal G$ is connected, this is isomorphic to a product of copies of $\ast/G$. In general, it can be analyzed via resolutions. In fact, the $\infty$-category of anima $S$ over $\ast/W$ admits compact projective generators, given by finite sets equipped with a $W$-torsor, so the whole $\infty$-category of such $S$ is obtained by animating the category of finite sets equipped with a $W$-torsor.

Now we look at the functor taking any such anima $S$ over $\ast/W$ to the $L$-linear symmetric monoidal presentable stable $\infty$-category
\[
\Ind\Perf(\Map_{\ast/W}(S,\mathcal G)).
\]
We can also restrict this functor to compact projective objects; then it commutes with all finite coproducts. This commutation with finite coproducts reduces to the assertion that for $G$ and $H$ linear-algebraic groups over $L$ with $G^\circ$ and $H^\circ$ reductive and $\pi_0 G$ and $\pi_0 H$ of order prime to $\ell$, one has
\[
\Ind\Perf(\ast/G)\otimes_{D(L)} \Ind\Perf(\ast/H)\cong \Ind\Perf(\ast/G\times H),
\]
which follows from highest weight theory (which gives infinite-semiorthogonal decompositions of all intervening categories whose associated gradeds are representations of finite groups of order prime to $\ell$, where the result is simple). We may then extend back to all $S$ via left Kan extension, to get a colimit-preserving functor
\[
S\mapsto \Ind\Perf(\Map^\Sigma_{\ast/W}(S,\mathcal G))
\]
from anima $S$ over $\ast/W$ to $L$-linear symmetric monoidal presentable stable $\infty$-categories. (The notation here is purely symbolic; we like to think of this symmetric monoidal category as the $\Ind$-category of perfect complexes on some hypothetical stack $\Map^\Sigma_{\ast/W}(S,\mathcal G)$. The symbol $\Sigma$ here is in reference to the use in \cite[Section 5.5.8]{LurieHTT} in relation to sifted colimits.) By the universal property of left Kan extensions, there is a functorial comparison map
\[
\Ind\Perf(\Map^\Sigma_{\ast/W}(S,\mathcal G))\to \Ind\Perf(\Map_{\ast/W}(S,\mathcal G)).
\]

If $S=\ast/F_n\to \ast/W$ is given by the classifying space of a finite free group $F_n$ equipped with a map $F_n\to W$, and $\mathcal G$ is connected, then $\Map_{\ast/W}(\ast/F_n,\mathcal G)$ is given by $\prod_{i=1}^n \pi^{-1}(\gamma_i)/G$ where $G$ acts by simultaneous conjugation and $\gamma_i\in W$ is the image of the $i$-th generator of $F_n$ (and $\pi: \mathcal E_{\mathcal G}\to W$ is the projection). In general, one gets a finite disjoint union of such.

\begin{proposition}\label{prop:firststep} In this situation, the functor
\[
\Ind\Perf(\Map_{\ast/W}^\Sigma(\ast/F_n,\mathcal G))\to \Ind\Perf(\Map_{\ast/W}(\ast/F_n,\mathcal G))
\]
is fully faithful and the essential image is the full subcategory generated under colimits by $\Ind\Perf(\ast/G)$, which is equivalent to the category of $\mathcal O(\prod_{i=1}^n \pi^{-1}(\gamma_i))$-modules in $\Ind\Perf(\ast/G)$. If $\pi_1(G^\circ)_{\mathrm{tor}}$ is of order prime to $\ell$, the displayed functor is an equivalence.

An obvious variant holds true if $\mathcal G$ is a finite disjoint union of such gerbes, or even if $\mathcal G$ is a stack such that $\mathcal G\times_{\ast/W} \ast$ is a finite disjoint union of such gerbes.
\end{proposition}

\begin{proof} For the identification of $\Ind\Perf(\Map_{\ast/W}^\Sigma(\ast/F_n,\mathcal G))$ with $\mathcal O(\prod_{i=1}^n \pi^{-1}(\gamma_i))$-modules in $\Ind\Perf(\ast/G)$, we may use that the functor commutes with colimits to reduce to $n=1$. Then we can write $\ast/\mathbb Z$ as the pushout of $\ast\leftarrow \ast\sqcup \ast\to \ast$ in anima and use that $\Ind\Perf(\ast/G)$ is equivalent to $\mathcal O(G)$-modules in $\Ind\Perf(\ast/G\times G)$ by Barr--Beck, and that module categories base change. This description of $\Ind\Perf(\Map_{\ast/W}^\Sigma(\ast/F_n,\mathcal G))$ shows that the comparison functor is fully faithful. By Proposition~\ref{prop:generateperfectcomplexespowerofG}, this is an equivalence when $\pi_1(G^\circ)_{\mathrm{tor}}$ is of order prime to $\ell$.
\end{proof}

Assume that $\mathcal G$ is connected, banded by $G$. Writing $W$ as a sifted colimit of finite free groups, we see that
\[
\Ind\Perf(\Map_{\ast/W}^\Sigma(\ast/W,\mathcal G))
\]
is given by the category of modules over 
\[
\colim_{(n,F_n\to W)} \mathcal O(\prod_{i=1}^n \pi^{-1}(\gamma_i))
\]
in $\Ind\Perf(\ast/G)$. Thus, to prove Theorem~\ref{thm:precisecoliminIndPerfreducechar}, we have to see that if $W$ is a discretization of $W_E$ as in the formulation of that theorem, then the functor
\[
\Ind\Perf(\Map_{\ast/W}^\Sigma(\ast/W,\mathcal G))\to \Ind\Perf(\Map_{\ast/W}(\ast/W,\mathcal G))
\]
is an equivalence. Indeed, $\Map_{\ast/W}(\ast/W,\mathcal G)=Z^1(W,G)/G$, so the fully faithfulness is equivalent to the isomorphism
\[
\colim_{(n,F_n\to W)} \mathcal O(\prod_{i=1}^n \pi^{-1}(\gamma_i))\to \mathcal O(Z^1(W,G))
\]
of algebras in $\Ind\Perf(\ast/G)$ (as the composite of the functor and its right adjoint is given by tensoring along this map), while the essential surjectivity amounts to the assertion that $\Perf(Z^1(W,G)/G)$ is generated by $\Perf(\ast/G)$ under cones and retracts.

We will now show more generally that if $W$ is a discretization of $W_E$ as in Theorem~\ref{thm:precisecoliminIndPerfreducechar} and $\mathcal G$ is a stack over $\ast/W$ such that $\mathcal G\times_{\ast/W} \ast$ is a finite disjoint union of classifying stacks of linear-algebraic groups $G$ with $G^\circ$ reductive and $\pi_0 G$ and $\pi_1(G^\circ)_{\mathrm{tor}}$ of order prime to $\ell$, then
\[
\Ind\Perf(\Map_{\ast/W}^\Sigma(\ast/W,\mathcal G))\to \Ind\Perf(\Map_{\ast/W}(\ast/W,\mathcal G))
\]
is an equivalence.

We will prove this first when $W$ is a discretization of the tame part, i.e.~$W=\mathbb Z[\tfrac 1p]\rtimes \mathbb Z$. We can assume that $\mathcal G$ is connected. If $\mathcal G\times_{\ast/W} \ast$ does not stay connected, then there are no sections; and also $\Ind\Perf(\Map_{\ast/W}^\Sigma(\ast/W,\mathcal G))=0$. We can also assume that the corresponding extension
\[
1\to G\to \mathcal E_{\mathcal G}\to W\to 1
\]
splits as otherwise both sides are zero. We can write $\ast/W$ as a pushout
\[
\ast\leftarrow \ast/\mathbb Z\to \ast/F_2
\]
in anima, where $F_2$ is the free group on generators $\sigma$ and $\tau$, and the map $\mathbb Z\to F_2$ picks out the relation $\sigma^{-1}\tau\sigma\tau^{-q}$ in $W$. Using Proposition~\ref{prop:firststep} it remains to see that
\[
\Ind\Perf(\ast/G)\otimes_{\Ind\Perf(G/_{\mathrm{ad}} G)} \Ind\Perf((G\times G)/_{(\sigma,\tau)} G)\to \Ind\Perf(Z^1(W,G)/G)
\]
is an equivalence. The tensor product here is dual to writing $Z^1(W,G)$ as the fibre product
\[\xymatrix{
Z^1(W,G)\ar[r]\ar[d] & G\times G\ar[d]\\
\ast\ar[r] & G
}\]
where the right vertical map takes $(g,h)$ to $g^{-1}hgh^{-q}$. The result follows from Proposition~\ref{prop:generationclosedsubscheme}, Proposition~\ref{prop:generateperfectcomplexespowerofG} and Proposition~\ref{prop:generateskyscraper}.

It remains to analyze the wild part. This uses the full force of our results on fixed point groups.

\begin{proposition}\label{prop:wildpart} Let $P$ be a finite solvable group of order prime to $\ell$ with a map $P\to W$. Assume that $\mathcal G$ is a stack over $\ast/W$ such that $\mathcal G\times_{\ast/W} \ast$ is a finite disjoint union of classifying stacks of linear-algebraic groups $G$ with $G^\circ$ reductive and $\pi_0 G$ of order prime to $\ell$. Then the pushforward $\mathcal G^t$ of $\mathcal G$ along $\ast/W\to \ast/W^t$, where $W^t=W/P$, has pullback along $\ast\to \ast/W^t$ given by
\[
\Map_{\ast/W}(\ast/P,\mathcal G),
\]
and this is again a finite disjoint union of classifying stacks of linear-algebraic groups $H$ with $H^\circ$ reductive and $\pi_0 H$ of order prime to $\ell$. Moreover, the functor
\[
\Ind\Perf(\Map_{\ast/W}^\Sigma(\ast/P,\mathcal G))\to \Ind\Perf(\Map_{\ast/W}(\ast/P,\mathcal G))
\]
is an equivalence.

If all $G$ appearing in $\mathcal G\times_{\ast/W}\ast$ have the property that $\pi_1(G^\circ)_{\mathrm{tor}}$ is of order prime to $\ell$, then the same holds true for the groups $H$ appearing in $\mathcal G^t\times_{\ast/W^t} \ast$.
\end{proposition}

\begin{proof} We can assume that $W=P$. We can also assume that $\mathcal G$ is connected and that the associated extension
\[
1\to G\to \mathcal E_{\mathcal G}\to P\to 1
\]
splits, as otherwise the stack $\Map_{\ast/P}(\ast/P,\mathcal G)$ is empty (in which case Haboush's theorem implies that also $\colim_{(n,F_n\to P)} \mathcal O(\prod_{i=1}^n \pi^{-1}(\gamma_i))^G=0$ and hence also $\Ind\Perf(\Map_{\ast/P}^\Sigma(\ast/P,\mathcal G))=0$). Fixing a section, we get an action of $P$ on $G$ and $\Map_{\ast/P}(\ast/P,\mathcal G)=Z^1(P,G)/G$. If $P=\mathbb Z/p\mathbb Z$ is cyclic of prime order $p\neq \ell$, the proposition follows from the results of Section~\ref{sec:fixedpoints}, noting that the colimit in Lemma~\ref{lem:fixedpointsIndPerf} is computing
\[
\colim_{(n,F_n\to \mathbb Z/p\mathbb Z)} \mathcal O(\prod_{i=1}^n \pi^{-1}(\gamma_i))
\]
in $\Ind\Perf(\ast/G)$, using a specific presentation (and that $X^P=Z^1(P,G)$). Namely, inside pointed anima over $\ast/P$, one can write $\ast/P$ as the quotient of $\ast\sqcup_P \ast$ by the action of $P$; this writes $\ast/P$ as a colimit of connected pointed anima over $\ast/P$ of the form $\ast/F_n$, i.e.~$P$ as a geometric realization of free groups inside animated groups, which can be used to compute the displayed colimit. Dually, this writes $Z^1(P,G)$ as the $P$-fixed points inside the $X$ from Lemma~\ref{lem:fixedpointsIndPerf}. Proposition~\ref{prop:strongdonkinsubgroup} shows that the resulting fully faithful functor
\[
\Ind\Perf(\Map_{\ast/P}^\Sigma(\ast/P,\mathcal G))\to \Ind\Perf(\Map_{\ast/P}(\ast/P,\mathcal G))
\]
is essentially surjective. Proposition~\ref{prop:fixedpointsfundamentalgroup} ensures the statement about $\pi_1$'s.

For the general case, we can find a normal subgroup $P'\subset P$ and quotient $\overline{P}=\mathbb Z/p\mathbb Z$ and we can assume by induction that the result holds for $P'$. We can then consider the functor taking an anima $\overline{S}$ with a map $\overline{S}\to \ast/\overline{P}$ to
\[
\Ind\Perf(\Map_{\ast/P}(\overline{S}\times_{\ast/\overline{P}} \ast/P,\mathcal G)).
\]
If $\overline{S}$ is a finite set $I$ (equipped with a $\overline{P}$-torsor), this is given by
\[
\Ind\Perf(\Map_{\ast/P'}(\ast/P',\mathcal G)^I).
\]
By the structure of $\Map_{\ast/P'}(\ast/P',\mathcal G)$, this is the $I$-fold tensor product of the value on a point; it follows that in this case
\[
\Ind\Perf(\Map^\Sigma_{\ast/P}(\overline{S}\times_{\ast/\overline{P}} \ast/P,\mathcal G))\to \Ind\Perf(\Map_{\ast/P}(\overline{S}\times_{\ast/\overline{P}} \ast/P,\mathcal G))
\]
is an equivalence. It follows that the functor
\[
\overline{S}\mapsto \Ind\Perf(\Map^\Sigma_{\ast/P}(\overline{S}\times_{\ast/\overline{P}} \ast/P,\mathcal G))
\]
is the left Kan extension of its restriction to finite sets with a $\overline{P}$-torsor, where it agrees with the functor associated to the stack $\overline{\mathcal G}$ over $\ast/\overline{P}$ which is the pushforward of $\mathcal G$ along $\ast/P\to \ast/\overline{P}$. After pullback along $\ast\to \ast/\overline{P}$, this is
\[
\Map_{\ast/P}(\ast/P',\mathcal G),
\]
which we have seen to be a finite disjoint union of classifying stacks of the desired form. Thus, we conclude by appealing to the case of $\overline{P}=\mathbb Z/p\mathbb Z$ already established, for the stack $\overline{\mathcal G}$.
\end{proof}

If now $W$ is a discretization of the Weil group $W_E$, we take $P\subset W$ to be the wild part (a finite $p$-group, normal in $W$). Using Proposition~\ref{prop:wildpart}, let $\mathcal G^t$ be the pushforward of $\mathcal G$ along $\ast/W\to \ast/W^t$, where $W^t=W/P$ is the tame quotient of $W$; thus, its pullback along $\ast\to \ast/W^t$ is
\[
\Map_{\ast/W}(\ast/P,\mathcal G).
\]
Arguing as in the inductive part of the proof of Proposition~\ref{prop:wildpart}, we can replace $W$ by $W^t$ and $\mathcal G$ by $\mathcal G^t$. This reduces us to the tame case already handled.

\chapter{The Hecke action}

The time has come to put everything together. As before, let $E$ be any nonarchimedean local field with residue field $\Fq$ of residue characteristic $p$, and let $G$ be a reductive group over $E$. For any $\mathbb Z_\ell$-algebra $\Lambda$, we have defined $D_\lis(\Bun_G,\Lambda)$, we have the geometric Satake equivalence relating $\hat{G}$ to perverse sheaves on the Hecke stack, and we have studied the stack of $L$-parameters.

Our first task is to use the geometric Satake equivalence to define the Hecke operators on $D_\lis(\Bun_G,\Lambda)$. As in the last chapter, we work over a $\mathbb Z_\ell[\sqrt{q}]$-algebra $\Lambda$ in order to trivialize the cyclotomic twist in the geometric Satake equivalence; let $Q$ be a finite quotient of $W_E$ over which the action on $\hat{G}$ factors. If $\Lambda$ is killed by a power of $\ell$, then we can define Hecke operators in the following standard way. For any finite set $I$ and $V\in \Rep_\Lambda(\hat{G}\rtimes Q)^I$, we get a perverse sheaf $\mathcal S_V$ on $\Hloc_G^I$, which we can pull back to the global Hecke stack $\Hglob_G^I$; we denote its pullback still by $\mathcal S_V$. Using the correspondence
\[\xymatrix{
& \Hglob_G^I\ar_{p_1}[dl]\ar^{p_2}[dr]\\
\Bun_G && \Bun_G\times (\Div^1)^I
}\]
we get the Hecke operator
\[
T_V: D_\et(\Bun_G,\Lambda)\to D_\et(\Bun_G\times (\Div^1)^I,\Lambda): A\mapsto Rp_{2\ast}(p_1^\ast A\dotimes_\Lambda \mathcal S_V).
\]
By Corollary~\ref{cor:drinfeldfullyfaithful}, the target has $D_\et(\Bun_G\times [\ast/\underline{W_E^I}],\Lambda)$ as a full subcategory, and we will see below that $T_V$ will factor over this subcategory. Working $\infty$-categorically in order to have descent, and using a little bit of condensed formalism in order to deal with $W_E^I$ not being discrete, we can in fact rewrite
\[
\mathcal D_\et(\Bun_G\times [\ast/\underline{W_E^I}],\Lambda)\cong \mathcal D_\et(\Bun_G,\Lambda)^{BW_E^I}
\]
as the $W_E^I$-equivariant objects of the condensed $\infty$-category $\mathcal D_\et(\Bun_G,\Lambda)$; we will discuss the condensed structure below.

The following theorem summarizes the properties of the Hecke operators. In particular, it asserts that these functors are defined even when $\Lambda$ is not torsion.

\begin{theorem}[Theorem~\ref{thm:heckefiniteness}; Corollary~\ref{cor:heckeaction}, Proposition~\ref{prop:actionofWEreasonable}] For any $\mathbb Z_\ell[\sqrt{q}]$-algebra $\Lambda$, any finite set $I$, and any $V\in \Rep_\Lambda(\hat{G}\rtimes Q)^I$, there is a natural Hecke operator
\[
T_V: \mathcal D_\lis(\Bun_G,\Lambda)\to \mathcal D_\lis(\Bun_G,\Lambda)^{BW_E^I}.
\]
\begin{altenumerate}
\item[{\rm (i)}] Forgetting the $W_E^I$-action, i.e.~as an endofunctor of $\mathcal D_\lis(\Bun_G,\Lambda)$, the functor $T_V$ commutes with all limits and colimits, and preserves compact objects and universally locally acyclic objects. Letting $\mathrm{sw}^\ast: \Rep_\Lambda \hat{G}^I\to \Rep_\Lambda \hat{G}^I$ be the involution of Proposition~\ref{prop:chevalleyinvolution}, there are natural isomorphisms
\[
\mathbb D_{\mathrm{BZ}}(T_V(A))\cong T_{\mathrm{sw}^\ast V^\vee}(\mathbb D_{\mathrm{BZ}}(A))\ ,\ R\sHom_\lis(T_V(A),\Lambda)\cong T_{\mathrm{sw}^\ast V^\vee} R\sHom_\lis(A,\Lambda).
\]
\item[{\rm (ii)}] As a functor of $V$, it induces an exact $\Rep_\Lambda(Q^I)$-linear monoidal functor
\[
\Rep_\Lambda(\hat{G}\rtimes Q)^I\to \End_\Lambda(\mathcal D_\lis(\Bun_G,\Lambda)^\omega)^{BW_E^I},
\]
where the target denotes $W_E^I$-equivariant objects inside the condensed $\infty$-category $\End_\Lambda(\mathcal D_\lis(\Bun_G,\Lambda)^\omega)$ (equipped with the trivial $W_E^I$-action). Moreover, for any compact object $X\in \mathcal D_\lis(\Bun_G,\Lambda)^\omega$, there is some open subgroup $P$ of the wild inertia subgroup of $W_E$ such that for all $I$ and $V$, the $P^I$-action on $T_V(X)$ is trivial. In particular, one can write $\mathcal D_\lis(\Bun_G,\Lambda)^\omega$ as an increasing union of full stable $\infty$-subcategories $\mathcal D_\lis^P(\Bun_G,\Lambda)^\omega$ such that the Hecke action defines functors
\[
\Rep_\Lambda(\hat{G}\rtimes Q)^I\to \End_\Lambda(\mathcal D_\lis^P(\Bun_G,\Lambda)^\omega)^{B(W_E/P)^I}.
\]
\item[{\rm (iii)}] Varying also $I$, the functors of (ii) are functorial in $I$.
\end{altenumerate}
\end{theorem}

Here, functoriality in $I$ means, more precisely, that one treats the preceding objects as coCartesian fibrations over the category of finite sets, and the functors are then required to lift to the total space of these coCartesian fibrations.

In particular, the categories $D_\lis^P(\Bun_G,\Lambda)^\omega$ fit the bill of the discussion of Section~\ref{sec:excursionoperators}, so Theorem~\ref{thm:constructionexcursionoperators} gives a construction of excursion operators. To state the outcome, we make the following definitions as in the introduction.

\begin{definition}\leavevmode
\begin{altenumerate}
\item[{\rm (i)}] The Bernstein center of $G(E)$ is
\[
\mathcal Z(G(E),\Lambda) = \pi_0\mathrm{End}(\mathrm{id}_{\mathcal D(G(E),\Lambda)}) = \varprojlim_{K\subset G(E)} \mathcal Z(\Lambda[K\backslash G(E)/K])
\]
where $K$ runs over open pro-$p$ subgroups of $G(E)$, and $\Lambda[K\backslash G(E)/K]=\mathrm{End}_{G(E)}(c\text-\mathrm{Ind}_K^{G(E)} \Lambda)$ is the Hecke algebra of level $K$.
\item[{\rm (ii)}] The geometric Bernstein center of $G$ is
\[
\mathcal Z^{\mathrm{geom}}(G,\Lambda) = \pi_0\mathrm{End}(\mathrm{id}_{\mathcal D_\lis(\Bun_G,\Lambda)}).
\]
Inside $\mathcal Z^{\mathrm{geom}}(G,\Lambda)$, we let $\mathcal Z^{\mathrm{geom}}_{\mathrm{Hecke}}(G,\Lambda)$ be the subring of all endomorphisms $f: \mathrm{id}\to \mathrm{id}$ commuting with Hecke operators, in the sense that for all $V\in \mathrm{Rep}(\hat{G}^I)$ and $A\in \mathcal D_\lis(\Bun_G,\Lambda)$, one has $T_V(f(A))=f(T_V(A))\in \mathrm{End}(T_V(A))$.
\item[{\rm (iii)}] The spectral Bernstein center of $G$ is
\[
\mathcal Z^{\mathrm{spec}}(G,\Lambda) = \mathcal O(Z^1(W_E,\hat{G})_\Lambda)^{\hat{G}},
\]
the ring of global functions on $Z^1(W_E,\hat{G})_\Lambda\sslash \hat{G}$.
\end{altenumerate}
\end{definition}

The inclusion $\mathcal D(G(E),\Lambda)\hookrightarrow \mathcal D_\lis(\Bun_G,\Lambda)$ induces a map of algebra $\mathcal Z^{\mathrm{geom}}(G,\Lambda)\to \mathcal Z(G(E),\Lambda)$. This discussion will lead to the following corollary.

\begin{corollary} Assume that the order of $\pi_0 Z(G)$ is invertible in $\Lambda$. There is a canonical map
\[
\mathcal Z^{\mathrm{spec}}(G,\Lambda)\to \mathcal Z^{\mathrm{geom}}_{\mathrm{Hecke}}(G,\Lambda)\subset \mathcal Z^{\mathrm{geom}}(G,\Lambda),
\]
and in particular a map
\[
\Psi_G: \mathcal Z^{\mathrm{spec}}(G,\Lambda)\to\mathcal Z(G(E),\Lambda).
\]
\end{corollary}

In particular, if $\Lambda=L$ is an algebraically closed field over $\mathbb Z_\ell[\sqrt{q}]$, we get the following construction of $L$-parameters. (This works even if $\ell$ does divide the order of $\pi_0 Z(G)$.)

\begin{definition} Let $L$ be an algebraically closed field over $\mathbb Z_\ell[\sqrt{q}]$, and let $A\in D_\lis(\Bun_G,L)$ be a Schur-irreducible object, i.e.~$\mathrm{End}(A)=L$ as condensed algebras. Then there is a unique semisimple $L$-parameter
\[
\phi_A: W_E\to \hat{G}(L)\rtimes Q
\]
such that for all excursion data $(I,V,\alpha,\beta,(\gamma_i)_{i\in I})$ consisting of a finite set $I$, $V\in \Rep((\hat{G}\rtimes Q)^I)$, $\alpha: 1\to V|_{\hat{G}}$, $\beta: V|_{\hat{G}}\to 1$ and $\gamma_i\in W_E$ for $i\in I$, the endomorphism
\[
A=T_1(A)\xrightarrow{\alpha} T_V(A)\xrightarrow{(\gamma_i)_{i\in I}} T_V(A)\xrightarrow{\beta} T_1(A)=A
\]
is given by the scalar
\[
L\xrightarrow{\alpha} V\xrightarrow{(\phi_A(\gamma_i))_{i\in I}} V\xrightarrow{\beta} L.
\]
\end{definition}

We can apply this in particular in the case of irreducible smooth representations $\pi$ of $G(E)$. Concerning the $L$-parameters we construct, we can prove the following basic results. (In fact, we prove slightly finer results on the level of Bernstein centers.)

\begin{theorem}[Sections~\ref{sec:propertiescorrespondence}, \ref{sec:reprGE}]\leavevmode
\begin{altenumerate}
\item[{\rm (i)}] If $G=T$ is a torus, then $\pi\mapsto \varphi_\pi$ is the usual Langlands correspondence.
\item[{\rm (ii)}] The correspondence $\pi\mapsto \varphi_\pi$ is compatible with twisting.
\item[{\rm (iii)}] The correspondence $\pi\mapsto \varphi_\pi$ is compatible with central characters (cf.~\cite[10.1]{Borel79}).
\item[{\rm (iv)}] The correspondence $\pi\mapsto \varphi_\pi$ is compatible with passage to congradients (cf.~\cite{AdamsVoganContragredient}).
\item[{\rm (v)}] If $G'\to G$ is a map of reductive groups inducing an isomorphism of adjoint groups, $\pi$ is an irreducible smooth representation of $G(E)$ and $\pi'$ is an irreducible constitutent of $\pi|_{G'(E)}$, then $\varphi_{\pi'}$ is the image of $\varphi_\pi$ under the induced map $\hat{G}\to \hat{G'}$.
\item[{\rm (vi)}] If $G=G_1\times G_2$ is a product of two groups and $\pi$ is an irreducible smooth representation of $G(E)$, then $\pi=\pi_1\boxtimes \pi_2$ for irreducible smooth representations $\pi_i$ of $G_i(E)$, and $\varphi_\pi=\varphi_{\pi_1}\times \varphi_{\pi_2}$ under $\hat{G}=\hat{G}_1\times \hat{G}_2$.
\item[{\rm (vii)}] If $G=\mathrm{Res}_{E'|E} G'$ is the Weil restriction of scalars of a reductive group $G'$ over some finite separable extension $E'|E$, so that $G(E)=G'(E')$, then $L$-parameters for $G|E$ agree with $L$-parameters for $G'|E'$.
\item[{\rm (viii)}] The correspondence $\pi\mapsto \varphi_\pi$ is compatible with parabolic induction.
\item[{\rm (ix)}] For $G=\GL_n$ and supercuspidal $\pi$, the correspondence $\pi\mapsto \varphi_\pi$ agrees with the usual local Langlands correspondence \cite{LaumonRapoportStuhler}, \cite{HarrisTaylor}, \cite{HenniartGLn}.
\end{altenumerate}
\end{theorem}

\section{Condensed $\infty$-categories}

In order to meaningfully talk about $W_E^I$-equivariant objects in $D_\lis(\Bun_G,\Lambda)$, we need to give $\mathcal D_\lis(\Bun_G,\Lambda)$ the structure of a condensed $\infty$-category. This is in fact easy to do: We can associate to any extremally disconnected profinite set $S$ the $\infty$-category $\mathcal D_\lis(\Bun_G\times S,\Lambda)$. This is a full condensed $\infty$-subcategory of the condensed $\infty$-category $\mathcal D_\solid(\Bun_G,\Lambda)$, taking any profinite $S$ to $\mathcal D_\solid(\Bun_G\times S,\Lambda)$. The latter defines a hypersheaf in $S$, by v-hyperdescent of $\mathcal D_\solid(X,\Lambda)$ (as follows from the case of $\mathcal D(X_v,\Lambda)$). With this definition, it becomes a direct consequence of descent that
\[
\mathcal D_\solid(\Bun_G\times [\ast/\underline{W_E^I}],\Lambda)\cong \mathcal D_\solid(\Bun_G,\Lambda)^{BW_E^I},
\]
where the latter is the evaluation of the condensed $\infty$-category $\mathcal D_\solid(\Bun_G,\Lambda)$ on the condensed anima $BW_E^I$. More concretely, this is the $\infty$-category of objects $A\in \mathcal D_\solid(\Bun_G,\Lambda)$ together with a map of condensed animated groups $W_E^I\to \Aut(A)$. We see in particular that to define $\mathcal D_\solid(\Bun_G,\Lambda)^{BW_E^I}$, we do not need to know the full structure as a condensed $\infty$-category. Rather, we only need the structure as an $\infty$-category enriched in condensed anima. This structure on $\mathcal D_\solid(\Bun_G,\Lambda)$ induces a similar structure on $\mathcal D_\lis(\Bun_G,\Lambda)$.

For the discussion of Hecke operators, we observe in particular the following result, that follows directly from the discussion above.

\begin{proposition}\label{prop:imageheckeoperator} Pullback under $\Bun_G\times (\Div^1)^I\to \Bun_G\times [\ast/\underline{W_E^I}]$ induces a fully faithful functor
\[
\mathcal D_\lis(\Bun_G,\Lambda)^{BW_E^I}\hookrightarrow \mathcal D_\solid(\Bun_G,\Lambda)^{BW_E^I}\cong \mathcal D_\solid(\Bun_G\times [\ast/\underline{W_E^I}],\Lambda)\hookrightarrow \mathcal D_\solid(\Bun_G\times (\Div^1)^I,\Lambda).
\]
The essential image of the first functor consists of all objects $A\in \mathcal D_\solid(\Bun_G\times [\ast/\underline{W_E^I}],\Lambda)$ whose pullback to $\Bun_G$ lies in $\mathcal D_\lis(\Bun_G,\Lambda)$.$\hfill \Box$
\end{proposition}

In fact, this structure of $\mathcal D_\lis(\Bun_G,\Lambda)$ as an $\infty$-category enriched in condensed anima, in fact condensed animated $\Lambda$-modules, can be obtained in the following way from its structure as a $\Lambda$-linear stable $\infty$-category.

\begin{proposition}\label{prop:condensedstructureDlis} For $A\in \mathcal D_\lis(\Bun_G,\Lambda)^\omega$ and $B\in \mathcal D_\lis(\Bun_G,\Lambda)$, the condensed animated $\Lambda$-module $\Hom_{\mathcal D_\lis(\Bun_G,\Lambda)}(A,B)$ is relatively discrete over $\mathbb Z_\ell$.
\end{proposition}

In other words, the condensed structure on $\mathcal D_\lis(\Bun_G,\Lambda)$ can also be defined as the relatively discrete condensed structure when restricted to compact objects, and in general induced from this. In particular, when restricting attention to the compact objects $\mathcal D_\lis(\Bun_G,\Lambda)^\omega$, it is simply the relatively discrete condensed structure.

\begin{proof} Take some $b\in B(G)$ and $K\subset G_b(E)$ an open pro-$p$-subgroup, and let $f_K: \tilde{\mathcal M}^b/\underline{K}\to \Bun_G$ be the local chart. We can assume $A=f_{K\natural} \mathbb Z_\ell$, as these form a family of generators. By adjunction, it is enough to show that for any $B'\in \mathcal D_\lis(\tilde{\mathcal M}^b/\underline{K},\Lambda)$, the global sections $R\Gamma(\tilde{\mathcal M}^b/\underline{K},B')$ have the relatively discrete condensed $\mathbb Z_\ell$-module structure. We claim that the restriction map $R\Gamma(\tilde{\mathcal M}^b/\underline{K},B')\to R\Gamma([\ast/\underline{K}],B')$ is an isomorphism, where $[\ast/\underline{K}]\subset \tilde{\mathcal M}^b/\underline{K}$ is the base point. Without the condensed structure, this was proved in the proof of Proposition~\ref{prop:compactobjectsDlis}, but actually the proof applies with condensed structure (as Theorem~\ref{thm:partialcompactsupportvanishingsolid} remembers the condensed structure). But $R\Gamma([\ast/\underline{K}],B')$ is a direct summand of the stalk of $B'$ at $\ast$, which has the relatively discrete condensed $\mathbb Z_\ell$-module structure (as this is true for all objects of $\mathcal D_\lis(\ast,\Lambda)$).
\end{proof}

\section{Hecke operators}

The geometric Satake equivalence gives exact $\Rep_{\mathbb Z_\ell[\sqrt{q}]}(Q^I)$-linear monoidal functors
\[
\Rep_{\mathbb Z_\ell[\sqrt{q}]}(\hat{G}\rtimes Q)^I\to \mathrm{Sat}(\Hloc_G^I,\mathbb Z_\ell[\sqrt{q}]): V\mapsto \mathcal S_V,
\]
where the target category is defined as the inverse limit over $n$ of the similar categories with $\mathbb Z/\ell^n[\sqrt{q}]$-coefficients, see Section VI.7. Moreover, this association is functorial in $I$. We can compose with the functor $A\mapsto \mathbb D(A)^\vee$ (where the Verdier duality is relative to the projection $\Hloc_G^I\to [(\Div^1)^I/L^+ G]$) to get exact $\Rep_{\mathbb Z_\ell[\sqrt{q}]}(Q^I)$-linear monoidal functors
\[
\Rep_{\mathbb Z_\ell[\sqrt{q}]}(\hat{G}\rtimes Q)^I\to \mathcal D_\solid(\Hloc_G^I,\mathbb Z_\ell[\sqrt{q}]),
\]
functorially in $I$. Here, the functor $A\mapsto \mathbb D(A)^\vee$ is monoidal with respect to the usual convolution on perverse sheaves, and the convolution of Section~\ref{sec:soliddualizability} on the right. We note that as the convolution on $\mathcal D_\solid$ makes use only of pullback, tensor product, and $\pi_\natural$-functors, all of which are defined naturally on $\infty$-categories, this monoidal structure is actually a monoidal structure on the $\infty$-category $\mathcal D_\solid(\Hloc_G^I,\mathbb Z_\ell[\sqrt{q}])$. (We would have to work harder to obtain this structure when employing lower-$!$-functors, as we have not defined them in a sufficiently structured way.) Also, the functor from $\Rep_{\mathbb Z_\ell[\sqrt{q}]}(\hat{G}\rtimes Q)^I$ is monoidal in this setting, as on perverse sheaves there are no higher coherences to take care of.

This extends by linearity uniquely to an exact $\Rep_\Lambda(Q^I)$-linear monoidal functor
\[
\Rep_\Lambda(\hat{G}\rtimes Q)^I\to \mathcal D_\solid(\Hloc_G^I,\Lambda): V\mapsto \mathcal S_V';
\]
here, we implicitly use highest weight theory to show
\[
\Perf(B(\hat{G}\rtimes Q)^I_{\mathbb Z_\ell[\sqrt{q}]})\otimes_{\Perf(BQ^I_{\mathbb Z_\ell[\sqrt{q}]})} \Perf(\ast/Q^I_{\Lambda})\cong \Perf(\ast/(\hat{G}\rtimes Q)^I_\Lambda),
\]
and that the free stable $\infty$-category with an exact functor from $\Rep_\Lambda(\hat{G}\rtimes Q)^I$ is $\Perf(\ast/(\hat{G}\rtimes Q)^I_\Lambda)$.

Pulling back to the global Hecke stack, we get exact $\Rep_\Lambda(Q^I)$-linear monoidal functors
\[
\Rep_\Lambda(\hat{G}\rtimes Q)^I\to \mathcal D_\solid(\Hglob_G^I,\Lambda).
\]
On the other hand, there is a natural exact $\Rep_\Lambda(Q^I)$-linear monoidal functor
\[
\mathcal D_\solid(\Hglob_G^I,\Lambda)\to \mathrm{End}_{\mathcal D_\solid((\Div^1)^I,\Lambda)}(\mathcal D_\solid(\Bun_G\times (\Div^1)^I,\Lambda)),
\]
where the right-hand side denotes the $\mathcal D_\solid((\Div^1)^I,\Lambda)$-linear endofunctors. In particular, any $V\in \Rep_\Lambda(\hat{G}\rtimes Q)^I$ gives rise to a functor
\[
T_V: \mathcal D_\lis(\Bun_G,\Lambda)\to \mathcal D_\solid(\Bun_G\times (\Div^1)^I,\Lambda)
\]
via
\[
T_V(A) = p_{2\natural}(p_1^\ast A\soliddotimesLambda \mathcal S_V')
\]
where we consider the usual diagram
\[\xymatrix{
& \Hglob_G^I\ar_{p_1}[dl]\ar^{p_2}[dr]\\
\Bun_G && \Bun_G\times (\Div^1)^I.
}\]
Note that we have thus essentially used the translation of Proposition~\ref{prop:generalizeRfshriek} to extend the Hecke operators from the case of torsion rings $\Lambda$ to all $\Lambda$.

We note that if we pull back to the diagonal geometric point $\Spd C\to (\Div^1)^I$, where $C=\hat{\overline{E}}$, then this functor depends only on the composite
\[
\Rep_\Lambda(\hat{G}\rtimes Q)^I\to \mathcal D_\solid(\Hglob_G^I,\Lambda)\to \mathcal D_\solid(\Hglob_G^I\times_{(\Div^1)^I} \Spd C,\Lambda),
\]
and this composite factors naturally over $\Rep_\Lambda(\hat{G}^I)$.

\begin{proposition}\label{prop:heckelisse} For any $V\in \Rep_\Lambda(\hat{G}^I)$, the functor
\[
T_V: \mathcal D_\solid(\Bun_G\times \Spd C,\Lambda)\to \mathcal D_\solid(\Bun_G\times \Spd C,\Lambda)
\]
restricts to a functor
\[
T_V: \mathcal D_\lis(\Bun_G,\Lambda)\to \mathcal D_\lis(\Bun_G,\Lambda).
\]
\end{proposition}

\begin{proof} By highest weight theory, one can reduce to the case that $V$ is an exterior tensor product of representations of $\hat{G}$, and then by using that $V\mapsto T_V$ is monoidal, we can reduce to the tensor factors, which reduces us to the case $I=\{\ast\}$. Consider the Hecke diagram
\[
\Bun_{G,C}\xleftarrow{h_1} \Hglob_{G,C}\xrightarrow{h_2} \Bun_{G,C}
\]
where $\Hglob_{G,C}$ parametrizes over $S\in \Perf_C$ pairs of $G$-torsors $\mathcal E_1,\mathcal E_2$ on $X_S$ together with an isomorphism over $X_S\setminus S^\sharp$ meromorphic along $S^\sharp$. It suffices to see that for all $B\in D^\ULA(\Hloc_{G,\Spd C/\Div^1_X},\mathbb Z_\ell)$, the object
\[
h_{2\natural}(h_1^\ast A\soliddotimes q^\ast B^\vee)\in D_\lis(\Bun_{G,C},\Lambda).
\]
Now the category of such $B$ is generated (under colimits) by the objects $Rf_{\dot{w}\ast} \mathbb Z_\ell$ for
\[
f_{\dot{w}}: L^+\mathcal I\backslash \mathrm{Dem}_{\dot{w}}\to \Hloc_{G,\Spd C/\Div^1_X}
\]
a Demazure resolution (modulo action of Iwahori) of some Schubert variety in the affine flag variety. Using Proposition~\ref{prop:dualembeddingpushforward}, it thus suffices to see that for the corresponding push-pull correspondence on $\Bun_{G,C}$ with kernel given by the Demazure resolution, one has preservation of $D_\lis(\Bun_{G,C},\Lambda)$. But this is a proper and cohomologically smooth correspondence.
\end{proof}

\begin{theorem}\label{thm:heckefiniteness} For any $V\in \Rep_\Lambda(\hat{G}^I)$, the action of $T_V$ on $\mathcal D_\lis(\Bun_G,\Lambda)$ preserves all limits and colimits, and the full subcategories of compact objects, and of universally locally acyclic objects. Moreover, for the automorphism $\mathrm{sw}^\ast$ of $\Rep_\Lambda(\hat{G}^I)$ given by Proposition~\ref{prop:chevalleyinvolution}, there are natural isomorphisms
\[
\mathbb D_{\mathrm{BZ}}(T_V(A))\cong T_{\mathrm{sw}^\ast V^\vee}(\mathbb D_{\mathrm{BZ}}(A))\ ,\ R\sHom_\lis(T_V(A),\Lambda)\cong T_{\mathrm{sw}^\ast V^\vee} R\sHom_\lis(A,\Lambda).
\]
\end{theorem}

\begin{proof} The functor $V\mapsto T_V$ is monoidal. As $V$ is dualizable in the Satake category, with dual $V^\vee$, it follows that $T_V$ has a left and a right adjoint, given by $T_{V^\vee}$, and hence it follows formally that it preserves all limits and colimits, and compact objects. Now recall that $A\in D_\lis(\Bun_G,\Lambda)$ is universally locally acyclic if and only if for all compact $B\in D_\lis(\Bun_G,\Lambda)$, the object $R\Hom_\Lambda(B,A)\in D(\Lambda)$ is perfect, by Proposition~\ref{prop:lisseULAbunG}. Thus, the preservation of universally locally acyclic objects follows by adjointness from the preservation of compact objects.

For the duality statements, we note that, for $\pi:\Bun_G\to \ast$ the projection, there are natural isomorphisms
\[
\pi_\natural(T_V(A)\soliddotimesLambda B)\cong \pi_\natural(A\soliddotimesLambda T_{\mathrm{sw}^\ast V}(B)),
\]
as follows from the definition of the Hecke operator, and Proposition~\ref{prop:chevalleyinvolution}: Both sides identify with the homology of $\Hglob_G^I\times_{(\Div^1)^I} \Spd C$ with coefficients in $h_1^\ast A\soliddotimesLambda h_2^\ast B\soliddotimesLambda \mathcal S_V'$. The displayed equation implies the statement for Bernstein--Zelevinsky duals by also using that $T_{\mathrm{sw}^\ast V^\vee}$ is right adjoint to $T_{\mathrm{sw}^\ast V}$, and the statement for naive duals by using that $T_{\mathrm{sw}^\ast V^\vee}$ is left adjoint to $T_{\mathrm{sw}^\ast V}$.
\end{proof}

Composing Hecke operators, we get the following corollary.

\begin{corollary} For any $V\in \Rep_\Lambda(\hat{G}\rtimes Q)^I$, the functor
\[
T_V: \mathcal D_\lis(\Bun_G,\Lambda)\to \mathcal D_\solid(\Bun_G\times (\Div^1_X)^I,\Lambda)
\]
takes image in the full subcategory $\mathcal D_\solid(\Bun_G\times [\ast/\underline{W_E^I}],\Lambda)$; moreover, all objects in the image have the property that their pullback to $\mathcal D_\solid(\Bun_G,\Lambda)$ lies in $\mathcal D_\lis(\Bun_G,\Lambda)$, so by Proposition~\ref{prop:imageheckeoperator} the functor $T_V$ induces a functor
\[
\mathcal D_\lis(\Bun_G,\Lambda)\to \mathcal D_\lis(\Bun_G,\Lambda)^{BW_E^I}.
\]
\end{corollary}

\begin{proof} We only need to see that the image lands in $\mathcal D_\solid(\Bun_G\times [\ast/\underline{W_E^I}],\Lambda)$; the rest follows from Proposition~\ref{prop:heckelisse}. One can reduce to the case that $V$ is an exterior tensor product of $|I|$ representations $V_i\in \Rep_\Lambda(\hat{G}\rtimes Q)$ --- one can always find a, possibly infinite, resolution by such exterior tensor products that involves only finitely many weights of $\hat{G}^I$, and thus induces a resolution in $\mathcal D_\solid(\Hloc_G^I,\Lambda)$ --- and thus reduce to $I=\{\ast\}$. By Corollary~\ref{cor:drinfeldfullyfaithfulsolid0}, it suffices to see that the pullback to $\mathcal D_\solid(\Bun_G\times \Spd C,\Lambda)$ lies in $\mathcal D_\solid(\Bun_G,\Lambda)$. But by Proposition~\ref{prop:heckelisse}, we know that it lies in $\mathcal D_\lis(\Bun_G\times \Spd C,\Lambda)$, and $\mathcal D_\lis(\Bun_G,\Lambda)\to \mathcal D_\lis(\Bun_G\times \Spd C,\Lambda)$ is an equivalence by Proposition~\ref{prop:DlisBunG}.
\end{proof}

Finally, we get the following Hecke action.

\begin{corollary}\label{cor:heckeaction} Endowing the stable $\mathbb Z_\ell$-linear $\infty$-category $\mathcal D_\lis(\Bun_G,\Lambda)^\omega$ with the relatively discrete condensed structure, the Hecke action defines exact $\Rep_\Lambda(Q^I)$-linear monoidal functors
\[
\Rep_\Lambda (\hat{G}\rtimes Q)^I\to \End_\Lambda(\mathcal D_\lis(\Bun_G,\Lambda)^\omega)^{BW_E^I},
\]
functorially in $I$.
\end{corollary}

\section{Cohomology of local Shimura varieties}\label{sec:localShimura}

Theorem~\ref{thm:heckefiniteness} encodes strong finiteness properties for the cohomology of local Shimura varieties, giving unconditional proofs, and refinements, of the results of \cite[Section 6]{RapoportViehmann}. For this, we first specialize to $E=\mathbb Q_p$ as this is the standard setting of local Shimura varieties. Consider any local Shimura datum, consisting of a reductive group $G$ over $\mathbb Q_p$, a conjugacy class of minuscule cocharacters $\mu: \mathbb G_m\to G_{\overline{\mathbb Q}_p}$ with field of definition $E|\mathbb Q_p$ and some element $b\in B(G,\mu)\subset B(G)$. (Beware that we are making a small sin here in changing the meaning of the letter $E$, using it now in its usual meaning as a reflex field.) In \cite[Lecture 24]{Berkeley}, we construct a tower of partially proper smooth rigid-analytic spaces
\[
(\mathcal M_{(G,b,\mu),K})_{K\subset G(\mathbb Q_p)}
\]
over $\breve{E}$, equipped with a Weil descent datum. Each object in the tower carries an action of $G_b(\mathbb Q_p)$, and the tower carries an action of $G(\mathbb Q_p)$. Following Huber \cite{Huberladic}, one defines
\[
R\Gamma_c(\mathcal M_{(G,b,\mu),K,C},\mathbb Z_\ell) = \varinjlim_U R\Gamma_c(U,\mathbb Z_\ell)
\]
where $U\subset \mathcal M_{(G,b,\mu),K,C}$ runs through quasicompact open subsets, and one defines $R\Gamma_c(U,\mathbb Z_\ell)=\varprojlim_m R\Gamma_c(U,\mathbb Z/\ell^m \mathbb Z)$. This carries an action of $G_b(\mathbb Q_p)$ as well as an action of the Weil group $W_E$.

\begin{theorem}\label{thm:cohomologylocalshimura} The complex $R\Gamma_c(\mathcal M_{(G,b,\mu),K,C},\mathbb Z_\ell)$ is naturally a complex of smooth $G_b(\mathbb Q_p)$-representations, and, if $K$ is pro-$p$, a compact object in $D(G_b(\mathbb Q_p),\mathbb Z_\ell)$. Moreover, the action of $W_E$ is continuous.
\end{theorem}

In particular, each $H^i_c(\mathcal M_{(G,b,\mu),K,C},\mathbb Z_\ell)$ is a finitely generated smooth $G_b(\mathbb Q_p)$-representation. By descent, this is true even for all $K$ (not necessarily pro-$p$).

\begin{proof} Let $f_K: \mathcal M_{(G,b,\mu),K,C}\to \Spa C$ be the projection. Up to shift, we can replace $\mathbb Z_\ell$ by the dualizing complex $Rf_K^! \mathbb Z_\ell$. Now by Proposition~\ref{prop:generalizeRfshriek}, one has
\[
Rf_{K!} Rf_K^! \mathbb Z_\ell|_U\cong f_{K\natural} \mathbb Z_\ell|_U
\]
for any quasicompact $U\subset \mathcal M_{(G,b,\mu),K,C}$. As the left-hand side is perfect, it is given by its limit over reductions modulo $\ell^m$. We see that $H^i_c(\mathcal M_{(G,b,\mu),K,C},\mathbb Z_\ell)$ can be identified with $H^i(f_{K\natural} \mathbb Z_\ell)$ up to shift.

Now $\mu$ gives rise to a Hecke operator $T_\mu=T_{V_\mu}$ where $V_\mu$ is the highest weight representation of weight $\mu$. It corresponds to the Hecke correspondence on $\Bun_{G,C}$ parametrizing modifications of type $\mu$; this Hecke correspondence is proper and smooth over both factors. We apply $T_\mu$ to the compact object
\[
A=j_! c\text-\mathrm{Ind}_K^{G(\mathbb Q_p)} \mathbb Z_\ell\in D_\lis(\Bun_G,\mathbb Z_\ell)
\]
where $j: \Bun_G^1\cong [\ast/\underline{G(\mathbb Q_p)}]\hookrightarrow \Bun_G$ is the open immersion. By Theorem~\ref{thm:heckefiniteness}, also $T_\mu(A)$ is compact. By Proposition~\ref{prop:compactlisseBunG}, it follows that also $i^{b\ast} T_\mu(A)\in D_\lis(\Bun_G^b,\mathbb Z_\ell)\cong D(G_b(\mathbb Q_p),\mathbb Z_\ell)$ is compact. But this is, up to shift again, precisely $f_{K\natural} \mathbb Z_\ell$, by the identification of $\mathcal M_{(G,b,\mu),K,C}$ with the space of modifications of $G$-torsors of type $\mu$ from the $G$-bundle $\mathcal E_b$ to the $G$-bundle $\mathcal E_1$, up to the action of $K$ (cf.~\cite[Lecture 23, 24]{Berkeley}).

Descending to $E$, note that $T_\mu$ can be defined with values in $D_\solid(\Bun_G\times \Spd E/\phi^{\Z},\mathbb Z_\ell)$, and takes values in those sheaves whose pullback to $\Bun_{G,C}$ lies in $D_\lis(\Bun_G,\mathbb Z_\ell)$. Thus $W_E$, as a condensed group, acts on $i^{b\ast} T_\mu(A)\in D_\lis(\Bun_G^b,\mathbb Z_\ell)\cong D(G_b(\mathbb Q_p),\mathbb Z_\ell)$ considered as representations on condensed $\mathbb Z_\ell$-modules. In classical language, this means that the action is continuous.
\end{proof}

In particular, for each admissible representation $\rho$ of $G_b(\mathbb Q_p)$ on a $\mathbb Z_\ell$-algebra $\Lambda$, the complex
\[
R\Hom_{G_b(\mathbb Q_p)}(R\Gamma_c(\mathcal M_{(G,b,\mu),K,C},\mathbb Z_\ell),\rho)
\]
is a perfect complex of $\Lambda$-modules. Passing to the colimit over $K$, one obtains at least on each cohomology group an admissible $G(\mathbb Q_p)$-representation. In fact, as $T_\mu$ is left adjoint to $T_{\mu^\vee}$, we see that this is (up to shift) given by
\[
i^{1\ast} T_{\mu^\vee}(Ri^b_\ast[\rho])\in D_\lis(\Bun_G^1,\Lambda)\cong D(G(\mathbb Q_p),\Lambda).
\]
Here $i^b: \Bun_G^b\hookrightarrow \Bun_G$ is the inclusion, and $[\rho]\in D(G_b(\mathbb Q_p),\Lambda)\cong D_\lis(\Bun_G^b,\Lambda)$ can be a complex of smooth $G_b(\mathbb Q_p)$-representations. This shows in particular that there is in fact a natural complex of admissible $G(\mathbb Q_p)$-representations underlying
\[
\mathrm{colim}_K R\Hom_{G_b(\mathbb Q_p)}(R\Gamma_c(\mathcal M_{(G,b,\mu),K,C},\mathbb Z_\ell),\rho).
\]

Assuming again that $\rho$ is admissible, one can pull through Verdier duality,
\[\begin{aligned}
i^{1\ast} T_{\mu^\vee}(Ri^b_\ast[\rho])&\cong i^{1\ast} T_{\mu^\vee}(Ri^b_\ast\mathbb D([\rho^\vee]))\\
&\cong i^{1\ast} T_{\mu^\vee}(\mathbb D(i^b_! [\rho^\vee]))\\
&\cong i^{1\ast} \mathbb D(T_{\mathrm{sw}^\ast \mu}(i^b_! [\rho^\vee]))\\
&\cong \mathbb D(i^{1\ast} T_{\mathrm{sw}^\ast \mu}(i^b_! [\rho^\vee])).
\end{aligned}\]
As $T_{\mathrm{sw}^\ast \mu}$ also preserves compact objects, it follows that \cite[Remark 6.2 (iii)]{RapoportViehmann} has a positive answer: If $\Lambda=\overline{\mathbb Q}_\ell$ and $\rho$ has finite length, then also each cohomology group of
\[
\mathrm{colim}_K R\Hom_{G_b(\mathbb Q_p)}(R\Gamma_c(\mathcal M_{(G,b,\mu),K,C},\mathbb Z_\ell),\rho)
\]
has finite length as $G(\mathbb Q_p)$-representation. Indeed, with $\overline{\mathbb Q}_\ell$-coefficients, the category of smooth representations has finite global dimension, and hence being compact is equivalent to each cohomology group being finitely generated. Compact objects are preserved under the Hecke operators, and so we see that each cohomology group is finitely generated. Being also admissible, it is then of finite length by Howe's theorem \cite[VI.6.3]{Renard}.

The same arguments apply to prove Corollary~\ref{cor:finitenessmain}. Let us recall the setup. We start with a general $E$ now. As in \cite[Lecture XXIII]{Berkeley}, for any collection $\{\mu_i\}_i$ of conjugacy classes of cocharacters with fields of definition $E_i/E$ and $b\in B(G)$, there is a tower of moduli spaces of local shtukas
\[
f_K: (\Sht_{(G,b,\mu_\bullet),K})_{K\subset G(E)}\to \prod_{i\in I} \Spd \breve{E}_i
\]
as $K$ ranges over compact open subgroups of $G(E)$, equipped with compatible \'etale period maps
\[
\pi_K: \Sht_{(G,b,\mu_\bullet),K}\to \Gr^\tw_{G,\prod_{i\in I} \Spd \breve{E}_i,\leq \mu_\bullet}.
\]
Here, $\Gr^\tw_{G,\prod_{i\in I} \Spd \breve{E}_i}\to \prod_{i\in I} \Spd \breve{E}$ is a certain twisted form of the convolution affine Grassmannian, cf.~\cite[Section 23.5]{Berkeley}. Let $W$ be the exterior tensor product $\boxtimes_{i\in I} V_{\mu_i}$ of highest weight representations, and $\mathcal S_W$ the corresponding sheaf on $\Gr^\tw_{G,\prod_{i\in I} \Spd \breve{E}_i}$. More precisely, away from Frobenius-twisted partial diagonals, $\Gr^\tw_{G,\prod_{i\in I} \Spd \breve{E}_i}$ is isomorphic to the Beilinson--Drinfeld Grassmannian $\Gr_{G,\prod_{i\in I} \Spd \breve{E}_i}$, and we have defined $\mathcal S_W$ on this locus. One can uniquely extend over these Frobenius-twisted partial diagonals to universally locally acyclic, necessarily perverse, sheaves, as in the discussion of the fusion product. We continue to write $\mathcal S_W$ for its pullback to $\Sht_{(G,b,\mu_\bullet),K}$. Let $\mathcal S_W' = \mathbb D(\mathcal S_W)^\vee$ be the corresponding solid sheaf. By Proposition~\ref{prop:generalizeRfshriek}, with torsion coefficients $f_{K\natural} \mathcal S_W'$ agrees with $Rf_{K!} \mathcal S_W$, but $f_{K\natural} \mathcal S_W'$ is well-defined in general.

\begin{proposition} The sheaf
\[
f_{K\natural} \mathcal S_W'\in \mathcal D_\solid([\ast/\underline{G_b(E)}]\times \prod_{i\in I} \Spd \breve{E}_i,\Lambda)
\]
is equipped with partial Frobenii, thus descends to an object of
\[
\mathcal D_\solid([\ast/\underline{G_b(E)}]\times \prod_{i\in I} \Spd \breve{E}_i/\varphi_i^{\mathbb Z},\Lambda).
\]
This object lives in the full $\infty$-subcategory
\[
\mathcal D(G_b(E),\Lambda)^{B\prod_{i\in I} W_{E_i}}\subset \mathcal D_\solid([\ast/\underline{G_b(E))}]\times \prod_{i\in I} \Spd \breve{E}_i/\varphi_i^{\mathbb Z},\Lambda),
\]
and its restriction to $\mathcal D(G_b(E),\Lambda)$ is compact. In particular, for any admissible representation $\rho$ of $G_b(E)$, the object
\[
R\Hom_{G_b(E)}(f_{K\natural} \mathcal S_W',\rho)\in \mathcal D(\Lambda)^{B\prod_{i\in I} W_{E_i}}
\]
is a representation of $\prod_{i\in I} W_{E_i}$ on a perfect complex of $\Lambda$-modules. Taking the colimit over $K$, this gives rise to a complex of admissible $G(E)$-representations
\[
\varinjlim_K R\Hom_{G_b(E)}(f_{K\natural} \mathcal S_W',\rho)
\]
equipped with a $\prod_{i\in I} W_{E_i}$-action.

If $\rho$ is compact, then so is
\[
\varinjlim_K R\Hom_{G_b(E)}(f_{K\natural} \mathcal S_W',\rho)
\]
as a complex of $G(E)$-representations.
\end{proposition}

\begin{proof} The key observation is that $f_{K\natural} \mathcal S_W'$ can be identified with $T_W(j_![c\text-\mathrm{Ind}_K^{G(E)} \Lambda])|_{\Bun_G^b}$. A priori, for the latter, we have to look at the moduli space $\mathcal M$ of modifications of type bounded by $\mu_\bullet$ from $\mathcal E_b$ to the trivial vector bundle, up to the action of $K$, and take the homology of $\mathcal M$ with coefficients in $\mathcal S_W'$; more precisely, the relative homology of $\mathcal M\to \prod_{i\in I} \Spd \breve{E}_i/\phi_i^\Z$. After pull back to $\prod_{i\in I} \Spd \breve{E}_i$, there is a natural map from $\mathcal M$ to $\Sht_{(G,b,\mu_\bullet,K)}$ that is an isomorphism away from Frobenius-twisted partial diagonals. Indeed, $\Sht_{(G,b,\mu_\bullet,K)}$ parametrizes $G$-torsors over $Y_S$ together with an isomorphism with their Frobenius pullback away from the given points, together with a level-$K$-trivialization of the $G$-bundle near $\{\pi=0\}$. This induces two vector bundles on $X_S$, given by the bundles near $\{\pi=0\}$ and near $\{[\varpi]=0\}$, and these are identified away from the images of the punctures in $X_S$. As long as their images in $X_S$ are disjoint, one can reverse this procedure. Now the fusion compatibility of $\mathcal S_W$ (and thus $\mathcal S_W'$) implies the desired result.

In particular, this shows that $f_{K\natural} \mathcal S_W'$ admits natural partial Frobenius operators. The rest of the proof is now as before.
\end{proof}

\section{$L$-parameter}

We can now define $L$-parameter.

\begin{defprop}\label{defprop:Lparameters} Let $L$ be an algebraically closed field over $\mathbb Z_\ell[\sqrt{q}]$, and let $A\in D_\lis(\Bun_G,L)$ be a Schur-irreducible object, i.e.~$\mathrm{End}(A)=L$ as condensed algebras. Then there is a unique semisimple $L$-parameter
\[
\phi_A: W_E\to \hat{G}(L)\rtimes Q
\]
such that for all excursion data $(I,V,\alpha,\beta,(\gamma_i)_{i\in I})$ consisting of a finite set $I$, $V\in \Rep((\hat{G}\rtimes Q)^I)$, $\alpha: 1\to V|_{\hat{G}}$, $\beta: V|_{\hat{G}}\to 1$ and $\gamma_i\in W_E$ for $i\in I$, the endomorphism
\[
A=T_1(A)\xrightarrow{\alpha} T_V(A)\xrightarrow{(\gamma_i)_{i\in I}} T_V(A)\xrightarrow{\beta} T_1(A)=A
\]
is given by the scalar
\[
L\xrightarrow{\alpha} V\xrightarrow{(\phi_A(\gamma_i))_{i\in I}} V\xrightarrow{\beta} L.
\]
\end{defprop}

\begin{proof} By the arguments of Section~\ref{sec:excursionoperators}, we can build excursion data as required for Proposition~\ref{prop:excursionoperatorsgiveLparameters}.
\end{proof}

\section{The Bernstein center}

As before, there is the problem that the stack $Z^1(W_E,\hat{G})/\hat{G}$ of $L$-parameters is not quasicompact, but an infinite disjoint union. We can now actually decompose $\mathcal D_\lis(\Bun_G,\Lambda)$ into a direct product according to the connected components of $Z^1(W_E,\hat{G})$. We start with the following observation.

\begin{proposition}\label{prop:actionofWEreasonable} Let $A\in \mathcal D_\lis(\Bun_G,\Lambda)^\omega$ be any compact object. Then there is an open subgroup $P\subset W_E$ of the wild inertia subgroup such that for all finite sets $I$ and all $V\in \Rep((\hat{G}\rtimes Q)^I)$, the object
\[
T_V(A)\in \mathcal D_\lis(\Bun_G,\Lambda)^{BW_E^I}
\]
lies in the full $\infty$-subcategory
\[
\mathcal D_\lis(\Bun_G,\Lambda)^{B(W_E/P)^I}\subset \mathcal D_\lis(\Bun_G,\Lambda)^{BW_E^I}.
\]
\end{proposition}

\begin{proof} First, note that indeed the functor
\[
\mathcal D_\lis(\Bun_G,\Lambda)^{B(W_E/P)^I}\to \mathcal D_\lis(\Bun_G,\Lambda)^{BW_E^I}.
\]
is fully faithful; this follows from fully faithfulness of the pullback functor
\[
f^\ast: \mathcal D_\solid(\Bun_G\times [\ast/\underline{(W_E/P)^I}],\Lambda)\to \mathcal D_\solid(\Bun_G\times [\ast/\underline{W_E^I}],\Lambda),
\]
which in turn follows from $f_\natural \Lambda\cong \Lambda$ (and the projection formula for $f_\natural$), which can be deduced via base change from the case of $[\ast/\underline{W_E^I}]\to [\ast/\underline{(W_E/P)^I}]$, or after pullback to a v-cover $\Spa C\to [\ast/\underline{(W_E/P)^I}]$, for $[\Spa C/\underline{P^I}]\to \Spa C$, where it amounts to the vanishing of the $\Lambda$-homology of $P^I$.

Now note that if $P^I$ acts trivially on $T_V(A)$ and on $T_W(A)$ for two $V,W\in \Rep_{\mathbb Z_\ell}((\hat{G}\rtimes Q)^I)$, then it also acts trivially on $T_{V\otimes W}(A)=T_V(T_W(A))=T_W(T_V(A))$: Indeed, the $W_E^{I\sqcup I}$-action on $T_V(T_W(A))\cong T_{V\boxtimes W}(A)\cong T_W(T_V(A))$ is trivial on $P^{I\sqcup \emptyset}$ and $P^{\emptyset\sqcup I}$, thus on $P^{I\sqcup I}$, and hence the diagonal $W_E^I$-action is trivial on $P^I$. Using reductions to exterior tensor products, we can also reduce to $I=\{\ast\}$. Then if $V\in \Rep_{\mathbb Z_\ell}(\hat{G}\rtimes Q)$ is a $\otimes$-generator, it follows that it suffices that $P$ acts trivially on $T_V(A)$. But
\[
(\mathcal D_\lis(\Bun_G,\Lambda)^\omega)^{BW_E}=\bigcup_P (\mathcal D_\lis(\Bun_G,\Lambda)^\omega)^{B(W_E/P)}
\]
as for any relatively discrete condensed animated $\mathbb Z_\ell$-algebra $R$ with a map $\mathbb Z_\ell[W_E]\to R$, the map factors over $\mathbb Z_\ell[W_E/P]$ for some $P$. Indeed, we may restrict to $\mathbb Z_\ell[I_E]$, and then (as $I_E$ is compact) the image is contained in some finitely generated $\mathbb Z_\ell$-submodule $R_0\subset R$, so we can assume that $R$ is finite over $\mathbb Z_\ell$; but then $\mathrm{Aut}_{\mathbb Z_\ell}(R)$ is profinite, and locally pro-$\ell$, so the map $I_E\to R^\times\subset \mathrm{Aut}_{\mathbb Z_\ell}(R)$ factors over $I_E/P$ for some $P$.
\end{proof}

Fix some open subgroup $P$ of the wild inertia subgroup of $W_E$, and let
\[
\mathcal D_\lis^P(\Bun_G,\Lambda)^{\omega}\subset \mathcal D_\lis(\Bun_G,\Lambda)^\omega
\]
be the full $\infty$-subcategory of all $A$ such that $P^I$ acts trivially on $T_V(A)$ for all $V\in \Rep((\hat{G}\rtimes Q)^I)$. Pick $W\subset W_E/P$ a discrete dense subgroup, by discretizing the tame inertia, as before. Then Theorem~\ref{thm:constructionexcursionoperators} gives a canonical map of algebras
\[
\mathrm{Exc}(W,\hat{G})\to \mathcal Z(\mathcal D_\lis^P(\Bun_G,\Lambda)^\omega) = \pi_0\mathrm{End}(\mathrm{id}_{\mathcal D_\lis^P(\Bun_G,\Lambda)^\omega}).
\]

As $\mathrm{Exc}(W,\hat{G})\otimes \Lambda\to \O(Z^1(W_E/P,\hat{G})_\Lambda)^{\hat{G}}$ is a universal homeomorphism, there are in particular idempotents corresponding to the connected components of $Z^1(W_E/P,\hat{G})_\Lambda$. Their action on $\mathcal D_\lis^P(\Bun_G,\Lambda)^\omega$ then induces a direct sum decomposition
\[
\mathcal D_\lis^P(\Bun_G,\Lambda)^\omega = \bigoplus_{c\in \pi_0 Z^1(W_E/P,\hat{G})_\Lambda} \mathcal D_\lis^c(\Bun_G,\Lambda)^\omega.
\]

Taking now a union over all $P$, we get a direct sum decomposition
\[
\mathcal D_\lis(\Bun_G,\Lambda)^\omega = \bigoplus_{c\in \pi_0 Z^1(W_E,\hat{G})_\Lambda} \mathcal D_\lis^c(\Bun_G,\Lambda)^\omega.
\]
On the level of $\Ind$-categories, this gives a direct product
\[
\mathcal D_\lis(\Bun_G,\Lambda) = \prod_{c\in \pi_0 Z^1(W_E,\hat{G})_\Lambda} \mathcal D_\lis^c(\Bun_G,\Lambda).
\]

Note in particular that any Schur-irreducible object $A\in \mathcal D_\lis(\Bun_G,\Lambda)$ necessarily lies in one of these factors, given by some connected component $c$ of $Z^1(W_E,\hat{G})_\Lambda$; and then the $L$-parameter $\phi_A$ of $A$ necessarily lies in this connected component.

Using excursion operators, we get the following result on the ``Bernstein center''. This is a generalization of results of Helm--Moss, \cite{HelmMoss}, noting that by the fully faithful functor $\mathcal D(G(E),\Lambda)\hookrightarrow \mathcal D_\lis(\Bun_G,\Lambda)$, there is a map of algebras
\[
\mathcal Z^{\mathrm{geom}}(G,\Lambda)\to \mathcal Z(G(E),\Lambda)
\]
to the usual Bernstein center of smooth $G(E)$-representations on $\Lambda$-modules. {\it From now on we assume that the order of $\pi_0 Z(G)$ is invertible in $\Lambda$.}

\begin{theorem}\label{thm:bernsteincenter} Assume that the order of $\pi_0 Z(G)$ is invertible in $\Lambda$. There is a natural map
\[
\mathcal Z^{\mathrm{spec}}(G,\Lambda)\to \mathcal Z^{\mathrm{geom}}(G,\Lambda)
\]
compatible with the above decomposition into connected components. Moreover, for all finite sets $I$, all $V\in \Rep_\Lambda(\hat{G}^I)$, and all $A\in D_\lis(\Bun_G,\Lambda)$, the diagram
\[\xymatrix{
\mathcal Z^{\mathrm{spec}}(G,\Lambda)\ar[r]\ar[dr] & \mathrm{End}(A)\ar[d]\\
& \mathrm{End}(T_V(A))
}\]
commutes, so the map factors over $\mathcal Z^{\mathrm{geom}}_{\mathrm{Hecke}}(G,\Lambda)\subset \mathcal Z^{\mathrm{geom}}(G,\Lambda)$.
\end{theorem}

\begin{proof} This follows from the decomposition into connected components, the map $\mathrm{Exc}(W,\hat{G})_\Lambda\to \mathcal Z(\mathcal D_\lis^P(\Bun_G,\Lambda)^\omega)$ above, and Theorem~\ref{thm:coliminIndPerf}. The statement about commutation with Hecke operators follows from the construction of excursion operators and the commutation of Hecke operators.
\end{proof}

Before going on, we make the following observation regarding duality. The Bernstein--Zelevinsky duality functor $\mathbb D_{\mathrm{BZ}}$ on $\mathcal D_\lis(\Bun_G,\Lambda)$ induces an involution $D^{\mathrm{geom}}$ of $\mathcal Z^{\mathrm{geom}}(G,\Lambda)$. On the other hand, on $Z^1(W_E,\hat{G})$, the Chevalley involution of $\hat{G}$ induces an involution; after passing to the quotient by the conjugation action of $\hat{G}$, we can also forget about the inner automorphism appearing in Proposition~\ref{prop:chevalleyinvolution}. Let $D^{\mathrm{spec}}$ be the induced involution of $\mathcal Z^{\mathrm{spec}}(G,\Lambda)$.

\begin{proposition} The diagram
\[\xymatrix{
\mathcal Z^{\mathrm{spec}}(G,\Lambda)\ar[d]^{D^{\mathrm{spec}}}\ar[r] & \mathcal Z^{\mathrm{geom}}(G,\Lambda)\ar[d]^{D^{\mathrm{geom}}}\\
\mathcal Z^{\mathrm{spec}}(G,\Lambda)\ar[r] &\mathcal Z^{\mathrm{geom}}(G,\Lambda)
}\]
commutes.

The formation of $L$-parameters for irreducible smooth representations of $G(E)$ is compatible with passage to Bernstein--Zelevinsky duals, and to smooth duals.
\end{proposition}

\begin{proof} The commutation follows easily from the construction of excursion operators and Proposition~\ref{prop:chevalleyinvolution}. For the final part, it now follows that the formation of $L$-parameters is compatible with passage to Bernstein--Zelevinsky duals. For supercuspidal representations, this agrees with the smooth dual. In general, the claim for smooth duals follows from the compatibility with parabolic induction proved below.
\end{proof}

\section{Properties of the correspondence}\label{sec:propertiescorrespondence}

In this section, we check various basic properties of the correspondence. Throughout, we assume for simplicity that the order of $\pi_0 Z(G)$ is invertible in $\Lambda$. All results admit an obvious variant replacing the spectral Bernstein center by an excursion algebra when this assumption is omitted, and in particular the claims about $L$-parameters of Schur-irreducible objects work in any characteristic ($\neq p$, of course).

\subsection{Isogenies}

\begin{theorem}\label{thm:Lparamisogenies} Let $G'\to G$ be a map of reductive groups inducing an isomorphism of adjoint groups, inducing a dual map $\hat{G}\to \hat{G'}$, and $\pi: \Bun_{G'}\to \Bun_G$. Then for any $A\in D_\lis(\Bun_G,\Lambda)$ the diagram
\[\xymatrix{
\mathcal Z^{\mathrm{spec}}(\hat{G'},\Lambda)\ar[r]\ar[d] & \End(\pi^\ast A)\\
\mathcal Z^{\mathrm{spec}}(\hat{G},\Lambda)\ar[r] & \End(A)\ar[u]
}\]
commutes. In particular, if $\Lambda=L$ is an algebraically closed field, $A$ is Schur-irreducible and $A'$ is a Schur-irreducible constituent of $\pi^\ast A$, then $\varphi_{A'}$ is the composite of $\varphi_A$ with $\hat{G}\to \hat{G'}$.
\end{theorem}

\begin{proof} Consider any excursion operator for $G'$, given by some finite set $I$, a representation $V'\in \Rep_\Lambda((\hat{G'}\rtimes Q)^I)$, maps $\alpha: 1\to V'|_{\hat{G'}}$, $\beta: V'|_{\hat{G'}}\to 1$ and elements $\gamma_i\in \Gamma$ as usual. Consider the diagram
\[\xymatrix{
\Bun_{G'}\ar[d]^{\pi} & \Hglob_{G'}^I\ar[l]_{h_1'}\ar[d]^{\pi_H}\ar[r]^-{h_2'} & \Bun_{G'}\times (\Div^1)^I\ar[d]^{\pi}\\
\Bun_G & \Hglob_G^I\ar[l]_{h_1}\ar[r]^-{h_2} & \Bun_G\times (\Div^1)^I.
}\]
Then
\[
T_{V'}(\pi^\ast A)=h_{2\natural}'(h_1^{\prime\ast}\pi^\ast A\soliddotimesLambda \mathcal S_{V'}').
\]
We are interested in computing an endomorphism of $\pi^\ast A$; in particular, it is enough to compute $\pi_\natural T_{V'}(\pi^\ast A)$. But
\[\begin{aligned}
\pi_\natural T_{V'}(\pi^\ast A)&=\pi_\natural h_{2\natural}'(h_1^{\prime\ast}\pi^\ast A\soliddotimesLambda \mathcal S_{V'}')\\
&\cong h_{2\natural} \pi_{H\natural}(\pi_H^\ast h_1^\ast A\soliddotimesLambda \mathcal S_{V'}')\\
&\cong h_{2\natural}(h_1^\ast A\soliddotimesLambda \pi_{H\natural} \mathcal S_{V'}')\\
&\cong h_{2\natural}(h_1^\ast A\soliddotimesLambda h_1^\ast \pi_\natural \Lambda\soliddotimesLambda \mathcal S_V')=T_V(A\otimes \pi_\natural \Lambda).
\end{aligned}\]
This identification is functorial in $V'$ and $I$, and is over $\Bun_G\times (\Div^1)^I$, hence implies the desired equality of excursion operators. Here, to identify $\pi_{H\natural}\mathcal S_{V'}$, we write $\pi_H$ as the composite
\[
\Hglob_{G'}^I\to \Hglob_G^I\times_{\Bun_G} \Bun_{G'}\to \Hglob_G^I.
\]
The first map is locally (over $\Bun_{G'}$) isomorphic to the map $\Gr_{G'}^I\to \Gr_G^I$ and hence pushforward takes $\mathcal S_{V'}$ to the pullback of $\mathcal S_V$, by the compatibility of the geometric Satake equivalence with the map $G\to G'$ inducing isomorphisms of adjoint groups, as in the proof of Theorem~\ref{thm:geometricsatake}. Now the projection formula shows
\[
\pi_{H\natural} \mathcal S_{V'}'\cong h_1^\ast \pi_\natural \Lambda\soliddotimesLambda \mathcal S_V'.\qedhere
\]
\end{proof}

\subsection{Products}

\begin{proposition}\label{prop:Lparamproduct} If $G=G_1\times G_2$ is a product of two groups, then the diagram
\[\xymatrix{
\mathcal Z^{\mathrm{spec}}(G_1,\Lambda)\otimes_\Lambda \mathcal Z^{\mathrm{spec}}(G_2,\Lambda)\ar[d]^{\cong}\ar[r] & \mathcal Z^{\mathrm{geom}}(G_1,\Lambda)\otimes_\Lambda\mathcal Z^{\mathrm{geom}}(G_1,\Lambda)\ar[d]\\
\mathcal Z^{\mathrm{spec}}(G,\Lambda)\ar[r] & \mathcal Z^{\mathrm{geom}}(G,\Lambda)
}\]
commutes.

In particular, if $\Lambda=L$ is an algebraically closed field and $A_1,A_2\in \mathcal D_\lis(\Bun_G,L)$ are Schur-irreducible, and $A$ is a Schur-irreducible constituent of $A_1\boxtimes A_2$, then
\[
\varphi_A=(\varphi_{A_1},\varphi_{A_2}): W_E\to \hat{G}(L)\cong \hat{G_1}(L)\times \hat{G_2}(L).
\]
\end{proposition}

\begin{proof} The statement can be checked using excursion operators, and the proof is a straightforward diagram chase, noting that everything decomposes into products.
\end{proof}

\subsection{Weil restriction}

\begin{proposition}\label{prop:LparamWeilrestriction} If $G=\mathrm{Res}_{E'|E} G'$ is a Weil restriction of scalars of some reductive group $G'$ over some finite separable extension $E'$ of $E$. Choose $P$ to be an open subgroup of the wild inertia of $W_{E'}\subset W_E$, and let $W'\subset W_{E'}/P$ be the preimage of $W\subset W_E/P$. Then there are canonical identifications $\Bun_{G'}\cong \Bun_G$, $Z^1(W_E,\hat{G})/\hat{G}\cong Z^1(W_{E'},\hat{G'})/\hat{G'}$ and $\mathrm{Exc}(W,\hat{G})\cong \mathrm{Exc}(W',\hat{G'})$, and the diagram
\[\xymatrix{
\mathcal Z^{\mathrm{spec}}(G',\Lambda)\ar[d]^\cong\ar[r] & \mathcal Z^{\mathrm{geom}}(G',\Lambda)\ar[d]^\cong\\
\mathcal Z^{\mathrm{spec}}(G,\Lambda)\ar[r] & \mathcal Z^{\mathrm{geom}}(G,\Lambda)
}\]
commutes. In particular, $L$-parameters are compatible with Weil restriction.
\end{proposition}

\begin{proof} The most nontrivial of these identifications is the identification
\[
\mathrm{Exc}(W,\hat{G})\cong \mathrm{Exc}(W',\hat{G'}).
\]
One way to understand this is to use the presentation
\[
\mathrm{Exc}(W,\hat{G}) = \colim_{(n,F_n\to W)} \O(Z^1(F_n,\hat{G}))^{\hat{G}}
\]
(and the similar presentation for $\mathrm{Exc}(W',\hat{G'})$) and the natural isomorphism $Z^1(F_n,\hat{G})\sslash\hat{G}\cong Z^1(F_n\times_W W',\hat{G'})\sslash\hat{G'}$ of affine schemes (and then passing to global sections), noting that $F_n\times_W W'\subset F_n$ is a subgroup of finite index, and thus itself a finitely generated free group. This shows in fact that restricting to those maps $F_n\to W$ factoring over $W'$ produces the same colimit, and so
\[\begin{aligned}
\mathrm{Exc}(W,\hat{G}) &= \colim_{(n,F_n\to W)} \O(Z^1(F_n,\hat{G}))^{\hat{G}}\\
&\xleftarrow{\sim} \colim_{(n,F_n\to W')} \O(Z^1(F_n,\hat{G}))^{\hat{G}}\\
&\cong \colim_{(n,F_n\to W')} \O(Z^1(F_n,\hat{G'}))^{\hat{G'}}\\
&=\mathrm{Exc}(W',\hat{G'}).
\end{aligned}\]

Now consider an excursion operator for $G'$, including a representation $V'$ of $(\hat{G'}\rtimes W_{E'})^I$. Note that $\hat{G}\rtimes W_E$ contains $\hat{G}\rtimes W_{E'}$ as a subgroup, and this admits a surjection onto $\hat{G'}\rtimes W_{E'}$ (noting that $\hat{G}=\prod_{E'\hookrightarrow \overline{E}} \hat{G'}$, where we picked out an embedding $E'\hookrightarrow \overline{E}$ and hence a projection $\hat{G}\to \hat{G'}$ when we regarded $W_{E'}\subset W_E$ as a subgroup). In this way, one can inflate $V'$ to a representation of $(\hat{G}\rtimes W_{E'})^I$ and then induce to $(\hat{G}\rtimes W_E)^I$ to get a representation $V$ of $(\hat{G}\rtimes W_E)^I$. Geometrically, this procedure amounts to the commutative diagram
\[\xymatrix{
\Bun_{G'}\ar[d]^{\cong} & \Hglob_{G'}^I\ar[l]_{h_1'}\ar[r]^-{h_2'}\ar[d]^{\psi} & \Bun_{G'}\times (\Div^{\prime 1})^I\ar[d]\\
\Bun_G & \Hglob_G^I\ar[l]_{h_1}\ar[r]^-{h_2} & \Bun_G\times (\Div^1)^I
}\]
and taking $\psi_\ast$ on sheaves. More precisely, we note that
\[
\Hglob_{G'}^I\to \Hglob_G^I\times_{(\Div^1)^I}(\Div^{\prime 1})^I
\]
is a closed immersion (compatibly with a similar closed immersion of Beilinson--Drinfeld Grassmannians). Now the claim follows from a diagram chase.
\end{proof}

\subsection{Tori}

If $G=T$ is a torus, then
\[
\mathcal D_\lis(\Bun_T,\Lambda)\cong \prod_{b\in B(T)=\pi_1(T)_\Gamma} \mathcal D(T(E),\Lambda)
\]
and in particular
\[
\mathcal Z^{\mathrm{geom}}(T,\Lambda) = \prod_{b\in B(T)} \mathcal Z(T(E),\Lambda)
\]
where $\mathcal Z(T(E),\Lambda)$ is the Bernstein center of $T(E)$; explicitly, this is
\[
\mathcal Z(T(E),\Lambda) = \varprojlim_{K\subset T(E)} \Lambda[T(E)/K]
\]
where $K$ runs over open subgroups of $T(E)$.

\begin{proposition}\label{prop:stackLparamtori} There is a natural isomorphism
\[
\mathcal Z^{\mathrm{spec}}(T,\Lambda)\cong \varprojlim_{K\subset T(E)} \Lambda[T(E)/K].
\]
\end{proposition}

\begin{proof} One can resolve $T$ by products of induced tori and then reduce to the case that $T$ is induced, and then by Weil restrictions of scalars to $T=\mathbb G_m$. In that case $Z^1(W_E,\mathbb G_m)=\mathrm{Hom}(E^\ast,\mathbb G_m)$ by local class field theory, giving the result.
\end{proof}

\begin{proposition}\label{prop:Lparamtori} Under the above identifications
\[
\mathcal Z^{\mathrm{spec}}(T,\Lambda) = \varprojlim_{K\subset T(E)}\Lambda[T(E)/K]
\]
and
\[
\mathcal Z^{\mathrm{geom}}(T,\Lambda)= \prod_{b\in B(T)} \varprojlim_{K\subset T(E)} \Lambda[T(E)/K],
\]
the map
\[
\mathcal Z^{\mathrm{spec}}(T,\Lambda)\to \mathcal Z^{\mathrm{geom}}(T,\Lambda)
\]
is the diagonal embedding.
\end{proposition}

\begin{proof} We may resolve $T$ by induced tori and use Theorem~\ref{thm:Lparamisogenies}, Proposition~\ref{prop:Lparamproduct} and Proposition~\ref{prop:LparamWeilrestriction} to reduce to the case of $T=\mathbb G_m$. It is enough to compute the excursion operators corresponding to $I=\{1,2\}$, $V=\mathrm{std}\boxtimes \mathrm{std}^\vee$ and the tautological maps $\alpha: 1\to \mathrm{std}\otimes \mathrm{std}^\vee$ and $\beta: \mathrm{std}\otimes \mathrm{std}^\vee\to 1$. It is then an easy consequence of Section~\ref{sec:LubinTate}.
\end{proof}

Proposition~\ref{prop:Lparamtori} in particular shows that the $L$-parameters we construct for tori are the usual $L$-parameters, and together with Theorem~\ref{thm:Lparamisogenies} and Proposition~\ref{prop:Lparamproduct} implies that $L$-parameters are compatible with central characters (in case of connected center) and twisting, by applying Theorem~\ref{thm:Lparamisogenies} to the maps $Z\times G\to G$ and $G\to G\times D$ where $Z\subset G$ is the center and $G\to D$ is the quotient by the derived group. To deduce compatibility with central characters in general, one can reduce to the case of connected center using z-extensions \cite[Section 5]{KalethaZExt}.

\section{Applications to representations of $G(E)$}\label{sec:reprGE}

Finally, we apply the preceding results to representations of $G(E)$. We get the following map to the Bernstein center.

\begin{definition}\label{def:bernsteinmap} The map
\[
\Psi_G: \mathcal Z^{\mathrm{spec}}(G,\Lambda)\to \mathcal Z(G(E),\Lambda)
\]
is the composite
\[
\mathcal Z^{\mathrm{spec}}(G,\Lambda)\to \mathcal Z^{\mathrm{geom}}(G,\Lambda)\to \mathcal Z(G(E),\Lambda)
\]
induced by the fully faithful functor
\[
j_!: \mathcal D(G(E),\Lambda)\cong \mathcal D_\lis(\Bun_G^1,\Lambda)\to \mathcal D_\lis(\Bun_G,\Lambda).
\]
\end{definition}

More generally, for any $b\in B(G)$, we can define a map
\[
\Psi_G^b: \mathcal Z^{\mathrm{spec}}(G,\Lambda)\to \mathcal Z(G_b(E),\Lambda)
\]
to the Bernstein center for $G_b(E)$ by using the fully faithful embedding
\[
\mathcal D(G_b(E),\Lambda)\cong \mathcal D_\lis(\Bun_G^b,\Lambda)\to \mathcal D_\lis(\Bun_G,\Lambda)
\]
determined for example by the left adjoint to $i^{b\ast}$, where $i^b: \Bun_G^b\hookrightarrow \Bun_G$ is the locally closed embedding (see Proposition~\ref{prop:compactobjectsDlis}). (Recall that in the $D_\lis$-setting, we do not have a general $i^b_!$-functor, although it can be defined in the present situation. All these maps will induce the same map to the Bernstein center.)

\subsection{Compatibility with $G_b$}

One can describe the maps $\Psi_G^b$ for $b\neq 1$ in terms of the maps $\Psi_{G_b}$. Note that $\hat{G_b}$ is naturally a Levi subgroup of $\hat{G}$, as $G_{b,\breve{E}}\subset G_{\breve{E}}$ is the centralizer of the slope morphism $\nu_b: \mathbb D\to G_{\breve{E}}$. This extends naturally to a morphism of $L$-groups
\[
\hat{G_b}\rtimes Q\to \hat{G}\rtimes Q
\]
where as usual $Q$ is a finite quotient of $W_E$ over which the action on $\hat{G}$ factors. However, from geometric Satake we rather get the natural inclusion
\[
\widecheck{G_b}\rtimes W_E\to \widecheck{G}\rtimes W_E
\]
where the $W_E$-actions include the cyclotomic twist. The latter induces a map
\[
Z^1(W_E,\hat{G_b})\to Z^1(W_E,\hat{G})
\]
that in terms of the usual $W_E$-action is given by sending a $1$-cocycle $\varphi: W_E\to \hat{G_b}(A)$ to the $1$-cocycle
\[
W_E\to \hat{G}(A): w\mapsto (2\rho_{\hat{G}}-2\rho_{\hat{G_b}})(\sqrt{q})^{|w|}\varphi(w)
\]
where $|\cdot|: W_E\to W_E/I_E\cong \mathbb Z$ is normalized as usual by sending a geometric Frobenius to $1$.

\begin{theorem}\label{thm:generalGbcenter} For all $G$ and $b\in B(G)$, the diagram
\[\xymatrix{
\mathcal Z^{\mathrm{spec}}(G,\Lambda)\ar[r]^{\Psi_G^b}\ar[d] & \mathcal Z(\mathcal D(G_b(E),\Lambda))\\
\mathcal Z^{\mathrm{spec}}(G_b,\Lambda)\ar[ur]_{\Psi_{G_b}}
}\]
commutes.
\end{theorem}

\begin{proof} We note that to prove the theorem, we can assume that $\Lambda$ is killed by power of $\ell$ (if $\ell$ divides the order of $\pi_0 Z(G)$, replacing the left-hand side with an algebra of excursion operators), as the result for $\Lambda=\mathbb Z_\ell[\sqrt{q}]$ implies it in general, and the right-hand side
\[
\mathcal Z(\mathcal D(G_b(E),\Lambda))=\varprojlim_{K\subset G_b(E)} \mathcal Z(\Lambda[G_b(E)\sslash K])
\]
is $\ell$-adically separated in that case. This means we can avoid the subtleties of $\mathcal D_\lis$ in place of $\mathcal D_\et$.

If $b$ is basic, the theorem follows from the identification $\Bun_G\cong \Bun_{G_b}$ of Corollary~\ref{cor:basictwistBunG}, which is equivariant for the Hecke action.

In general, we first reduce to the case that $G$ is quasisplit. Take a z-embedding $G\hookrightarrow G'$ as in \cite[Section 5]{KalethaZExt}, with quotient a torus $D$, so that the center $Z(G')$ is connected. Then $\Bun_G = \Bun_{G'}\times_{\Bun_D} \{\ast\}$ and the map $B(G)\to B(G')$ is injective. To see the latter, by the description of the stacks, it suffices to see that for all $b\in B(G)$ with image $b'\in B(G')$, the map $G'_{b'}(E)\to D(E)$ is surjective. But for any $b\in B(G)$, the map $G_b\to G'_{b'}$ is a z-embedding with quotient $D$, and $Z'(E)\to D(E)$ is surjective by \cite[Fact 5.5]{KalethaZExt}, where also $Z'\subset G'_{b'}$, so in particular $G'_{b'}(E)\to D(E)$ is surjective. An element of
\[
\mathcal Z(\mathcal D(G_b(E),\Lambda))=\varprojlim_{K\subset G_b(E)} \mathcal Z(\Lambda[G_b(E)\sslash K])
\]
of the Bernstein center of $G_b(E)$ is determined by its action on $\pi'|_{G_b(E)}$ for representations $\pi'$ of $G'_{b'}(E)$. By Theorem~\ref{thm:Lparamisogenies}, we can thus reduce to $G'$ in place of $G$, i.e.~that the center of $G$ is connected. When $Z(G)$ is connected, there is some basic $b_0\in B(G)$ such that $G_{b_0}$ is quasisplit. Using the Hecke-equivariant isomorphism $\Bun_G\cong \Bun_{G_{b_0}}$ we can thus assume that $G$ is quasisplit.

Now if $G$ is quasisplit, fix a Borel $B\subset G$. Any $b\in B(G)$ then admits a reduction to a canonical parabolic $P=P_b\subset G$ containing $B$. Pick a cocharacter $\mu: \mathbb G_m\to G$ with dynamical parabolic $P$. For any $N\geq 0$, let $b_N = b\mu(\pi^N)$: This is a sequence of elements of $B(G)$ associated to the same parabolic $P$ but increasingly instable. Moreover, $G_b=G_{b_N}$. We note that the diagram
\[\xymatrix{
\mathcal Z^{\mathrm{spec}}(G,\Lambda)\ar[r]^{\Psi_G^b}\ar[d]^{=} & \mathcal Z(G_b(E),\Lambda)\ar[d]^{=}\\
\mathcal Z^{\mathrm{spec}}(G,\Lambda)\ar[r]^{\Psi_G^{b_N}}& \mathcal Z(G_{b_N}(E),\Lambda)
}\]
commutes. For this, take any representation $\sigma$ of $G_b(E)$ and consider the sheaf $A_N\in \mathcal D_\et(\Bun_G,\Lambda)$ concentrated on $\Bun_G^{b_N}$, corresponding to the representation $\sigma$. Let $V\in \Rep \hat{G}$ be the highest weight representation with weight $\mu^N$. We claim that $T_V(A_N)|_{\Bun_G^b}$ is given by the representation $\sigma$. As Hecke operators commute with excursion operators, this implies the desired result. To compute $T_V(A_N)|_{\Bun_G^b}$, we have to analyze the moduli space of modifications of $\mathcal E_b$ of type bounded by $\mu^N$ that are isomorphic to $\mathcal E_{b_N}$. There is in fact precisely one such modification, given by pushout of the standard modification of line bundles from $\mathcal O$ to $\mathcal O(1)$ via $\mu^N: \mathbb G_m\to G$; its type is exactly $\mu^N$. This gives the claim.

Now to prove the theorem, we have to prove the commutativity of the diagram for any excursion operator, given by excursion data $(I,V,\alpha,\beta,(\gamma_i)_{i\in I})$. For any such excursion data, we can pick $N$ large enough so that any modification of $\mathcal E_{b_N}$ to itself, of type bounded by $V$, is automatically compatible with the Harder--Narasimhan reduction to $P$. In that case, for $\sigma$ and $A_N$ as above, to analyze the excursion operators
\[
A_N=T_1(A_N)\xrightarrow{\alpha} T_V(A_N)\xrightarrow{(\gamma_i)_{i\in I}} T_V(A_N)\xrightarrow{\beta} T_1(A_N)=A_N,
\]
we have to analyze the moduli space of modifications of $\mathcal E_{b_N}$, at $I$ varying points, of type bounded by $V$, and that are isomorphic to $\mathcal E_{b_N}$. By assumption on $N$, this is the same as the moduli space of such modifications as $P$-bundles. This maps to the similar moduli space parametrizing modifications as $M$-bundles, where $M$ is the Levi of $P$ (and $G_b=M_{b_M}$ for a basic $b_M\in B(M)$). We want to compute $T_V(A_N)$. Note that $A_N$ comes from $A_N'\in D_\et(\Bun_P,\Lambda)$ as it is concentrated on $\Bun_G^{b_N}\cong \Bun_P^{b_N}\subset \Bun_P$, by the Harder--Narasimhan reduction. 

Consider the diagram
\[\xymatrix{
\Bun_M\ar[d]^\psi & \Hglob_{M,P}^I\ar[d]^{\psi_H}\ar[l]_{h_1''}\ar[r]^-{h_2''} & \Bun_P\times (\Div^1)^I\ar[d]^{=}\\
\Bun_{P}\ar[d]^{\pi} & \Hglob_{P}^I\ar[l]_{h_1'}\ar[d]^{\pi_H}\ar[r]^-{h_2'} & \Bun_{P}\times (\Div^1)^I\ar[d]^{\pi}\\
\Bun_{M} & \Hglob_{M}^I\ar[l]_{h_1}\ar[r]^-{h_2} & \Bun_{M}\times (\Div^1)^I
}\]
where $\Hglob_{M,P}^I$ is defined as the fibre product $\Hglob_P^I\times_{\Bun_P} \Bun_M$, and thus parametrizes modifications from an $M$-bundle to a $P$-bundle.

We need to compute $T_V(A_N)|_{\Bun_G^{b_N}}$, which by the above argument that any modification of $\mathcal E_{b_N}$ to itself of type bounded by $V$ is a modification as $P$-bundles, can be computed in terms of the middle diagram, as
\[
Rh_{2!}'(h_1^{\prime\ast} A_N'\dotimes_\Lambda \mathcal S_V)|_{\Bun_P^{b_N}}
\]
where $\mathcal S_V\in \Hloc_G^I$ is the perverse sheaf determined by $V$ under the geometric Satake equivalence, and we continue to denote by $\mathcal S_V$ any of its pullbacks.

There is some $B_N\in D_\et(\Bun_M,\Lambda)$ such that $A_N'=R\psi_! B_N$. In fact, one can take $B_N=R\pi_! A_N'$, noting that on the support of $A_N'$, the map $\pi: \Bun_P\to \Bun_M$ is (cohomologically) smooth, so $R\pi_!$ is defined on $A_N'$ (although $\pi$ is a stacky map). (Indeed, everything is concentrated on one stratum, and the relevant categories are all equivalent to $D(G_b(E),\Lambda)$.) Moreover, to compute the restriction to $\Bun_P^{b_N}$ it is enough to do the computation after applying $R\pi_!$. We compute:
\[\begin{aligned}
R\pi_!Rh_{2!}'(h_1^{\prime\ast} A_N'\dotimes_\Lambda \mathcal S_V)&= R\pi_!Rh_{2!}'R\psi_{\mathcal H!}(h_1^{\prime\prime\ast} B_N\dotimes_\Lambda \mathcal S_V)\\
&=Rh_{2!}(h_1^\ast B_N\dotimes_\Lambda Rg_!\mathcal S_V)
\end{aligned}\]
where $g: \Hglob_{M,P}^I\to \Hglob_M^I$ is the projection. But this is the pullback of the map $L^+M\backslash \Gr_P^I\to L^+M\backslash \Gr_M^I=\Hloc_M^I$ under $\Hglob_M^I\to \Hloc_M^I$. This means that $Rg_!\mathcal S_V$ arises via pullback from $\mathrm{CT}_P(\mathcal S_V)\in D_\et(\Hloc_M^I,\Lambda)$. Up to the shift $[\mathrm{deg}_P]$, this agrees with $\mathcal S_{V|_{(\hat{M}\rtimes Q)^I}}$, where the restriction involves a cyclotomic twist, as above. (It is the canonical restriction along $\widecheck{M}\to \widecheck{G}$ for the canonical $W_E$-actions arising geometrically.) Now the excursion operators, which involve maps from and to the sheaf corresponding to $V=1$, require only the connected component where $\mathrm{deg}_P=0$, so we can ignore the shift.

With these translations, we see that the excursion operators on $\Bun_G^{b_N}$ and on $\Bun_M^{b_N}$ agree, giving the desired result.
\end{proof}

\subsection{Parabolic induction}

A corollary of this result is compatibility with parabolic induction.

\begin{corollary}\label{cor:Lparamparabinduction} Let $G$ be a reductive group with a parabolic $P\subset G$ and Levi $P\to M$. Then for all representations $\sigma$ of $M(E)$ with (unnormalized) parabolic induction $\mathrm{Ind}_{P(E)}^{G(E)} \sigma$, the diagram
\[\xymatrix{
\mathcal Z^{\mathrm{spec}}(G,\Lambda)\ar[r]\ar[d] & \mathrm{End}(\mathrm{Ind}_{P(E)}^{G(E)} \sigma)\\
\mathcal Z^{\mathrm{spec}}(M,\Lambda)\ar[r] & \mathrm{End}(\sigma)\ar[u]
}\]
commutes. In particular, the formation of $L$-parameters is compatible with parabolic induction: If $\Lambda=L$ is an algebraically closed field, $\sigma$ is irreducible and $\tilde{\sigma}$ is an irreducible subquotient of $\mathrm{Ind}_{P(E)}^{G(E)} \sigma$, then $\varphi_{\tilde{\sigma}}$ is conjugate to the composite
\[
W_E\xrightarrow{\varphi_\sigma} \hat{M}(L)\rtimes W_E\to \hat{G}(L)\rtimes W_E
\]
where the map $\hat{M}\rtimes W_E\to \hat{G}\rtimes W_E$ is defined as above, involving the cyclotomic twist.
\end{corollary}

\begin{proof} It suffices to prove the result for $\sigma=c\text-\mathrm{Ind}_K^{M(E)} \Lambda$ for $K\subset M(E)$ an open pro-$p$-subgroup, and then one can assume $\Lambda=\mathbb Z_\ell[\sqrt{q}]$, where one can further by $\ell$-adic separatedness reduce to torsion coefficients.

Let $\mu: \mathbb G_m\to G$ be a cocharacter with dynamical parabolic $P$ and let $b=\mu(\pi^{-1})\in B(G)$. Then $G_b=M$, and we can build a sheaf $A\in \mathcal D_\et(\Bun_G,\Lambda)$ concentrated on $\Bun_G^b$, given by the representation $\sigma$. Then $T_{\mu^{-1}}(A)|_{\Bun_G^1}$ is given by a parabolic induction of $\sigma$, more precisely $\mathrm{Ind}_{P(E)}^{G(E)}\sigma(-\tfrac d2)[-d]$ where $d=\langle 2\rho,\mu\rangle$: To see this, we have to understand the moduli space of modifications of the trivial $G$-torsor of type bounded by $\mu$ that are isomorphic to $\mathcal E_b$. This is in fact given by $\underline{G(E)/P(E)}$, the $G(E)$-orbit of the pushout of the modification from $\mathcal O$ to $\mathcal O(1)$ via $\mu$. All of these modifications are of type exactly $\mu$ (so the Satake sheaf is simply a twist $(-\tfrac d2)[-d]$ of the constant sheaf). This easily gives the claim on $T_{\mu^{-1}}(A)|_{\Bun_G^1}$.\footnote{See \cite[Theorem 4.26]{GaisinImai} for more details on the analysis of Hodge--Newton reducible local Shimura varieties.}
Now as Hecke operators commute with excursion operators, the excursion operators on $\mathrm{Ind}_{P(E)}^{G(E)}\sigma$ agree with those on $A$, and these are determined by Theorem~\ref{thm:generalGbcenter}, giving the result.
\end{proof}

\subsection{The case $G=\mathrm{GL}_n$}

For the group $G=\GL_n$, we can identify the $L$-parameters with the usual $L$-parameters of \cite{LaumonRapoportStuhler}, \cite{HarrisTaylor}, \cite{HenniartGLn}. This is the only place of this paper where we rely on previous work on the local Langlands correspondence, or (implicitly) rely on global arguments. More precisely, we use the identification of the cohomology of the Lubin--Tate and Drinfeld tower, see \cite{BoyerMauvais}, \cite{Harris1}, \cite{HarrisTaylor}, \cite{Hausberger}, \cite{DatElliptic}. In the proof, we use the translation between Hecke operators and local Shimura varieties as Section~\ref{sec:localShimura}, together with the description of these as the Lubin--Tate tower and Drinfeld tower in special cases, see \cite{ScholzeWeinstein}.

\begin{theorem}\label{thm:LparamGLn} Let $\pi$ be any irreducible smooth $\overline{\mathbb Q}_\ell$-representation of $\GL_n(E)$. Then the $L$-parameter $\varphi_\pi$ agrees with the usual (semisimplified) $L$-parameter.
\end{theorem}

\begin{proof} By Corollary~\ref{cor:Lparamparabinduction}, we can assume that $\pi$ is supercuspidal. We only need to evaluate the excursion operators for the excursion data given by $I=\{1,2\}$, the representation $V=\mathrm{std}\boxtimes \mathrm{std}^\vee$ of $\widehat{\GL_n}^2$, and the unit/counit maps $\alpha: 1\to \mathrm{std}\otimes \mathrm{std}^\vee$ and $\beta: \mathrm{std}\otimes \mathrm{std}^\vee\to 1$, as these excursion operators determine the trace of the representation (and thus the semisimplified representation).

First, we analyze these excursion operators on the sheaf $B$ which is the sheaf on $\Bun_{\GL_n}^b$ for $b$ corresponding to the bundle $\mathcal O(-\tfrac 1n)$, given by the representation $\sigma=\mathrm{JL}(\pi)$ of $D^\times$; here $D$ is the division algebra of invariant $\tfrac 1n$. The Hecke operator $T_V$ is the composite of two operators. The first Hecke operator, corresponding to $\mathrm{std}$, takes minuscule modifications $\mathcal O(-\tfrac 1n)\subset \mathcal E$ with cokernel a skyscraper sheaf of rank $1$. Such an $\mathcal E$ is necessarily isomorphic to $\mathcal O^n$, and the Hecke operator will then produce the $\sigma$-isotypic part of the cohomology of the Lubin--Tate tower, which is $\pi\otimes \rho_\pi$, where $\rho_\pi$ is the irreducible $n$-dimensional $W_E$-representation associated to $\pi$ by the local Langlands correspondence. (Note that the shift $[n-1]$, as well as the cyclotomic twist $(\frac{n-1}2)$ that usually appears, is hidden inside the normalization of the perverse sheaf corresponding to the standard representation.) Now the second Hecke operator, when restricted to $\Bun_{\GL_n}^b$, produces the $\pi$-isotypic component of the cohomology of the Drinfeld tower, which is $\sigma\otimes \rho_\pi^\ast$. In total, we see that $T_V(B)|_{\Bun_{\GL_n}^b}$ is given by $\sigma\otimes \rho_\pi\otimes \rho_\pi^\ast$ as representation of $D^\times\times W_E\times W_E$. By irreducibility of $\rho_\pi$, the $W_E$-equivariant map $\sigma\to \sigma\otimes \rho_\pi\otimes \rho_\pi^\ast$ induced by $\alpha$ (and the similar map backwards induced by $\beta$) must agree up to scalar with the obvious map. The scalar of the total composite can be identified by taking both elements of $W_E$ to be equal to $1$. This shows that $B$ has the correct $L$-parameter. Now use that the sheaf corresponding to $\pi$ appears as a summand of $T_{\mathrm{std}}(B)$ (after forgetting the $W_E$-action) to conclude the same for $\pi$.
\end{proof}

In particular, it follows that the map $\mathcal Z^{\mathrm{spec}}(\GL_n,\overline{\mathbb Q}_\ell)\to \mathcal{Z}(\GL_n(E),\overline{\mathbb Q}_\ell)$ to the Bernstein center agrees with the usual map. But this refines to a map
\[
\mathcal Z^{\mathrm{spec}}(\GL_n,\mathbb Z_\ell[\sqrt{q}])\to \mathcal{Z}(\GL_n(E),\mathbb Z_\ell[\sqrt{q}])
\]
to the integral Bernstein center, recovering a result of Helm--Moss \cite{HelmMoss}.

\chapter{The spectral action}\label{ch:spectral}

As a final topic, we construct the spectral action. We will first construct it with characteristic $0$ coefficients, and then explain refinements with integral coefficients.

Let $\Lambda$ be the ring of integers in a finite extension of $\mathbb Q_\ell(\sqrt{q})$. We have the stable $\infty$-category $\mathcal C=\mathcal D_\lis(\Bun_G,\Lambda)^\omega$ of compact objects, which is linear over $\Lambda$, and functorially in the finite set $I$ an exact monoidal functor $\Rep_\Lambda(\hat{G}\rtimes Q)^I\to \End_\Lambda(\mathcal C)^{BW_E^I}$ that is linear over $\Rep_\Lambda(Q^I)$. A first version of the following theorem is due to Nadler--Yun \cite{NadlerYun} in the context of Betti geometric Langlands, and a more general version appeared in the work of Gaitsgory--Kazhdan--Rozenblyum--Varshavsky \cite{GaitsgoryKazhdanRozenblyumVarshavsky}. Both references, however, effectively assume that $G$ is split, work only with characteristic $0$ coefficients, and work with a discrete group in place of $W_E$. At least the extension to $\mathbb Z_\ell$-coefficients is a nontrivial matter.

Note that $Z^1(W_E,\hat{G})$ is not quasicompact, as it has infinitely many connected components; it can be written as the increasing union of open and closed quasicompact subschemes $Z^1(W_E/P,\hat{G})$. We say that an action of $\Perf(Z^1(W_E,\hat{G})/\hat{G})$ on a stable $\infty$-category $\mathcal C$ is compactly supported if for all $X\in \mathcal C$ the functor $\Perf(Z^1(W_E,\hat{G})/\hat{G})\to \mathcal C$ (induced by acting on $X$) factors over some $\Perf(Z^1(W_E/P,\hat{G})/\hat{G})$.

The goal of this chapter is to prove the following theorem. Below, ``functorially in the finite set $I$'' means a map on total spaces over $\mathrm{Fin}$ of the corresponding coCartesian fibrations.

\begin{theorem}\label{thm:constructionspectralaction} Assume that $\ell$ does not divide the order of $\pi_1(\hat{G})_{\mathrm{tor}}$. Let $\mathcal C$ be a small idempotent-complete $\Lambda$-linear stable $\infty$-category. Then giving, functorially in the finite set $I$, an exact $\Rep_\Lambda(Q^I)$-linear monoidal functor
\[
\Rep_\Lambda(\hat{G}\rtimes Q)^I\to \End_\Lambda(\mathcal C)^{BW_E^I}
\]
is equivalent to giving a compactly supported $\Lambda$-linear action of
\[
\Perf(Z^1(W_E,\hat{G})_\Lambda/\hat{G}).
\]
Here, given a $\Lambda$-linear action of $\Perf(Z^1(W_E,\hat{G})_\Lambda/\hat{G})$, one can produce such an exact $\Rep_\Lambda(Q^I)$-linear monoidal functor
\[
\Rep_\Lambda(\hat{G}\rtimes Q)^I\to \End_\Lambda(\mathcal C)^{BW_E^I}
\]
functorially in $I$ by composing the exact $\Rep_\Lambda(Q^I)$-linear symmetric monoidal functor
\[
\Rep_\Lambda(\hat{G}\rtimes Q)^I\to \Perf(Z^1(W_E,\hat{G})_\Lambda/\hat{G})^{BW_E^I}
\]
with the action of $\Perf(Z^1(W_E,\hat{G})_\Lambda/\hat{G})$.

The same result holds true if $\Lambda$ is a field over $\mathbb Q_\ell(\sqrt{q})$, for any prime $\ell$.
\end{theorem}

Here, the exact $\Rep_\Lambda(Q^I)$-linear symmetric monoidal functor
\[
\Rep_\Lambda(\hat{G}\rtimes Q)^I\to \Perf(Z^1(W_E,\hat{G})_\Lambda/\hat{G})^{BW_E^I}
\]
is induced by tensor products and the exact $\Rep_\Lambda(Q)$-linear symmetric monoidal functor
\[
\Rep_\Lambda(\hat{G}\rtimes Q)\to \Perf(Z^1(W_E,\hat{G})_\Lambda/\hat{G})^{BW_E}
\]
corresponding to the universal $\hat{G}\rtimes Q$-torsor, with the universal $W_E$-equivariance as parametrized by $Z^1(W_E,\hat{G})/\hat{G}$.

Before starting the proof, we note that the proof of Proposition~\ref{prop:actionofWEreasonable} shows that we may replace $W_E$ by $W_E/P$ in the statement of Theorem~\ref{thm:constructionspectralaction}. Choosing moreover a discretization $W\subset W_E/P$, we reduce to the following variant.

\begin{theorem}\label{thm:constructionspectralaction2} Assume that $\ell$ does not divide the order of $\pi_1(\hat{G})_{\mathrm{tor}}$. Let $\mathcal C$ be a small idempotent-complete $\Lambda$-linear stable $\infty$-category. Then giving, functorially in the finite set $I$, an exact $\Rep_\Lambda(Q^I)$-linear monoidal functor
\[
\Rep_\Lambda(\hat{G}\rtimes Q)^I\to \End_\Lambda(\mathcal C)^{BW^I}
\]
is equivalent to giving a $\Lambda$-linear action of
\[
\Perf(Z^1(W,\hat{G})_\Lambda/\hat{G}),
\]
with the same compatibility as above. The same result holds true if $\Lambda$ is a field over $\mathbb Q_\ell$, for any prime $\ell$.
\end{theorem}

\section{Rational coefficients}

With rational coefficients, we can prove a much more general result, following \cite{GaitsgoryKazhdanRozenblyumVarshavsky}. Consider a reductive group $H$ over a field $L$ of characteristic $0$ (like $\hat{G}$ over $\overline{\mathbb Q}_\ell$) with an action of a finite group $Q$. Let $S$ be any anima over $\ast/Q$ (like $\ast/W$, where $W\subset W_E/P$ is a discretization of $W_E/P$ for an open subgroup of the wild inertia, as usual). We can then consider the (derived) stack $\Map_{\ast/Q}(S,\ast/(H\rtimes Q))$ over $L$, whose values in an animated $L$-algebra $A$ are the maps of anima $S\to \ast/(H\rtimes Q)(A)$ over $\ast/Q$. This recovers the stack $[Z^1(W,\hat{G})_{\overline{\mathbb Q}_\ell}/\hat{G}]$ in the above example, using Proposition~\ref{prop:colimunderlying}.

In general, $\Map_{\ast/Q}(S,\ast/(H\rtimes Q))$ is the fpqc quotient of an affine derived scheme by a power of $H$. Indeed, pick a surjection $S'\to S\times_{\ast/Q}\ast$ from a set $S'$. Then $\Map_{\ast/Q}(S,\ast/(H\rtimes Q))$ maps to $\ast/H^{S'}$; we claim that the fibre is an affine derived scheme, i.e.~representable by an animated $L$-algebra. For this, note that
\[
\Map_{\ast/Q}(S,\ast/(H\rtimes Q))\to \Map(S\times_{\ast/Q}\ast,\ast/H)
\]
is relatively representable, as it is given by the $Q$-fixed points. To show that the right-hand side is relatively representable over $\ast/H^{S'}$, we can replace $S\times_{\ast/Q}\ast$ by a connected anima $T$, and $S'$ by a point. Then $\Map(T,\ast/H)\times_{\ast/H} \ast$ parametrizes pointed maps $T\to \ast/H$, which are equivalent to maps of $\mathbb E_1$-groups $\Omega(T)\to H$ , cf.~\cite[Lemma 7.2.2.11 (1)]{LurieHTT}. Writing $\Omega(T)$ as a sifted colimit of finite free groups $F_n$, one reduces to representability of maps of groups $F_n\to H$, which is representable by $H^n$.

\begin{theorem}\label{thm:spectralactionchar0} Let $\mathcal C$ be an idempotent-complete small stable $L$-linear $\infty$-category. Giving, functorially in finite sets $I$, an exact $\Rep_L(Q^I)$-linear monoidal functor
\[
\Rep_L((H\rtimes Q)^I)\to \End_L(\mathcal C)^{S^I}
\]
is equivalent to giving an $L$-linear action of $\Perf(\Map_{\ast/Q}(S,\ast/(H\rtimes Q)))$ on $\mathcal C$. Here, given such an action of $\Perf(\Map_{\ast/Q}(S,\ast/(H\rtimes Q)))$, one gets exact $\Rep_L(Q^I)$-linear monoidal functors
\[
\Rep_L((H\rtimes Q)^I)\to \End_L(\mathcal C)^{S^I}
\]
by precomposing the exact monoidal functor $\Perf(\Map_{\ast/Q}(S,\ast/(H\rtimes Q)))\to \End_L(\mathcal C)$ with the natural exact $\Rep_L(Q^I)$-linear symmetric monoidal functor
\[
\Rep_L((H\rtimes Q)^I)\to \Perf(\Map_{\ast/Q}(S,\ast/(H\rtimes Q)))^{S^I}
\]
given by $I$-fold tensor product of the exact $\Rep_L(Q)$-linear symmetric monoidal functor
\[
\Rep_L(H\rtimes Q)\to \Perf(\Map_{\ast/Q}(S,\ast/(H\rtimes Q)))^S
\]
assigning to each $s\in S$ pullback along evaluation at $s$, $\Map_{\ast/Q}(S,\ast/(H\rtimes Q))\to \ast/(H\rtimes Q)$.
\end{theorem}

\begin{proof} Note first that, for any $L$-linear idempotent-complete small stable $\infty$-category $\mathcal C$, giving an exact $L$-linear functor $\Rep_L((H\rtimes Q)^I)\to \mathcal C$ is equivalent to giving an exact $L$-linear functor of stable $\infty$-categories $\Perf((H\rtimes Q)^I)\to \mathcal C$, as the $\infty$-category of perfect complexes is freely generated by the exact category of representations. Indeed, such functors extend to the $\infty$-category obtained by inverting quasi-isomorphisms in $\mathrm{Ch}^b(\Rep_L(H\rtimes Q)^I)$, and this is $\Perf(\ast/(H\rtimes Q)^I)$.

For any $S$, we have the anima $F_1(S)$ of $L$-linear actions of $\Perf(\Map_{\ast/Q}(S,\ast/(H\rtimes Q)))$ on $\mathcal C$,\footnote{Indeed, any $\infty$-categorical datum will naturally produce an anima, or $\infty$-groupoid, via incorporating all isomorphisms and higher isomorphisms.} and the anima $F_2(S)$ of functorial exact monoidal functors
\[
\Rep_L((H\rtimes Q)^I)\to \End_L(\mathcal C)^{S^I}
\]
linear over $\Rep_L(Q^I)$, and a natural map $F_1(S)\to F_2(S)$ functorial in $S$ (where both $F_1$ and $F_2$ are contravariant functors of $S$). Both functors take sifted colimits in $S$ to limits. This is clear for $F_2$ (as $S\mapsto S^I$ commutes with sifted colimits). For $F_1$, it is enough to see that taking $S$ to $\Perf(\Map_{\ast/Q}(S,\ast/(H\rtimes Q)))$ commutes with sifted colimits (taken in idempotent-complete stable $\infty$-categories), which is Lemma~\ref{lem:perfectcomplexeslimit} below.

Therefore it suffices to handle the case that $S$ is a finite set, for which the map $S\to \ast/Q$ can be factored over $\ast$. Then $\Map_{\ast/Q}(S,\ast/(H\rtimes Q))\cong \ast/H^S$. Similarly, exact monoidal functors
\[
\Rep_L((H\rtimes Q)^I)\to \End_L(\mathcal C)^{S^I}
\]
linear over $\Rep_L(Q^I)$ are equivalent to exact monoidal functors
\[
\Rep_L(H^I)\to \End_L(\mathcal C)^{S^I}
\]
linear over $L$. Here, we use $\Perf(\ast/(H\rtimes Q)^I)\otimes_{\Perf(\ast/Q^I)} \Perf(L)\cong \Perf(\ast/H^I)$, which follows easily from highest weight theory.

The latter data is equivalent to maps
\[
\Hom(I,S)=S^I\to \Fun^{\mathrm{mon}}_{\mathrm{ex},L}(\Rep_L(H^I),\End_L(\mathcal C))
\]
functorially in $I$, where $\Fun^{\mathrm{mon}}_{\mathrm{ex},L}$ denotes the exact $L$-linear monoidal functors. Both sides here are functors in $I$, and on the left-hand side we have a representable functor. By the Yoneda lemma, it follows that this data is equivalent to $L$-linear exact monoidal functors
\[
\Rep_L(H^S)\to \End_L(\mathcal C).
\]
Such actions extend uniquely to $\Perf(\ast/H^S)$ by the observation explained in the first paragraph of this proof, giving the desired result.
\end{proof}

\begin{lemma}\label{lem:perfectcomplexeslimit} The functor taking an anima $S$ over $\ast/Q$ to $\Perf(\Map_{\ast/Q}(S,\ast/(H\rtimes Q)))$, regarded as an idempotent-complete stable $\infty$-category, commutes with sifted colimits. More precisely, as a functor into $L$-linear symmetric monoidal idempotent-complete stable $\infty$-categories, it commutes with all colimits.
\end{lemma}

We use here \cite[Corollary 3.2.3.2]{LurieHA} to see that sifted colimits agree with or without symmetric monoidal structure.

\begin{proof} We first check that it commutes with filtered colimits. For this, let $S_i$, $i\in I$, be a filtered diagram of anima over $\ast/Q$, and choose compatible surjections $S_i'\to S_i\times_{\ast/Q}\ast$ from sets $S_i'$. Let $S=\colim_i S_i$ and $S'=\colim_i S_i'$, which is a set surjecting onto $S\times_{\ast/Q}\ast$. Letting $G_i=H^{S_i'}$ and $G=H^{S'}$, we get presentations $\Map_{\ast/Q}(S_i,\ast/(H\rtimes Q))=X_i/G_i$ as quotients of affine derived $L$-schemes $X_i$ by the pro-reductive group $G_i$, and similarly $\Map_{\ast/Q}(S,\ast/(H\rtimes Q))=X/G$, with $X=\varprojlim_i X_i$ and $G=\varprojlim_i G_i$. We claim that in this generality
\[
\varinjlim_i \Perf(X_i/G_i)\to \Perf(X/G)
\]
is an isomorphism of idempotent-complete stable $\infty$-categories.

Assume first that all $X=\Spec L$ are a point. Then note that $\Perf(\ast/G)$ is generated by $\Rep(G)$, which is easily seen to be the filtered colimit $\varinjlim_i \Rep(G_i)$, and (by writing it as limit of reductive groups) is seen to be semisimple. The claim is easily checked in this case.

In general, $\Perf(X/G)$ is generated by $\Rep(G)$ as an idempotent complete stable $\infty$-category. Indeed, given any perfect complex $A\in \Perf(X/G)$, we can look at the largest $n$ for which the cohomology sheaf $H^n(A)$ is nonzero; after shift, $n=0$. Pick $V\in \Rep(G)$ with a map $V\to H^0(A)$ such that $V\otimes_L \mathcal O_{X/G}\to H^0(A)$ is surjective. By semisimplicity of $\Rep(G)$, we can lift $V\to H^0(A)$ to $V\to A$, and then pass to the cone of $V\otimes_L \mathcal O_{X/G}\to A$ to reduce the projective amplitude until $A$ is a vector bundle. In that case the homotopy fibre $B$ of $V\otimes_L \mathcal O_{X/G}\to A$ is again a vector bundle, and the map $V\otimes_L \mathcal O_{X/G}\to A$ splits, as the obstruction is $H^1(X/G,A^\vee\otimes_{\mathcal O_{X/G}} B)$, which vanishes by semisimplicity of $\Rep(G)$.

This already proves essential surjectivity. For fully faithfulness, it suffices by passage to internal Hom's to show that for all $A_{i_0}\in \Perf(X_{i_0}/G_{i_0})$ (for some chosen $i_0$) with pullbacks $A_i\in \Perf(X_i/G_i)$ for $i\to i_0$ and $A\in \Perf(X/G)$, the map $\varinjlim_i R\Gamma(X_i/G_i,A_i)\to R\Gamma(X/G,A)$ is an isomorphism. By semisimplicity of $\Rep(G_i)$ and $\Rep(G)$, it suffices to see that $\varinjlim_i R\Gamma(X_i,A_i)\to R\Gamma(X,A)$ is an isomorphism, which is clear by affineness.

This handles the case of filtered colimits. For the more precise claim, it is also easy to see that it commutes with disjoint unions. It is now enough to handle pushouts, so consider a diagram $S_1\leftarrow S_0\to S_2$ of anima over $\ast/Q$, with pushout $S$. We can assume that the maps $S_0\to S_1$ and $S_0\to S_2$ are surjective, as otherwise we can use compatibility with disjoint unions (replacing $S_2$ by the disjoint union of the image of $S_0$ and its complement). Then choose a surjection $S'\to S_0\times_{\ast/Q}\ast$, which induces similar surjections in the other cases. Thus, we get affine derived $L$-schemes $X_1\to X_0\leftarrow X_2$ with actions by $G=H^{S'}$, and $X=X_1\times_{X_0} X_2$, and we want to see that the functor
\[
\Perf(X_1/G)\otimes_{\Perf(X_0/G)} \Perf(X_2/G)\to \Perf(X_1\times_X X_2/G)
\]
is an equivalence. On the level of $\Ind$-categories, $\Ind\Perf(X_i/G)$ is the $\infty$-category of $\O(X_i)$-modules in $\Ind\Perf(\ast/G)$: This is a consequence of Barr--Beck--Lurie \cite[Theorem 4.7.4.5]{LurieHA} and the fact observed above that $\Perf(\ast/G)$ generates $\Perf(X_i/G)$, so that the forgetful functor $\Ind\Perf(X_i/G)\to \Ind\Perf(\ast/G)$ is conservative. It follows that the tensor product is the $\infty$-category of $\O(X_1)\otimes_{\O(X_0)} \O(X_2)$-modules in $\Ind\Perf(\ast/G)$, the tensor product $\O(X_1)\otimes_{\O(X_0)} \O(X_2)$ taken in the symmetric monoidal stable $\infty$-category $\Ind\Perf(\ast/G)$. The map $\O(X_1)\otimes_{\O(X_0)}\O(X_2)\to \O(X)$ is an isomorphism in $\Ind\Perf(\ast/G)$: This can be checked after the forgetful functor $\Ind\Perf(\ast/G)\to \mathcal D(L)$ as it is conservative (using that $G$ is pro-reductive, hence $\Rep(G)$ is semisimple), and then it amounts to $X=X_1\times_{X_0} X_2$.
\end{proof}

In particular, we get the following corollary.

\begin{corollary}\label{cor:spectralactionchar0} Let $L$ be a field over $\mathbb Q_\ell(\sqrt{q})$. There is a natural compactly supported $L$-linear action of $\Perf(Z^1(W_E,\hat{G})_L/\hat{G})$ on $\mathcal D_\lis(\Bun_G,L)^\omega$, uniquely characterized by the requirement that by restricting along the $\Rep_L(Q^I)$-linear maps
\[
\Rep_L((\hat{G}\rtimes Q)^I)\to \Perf(Z^1(W_E,\hat{G})_L/\hat{G})^{BW_E^I}
\]
it induces the Hecke action, which gives functorially in the finite set $I$ exact $\Rep_L(Q^I)$-linear functors
\[
\Rep_L((\hat{G}\rtimes Q)^I)\to \mathrm{End}_L(\mathcal D_\lis(\Bun_G,L)^\omega)^{BW_E^I}.
\]
\end{corollary}

\begin{proof} We can reduce to the subcategories $\mathcal D_\lis^P(\Bun_G,L)^\omega\subset \mathcal D_\lis(\Bun_G,L)$ for open subgroups $P$ of the wild inertia of $W_E$, acting trivially on $\hat{G}$. Then we can replace $W_E$ by $W_E/P$ throughout. In that case, restricting the given Hecke action to $W\subset W_E/P$, Theorem~\ref{thm:spectralactionchar0} gives an action of $\Perf(Z^1(W,\hat{G})_L/\hat{G})$, and $Z^1(W,\hat{G})=Z^1(W_E/P,\hat{G})$, so we get the desired action of $\Perf(Z^1(W_E/P,\hat{G})_L/\hat{G})$.
\end{proof}

With this action, we can formulate the main conjecture, ``the categorical form of the geometric Langlands conjecture on the Fargues--Fontaine curve''. Recall that for a quasisplit reductive group $G$ over $E$, Whittaker data consist of a choice of a Borel $B\subset G$ with unipotent radical $U\subset B$, together with a generic character $\psi: U(E)\to \overline{\mathbb Q}_\ell^\times$. As usual, we also fix $\sqrt{q}\in \overline{\mathbb Q}_\ell$. Let
\[
\mathcal W_\psi\in \mathcal D_\lis(\Bun_G,\overline{\mathbb Q}_\ell)
\]
be the Whittaker sheaf, which is the sheaf concentrated on $\Bun_G^1\subset \Bun_G$ corresponding to the representation $c\text-\mathrm{Ind}_{U(E)}^{G(E)} \psi$ of $G(E)$.

\begin{conj}\label{conj:mainconjecturerational} Consider the functor
\[
\mathcal D_{\mathrm{qcoh}}(Z^1(W_E,\hat{G})_{\overline{\mathbb Q}_\ell}/\hat{G}) = \Ind\mathrm{Perf}^{\mathrm{qc}}(Z^1(W_E,\hat{G})_{\overline{\mathbb Q}_\ell}/\hat{G})\to \mathcal D_\lis(\Bun_G,\overline{\mathbb Q}_\ell): M \mapsto \mathrm{Act}_M(\mathcal W_\psi)
\]
given as the colimit-preserving extension of the spectral action $\mathrm{Act}$ on the Whittaker sheaf. The corresponding right adjoint functor is fully faithful when restricted to the compact objects, and induces an equivalence of ($\Perf(Z^1(W_E,\hat{G})_{\overline{\mathbb Q}_\ell}/\hat{G})$-linear small stable) $\infty$-categories
\[
\mathcal D_\lis(\Bun_G,\overline{\mathbb Q}_\ell)^\omega\cong \mathcal D^{b,\mathrm{qc}}_{\mathrm{coh}}(Z^1(W_E,\hat{G})_{\overline{\mathbb Q}_\ell}/\hat{G}).
\]
\end{conj}

To be precise, $\mathcal D^{b,\mathrm{qc}}_{\mathrm{coh}}(Z^1(W_E,\hat{G})_{\overline{\mathbb Q}_\ell}/\hat{G})$ refers here to the $\infty$-category of those bounded complexes with coherent cohomology that also have quasicompact support, i.e.~only live on finitely many connected components.

\begin{remark}
There is an orthogonal decomposition $$\mathcal D_{\lis} (\Bun_G,\Qlb)^\omega = \bigoplus_{\alpha\in \pi_1 (G)_\Gamma} \mathcal D_{\lis} ( \Bun_G^{c_1=\alpha},\Qlb)^\omega$$ given by the connected components of $\Bun_G$. There is a morphism $Z(\hat{G})^\Gamma \to \mathrm{Aut}(\mathrm{Id}_{Z^1(W_E,\hat{G})_{\overline{\mathbb Q}_\ell}/\hat{G}})$, as $Z(\hat{G})^\Gamma\subset \hat{G}$ acts trivially on $Z^1(W_E,\hat{G})$. There is an associated ``eigenspace'' decomposition $$\mathcal D^{b,\mathrm{qc}}_{\mathrm{coh}}( Z^1(W_E,\hat{G})_{\overline{\mathbb Q}_\ell}/\hat{G})= \bigoplus_{\chi \in X^*(Z(\hat{G})^\Gamma)} \mathcal D^{b,\mathrm{qc}}_{\mathrm{coh}}( Z^1(W_E,\hat{G})_{\overline{\mathbb Q}_\ell}/\hat{G})_\chi.$$ Compatibility with the spectral action implies that via the identification $\pi_1(G)_\Gamma= X^*(Z(\hat{G})^\Gamma)$ those two decompositions should match.
\end{remark}

Another way to phrase the preceding conjecture is to say that, noting $\ast$ the spectral action, the ``non-abelian Fourier transform''\footnote{No precise meaning is implied by these words.}
\begin{align*}
\Perf^{\mathrm{qc}}(Z^1(W_E,\hat{G})_{\overline{\mathbb Q}_\ell}/\hat{G}) & \lto \mathcal D_\lis(\Bun_G,\overline{\mathbb{Q}}_\ell)
 \\
M & \longmapsto M\ast \mathcal W_\psi
\end{align*}
is fully faithful and extends to an equivalence of $\overline{\mathbb Q}_\ell$-linear small stable $\infty$-categories
$$
\mathcal D^{b,\mathrm{qc}}_{\mathrm{coh}}(Z^1(W_E,\hat{G})_{\overline{\mathbb Q}_\ell}/\hat{G})
 \xrightarrow{\sim} \mathcal D_\lis(\Bun_G,\overline{\mathbb{Q}}_\ell)^\omega.
$$

\begin{example}\label{rem:stable berstein center egal}
Fully faithfulness in the categorical conjecture, applied to the structure sheaf, implies that 
$$
\mathcal{Z}^{\mathrm{spec}} (G,\Qlb) \iso \End (c\text-\mathrm{Ind}_{U(E)}^{G(E)} \psi).
$$
\end{example}

\begin{example}[Kernel of functoriality]
Conjecture \ref{conj:mainconjecturerational} implies the existence of a kernel of functoriality for the local Langlands correspondence in the following way. Let $$f:{}^L H \to {}^L G$$ be an $L$-morphism between the $L$-groups of two quasi-split reductive groups $H$ and $G$ over $E$. This defines a morphism of stacks 
$$
Z^1 (W_E, \hat{H} )_{\Qlb} / \hat{H} \lto Z^1 (W_E,\hat{G})_{\Qlb}/\hat{G},
$$
and pushforward along this map induces a functor
\[
\Ind\mathcal D^{b,\mathrm{qc}}_{\mathrm{coh}}(Z^1 (W_E, \hat{H} )_{\Qlb} / \hat{H})\to \Ind\mathcal D^{b,\mathrm{qc}}_{\mathrm{coh}}(Z^1 (W_E,\hat{G})_{\Qlb}/\hat{G}).
\]
(There may be slightly different ways of handling the singularities here. One way to argue is to observe that pushforward is naturally a functor
\[
\mathcal D_{\mathrm{qcoh}}^{\geq 0}(Z^1 (W_E, \hat{H} )_{\Qlb} / \hat{H})\to \mathcal D_{\mathrm{qcoh}}^{\geq 0}(Z^1 (W_E, \hat{G} )_{\Qlb} / \hat{G}),
\]
and $\mathcal D_{\mathrm{qcoh}}^{\geq 0} = \mathrm{Ind}\mathcal D^{b,\mathrm{qc},\geq 0}_{\mathrm{coh}}$, and then extend by shifts.) The categorical equivalence then leads to a canonical functor
\[
\mathcal D_\lis(\Bun_H,\Qlb)\to \mathcal D_\lis(\Bun_G,\Qlb).
\]
By the self-duality of $\mathcal D_\lis$ coming from Bernstein--Zelevinsky duality, and Proposition~\ref{prop:lisseBunGproduct}, any such functor is given by a kernel
\[
A_f\in \mathcal D_\lis(\Bun_H\times \Bun_G,\Qlb).
\]
One could, in fact, identify the image of $A_f$ under the categorical equivalence for $H\times G$; up to minor twists, it should be given by the structure sheaf of the graph of $Z^1 (W_E, \hat{H} )_{\Qlb} / \hat{H} \lto Z^1 (W_E,\hat{G})_{\Qlb}/\hat{G}$. It would be very interesting if some examples of such kernels $A_f$ can be constructed explicitly.

Since $D(H(E),\Qlb)$, resp. $D(G(E),\Qlb)$, are direct factors of $D_\lis(\Bun_H,\Qlb)$, resp.~$D_\lis(\Bun_G,\Qlb)$, this should give rise to the ``classical'' Langlands functoriality $D(H(E),\Qlb)\to D(G(E),\Qlb)$.
\end{example}

\begin{remark}
Above, we assumed $G$ and $H$ are quasisplit, in order to define the Whittaker sheaf. To some extent, this is necessary, as the Jacquet--Langlands correspondence cannot be given by a completely canonical functor (in particular, one defined over $\mathbb Q_\ell$): In fact, as is well-known, any discrete series representation of $\GL_n(E)$ is defined over its field of moduli but this is not the case for smooth irreducible representations of $D^\times$ if $D$ is a division algebra over $E$.
\end{remark}

Let us now explain how {\it Fargues's original conjecture} fits into this context. Let $\varphi:W_E\to \hat{G}(\Qlb)$ be a Langlands parameter. Consider the map $i: \Spec \Qlb\to Z^1(W_E,\hat{G})_{\Qlb}/\hat{G}$ corresponding to $\varphi$, and let
\[
\mathcal E_\phi = i_\ast \Qlb\in \mathcal D_{\mathrm{qcoh}}(Z^1(W_E,\hat{G})_{\Qlb}/\hat{G})=\Ind\Perf^{\mathrm{qc}}(Z^1(W_E,\hat{G})_{\Qlb}/\hat{G}).
\]
(We take the pushforward here in the sense of $\mathcal D_{\mathrm{qcoh}}$. One could a priori produce a more refined object of $\Ind\mathcal D^{b,\mathrm{qc}}_{\mathrm{coh}}$, but we do not consider this here.) Factoring the map $i$ via $[\Spec \Qlb/S_\varphi]$, one actually sees that $\mathcal E_\phi$ carries naturally an action of $S_\varphi$. Moreover, if one acts via tensoring with a representation $V$ of $\hat{G}\rtimes Q$, then by the projection formula the sheaf $\mathcal E_\phi$ gets taken to itself, tensored with the $W_E$-representation $V\circ \phi$. Using the spectral action, we find an $S_{\varphi}$-equivariant ``automorphic complex''
\[
\mathrm{Aut}_{\varphi}=\mathcal E_\phi\ast \mathcal W_\psi\in \mathcal D_\lis(\Bun_G, \Qlb).
\]
It already follows that $\mathrm{Aut}_\phi\in \mathcal D_\lis(\Bun_G,\Qlb)$ {\it is} a Hecke eigensheaf, with eigenvalue $\phi$, so the spectral action produces Hecke eigensheaves. Except, it is not clear whether $\mathrm{Aut}_\phi\neq 0$. Under the fully faithfulness part of the categorical conjecture, one sees that it must be nonzero, and moreover have some of the properties stated in \cite{FarguesGeom}, in particular regarding the relation to $L$-packets. The particular case of elliptic parameters is further spelled out in the next section.

\section{Elliptic parameters}

Let us make explicit what the spectral action, and Conjecture~\ref{conj:mainconjecturerational}, entails in the case of elliptic parameters. As coefficients, we take $L=\overline{\mathbb Q}_\ell$ for simplicity.

\begin{definition} An $L$-parameter $\phi: W_E\to \hat{G}(\overline{\mathbb Q}_\ell)$ is elliptic if it is semisimple and the centralizer $S_\phi\subset \hat{G}_{\overline{\mathbb Q}_\ell}$ has the property that $S_\phi/Z(\hat{G})_{\overline{\mathbb Q}_\ell}^\Gamma$ is finite.
\end{definition}

By deformation theory,\footnote{One has $H^2(W_E,\mathrm{ad}\, \phi)=0$ using Tate duality, and the $H^0$ reduces to the Lie algebra of $Z(\hat{G})^\Gamma$. The $H^1$ must thus be of the same dimension and be accounted for by the unramified twists.} it follows that the unramified twists of $\phi$ define a connected component
\[
C_\phi\hookrightarrow [Z^1(W_E,\hat{G})_{\overline{\mathbb Q}_\ell}/\hat{G}].
\]
Thus, the spectral action (in fact, the excursion operators are enough for this, see the discussion around Theorem~\ref{thm:bernsteincenter}) implies that there is a corresponding direct summand
\[
\mathcal D_\lis^{C_\phi}(\Bun_G,\overline{\mathbb Q}_\ell)^\omega\subset \mathcal D_\lis(\Bun_G,\overline{\mathbb Q}_\ell)^\omega,
\]
explicitly given as those objects on which the excursion operator corresponding to the function that is $1$ on $C_\phi$ and $0$ elsewhere acts via the identity. For any Schur-irreducible $A\in \mathcal D_\lis^{C_\phi}(\Bun_G,\overline{\mathbb Q}_\ell)^\omega$, the excursion operators act via scalars on $A$, as determined by an unramified twist of $\phi$. In particular, they act in this way on $i^{b\ast} A$ for any $b\in B(G)$. By compatibility with parabolic induction, it follows that for any $A\in \mathcal D_\lis^{C_\phi}(\Bun_G,\overline{\mathbb Q}_\ell)^\omega$, the restriction $i^{b\ast} A$ is equal to $0$ if $b$ is not basic (if it was not zero, one could find an irreducible subquotient to which this argument applies). Thus,
\[
\mathcal D_\lis^{C_\phi}(\Bun_G,\overline{\mathbb Q}_\ell)^\omega\cong \bigoplus_{b\in B(G)_{\mathrm{basic}}} \mathcal D^{C_\phi}(G_b(E),\overline{\mathbb Q}_\ell)^\omega.
\]
Moreover, all $A\in \mathcal D^{C_\phi}(G_b(E),\overline{\mathbb Q}_\ell)^\omega$ must lie in only supercuspidal components of the Bernstein center, again by compatibility with parabolic induction. If $Z(\hat{G})^\Gamma$ is finite (equivalently, if the connected split center of $G$ is trivial), then $C_\phi=[\ast/S_\phi]$ is a point and it follows that all $A$ are finite direct sums of shifts supercuspidal representations of $G_b(E)$, and so
\[
\mathcal D_\lis^{C_\phi}(\Bun_G,\overline{\mathbb Q}_\ell)^\omega\cong \bigoplus_{b\in B(G)_{\mathrm{basic}}} \bigoplus_\pi \Perf(\overline{\mathbb Q}_\ell)\otimes \pi,
\]
where $\pi$ runs over supercuspidal $\overline{\mathbb Q}_\ell$-representations of $G_b(E)$ with $L$-parameter $\varphi_\pi=\phi$.

In general, acting on $\mathcal D_\lis^{C_\phi}(\Bun_G,\overline{\mathbb Q}_\ell)^\omega$, we have the direct summand
\[
\Perf(C_\phi)
\]
of
\[
\Perf([Z^1(W_E,\hat{G})_{\overline{\mathbb Q}_\ell}/\hat{G}]).
\]
If $Z(\hat{G})^\Gamma$ is finite, one has $C_\phi=[\ast/S_\phi]$, and hence we get an action of $\Rep(S_\phi)$ on $\mathcal D_\lis^{C_\phi}(\Bun_G,\overline{\mathbb Q}_\ell)^\omega$. In general, one can get a similar picture by fixing central characters; let us for simplicity only spell out the case when $Z(\hat{G})^\Gamma$ is finite, i.e.~the connected split center of $G$ is trivial.

If $\pi_b$ is a supercuspidal representation of some $G_b(E)$ with $\varphi_{\pi_b}=\phi$, and $W\in \Rep(S_\phi)$ then acting via $W$ on $\pi_b$ we get some object
\[
\mathrm{Act}_W(\pi_b)\in \bigoplus_{b'\in B(G)_{\mathrm{basic}}} \bigoplus_\pi \Perf(\overline{\mathbb Q}_\ell)\otimes \pi.
\]
Assume that $W|_{Z(\hat{G})^\Gamma}$ is isotypic, given by some character $\chi: Z(\hat{G})^\Gamma\to \overline{\mathbb Q}_\ell^\times$. As $Z(\hat{G})^\Gamma$ is the diagonalizable group with characters $\pi_1(G)_\Gamma$, it follows that we get an element $b_\chi\in \pi_1(G)_\Gamma = B(G)_{\mathrm{basic}}$. Then $\mathrm{Act}_W(\pi_b)$ is concentrated on $b'=b+b_\chi$, and so
\[
\mathrm{Act}_W(\pi_b)\cong \bigoplus_{\pi_{b'}} V_{\pi_{b'}}\otimes \pi_{b'}
\]
for a certain multiplicity space $V_{\pi_{b'}}\in \Perf(\overline{\mathbb Q}_\ell)$, where $\pi_{b'}$ runs over supercuspidal representations of $G_{b'}(E)$, $b'=b+b_\chi$, with $L$-parameter $\varphi_{\pi_{b'}}=\phi$.

The conjectural description of $L$-packets \cite{KalethaLocalLanglands} then suggests the following conjecture, which is (up to the added $t$-exactness) the specialization of Conjecture~\ref{conj:mainconjecturerational} to the case of elliptic $L$-parameters. (If one projects the Whittaker sheaf to the $C_\phi$-component, a priori it could break into a direct sum of several supercuspidal representations; but then the functor would not have a chance of being an equivalence.)

\begin{conj}\label{conj:ellipticLpackets} Assume that $G$ is quasisplit, with a fixed Whittaker datum, and that the connected split center of $G$ is trivial. Then there is a unique generic supercuspidal representation $\pi$ of $G(E)$ with $L$-parameter $\varphi_\pi=\phi$, and the functor
\[
\Perf([\Spec \overline{\mathbb Q}_\ell/S_\phi])\to \mathcal D_\lis^{C_\phi}(\Bun_G,\overline{\mathbb Q}_\ell)^\omega: W\mapsto \mathrm{Act}_W(\pi)
\]
is an equivalence. In particular, the set of irreducible supercuspidal representations of some $G_b(E)$ with $L$-parameter $\phi$ is in bijection with the set of irreducible representations of $S_\phi$.

Moreover, the equivalence is $t$-exact for the standard $t$-structures on source and target. 
\end{conj}

Thus, the conjecture gives an explicit parametrization of $L$-packets.

Let us explain what the compatibility of the spectral action with Hecke operators entails in this case. Given $V\in \Rep(\hat{G}\rtimes Q)$, the restriction of $V$ to $S_\phi$ admits a commuting $W_E$-action given by $\phi$. This defines a functor
\[
\Rep(\hat{G}\rtimes Q)\to \Rep(S_\phi)^{BW_E}.
\]
Now the diagram of monoidal functors
\[\xymatrix{
\Rep(\hat{G}\rtimes Q)\ar[r]\ar[d] & \End_{\overline{\mathbb Q}_\ell}(\mathcal D_\lis^{C_\phi}(\Bun_G,\overline{\mathbb Q}_\ell))^{BW_E} \\
\Rep(S_\phi)^{BW_E}\ar[ur]
}\]
commutes; this follows from the compatibility of the spectral action with the Hecke action.

Concretely, given $\pi$ as above and $V\in \Rep(\hat{G}\rtimes Q)$, decompose the image of $V$ in $\Rep(S_\phi)^{BW_E}$ as a direct sum $\bigoplus_{i\in I} W_i\otimes \sigma_i$ where $W_i\in \Rep(S_\phi)$ is irreducible and $\sigma_i$ is some continuous representation of $W_E$ on a finite-dimensional $\overline{\mathbb Q}_\ell$-space. Then
\[
T_V(\pi)\cong \bigoplus_{i\in I} \mathrm{Act}_{W_i}(\pi)\otimes \sigma_i.
\]

Recall that $T_V(\pi)$ can be calculated concretely through the cohomology of local Shimura varieties, or in general moduli spaces of local shtukas. Noting that the functor $\mathrm{Act}_{W_i}$ is realizing a form of the Jacquet--Langlands correspondence relating different inner forms, the formula above is essentially the conjecture of Kottwitz \cite[Conjecture 7.3]{RapoportViehmann}. In fact, assuming Conjecture~\ref{conj:ellipticLpackets}, it is an easy exercise to deduce \cite[Conjecture 7.3]{RapoportViehmann}, assuming that the parametrization of the Conjecture~\ref{conj:ellipticLpackets} agrees with the parametrization implicit in \cite[Conjecture 7.3]{RapoportViehmann}.

\section{Integral coefficients}\label{sec:integralcoeff}

We want to construct the spectral action with integral coefficients. Unfortunately, the naive analogue of Theorem~\ref{thm:spectralactionchar0} is not true, the problem being that the analogue of Lemma~\ref{lem:perfectcomplexeslimit} fails. However, the rest of the argument still works, and gives the following result.

Consider a split reductive group $H$ over a discrete valuation ring $R$ with an action of a finite group $Q$. Let $S$ be any anima over $\ast/Q$. As before, we can define a derived stack $\Map_{\ast/Q}(S,\ast/(H\rtimes Q))$ over $R$, whose values in an animated $R$-algebra $A$ are the maps of anima $S\to \ast/(H\rtimes Q)(A)$ over $\ast/Q$. In general, the functor $S\mapsto \Perf(\Map_{\ast/Q}(S,\ast/(H\rtimes Q)))$ does not commute with sifted colimits in $S$.

However, we can consider the best approximation to it that does commute with sifted colimits. Note that the $\infty$-category of anima over $\ast/Q$ is the animation of the category of sets equipped with a $Q$-torsor; it is freely generated under sifted colimits by the category of finite sets equipped with a $Q$-torsor. Thus, the sifted-colimit approximation to $S\mapsto \Perf(\Map_{\ast/Q}(S,\ast/(H\rtimes Q)))$ is the animation of its restriction to finite sets with $Q$-torsors; we denote it by
\[
S\mapsto \Perf(\Map_{\ast/Q}^{\Sigma}(S,\ast/(H\rtimes Q))),
\]
with the idea in mind that it is like the $\infty$-category of perfect complexes on some (nonexistent) derived stack $\Map_{\ast/Q}^{\Sigma}(S,\ast/(H\rtimes Q))$, gotten as a (co-)sifted limit approximation to $\Map_{\ast/Q}(S,\ast/(H\rtimes Q))$. The symbol $\Sigma$ here is in reference to the notation used in \cite[Section 5.5.8]{LurieHTT} in relation to sifted colimits. Thus $\Perf(\Map_{\ast/Q}^{\Sigma}(S,\ast/(H\rtimes Q)))$ is an $R$-linear idempotent-complete small stable $\infty$-category, mapping to $\Perf(\Map_{\ast/Q}(S,\ast/(H\rtimes Q)))$.

\begin{proposition}\label{prop:spectralactiongeneral} Let $\mathcal C$ be an $R$-linear idempotent-complete small stable  $\infty$-category. Giving, functorially in finite sets $I$, an exact $\Rep_R(Q^I)$-linear monoidal functor
\[
\Rep_R((H\rtimes Q)^I)\to \End_R(\mathcal C)^{S^I}
\]
is equivalent to giving an $R$-linear action of $\Perf(\Map_{\ast/Q}^{\Sigma}(S,\ast/(H\rtimes Q)))$ on $\mathcal C$. Here, given such an action of $\Perf(\Map_{\ast/Q}^{\Sigma}(S,\ast/(H\rtimes Q)))$, one gets exact $\Rep_R(Q^I)$-linear monoidal functors
\[
\Rep_R((H\rtimes Q)^I)\to \End_R(\mathcal C)^{S^I}
\]
by composing the exact monoidal functor $\Perf(\Map_{\ast/Q}^{\Sigma}(S,\ast/(H\rtimes Q)))\to \End_R(\mathcal C)$ with the natural exact $\Rep_R(Q^I)$-linear symmetric monoidal functor
\[
\Rep_R((H\rtimes Q)^I)\to \Perf(\Map_{\ast/Q}^{\Sigma}(S,\ast/(H\rtimes Q)))^{S^I}
\]
given by $I$-fold tensor product of the exact $\Rep_R(Q)$-linear symmetric monoidal functor
\[
\Rep_R(H\rtimes Q)\to \Perf(\Map_{\ast/Q}^{\Sigma}(S,\ast/(H\rtimes Q)))^S
\]
assigning to each $s\in S$ pullback along evaluation at $s$, $\Map_{\ast/Q}(S,\ast/(H\rtimes Q))\to \ast/(H\rtimes Q)$; more precisely, it is defined in this way if $S$ is a finite set, and in general by animation.
\end{proposition}

\begin{proof} This follows from the proof of Theorem~\ref{thm:spectralactionchar0}.
\end{proof}

To make use of Proposition~\ref{prop:spectralactiongeneral}, we need to find sufficiently many situations in which the functor
\[
\Perf(\Map_{\ast/Q}^{\Sigma}(S,\ast/(H\rtimes Q)))\to \Perf(\Map_{\ast/Q}(S,\ast/(H\rtimes Q)))
\]
is an equivalence, and specifically we need to prove this for $\Map_{\ast/Q}(\ast/W,\ast/(\hat{G}\rtimes Q))=Z^1(W,\hat{G})/\hat{G}$.

First, we have the following result.

\begin{proposition} The functor $S\mapsto \Perf(\Map_{\ast/Q}^\Sigma(S,\ast/(H\rtimes Q)))$ from anima over $\ast/Q$ to symmetric monoidal idempotent-complete stable $R$-linear $\infty$-categories commutes with all colimits.
\end{proposition}

\begin{proof} As the functor commutes with sifted colimits by definition, it suffices to show that when restricted to finite sets $S$ equipped with $Q$-torsors, it commutes with disjoint unions. But for such $S$, the map $S\to \ast/Q$ can be factored over a point, and then $\Map_{\ast/Q}(S,\ast/(H\rtimes Q))=\ast/H^S$. Thus, one has to see that for two finite sets $S_1$, $S_2$, the functor
\[
\Perf(\ast/H^{S_1})\otimes_{\Perf(R)} \Perf(\ast/H^{S_2})\to \Perf(\ast/H^{S_1\sqcup S_2})
\]
is an equivalence. But this follows easily from highest weight theory, which for any split reductive group $H$ filters $\Perf(\ast/H)$ in terms of copies of $\Perf(R)$ enumerated by highest weights.
\end{proof}

\begin{proposition} Assume that $S=\ast/F_n$ is the classifying space of a free group. Then the functor
\[
\Perf(\Map_{\ast/Q}^{\Sigma}(S,\ast/(H\rtimes Q)))\to \Perf(\Map_{\ast/Q}(S,\ast/(H\rtimes Q)))
\]
is fully faithful, and the essential image is the idempotent-complete stable $\infty$-subcategory generated by the image of $\Rep_R(H)$.
\end{proposition}

\begin{proof} Represent $\ast/F_n\to \ast/Q$ by a map $F_n\to Q$, and let $\sigma_1,\ldots,\sigma_n\in Q$ be the images of the generators. Then $\Map_{\ast/Q}(S,\ast/(H\rtimes Q))$ can be identified with $[H^n/H]$, where $H$ acts on $H^n$ via the $(\sigma_1,\ldots,\sigma_n)$-twisted diagonal conjugation action. We claim that
\[
\Perf(\Map_{\ast/Q}^\Sigma(S,\ast/(H\rtimes Q)))
\]
is the $\infty$-category of compact objects in the $\infty$-category of modules over $\O(H^n)$ in $\Ind\Perf(\ast/H)$; in fact, this is equivalent to the claim, as by Barr--Beck--Lurie \cite[Theorem 4.7.4.5]{LurieHA} this gives a description of the full $\infty$-subcategory of $\Perf([H^n/H])$ generated by $\Perf(\ast/H)$.

As $\O(H^n) = \O(H)\otimes\ldots\otimes \O(H)$ in $\Ind\Perf(\ast/H)$, one reduces to the case $n=1$. In that case $S=\ast/F_1$ is a circle, which we can present as a pushout of $\ast\sqcup\ast\rightrightarrows \ast$. Thus, we have to compute
\[
\Perf(\ast/H)\otimes_{\Perf(\ast/H^2)}\Perf(\ast/H)
\]
where the two implicit maps $H\to H^2$ are given by the diagonal and the $\sigma_1$-twisted diagonal, respectively. As the pullback functors $\Perf(\ast/H^2)\to \Perf(\ast/H)$ generate the image, we can write $\Perf(\ast/H)=\Perf(H/H^2)$ as the compact objects in the $\infty$-category of $\O(H)$-modules in $\Ind\Perf(\ast/H^2)$. Similarly, the expected answer $\Perf([H/H])=\Perf([H\times H/H^2])$ is given by the compact objects in the $\infty$-category of modules over $\O(H\times H)=\O(H)\otimes \O(H)$ in $\Ind\Perf(\ast/H^2)$, thus implying the result.
\end{proof}

\begin{proposition} Let $S=\ast/\Gamma$, where $\Gamma$ is any discrete group, and lift the map $S\to \ast/Q$ to a map $\Gamma\to Q$. One can write
\[
\ast/\Gamma=\colim_{(n,F_n\to \Gamma)} \ast/F_n
\]
as a sifted colimit (in anima). Then $\Perf(\Map_{\ast/Q}^{\Sigma}(S,\ast/(H\rtimes Q)))$ is the $\infty$-category of compact objects in the $\infty$-category of modules over
\[
\colim_{(n,F_n\to \Gamma)} \O(H^n)
\]
in $\Ind\Perf(\ast/H)$, where $\O(H^n)$ is equipped with the twisted (via the map $F_n\to \Gamma\to Q$) diagonal conjugation of $H$.
\end{proposition}

\begin{proof} As $\mathbb E_1$-groups in anima are equivalent to animated groups, with compact projective generators the free groups $F_n$, it follows that $\ast/\Gamma$ is the sifted colimit $\colim_{(n,F_n\to \Gamma)} \ast/F_n$. Now the result follows from the previous proposition (and its proof), together with the commutation with sifted colimits.
\end{proof}

Combining this with Theorem~\ref{thm:precisecoliminIndPerf}, we have finished the proof of Theorem~\ref{thm:constructionspectralaction2}. In particular, this gives the spectral action on $\mathcal D_\lis(\Bun_G)$.

Let us end by stating again the main conjecture with integral coefficients. The formulation may have to be adapted at very small primes (i.e., bad primes) as then the nilpotent cone implicit in the formulation of nilpotent singular support of Section~\ref{sec:singularsupport} may not be well-behaved.

\begin{conj} Assume that $G$ is quasisplit and choose Whittaker data consisting of a Borel $B\subset G$ and generic character $\psi: U(E)\to \mathcal O_L^\times$ of the unipotent radical $U\subset B$, where $L/\mathbb Q_\ell$ is some algebraic extension; also fix $\sqrt{q}\in \mathcal O_L$. Let $n$ be the order of $\pi_0 Z(G)$ and let $\Lambda=\mathcal O_L[\tfrac 1n]$. Let
\[
\mathcal W_\psi\in \mathcal D_\lis(\Bun_G,\Lambda)
\]
be the Whittaker sheaf, which is the sheaf concentrated on $\Bun_G^1$ corresponding to the Whittaker representation $c\text-\mathrm{Ind}_{U(E)}^{G(E)} \psi$, and let
\[
\Ind\Perf^{\mathrm{qc}}(Z^1(W_E,\hat{G})_\Lambda/\hat{G})\to \mathcal D_\lis(\Bun_G,\Lambda): M\mapsto \mathrm{Act}_M(\mathcal W_\psi)
\]
be defined as the colimit-preserving extension of the spectral action on $\mathcal W_\psi$. Then the corresponding right adjoint functor is fully faithful when restricted to the compact objects, and induces an equivalence of ($\Perf(Z^1(W_E,\hat{G})_\Lambda/\hat{G})$-linear small stable) $\infty$-categories
\[
\mathcal D(\Bun_G,\Lambda)^\omega\cong \mathcal D^{b,\mathrm{qc}}_{\mathrm{coh},\mathrm{Nilp}}(Z^1(W_E,\hat{G})_\Lambda/\hat{G}).
\]
\end{conj}

\bibliographystyle{amsalpha}
\bibliography{Geometrization1}

\end{document}